\documentclass{article}
\usepackage{amsmath}
\usepackage{amssymb}
\usepackage{amsthm}
\usepackage{cite}
\usepackage{euscript}
\usepackage[arrow,curve,matrix,arc,2cell]{xy}
\UseAllTwocells
\usepackage{makeidx}
\usepackage[refpage,intoc]{nomencl}
\makenomenclature
\makeindex
\DeclareFontFamily{U}{rsfs}{} \DeclareFontShape{U}{rsfs}{n}{it}{<->
rsfs10}{} \DeclareSymbolFont{mscr}{U}{rsfs}{n}{it}
\DeclareSymbolFontAlphabet{\scr}{mscr}
\def\mathscr{\scr}
\begin{document}
\renewcommand{\nomname}{Glossary of Notation}
\renewcommand{\pagedeclaration}[1]{, #1}
\setlength{\nomitemsep}{0pt}
\def\G[#1]#2#3{\nomenclature[#1]{#2}{#3}}
\def\I#1{\index{#1}}
\def\e#1\e{\begin{equation}#1\end{equation}}
\def\ea#1\ea{\begin{align}#1\end{align}}
\def\eq#1{{\rm(\ref{#1})}}
\theoremstyle{plain}
\newtheorem{thm}{Theorem}[section]
\newtheorem{lem}[thm]{Lemma}
\newtheorem{prop}[thm]{Proposition}
\newtheorem{cor}[thm]{Corollary}
\theoremstyle{definition}
\newtheorem{dfn}[thm]{Definition}
\newtheorem{ex}[thm]{Example}
\newtheorem{rem}[thm]{Remark}
\numberwithin{figure}{section}
\numberwithin{equation}{section}
\def\dim{\mathop{\rm dim}\nolimits}
\def\vdim{\mathop{\rm vdim}\nolimits}
\def\Rep{\mathop{\rm Rep}\nolimits}
\def\Tot{\mathop{\rm Tot}}
\def\Totc{\mathop{\rm Tot^c}}
\def\Ker{\mathop{\rm Ker}}
\def\Coker{\mathop{\rm Coker}}
\def\GL{\mathop{\rm GL}}
\def\Stab{\mathop{\rm Stab}\nolimits}
\def\supp{\mathop{\rm supp}}
\def\rank{\mathop{\rm rank}\nolimits}
\def\Hom{\mathop{\rm Hom}\nolimits}
\def\Ho{\mathop{\rm Ho}\nolimits}
\def\Ext{\mathop{\rm Ext}\nolimits}
\def\Ad{\mathop{\rm Ad}\nolimits}
\def\Aut{\mathop{\rm Aut}}
\def\coh{\mathop{\rm coh}}
\def\qcoh{\mathop{\rm qcoh}}
\def\vect{\mathop{\rm vect}}
\def\vqcoh{\mathop{\rm vqcoh}}
\def\vvect{\mathop{\rm vvect}}
\def\id{{\mathop{\rm id}\nolimits}}
\def\top{{\kern.05em\rm top}}
\def\Iso{{\rm Iso}}
\def\inc{{\rm inc}}
\def\tr{{\rm tr}}
\def\nt{{\rm nt}}
\def\ev{{\rm ev}}
\def\Obj{{\rm Obj}}
\def\Spec{\mathop{\rm Spec}\nolimits}
\def\MSpec{\mathop{\rm MSpec}\nolimits}
\def\Top{{\mathop{\bf Top}}}
\def\CRings{{\mathop{\bf C^{\bs\iy}Rings}}}
\def\CRS{{\mathop{\bf C^{\bs\iy}RS}}}
\def\CSch{{\mathop{\bf C^{\bs\iy}Sch}}}
\def\CSchlf{{\mathop{\bf C^{\bs\iy}Sch^{lf}}}}
\def\CSchlfssc{{\mathop{\bf C^{\bs\iy}Sch^{lf}_{ssc}}}}
\def\hCSchlfssc{{\mathop{\bf \hat{C}^{\bs\iy}Sch^{lf}_{ssc}}}}
\def\bCSch{{\mathop{\bf{\bar C}^{\bs\iy}Sch}}}
\def\CSta{{\mathop{\bf C^{\bs\iy}Sta}}}
\def\DMCSta{{\mathop{\bf DMC^{\bs\iy}Sta}}}
\def\DMCStalf{{\mathop{\bf DMC^{\bs\iy}Sta^{lf}}}}
\def\DMCStalfssc{{\mathop{\bf DMC^{\bs\iy}Sta^{lf}_{ssc}}}}
\def\hDMCStalfssc{{\mathop{\bf D\hat{M}C^{\bs\iy}Sta^{lf}_{ssc}}}}
\def\Sta{{\mathop{\bf Sta}\nolimits}}
\def\Man{{\mathop{\bf Man}}}
\def\Manb{{\mathop{\bf Man^b}}}
\def\Manc{{\mathop{\bf Man^c}}}
\def\cManc{{\mathop{\bf\check{M}an^c}}}
\def\dMan{{\mathop{\bf dMan}}}
\def\bdMan{{\mathop{\bf d\bar{M}an}}}
\def\hdMan{{\mathop{\bf d\hat{M}an}}}
\def\dManb{{\mathop{\bf dMan^b}}}
\def\dManc{{\mathop{\bf dMan^c}}}
\def\dtManc{{\mathop{\bf d\widetilde{Ma}n{}^c}}}
\def\dcManc{{\mathop{\bf d\check{M}an{}^c}}}
\def\hdManc{{\mathop{\bf d\hat{M}an{}^c}}}
\def\dOrb{{\mathop{\bf dOrb}}}
\def\bdOrb{{\mathop{\bf d\bar{O}rb}}}
\def\dOrbb{{\mathop{\bf dOrb^b}}}
\def\dOrbc{{\mathop{\bf dOrb^c}}}
\def\dcOrbc{{\mathop{\bf d\check{O}rb{}^c}}}
\def\dSpa{{\mathop{\bf dSpa}}}
\def\dSta{{\mathop{\bf dSta}}}
\def\dSpab{{\mathop{\bf dSpa^b}}}
\def\dStab{{\mathop{\bf dSta^b}}}
\def\dSpac{{\mathop{\bf dSpa^c}}}
\def\dStac{{\mathop{\bf dSta^c}}}
\def\hdSpa{{\mathop{\bf d\hat{S}pa}}}
\def\bdSpa{{\mathop{\bf d\bar{S}pa}}}
\def\bdSta{{\mathop{\bf d\bar{S}ta}}}
\def\hdSpac{{\mathop{\bf d\hat{S}pa{}^c}}}
\def\hMan{{\mathop{\bf \hat{M}an}}}
\def\bMan{{\mathop{\bf \bar{M}an}}}
\def\bManb{{\mathop{\bf \bar{M}an^b}}}
\def\bManc{{\mathop{\bf \bar{M}an^c}}}
\def\Orb{{\mathop{\bf Orb}}}
\def\cOrb{{\mathop{\bf\check{O}rb}}}
\def\dotOrb{{\mathop{\bf \dot{O}rb}}}
\def\hOrb{{\mathop{\bf \hat{O}rb}}}
\def\bOrb{{\mathop{\bf \bar{O}rb}}}
\def\bOrbb{{\mathop{\bf \bar{O}rb^b}}}
\def\bOrbc{{\mathop{\bf \bar{O}rb^c}}}
\def\Orbb{{\mathop{\bf Orb^b}}}
\def\Orbc{{\mathop{\bf Orb^c}}}
\def\cOrbc{{\mathop{\bf\check{O}rb^c}}}
\def\DerMan{{\mathop{\bf DerMan}}}
\def\bo{{\rm bo}}
\def\hom{{\rm hom}}
\def\dbo{{\rm dbo}}
\def\dorb{{\rm dorb}}
\def\orb{{\rm orb}}
\def\sef{{\rm sef}}
\def\eff{{\rm eff}}
\def\deff{{\rm deff}}
\def\obo{{\rm obo}}
\def\dobo{{\rm dobo}}
\def\ul{\underline}
\def\bs{\boldsymbol}
\def\ge{\geqslant}
\def\le{\leqslant\nobreak}
\def\O{{\mathcal O}}
\def\ES{{\mathcal E}_S}
\def\FS{{\mathcal F}_S}
\def\IS{{\mathcal I}_S}
\def\OS{{{\mathcal O}_{\kern -.1em S}}}
\def\OSp{{{\mathcal O}^{\smash{\prime}}_{\kern -.1em S}}}
\def\EW{{\mathcal E}_W}
\def\FW{{\mathcal F}_W}
\def\OW{{{\mathcal O}_{\kern -.1em W}}}
\def\OWp{{{\mathcal O}^{\smash{\prime}}_{\kern -.1em W}}}
\def\EX{{\mathcal E}_X}
\def\FX{{\mathcal F}_X}
\def\IX{{\mathcal I}_X}
\def\OX{{{\mathcal O}_{\kern -.1em X}}}
\def\OXp{{{\mathcal O}^{\smash{\prime}}_{\kern -.1em X}}}
\def\EY{{\mathcal E}_Y}
\def\FY{{\mathcal F}_Y}
\def\IY{{\mathcal I}_Y}
\def\OY{{{\mathcal O}_{\kern -.1em Y}}}
\def\OYp{{{\mathcal O}^{\smash{\prime}}_{\kern -.1em Y}}}
\def\EZ{{\mathcal E}_Z}
\def\FZ{{\mathcal F}_Z}
\def\IZ{{\mathcal I}_Z}
\def\OZ{{{\mathcal O}_{\kern -.1em Z}}}
\def\OZp{{{\mathcal O}^{\smash{\prime}}_{\kern -.1em Z}}}
\def\EcS{{\mathcal E}_\cS}
\def\FcS{{\mathcal F}_\cS}
\def\OcS{{{\mathcal O}_{\kern -.1em\cS}}}
\def\OcSp{{{\mathcal O}^{\smash{\prime}}_{\kern -.1em\cS}}}
\def\OcV{{{\mathcal O}_{\kern -.1em\cV}}}
\def\OcW{{{\mathcal O}_{\kern -.1em\cW}}}
\def\EcX{{\mathcal E}_\cX}
\def\FcX{{\mathcal F}_\cX}
\def\IcX{{\mathcal I}_\cX}
\def\OcX{{{\mathcal O}_{\kern -.1em\cX}}}
\def\OcXp{{{\mathcal O}^{\smash{\prime}}_{\kern -.1em\cX}}}
\def\EcY{{\mathcal E}_\cY}
\def\FcY{{\mathcal F}_\cY}
\def\IcY{{\mathcal I}_\cY}
\def\OcY{{{\mathcal O}_{\kern -.1em\cY}}}
\def\OcYp{{{\mathcal O}^{\smash{\prime}}_{\kern -.1em\cY}}}
\def\EcZ{{\mathcal E}_\cZ}
\def\FcZ{{\mathcal F}_\cZ}
\def\IcZ{{\mathcal I}_\cZ}
\def\OcZ{{{\mathcal O}_{\kern -.1em\cZ}}}
\def\OcZp{{{\mathcal O}^{\smash{\prime}}_{\kern -.1em\cZ}}}
\def\im{\imath}
\def\jm{\jmath}
\def\bR{{\mathbin{\pmb{\mathbb R}}}}
\def\K{{\mathbin{\mathbb K}}}
\def\R{{\mathbin{\mathbb R}}}
\def\Z{{\mathbin{\mathbb Z}}}
\def\Q{{\mathbin{\mathbb Q}}}
\def\N{{\mathbin{\mathbb N}}}
\def\C{{\mathbin{\mathbb C}}}
\def\CP{{\mathbin{\mathbb{CP}}}}
\def\RP{{\mathbin{\mathbb{RP}}}}
\def\fC{{\mathbin{\mathfrak C}\kern.05em}}
\def\fD{{\mathbin{\mathfrak D}}}
\def\fE{{\mathbin{\mathfrak E}}}
\def\fF{{\mathbin{\mathfrak F}}}
\def\fn{{\mathbin{\mathfrak n}}}
\def\cA{{\mathbin{\cal A}}}
\def\cB{{\mathbin{\cal B}}}
\def\cC{{\mathbin{\cal C}}}
\def\cD{{\mathbin{\cal D}}}
\def\cE{{\mathbin{\cal E}}}
\def\cF{{\mathbin{\cal F}}}
\def\cG{{\mathbin{\cal G}}}
\def\cH{{\mathbin{\cal H}}}
\def\cI{{\mathbin{\cal I}}}
\def\cJ{{\mathbin{\cal J}}}
\def\cK{{\mathbin{\cal K}}}
\def\cL{{\mathbin{\cal L}}}
\def\cM{{\mathbin{\cal M}}}
\def\cN{{\mathbin{\cal N}}}
\def\cP{{\mathbin{\cal P}}}
\def\cQ{{\mathbin{\cal Q}}}
\def\cR{{\mathbin{\cal R}}}
\def\cS{{\mathbin{\cal S}\kern -0.1em}}
\def\cT{{\mathbin{\cal T}\kern -0.1em}}
\def\cW{{\mathbin{\cal W}}}
\def\oM{{\mathbin{\smash{\,\,\overline{\!\!\mathcal M\!}\,}}}}
\def\fCmod{{\mathbin{{\mathfrak C}\text{\rm -mod}}}}
\def\fDmod{{\mathbin{{\mathfrak D}\text{\rm -mod}}}}
\def\OXmod{{\mathbin{\O_X\text{\rm -mod}}}}
\def\OYmod{{\mathbin{\O_Y\text{\rm -mod}}}}
\def\OZmod{{\mathbin{\O_Z\text{\rm -mod}}}}
\def\OcXmod{{\mathbin{{\O\kern -0.1em}_{\cal X}\text{\rm -mod}}}}
\def\OcYmod{{\mathbin{\O_{\cal Y}\text{\rm -mod}}}}
\def\OcZmod{{\mathbin{{\O\kern -0.1em}_{\cal Z}\text{\rm -mod}}}}
\def\m{{\mathfrak m}}
\def\ua{{\underline{a}}{}}
\def\ub{{\underline{b\kern -0.1em}\kern 0.1em}{}}
\def\uc{{\underline{c}}{}}
\def\ud{{\underline{d}}{}}
\def\ue{{\underline{e}}{}}
\def\uf{{\underline{f\!}\,}{}}
\def\ug{{\underline{g\!}\,}{}}
\def\uh{{\underline{h\kern -0.1em}\kern 0.1em}{}}
\def\ui{{\underline{i\kern -0.07em}\kern 0.07em}{}}
\def\uim{{\underline{\imath\kern -0.07em}\kern 0.07em}{}}
\def\uj{{\underline{j\kern -0.1em}\kern 0.1em}{}}
\def\uk{{\underline{k\kern -0.1em}\kern 0.1em}{}}
\def\um{{\underline{m\kern -0.1em}\kern 0.1em}{}}
\def\uq{{\underline{q\kern -0.15em}\kern 0.15em}{}}
\def\ur{{\underline{r\kern -0.15em}\kern 0.15em}{}}
\def\us{{\underline{s\kern -0.15em}\kern 0.15em}{}}
\def\ut{{\underline{t\kern -0.1em}\kern 0.1em}{}}
\def\uu{{\underline{u\kern -0.1em}\kern 0.1em}{}}
\def\uv{{\underline{v\!}\,}{}}
\def\uw{{\underline{w\!}\,}{}}
\def\ux{{\underline{x\!}\,}{}}
\def\uy{{\underline{y\!}\,}{}}
\def\uz{{\underline{z\!}\,}{}}
\def\uA{{{\underline{A\!}\,}}{}}
\def\uB{{{\underline{B\!}\,}}{}}
\def\uC{{{\underline{C\!}\,}}{}}
\def\uD{{{\underline{D\!}\,}}{}}
\def\uE{{{\underline{E\!}\,}}{}}
\def\uF{{{\underline{F\!}\,}}{}}
\def\uG{{{\underline{G\!}\,}}{}}
\def\uH{{{\underline{H\!}\,}}{}}
\def\uI{{{\underline{I\kern -0.15em}\kern 0.1em}{}}}
\def\uK{{{\underline{K\!}\,}{}}}
\def\uP{{{\underline{P\!}\,}{}}}
\def\uQ{{{\underline{Q\!}\,}{}}}
\def\uR{{{\underline{R\kern -0.1em}\kern 0.1em}{}}}
\def\uS{{{\underline{S\!}\,}{}}}
\def\uT{{{\underline{T\!}\,}{}}}
\def\uU{{{\underline{U\kern -0.25em}\kern 0.2em}{}}}
\def\utU{{{\underline{\ti U\kern -0.25em}\kern 0.2em}{}}}
\def\uhU{{{\underline{\hat U\kern -0.25em}\kern 0.2em}{}}}
\def\uV{{{\underline{V\kern -0.25em}\kern 0.2em}{}}}
\def\utV{{{\underline{\ti V\kern -0.25em}\kern 0.2em}{}}}
\def\uhV{{{\underline{\hat V\kern -0.25em}\kern 0.2em}{}}}
\def\utW{{{\underline{\ti W\!\!}\,\,}{}}}
\def\uW{{{\underline{W\!\!}\,\,}{}}}
\def\uX{{{\underline{X\!}\,}{}}}
\def\uY{{{\underline{Y\!\!}\,\,}{}}}
\def\uZ{{{\underline{Z\!}\,}{}}}
\def\cU{{\cal U}}
\def\cV{{\cal V}}
\def\cW{{\cal W}}
\def\cX{{\cal X}}
\def\hcX{{\mathbin{\hat{\cal X}}}{}}
\def\tcX{{\mathbin{\tilde{\cal X}}}{}}
\def\cY{{\cal Y}}
\def\cZ{{\cal Z}}
\def\bA{{\bs A}}
\def\bB{{\bs B}}
\def\bC{{\bs C}}
\def\bD{{\bs D}}
\def\bE{{\bs E}}
\def\bF{{\bs F}}
\def\bG{{\bs G}}
\def\bH{{\bs H}}
\def\bM{{\bs M}}
\def\bN{{\bs N}}
\def\bO{{\bs O}}
\def\bP{{\bs P}}
\def\bQ{{\bs Q}}
\def\bS{{\bs S}}
\def\bT{{\bs T}}
\def\bU{{\bs U}}
\def\bV{{\bs V}}
\def\bW{{\bs W}\kern -0.1em}
\def\bX{{\bs X}}
\def\bY{{\bs Y}\kern -0.1em}
\def\bZ{{\bs Z}}
\def\bpS{{\bs{\pd S}}}
\def\bpW{{\bs{\pd W}}}
\def\bpX{{\bs{\pd X}}}
\def\bpY{{\bs{\pd Y}}}
\def\bpZ{{\bs{\pd Z}}}
\def\rE{{\bs{\rm E}}}
\def\rS{{\bs{\rm S}}}
\def\rT{{\bs{\rm T}}}
\def\rU{{\bs{\rm U}}}
\def\rV{{\bs{\rm V}}}
\def\rW{{\bs{\rm W}}}
\def\rX{{\bs{\rm X}}}
\def\rY{{\bs{\rm Y}}}
\def\rZ{{\bs{\rm Z}}}
\def\oE{{\EuScript E}}
\def\oU{{\EuScript U}}
\def\oV{{\EuScript V}}
\def\oW{{\EuScript W}}
\def\oX{{\EuScript X}}
\def\oY{{\EuScript Y}}
\def\oZ{{\EuScript Z}}
\def\eE{{\bs{\EuScript E}}}
\def\eS{{\bs{\EuScript S}}}
\def\eT{{\bs{\EuScript T}}}
\def\eU{{\bs{\EuScript U}}}
\def\eV{{\bs{\EuScript V}}}
\def\eW{{\bs{\EuScript W}}}
\def\eX{{\bs{\EuScript X}}}
\def\eY{{\bs{\EuScript Y}}}
\def\eZ{{\bs{\EuScript Z}}}
\def\ehX{{\bs{\hat{\EuScript X}}}{}}
\def\ehY{{\bs{\hat{\EuScript Y}}}{}}
\def\ehZ{{\bs{\hat{\EuScript Z}}}{}}
\def\etX{{\bs{\ti{\EuScript X}}}{}}
\def\etY{{\bs{\ti{\EuScript Y}}}{}}
\def\etZ{{\bs{\ti{\EuScript Z}}}{}}
\def\bcE{{\bs{\cal E}}}
\def\bcM{{\bs{\cal M}}}
\def\bcS{{\bs{\cal S}}}
\def\bcT{{\bs{\cal T}}}
\def\bcU{{\bs{\cal U}}}
\def\bcV{{\bs{\cal V}}}
\def\bcW{{\bs{\cal W}}{}}
\def\bcX{{\bs{\cal X}}{}}
\def\bcY{{\bs{\cal Y}}{}}
\def\bcZ{{\bs{\cal Z}}{}}
\def\bcpS{{\bs{\pd{\cal S}}}}
\def\bcpX{{\bs{\pd{\cal X}}}}
\def\bcpY{{\bs{\pd{\cal Y}}}}
\def\bcpZ{{\bs{\pd{\cal Z}}}}
\def\bcptX{{\bs{\pd^2{\cal X}}}}
\def\cpS{{\pd{\cal S}}}
\def\cpT{{\pd{\cal T}}}
\def\cpU{{\pd{\cal U}}}
\def\cpV{{\pd{\cal V}}}
\def\cpW{{\pd{\cal W}}}
\def\cpX{{\pd{\cal X}}}
\def\cpY{{\pd{\cal Y}}}
\def\cpZ{{\pd{\cal Z}}}
\def\cptU{{\pd^2{\cal U}}}
\def\cptV{{\pd^2{\cal V}}}
\def\cptW{{\pd^2{\cal W}}}
\def\cptX{{\pd^2{\cal X}}}
\def\cptY{{\pd^2{\cal Y}}}
\def\cptZ{{\pd^2{\cal Z}}}
\def\bcpX{{\bs{\pd{\cal X}}}}
\def\bcpY{{\bs{\pd{\cal Y}}}}
\def\bcpZ{{\bs{\pd{\cal Z}}}}
\def\btcX{{\bs{\ti{\cal X}}}{}}
\def\bhcX{{\bs{\hat{\cal X}}}{}}
\def\btcY{{\bs{\ti{\cal Y}}}{}}
\def\bhcY{{\bs{\hat{\cal Y}}}{}}
\def\upS{{\underline{\pd S\!}\,}}
\def\upU{{\underline{\pd U\!}\,}}
\def\upV{{\underline{\pd V\!\!}\,\,}}
\def\upX{{\underline{\pd X\!}\,}}
\def\upY{{\underline{\pd Y\!\!}\,\,}}
\def\upZ{{\underline{\pd Z\!}\,}}
\def\uptY{{\underline{\pd^2Y\!\!}\,\,}}
\def\uptZ{{\underline{\pd^2Z\!}\,}}
\def\upi{{\underline{\pi\!}\,}}
\def\uid{{\underline{\id\kern -0.1em}\kern 0.1em}}
\def\al{\alpha}
\def\be{\beta}
\def\ga{\gamma}
\def\de{\delta}
\def\io{\iota}
\def\ep{\epsilon}
\def\la{\lambda}
\def\ka{\kappa}
\def\th{\theta}
\def\ze{\zeta}
\def\up{\upsilon}
\def\vp{\varphi}
\def\si{\sigma}
\def\om{\omega}
\def\De{\Delta}
\def\La{\Lambda}
\def\Om{\Omega}
\def\Up{\Upsilon}
\def\Ga{\Gamma}
\def\Si{\Sigma}
\def\Th{\Theta}
\def\pd{\partial}
\def\ts{\textstyle}
\def\st{\scriptstyle}
\def\sst{\scriptscriptstyle}
\def\w{\wedge}
\def\sm{\setminus}
\def\lt{\ltimes}
\def\bu{\bullet}
\def\sh{\sharp}
\def\op{\oplus}
\def\od{\odot}
\def\op{\oplus}
\def\ot{\otimes}
\def\ov{\overline}
\def\bigop{\bigoplus}
\def\bigot{\bigotimes}
\def\iy{\infty}
\def\es{\emptyset}
\def\ra{\rightarrow}
\def\rra{\rightrightarrows}
\def\Ra{\Rightarrow}
\def\Longra{\Longrightarrow}
\def\ab{\allowbreak}
\def\longra{\longrightarrow}
\def\hookra{\hookrightarrow}
\def\dashra{\dashrightarrow}
\def\t{\times}
\def\ci{\circ}
\def\ti{\tilde}
\def\d{{\rm d}}
\def\md#1{\vert #1 \vert}
\def\bmd#1{\big\vert #1 \big\vert}
\def\an#1{\langle #1 \rangle}
\title{D-manifolds, d-orbifolds and derived differential geometry:
a detailed summary}
\author{Dominic Joyce}
\date{}
\maketitle

\begin{abstract} This is a detailed summary of the author's (rather
longer) book \cite{Joyc6}. We introduce a 2-category $\dMan$ of {\it
d-manifolds}, new geometric objects which are `derived' smooth
manifolds, in the sense of the `derived algebraic
geometry'\I{derived algebraic geometry} of To\"en and Lurie.
Manifolds $\Man$ embed in $\dMan$ as a full (2-)subcategory. There
are also 2-categories $\dManb,\dManc$ of {\it d-manifolds with
boundary\/} and {\it with corners}, and orbifold versions
$\dOrb,\dOrbb,\dOrbc$ of all of these, {\it d-orbifolds}.

Much of differential geometry extends very nicely to d-manifolds and
d-orbifolds --- immersions, submersions, submanifolds, transverse
fibre products, orientations, bordism groups, etc. Compact oriented
d-manifolds and d-orbifolds have virtual classes. Boundaries of
d-manifolds and d-orbifolds with corners behave in a functorial way.

Many important areas of geometry involve forming moduli
spaces\I{moduli space} $\cM$ of geometric objects, and `counting'
them to get an enumerative invariant, or a more general structure in
homological algebra, such as a Floer homology theory. These areas
include Donaldson invariants and Seiberg--Witten invariants of
4-manifolds, Donaldson--Thomas invariants\I{Donaldson--Thomas
invariants} of Calabi--Yau 3-folds,\I{Calabi--Yau 3-fold}
Gromov--Witten invariants\I{Gromov--Witten invariants} in both
algebraic and symplectic geometry,\I{symplectic geometry} and
Lagrangian Floer cohomology,\I{Lagrangian Floer cohomology} Fukaya
categories,\I{Fukaya categories} contact homology,\I{contact
homology} and Symplectic Field Theory\I{Symplectic Field Theory} in
symplectic geometry.

In all these areas, one first defines an appropriate geometric
structure on $\cM$, and then applies a `virtual class'\I{virtual
class} or `virtual chain'\I{virtual chain} construction to do the
`counting' and define the invariants. The geometric structures used
for this purpose include $\C$-schemes\I{scheme with obstruction
theory} and Deligne--Mumford $\C$-stacks with perfect obstruction
theories\I{Deligne--Mumford stack with \\ obstruction
theory}\I{stack} in complex algebraic geometry, and
polyfolds\I{polyfold} and Kuranishi spaces\I{Kuranishi space} in
symplectic geometry.

There are truncation functors\I{truncation
functor}\I{functor!truncation} from all of these classes of
geometric structures on $\cM$ to d-manifolds or d-orbifolds, with or
without corners. There are also truncation functors from
quasi-smooth derived $\C$-schemes\I{derived scheme!quasi-smooth} and
Spivak's derived manifolds\I{Spivak's derived manifolds} to
d-manifolds. As a result, all the areas of geometry above involving
`counting' moduli spaces can be rewritten in terms of d-manifolds
and d-orbifolds. This will lead to new results and simplifications
of existing proofs, particularly in areas involving moduli spaces
with boundary and corners.

A (rather shorter) survey paper on the book, focussing on
d-manifolds without boundary, is~\cite{Joyc7}.
\end{abstract}

\setcounter{tocdepth}{2}
\tableofcontents

\section{Introduction}
\label{ds1}
\G[Man]{$\Man$}{category of manifolds\nomnorefpage}

This is a summary of the author's (rather longer) book \cite{Joyc6}
on `D-manifolds and d-orbifolds: a theory of derived differential
geometry'. A (rather shorter) survey paper on the book, focussing on
d-manifolds without boundary, is \cite{Joyc7}, and readers just
wanting a general overview of the theory are advised to read
\cite{Joyc7}.

In this paper we aim to provide a fairly complete coverage of the
main definitions and results of \cite{Joyc6}, omitting almost all
proofs and some more abstruse technical details, which is suitable
to be the primary reference for those wanting to use d-manifolds and
d-orbifolds in their own research.

We develop a new theory of `derived differential geometry'. The
objects in this theory are {\it d-manifolds}, `derived' versions of
smooth manifolds, which form a (strict) 2-category $\dMan$. There
are also 2-categories of {\it d-manifolds with boundary\/} $\dManb$
and {\it d-manifolds with corners\/} $\dManc$, and orbifold versions
of all these, {\it d-orbifolds\/} $\dOrb,\dOrbb,\dOrbc$.

Here `derived' is intended in the sense of {\it derived algebraic
geometry}.\I{derived algebraic geometry|(} The original motivating
idea for derived algebraic geometry, as in Kontsevich \cite{Kont}
for instance, was that certain moduli schemes $\cM$ appearing in
enumerative invariant problems may be very singular as schemes.
However, it may be natural to realize $\cM$ as a truncation of some
`derived' moduli space $\bcM$, a new kind of geometric object living
in a higher category. The geometric structure on $\bcM$ should
encode the full deformation theory of the moduli problem, the
obstructions as well as the deformations. It was hoped that $\bcM$
would be `smooth', and so in some sense simpler than its
truncation~$\cM$.

Early work in derived algebraic geometry focussed on {\it
dg-schemes},\I{dg-scheme} as in Ciocan-Fontanine and Kapranov
\cite{CiKa}. These have largely been replaced by the {\it derived
stacks\/}\I{stack} of To\"en and Vezzosi \cite{Toen,ToVe1,ToVe2},
and the {\it structured spaces\/} of Lurie \cite{Luri1,Luri2,Luri3}.
{\it Derived differential geometry\/} aims to generalize these ideas
to differential geometry and smooth manifolds. A brief note about it
can be found in Lurie \cite[\S 4.5]{Luri3}; the ideas are worked out
in detail by Lurie's student David Spivak \cite{Spiv}, who defines
an $\iy$-category\I{$\iy$-category} of {\it derived
manifolds}.\I{Spivak's derived manifolds}

The author came to these questions from a different direction,
symplectic geometry. Many important areas in symplectic
geometry\I{symplectic geometry|(} involve forming moduli spaces
$\oM_{g,m}(X,J,\be)$ of $J$-holomorphic curves in some symplectic
manifold $(X,\om)$, possibly with boundary in a Lagrangian
$Y$,\I{Lagrangian submanifold} and then `counting' these moduli
spaces to get `invariants' with interesting properties. Such areas
include Gromov--Witten invariants\I{Gromov--Witten invariants} (open
and closed), Lagrangian Floer cohomology,\I{Lagrangian Floer
cohomology} Symplectic Field Theory,\I{Symplectic Field Theory}
contact homology,\I{contact homology} and Fukaya
categories.\I{Fukaya categories}

To do this `counting', one needs to put a suitable geometric
structure on $\oM_{g,m}(X,J,\be)$ --- something like the `derived'
moduli spaces $\bcM$ above --- and use this to define a `virtual
class' or `virtual chain'\I{virtual chain} in $\Z,\Q$ or some
homology theory. Two alternative theories for geometric structures
to put on moduli spaces $\oM_{g,m}(X,J,\be)$ are the {\it Kuranishi
spaces\/}\I{Kuranishi space} of Fukaya, Oh, Ohta and Ono
\cite{FuOn,FOOO} and the {\it polyfolds\/}\I{polyfold} of Hofer,
Wysocki and Zehnder~\cite{HWZ1,HWZ2,HWZ3,HWZ4,HWZ5,HWZ6}.

The philosophies of Kuranishi spaces and of polyfolds are in a sense
opposite: Kuranishi spaces remember only the minimal information
needed to form virtual chains, but polyfolds remember a huge amount
more information, essentially a complete description of the
functional-analytic problem which gives rise to the moduli space.
There is a truncation functor\I{truncation
functor}\I{functor!truncation} from polyfolds to Kuranishi spaces.

The theory of Kuranishi spaces in \cite{FuOn,FOOO} does not go far
--- they define Kuranishi spaces, and construct virtual cycles upon
them, but they do not define morphisms between Kuranishi spaces, for
instance. The author tried to study and work with Kuranishi spaces
as geometric spaces in their own right, but ran into problems, and
became convinced that a new definition was needed. Upon reading
Spivak's theory of derived manifolds \cite{Spiv},\I{Spivak's derived
manifolds} it became clear that some form of `derived differential
geometry' was required: {\it Kuranishi spaces in the sense of\/
{\rm\cite[\S A]{FOOO}} ought to be defined to be `derived orbifolds
with corners'}.

The purpose of \cite{Joyc6}, summarized here, is to build a
comprehensive, rigorous theory of derived differential geometry
designed for applications in symplectic geometry, and other areas of
mathematics such as String Topology.\I{String Topology}

As the moduli spaces of interest in the symplectic geometry of
Lagrangian submanifolds\I{Lagrangian submanifold} should be `derived
orbifolds with corners', it was necessary that this theory should
cover not just derived manifolds without boundary, but also derived
manifolds and derived orbifolds with boundary and with corners. It
turns out that doing `things with corners' properly is a complex,
fascinating, and hitherto almost unexplored area. This has added
considerably to the length of the project: the parts (sections
\ref{ds2}--\ref{ds4} and Appendix \ref{dsA}) dealing with
d-manifolds without boundary are only a quarter of the whole.

The author wants the theory to be easily usable by symplectic
geometers, and others who are not specialists in derived algebraic
geometry. In applications, much of the theory can be treated as a
`black box', as they do not require a detailed understanding of what
a d-manifold or d-orbifold really is, but only a general idea, plus
a list of useful properties of the
2-categories~$\dMan,\dOrb$.\I{symplectic geometry|)}

Our theory of derived differential geometry has a major
simplification compared to the derived algebraic geometry\I{derived
algebraic geometry} of To\"en and Vezzosi \cite{Toen,ToVe1,ToVe2}
and Lurie \cite{Luri1,Luri2,Luri3}, and the derived manifolds of
Spivak \cite{Spiv}.\I{Spivak's derived manifolds} All of the
`derived' spaces in \cite{Luri1,Luri2,Luri3,Toen,ToVe1,ToVe2,Spiv}
form some kind of $\iy$-category\I{$\iy$-category} (simplicial
category, model category, Segal category, quasicategory, \ldots). In
contrast, our d-manifolds and d-orbifolds form (strict) 2-categories
$\dMan,\ldots,\dOrbc$, which are the simplest and most friendly kind
of higher category.

Furthermore, the $\iy$-categories in \cite{Luri1,Luri2,Luri3,Toen,
ToVe1,ToVe2,Spiv} are usually formed by localization (inversion of
some class of morphisms), so the (higher) morphisms in the resulting
$\iy$-category are difficult to describe and work with. But the 1-
and 2-morphisms in $\dMan,\ldots,\dOrbc$ are defined explicitly,
without localization.

The essence of our simplification is this. Consider a `derived'
moduli space $\bcM$ of some objects $E$, e.g. vector bundles on some
$\C$-scheme $X$. One expects $\bcM$ to have a `cotangent
complex'\I{cotangent complex} ${\mathbb L}_\bcM$, a complex in some
derived category with cohomology $h^i({\mathbb
L}_\bcM)\vert_E\cong\Ext^{1-i}(E,E)^*$ for $i\in\Z$. In general,
${\mathbb L}_\bcM$ can have nontrivial cohomology in many negative
degrees, and because of this such objects $\bcM$ must form an
$\iy$-category to properly describe their geometry.

However, the moduli spaces relevant to enumerative invariant
problems are of a restricted kind: one considers only $\bcM$ such
that ${\mathbb L}_\bcM$ has nontrivial cohomology only in degrees
$-1,0$, where $h^0({\mathbb L}_\bcM)$ encodes the (dual of the)
deformations $\Ext^1(E,E)^*$, and $h^{-1}({\mathbb L}_\bcM)$ the
(dual of the) obstructions $\Ext^2(E,E)^*$. As in To\"en \cite[\S
4.4.3]{Toen}, such derived spaces are called {\it
quasi-smooth},\I{quasi-smooth}\I{derived scheme!quasi-smooth} and
this is a necessary condition on $\bcM$ for the construction of a
virtual fundamental class.

Our construction of d-manifolds replaces complexes in a derived
category $D^b\coh(\cM)$ with a 2-category of complexes in degrees
$-1,0$ only. For general $\bcM$ this loses a lot of information, but
for quasi-smooth $\bcM$, since ${\mathbb L}_\bcM$ is concentrated in
degrees $-1,0$, the important information is retained. In the
language of dg-schemes,\I{dg-scheme} this corresponds to working
with a subclass of derived schemes whose dg-algebras are of a
special kind: they are 2-step supercommutative dg-algebras
$A^{-1}\,{\buildrel\d\over \longra}\,A^0$ such that $\d(A^{-1})\cdot
A^{-1}=0$. Then $\d(A^{-1})$ is a square zero ideal in $A^0$, and
$A^{-1}$ is a module over~$H^0\bigl(A^{-1}\,{\buildrel
\d\over\longra}\,A^0\bigr)$.

An important reason why this 2-category style derived geometry works
successfully in our differential-geometric context is the existence
of {\it partitions of unity\/}\I{partition of unity} on smooth
manifolds, and on nice $C^\iy$-schemes. This means that (derived)
structure sheaves are `fine' or `soft', which simplifies their
behaviour. Partitions of unity are also essential for constructions
such as gluing d-manifolds by equivalences on open d-subspaces in
$\dMan$. In conventional derived algebraic geometry, where
partitions of unity do not exist, one needs the extra freedom of an
$\iy$-category\I{$\iy$-category} to glue by equivalences.\I{derived
algebraic geometry|)}

Throughout the paper, following \cite{Joyc6}, we will consistently
use different typefaces to indicate different classes of geometrical
objects. In particular:
\begin{itemize}
\setlength{\itemsep}{0pt}
\setlength{\parsep}{0pt}
\item $W,X,Y,\ldots$ will denote manifolds (of any kind), or
topological spaces.
\item $\uW,\uX,\uY,\ldots$ will denote $C^\iy$-schemes.
\item $\bW,\bX,\bY,\ldots$ will denote d-spaces, including d-manifolds.
\item $\cW,\cX,\cY,\ldots$ will denote Deligne--Mumford
$C^\iy$-stacks, including orbifolds.
\item $\bcW,\bcX,\bcY,\ldots$ will denote d-stacks, including d-orbifolds.
\item $\rW,\rX,\rY,\ldots$ will denote d-spaces with corners,
including d-manifolds with corners.
\item $\oW,\oX,\oY,\ldots$ will denote orbifolds with corners.
\item $\eW,\eX,\eY,\ldots$ will denote d-stacks with corners,
including d-orbifolds with corners.
\end{itemize}

\smallskip

\noindent{\it Acknowledgements.} My particular thanks to Dennis
Borisov, Jacob Lurie and Bertrand To\"en for help with derived
manifolds. I would also like to thank Manabu Akaho, Tom Bridgeland,
Kenji Fukaya, Hiroshi Ohta, Kauru Ono, and Timo Sch\"urg for useful
conversations.

\section{$C^\iy$-rings and $C^\iy$-schemes}
\label{ds2}
\I{C-algebraic geometry@$C^\iy$-algebraic geometry|(}

If $X$ is a manifold then the $\R$-algebra $C^\iy(X)$ of smooth
functions $c:X\ra\R$ is a $C^\iy$-{\it
ring}.\I{C-ring@$C^\iy$-ring|(} That is, for each smooth function
$f:\R^n\ra\R$ there is an $n$-fold operation $\Phi_f:C^\iy(X)^n\ra
C^\iy(X)$ acting by $\Phi_f:c_1,\ldots,c_n\mapsto
f(c_1,\ldots,c_n)$, and these operations $\Phi_f$ satisfy many
natural identities. Thus, $C^\iy(X)$ actually has a far richer
algebraic structure than the obvious $\R$-algebra structure.

$C^\iy$-{\it algebraic geometry\/} is a version of algebraic
geometry in which rings or algebras are replaced by $C^\iy$-rings.
The basic objects are $C^\iy$-{\it
schemes},\I{C-scheme@$C^\iy$-scheme} a category of
differential-geometric spaces including smooth manifolds, and also
many singular spaces. They were introduced in synthetic differential
geometry\I{synthetic differential geometry} (see for instance Dubuc
\cite{Dubu} and Moerdijk and Reyes \cite{MoRe}), and developed
further by the author in \cite{Joyc4} (surveyed in \cite{Joyc5}
and~\cite[App.~B]{Joyc6}).

This section briefly discusses $C^\iy$-rings, $C^\iy$-schemes, and
quasicoherent sheaves on $C^\iy$-schemes, following the author's
treatment~\cite[\S 2--\S 6]{Joyc4}.

\subsection{$C^\iy$-rings}
\label{ds21}

\begin{dfn} A $C^\iy$-{\it ring\/} is a set $\fC$ together with
operations $\Phi_f:\fC^n\ra\fC$\G[CDE]{$\fC,{\mathfrak D},{\mathfrak
E},\ldots$}{$C^\iy$-rings} \G[Phif]{$\Phi_f:\fC^n\ra\fC$}{operations
on $C^\iy$-ring $\fC$, for smooth $f:\R^n\ra\R$} for all $n\ge 0$
and smooth maps $f:\R^n\ra\R$, where by convention when $n=0$ we
define $\fC^0$ to be the single point $\{\es\}$. These operations
must satisfy the following relations: suppose $m,n\ge 0$, and
$f_i:\R^n\ra\R$ for $i=1,\ldots,m$ and $g:\R^m\ra\R$ are smooth
functions. Define a smooth function $h:\R^n\ra\R$ by
\begin{equation*}
h(x_1,\ldots,x_n)=g\bigl(f_1(x_1,\ldots,x_n),\ldots,f_m(x_1
\ldots,x_n)\bigr),
\end{equation*}
for all $(x_1,\ldots,x_n)\in\R^n$. Then for all
$(c_1,\ldots,c_n)\in\fC^n$ we have
\begin{equation*}
\Phi_h(c_1,\ldots,c_n)=\Phi_g\bigl(\Phi_{f_1}(c_1,\ldots,c_n),
\ldots,\Phi_{f_m}(c_1,\ldots,c_n)\bigr).
\end{equation*}
We also require that for all $1\le j\le n$, defining
$\pi_j:\R^n\ra\R$ by $\pi_j:(x_1,\ldots,x_n)\mapsto x_j$, we have
$\Phi_{\pi_j}(c_1,\ldots,c_n)=c_j$ for
all~$(c_1,\ldots,c_n)\in\fC^n$.

Usually we refer to $\fC$ as the $C^\iy$-ring, leaving the
operations $\Phi_f$ implicit.

A {\it morphism\/} between $C^\iy$-rings $\bigl(\fC,(\Phi_f)_{
f:\R^n\ra\R\,\,C^\iy}\bigr)$, $\bigl({\mathfrak
D},(\Psi_f)_{f:\R^n\ra\R\,\,C^\iy}\bigr)$ is a map
$\phi:\fC\ra{\mathfrak D}$ such that $\Psi_f\bigl(\phi
(c_1),\ldots,\phi(c_n)\bigr)=\phi\ci\Phi_f(c_1,\ldots,c_n)$ for all
smooth $f:\R^n\ra\R$ and $c_1,\ldots,c_n\in\fC$. We will write
$\CRings$ for the category of
$C^\iy$-rings.\G[CRings]{$\CRings$}{category of $C^\iy$-rings}
\label{ds2def1}
\end{dfn}

Here is the motivating example:

\begin{ex} Let $X$ be a manifold. Write $C^\iy(X)$ for the set of
smooth functions $c:X\ra\R$. For $n\ge 0$ and $f:\R^n\ra\R$ smooth,
define $\Phi_f:C^\iy(X)^n\ra C^\iy(X)$ by
\e
\bigl(\Phi_f(c_1,\ldots,c_n)\bigr)(x)=f\bigl(c_1(x),\ldots,c_n(x)\bigr),
\label{ds2eq1}
\e
for all $c_1,\ldots,c_n\in C^\iy(X)$ and $x\in X$. It is easy to see
that $C^\iy(X)$ and the operations $\Phi_f$ form a~$C^\iy$-ring.

Now let $f:X\ra Y$ be a smooth map of manifolds. Then pullback
$f^*:C^\iy(Y)\ra C^\iy(X)$ mapping $f^*:c\mapsto c\ci f$ is a
morphism of $C^\iy$-rings. Furthermore (at least for $Y$ without
boundary), every $C^\iy$-ring morphism $\phi:C^\iy(Y)\ra C^\iy(X)$
is of the form $\phi=f^*$ for a unique smooth map~$f:X\ra Y$.

Write $\CRings^{\rm op}$ for the opposite category of $\CRings$,
with directions of morphisms reversed, and $\Man$ for the category
of manifolds without boundary. Then we have a full and faithful
functor\I{functor!full}\I{functor!faithful}
$F_\Man^\CRings:\Man\ra\CRings^{\rm op}$ acting by
$F_\Man^\CRings(X)=C^\iy(X)$ on objects and $F_\Man^\CRings(f)=f^*$
on morphisms. This embeds $\Man$ as a full subcategory
of~$\CRings^{\rm op}$.
\label{ds2ex1}
\end{ex}

Note that $C^\iy$-rings are far more general than those coming from
manifolds. For example, if $X$ is any topological space we could
define a $C^\iy$-ring $C^0(X)$ to be the set of {\it continuous\/}
$c:X\ra\R$, with operations $\Phi_f$ defined as in \eq{ds2eq1}. For
$X$ a manifold with $\dim X>0$, the $C^\iy$-rings $C^\iy(X)$ and
$C^0(X)$ are different.

\begin{dfn} Let $\fC$ be a $C^\iy$-ring. Then we may give $\fC$ the
structure of a {\it commutative\/ $\R$-algebra}. Define addition
`$+$' on $\fC$ by $c+c'=\Phi_f(c,c')$ for $c,c'\in\fC$, where
$f:\R^2\ra\R$ is $f(x,y)=x+y$. Define multiplication `$\,\cdot\,$'
on $\fC$ by $c\cdot c'=\Phi_g(c,c')$, where $g:\R^2\ra\R$ is
$g(x,y)=xy$. Define scalar multiplication by $\la\in\R$ by $\la
c=\Phi_{\la'}(c)$, where $\la':\R\ra\R$ is $\la'(x)=\la x$. Define
elements $0,1\in\fC$ by $0=\Phi_{0'}(\es)$ and $1=\Phi_{1'}(\es)$,
where $0':\R^0\ra\R$ and $1':\R^0\ra\R$ are the maps $0':\es\mapsto
0$ and $1':\es\mapsto 1$. One can show using the relations on the
$\Phi_f$ that the axioms of a commutative $\R$-algebra are
satisfied. In Example \ref{ds2ex1}, this yields the obvious
$\R$-algebra structure on the smooth functions~$c:X\ra\R$.

An {\it ideal\/} $I$ in $\fC$ is an ideal $I\subset\fC$ in $\fC$
regarded as a commutative $\R$-algebra. Then we make the quotient
$\fC/I$ into a $C^\iy$-ring as follows. If $f:\R^n\ra\R$ is smooth,
define $\Phi_f^I:(\fC/I)^n\ra\fC/I$ by
\begin{equation*}
\bigl(\Phi_f^I(c_1+I,\ldots,c_n+I)\bigr)(x)=f\bigl(c_1(x),\ldots,
c_n(x)\bigr)+I.
\end{equation*}
Using Hadamard's Lemma,\I{Hadamard's Lemma} one can show that this
is independent of the choice of representatives $c_1,\ldots,c_n$.
Then $\bigl(\fC/I,(\Phi_f^I)_{ f:\R^n\ra\R\,\,C^\iy}\bigr)$ is a
$C^\iy$-ring.

A $C^\iy$-ring $\fC$ is called {\it finitely
generated\/}\I{C-ring@$C^\iy$-ring!finitely generated} if there
exist $c_1,\ldots,c_n$ in $\fC$ which generate $\fC$ over all
$C^\iy$-operations. That is, for each $c\in\fC$ there exists smooth
$f:\R^n\ra\R$ with $c=\Phi_f(c_1,\ldots,c_n)$. Given such
$\fC,c_1,\ldots,c_n$, define $\phi:C^\iy(\R^n)\ra\fC$ by
$\phi(f)=\Phi_f(c_1,\ldots,c_n)$ for smooth $f:\R^n\ra\R$, where
$C^\iy(\R^n)$ is as in Example \ref{ds2ex1} with $X=\R^n$. Then
$\phi$ is a surjective morphism of $C^\iy$-rings, so $I=\Ker\phi$ is
an ideal in $C^\iy(\R^n)$, and $\fC\cong C^\iy(\R^n)/I$ as a
$C^\iy$-ring. Thus, $\fC$ is finitely generated if and only if
$\fC\cong C^\iy(\R^n)/I$ for some $n\ge 0$ and some ideal $I$
in~$C^\iy(\R^n)$.
\label{ds2def2}
\end{dfn}
\I{C-ring@$C^\iy$-ring|)}

\subsection{$C^\iy$-schemes}
\label{ds22}
\I{C-scheme@$C^\iy$-scheme|(}

Next we summarize material in \cite[\S 4]{Joyc4} on $C^\iy$-schemes.

\begin{dfn} A {\it $C^\iy$-ringed space\/}\I{C-ringed
space@$C^\iy$-ringed space} $\uX=(X,\O_X)$ is a topological space
$X$ with a sheaf $\O_X$ of $C^\iy$-rings on $X$.

A {\it morphism\/} $\uf=(f,f^\sh):(X,\O_X)\ra (Y,\O_Y)$ of $C^\iy$
ringed spaces is a continuous map $f:X\ra Y$ and a morphism
$f^\sh:f^{-1}(\OY)\ra\OX$ of sheaves of $C^\iy$-rings on $X$, where
$f^{-1}(\OY)$ is the inverse image sheaf. There is another way to
write the data $f^\sh$: since direct image of sheaves $f_*$ is right
adjoint to inverse image $f^{-1}$, there is a natural bijection
\e
\Hom_X\bigl(f^{-1}(\O_Y),\O_X\bigr)\cong\Hom_Y\bigl(\O_Y,f_*(\O_X)\bigr).
\label{ds2eq2}
\e
Write $f_\sh:\O_Y\ra f_*(\O_X)$ for the morphism of sheaves of
$C^\iy$-rings on $Y$ corresponding to $f^\sh$ under \eq{ds2eq2}, so
that
\e
f^\sh:f^{-1}(\O_Y)\longra\O_X\quad \leftrightsquigarrow\quad
f_\sh:\O_Y\longra f_*(\O_X).
\label{ds2eq3}
\e
Depending on the application, either $f^\sh$ or $f_\sh$ may be more
useful. We choose to regard $f^\sh$ as primary and write morphisms
as $\uf=(f,f^\sh)$ rather than $(f,f_\sh)$, because we find it
convenient to work uniformly using pullbacks, rather than mixing
pullbacks and pushforwards.

Write $\CRS$ for the category of $C^\iy$-ringed
spaces.\G[WXYZa]{$\uW,\uX,\uY,\uZ,\ldots$}{$C^\iy$-schemes} As in
\cite[Th.~8]{Dubu} there is a {\it spectrum
functor\/}\I{C-scheme@$C^\iy$-scheme!spectrum functor}
$\Spec:\CRings^{\rm op}\ra\CRS$, defined explicitly in
\cite[Def.~4.12]{Joyc4}. A $C^\iy$-ringed space $\uX$ is called an
{\it affine\/ $C^\iy$-scheme\/}\I{C-scheme@$C^\iy$-scheme!affine} if
it is isomorphic in $\CRS$ to $\Spec\fC$ for some $C^\iy$-ring
$\fC$. A $C^\iy$-ringed space $\uX=(X,\OX)$ is called a $C^\iy$-{\it
scheme\/} if $X$ can be covered by open sets $U\subseteq X$ such
that $(U,\O_X\vert_U)$ is an affine $C^\iy$-scheme. Write $\CSch$
for the full subcategory of $C^\iy$-schemes
in~$\CRS$.\G[CSch]{$\CSch$}{category of $C^\iy$-schemes}

A $C^\iy$-scheme $\uX=(X,\OX)$ is called {\it locally
fair}\I{C-scheme@$C^\iy$-scheme!locally fair} if $X$ can be covered
by open $U\subseteq X$ with $(U,\O_X\vert_U)\cong\Spec\fC$ for some
finitely generated $C^\iy$-ring $\fC$. Roughly speaking this means
that $\uX$ is locally finite-dimensional. Write $\CSchlf$ for the
full subcategory of locally fair $C^\iy$-schemes in~$\CSch$.

We call a $C^\iy$-scheme $\uX$ {\it separated, second countable,
compact, locally compact}, or {\it paracompact}, if the underlying
topological space $X$ is Hausdorff, second countable, compact,
locally compact, or paracompact, respectively.
\label{ds2def3}
\end{dfn}

We define a $C^\iy$-scheme $\uX$ for each manifold~$X$.

\begin{ex} Let $X$ be a manifold. Define a $C^\iy$-ringed space
$\uX=(X,\O_X)$ to have topological space $X$ and $\O_X(U)=C^\iy(U)$
for each open $U\subseteq X$, where $C^\iy(U)$ is the $C^\iy$-ring
of smooth maps $c:U\ra\R$, and if $V\subseteq U\subseteq X$ are open
define $\rho_{UV}:C^\iy(U)\ra C^\iy(V)$ by $\rho_{UV}:c\mapsto
c\vert_V$. Then $\uX=(X,\O_X)$ is a local $C^\iy$-ringed space. It
is canonically isomorphic to $\Spec C^\iy(X)$, and so is an affine
$C^\iy$-scheme. It is locally fair.

Define a functor $F_\Man^\CSch:\Man\ra\CSchlf\subset\CSch$ by
$F_\Man^\CSch=\Spec\ci F_\Man^\CRings$. Then $F_\Man^\CSch$ is full
and faithful,\I{functor!full}\I{functor!faithful} and embeds $\Man$
as a full subcategory of $\CSch$.
\label{ds2ex2}
\end{ex}

By \cite[Cor.~4.21 \& Th.~4.33]{Joyc4} we have:
\I{C-scheme@$C^\iy$-scheme!fibre products|(}\I{fibre product!of
C-schemes@of $C^\iy$-schemes}

\begin{thm} Fibre products and all finite limits exist in\/
$\CSch$. The subcategory $\CSchlf$ is closed under fibre products
and finite limits. The functor\/ $F_\Man^\CSch$ takes transverse
fibre products in $\Man$ to fibre products in\/~$\CSch$.
\label{ds2thm1}
\end{thm}
\I{manifold!transverse fibre products}

The proof of the existence of fibre products in $\CSch$ follows that
for fibre products of schemes in Hartshorne \cite[Th.~II.3.3]{Hart},
together with the existence of $C^\iy$-scheme products $\uX\t\uY$ of
affine $C^\iy$-schemes $\uX,\uY$. The latter follows from the
existence of coproducts $\fC\hat\ot\fD$ in $\CRings$ of
$C^\iy$-rings $\fC,\fD$. Here $\fC\hat\ot\fD$ may be thought of as a
`completed tensor product' of $\fC,\fD$. The actual tensor product
$\fC\ot_\R\fD$ is naturally an $\R$-algebra but not a $C^\iy$-ring,
with an inclusion of $\R$-algebras $\fC\ot_\R\fD\hookra
\fC\hat\ot\fD$, but $\fC\hat\ot\fD$ is often much larger than
$\fC\ot_\R\fD$. For free $C^\iy$-rings we have~$C^\iy(\R^m)\hat\ot
C^\iy(\R^n)\cong C^\iy(\R^{m+n})$.\I{C-scheme@$C^\iy$-scheme!fibre
products|)}

In \cite[Def.~4.34 \& Prop.~4.35]{Joyc4} we discuss {\it partitions
of unity\/} on $C^\iy$-schemes.\I{partition of unity|(}

\begin{dfn} Let $\uX=(X,\O_X)$ be a $C^\iy$-scheme. Consider a
formal sum $\sum_{a\in A}c_a$, where $A$ is an indexing set and
$c_a\in\O_X(X)$ for $a\in A$. We say $\sum_{a\in A}c_a$ is a {\it
locally finite sum on\/} $\uX$ if $X$ can be covered by open
$U\subseteq X$ such that for all but finitely many $a\in A$ we have
$\rho_{XU}(c_a)=0$ in~$\O_X(U)$.

By the sheaf axioms for $\O_X$, if $\sum_{a\in A}c_a$ is a locally
finite sum there exists a unique $c\in\O_X(X)$ such that for all
open $U\subseteq X$ with $\rho_{XU}(c_a)=0$ in $\O_X(U)$ for all but
finitely many $a\in A$, we have $\rho_{XU}(c)=\sum_{a\in
A}\rho_{XU}(c_a)$ in $\O_X(U)$, where the sum makes sense as there
are only finitely many nonzero terms. We call $c$ the {\it limit\/}
of $\sum_{a\in A}c_a$, written~$\sum_{a\in A}c_a=c$.

Let $c\in\O_X(X)$. Then there is a unique maximal open set
$V\subseteq X$ with $\rho_{XV}(c)=0$ in $\O_X(V)$. Define the {\it
support\/} $\supp c$ to be $X\sm V$, so that $\supp c$ is closed in
$X$. If $U\subseteq X$ is open, we say that $c$ {\it is supported
in\/} $U$ if~$\supp c\subseteq U$.

Let $\{U_a:a\in A\}$ be an open cover of $X$. A {\it partition of
unity on\/ $\uX$ subordinate to\/} $\{U_a:a\in A\}$ is
$\{\eta_a:a\in A\}$ with $\eta_a\in\O_X(X)$ supported on $U_a$ for
$a\in A$, such that $\sum_{a\in A}\eta_a$ is a locally finite sum on
$\uX$ with~$\sum_{a\in A}\eta_a=1$.
\label{ds2def4}
\end{dfn}

\begin{prop} Suppose $\uX$ is a separated, paracompact, locally fair
$C^\iy$-scheme, and\/ $\{\uU_a:a\in A\}$ an open cover of\/ $\uX$.
Then there exists a partition of unity $\{\eta_a:a\in A\}$ on $\uX$
subordinate to\/~$\{\uU_a:a\in A\}$.\I{partition of unity|)}
\label{ds2prop1}
\end{prop}

Here are some differences between ordinary schemes and
$C^\iy$-schemes:

\begin{rem}{\bf(i)} If $A$ is a ring or algebra, then points of the
corresponding scheme $\Spec A$ are prime ideals in $A$. However, if
$\fC$ is a $C^\iy$-ring then (by definition) points of $\Spec\fC$
are maximal ideals in $\fC$ with residue field $\R$, or
equivalently, $\R$-algebra morphisms $x:\fC\ra\R$. This has the
effect that if $X$ is a manifold then points of $\Spec C^\iy(X)$ are
just points of~$X$.
\smallskip

\noindent{\bf(ii)} In conventional algebraic geometry, affine
schemes are a restrictive class. Central examples such as $\CP^n$
are not affine, and affine schemes are not closed under open
subsets, so that $\C^2$ is affine but $\C^2\sm\{0\}$ is not. In
contrast, affine $C^\iy$-schemes are already general enough for many
purposes. For example:
\begin{itemize}
\setlength{\itemsep}{0pt}
\setlength{\parsep}{0pt}
\item All manifolds are fair affine $C^\iy$-schemes.
\item Open $C^\iy$-subschemes of fair affine $C^\iy$-schemes are
fair and affine.
\item Separated, second countable, locally fair $C^\iy$-schemes
are affine.
\end{itemize}
Affine $C^\iy$-schemes are always separated (Hausdorff), so we need
general $C^\iy$-schemes to include non-Hausdorff behaviour.
\smallskip

\noindent{\bf(iii)} In conventional algebraic geometry the Zariski
topology\I{Zariski topology} is too coarse for many purposes, so one
has to introduce the \'etale topology.\I{etale topology@\'etale
topology} In $C^\iy$-algebraic geometry there is no need for this,
as affine $C^\iy$-schemes are Hausdorff.
\smallskip

\noindent{\bf(iv)} Even very basic $C^\iy$-rings such as
$C^\iy(\R^n)$ for $n>0$ are not noetherian as $\R$-algebras. So
$C^\iy$-schemes should be compared to non-noetherian schemes in
conventional algebraic geometry.
\smallskip

\noindent{\bf(v)} The existence of partitions of unity,\I{partition
of unity} as in Proposition \ref{ds2prop1}, makes some things easier
in $C^\iy$-algebraic geometry than in conventional algebraic
geometry. For example, geometric objects can often be `glued
together' over the subsets of an open cover using partitions of
unity, and if $\cE$ is a quasicoherent sheaf on a separated,
paracompact, locally fair $C^\iy$-scheme $\uX$ then $H^i(\cE)=0$
for~$i>0$.
\label{ds2rem1}
\end{rem}
\I{C-scheme@$C^\iy$-scheme|)}

\subsection{Modules over $C^\iy$-rings, and cotangent modules}
\label{ds23}
\I{module over C-ring@module over $C^\iy$-ring|(}

In \cite[\S 5]{Joyc4} we discuss modules over $C^\iy$-rings.

\begin{dfn} Let $\fC$ be a $C^\iy$-ring. A $\fC$-{\it module} $M$
is a module over $\fC$ regarded as a commutative $\R$-algebra as in
Definition \ref{ds2def2}. $\fC$-modules form an abelian
category,\I{abelian category} which we write as $\fCmod$. For
example, $\fC$ is a $\fC$-module, and more generally $\fC\ot_\R V$
is a $\fC$-module for any real vector space $V$. Let
$\phi:\fC\ra\fD$ be a morphism of $C^\iy$-rings. If $M$ is a
$\fC$-module then $\phi_*(M)=M\ot_\fC\fD$ is a $\fD$-module. This
induces a functor~$\phi_*:\fCmod\ra\fDmod$.
\label{ds2def5}
\end{dfn}

\begin{ex} Let $X$ be a manifold, and $E\ra X$  a vector bundle.
Write $C^\iy(E)$ for the vector space of smooth sections $e$ of $E$.
Then $C^\iy(X)$ acts on $C^\iy(E)$ by multiplication, so $C^\iy(E)$
is a $C^\iy(X)$-module.
\label{ds2ex3}
\end{ex}

In \cite[\S 5.3]{Joyc4} we define the {\it cotangent module\/}
$\Om_\fC$ of a $C^\iy$-ring~$\fC$.\I{C-ring@$C^\iy$-ring!cotangent
module $\Om_\fC$|(}

\begin{dfn} Let $\fC$ be a $C^\iy$-ring, and $M$ a
$\fC$-module. A $C^\iy$-{\it derivation} is an $\R$-linear map
$\d:\fC\ra M$ such that whenever $f:\R^n\ra\R$ is a smooth map and
$c_1,\ldots,c_n\in\fC$, we have
\begin{equation*}
\d\Phi_f(c_1,\ldots,c_n)=\ts\sum_{i=1}^n\Phi_{\frac{\pd f}{\pd
x_i}}(c_1,\ldots,c_n)\cdot\d c_i.
\end{equation*}
We call such a pair $M,\d$ a {\it cotangent module\/} for $\fC$ if
it has the universal property that for any $\fC$-module $M'$ and
$C^\iy$-derivation $\d':\fC\ra M'$, there exists a unique morphism
of $\fC$-modules $\phi:M\ra M'$ with~$\d'=\phi\ci\d$.

Define $\Om_\fC$ to be the quotient of the free $\fC$-module with
basis of symbols $\d c$ for $c\in\fC$ by the $\fC$-submodule spanned
by all expressions of the form
$\d\bigl(\Phi_f(c_1,\ldots,c_n)\bigr)-\sum_{i=1}^n \Phi_{\frac{\pd
f}{\pd x_i}}(c_1,\ldots,c_n)\cdot\d c_i$ for $f:\R^n\ra\R$ smooth
and $c_1,\ldots,c_n\in\fC$, and define $\d_\fC:\fC\ra \Om_\fC$ by
$\d_\fC:c\mapsto\d c$. Then $\Om_\fC,\d_\fC$ is a cotangent module
for $\fC$. Thus cotangent modules always exist, and are unique up to
unique isomorphism.

Let $\fC,\fD$ be $C^\iy$-rings with cotangent modules
$\Om_\fC,\d_\fC$, $\Om_\fD,\d_\fD$, and $\phi:\fC\ra\fD$ be a
morphism of $C^\iy$-rings. Then $\phi$ makes $\Om_\fD$ into a
$\fC$-module, and there is a unique morphism
$\Om_\phi:\Om_\fC\ra\Om_\fD$ in $\fC$-mod with
$\d_\fD\ci\phi=\Om_\phi\ci\d_\fC$. This induces a morphism
$(\Om_\phi)_*:\Om_\fC\ot_\fC\fD\ra\Om_\fD$ in $\fD$-mod
with~$(\Om_\phi)_*\ci (\d_\fC\ot\id_\fD)=\d_\fD$.
\label{ds2def6}
\end{dfn}

\begin{ex} Let $X$ be a manifold. Then the cotangent bundle $T^*X$
is a vector bundle over $X$, so as in Example \ref{ds2ex3} it yields
a $C^\iy(X)$-module $C^\iy(T^*X)$. The exterior derivative
$\d:C^\iy(X)\ra C^\iy(T^*X)$ is a $C^\iy$-derivation. These
$C^\iy(T^*X),\d$ have the universal property in Definition
\ref{ds2def6}, and so form a {\it cotangent module\/}
for~$C^\iy(X)$.

Now let $X,Y$ be manifolds, and $f:X\ra Y$ be smooth. Then
$f^*(TY),TX$ are vector bundles over $X$, and the derivative of $f$
is a vector bundle morphism $\d f:TX\ra f^*(TY)$. The dual of this
morphism is $\d f^*:f^*(T^*Y)\ra T^*X$. This induces a morphism of
$C^\iy(X)$-modules $(\d f^*)_*:C^\iy\bigl(f^*(T^*Y)\bigr)\ra
C^\iy(T^*X)$. This $(\d f^*)_*$ is identified with $(\Om_{f^*})_*$
in Definition \ref{ds2def6} under the natural isomorphism
$C^\iy\bigl(f^*(T^*Y)\bigr)\cong C^\iy(T^*Y)\ot_{C^\iy(Y)}C^\iy(X)$.
\label{ds2ex4}
\end{ex}

Definition \ref{ds2def6} abstracts the notion of cotangent bundle of
a manifold in a way that makes sense for any
$C^\iy$-ring.\I{C-ring@$C^\iy$-ring!cotangent module $\Om_\fC$|)}
\I{module over C-ring@module over $C^\iy$-ring|)}

\subsection{Quasicoherent sheaves on $C^\iy$-schemes}
\label{ds24}
\I{C-scheme@$C^\iy$-scheme!quasicoherent sheaves on|(}

In \cite[\S 6]{Joyc4} we discuss sheaves of modules on
$C^\iy$-schemes.

\begin{dfn} Let $\uX=(X,\O_X)$ be a $C^\iy$-scheme. An $\O_X$-{\it
module\/} $\cE$ on $\uX$ assigns a module $\cE(U)$ over $\O_X(U)$
for each open set $U\subseteq X$, with $\O_X(U)$-action
$\mu_U:\O_X(U)\t\cE(U)\ra\cE(U)$, and a linear map
$\cE_{UV}:\cE(U)\ra\cE(V)$ for each inclusion of open sets
$V\subseteq U\subseteq X$, such that the following commutes:
\begin{equation*}
\xymatrix@R=10pt@C=60pt{ \O_X(U)\t \cE(U) \ar[d]^{\rho_{UV}\t
\cE_{UV}} \ar[r]_{\mu_U} & \cE(U) \ar[d]_{\cE_{UV}} \\
\O_X(V)\t \cE(V)\ar[r]^{\mu_V} & \cE(V),}
\end{equation*}
and all this data $\cE(U),\cE_{UV}$ satisfies the usual sheaf
axioms~\cite[\S II.1]{Hart} .

A {\it morphism of\/ $\O_X$-modules\/} $\phi:\cE\ra\cF$ assigns a
morphism of $\O_X(U)$-modules $\phi(U):\cE(U)\ra\cF(U)$ for each
open set $U\subseteq X$, such that $\phi(V)\ci\cE_{UV}=
\cF_{UV}\ci\phi(U)$ for each inclusion of open sets $V\subseteq
U\subseteq X$. Then $\O_X$-modules form an abelian
category,\I{abelian category} which we write as~$\OXmod$.

As in \cite[\S 6.2]{Joyc4}, the spectrum
functor\I{C-scheme@$C^\iy$-scheme!spectrum functor}
$\Spec:\CRings^{\rm op}\ab\ra\CSch$ has a counterpart for modules:
if $\fC$ is a $C^\iy$-ring and $(X,\O_X)=\Spec\fC$ we can define a
functor $\MSpec:\fCmod\ra\OXmod$. If $\fC$ is a {\it fair\/}
$C^\iy$-ring, there is a full abelian subcategory $\fCmod^{\rm co}$
of {\it complete\/}\I{module over C-ring@module over
$C^\iy$-ring!complete} $\fC$-modules in $\fCmod$, such that
$\MSpec\vert_{\fCmod^{\rm co}}:\fCmod^{\rm co}\ra\OXmod$ is an
equivalence of categories, with quasi-inverse the global sections
functor $\Ga:\OXmod\ra\fCmod^{\rm co}$. Let $\uX=(X,\O_X)$ be a
$C^\iy$-scheme, and $\cE$ an $\O_X$-module. We call $\cE$ {\it
quasicoherent\/} if $\uX$ can be covered by open $\uU$ with
$\uU\cong\Spec\fC$ for some $C^\iy$-ring $\fC$, and under this
identification $\cE\vert_U\cong\MSpec M$ for some $\fC$-module $M$.
We call $\cE$ a {\it vector bundle of rank\/} $n\ge 0$ if $\uX$ may
be covered by open $\uU$ such
that~$\cE\vert_\uU\cong\O_U\ot_\R\R^n$.\I{C-scheme@$C^\iy$-scheme!vector
bundles on}

Write $\qcoh(\uX),\vect(\uX)$ for the full subcategories of
quasicoherent sheaves and vector bundles in $\OXmod$. Then
$\qcoh(\uX)$ is an abelian category.\I{abelian category} Since
$\MSpec:\fCmod^{\rm co}\ra\OXmod$ is an equivalence for $\fC$ fair
and $(X,\O_X)=\Spec\fC$, as in \cite[Cor.~6.11]{Joyc4} we see that
if $\uX$ is a locally fair $C^\iy$-scheme then every $\O_X$-module
$\cE$ on $\uX$ is quasicoherent, that is,~$\qcoh(\uX)=\OXmod$.
\label{ds2def7}
\end{dfn}

\begin{rem}{\bf(a)} If $\uX$ is a separated, paracompact, locally
fair $C^\iy$-scheme then vector bundles on $\uX$ are projective
objects in the abelian category~$\qcoh(\uX)$.

\smallskip

\noindent{\bf(b)} In \cite[\S 6.3]{Joyc4} we also define a
subcategory $\coh(\uX)$ of {\it coherent
sheaves\/}\I{C-scheme@$C^\iy$-scheme!coherent sheaves on} in
$\qcoh(\uX)$. But we will not use them in this paper, as they do not
have all the good properties we want. In conventional algebraic
geometry, one usually restricts to noetherian schemes, where
coherent sheaves are well behaved, and form an abelian category.
However, as in Remark \ref{ds2rem1}(iv), even very basic
$C^\iy$-schemes $\uX$ such as $\ul\R^n$ for $n>0$ are
non-noetherian. Because of this, $\coh(\uX)$ is not closed under
kernels in $\qcoh(\uX)$, and is not an abelian category.\I{abelian
category}
\label{ds2rem2}
\end{rem}

\begin{dfn} Let\I{C-scheme@$C^\iy$-scheme!quasicoherent sheaves on!pullback} $\uf:\uX\ra\uY$ be a morphism of $C^\iy$-schemes,
and let $\cE$ be an $\OY$-module. Define the {\it pullback\/}
$\uf^*(\cE)$, an $\OX$-module, by $\uf^*(\cE)=f^{-1}(\cE)
\ot_{f^{-1}(\OY)}\OX$, where $f^{-1}(\cE),f^{-1}(\OY)$ are inverse
image sheaves, and the tensor product uses the morphism
$f^\sh:f^{-1}(\O_Y)\ra\O_X$. If $\phi:\cE\ra\cF$ is a morphism in
$\OYmod$ we have an induced morphism $\uf^*(\phi)=f^{-1}(\phi)
\ot\id_{\OX}:\uf^*(\cE)\ra\uf^*(\cF)$ in $\OXmod$. Then
$\uf^*:\OYmod\ra\OXmod$ is a right exact functor, which restricts to
a right exact functor~$\uf^*:\qcoh(\uY)\ra\qcoh(\uX)$.
\label{ds2def8}
\end{dfn}

\begin{rem} Pullbacks $\uf^*(\cE)$ are characterized by a universal
property, and so are {\it unique up to canonical isomorphism},
rather than unique. Our definition of $\uf^*(\cE)$ is not functorial
in $\uf$. That is, if $\uf:\uX\ra\uY$, $\ug:\uY\ra\uZ$ are morphisms
and $\cE\in\OZmod$ then $(\ug\ci\uf)^*(\cE)$ and $\uf^*(\ug^*(\cE))$
are canonically isomorphic in $\OXmod$, but may not be equal. We
will write $I_{\uf,\ug}(\cE):(\ug\ci\uf)^*(\cE)\ra\uf^*(\ug^*(\cE))$
for these canonical isomorphisms. Then $I_{\uf,\ug}:(\ug\ci\uf)^*\Ra
\uf^*\ci\ug^*$ is a natural isomorphism of functors.

Similarly, when $\uf$ is the identity $\uid_{\smash{\uX}}:\uX\ra\uX$
and $\cE\in\OXmod$ we may not have $\uid^*_{\smash{\uX}}(\cE)=\cE$,
but there is a canonical isomorphism
$\de_\uX(\cE):\uid^*_{\smash{\uX}}(\cE)\ra\cE$, and
$\de_\uX:\uid^*_{\smash{\uX}}\Ra\id_\OXmod$ is a natural isomorphism
of functors.

In fact it is a common abuse of notation in algebraic geometry to
omit these isomorphisms $I_{\uf,\ug}(\cE),
\uid^*_{\smash{\uX}}(\cE)$, and just assume that
$(\ug\ci\uf)^*(\cE)=\uf^*(\ug^*(\cE))$ and
$\uid^*_{\smash{\uX}}(\cE)=\cE$. An author who treats them
rigorously is Vistoli \cite{Vist}, see in particular
\cite[Introduction \& \S 3.2.1]{Vist}. One reason we decided to
include them is to be sure that $\dSpa,\dMan,\ldots$ defined below
are strict 2-categories,\I{2-category!strict} rather than weak
2-categories\I{2-category!weak} or some other structure.
\label{ds2rem3}
\end{rem}

\begin{ex} Let $X$ be a manifold, and $\uX$ the associated
$C^\iy$-scheme from Example \ref{ds2ex2}, so that $\O_X(U)=C^\iy(U)$
for all open $U\subseteq X$. Let $E\ra X$ be a vector bundle. Define
an $\OX$-module $\cE$ on $\uX$ by $\cE(U)=C^\iy(E\vert_U)$, the
smooth sections of the vector bundle $E\vert_U\ra U$, and for open
$V\subseteq U\subseteq X$ define $\cE_{UV}:\cE(U)\ra\cE(V)$ by
$\cE_{UV}:e_U\mapsto e_U\vert_V$. Then $\cE\in\vect(\uX)$ is a
vector bundle on $\uX$, which we think of as a lift of $E$ from
manifolds to $C^\iy$-schemes.

Let $f:X\ra Y$ be a smooth map of manifolds, and $\uf:\uX\ra\uY$ the
corresponding morphism of $C^\iy$-schemes. Let $F\ra Y$ be a vector
bundle over $Y$, so that $f^*(F)\ra X$ is a vector bundle over $X$.
Let $\cF\in\vect(\uY)$ be the vector bundle over $\uY$ lifting $F$.
Then $\uf^*(\cF)$ is canonically isomorphic to the vector bundle
over $\uX$ lifting~$f^*(F)$.
\label{ds2ex5}
\end{ex}

We define {\it cotangent sheaves}, the sheaf version of cotangent
modules in~\S\ref{ds23}.\I{C-scheme@$C^\iy$-scheme!cotangent
sheaf|(}

\begin{dfn} Let $\uX$ be a $C^\iy$-scheme. Define
${\cal P}T^*\uX$ to associate to each open $U\subseteq X$ the
cotangent module $\Om_{\O_X(U)}$, and to each inclusion of open sets
$V\subseteq U\subseteq X$ the morphism of $\O_X(U)$-modules
$\Om_{\rho_{UV}}:\Om_{\O_X(U)}\ra\Om_{\O_X(V)}$ associated to the
morphism of $C^\iy$-rings $\rho_{UV}:\O_X(U)\ra\O_X(V)$. Then ${\cal
P}T^*\uX$ is a {\it presheaf of\/ $\O_X$-modules on\/} $\uX$. Define
the {\it cotangent sheaf\/ $T^*\uX$ of\/} $\uX$ to be the
sheafification of ${\cal P}T^*\uX$, as an $\OX$-module.

Let $\uf:\uX\ra\uY$ be a morphism of $C^\iy$-schemes. Then by
Definition \ref{ds2def8}, $\uf^*\bigl(T^*\uY\bigr)=f^{-1}(T^*\uY)
\ot_{f^{-1}(\O_Y)}\O_X,$ where $T^*\uY$ is the sheafification of the
presheaf $V\mapsto\Om_{\O_Y(V)}$, and $f^{-1}(T^*\uY)$ the
sheafification of the presheaf $U\mapsto\lim_{V\supseteq f(U)}
(T^*\uY)(V)$, and $f^{-1}(\O_Y)$ the sheafification of the presheaf
$U\mapsto\lim_{V\supseteq f(U)}\O_Y(V)$. The three sheafifications
combine into one, so that $\uf^*\bigl(T^*\uY\bigr)$ is the
sheafification of the presheaf ${\cal P}(\uf^*(T^*\uY))$ acting by
\begin{equation*}
U\longmapsto{\cal P}(\uf^*(T^*\uY))(U)=
\ts\lim_{V\supseteq f(U)}\Om_{\O_Y(V)}\ot_{\O_Y(V)}\O_X(U).
\end{equation*}

Define a morphism of presheaves ${\cal P}\Om_\uf:{\cal
P}(\uf^*(T^*\uY))\ra{\cal P}T^*\uX$ on $X$ by
\begin{equation*}
({\cal P}\Om_\uf)(U)=\ts\lim_{V\supseteq f(U)}
(\Om_{\rho_{f^{-1}(V)\,U}\ci f_\sh(V)})_*,
\end{equation*}
where $(\Om_{\rho_{f^{-1}(V)\,U}\ci f_\sh(V)})_*:\Om_{\O_Y(V)}
\ot_{\O_Y(V)}\O_X(U)\ra\Om_{\O_X(U)}=({\cal P} T^*\uX)(U)$ is
constructed as in Definition \ref{ds2def6} from the $C^\iy$-ring
morphisms $f_\sh(V):\O_Y(V)\ra\O_X(f^{-1}(V))$ from $f_\sh:\O_Y\ra
f_*(\O_X)$ corresponding to $f^\sh$ in $\uf$ as in \eq{ds2eq3}, and
$\rho_{f^{-1}(V)\,U}:\O_X(f^{-1}(V))\ra\O_X(U)$ in $\O_X$. Define
$\Om_\uf:\uf^*\bigl(T^*\uY\bigr)\ra T^*\uX$ to be the induced
morphism of the associated sheaves.
\label{ds2def9}
\end{dfn}

\begin{ex} Let $X$ be a manifold, and $\uX$ the associated
$C^\iy$-scheme. Then $T^*\uX$ is a vector bundle on $\uX$, and is
canonically isomorphic to the lift to $C^\iy$-schemes from Example
\ref{ds2ex5} of the cotangent vector bundle $T^*X$ of~$X$.
\label{ds2ex6}
\end{ex}

Here \cite[Th.~6.17]{Joyc4} are some properties of cotangent
sheaves.

\begin{thm}{\bf(a)} Let\/ $\uf:\uX\ra\uY$ and\/ $\ug:\uY\ra\uZ$ be
morphisms of\/ $C^\iy$-schemes. Then
\begin{equation*}
\Om_{\ug\ci\uf}=\Om_\uf\ci \uf^*(\Om_\ug)\ci I_{\uf,\ug}(T^*\uZ)
\end{equation*}
as morphisms $(\ug\ci\uf)^*(T^*\uZ)\ra T^*\uX$. Here
$\Om_\ug:\ug^*(T^*\uZ)\ra T^*\uY$ is a morphism in $\OYmod,$ so
applying $\uf^*$ gives $\uf^*(\Om_\ug):\uf^*(\ug^*(T^*\uZ))\ra
\uf^*(T^*\uY)$ in\/ $\OXmod,$ and\/
$I_{\uf,\ug}(T^*\uZ):(\ug\ci\uf)^*(T^*\uZ) \ra\uf^*(\ug^*(T^*\uZ))$
is as in Remark\/~{\rm\ref{ds2rem3}}.
\smallskip

\noindent{\bf(b)} Suppose\/ $\uW,\uX,\uY,\uZ$ are locally fair\/
$C^\iy$-schemes with a Cartesian square
\begin{equation*}
\xymatrix@C=60pt@R=12pt{ \uW \ar[r]_\uf \ar[d]^\ue & \uY \ar[d]_\uh \\
\uX \ar[r]^\ug & \uZ}
\end{equation*}
in $\CSchlf,$ so that\/ $\uW=\uX\t_\uZ\uY$. Then the following is
exact in~$\qcoh(\uW)\!:$\I{C-scheme@$C^\iy$-scheme!cotangent
sheaf|)} \I{C-scheme@$C^\iy$-scheme!quasicoherent sheaves
on|)}\I{C-algebraic geometry@$C^\iy$-algebraic geometry|)}
\begin{equation*}
\xymatrix@C=15pt{ (\ug\ci\ue)^*(T^*\uZ)
\ar[rrrr]^(0.45){\begin{subarray}{l}\ue^*(\Om_\ug)\ci
I_{\ue,\ug}(T^*\uZ)\op\\ -\uf^*(\Om_\uh)\ci
I_{\uf,\uh}(T^*\uZ)\end{subarray}} &&&&
\ue^*(T^*\uX)\!\op\!\uf^*(T^*\uY) \ar[rr]^(0.63){\Om_\ue\op\Om_\uf}
&& T^*\uW \ar[r] & 0.}
\end{equation*}
\label{ds2thm2}
\end{thm}

\section{The 2-category of d-spaces}
\label{ds3}
\I{d-space|(}

We will now define the 2-category of {\it d-spaces\/} $\dSpa$,
following \cite[Chap.~2]{Joyc6}. D-spaces are `derived' versions of
$C^\iy$-schemes. In \S\ref{ds4} we will define the 2-category of
d-manifolds $\dMan$ as a 2-subcategory of $\dSpa$. For an
introduction to 2-categories, see~\S\ref{dsA3}--\S\ref{dsA4}.

\subsection{The definition of d-spaces}
\label{ds31}
\I{d-space!definition}\I{2-category|(}

\begin{dfn} A {\it d-space\/} $\bX$ is a quintuple
$\bX=(\uX,\OXp,\EX,\im_X,\jm_X)$ such that $\uX=(X,\OX)$ is a
separated, second countable, locally fair $C^\iy$-scheme, and
$\OXp,\EX,\im_X,\jm_X$ fit into an exact sequence of sheaves on~$X$
\begin{equation*}
\smash{\xymatrix@C=25pt{ \EX \ar[rr]^(0.45){\jm_X} && \OXp
\ar[rr]^(0.55){\im_X} && \OX \ar[r] & 0,}}
\end{equation*}
satisfying the conditions:
\begin{itemize}
\setlength{\itemsep}{0pt}
\setlength{\parsep}{0pt}
\item[(a)] $\OXp$ is a sheaf of $C^\iy$-rings on $X$, with
$\uX'=(X,\OXp)$ a $C^\iy$-scheme.
\item[(b)] $\im_X:\OXp\ra\OX$ is a surjective morphism of sheaves of
$C^\iy$-rings on $X$. Its kernel $\ka_X:\IX\ra\OXp$ is a sheaf
of ideals $\IX$ in $\OXp$, which should be a sheaf of square
zero ideals. Here a {\it square zero ideal\/}\I{square zero
ideal} in a commutative $\R$-algebra $A$ is an ideal $I$ with
$i\cdot j=0$ for all $i,j\in I$. Then $\IX$ is an $\OXp$-module,
but as $\IX$ consists of square zero ideals and $\im_X$ is
surjective, the $\OXp$-action factors through an $\OX$-action.
Hence $\IX$ is an $\OX$-module, and thus a quasicoherent sheaf
on $\uX$, as $\uX$ is locally fair.
\item[(c)] $\EX$ is a quasicoherent sheaf on $\uX$, and
$\jm_X:\EX\ra\IX$ is a surjective morphism in~$\qcoh(\uX)$.
\end{itemize}
As $\uX$ is locally fair, the underlying topological space $X$ is
locally homeomorphic to a closed subset of $\R^n$, so it is {\it
locally compact}. But Hausdorff, second countable and locally
compact imply paracompact, and thus $\uX$ is {\it paracompact}.
\G[WXYZb]{$\bW,\bX,\bY,\bZ,\ldots$}{d-spaces, including d-manifolds}

The sheaf of $C^\iy$-rings $\OXp$ has a sheaf of cotangent modules
$\Om_{\OXp}$, which is an $\OXp$-module with exterior derivative
$\d:\OXp\ra\Om_{\OXp}$. Define $\FX=\Om_{\OXp}\ot_\OXp\OX$ to be the
associated $\OX$-module, a quasicoherent sheaf on $\uX$, and set
$\psi_X=\Om_{\im_X}\ot\id:\FX\ra T^*\uX$, a morphism in
$\qcoh(\uX)$. Define $\phi_X:\EX\ra\FX$ to be the composition of
morphisms of sheaves of abelian groups on~$X$:
\begin{equation*}
\smash{\xymatrix@C=7pt{ \EX  \ar[rr]^{\jm_X} && \IX
\ar[rr]^{\d\vert_{\IX}} && \Om_{\OXp} \ar@{=}[r]^(0.31)\sim &
\Om_{\OXp}\ot_\OXp\OXp \ar[rrr]^{\id\ot\im_X} &&&
\Om_{\OXp}\ot_\OXp\OX \ar@{=}[r] & \FX. }}
\end{equation*}
It turns out that $\phi_X$ is actually a morphism of $\OX$-modules,
and the following sequence is exact in~$\qcoh(\uX)\!:$
\begin{equation*}
\xymatrix@C=20pt{ \EX \ar[rr]^{\phi_X} && \FX \ar[rr]^{\psi_X} &&
T^*\uX \ar[r] & 0.}
\end{equation*}
The morphism $\phi_X:\EX\ra\FX$ will be called the {\it virtual
cotangent sheaf\/} of $\bX$, for reasons we explain
in~\S\ref{ds43}.\I{d-space!virtual cotangent sheaf}
\I{d-space!1-morphism}

Let $\bX,\bY$ be d-spaces. A 1-{\it
morphism\/}\I{2-category!1-morphism} $\bs f:\bX\ra\bY$ is a triple
$\bs f=(\uf,f',f'')$, where $\uf=(f,f^\sh):\uX\ra\uY$ is a morphism
of $C^\iy$-schemes, $f':f^{-1}(\OYp)\ra\OXp$ a morphism of sheaves
of $C^\iy$-rings on $X$, and $f'':\uf^*(\EY)\ra\EX$ a morphism in
$\qcoh(\uX)$, such that the following diagram of sheaves on $X$
commutes:
\begin{equation*}
\xymatrix@C=11pt@R=1pt{
f^{-1}\!(\EY)\ot_{f^{-1}(\OY)}^{\id}\!f^{-1}\!(\OY) \ar@{=}[r]
\ar[dd]^(0.4){{}\,\id\ot f^\sh} & f^{-1}\!(\EY)
\ar[rr]_(0.45){\raisebox{-9pt}{$\scriptstyle f^{-1}(\jm_Y)$}}
 && f^{-1}\!(\OYp)
\ar[rr]_(0.52){\raisebox{-9pt}{$\scriptstyle f^{-1}(\im_Y)$}}
\ar[ddd]^{f'} && f^{-1}\!(\OY) \ar[r] \ar[ddd]_{f^\sh} & 0 \\ \\
{\begin{subarray}{l}\ts \uf^*(\EY)=\\
\ts f^{-1}(\EY) \ot_{f^{-1}(\OY)}^{f^\sh}\OX\end{subarray}}
\ar[dr]^(0.7){f''} \\  & \EX \ar[rr]^{\jm_X} && \OXp
\ar[rr]^(0.55){\im_X} && \OX \ar[r] &  {0.\kern -.28em} }
\end{equation*}
Define morphisms $f^2=\Om_{f'}\ot\id:\uf^*(\FY)\ra\FX$ and
$f^3=\Om_\uf:\uf^*(T^*\uY)\ra T^*\uX$ in $\qcoh(\uX)$. Then the
following commutes in $\qcoh(\uX)$, with exact rows:
\e
\begin{gathered}
\xymatrix@C=20pt@R=12pt{ \uf^*(\EY) \ar[rr]_{\uf^*(\phi_Y)}
\ar[d]^{f''} && \uf^*(\FY) \ar[rr]_{\uf^*(\psi_Y)} \ar[d]^{f^2} &&
\uf^*(T^*\uY) \ar[r] \ar[d]^{f^3} & 0 \\
\EX \ar[rr]^{\phi_X} && \FX \ar[rr]^{\psi_X} && T^*\uX \ar[r] &
{0.\kern -.28em} }
\end{gathered}
\label{ds3eq1}
\e

If $\bX$ is a d-space, the {\it identity $1$-morphism\/}
$\bs\id_\bX:\bX\ra\bX$ is $\bs\id_\bX=\bigl(\uid_\uX,
\de_X(\OXp),\de_\uX(\EX)\bigr)$, where $\de_X(*)$ are the canonical
isomorphisms of Remark \ref{ds2rem3}. Let $\bX,\bY,\bZ$ be d-spaces,
and $\bs f:\bX\ra\bY$, $\bs g:\bY\ra\bZ$ be 1-morphisms. Define the
{\it composition\I{2-category!1-morphism!composition} of\/
$1$-morphisms\/} $\bs g\ci\bs f:\bX\ra\bZ$ to be
\e
\bs g\ci\bs f=\bigl(\ug\ci\uf,f'\ci f^{-1}(g')\ci I_{f,g}(\OZp),
f''\ci\uf^*(g'')\ci I_{\uf,\ug}(\EZ)\bigr),
\label{ds3eq2}
\e
where $I_{*,*}(*)$ are the canonical isomorphisms of
Remark~\ref{ds2rem3}.

Let $\bs f,\bs g:\bX\ra\bY$ be 1-morphisms of d-spaces, where $\bs
f=(\uf,f',f'')$ and $\bs g=(\ug,g',g'')$. Suppose $\uf=\ug$. A
2-{\it morphism\/}\I{d-space!2-morphism}\I{2-category!2-morphism}
$\eta:\bs f\Ra\bs g$ is a morphism $\eta:\uf^*(\FY)\ra\EX$ in
$\qcoh(\uX)$, such that
\begin{gather*}
g'=f'+\jm_X\ci\eta\ci\bigl(\id\ot(f^\sh\ci
f^{-1}(\im_Y))\bigr) \ci\bigl(f^{-1}(\d)\bigr)\\
\text{and}\qquad g''=f''+\eta\ci\uf^*(\phi_Y).
\end{gather*}
Then $g^2=f^2+\phi_X\ci\eta$ and $g^3=f^3$, so \eq{ds3eq1} for $\bs
f,\bs g$ combine to give a diagram
\e
\begin{gathered}
\xymatrix@C=30pt@R=17pt{ \uf^*(\EY) \ar[rr]^{\uf^*(\phi_Y)}
\ar@<-.8ex>[d]_(0.3){f''}
\ar@<.2ex>[d]^(0.3){g''=f''+\eta\ci\uf^*(\phi_Y)} && \uf^*(\FY)
\ar[rr]^{\uf^*(\psi_Y)} \ar@<-.5ex>[d]_(0.45){f^2}
\ar@<.5ex>[d]^(0.45){g^2=f^2+\phi_X\ci\eta}
\ar@<.5ex>@{.>}[dll]^(0.3)\eta
 && \uf^*(T^*\uY) \ar[r] \ar[d]^{f^3=g^3} & 0 \\
\EX \ar@<-.2ex>[rr]^(0.6){\phi_X} && \FX \ar@<-.2ex>[rr]^{\psi_X} &&
T^*\uX \ar@<-.2ex>[r] & {0.\!\!} }
\end{gathered}\!\!\!\!\!
\label{ds3eq3}
\e
That is, $\eta$ is a homotopy between the morphisms of complexes
\eq{ds3eq1} from~$\bs f,\bs g$.

If $\bs f:\bX\ra\bY$ is a 1-morphism, the {\it identity\/
$2$-morphism\/} $\id_{\bs f}:\bs f\Ra\bs f$ is the zero morphism
$0:\uf^*(\FY)\ra\EX$. Suppose $\bX,\bY$ are d-spaces, $\bs f,\bs
g,\bs h:\bX\ra\bY$ are 1-morphisms and $\eta:\bs f\Ra\bs g$,
$\ze:\bs g\Ra\bs h$ are 2-morphisms. The {\it vertical composition
of\/ $2$-morphisms\/}\I{2-category!2-morphism!vertical composition}
$\ze\od\eta:\bs f\Ra\bs h$ as in \eq{dsAeq1}
is~$\ze\od\eta=\ze+\eta$.

Let $\bX,\bY,\bZ$ be d-spaces, $\bs f,\bs{\ti f}:\bX\ra\bY$ and $\bs
g,\bs{\ti g}:\bY\ra\bZ$ be 1-morphisms, and $\eta:\bs f\Ra\bs{\ti
f}$, $\ze:\bs g\Ra\bs{\ti g}$ be 2-morphisms. The {\it horizontal
composition\I{2-category!2-morphism!horizontal composition} of\/
$2$-morphisms\/} $\ze*\eta:\bs g\ci\bs f\Ra\bs{\ti g}\ci\bs{\ti f}$
as in \eq{dsAeq2} is
\begin{equation*}
\ze*\eta=\bigl(\eta\ci\uf^*(g^2)+f''\ci\uf^*(\ze)+\eta\ci\uf^*(\phi_Y)
\ci\uf^*(\ze)\bigr)\ci I_{\uf,\ug}(\FZ).
\end{equation*}
This completes the definition of the 2-category of d-spaces~$\dSpa$.
\G[dSpa]{$\dSpa$}{2-category of d-spaces}\I{2-category|)}

Regard the category $\CSchlfssc$\G[CSchlfssc]{$\CSchlfssc$}{category
of separated, second countable, locally fair $C^\iy$-schemes} of
separated, second countable, locally fair $C^\iy$-schemes as a
2-category with only identity 2-morphisms $\id_{\uf}$ for
(1-)mor\-phisms $\uf:\uX\ra\uY$. Define a 2-functor
$F_\CSch^\dSpa:\CSchlfssc\ra\dSpa$ to map $\uX$ to
$\bX=(\uX,\OX,0,\id_{\OX},0)$ on objects $\uX$, to map $\uf$ to $\bs
f=(\uf,f^\sh,0)$ on (1-)morphisms $\uf:\uX\ra\uY$, and to map
identity 2-morphisms $\id_\uf:\uf\Ra\uf$ to identity 2-morphisms
$\id_{\bs f}:\bs f\Ra\bs f$. Define a 2-functor
$F_\Man^\dSpa:\Man\ra\dSpa$ by~$F_{\Man}^\dSpa=F_\CSch^\dSpa\ci
F_\Man^\CSch$.

Write $\hCSchlfssc$ for the full 2-subcategory of objects $\bX$ in
$\dSpa$ equivalent to $F_\CSch^\dSpa(\uX)$ for some $\uX$ in
$\CSchlfssc$, and $\hMan$\G[Man']{$\hMan$}{2-subcategory of d-spaces
equivalent to manifolds} for the full 2-subcategory of objects $\bX$
in $\dSpa$ equivalent to $F_\Man^\dSpa(X)$ for some manifold $X$.
When we say that a d-space $\bX$ {\it is a $C^\iy$-scheme}, or {\it
is a manifold}, we mean that $\bX\in\hCSchlfssc$, or $\bX\in\hMan$,
respectively.\I{d-space!is a $C^\iy$-scheme}\I{d-space!is a
manifold}
\label{ds3def1}
\end{dfn}

In \cite[\S 2.2]{Joyc6} we prove:

\begin{thm}{\bf(a)} Definition\/ {\rm\ref{ds3def1}} defines a
strict\/ $2$-category $\dSpa,$ in which all\/ $2$-morphisms are
$2$-isomorphisms.\I{2-category}
\smallskip

\noindent{\bf(b)} For any $1$-morphism $\bs f:\bX\ra\bY$ in\/
$\dSpa$ the\/ $2$-morphisms $\eta:\bs f\Ra\bs f$ form an abelian
group under vertical composition, and in fact a real vector space.
\smallskip

\noindent{\bf(c)} $F_\CSch^\dSpa$ and\/ $F_\Man^\dSpa$ in
Definition\/ {\rm\ref{ds3def1}} are full and
faithful\I{functor!full}\I{functor!faithful} strict\/ $2$-functors.
Hence $\CSchlfssc,\Man$ and\/ $\hCSchlfssc,\hMan$ are equivalent\/
$2$-categories.
\label{ds3thm1}
\end{thm}

\begin{rem}{\bf(i)} One should think of a d-space
$\bX=(\uX,\OXp,\EX,\im_X, \jm_X)$ as being a $C^\iy$-scheme $\uX$,
which is the `classical' part of $\bX$ and lives in a category
rather than a 2-category, together with some extra `derived'
information $\OXp,\EX,\im_X,\jm_X$. 2-morphisms in $\dSpa$ are
wholly to do with this derived part. The sheaf $\EX$ may be thought
of as a (dual) `obstruction sheaf' on~$\uX$.

\smallskip

\noindent{\bf(ii)} Readers familiar with derived algebraic
geometry\I{derived algebraic geometry|(} may find the following
(oversimplified) explanation of d-spaces helpful; more details are
given in \cite[\S 14.4]{Joyc6}. In conventional algebraic geometry,
a $\K$-scheme $(X,\O_X)$ is a topological space $X$ equipped with a
sheaf of $\K$-algebras $\O_X$. In derived algebraic geometry, as in
To\"en and Vezzosi \cite{ToVe1,ToVe2} and Lurie
\cite{Luri1,Luri2,Luri3}, a {\it derived\/ $\K$-scheme\/}\I{derived
scheme|(} $(X,\O_X)$ is (roughly) a topological space $X$ with a
(homotopy) sheaf of (commutative) dg-algebras over $\K$. Here a
({\it commutative\/}) {\it dg-algebra\/} $(A_*,\d)$ is a
nonpositively graded $\K$-algebra $\bigop_{k\le 0}A_k$, with
differentials $\d:A_k\ra A_{k+1}$ satisfying $\d^2=0$ and
$ab=(-1)^{kl}ba$, $\d(ab)=(\d a)b+(-1)^ka(\d b)$ for all $a\in A_k$
and~$b\in A_l$.\I{dg-algebra}

We call a dg-algebra $(A_*,\d)$ {\it square
zero\/}\I{dg-algebra!square zero} if $A_k=0$ for $k\ne 0,-1$ and
$A_{-1}\cdot \d(A_{-1})=0$. This implies that $\d(A_{-1})$ is a
square zero ideal\I{square zero ideal} in $A_0$. General dg-algebras
form an $\iy$-category,\I{$\iy$-category} but square zero
dg-algebras form a 2-category. Ignoring $C^\iy$-rings for the
moment, we can think of the data
$\EX\,{\buildrel\jm_X\over\longra}\,\OXp$ in a d-space $\bX$ as a
sheaf of square zero dg-algebras $\smash{A_{-1}\,{\buildrel\d\over
\longra}\,A_0}$ on $X$. The remaining data $\OX,\im_X$ can be
recovered from this, since $\OXp\,{\buildrel\im_X\over\longra}\,\OX$
is the cokernel of $\EX\,{\buildrel\jm_X\over\longra}\,\OXp$. Thus,
a d-space $\bX$ is like a special kind of derived $\R$-scheme, in
which the dg-algebras are all square zero.\I{derived algebraic
geometry|)}\I{derived scheme|)}
\label{ds3rem1}
\end{rem}

\subsection{Gluing d-spaces by equivalences}
\label{ds32}

Next we discuss gluing of d-spaces and 1-morphisms on open
d-subspaces.

\begin{dfn} Let $\bX=(\uX,\OXp,\EX,\im_X,\jm_X)$ be a d-space.
Suppose $\uU\subseteq\uX$ is an open $C^\iy$-subscheme. Then $\bU=
\bigl(\uU,\OXp\vert_U,\EX\vert_U,\im_X\vert_U,\jm_X\vert_U\bigr)$ is
a d-space. We call $\bU$ an {\it open d-subspace\/} of $\bX$. An
{\it open cover\/} of a d-space $\bX$ is a family $\{\bU_a:a\in A\}$
of open d-subspaces $\bU_a$ of $\bX$ with~$\uX=\bigcup_{a\in
A}\uU_a$.\I{d-space!open d-subspace}\I{d-space!open cover}
\label{ds3def2}
\end{dfn}

As in \cite[\S 2.4]{Joyc6}, we can glue 1-morphisms on open
d-subspaces which are 2-isomorphic on the overlap. The proof uses
partitions of unity, as in~\S\ref{ds22}.

\begin{prop} Let\/ $\bX,\bY$ be d-spaces, $\bU,\bV\subseteq\bX$
be open d-subspaces with\/ $\bX=\bU\cup\bV,$ $\bs f:\bU\ra\bY$ and\/
$\bs g:\bV\ra\bY$ be $1$-morphisms, and\/ $\eta:\bs
f\vert_{\bU\cap\bV}\Ra\bs g\vert_{\bU\cap\bV}$ a $2$-morphism. Then
there exist a $1$-morphism $\bs h:\bX\ra\bY$ and\/ $2$-morphisms
$\ze:\bs h\vert_\bU\Ra\bs f,$ $\th:\bs h\vert_\bV\Ra\bs g$ such
that\/ $\th\vert_{\bU\cap\bV}=
\eta\od\ze\vert_{\bU\cap\bV}:\bs h\vert_{\bU\cap\bV}\Ra \bs
g\vert_{\bU\cap\bV}$. This $\bs h$ is unique up to $2$-isomorphism,
and independent up to $2$-isomorphism of the choice of\/~$\eta$.
\label{ds3prop}
\end{prop}

{\it Equivalences\/}\I{d-space!equivalence}\I{d-space!gluing by
equivalences|(} $\bs f:\bX\ra\bY$ in a 2-category are defined in
\S\ref{dsA3}, and are the natural notion of when two objects
$\bX,\bY$ are `the same'. In \cite[\S 2.4]{Joyc6} we prove theorems
on gluing d-spaces by equivalences. See Spivak \cite[Lem.~6.8 \&
Prop.~6.9]{Spiv} for results similar to Theorem \ref{ds3thm2} for
his `local $C^\iy$-ringed spaces', an
$\iy$-categorical\I{$\iy$-category} analogue of our d-spaces.

\begin{thm} Suppose $\bX,\bY$ are d-spaces, $\bU\subseteq\bX,$
$\bV\subseteq\bY$ are open d-subspaces, and\/ $\bs f:\bU\ra\bV$ is
an equivalence in $\dSpa$. At the level of topological spaces, we
have open $U\subseteq X,$ $V\subseteq Y$ with a homeomorphism
$f:U\ra V,$ so we can form the quotient topological space
$Z:=X\amalg_fY=(X\amalg Y)/\sim,$ where the equivalence relation
$\sim$ on $X\amalg Y$ identifies $u\in U\subseteq X$ with\/~$f(u)\in
V\subseteq Y$.

Suppose $Z$ is Hausdorff. Then there exist a d-space $\bZ$ with
topological space $Z,$ open d-subspaces $\bs{\hat X},\bs{\hat Y}$ in
$\bZ$ with\/ $\bZ=\bs{\hat X}\cup\bs{\hat Y},$ equivalences $\bs
g:\bX\ra\bs{\hat X}$ and\/ $\bs h:\bY\ra\bs{\hat Y}$ in $\dSpa$ such
that\/ $\bs g\vert_\bU$ and\/ $\bs h\vert_\bV$ are both equivalences
with\/ $\bs{\hat X}\cap\bs{\hat Y},$ and a $2$-morphism $\eta:\bs
g\vert_\bU\Ra\bs h\ci\bs f:\bU\ra\bs{\hat X}\cap\bs{\hat Y}$.
Furthermore, $\bZ$ is independent of choices up to equivalence.
\label{ds3thm2}
\end{thm}

In Theorem \ref{ds3thm2}, $\bZ$ is a {\it
pushout\/}\I{2-category!pushout}\I{pushout}
$\bX\amalg_{\bs\id_\bU,\bU,\bs f}\bY$ in the 2-category~$\dSpa$.

\begin{thm} Suppose\/ $I$ is an indexing set, and\/ $<$ is a total order
on $I,$ and\/ $\bX_i$ for $i\in I$ are d-spaces, and for all\/ $i<j$
in $I$ we are given open d-subspaces $\bU_{ij}\subseteq\bX_i,$
$\bU_{ji}\subseteq\bX_j$ and an equivalence $\bs
e_{ij}:\bU_{ij}\ra\bU_{ji},$ such that for all\/ $i<j<k$ in $I$ we
have a $2$-commutative diagram
\begin{equation*}
\xymatrix@C=70pt@R=10pt{ & \bU_{ji}\cap\bU_{jk} \ar@<.5ex>[dr]^{\bs
e_{jk}\vert_{\bU_{ji}\cap\bU_{jk}}} \ar@{=>}[d]^{\eta_{ijk}} \\
\bU_{ij}\cap\bU_{ik} \ar@<.5ex>[ur]^{\bs
e_{ij}\vert_{\bU_{ij}\cap\bU_{ik}}} \ar@<-.25ex>[rr]^(0.37){\bs
e_{ik}\vert_{\bU_{ij}\cap\bU_{ik}}} && \bU_{ki}\cap\bU_{kj}}
\end{equation*}
for some $\eta_{ijk},$ where all three $1$-morphisms are
equivalences.

On the level of topological spaces, define the quotient topological
space $Y=(\coprod_{i\in I}X_i)/\sim,$ where $\sim$ is the
equivalence relation generated by $x_i\sim x_j$ if\/ $i<j,$ $x_i\in
U_{ij}\subseteq X_i$ and\/ $x_j\in U_{ji}\subseteq X_j$ with\/
$e_{ij}(x_i)=x_j$. Suppose $Y$ is Hausdorff and second countable.
Then there exist a d-space $\bY$ and a $1$-morphism $\bs
f_i:\bX_i\ra\bY$ which is an equivalence with an open d-subspace
$\bs{\hat X}_i\subseteq\bY$ for all\/ $i\in I,$ where
$\bY=\bigcup_{i\in I}\bs{\hat X}_i,$ such that\/ $\bs
f_i\vert_{\bU_{ij}}$ is an equivalence $\bU_{ij}\ra\bs{\hat
X}_i\cap\bs{\hat X}_j$ for all\/ $i<j$ in $I,$ and there exists a
$2$-morphism\/ $\eta_{ij}:\bs f_j\ci\bs e_{ij}\Ra\bs
f_i\vert_{\bU_{ij}}$. The d-space $\bY$ is unique up to equivalence,
and is independent of choice of\/ $2$-morphisms\/~$\eta_{ijk}$.

Suppose also that\/ $\bZ$ is a d-space, and\/ $\bs g_i:\bX_i\ra\bZ$
are $1$-morphisms for all\/ $i\in I,$ and there exist\/
$2$-morphisms $\ze_{ij}:\bs g_j\ci\bs e_{ij}\Ra\bs
g_i\vert_{\bU_{ij}}$ for all\/ $i<j$ in $I$. Then there exist a
$1$-morphism $\bs h:\bY\ra\bZ$ and\/ $2$-morphisms $\ze_i:\bs
h\ci\bs f_i\Ra\bs g_i$ for all\/ $i\in I$. The $1$-morphism $\bs h$
is unique up to $2$-isomorphism, and is independent of the choice
of\/ $2$-morphisms~$\ze_{ij}$.
\label{ds3thm3}
\end{thm}

\begin{rem} In Proposition \ref{ds3prop}, it is surprising that
$\bs h$ is independent of $\eta$ up to $2$-isomorphism. It holds
because of the existence of partitions of unity on nice
$C^\iy$-schemes, as in Proposition \ref{ds2prop1}. Here is a sketch
proof: suppose $\eta,\bs h,\ze,\th$ and $\eta',\bs h',\ze',\th'$ are
alternative choices in Proposition \ref{ds3prop}. Then we have
2-morphisms $(\ze')^{-1}\od\ze:\bs h\vert_\bU\Ra \bs h'\vert_\bU$
and $(\th')^{-1}\od\th:\bs h\vert_\bV\Ra \bs h'\vert_\bV$. Choose a
partition of unity\I{partition of unity} $\{\al,1-\al\}$ on $\uX$
subordinate to $\{\uU,\uV\}$, so that $\al:\uX\ra\R$ is smooth with
$\al$ supported on $\uU\subseteq\uX$ and $1-\al$ supported on
$\uV\subseteq\uX$. Then
$\al\cdot\bigl((\ze')^{-1}\od\ze\bigr)+(1-\al)
\cdot\bigl((\th')^{-1} \od\th\bigr)$ is a 2-morphism $\bs h\Ra\bs
h'$, where $\al\cdot\bigl((\ze')^{-1}\od\ze\bigr)$ makes sense on
all of $\uX$ (rather than just on $\uU$ where $(\ze')^{-1}\od\ze$ is
defined) as $\al$ is supported on $\uU$, so we extend by zero
on~$\uX\sm\uU$.

Similarly, in Theorem \ref{ds3thm3}, the compatibility conditions on
the gluing data $\bX_i,\bU_{ij},\bs e_{ij}$ are significantly weaker
than you might expect, because of the existence of partitions of
unity. The 2-morphisms $\eta_{ijk}$ on overlaps $\bs{\hat
X}_i\cap\bs{\hat X}_j\cap\bs{\hat X}_k$ are only required to exist,
not to satisfy any further conditions. In particular, one might
think that on overlaps $\bs{\hat X}_i\cap\bs{\hat X}_j\cap\bs{\hat
X}_k \cap\bs{\hat X}_l$ we should require
\e
\eta_{ikl}\od(\id_{{\bs f}_{kl}}*\eta_{ijk})\vert_{\bU_{ij}\cap
\bU_{ik}\cap\bU_{il}}=\eta_{ijl}\od (\eta_{jkl}*\id_{{\bs f}_{ij}})
\vert_{\bU_{ij}\cap \bU_{ik}\cap\bU_{il}},
\label{ds3eq4}
\e
but we do not. Also, one might expect the $\ze_{ij}$ should satisfy
conditions on triple overlaps $\bs{\hat X}_i\cap\bs{\hat
X}_j\cap\bs{\hat X}_k$, but they need not.

The moral is that constructing d-spaces by gluing together patches
$\bX_i$ is straightforward, as one only has to verify mild
conditions on triple overlaps $\bX_i\cap\bX_j\cap\bX_k$. Again, this
works because of the existence of partitions of unity on nice
$C^\iy$-schemes, which are used to construct the glued d-spaces
$\bZ$ and 1- and 2-morphisms in Theorems \ref{ds3thm2}
and~\ref{ds3thm3}.

In contrast, for gluing d-stacks in \S\ref{ds103}, we do need
compatibility conditions of the form \eq{ds3eq4}.\I{d-stack!gluing
by equivalences} The problem of gluing geometric spaces in an
$\iy$-category\I{$\iy$-category} $\bs{\mathcal{C}}$ by equivalences,
such as Spivak's derived manifolds \cite{Spiv},\I{Spivak's derived
manifolds} is discussed by To\"en and Vezzosi \cite[\S 1.3.4]{ToVe1}
and Lurie \cite[\S 6.1.2]{Luri1}. It requires nontrivial conditions
on overlaps $\bcX_{i_1}\cap\cdots\cap\bcX_{i_n}$ for
all~$n=2,3,\ldots.$\I{d-space!gluing by equivalences|)}
\label{ds3rem2}
\end{rem}

\subsection{Fibre products in $\dSpa$}
\label{ds33}
\I{d-space!fibre products|(}\I{fibre product!of
d-spaces}\I{2-category!fibre products in}

Fibre products in 2-categories are explained in \S\ref{dsA4}. In
\cite[\S 2.5--\S 2.6]{Joyc6} we discuss fibre products in $\dSpa$,
and their relation to transverse fibre products in~$\Man$.

\begin{thm}{\bf(a)} All fibre products exist in the
$2$-category~$\dSpa$.
\smallskip

\noindent{\bf(b)} Let\/ $g:X\ra Z$ and\/ $h:Y\ra Z$ be smooth maps
of manifolds, and write\/ $\bX=F_\Man^\dSpa(X),$ and similarly for
$\bY,\bZ,\bs g,\bs h$. If\/ $g,h$ are transverse, so that a fibre
product $X\t_{g,Z,h}Y$ exists in $\Man,$ then the fibre product\/
$\bX\t_{\bs g,\bZ,\bs h}\bY$ in $\dSpa$ is equivalent in $\dSpa$
to\/ $F_\Man^\dSpa(X\t_{g,Z,h}Y)$. If\/ $g,h$ are not transverse
then $\bX\t_{\bs g,\bZ,\bs h}\bY$ exists in $\dSpa,$ but is not a
manifold.\I{manifold!transverse fibre products}
\label{ds3thm4}
\end{thm}

To prove (a), given 1-morphisms $\bs g:\bX\ra\bZ$ and $\bs
h:\bY\ra\bZ$, we write down an explicit d-space
$\bW=(\uW,\OWp,\EW,\im_W,\jm_W)$, 1-morphisms $\bs
e=(\ue,e',e''):\bW\ra\bX$ and $\bs f=(\uf,f',f''):\bW\ra\bY$ and a
2-morphism $\eta:\bs g\ci\bs e\Ra\bs h\ci\bs f$, and verify the
universal property for
\begin{equation*}
\xymatrix@C=60pt@R=10pt{ \bW \ar[r]_(0.25){\bs f} \ar[d]^{\bs e}
\drtwocell_{}\omit^{}\omit{^{\eta}}
 & \bY \ar[d]_{\bs h} \\ \bX \ar[r]^(0.7){\bs g} & \bZ}
\end{equation*}
to be a 2-Cartesian square\I{2-category!2-Cartesian square} in
$\dSpa$. The underlying $C^\iy$-scheme $\uW$ is the fibre product
$\uW=\uX\t_{\ug,\uZ,\uh}\uY$ in $\CSch$, and $\ue:\uW\ra\uX$,
$\uf:\uW\ra\uY$ are the projections from the fibre product. The
definitions of $\OWp,\im_W,\jm_W,e',f'$ in \cite[\S 2.5]{Joyc6} are
complex, and we will not give them here. The remaining data
$\EW,e'',f'',\eta$, as well as the virtual cotangent sheaf
$\phi_W:\EW\ra\FW$, is characterized by the following commutative
diagram in $\qcoh(\uW)$, with exact top row:
\begin{equation*}
\xymatrix@C=15pt@R=30pt{ (\ug\ci\ue)^*(\EZ)
\ar[rrrr]^(0.47){\begin{pmatrix}\st \ue^*(g'')\ci
I_{\ue,\ug}(\EZ) \\ \st -\uf^*(h'')\ci I_{\uf,\uh}(\EZ) \\
\st (\ug\ci\ue)^*(\phi_Z)\end{pmatrix}} &&&&
{\vphantom{\EW}\smash{\raisebox{11pt}{$\begin{subarray}{l}
\ts \ue^*(\EX)\op \\ \ts \uf^*(\EY)\op \\ \ts (\ug\ci\ue)^*(\FZ)
\end{subarray}$}}} \ar[d]_(0.45){\begin{pmatrix}
\st -\ue^*(\phi_X)\!{} & \st 0\!{} & \st \ue^*(g^2)\ci I_{\ue,\ug}(\FZ) \\
\st 0\!{} & \st -\uf^*(\phi_Y)\!{} & \st -\uf^*(h^2)\ci
I_{\uf,\uh}(\FZ)
\end{pmatrix}}
\ar@{.>}[rrr]^(0.55){\begin{pmatrix} \st e''\!\!\!\!{} & \st f''\!\!\!\!{}
& \st \eta \end{pmatrix}} &&& \cE_W \ar[r] \ar@{.>}[d]^{\phi_W} & 0
\\
&&&& {\begin{subarray}{l} \ts\ue^*(\FX)\op \\
\ts \uf^*(\FY)\end{subarray}}
\ar@{.>}[rrr]^(0.55){\begin{pmatrix} \st e^2\!\!\!\!{} & \st f^2
\end{pmatrix}}_(0.55)\cong
&&& \FW. }
\end{equation*}
\I{d-space!fibre products|)}

\subsection{Fixed point loci of finite groups in d-spaces}
\label{ds34}
\I{d-space!fixed point loci|(}

If a finite group $\Ga$ acts on a manifold $X$ by diffeomorphisms,
then the fixed point locus $X^\Ga$ is a disjoint union of closed,
embedded submanifolds of $X$. In a similar way, if $\Ga$ acts on a
d-space $\bX$ by 1-isomorphisms, in \cite[\S 2.7]{Joyc6} we define a
d-space $\bX^\Ga$ called the {\it fixed d-subspace of\/} $\Ga$ in
$\bX$, with an inclusion 1-morphism $\bs
j_{\bX,\Ga}:\bX^\Ga\hookra\bX$, whose topological space $X^\Ga$ is
the fixed point locus of $\Ga$ in~$X$.\G[XGaa]{$\bX^\Ga$}{fixed
d-subspace of group $\Ga$ acting on a d-space $\bX$}\G[jXGaa]{$\bs
j_{\bX,\Ga}:\bX^\Ga\hookra \bX$}{inclusion of $\Ga$-fixed d-subspace
$\bX^\Ga$ in a d-space $\bX$}

Note that by an {\it action\/ $\bs r:\Ga\ra\Aut(\bX)$ of\/ $\Ga$
on\/} $\bX$ we shall always mean a {\it strict\/} action, that is,
$\bs r(\ga):\bX\ra\bX$ is a 1-isomorphism for all $\ga\in\Ga$ and
$\bs r(\ga\de)=\bs r(\ga)\bs r(\de)$ for all $\ga,\de\in\Ga$, rather
than $\bs r(\ga\de)$ only being 2-isomorphic to $\bs r(\ga)\bs
r(\de)$. The next theorem summarizes our results.

\begin{thm} Let\/ $\bX$ be a d-space, $\Ga$ a finite group, and\/
$\bs r:\Ga\ra\Aut(\bX)$ an action of\/ $\Ga$ on\/ $\bX$ by
$1$-isomorphisms. Then we can define a d-space\/ $\bX^\Ga$ called
the \begin{bfseries}fixed d-subspace of\/ $\Ga$ in\/\end{bfseries}
$\bX,$ with an inclusion $1$-morphism\/ $\bs j_{\bX,\Ga}:\bX^\Ga
\ra\bX$. It has the following properties:
\begin{itemize}
\setlength{\itemsep}{0pt}
\setlength{\parsep}{0pt}
\item[{\bf(a)}] Let\/ $\bX,\Ga,\bs r$ and\/ $\bs
j_{\bX,\Ga}:\bX^\Ga\ra\bX$ be as above. Suppose $\bs
f:\bW\ra\bX$ is a $1$-morphism in $\dSpa$. Then $\bs f$
factorizes as $\bs f=\bs j_{\bX,\Ga}\ci\bs g$ for some
$1$-morphism $\bs g:\bW\ra\bX^\Ga$ in $\dSpa,$ which must be
unique, if and only if\/ $\bs r(\ga)\ci\bs f=\bs f$ for
all\/~$\ga\in\Ga$.

\item[{\bf(b)}] Suppose $\bX,\bY$ are d-spaces, $\Ga$ is a finite
group, $\bs r:\Ga\ra\Aut(\bX),$ $\bs s:\Ga\ra\Aut(\bY)$ are
actions of\/ $\Ga$ on $\bX,\bY,$ and\/ $\bs f:\bX\ra\bY$ is a
$\Ga$-equivariant\/ $1$-morphism in $\dSpa,$ that is, $\bs
f\ci\bs r(\ga)=\bs s(\ga)\ci\bs f$ for\/ $\ga\in\Ga$. Then there
exists a unique\/ $1$-morphism $\bs f^\Ga:\bX^\Ga\ra\bY^\Ga$
such that\/~$\bs j_{\bY,\Ga}\ci\bs f^\Ga=\bs f\ci\bs
j_{\bX,\Ga}$.
\item[{\bf(c)}] Let\/ $\bs f,\bs g:\bX\ra\bY$ be
$\Ga$-equivariant\/ $1$-morphisms as in {\bf(b)\rm,} and\/
$\eta:\bs f\Ra\bs g$ be a $\Ga$-equivariant\/ $2$-morphism, that
is, $\eta*\id_{\bs r(\ga)}=\id_{\bs s(\ga)}*\eta$ for\/
$\ga\in\Ga$. Then there exists a unique $2$-morphism
$\eta^\Ga:\bs f^\Ga\Ra\bs g^\Ga$ such that\/~$\id_{\bs
j_{\bY,\Ga}}*\eta^\Ga=\eta*\id_{\bs j_{\bX,\Ga}}$.
\end{itemize}
Note that\/ {\bf(a)} is a universal property that determines
$\bX^\Ga,\bs j_{\bX,\Ga}$ up to canonical\/ $1$-isomorphism.
\label{ds3thm5}
\end{thm}

We will use fixed d-subspaces $\bX^\Ga$ in Theorem \ref{ds10thm9}
below to describe orbifold strata\I{orbifold strata!of
d-stacks}\I{d-stack!orbifold strata}\I{quotient
d-stack}\I{d-stack!quotients $[\bX/G]$} $\bcX^\Ga$ of quotient
d-stacks $\bcX=[\bX/G]$. If $\bX$ is a d-manifold, as in
\S\ref{ds4}, then in general the fixed d-subspaces $\bX^\Ga$ are
disjoint unions of d-manifolds of different
dimensions.\I{d-space|)}\I{d-space!fixed point loci|)}

\section{The 2-category of d-manifolds}
\label{ds4}
\I{d-manifold|(}

We can now define and discuss {\it d-manifolds}, our derived version
of smooth manifolds (without boundary), following~\cite[Chap.s 3 \&
4]{Joyc6}.

\subsection{The definition of d-manifolds}
\label{ds41}
\I{principal d-manifold|(}\I{d-manifold!principal|(}

\begin{dfn} A d-space $\bU$ is called a {\it principal
d-manifold\/} if is equivalent in $\dSpa$ to a fibre product
$\bX\t_{\bs g,\bZ,\bs h}\bY$ with $\bX,\bY,\bZ\in\hMan$. That is,
\begin{equation*}
\bU\simeq F_\Man^\dSpa(X)\t_{F_\Man^\dSpa(g),F_\Man^\dSpa(Z),
F_\Man^\dSpa(h)}F_\Man^\dSpa(Y)
\end{equation*}
for manifolds $X,Y,Z$ and smooth maps $g:X\ra Z$ and $h:Y\ra Z$. The
{\it virtual dimension\/} $\vdim\bU$ of $\bU$ is defined to be
$\vdim\bU=\dim X+\dim Y-\dim Z$. Proposition \ref{ds4prop2}(b) below
shows that if $\bU\ne\bs\es$ then $\vdim\bU$ depends only on the
d-space $\bU$, and not on the choice of $X,Y,Z,g,h$, and so is well
defined.\I{d-manifold!virtual dimension}

A d-space $\bW$ is called a {\it d-manifold of virtual dimension\/}
$n\in\Z$, written $\vdim\bW=n$, if $\bW$ can be covered by nonempty
open d-subspaces $\bU$ which are principal d-manifolds
with~$\vdim\bU=n$.

Write $\dMan$\G[dMan]{$\dMan$}{2-category of d-manifolds} for the
full 2-subcategory of d-manifolds in $\dSpa$. If $\bX\in\hMan$ then
$\bX\simeq\bX\t_{\bs *}\bs *$, so $\bX$ is a principal d-manifold,
and thus a d-manifold. Therefore $\hMan$ in \S\ref{ds31} is a
2-subcategory of $\dMan$. We say that a d-manifold $\bX$ {\it is a
manifold\/}\I{d-manifold!is a manifold} if it lies in $\hMan$. The
2-functor $F_\Man^\dSpa:\Man\ra\dSpa$ maps into $\dMan$, and we will
write $F_\Man^\dMan=F_\Man^\dSpa:\Man\ra\dMan$.
\label{ds4def1}
\end{dfn}

Here, as in \cite[\S 3.2]{Joyc6}, are alternative descriptions of
principal d-manifolds:

\begin{prop} The following are equivalent characterizations of when
a d-space $\bW$ is a principal d-manifold:
\begin{itemize}
\setlength{\itemsep}{0pt}
\setlength{\parsep}{0pt}
\item[{\bf(a)}] $\bW\simeq \bX\t_{\bs g,\bZ,\bs
h}\bY$ for $\bX,\bY,\bZ\in\hMan$.
\item[{\bf(b)}] $\bW\simeq \bX\t_{\bs i,\bZ,\bs
j}\bY,$ where $X,Y,Z$ are manifolds, $i:X\ra Z,$ $j:Y\ra Z$ are
embeddings, $\bX=F_\Man^\dSpa(X),$ and similarly for $Y,Z,i,j$.
That is, $\bW$ is an intersection of two submanifolds $X,Y$ in
$Z,$ in the sense of d-spaces.
\item[{\bf(c)}] $\bW\simeq \bV\t_{\bs s,\bs E,\bs 0}\bV,$ where $V$
is a manifold, $E\ra V$ is a vector bundle, $s:V\ra E$ is a
smooth section, $0:V\ra E$ is the zero section,
$\bV=F_\Man^\dSpa(V),$ and similarly for $E,s,0$. That is, $\bW$
is the zeroes $s^{-1}(0)$ of a smooth section $s$ of a vector
bundle $E,$ in the sense of d-spaces.\I{principal
d-manifold|)}\I{d-manifold!principal|)}
\end{itemize}
\label{ds4prop1}
\end{prop}

\begin{ex} Let $X\subseteq\R^n$ be any closed subset. By a lemma of
Whitney's, we can write $X$ as the zero set of a smooth function
$f:\R^n\ra\R$. Then $\bX=\bR^{\bs n}\t_{\bs f,\bR,\bs 0}\bs *$ is a
principal d-manifold, with topological space~$X$.
\label{ds4ex1}
\end{ex}

This example shows that the topological spaces $X$ underlying
d-manifolds $\bX$ can be fairly wild, for example, $X$ could be a
fractal\I{fractal} such as the Cantor set.\I{Cantor set}

\subsection{`Standard model' d-manifolds, 1- and 2-morphisms}
\label{ds42}
\I{d-manifold!standard model|(}

The next three examples, taken from \cite[\S 3.2 \& \S 3.4]{Joyc6},
give explicit models for principal d-manifolds in the form
$\bV\t_{\bs s,\bs E,\bs 0}\bV$ from Proposition \ref{ds4prop1}(c)
and their 1- and 2-morphisms, which we call {\it standard models}.

\begin{ex} Let $V$ be a manifold, $E\ra V$ a vector bundle (which we
sometimes call the {\it obstruction bundle\/}), and $s\in C^\iy(E)$.
We will write down an explicit principal d-manifold
$\bS=(\uS,\OSp,\ES,\im_S,\jm_S)$ which is equivalent to $\bV\t_{\bs
s,\bs E,\bs 0}\bV$ in Proposition \ref{ds4prop1}(c). We call $\bS$
the {\it standard model\/} of $(V,E,s)$, and also write it
$\bS_{V,E,s}$.\G[SVEsa]{$\bS_{V,E,s}$}{`standard model' d-manifold}
Proposition \ref{ds4prop1} shows that every principal d-manifold
$\bW$ is equivalent to $\bS_{V,E,s}$ for some~$V,E,s$.\I{principal
d-manifold}\I{d-manifold!principal}

Write $C^\iy(V)$ for the $C^\iy$-ring of smooth functions
$c:V\ra\R$, and $C^\iy(E),\ab C^\iy(E^*)$ for the vector spaces of
smooth sections of $E,E^*$ over $V$. Then $s$ lies in $C^\iy(E)$,
and $C^\iy(E),C^\iy(E^*)$ are modules over $C^\iy(V)$, and there is
a natural bilinear product $\cdot\,:C^\iy(E^*)\t C^\iy(E)\ra
C^\iy(V)$. Define $I_s\subseteq C^\iy(V)$ to be the ideal generated
by $s$. That is,
\e
I_s=\bigl\{\al\cdot s:\al\in C^\iy(E^*)\bigr\}\subseteq C^\iy(V).
\label{ds4eq1}
\e

Let $I_s^2=\an{fg:f,g\in I_s}_\R$ be the square of $I_s$. Then
$I_s^2$ is an ideal in $C^\iy(V)$, the ideal generated by $s\ot s\in
C^\iy(E\ot E)$. That is,
\begin{equation*}
I_s^2=\bigl\{\be\cdot (s\ot s):\be\in C^\iy(E^*\ot E^*)\bigr\}
\subseteq C^\iy(V).
\end{equation*}
Define $C^\iy$-rings $\fC=C^\iy(V)/I_s$, $\fC'=C^\iy(V)/I_s^2$, and
let $\pi:\fC'\ra\fC$ be the natural projection from the inclusion
$I_s^2\subseteq I_s$. Define a topological space $S=\{v\in
V:s(v)=0\}$, as a subspace of $V$. Now $s(v)=0$ if and only if
$(s\ot s)(v)=0$. Thus $S$ is the underlying topological space for
both $\Spec\fC$ and $\Spec\fC'$. So $\Spec\fC=\uS=(S,\OS)$,
$\Spec\fC'=\uS'=(S,\OSp)$, and $\Spec\pi=\uim_S=(\id_S,\im_S):
\uS'\ra\uS$, where $\uS,\uS'$ are fair affine $C^\iy$-schemes, and
$\OS,\OSp$ are sheaves of $C^\iy$-rings on $S$, and
$\im_S:\OSp\ra\OS$ is a morphism of sheaves of $C^\iy$-rings. Since
$\pi$ is surjective with kernel the square zero ideal\I{square zero
ideal} $I_s/I_s^2$, $\im_S$ is surjective, with kernel $\IS$ a sheaf
of square zero ideals in~$\OSp$.

From \eq{ds4eq1} we have a surjective $C^\iy(V)$-module morphism
$C^\iy(E^*)\ra I_s$ mapping $\al\mapsto \al\cdot s$. Applying
$\ot_{C^\iy(V)}\fC$ gives a surjective $\fC$-module morphism
\begin{equation*}
\si:C^\iy(E^*)/(I_s\cdot C^\iy(E^*))\longra I_s/I_s^2,\quad
\si: \al+(I_s\cdot C^\iy(E^*))\longmapsto \al\cdot s+I_s^2.
\end{equation*}
Define $\ES=\MSpec\bigl(C^\iy(E^*)/(I_s\cdot C^\iy(E^*))\bigr)$.
Also $\MSpec(I_s/I_s^2)=\IS$, so $\jm_S=\MSpec\si$ is a surjective
morphism $\jm_S:\ES\ra\IS$ in $\qcoh(\uS)$. Therefore
$\bS_{V,E,s}=\bS=(\uS,\OSp,\ES,\im_S,\jm_S)$ is a d-space.

In fact $\cE_S$ is a vector bundle on $\uS$ naturally isomorphic to
$\cE^*\vert_\uS$, where $\cE$ is the vector bundle on
$\uV=F_\Man^\CSch(V)$ corresponding to $E\ra V$. Also $\FS\cong
T^*\uV\vert_\uS$. The morphism $\phi_S:\ES\ra\FS$ can be interpreted
as follows: choose a connection $\nabla$ on $E\ra V$. Then $\nabla
s\in C^\iy(E\ot T^*V)$, so we can regard $\nabla s$ as a morphism of
vector bundles $E^*\ra T^*V$ on $V$. This lifts to a morphism of
vector bundles $\hat\nabla s:\cE^*\ra T^*\uV$ on the $C^\iy$-scheme
$\uV$, and $\phi_S$ is identified with $\hat\nabla
s\vert_\uS:\cE^*\vert_\uS\ra T^*\uV\vert_\uS$ under the isomorphisms
$\ES\cong\cE^*\vert_\uS$, $\FS\cong T^*\uV\vert_\uS$.
\label{ds4ex2}
\end{ex}

Proposition \ref{ds4prop1} implies that every principal d-manifold
$\bW$ is equivalent to $\bS_{V,E,s}$ for some $V,E,s$. The notation
$O(s)$\G[O(s)]{$O(s)$}{an error term in the ideal generated by a
section $s\in C^\iy(E)$} and $O(s^2)$\G[O(s)']{$O(s^2)$}{an error
term in the ideal generated by $s\ot s$ for $s\in C^\iy(E)$} used
below should be interpreted as follows. Let $V$ be a manifold, $E\ra
V$ a vector bundle, and $s\in C^\iy(E)$. If $F\ra V$ is another
vector bundle and $t\in C^\iy(F)$, then we write $t=O(s)$ if
$t=\al\cdot s$ for some $\al\in C^\iy(F\ot E^*)$, and $t=O(s^2)$ if
$t=\be\cdot (s\ot s)$ for some $\be\in C^\iy(F\ot E^*\ot E^*)$.
Similarly, if $W$ is a manifold and $f,g:V\ra W$ are smooth then we
write $f=g+O(s)$ if $c\ci f-c\ci g=O(s)$ for all smooth $c:W\ra\R$,
and $f=g+O(s^2)$ if $c\ci f-c\ci g=O(s^2)$ for all~$c$.

\begin{ex} Let $V,W$ be manifolds, $E\ra V$, $F\ra W$ be vector
bundles, and $s\in C^\iy(E)$, $t\in C^\iy(F)$. Write
$\bX=\bS_{V,E,s}$, $\bY=\bS_{W,F,t}$ for the `standard model'
principal d-manifolds from Example \ref{ds4ex2}. Suppose $f:V\ra W$
is a smooth map, and $\hat f:E\ra f^*(F)$ is a morphism of vector
bundles on $V$ satisfying
\e
\hat f\ci s= f^*(t)+O(s^2)\quad\text{in $C^\iy\bigl(f^*(F)\bigr)$.}
\label{ds4eq2}
\e

We will define a 1-morphism $\bs g=(\ug,g',g''):\bX\ra\bY$ in
$\dMan$ using $f,\hat f$. We will also write $\bs g:\bX\ra\bY$ as
$\bS_{\smash{f,\hat f}}:\bS_{V,E,s}\ra \bS_{W,F,t}$, and call it a
{\it standard model\/ $1$-morphism}.\I{d-manifold!standard
model!1-morphism} If $x\in X$ then $x\in V$ with $s(x)=0$, so
\eq{ds4eq2} implies that
\begin{equation*}
t\bigl(f(x)\bigr)=\bigl(f^*(t)\bigr)(x)=
\hat f\bigl(s(x)\bigr)+O\bigl(s(x)^2\bigr)=0,
\end{equation*}
so $f(x)\in Y\subseteq W$. Thus $g:=f\vert_X$ maps~$X\ra Y$.

Define morphisms of $C^\iy$-rings
\begin{align*}
\phi:C^\iy(W)/I_t&\longra C^\iy(V)/I_s,&
\phi':C^\iy(W)/I_t^2&\longra C^\iy(V)/I_s^2,\\
\text{by}\quad \phi:c+I_t&\longmapsto c\ci f+I_s,&
\phi':c+I_t^2&\longmapsto c\ci f+I_s^2.
\end{align*}
Here $\phi$ is well-defined since if $c\in I_t$ then $c=\ga\cdot t$
for some $\ga\in C^\iy(F^*)$, so
\begin{equation*}
c\ci f\!=\!(\ga\cdot t)\ci f\!=\!f^*(\ga)\cdot f^*(t)\!=\!f^*(\ga)
\cdot\bigl(\hat f\ci s+O(s^2)\bigr)\!=\!\bigl(\hat f\ci f^*(\ga)
\bigr)\cdot s+O(s^2)\in I_s.
\end{equation*}
Similarly if $c\in I_t^2$ then $c\ci f\in I_s^2$, so $\phi'$ is
well-defined. Thus we have $C^\iy$-scheme morphisms
$\ug=(g,g^\sh)=\Spec\phi:\uX\ra\uY$, and
$(g,g')=\Spec\phi':(X,\OXp)\ra(Y,\OYp)$, both with underlying map
$g$. Hence $g^\sh:g^{-1}(\OY)\ra\OX$ and $g':g^{-1}(\OYp)\ra\OXp$
are morphisms of sheaves of $C^\iy$-rings on~$X$.

Since $\ug^*(\EY)=\MSpec\bigl(C^\iy(f^*(F^*))/(I_s\cdot
C^\iy(f^*(F^*))\bigr)$, we may define $g'':\ug^*(\EY)\ra\EX$ by
$g''=\MSpec(G'')$, where
\begin{align*}
G'':C^\iy(f^*(F^*))/(I_s\cdot C^\iy(f^*(F^*))&\longra
C^\iy(E^*)/(I_s\cdot C^\iy(E^*))\\
\text{is defined by}\quad G'':\ga+I_s\cdot C^\iy(f^*(F^*))
&\longmapsto\ga\ci\hat f+I_s\cdot C^\iy(E^*).
\end{align*}
This defines $\bs g=(\ug,g',g'')$. One can show it is a 1-morphism
$\bs g:\bX\ra\bY$ in $\dMan$, which we also write
as~$\bS_{\smash{f,\hat f}}:\bS_{V,E,s}\ra
\bS_{W,F,t}$.\G[Sffa]{$\bS_{\smash{f,\hat f}}:\bS_{V,E,s}\ra
\bS_{W,F,t}$}{`standard model' 1-morphism in $\dMan$}

Suppose $\ti V\subseteq V$ is open, with inclusion $i_{\smash{\ti
V}}:\ti V\ra V$. Write $\ti E=E\vert_{\ti V}= i_{\ti V}^*(E)$ and
$\ti s=s\vert_{\ti V}$. Define $\bs i_{\smash{\ti
V,V}}=\bS_{i_{\smash{\ti V}},\id_{\ti E}}:\bS_{\smash{\ti V,\ti
E,\ti s}}\ra\bS_{V,E,s}$.\G[iVVa]{$\bs i_{\smash{\ti
V,V}}:\bS_{\smash{\ti V,\ti E,\ti s}}\ra\bS_{V,E,s}$}{inclusion of
open set in `standard model' d-manifold} If $s^{-1}(0)\subseteq\ti
V$ then $\bs i_{\smash{\ti V,V}}$ is a 1-isomorphism, with inverse
$\bs i_{\smash{\ti V,V}}^{-1}$. That is, making $V$ smaller without
making $s^{-1}(0)$ smaller does not really change $\bS_{V,E,s}$; the
d-manifold $\bS_{V,E,s}$ depends only on $E,s$ in an arbitrarily
small open neighbourhood of $s^{-1}(0)$ in~$V$.

\label{ds4ex3}
\end{ex}

\begin{ex} Let $V,W$ be manifolds, $E\ra V$, $F\ra W$ be vector
bundles, and $s\in C^\iy(E)$, $t\in C^\iy(F)$. Suppose $f,g:V\ra W$
are smooth and $\hat f:E\ra f^*(F),$ $\hat g:E\ra g^*(F)$ are vector
bundle morphisms with $\hat f\ci s= f^*(t)+O(s^2)$ and $\hat g\ci s=
g^*(t)+O(s^2),$ so we have 1-morphisms $\bS_{\smash{f,\hat
f}},\bS_{\smash{g,\hat g}}:\bS_{V,E,s}\ra \bS_{W,F,t}$. It is easy
to show that $\bS_{\smash{f,\hat f}}=\bS_{\smash{g,\hat g}}$ if and
only if $g=f+O(s^2)$ and~$\hat g=\hat f+O(s)$.

Now suppose $\La:E\ra f^*(TW)$ is a morphism of vector bundles on
$V$. Taking the dual of $\La$ and lifting to $\uV$ gives
$\La^*:\uf^*(T^*\uW)\ra\cE^*$. Restricting to the $C^\iy$-subscheme
$\uX=s^{-1}(0)$ in $\uV$ gives $\la=\La^*\vert_\uX:
\uf^*(\FY)\cong\uf^*(T^*\uW) \vert_\uX\ra\cE^*\vert_\uX=\EX$. One
can show that $\la$ is a 2-morphism $\bS_{\smash{f,\hat
f}}\Ra\bS_{\smash{g,\hat g}}$ if and only if
\begin{equation*}
g= f+\La\ci s+O(s^2)\quad\text{and}\quad \hat g=\hat f+f^*(\d
t)\ci\La+O(s).
\end{equation*}
Then we write $\la$ as $\bS_\La:\bS_{\smash{f,\hat f}}\Ra\bS_{g,\hat
g}$,\G[SLa]{$\bS_\La:\bS_{\smash{f,\hat f}}\Ra\bS_{g,\hat
g}$}{`standard model' 2-morphism in $\dMan$} and call it a {\it
standard model\/ $2$-morphism}.\I{d-manifold!standard
model!2-morphism} Every 2-morphism $\eta:\bS_{\smash{f,\hat
f}}\Ra\bS_{\smash{g,\hat g}}$ is $\bS_\La$ for some $\La$. Two
vector bundle morphisms $\La,\La':E\ra f^*(TW)$ have
$\bS_\La=\bS_{\smash{\La'}}$ if and only if~$\La=\La'+O(s)$.
\label{ds4ex4}
\end{ex}

If $\bX$ is a d-manifold and $x\in\bX$ then $x$ has an open
neighbourhood $\bU$ in $\bX$ equivalent in $\dSpa$ to $\bS_{V,E,s}$
for some manifold $V$, vector bundle $E\ra V$ and $s\in C^\iy(E)$.
In \cite[\S 3.3]{Joyc6} we investigate the extent to which $\bX$
determines $V,E,s$ near a point in $\bX$ and $V$, and prove:

\begin{thm} Let\/ $\bX$ be a d-manifold, and\/ $x\in\bX$. Then
there exists an open neighbourhood\/ $\bU$ of\/ $x$ in $\bX$ and an
equivalence $\bU\simeq\bS_{V,E,s}$ in $\dMan$ for some manifold\/
$V,$ vector bundle $E\ra V$ and\/ $s\in C^\iy(E)$ which identifies
$x\in\bU$ with a point\/ $v\in V$ such that\/ $s(v)=\d s(v)=0,$
where $\bS_{V,E,s}$ is as in Example\/ {\rm\ref{ds4ex2}}. These
$V,E,s$ are determined up to non-canonical isomorphism near $v$ by
$\bX$ near $x,$ and in fact they depend only on the underlying\/
$C^\iy$-scheme $\uX$ and the integer\/~$\vdim\bX$.
\label{ds4thm1}
\end{thm}

Thus, if we impose the extra condition $\d s(v)=0$, which is in fact
equivalent to choosing $V,E,s$ with $\dim V$ as small as possible,
then $V,E,s$ are determined uniquely near $v$ by $\bX$ near $x$
(that is, $V,E,s$ are determined locally up to isomorphism, but not
up to canonical isomorphism). If we drop the condition $\d s(v)=0$
then $V,E,s$ are determined uniquely near $v$ by $\bX$ near $x$
and~$\dim V$.

Theorem \ref{ds4thm1} shows that any d-manifold
$\bX=(\uX,\OXp,\EX,\im_X,\jm_X)$ is determined up to equivalence in
$\dSpa$ near any point $x\in\bX$ by the `classical' underlying
$C^\iy$-scheme $\uX$ and the integer $\vdim\bX$. So we can ask: what
extra information about $\bX$ is contained in the `derived' data
$\OXp,\EX,\im_X,\jm_X$? One can think of this extra information as
like a vector bundle $\cE$ over $\uX$. The only local information in
a vector bundle $\cE$ is $\rank\cE\in\Z$, but globally it also
contains nontrivial algebraic-topological information.

Suppose now that $\bs f:\bX\ra\bY$ is a 1-morphism in $\dMan$, and
$x\in\bX$ with $\bs f(x)=y\in\bY$. Then by Theorem \ref{ds4thm1} we
have $\bX\simeq\bS_{V,E,s}$ near $x$ and $\bY\simeq\bS_{W,F,t}$ near
$y$. So up to composition with equivalences, we can identify $\bs f$
near $x$ with a 1-morphism $\bs g:\bS_{V,E,s}\ra\bS_{W,F,t}$. Thus,
to understand arbitrary 1-morphisms $\bs f$ in $\dMan$ near a point,
it is enough to study 1-morphisms $\bs g:\bS_{V,E,s}\ra\bS_{W,F,t}$.
Our next theorem, proved in \cite[\S 3.4]{Joyc6}, shows that after
making $V$ smaller, every 1-morphism $\bs g:\bS_{V,E,s}
\ra\bS_{W,F,t}$ is of the form~$\bS_{\smash{f,\hat f}}$.

\begin{thm} Let\/ $V,W$ be manifolds, $E\ra V,$ $F\ra W$ be vector
bundles, and\/ $s\in C^\iy(E),$ $t\in C^\iy(F)$. Define principal
d-manifolds\/ $\bX=\bS_{V,E,s},$ $\bY=\bS_{W,F,t},$ with topological
spaces $X=\{v\in V:s(v)=0\}$ and\/ $Y=\{w\in W:t(w)=0\}$. Suppose
$\bs g:\bX\ra\bY$ is a $1$-morphism. Then there exist an open
neighbourhood\/ $\ti V$ of\/ $X$ in $V,$ a smooth map $f:\ti V\ra
W,$ and a morphism of vector bundles $\hat f:\ti E\ra f^*(F)$ with\/
$\hat f\ci\ti s= f^*(t),$ where $\ti E=E\vert_{\ti V},$ $\ti
s=s\vert_{\ti V},$ such that\/ $\bs g=\bS_{\smash{f,\hat f}}\ci\bs
i_{\ti V,V}^{-1},$ where $\bs i_{\ti V,V}=\bS_{\smash{\id_{\ti
V},\id_{\ti E}}}:\bS_{\smash{\ti V,\ti E,\ti s}}\ra \bS_{V,E,s}$ is
a $1$-isomorphism, and\/ $\bS_{\smash{f,\hat f}}:\bS_{\smash{\ti
V,\ti E,\ti s}}\ra \bS_{W,F,t}$ is as in
Example\/~{\rm\ref{ds4ex3}}.
\label{ds4thm2}
\end{thm}

These results give a good differential-geometric picture of
d-manifolds and their 1- and 2-morphisms near a point. The $O(s)$
and $O(s^2)$ notation helps keep track of what information from
$V,E,s$ and $f,\hat f$ and $\La$ is remembered and what forgotten by
the d-manifolds $\bS_{V,E,s}$, 1-morphisms $\bS_{\smash{f,\hat f}}$
and 2-morphisms~$\bS_\La$.\I{d-manifold!standard model|)}

\subsection{The 2-category of virtual vector bundles}
\label{ds43}
\I{virtual vector bundle|(}\I{virtual quasicoherent sheaf|(}

In our theory of derived differential geometry, it is a general
principle that categories in classical differential geometry should
often be replaced by 2-categories, and classical concepts be
replaced by 2-categorical analogues.

In classical differential geometry, if $X$ is a manifold, the vector
bundles $E\ra X$ and their morphisms form a category $\vect(X)$. The
cotangent bundle $T^*X$ is an important example of a vector bundle.
If $f:X\ra Y$ is smooth then pullback $f^*:\vect(Y)\ra\vect(X)$ is a
functor. There is a natural morphism $\d f^*:f^*(T^*Y)\ra T^*X$. We
now explain 2-categorical analogues of all this for d-manifolds,
following~\cite[\S 3.1--\S 3.2]{Joyc6}.

\begin{dfn} Let $\uX$ be a $C^\iy$-scheme, which will
usually be the $C^\iy$-scheme underlying a d-manifold $\bX$. We will
define a 2-category\I{2-category}
$\vqcoh(\uX)$\G[vqcoh(X)a]{$\vqcoh(\uX)$}{2-category of virtual
quasicoherent sheaves on a $C^\iy$-scheme $\uX$} of {\it virtual
quasicoherent sheaves\/} on $\uX$. {\it Objects\/} of $\vqcoh(\uX)$
are morphisms $\phi:\cE^1\ra\cE^2$ in $\qcoh(\uX)$, which we also
may write as $(\cE^1,\cE^2,\phi)$ or
$(\cE^\bu,\phi)$.\G[Ephi]{$(\cE^\bu,\phi)$}{virtual quasicoherent
sheaf, or virtual vector bundle} Given objects $\phi:\cE^1\ra\cE^2$
and $\psi:\cF^1\ra\cF^2$, a 1-{\it morphism\/}
$(f^1,f^2):(\cE^\bu,\phi)\ra(\cF^\bu,\psi)$ is a pair of morphisms
$f^1:\cE^1\ra\cF^1$, $f^2:\cE^2\ra\cF^2$ in $\qcoh(\uX)$ such that
$\psi\ci f^1=f^2\ci\phi$. We write $f^\bu$
for~$(f^1,f^2)$.\I{virtual vector bundle!on a
$C^\iy$-scheme}\I{virtual quasicoherent sheaf!on $C^\iy$-scheme}

The {\it identity\/ $1$-morphism\/} of $(\cE^\bu,\phi)$ is
$(\id_{\cE^1},\id_{\cE^2})$. The {\it composition\/} of 1-morphisms
$f^\bu:(\cE^\bu,\phi)\ra(\cF^\bu,\psi)$ and
$g^\bu:(\cF^\bu,\psi)\ra(\cG^\bu,\xi)$ is $g^\bu\ci f^\bu=(g^1\ci
f^1,g^2\ci f^2):(\cE^\bu,\phi)\ra(\cG^\bu,\xi)$.

Given $f^\bu,g^\bu:(\cE^\bu,\phi)\ra(\cF^\bu,\psi)$, a 2-{\it
morphism\/} $\eta:f^\bu\Ra g^\bu$ is a morphism $\eta:\cE^2\ra\cF^1$
in $\qcoh(\uX)$ such that $g^1=f^1+\eta\ci\phi$ and
$g^2=f^2+\psi\ci\eta$. The {\it identity\/ $2$-morphism\/} for
$f^\bu$ is $\id_{f^\bu}=0$. If $f^\bu,g^\bu,h^\bu:(\cE^\bu,\phi)\ra
(\cF^\bu,\psi)$ are 1-morphisms and $\eta:f^\bu\Ra g^\bu$,
$\ze:g^\bu\Ra h^\bu$ are 2-morphisms, the {\it vertical composition
of\/ $2$-morphisms\/} $\ze\od\eta:f^\bu\Ra h^\bu$ is
$\ze\od\eta=\ze+\eta$. If $f^\bu,\ti
f{}^\bu:(\cE^\bu,\phi)\ra(\cF^\bu, \psi)$ and $g^\bu,\ti
g{}^\bu:(\cF^\bu,\psi)\ra(\cG^\bu,\xi)$ are 1-morphisms and
$\eta:f^\bu\Ra\ti f{}^\bu$, $\ze:g^\bu\Ra\ti g{}^\bu$ are
2-morphisms, the {\it horizontal composition of\/ $2$-morphisms\/}
$\ze*\eta:g^\bu\ci f^\bu\Ra\ti g{}^\bu\ci\ti f{}^\bu$ is
$\ze*\eta=g^1\ci\eta+\ze\ci f^2+\ze\ci\psi\ci\eta$. This defines a
strict 2-category $\vqcoh(\uX)$, the obvious 2-category of 2-term
complexes in~$\qcoh(\uX)$.

If $\uU\subseteq\uX$ is an open $C^\iy$-subscheme then restriction
from $\uX$ to $\uU$ defines a strict 2-functor
$\vert_\uU:\vqcoh(\uX)\ra\vqcoh(\uU)$. An object $(\cE^\bu,\phi)$ in
$\vqcoh(\uX)$ is called a {\it virtual vector bundle of rank\/}
$d\in\Z$ if $\uX$ may be covered by open $\uU\subseteq\uX$ such that
$(\cE^\bu,\phi)\vert_\uU$ is equivalent in $\vqcoh(\uU)$ to some
$(\cF^\bu,\psi)$ for $\cF^1,\cF^2$ vector bundles on $\uU$ with
$\rank\cF^2-\rank\cF^1=d$. We write $\rank(\cE^\bu,\phi)=d$. If
$\uX\ne\es$ then $\rank(\cE^\bu,\phi)$ depends only on
$\cE^1,\cE^2,\phi$, so it is well-defined. Write
$\vvect(\uX)$\G[vvect(X)a]{$\vvect(\uX)$}{2-category of virtual
vector bundles on a $C^\iy$-scheme $\uX$} for the full 2-subcategory
of virtual vector bundles in~$\vqcoh(\uX)$.

If $\uf:\uX\ra\uY$ is a $C^\iy$-scheme morphism then pullback gives
a strict 2-functor\I{2-category!strict 2-functor}
$\uf^*:\vqcoh(\uY)\ra\vqcoh(\uX)$, which maps
$\vvect(\uY)\ra\vvect(\uX)$.
\label{ds4def2}
\end{dfn}

We apply these ideas to d-spaces.

\begin{dfn} Let $\bX=(\uX,\OXp,\EX,\im_X,\jm_X)$ be a d-space. Define
the {\it virtual cotangent sheaf\/} $T^*\bX$ of $\bX$ to be the
morphism $\phi_X:\EX\ra\FX$ in $\qcoh(\uX)$ from Definition
\ref{ds3def1}, regarded as a virtual quasicoherent sheaf
on~$\uX$.\G[T*Xa]{$T^*\bX$}{virtual cotangent sheaf of a d-space
$\bX$}

Let $\bs f\!=\!(\uf,f',f''):\bX\ra\bY$ be a 1-morphism in $\dSpa$.
Then $T^*\bX\!=\!(\EX,\FX,\phi_X)$ and
$\uf^*(T^*\bY)\!=\!\bigl(\uf^*(\EY), \uf^*(\FY),\uf^*(\phi_Y)\bigr)$
are virtual quasicoherent sheaves on $\uX$, and $\Om_{\bs
f}:=\!(f'',f^2)$ is a 1-morphism $\uf^*(T^*\bY)\!\ra\! T^*\bX$ in
$\vqcoh(\uX)$, as \eq{ds3eq1} commutes.

Let $\bs f,\bs g:\bX\ra\bY$ be 1-morphisms in $\dSpa$, and $\eta:\bs
f\Ra\bs g$ a 2-morphism. Then $\eta:\uf^*(\FY)\ra\EX$ with
$g''=f''+\eta\ci\uf^*(\phi_Y)$ and $g^2=f^2+\phi_X\ci\eta$, as in
\eq{ds3eq3}. It follows that $\eta$ is a 2-morphism $\Om_{\bs
f}\Ra\Om_{\bs g}$ in $\vqcoh(\uX)$. Thus, objects, 1-morphisms and
2-morphisms in $\dSpa$ lift to objects, 1-morphisms and 2-morphisms
in $\vqcoh(\uX)$.
\label{ds4def3}
\end{dfn}

The next proposition justifies the definition of virtual vector
bundle. Because of part (b), if $\bX$ is a d-manifold we call
$T^*\bX$ the {\it virtual cotangent bundle\/}\I{virtual cotangent
bundle} of $\bX$, rather than the virtual cotangent
sheaf.\I{d-manifold!virtual cotangent bundle}

\begin{prop}{\bf(a)} Let\/ $V$ be a manifold, $E\ra V$ a vector
bundle, and\/ $s\in C^\iy(E)$. Then Example\/ {\rm\ref{ds4ex2}}
defines a d-manifold\/ $\bS_{V,E,s}$. Its cotangent bundle
$T^*\bS_{V,E,s}$ is a virtual vector bundle on $\uS_{V,E,s}$ of
rank\/~$\dim V-\rank E$.
\smallskip

\noindent{\bf(b)} Let\/ $\bX$ be a d-manifold. Then $T^*\bX$ is a
virtual vector bundle on $\uX$ of rank\/ $\vdim\bX$. Hence if\/
$\bX\ne\bs\es$ then $\vdim\bX$ is well-defined.
\label{ds4prop2}
\end{prop}

The virtual cotangent bundle\I{virtual cotangent bundle} $T^*\bX$ of
a d-manifold $\bX$ is a d-space analogue of the {\it cotangent
complex\/}\I{cotangent complex} in algebraic geometry, as in Illusie
\cite{Illu}. It contains only a fraction of the information in
$\bX=(\uX,\OXp,\EX,\im_X, \jm_X)$, but many interesting properties
of d-manifolds $\bX$ and 1-morphisms $\bs f:\bX\ra\bY$ can be
expressed solely in terms of virtual cotangent bundles
$T^*\bX,T^*\bY$ and 1-morphisms $\Om_{\bs f}:\uf^*(T^*\bY)\ra
T^*\bX$. Here is an example of this.

\begin{dfn} Let $\uX$ be a $C^\iy$-scheme. We say that a
virtual vector bundle $(\cE^1,\cE^2,\phi)$ on $\uX$ {\it is a vector
bundle\/} if it is equivalent in $\vvect(\uX)$ to $(0,\cE,0)$ for
some vector bundle $\cE$ on $\uX$. One can show $(\cE^1,\cE^2,\phi)$
is a vector bundle if and only if $\phi$ has a left inverse
in~$\qcoh(\uX)$.\I{virtual vector bundle!is a vector bundle}
\label{ds4def4}
\end{dfn}

\begin{prop} Let\/ $\bX$ be a d-manifold. Then $\bX$ is a manifold
(that is, $\bX\in\hMan$) if and only if\/ $T^*\bX$ is a vector
bundle, or equivalently, if\/ $\phi_X:\EX\ra\FX$ has a left inverse
in\/~$\qcoh(\uX)$.\I{d-manifold!is a manifold}\I{virtual vector
bundle|)}\I{virtual quasicoherent sheaf|)}
\label{ds4prop3}
\end{prop}

\subsection{Equivalences in $\dMan$, and gluing by equivalences}
\label{ds44}
\I{d-manifold!equivalence|(}

Equivalences in a 2-category are defined in \S\ref{dsA3}.
Equivalences in $\dMan$ are the best derived analogue of
isomorphisms in $\Man$, that is, of diffeomorphisms of manifolds. A
smooth map of manifolds $f:X\ra Y$ is called {\it \'etale\/} if it
is a local diffeomorphism. Here is the derived analogue.

\begin{dfn} Let $\bs f:\bX\ra\bY$ be a 1-morphism in $\dMan$. We
call $\bs f$ {\it \'etale\/} if it is a {\it local equivalence},
that is, if for each $x\in\bX$ there exist open
$x\in\bU\subseteq\bX$ and $\bs f(x)\in\bV\subseteq\bY$ such that
$\bs f(\bU)=\bV$ and $\bs f\vert_\bU:\bU\ra\bV$ is an equivalence.
\label{ds4def5}
\end{dfn}

If $f:X\ra Y$\I{d-manifold!etale 1-morphism@\'etale 1-morphism|(} is
a smooth map of manifolds, then $f$ is \'etale if and only if $\d
f^*:f^*(T^*Y)\ra T^*X$ is an isomorphism of vector bundles. (The
analogue is false for schemes.) In \cite[\S 3.5]{Joyc6} we prove a
version of this for d-manifolds:

\begin{thm} Suppose $\bs f:\bX\ra\bY$ is a $1$-morphism of
d-manifolds. Then the following are equivalent:
\begin{itemize}
\setlength{\itemsep}{0pt}
\setlength{\parsep}{0pt}
\item[{\rm(i)}] $\bs f$ is \'etale;
\item[{\rm(ii)}] $\Om_{\bs f}:\uf^*(T^*\bY)\ra T^*\bX$ is an
equivalence in $\vqcoh(\uX);$ and
\item[{\rm(iii)}] the following is a split short
exact sequence in\/~$\qcoh(\uX)\!:$
\end{itemize}
\e
\smash{\xymatrix@C=25pt{ 0 \ar[r] & \uf^*(\EY) \ar[rr]^(0.45){f''
\op -\uf^*(\phi_Y)} && \EX\op \uf^*(\FY) \ar[rr]^(0.6){\phi_X\op
f^2} && \FX \ar[r] & 0.}}
\label{ds4eq3}
\e
If in addition $f:X\ra Y$ is a bijection, then $\bs f$ is an
equivalence in\/~$\dMan$.
\label{ds4thm3}
\end{thm}

Here a complex $0\ra E\ra F\ra G\ra 0$ in an abelian category $\cA$
is called a {\it split short exact sequence\/} if there exists an
isomorphism $F\cong E\op G$ in $\cA$ identifying the complex
with~$0\ra E \,{\buildrel \id\op 0\over\longra}\, E\op G\,{\buildrel
0\op \id\over\longra}\,G\ra 0$.\I{split short exact
sequence}\I{abelian category!split short exact sequence}

The analogue of Theorem \ref{ds4thm3} for d-spaces is false. When
$\bs f:\bX\ra\bY$ is a `standard model' 1-morphism
$\bS_{\smash{f,\hat f}}:\bS_{V,E,s}\ra\bS_{W,F,t}$, as in
\S\ref{ds42}, we can express the conditions for $\bS_{\smash{f,\hat
f}}$ to be \'etale or an equivalence in terms of~$f,\hat
f$.\I{d-manifold!standard model|(}

\begin{thm} Let\/ $V,W$ be manifolds, $E\ra V,$ $F\ra W$ be vector
bundles, $s\in C^\iy(E),$ $t\in C^\iy(F),$ $f:V\ra W$ be smooth,
and\/ $\hat f:E\ra f^*(F)$ be a morphism of vector bundles on $V$
with\/ $\hat f\ci s= f^*(t)+O(s^2)$. Then Example\/
{\rm\ref{ds4ex3}} defines a $1$-morphism\/ $\bS_{\smash{f,\hat
f}}:\bS_{V,E,s}\ra\bS_{W,F,t}$ in $\dMan$. This $\bS_{\smash{f,\hat
f}}$ is \'etale if and only if for each\/ $v\in V$ with\/ $s(v)=0$
and\/ $w=f(v)\in W,$ the following sequence of vector spaces is
exact:
\e
\smash{\xymatrix@C=18pt{ 0 \ar[r] & T_vV \ar[rrr]^(0.42){\d s(v)\op
\,\d f(v)} &&& E_v\op T_wW \ar[rrr]^(0.57){\hat f(v)\op\, -\d t(w)}
&&& F_w \ar[r] & 0.}}
\label{ds4eq4}
\e
Also $\bS_{\smash{f,\hat f}}$ is an equivalence if and only if in
addition\/ $f\vert_{s^{-1}(0)}:s^{-1}(0)\!\ra\! t^{-1}(0)$ is a
bijection, where $s^{-1}(0)\!=\!\{v\in V:s(v)\!=\!0\},$
$t^{-1}(0)\!=\!\{w\in W:t(w)\!=\!0\}$.
\label{ds4thm4}
\end{thm}
\I{d-manifold!etale 1-morphism@\'etale 1-morphism|)}

Section \ref{ds32} discussed gluing d-spaces by equivalences on open
d-subspaces. It generalizes immediately to d-manifolds: if in
Theorem \ref{ds3thm3} we fix $n\in\Z$ and take the initial d-spaces
$\bX_i$ to be d-manifolds with $\vdim\bX_i=n$, then the glued
d-space $\bY$ is also a d-manifold with $\vdim\bY=n$.

Here is an analogue of Theorem \ref{ds3thm3}, taken from \cite[\S
3.6]{Joyc6}, in which we take the d-spaces $\bX_i$ to be `standard
model' d-manifolds $\bS_{V_i,E_i,s_i}$, and the 1-morphisms $\bs
e_{ij}$ to be `standard model' 1-morphisms $\bS_{e_{ij},\hat
e_{ij}}$. We also use Theorem \ref{ds4thm4} in (ii) to characterize
when $\bs e_{ij}$ is an equivalence.\I{d-manifold!gluing by
equivalences|(}

\begin{thm} Suppose we are given the following data:
\begin{itemize}
\setlength{\itemsep}{0pt}
\setlength{\parsep}{0pt}
\item[{\rm(a)}] an integer $n;$
\item[{\rm(b)}] a Hausdorff, second countable topological space $X;$
\item[{\rm(c)}] an indexing set\/ $I,$ and a total order $<$ on $I;$
\item[{\rm(d)}] for each\/ $i$ in $I,$ a manifold\/ $V_i,$ a vector
bundle $E_i\ra V_i$ with\/ $\dim V_i-\rank E_i=n,$ a smooth
section $s_i:V_i\ra E_i,$ and a homeomorphism $\psi_i:X_i\ra\hat
X_i,$ where $X_i=\{v_i\in V_i:s_i(v_i)=0\}$ and\/ $\hat
X_i\subseteq X$ is open; and
\item[{\rm(e)}] for all\/ $i<j$ in $I,$ an open submanifold
$V_{ij}\subseteq V_i,$ a smooth map $e_{ij}:V_{ij}\ra V_j,$ and
a morphism of vector bundles $\hat e_{ij}:E_i\vert_{V_{ij}}\ra
e_{ij}^*(E_j)$.
\end{itemize}
Using notation $O(s_i),O(s_i^2)$ as in\/ {\rm\S\ref{ds42},} let this
data satisfy the conditions:
\begin{itemize}
\setlength{\itemsep}{0pt}
\setlength{\parsep}{0pt}
\item[{\rm(i)}] $X=\bigcup_{i\in I}\hat X_i;$
\item[{\rm(ii)}] if\/ $i<j$ in $I$ then $\hat e_{ij}\ci
s_i\vert_{V_{ij}}= e_{ij}^*(s_j)+O(s_i^2),$ $\psi_i(X_i\cap
V_{ij})=\hat X_i\cap\hat X_j,$ and\/ $\psi_i\vert_{X_i\cap
V_{ij}}=\psi_j\ci e_{ij}\vert_{X_i\cap V_{ij}},$ and if\/
$v_i\in V_{ij}$ with\/ $s_i(v_i)=0$ and\/ $v_j=e_{ij}(v_i)$ then
the following is exact:
\begin{equation*}
\smash{\xymatrix@C=18pt{ 0 \ar[r] & T_{v_i}V_i \ar[rrr]^(0.42){\d
s_i(v_i)\op \,\d e_{ij}(v_i)} &&& E_i\vert_{v_i}\!\op\! T_{v_j}V_j
\ar[rrr]^(0.57){\hat e_{ij}(v_i)\op\, -\d s_j(v_j)} &&&
E_j\vert_{v_j} \ar[r] & 0;}}
\end{equation*}
\item[{\rm(iii)}] if\/ $i<j<k$ in $I$ then
\begin{align*}
e_{ik}\vert_{V_{ij}\cap V_{ik}}&= e_{jk}\ci
e_{ij}\vert_{V_{ij}\cap V_{ik}}+O(s_i^2)\qquad\text{and}\\
\hat e_{ik}\vert_{V_{ij}\cap V_{ik}}&= e_{ij}\vert_{V_{ij}\cap
V_{ik}}^*(\hat e_{jk})\ci \hat e_{ij}\vert_{V_{ij}\cap
V_{ik}}+O(s_i).
\end{align*}
\end{itemize}

Then there exist a d-manifold\/ $\bX$ with\/ $\vdim\bX=n$ and
underlying topological space $X,$ and a $1$-morphism
$\bs\psi_i:\bS_{V_i,E_i,s_i}\ra\bX$ with underlying continuous map
$\psi_i,$ which is an equivalence with the open d-submanifold\/
$\bs{\hat X}_i\subseteq\bX$ corresponding to $\hat X_i\subseteq X$
for all\/ $i\in I,$ such that for all\/ $i<j$ in $I$ there exists a
$2$-morphism\/ $\eta_{ij}:\bs\psi_j\ci\bS_{e_{ij},\hat
e_{ij}}\Ra\bs\psi_i\ci\bs i_{V_{ij},V_i},$ where $\bS_{e_{ij},\hat
e_{ij}}:\bS_{V_{ij},E_i\vert_{V_{ij}},s_i\vert_{V_{ij}}}\ra
\bS_{V_j,E_j,s_j}$ and\/ $\bs i_{V_{ij},V_i}:\bS_{V_{ij},
E_i\vert_{V_{ij}},s_i\vert_{V_{ij}}}\ra \bS_{V_i,E_i,s_i}$ are as in
Example\/ {\rm\ref{ds4ex2}}. This d-manifold\/ $\bX$ is unique up to
equivalence in~$\dMan$.

Suppose also that\/ $Y$ is a manifold, and\/ $g_i:V_i\ra Y$ are
smooth maps for all\/ $i\in I,$ and\/ $g_j\ci
e_{ij}=g_i\vert_{V_{ij}}+O(s_i)$ for all\/ $i<j$ in $I$. Then there
exist a $1$-morphism $\bs h:\bX\ra\bY$ unique up to $2$-isomorphism,
where $\bY=F_\Man^\dMan(Y)=\bS_{Y,0,0},$ and\/ $2$-morphisms
$\ze_i:\bs h\ci\bs\psi_i\Ra\bS_{g_i,0}$ for all\/ $i\in I$. Here
$\bS_{Y,0,0}$ is from Example\/ {\rm\ref{ds4ex2}} with vector bundle
$E$ and section $s$ both zero, and\/
$\bS_{g_i,0}:\bS_{V_i,E_i,s_i}\ra \bS_{Y,0,0}=\bY$ is from Example
{\rm\ref{ds4ex3}} with\/~$\hat g_i=0$.
\label{ds4thm5}
\end{thm}
\I{d-manifold!standard model|)}

The hypotheses of Theorem \ref{ds4thm5} are similar to the notion of
{\it good coordinate system\/}\I{good coordinate
system}\I{d-orbifold!good coordinate system} in the theory of
Kuranishi spaces\I{Kuranishi space} of Fukaya and Ono
\cite[Def.~6.1]{FuOn}, as discussed in \S\ref{ds119}. The importance
of Theorem \ref{ds4thm5} is that all the ingredients are described
wholly in differential-geometric or topological terms. So we can use
the theorem as a tool to prove the existence of d-manifold
structures on spaces coming from other areas of geometry, for
instance, on moduli spaces.\I{moduli space}\I{d-manifold!gluing by
equivalences|)}\I{d-manifold!equivalence|)}

\subsection{Submersions, immersions and embeddings}
\label{ds45}

Let $f:X\ra Y$ be a smooth map of manifolds. Then $\d
f^*:f^*(T^*Y)\ra T^*X$ is a morphism of vector bundles on $X$, and
$f$ is a {\it submersion\/} if $\d f^*$ is injective, and $f$ is an
{\it immersion\/} if $\d f^*$ is surjective. Here the appropriate
notions of injective and surjective for morphisms of vector bundles
are stronger than the corresponding notions for sheaves: $\d f^*$ is
{\it injective\/} if it has a left inverse, and {\it surjective\/}
if it has a right inverse.

In a similar way, if $\bs f:\bX\ra\bY$ is a 1-morphism of
d-manifolds, we would like to define $\bs f$ to be a submersion or
immersion if the 1-morphism $\Om_{\bs f}:\uf^*(T^*\bY)\ra T^*\bX$ in
$\vvect(\uX)$ is injective or surjective in some suitable sense. It
turns out that there are two different notions of injective and
surjective 1-morphisms in the 2-category $\vvect(\uX)$, a weak and a
strong:

\begin{dfn} Let $\uX$ be a $C^\iy$-scheme, $(\cE^1,\cE^2,\phi)$ and
$(\cF^1,\cF^2,\psi)$ be virtual vector bundles on $\uX$, and
$(f^1,f^2):(\cE^\bu,\phi)\ra(\cF^\bu,\psi)$ be a 1-morphism in
$\vvect(\uX)$. Then we have a complex in~$\qcoh(\uX)$:
\e
\xymatrix@C=25pt{ 0 \ar[r] & \cE^1 \ar@<.5ex>[rr]^(0.45){f^1 \op
-\phi} && \cF^1\op \cE^2 \ar@<.5ex>[rr]^(0.6){\psi\op f^2}
\ar@<.5ex>@{.>}[ll]^(0.55)\ga && \cF^2 \ar@<.5ex>@{.>}[ll]^(0.4)\de
\ar[r] & 0.}
\label{ds4eq5}
\e
One can show that $f^\bu$ is an equivalence in $\vvect(\uX)$ if and
only if \eq{ds4eq5} is a {\it split short exact sequence\/}\I{split
short exact sequence}\I{abelian category!split short exact sequence}
in $\qcoh(\uX)$. That is, $f^\bu$ is an equivalence if and only if
there exist morphisms $\ga,\de$ as shown in \eq{ds4eq5} satisfying
the conditions:
\e
\begin{aligned}
\ga\ci\de&=0,& \ga\ci(f^1 \op -\phi)&=\id_{\cE^1},\\
(f^1 \op -\phi)\ci\ga+\de\ci(\psi\op f^2)&=\id_{\cF^1\op \cE^2}, &
(\psi\op f^2)\ci\de&=\id_{\cF^2}.
\end{aligned}
\label{ds4eq6}
\e

Our notions of $f^\bu$ injective or surjective impose some but not
all of \eq{ds4eq6}:
\begin{itemize}
\setlength{\itemsep}{0pt}
\setlength{\parsep}{0pt}
\item[(a)] We call $f^\bu$ {\it weakly injective\/} if there
exists $\ga:\cF^1\op\cE^2\ra\cE^1$ in $\qcoh(\uX)$
with~$\ga\ci(f^1 \op -\phi)=\id_{\cE^1}$.
\item[(b)] We call $f^\bu$ {\it injective\/} if there exist
$\ga:\cF^1\op\cE^2\ra\cE^1$ and $\de:\cF^2\ra\cF^1\op\cE^2$ with
$\ga\ci\de=0$, $\ga\ci(f^1 \op -\phi)= \id_{\cE^1}$ and~$(f^1
\op -\phi)\ci\ga+\de\ci(\psi\op f^2)=\id_{\cF^1\op\cE^2}$.
\item[(c)] We call $f^\bu$ {\it weakly surjective\/} if there exists
$\de:\cF^2\ra\cF^1\op\cE^2$ in $\qcoh(\uX)$ with~$(\psi\op
f^2)\ci\de=\id_{\cF^2}$.
\item[(d)] We call $f^\bu$ {\it surjective\/} if there exist
$\ga:\cF^1\op\cE^2\ra\cE^1$ and $\de:\cF^2\ra\cF^1\op\cE^2$ with
$\ga\ci\de=0$, $\ga\ci(f^1 \op -\phi)= \id_{\cE^1}$
and~$(\psi\op f^2)\ci\de=\id_{\cF^2}$.
\end{itemize}
If $\uX$ is separated, paracompact, and locally fair, these are
local conditions\I{virtual vector bundle!weakly injective
1-morphism}\I{virtual vector bundle!injective 1-morphism}\I{virtual
vector bundle!weakly surjective 1-morphism}\I{virtual vector
bundle!surjective 1-morphism} on~$\uX$.

\label{ds4def6}
\end{dfn}

Using these we define weak and strong forms of submersions,
immersions, and embeddings for
d-manifolds.\I{d-manifold!w-submersion|(}\I{d-manifold!submersion|(}
\I{d-manifold!w-immersion|(}\I{d-manifold!immersion|(}
\I{d-manifold!w-embedding|(}\I{d-manifold!embedding|(}

\begin{dfn} Let $\bs f:\bX\ra\bY$ be a 1-morphism of d-manifolds.
Definition \ref{ds4def3} defines a 1-morphism $\Om_{\bs
f}:\uf^*(T^*\bY)\ra T^*\bX$ in $\vvect(\uX)$. Then:
\begin{itemize}
\setlength{\itemsep}{0pt}
\setlength{\parsep}{0pt}
\item[(a)] We call $\bs f$ a {\it w-submersion\/} if $\Om_{\bs
f}$ is weakly injective.
\item[(b)] We call $\bs f$ a {\it submersion\/} if $\Om_{\bs
f}$ is injective.
\item[(c)] We call $\bs f$ a {\it w-immersion\/} if $\Om_{\bs
f}$ is weakly surjective.
\item[(d)] We call $\bs f$ an {\it immersion\/} if $\Om_{\bs
f}$ is surjective.
\item[(e)] We call $\bs f$ a {\it w-embedding\/} if it is a
w-immersion and $f:X\ra f(X)$ is a homeomorphism, so in
particular $f$ is injective.
\item[(f)] We call $\bs f$ an {\it embedding\/} if it is an
immersion and $f$ is a homeomorphism with its image.
\end{itemize}
Here w-submersion is short for {\it weak submersion}, etc.
Conditions (a)--(d) all concern the existence of morphisms $\ga,\de$
in the next equation satisfying identities:
\begin{equation*}
\xymatrix@C=25pt{ 0 \ar[r] & \uf^*(\EY) \ar@<.5ex>[rr]^(0.45){f''
\op -\uf^*(\phi_Y)} && \EX\op \uf^*(\FY)
\ar@<.5ex>@{.>}[ll]^(0.55)\ga \ar@<.5ex>[rr]^(0.6){\phi_X\op f^2} &&
\FX \ar@<.5ex>@{.>}[ll]^(0.4)\de \ar[r] & 0.}
\end{equation*}

Parts (c)--(f) enable us to define {\it d-submanifolds\/} of
d-manifolds. {\it Open d-submanifolds\/} are open d-subspaces of a
d-manifold. More generally, we call $\bs i:\bX\ra\bY$ a {\it
w-immersed}, or {\it immersed}, or {\it w-embedded}, or {\it
embedded d-submanifold\/} of $\bY$, if $\bX,\bY$ are d-manifolds and
$\bs i$ is a w-immersion, immersion, w-embedding, or embedding,
respectively.\I{d-manifold!d-submanifold}
\label{ds4def7}
\end{dfn}

Here are some properties of these, taken from~\cite[\S 4.1--\S
4.2]{Joyc6}:

\begin{thm}{\bf(i)} Any equivalence of d-manifolds is a
w-submersion, submersion, w-immersion, immersion, w-embedding and
embedding.\I{d-manifold!equivalence}
\smallskip

\noindent{\bf(ii)} If\/ $\bs f,\bs g:\bX\ra\bY$ are $2$-isomorphic
$1$-morphisms of d-manifolds then $\bs f$ is a w-submersion,
submersion, \ldots, embedding, if and only if\/ $\bs g$ is.
\smallskip

\noindent{\bf(iii)} Compositions of w-submersions, submersions,
w-immersions, immersions, w-embeddings, and embeddings are
$1$-morphisms of the same kind.
\smallskip

\noindent{\bf(iv)} The conditions that a $1$-morphism of
d-manifolds\/ $\bs f:\bX\ra\bY$ is a w-submersion, submersion,
w-immersion or immersion are local in $\bX$ and\/ $\bY$. That is,
for each\/ $x\in\bX$ with\/ $\bs f(x)=y\in\bY,$ it suffices to check
the conditions for $\bs f\vert_\bU:\bU\ra\bV$ with $\bV$ an open
neighbourhood of\/ $y$ in $\bY,$ and\/ $\bU$ an open neighbourhood
of\/ $x$ in\/ ${\bs f}^{-1}(\bV)\subseteq\bX$. The conditions that\/
$\bs f:\bX\ra\bY$ is a w-embedding or embedding are local in $\bY,$
but not in\/~$\bX$.
\smallskip

\noindent{\bf(v)} Let\/ $\bs f:\bX\ra\bY$ be a submersion of
d-manifolds. Then $\vdim\bX\ge\vdim\bY,$ and if\/
$\vdim\bX=\vdim\bY$ then\/ $\bs f$ is \'etale.
\smallskip

\noindent{\bf(vi)} Let\/ $\bs f:\bX\ra\bY$ be an immersion of
d-manifolds. Then $\vdim\bX\le\vdim\bY,$ and if\/
$\vdim\bX=\vdim\bY$ then\/ $\bs f$ is \'etale.
\smallskip

\noindent{\bf(vii)} Let\/ $f:X\ra Y$ be a smooth map of manifolds,
and\/ $\bs f=F_\Man^\dMan(f)$. Then $\bs f$ is a submersion,
immersion, or embedding in $\dMan$ if and only if\/ $f$ is a
submersion, immersion, or embedding in $\Man,$ respectively. Also
$\bs f$ is a w-immersion or w-embedding if and only if\/ $f$ is an
immersion or embedding.
\smallskip

\noindent{\bf(viii)} Let\/ $\bs f:\bX\ra\bY$ be a $1$-morphism of
d-manifolds, with\/ $\bY$ a manifold. Then $\bs f$ is a
w-submersion.
\smallskip

\noindent{\bf(ix)} Let\/ $\bX,\bY$ be d-manifolds, with\/ $\bY$ a
manifold. Then $\bs\pi_\bX:\bX\t\bY\ra\bX$ is a submersion.
\smallskip

\noindent{\bf(x)} Let\/ $\bs f:\bX\ra\bY$ be a submersion of
d-manifolds, and\/ $x\in\bX$ with\/ $\bs f(x)=y\in\bY$. Then there
exist open $x\in\bU\subseteq\bX$ and\/ $y\in\bV\subseteq\bY$ with\/
$\bs f(\bU)=\bV,$ a manifold\/ $\bZ,$ and an equivalence $\bs
i:\bU\ra\bV\t\bZ,$ such that\/ $\bs f\vert_\bU:\bU\ra\bV$ is
$2$-isomorphic to $\bs\pi_\bV\ci\bs i,$ where
$\bs\pi_\bV:\bV\t\bZ\ra\bV$ is the projection.
\smallskip

\noindent{\bf(xi)} Let\/ $\bs f:\bX\ra\bY$ be a submersion of
d-manifolds with\/ $\bY$ a manifold. Then\/ $\bX$ is a
manifold.\I{d-manifold!w-submersion|)}\I{d-manifold!submersion|)}
\I{d-manifold!w-immersion|)}\I{d-manifold!immersion|)}
\I{d-manifold!w-embedding|)}\I{d-manifold!embedding|)}
\label{ds4thm6}
\end{thm}

\subsection{D-transversality and fibre products}
\label{ds46}
\I{d-transversality|(} \I{d-manifold!d-transverse 1-morphisms|(}
\I{d-manifold!fibre products|(}

From \S\ref{ds33}, if $\bs g:\bX\ra\bZ$ and $\bs h:\bY\ra\bZ$ are
1-morphisms of d-manifolds then a fibre product $\bW=\bX_{\bs
g,\bZ,\bs h}\bY$ exists in $\dSpa$, and is unique up to equivalence.
We want to know whether $\bW$ is a d-manifold. We will define when
$\bs g,\bs h$ are {\it d-transverse}, which is a sufficient
condition for $\bW$ to be a d-manifold.

Recall that if $g:X\ra Z$, $h:Y\ra Z$ are smooth maps of manifolds,
then a fibre product $W=X\t_{g,Z,h}Y$ in $\Man$ exists if $g,h$ are
{\it transverse},\I{manifold!transverse fibre products} that is, if
$T_zZ=\d g\vert_x(T_xX)+\d h\vert_y(T_yY)$ for all $x\in X$ and
$y\in Y$ with $g(x)=h(y)=z\in Z$. Equivalently, $\d g\vert_x^*\op\d
h\vert_y^*:T_zZ^*\ra T_x^*X\op T_y^*Y$ should be injective. Writing
$W=X\t_ZY$ for the topological fibre product and $e:W\ra X$, $f:W\ra
Y$ for the projections, with $g\ci e=h\ci f$, we see that $g,h$ are
transverse if and only if
\e
e^*(\d g^*)\op f^*(\d h^*):(g\ci e)^*(T^*Z)\ra e^*(T^*X)\op
f^*(T^*Y)
\label{ds4eq7}
\e
is an injective morphism of vector bundles on the topological space
$W$, that is, it has a left inverse. The condition that \eq{ds4eq8}
has a left inverse is an analogue of this, but on (dual) obstruction
rather than cotangent bundles.

\begin{dfn} Let $\bX,\bY,\bZ$ be d-manifolds and $\bs g:\bX\ra\bZ,$
$\bs h:\bY\ra\bZ$ be 1-morphisms. Let $\uW=\uX\t_{\ug,\uZ,\uh}\uY$
be the $C^\iy$-scheme fibre product, and write $\ue:\uW\ra\uX$,
$\uf:\uW\ra\uY$ for the projections. Consider the morphism
\e
\al=\begin{pmatrix} \ue^*(g'')\ci I_{\ue,\ug}(\EZ) \\
-\uf^*(h'')\ci I_{\uf,\uh}(\EZ) \\
(\ug\ci\ue)^*(\phi_Z)\end{pmatrix}:\begin{aligned}[t]
&(\ug\ci\ue)^*(\EZ)\longra\\
&\quad\ue^*(\EX)\op\uf^*(\EY)\op(\ug\ci\ue)^*(\FZ)\end{aligned}
\label{ds4eq8}
\e
in $\qcoh(\uW)$. We call $\bs g,\bs h$ {\it d-transverse\/} if $\al$
has a left inverse. Note that this is a local condition in $\uW$,
since local choices of left inverse for $\al$ can be combined using
a partition of unity\I{partition of unity} on $\uW$ to make a global
left inverse.

In the notation of \S\ref{ds43} and \S\ref{ds45}, we have
1-morphisms $\Om_{\bs g}:\ug^*(T^*\bZ)\ra T^*\bX$ in $\vvect(\uX)$
and $\Om_{\bs h}:\uh^*(T^*\bZ)\ra T^*\bY$ in $\vvect(\uY)$. Pulling
these back to $\vvect(\uW)$ using $\ue^*,\uf^*$ we form the
1-morphism in~$\vvect(\uW)$:
\e
\begin{split}
\bigl(\ue^*(\Om_{\bs g})\ci I_{\ue,\ug}(T^*\bZ)\bigr) \op
\bigl(\uf^*(\Om_{\bs h})\ci I_{\uf,\uh}(T^*\bZ)\bigr):
(\ug\ci\ue)^*(T^*\bZ)&\\
\longra\ue^*(T^*\bX)\op\uf^*(T^*&\bY).
\end{split}
\label{ds4eq9}
\e
For \eq{ds4eq8} to have a left inverse is equivalent to \eq{ds4eq9}
being weakly injective, as in Definition \ref{ds4def6}. This is the
d-manifold analogue of \eq{ds4eq7} being injective.
\label{ds4def8}
\end{dfn}

Here are the main results of~\cite[\S 4.3]{Joyc6}:

\begin{thm} Suppose\/ $\bX,\bY,\bZ$ are d-manifolds and\/ $\bs
g:\bX\ra\bZ,$ $\bs h:\bY\ra\bZ$ are d-transverse $1$-morphisms, and
let\/ $\bW=\bX\t_{\bs g,\bZ,\bs h}\bY$ be the d-space fibre product.
Then $\bW$ is a d-manifold, with\I{d-manifold!fibre
products!d-transverse}\I{d-manifold!virtual dimension}\I{fibre
product!of d-manifolds}
\e
\vdim\bW=\vdim\bX+\vdim\bY-\vdim\bZ.
\label{ds4eq10}
\e
\label{ds4thm7}
\end{thm}

\begin{thm} Suppose\/ $\bs g:\bX\ra\bZ,$ $\bs h:\bY\ra\bZ$ are
$1$-morphisms of d-manifolds. The following are sufficient
conditions for $\bs g,\bs h$ to be d-transverse, so that\/
$\bW=\bX\t_{\bs g,\bZ,\bs h}\bY$ is a d-manifold of virtual
dimension\/~{\rm\eq{ds4eq10}:}
\begin{itemize}
\setlength{\itemsep}{0pt}
\setlength{\parsep}{0pt}
\item[{\bf(a)}] $\bZ$ is a manifold, that is, $\bZ\in\hMan;$ or
\item[{\bf(b)}] $\bs g$ or $\bs h$ is a
w-submersion.\I{d-manifold!is a
manifold}\I{d-manifold!w-submersion}
\end{itemize}
\label{ds4thm8}
\end{thm}

The point here is that roughly speaking, $\bs g,\bs h$ are
d-transverse if they map the direct sum of the obstruction spaces of
$\bX,\bY$ surjectively onto the obstruction spaces of $\bZ$. If
$\bZ$ is a manifold its obstruction spaces are zero. If $\bs g$ is a
w-submersion it maps the obstruction spaces of $\bX$ surjectively
onto the obstruction spaces of $\bZ$. In both cases,
d-transversality follows. See \cite[Th.~8.15]{Spiv} for the analogue
of Theorem \ref{ds4thm8}(a) for Spivak's derived
manifolds.\I{Spivak's derived manifolds}

\begin{thm} Let\/ $\bX,\bZ$ be d-manifolds, $\bY$ a manifold, and\/
$\bs g:\bX\ra\bZ,$ $\bs h:\bY\ra\bZ$ be $1$-morphisms with\/ $\bs g$
a submersion. Then\/ $\bW=\bX\t_{\bs g,\bZ,\bs h}\bY$ is a manifold,
with\/~$\dim\bW=\vdim\bX+\dim\bY-\vdim\bZ$.\I{d-manifold!submersion}
\label{ds4thm9}
\end{thm}

Theorem \ref{ds4thm9} shows that we may think of submersions as
`representable 1-morphisms' in $\dMan$. We can locally characterize
embeddings and immersions in $\dMan$ in terms of fibre products with
$\bR^{\bs n}$ in~$\dMan$.

\begin{thm}{\bf(i)} Let\/ $\bX$ be a d-manifold and\/ $\bs
g:\bX\ra\bR^{\bs n}$ a $1$-morphism in $\dMan$. Then the fibre
product\/ $\bW=\bX\t_{\bs g,\bR^{\bs n},\bs 0}\bs *$ exists in
$\dMan$ by Theorem\/ {\rm\ref{ds4thm8}(a),} and the projection
$\bs\pi_\bX:\bW\ra\bX$ is an embedding.
\smallskip

\noindent{\bf(ii)} Suppose\/ $\bs f:\bX\ra\bY$ is an immersion of
d-manifolds, and\/ $x\in\bX$ with\/ $\bs f(x)=y\in\bY$. Then there
exist open d-submanifolds $x\in\bU\subseteq\bX$ and\/
$y\in\bV\subseteq\bY$ with\/ $\bs f(\bU)\subseteq\bV,$ and a
$1$-morphism $\bs g:\bV\ra\bR^{\bs n}$ with\/ $\bs g(y)=0,$ where
$n=\vdim\bY-\vdim\bX\ge 0,$ fitting into a
$2$-Cartesian\I{2-category!2-Cartesian square} square in~$\dMan:$
\begin{equation*}
\xymatrix@C=60pt@R=10pt{\bU \ar[d]^{\bs f\vert_\bU} \ar[r]
\drtwocell_{}\omit^{}\omit{^{}} & \bs{*} \ar[d]_{\bs 0} \\
\bV \ar[r]^(0.7){\bs g} & \bR^{\bs n}. }
\end{equation*}
If\/ $\bs f$ is an embedding we may take\/~$\bU=\bs
f^{-1}(\bV)$.\I{d-manifold!immersion}\I{d-manifold!embedding}
\label{ds4thm10}
\end{thm}

\begin{rem} For the applications the author has in
mind, it will be crucial that if $\bs g:\bX\ra\bZ$ and $\bs
h:\bY\ra\bZ$ are 1-morphisms with $\bX,\bY$ d-manifolds and $\bZ$ a
manifold then $\bW=\bX\t_\bZ\bY$ is a d-manifold, with
$\vdim\bW=\vdim\bX+\vdim\bY-\dim\bZ$, as in Theorem
\ref{ds4thm8}(a). We will show by example, following Spivak
\cite[Prop.~1.7]{Spiv}, that if d-manifolds $\dMan$ were an ordinary
category containing manifolds as a full subcategory, then this would
be false.

Consider the fibre product $\bs{*}\t_{\bs 0,\bR,\bs 0}\bs{*}$ in
$\dMan$. If $\dMan$ were a category then as $\bs{*}$ is a terminal
object, the fibre product would be $\bs{*}$. But then
\begin{equation*}
\vdim(\bs{*}\t_{\bs 0,\bR,\bs 0}\bs{*})=\vdim\bs{*}=0\ne -1=\vdim\bs{*}+
\vdim\bs{*}-\vdim\bR,
\end{equation*}
so equation \eq{ds4eq10} and Theorem \ref{ds4thm8}(a) would be
false. Thus, if we want fibre products of d-manifolds over manifolds
to be well behaved, then $\dMan$ must be at least a 2-category. It
could be an $\iy$-category,\I{$\iy$-category} as for Spivak's
derived manifolds\I{Spivak's derived manifolds} \cite{Spiv}, or some
other kind of higher category. Making d-manifolds into a 2-category,
as we have done, is the simplest of the available
options.\I{d-manifold!why $\dMan$ is a 2-category}
\I{d-transversality|)} \I{d-manifold!d-transverse 1-morphisms|)}
\I{d-manifold!fibre products|)}
\label{ds4rem}
\end{rem}

\subsection{Embedding d-manifolds into manifolds}
\label{ds47}
\I{principal
d-manifold|(}\I{d-manifold!principal|(}\I{d-manifold!embedding!into
manifolds|(}

Let $V$ be a manifold, $E\ra V$ a vector bundle, and $s\in
C^\iy(E)$. Then Example \ref{ds4ex2} defines a `standard model'
principal d-manifold $\bS_{V,E,s}$. When $E$ and $s$ are zero, we
have $\bS_{V,0,0}=\bV=F_\Man^\dMan(V)$, so that $\bS_{V,0,0}$ is a
manifold. For general $V,E,s$, taking $f=\id_V:V\ra V$ and $\hat
f=0:E\ra 0$ in Example \ref{ds4ex3} gives a `standard model'
1-morphism $\bS_{\id_V,0}:\bS_{V,E,s}\ra \bS_{V,0,0}=\bV$. One can
show $\bS_{\id_V,0}$ is an embedding, in the sense of Definition
\ref{ds4def7}. Any principal d-manifold $\bU$ is equivalent to some
$\bS_{V,E,s}$. Thus we deduce:

\begin{lem} Any principal d-manifold\/ $\bU$ admits an embedding\/
$\bs i:\bU\ra\bV$ into a manifold\/~$\bV$.
\label{ds4lem1}
\end{lem}

Theorem \ref{ds4thm13} below is a converse to this: if a d-manifold
$\bX$ can be embedded into a manifold $\bY$, then $\bX$ is
principal. So it will be useful to study embeddings of d-manifolds
into manifolds. The following classical facts are due to
Whitney~\cite{Whit}.

\begin{thm}{\bf(a)} Let\/ $X$ be an $m$-manifold and\/ $n\ge 2m$.
Then a generic smooth map $f:X\ra\R^n$ is an
immersion.\I{manifold!immersion}
\smallskip

\noindent{\bf(b)} Let\/ $X$ be an $m$-manifold and\/ $n\ge 2m+1$.
Then there exists an embedding $f:X\ra\R^n,$ and we can choose
such\/ $f$ with\/ $f(X)$ closed in\/ $\R^n$. Generic smooth maps\/
$f:X\ra\R^n$ are embeddings.\I{manifold!embedding}
\label{ds4thm11}
\end{thm}

In \cite[\S 4.4]{Joyc6} we generalize Theorem \ref{ds4thm11} to
d-manifolds.

\begin{thm} Let\/ $\bX$ be a d-manifold. Then there exist
immersions and/or embeddings $\bs f:\bX\ra\bR^{\bs n}$ for some
$n\gg 0$ if and only if there is an upper bound for\/ $\dim
T^*_x\uX$ for all\/ $x\in\uX$. If there is such an upper bound, then
immersions $\bs f:\bX\ra\bR^{\bs n}$ exist provided\/ $n\ge 2\dim
T_x^*\uX$ for all\/ $x\in\uX,$ and embeddings $\bs f:\bX\ra\bR^{\bs
n}$ exist provided\/ $n\ge 2\dim T_x^*\uX+1$ for all\/ $x\in\uX$.
For embeddings we may also choose $\bs f$ with\/ $f(X)$ closed
in\/~$\R^n$.\I{d-manifold!immersion}
\label{ds4thm12}
\end{thm}

Here is an example in which the condition does not hold.

\begin{ex} $\bR^{\bs k}\t_{\bs 0,\bR^{\bs k},\bs 0}\bs *$ is a
principal d-manifold of virtual dimension 0, with $C^\iy$-scheme
$\ul{\R}^k$, and obstruction bundle $\R^k$. Thus $\bX=\coprod_{k\ge
0}\bR^{\bs k}\t_{\bs 0,\bR^{\bs k},\bs 0}\bs *$ is a d-manifold of
virtual dimension 0, with $C^\iy$-scheme $\uX=\coprod_{k\ge
0}\ul{\R}^k$. Since $T^*_x\uX\cong\R^n$ for
$x\in\R^n\subset\coprod_{k\ge 0}\R^k$, $\dim T^*_x\uX$ realizes all
values $n\ge 0$. Hence there cannot exist immersions or embeddings
$\bs f:\bX\ra\bR^{\bs n}$ for any~$n\ge 0$.
\label{ds4ex5}
\end{ex}\I{d-manifold!example which is not principal}

As $x\mapsto\dim T_x^*\uX$ is an upper semicontinuous map $X\ra\N$,
if $\bX$ is compact then $\dim T_x^*\uX$ is bounded above, giving:

\begin{cor} Let\/ $\bX$ be a compact d-manifold. Then there exists
an embedding\/ $\bs f:\bX\ra\bR^{\bs n}$ for some\/~$n\gg 0$.
\label{ds4cor1}
\end{cor}

If a d-manifold $\bX$ can be embedded into a manifold $Y$, we show
in \cite[\S 4.4]{Joyc6} that we can write $\bX$ as the zeroes of a
section of a vector bundle over $Y$ near its image. See
\cite[Prop.~9.5]{Spiv} for the analogue for Spivak's derived
manifolds.\I{Spivak's derived manifolds}

\begin{thm} Suppose $\bX$ is a d-manifold, $Y$ a manifold, and\/
$\bs f:\bX\ra\bY$ an embedding, in the sense of Definition\/
{\rm\ref{ds4def7}}. Then there exist an open subset $V$ in $Y$ with
$\bs f(\bX)\subseteq\bV,$ a vector bundle $E\ra V,$ and\/ $s\in
C^\iy(E)$ fitting into a $2$-Cartesian\I{2-category!2-Cartesian
square} diagram in $\dSpa\!:$
\begin{equation*}
\xymatrix@C=60pt@R=10pt{ \bX \ar[r]_(0.25){\bs f} \ar[d]^{\bs f}
\drtwocell_{}\omit^{}\omit{^{\eta}}
 & \bV \ar[d]_{\bs 0} \\ \bV \ar[r]^(0.7){\bs s} & \bE.}
\end{equation*}
Here \/ $\bY=F_\Man^\dMan(Y),$ and similarly for $\bV,\bE,\bs s,\bs
0,$ with\/ $0:V\ra E$ the zero section. Hence $\bX$ is equivalent to
the `standard model' d-manifold $\bS_{V,E,s}$ of Example\/
{\rm\ref{ds4ex2},} and is a principal d-manifold.
\label{ds4thm13}
\end{thm}\I{d-manifold!standard model}

Combining Theorems \ref{ds4thm12} and \ref{ds4thm13}, Lemma
\ref{ds4lem1}, and Corollary \ref{ds4cor1} yields:

\begin{cor} Let\/ $\bX$ be a d-manifold. Then $\bX$ is a principal
d-manifold if and only if\/ $\dim T^*_x\uX$ is bounded above for
all\/ $x\in\uX$. In particular, if\/ $\bX$ is compact, then $\bX$ is
principal.
\label{ds4cor2}
\end{cor}

Corollary \ref{ds4cor2} suggests that most interesting d-manifolds
are principal, in a similar way to most interesting $C^\iy$-schemes
being affine in Remark \ref{ds2rem1}(ii). Example \ref{ds4ex5} gives
a d-manifold which is not principal.\I{principal
d-manifold|)}\I{d-manifold!principal|)}\I{d-manifold!embedding!into
manifolds|)}

\subsection{Orientations on d-manifolds}
\label{ds48}
\I{d-manifold!orientations|(}

Let $X$ be an $n$-manifold. Then $T^*X$ is a rank $n$ vector bundle
on $X$, so its top exterior power $\La^nT^*X$ is a line bundle (rank
1 vector bundle) on $X$. In algebraic geometry, $\La^nT^*X$ would be
called the canonical bundle of $X$. We define an {\it orientation\/}
$\om$ on $X$ to be an {\it orientation on the fibres of\/}
$\La^nT^*X$.\I{manifold!orientation} That is, $\om$ is an
equivalence class $[\tau]$ of isomorphisms of line bundles
$\tau:O_X\ra\La^nT^*X$, where $O_X$ is the trivial line bundle $\R\t
X\ra X$, and $\tau,\tau'$ are equivalent if $\tau'=\tau\cdot c$ for
some smooth~$c:X\ra(0,\iy)$.

To generalize all this to d-manifolds, we will need a notion of the
`top exterior power' $\cL_{\smash{ (\cE^\bu,\phi)}}$ of a virtual
vector bundle $(\cE^\bu,\phi)$ in \S\ref{ds43}. As the definition in
\cite[\S 4.5]{Joyc6} is long, we will not give it, but just state
its important properties:\G[LE]{$\cL_{\smash{
(\cE^\bu,\phi)}}$}{orientation line bundle of a virtual vector
bundle $(\cE^\bu,\phi)$}

\begin{thm} Let\/ $\uX$ be a $C^\iy$-scheme, and\/ $(\cE^\bu,\phi)$
a virtual vector bundle on $\uX$. Then in {\rm\cite[\S 4.5]{Joyc6}}
we define a line bundle (rank\/ $1$ vector bundle)\I{virtual vector
bundle!orientation line bundle of} $\cL_{\smash{(\cE^\bu,\phi)}}$ on
$\uX,$ which we call the \begin{bfseries}orientation line
bundle\end{bfseries}\I{orientation line bundle} of\/
$(\cE^\bu,\phi)$. This satisfies:
\begin{itemize}
\setlength{\itemsep}{0pt}
\setlength{\parsep}{0pt}
\item[{\bf(a)}] Suppose $\cE^1,\cE^2$ are vector bundles on $\uX$ with
ranks $k_1,k_2,$ and\/ $\phi:\cE^1\ra\cE^2$ is a morphism. Then
$(\cE^\bu,\phi)$ is a virtual vector bundle of rank\/ $k_2-k_1,$
and there is a canonical isomorphism\/~$\cL_{\smash{
(\cE^\bu,\phi)}}\cong \La^{k_1}(\cE^1)^*\ot\La^{k_2}\cE^2$.
\item[{\bf(b)}] Let\/ $f^\bu: (\cE^\bu,\phi)\ra(\cF^\bu,\psi)$
be an equivalence in $\vvect(\uX)$. Then there is a canonical
isomorphism\/ $\cL_{f^\bu}:\cL_{\smash{(\cE^\bu,\phi)}}\ra
\cL_{\smash{ (\cF^\bu,\psi)}}$ in\/~$\qcoh(\uX)$.
\item[{\bf(c)}] If\/ $(\cE^\bu,\phi)\in\vvect(\uX)$
then\/~$\cL_{\id_\phi}=\id_{\cL_{\smash{(\cE^\bu,\phi)}}}:
\cL_{\smash{(\cE^\bu,\phi)}}\ra\cL_{\smash{(\cE^\bu,\phi)}}$.
\item[{\bf(d)}] If\/ $f^\bu:(\cE^\bu,\phi)\ra(\cF^\bu,\psi)$ and\/
$g^\bu:(\cF^\bu,\psi)\ra(\cG^\bu,\xi)$ are equivalences in
$\vvect(\uX)$ then $\cL_{g^\bu\ci f^\bu}=\cL_{g^\bu}\ci\cL_{
f^\bu}:\cL_{\smash{(\cE^\bu,\phi)}}\ra\cL_{\smash{(\cG^\bu,\xi)}}$.
\item[{\bf(e)}] If\/ $f^\bu,g^\bu:(\cE^\bu,\phi)\ra(\cF^\bu,\psi)$
are $2$-isomorphic equivalences in $\vvect(\uX)$ then\/~$\cL_{
f^\bu}=\cL_{g^\bu} :\cL_{\smash{(\cE^\bu,\phi)}}\ra\cL_{\smash{
(\cF^\bu,\psi)}}$.
\item[{\bf(f)}] Let\/ $\uf:\uX\ra\uY$ be a morphism of\/
$C^\iy$-schemes, and\/ $(\cE^\bu,\phi)\in\vvect(\uY)$. Then
there is a canonical isomorphism\/~$I_{\uf,(\cE^\bu,\phi)}:
\uf^*(\cL_{\smash{(\cE^\bu, \phi)}})\ra
\cL_{\smash{\uf^*(\cE^\bu,\phi)}}$.
\end{itemize}
\label{ds4thm14}
\end{thm}

Now we can define orientations on d-manifolds.

\begin{dfn} Let $\bX$ be a d-manifold. Then the virtual cotangent
bundle\I{virtual cotangent bundle} $T^*\bX$ is a virtual vector
bundle on $\uX$ by Proposition \ref{ds4prop2}(b), so Theorem
\ref{ds4thm14} gives a line bundle $\cL_{T^*\bX}$ on $\uX$. We call
$\cL_{T^*\bX}$\G[LT*Xa]{$\cL_{T^*\bX}$}{orientation line bundle of a
d-manifold $\bX$} the {\it orientation line bundle\/}\I{orientation
line bundle} of~$\bX$.\I{d-manifold!orientation line bundle}

An {\it orientation\/} $\om$ on $\bX$ is an orientation on
$\cL_{T^*\bX}$. That is, $\om$ is an equivalence class $[\tau]$ of
isomorphisms $\tau:\OX\ra\cL_{T^*\bX}$ in $\qcoh(\uX)$, where
$\tau,\tau'$ are equivalent if they are proportional by a smooth
positive function on~$\uX$.

If $\om=[\tau]$ is an orientation on $\bX$, the {\it opposite
orientation\/} is $-\om=[-\tau]$, which changes the sign of the
isomorphism $\tau:\O_X\ra\cL_{T^*\bX}$. When we refer to $\bX$ as an
oriented d-manifold, $-\bX$ will mean $\bX$ with the opposite
orientation, that is, $\bX$ is short for $(\bX,\om)$ and $-\bX$ is
short for~$(\bX,-\om)$.
\label{ds4def9}
\end{dfn}

\begin{ex}{\bf(a)} Let $X$ be an $n$-manifold, and
$\bX=F_\Man^\dMan(X)$ the associated d-manifold. Then
$\uX=F_\Man^\CSch(X)$, $\EX=0$ and $\FX=T^*\uX$. So $\EX,\FX$ are
vector bundles of ranks $0,n$. As $\La^0\EX\cong\O_X$, Theorem
\ref{ds4thm14}(a) gives a canonical isomorphism
$\cL_{T^*\bX}\cong\La^nT^*\uX$. That is, $\cL_{T^*\bX}$ is
isomorphic to the lift to $C^\iy$-schemes of the line bundle
$\La^nT^*X$ on the manifold~$X$.

As above, an orientation on $X$ is an orientation on the line bundle
$\La^nT^*X$. Hence orientations on the d-manifold
$\bX=F_\Man^\dMan(X)$ in the sense of Definition \ref{ds4def9} are
equivalent to orientations on the manifold $X$ in the usual sense.
\smallskip

\noindent{\bf(b)} Let $V$ be an $n$-manifold, $E\ra V$ a vector
bundle of rank $k$, and $s\in C^\iy(E)$. Then Example \ref{ds4ex2}
defines a `standard model' principal d-manifold $\bS=\bS_{V,E,s}$,
which has $\ES\cong\cE^*\vert_\uS$, $\FS\cong T^*\uV\vert_\uS$,
where $\cE,T^*\uV$ are the lifts of the vector bundles $E,T^*V$ on
$V$ to $\uV$. Hence $\ES,\FS$ are vector bundles on $\uS_{V,E,s}$ of
ranks $k,n$, so Theorem \ref{ds4thm14}(a) gives an
isomorphism~$\cL_{T^*\bS_{V,E,s}}\cong(\La^k\cE\ot
\La^nT^*\uV)\vert_\uS$.\I{d-manifold!standard model!orientations on}

Thus $\cL_{T^*\bS_{V,E,s}}$ is the lift to $\uS_{V,E,s}$ of the line
bundle $\La^kE\ot\La^nT^*V$ over the manifold $V$. Therefore we may
induce an orientation on the d-manifold $\bS_{V,E,s}$ from an
orientation on the line bundle $\La^kE\ot\La^nT^*V$ over $V$.
Equivalently, we can induce an orientation on $\bS_{V,E,s}$ from an
orientation on the total space of the vector bundle $E^*$ over $V$,
or from an orientation on the total space of~$E$.
\label{ds4ex6}
\end{ex}

We can construct orientations on d-transverse fibre products of
oriented d-manifolds. Note that \eq{ds4eq11} depends on an {\it
orientation convention\/}:\I{orientation convention|(} a different
choice would change \eq{ds4eq11} by a sign depending on
$\vdim\bX,\vdim\bY,\vdim\bZ$. Our conventions follow those of Fukaya
et al.\ \cite[\S 8.2]{FOOO} for Kuranishi spaces.\I{d-manifold!fibre
products!orientations on|(}\I{Kuranishi space}

\begin{thm} Work in the situation of Theorem\/ {\rm\ref{ds4thm7},} so
that\/ $\bW,\bX,\bY,\bZ$ are d-manifolds with\/ $\bW=\bX\t_{\bs
g,\bZ,\bs h}\bY$ for $\bs g,\bs h$ d-transverse, where $\bs
e:\bW\ra\bX,$ $\bs f:\bW\ra\bY$ are the projections. Then we have
orientation line bundles\I{orientation line bundle}
$\cL_{T^*\bW},\ldots,\cL_{T^*\bZ}$ on $\uW,\ldots,\uZ,$ so\/
$\cL_{T^*\bW},\ue^*(\cL_{T^*\bX}),
\uf^*(\cL_{T^*\bY}),(\ug\ci\ue)^*(\cL_{T^*\bZ})$ are line bundles on
$\uW$. With a suitable choice of orientation
convention,\I{orientation convention|)} there is a canonical
isomorphism
\e
\Phi:\cL_{T^*\bW}\longra\ue^*(\cL_{T^*\bX})\ot_{\O_W}\uf^*
(\cL_{T^*\bY})\ot_{\O_W}(\ug\ci\ue)^*(\cL_{T^*\bZ})^*.
\label{ds4eq11}
\e

Hence, if\/ $\bX,\bY,\bZ$ are oriented d-manifolds, then $\bW$ also
has a natural orientation, since trivializations of\/
$\cL_{T^*\bX},\cL_{T^*\bY},\cL_{T^*\bZ}$ induce a trivialization
of\/ $\cL_{T^*\bW}$ by\/~\eq{ds4eq11}.
\label{ds4thm15}
\end{thm}

Fibre products have natural commutativity and associativity
properties. When we include orientations, the orientations differ by
some sign. Here is an analogue of results of Fukaya et al.\
\cite[Lem.~8.2.3]{FOOO} for Kuranishi spaces.\I{Kuranishi space}

\begin{prop} Suppose $\bV,\ldots,\bZ$ are oriented d-manifolds,
$\bs e,\ldots,\bs h$ are $1$-morphisms, and all fibre products below
are d-transverse. Then the following hold, in oriented d-manifolds:
\begin{itemize}
\setlength{\itemsep}{0pt}
\setlength{\parsep}{0pt}
\item[{\bf(a)}] For\/ $\bs g:\bX\ra\bZ$ and\/
$\bs h:\bY\ra\bZ$ we have
\begin{equation*}
\bX\t_{\bs g,\bZ,\bs h}\bY\simeq(-1)^{(\vdim \bX-\vdim\bZ)(\vdim
\bY-\vdim \bZ)}\bY\t_{\bs h,\bZ,\bs g}\bX.
\end{equation*}
In particular, when $\bZ=\bs{*}$ so that\/
$\bX\t_\bZ\bY=\bX\t\bY$ we have
\begin{equation*}
\bX\t\bY\simeq(-1)^{\vdim \bX\vdim \bY}\bY\t\bX.
\end{equation*}
\item[{\bf(b)}] For\/ $\bs e:\bV\ra\bY,$ $\bs f:\bW\ra\bY,$ $\bs
g:\bW\ra\bZ,$ and\/ $\bs h:\bX\ra\bZ$ we have
\begin{equation*}
\bV\t_{\bs e,\bY,\bs f\ci\bs\pi_\bW}\bigl(\bW\t_{\bs g,\bZ,\bs
h}\bX\bigr)\simeq \bigl(\bV\t_{\bs e,\bY,\bs f}\bW\bigr)\t_{\bs
g\ci\bs\pi_\bW,\bZ,\bs h}\bX.
\end{equation*}
\item[{\bf(c)}] For\/ $\bs e:\bV\ra\bY,$ $\bs f:\bV\ra\bZ,$
$\bs g:\bW\ra\bY,$ and\/ $\bs h:\bX\ra\bZ$ we
have\I{d-manifold!fibre products!orientations on|)}
\I{d-manifold|)}\I{d-manifold!orientations|)}
\begin{align*}
&\bV\t_{(\bs e,\bs f),\bY\t\bZ,\bs g\t\bs h}(\bW\t\bX)\simeq \\
&\quad(-1)^{\vdim\bZ(\vdim\bY+\vdim\bW)} (\bV\t_{\bs e,\bY,\bs
g}\bW)\t_{\bs f\ci\bs\pi_\bV,\bZ,\bs h}\bX.
\end{align*}
\end{itemize}
\label{ds4prop4}
\end{prop}

\section[Manifolds with boundary and manifolds with
corners]{Manifolds with boundary and corners}
\label{ds5}
\I{manifold!with boundary|see{manifold with
boundary}}\I{manifold!with corners|see{manifold with \\ corners}}
\I{manifold with boundary|(}\I{manifold with corners|(}

So far we have discussed only manifolds {\it without boundary\/}
(locally modelled on $\R^n$). One can also consider {\it manifolds
with boundary\/} (locally modelled on $[0,\iy)\t\R^{n-1}$) and {\it
manifolds with corners} (locally modelled on $[0,\iy)^k\t\R^{n-k}$).
The author \cite{Joyc3} studied manifolds with corners, giving a new
definition of smooth map $f:X\ra Y$ between manifolds with corners
$X,Y$, satisfying extra conditions over $\pd^kX,\pd^lY$. This yields
categories $\Manb,\Manc$ of manifolds with boundary and with corners
with good properties {\it as categories}.

In \cite[Chap.~5]{Joyc6} we surveyed \cite{Joyc3}, changing some
notation, and including some new material. This section summarizes
\cite{Joyc3}, \cite[Chap.~5]{Joyc6}, following the notation of
\cite[Chap.~5]{Joyc6}. See \cite{Joyc3} and \cite[Chap.~5]{Joyc6}
for further references on manifolds with corners.

\subsection{Boundaries and smooth maps}
\label{ds51}

The definition of an $n$-manifold with corners $X$ in \cite[\S
2]{Joyc3} involves an atlas of charts $(U,\phi)$ on $X$ with
$U\subseteq[0,\iy)^k\t\R^{n-k}$ open and $\phi:U\hookra X$ a
homeomorphism with an open set in $X$. Apart from taking
$U\subseteq[0,\iy)^k\t\R^{n-k}$ rather than $U\subseteq\R^n$, there
is no difference with the usual definition of $n$-manifold. The
definitions of the {\it boundary\/} $\pd X$ of $X$ in \cite[\S
2]{Joyc3}, and of {\it smooth map\/} $f:X\ra Y$ between manifolds
with corners in \cite[\S 3]{Joyc3}, may be surprising for readers
who have not thought much about corners, so we give them here.

\begin{dfn} Let $X$ be a manifold with corners, of dimension $n$.
Then there is a natural stratification $X=\coprod_{k=0}^nS^k(X)$,
where $S^k(X)$ is the {\it depth\/ $k$ stratum} of $X$, that is, the
set of points $x\in X$ such that $X$ near $x$ is locally modelled on
$[0,\iy)^k\t\R^{n-k}$ near 0. Then $S^k(X)$ is an $(n-k)$-manifold
without boundary, and $\ov{S^k(X)}=\coprod_{l=k}^nS^l(X)$. The {\it
interior\/} of $X$ is~$X^\ci=S^0(X)$.\G[SkX]{$S^k(X)$}{depth $k$
stratum of a manifold with corners $X$}

A {\it local boundary component\/ $\be$ of\/ $X$ at\/} $x$ is a
local choice of connected component of $S^1(X)$ near $x$. That is,
for each sufficiently small open neighbourhood $V$ of $x$ in $X$,
$\be$ gives a choice of connected component $W$ of $V\cap S^1(X)$
with $x\in\overline W$, and any two such choices $V,W$ and $V',W'$
must be compatible in the sense that $x\in\overline{(W\cap W')}$. As
a set, define the {\it boundary\/}\I{manifold with
corners!boundary}\I{boundary!of a manifold with corners}\I{manifold
with corners!local boundary component}\G[dXa]{$\pd X$}{boundary of a
manifold with corners $X$}
\begin{equation*}
\pd X=\bigl\{(x,\be):\text{$x\in X$, $\be$ is a local boundary
component for $X$ at $x$}\bigr\}.
\end{equation*}
Then $\pd X$ is an $(n-1)$-manifold with corners if $n>0$, and $\pd
X=\es$ if~$n=0$.

Define a smooth map $i_X:\pd X\ra X$ by~$i_X:(x,\be)\mapsto
x$.\G[iXa]{$i_X:\pd X\ra X$}{inclusion of boundary $\pd X$ into a
manifold with corners $X$}
\label{ds5def1}
\end{dfn}

\begin{ex} The manifold with corners $X=[0,\iy)^2$ has strata
$S^0(X)=(0,\iy)^2$, $S^1(X)=\bigl(\{0\}\t(0,\iy)\bigr)\amalg
\bigl((0,\iy)\t\{0\}\bigr)$ and $S^2(X)=\{(0,0)\}$. A point $(a,b)$
in $X$ has local boundary components $\{x=0\}$ if $a=0$ and
$\{y=0\}$ if $b=0$. Thus
\begin{align*}
\pd X&=\bigl\{\bigl((x,0),\{y=0\}\bigr):x\in[0,\iy)\bigr\}\amalg
\bigl\{\bigl((0,y),\{x=0\}\bigr):y\in[0,\iy)\bigr\}\\
&\cong [0,\iy)\amalg[0,\iy).
\end{align*}

Note that $i_X:\pd X\ra X$  maps two points
$\bigl((0,0),\{x=0\}\bigr),\bigl((0,0),\{y=0\}\bigr)$ to $(0,0)$. In
general, if a manifold with corners $X$ has $\pd^2X\ne\es$ then
$i_X$ {\it is not injective}, so {\it the boundary $\pd X$ is not a
subset of\/}~$X$.
\label{ds5ex1}
\end{ex}

\begin{dfn} Let $X,Y$ be manifolds with corners of dimensions
$m,n$. A continuous map $f:X\ra Y$ is called {\it weakly
smooth\/}\I{manifold with corners!weakly smooth map} if whenever
$(U,\phi),(V,\psi)$ are charts on $X,Y$ then
\begin{equation*}
\psi^{-1}\ci f\ci\phi:(f\ci\phi)^{-1}(\psi(V))\longra V
\end{equation*}
is a smooth map from $(f\ci\phi)^{-1}(\psi(V))\subset\R^m$ to
$V\subset\R^n$.

Let $(x,\be)\in\pd X$. A {\it boundary defining function for\/ $X$
at\/}\I{manifold with corners!boundary defining function} $(x,\be)$
is a pair $(V,b)$, where $V$ is an open neighbourhood of $x$ in $X$
and $b:V\ra[0,\iy)$ is a weakly smooth map, such that $\d
b\vert_v:T_vV\ra T_{b(v)}[0,\iy)$ is nonzero for all $v\in V$, and
there exists an open neighbourhood $U$ of $(x,\be)$ in
$i_X^{-1}(V)\subseteq\pd X$, with $b\ci i_X\vert_U= 0$, and
$i_X\vert_U:U\longra\bigl\{v\in V:b(v)=0\bigr\}$ is a homeomorphism.

A weakly smooth map of manifolds with corners $f:X\ra Y$ is called
{\it smooth\/}\I{manifold with corners!smooth map} if it satisfies
the following additional condition over $\pd X,\pd Y$. Suppose $x\in
X$ with $f(x)=y\in Y$, and $\be$ is a local boundary component of
$Y$ at $y$. Let $(V,b)$ be a boundary defining function for $Y$ at
$(y,\be)$. We require that either:
\begin{itemize}
\setlength{\itemsep}{0pt}
\setlength{\parsep}{0pt}
\item[(i)] There exists an open $x\in\ti V\subseteq
f^{-1}(V)\subseteq X$ such that $(\ti V,b\ci f\vert_{\ti V})$ is
a boundary defining function for $X$ at $(x,\ti\be)$, for some
unique local boundary component $\ti\be$ of $X$ at $x$; or
\item[(ii)] There exists an open $x\in W\subseteq
f^{-1}(V)\subseteq X$ with~$b\ci f\vert_W= 0$.
\end{itemize}

Form the fibre products of topological spaces
\begin{align*}
\pd X\!\t_{f\ci i_X,Y,i_Y}\!\pd Y&=\bigl\{\bigl((x,\ti\be),
(y,\be)\bigr)\!\in\!\pd X\!\t\!\pd Y:f\!\ci\!
i_X(x,\ti\be)\!=\!y\!=\!i_Y(y,\be)\bigr\},\\
X\t_{f,Y,i_Y}\pd Y&=\bigl\{\bigl(x,(y,\be)\bigr)\in X\t\pd
Y:f(x)=y=i_Y(y,\be)\bigr\}.
\end{align*}
Define subsets $S_f\subseteq\pd X\t_Y\pd Y$ and $T_f\subseteq
X\t_Y\pd Y$\G[Sfa]{$S_f\subseteq\pd X\t_Y\pd Y$}{set associated to
smooth map $f:X\ra Y$ in $\Manc$}\G[Tfa]{$T_f\subseteq X\t_Y\pd
Y$}{set associated to smooth map $f:X\ra Y$ in $\Manc$} by
$\bigl((x,\ti\be),(y,\be)\bigr)\in S_f$ in case (i) above, and
$\bigl(x,(y,\be)\bigr)\in T_f$ in case (ii) above. Define maps
$s_f:S_f\ra\pd X$, $t_f:T_f\ra X$, $u_f:S_f\ra\pd Y$, $v_f:T_f\ra\pd
Y$ to be the projections from the fibre products. Then $S_f,T_f$ are
open and closed in $\pd X\t_Y\pd Y,X\t_Y\pd Y$ and have the
structure of manifolds with corners, with $\dim S_f=\dim X-1$ and
$\dim T_f=\dim X$, and $s_t,t_f,u_f,v_f$ are smooth maps with
$s_f,t_f$ \'etale.
\label{ds5def2}
\end{dfn}

We write $\Manc$ for the category of manifolds with corners, with
morphisms smooth maps, and $\Manb$ for the full subcategory of
manifolds with boundary.\G[Manb]{$\Manb$}{category of manifolds with
boundary}\G[Manc]{$\Manc$}{category of manifolds with corners}

\subsection[(Semi)simple maps, submersions, immersions,
embeddings]{(Semi)simple maps, submersions, immersions, \\ and
embeddings}
\label{ds52}
\I{manifold with corners!simple map|(} \I{manifold with
corners!semisimple map|(} \I{manifold with corners!flat map|(}
\I{manifold with corners!submersion|(}\I{manifold with
corners!s-submersion|(} \I{manifold with corners!immersion|(}
\I{manifold with corners!s-immersion|(}\I{manifold with
corners!sf-immersion|(} \I{manifold with corners!embedding|(}
\I{manifold with corners!s-embedding|(} \I{manifold with
corners!sf-embedding|(}

In \cite[\S 5.4 \& \S 5.7]{Joyc6} we define some interesting classes
of smooth maps:

\begin{dfn} Let $f:X\ra Y$ be a smooth map of manifolds with
corners.
\begin{itemize}
\setlength{\itemsep}{0pt}
\setlength{\parsep}{0pt}
\item[(a)] We call $f$ {\it simple\/} if $s_f:S_f\ra\pd X$ in
Definition \ref{ds5def2} is bijective.
\item[(b)] We call $f$ {\it semisimple\/} if $s_f:S_f\ra\pd X$ is
injective.
\item[(c)] We call $f$ {\it flat\/} if $T_f=\es$ in
Definition~\ref{ds5def2}.
\item[(d)] We call $f$ a {\it diffeomorphism\/} if it has a
smooth inverse $f^{-1}:Y\ra X$.\I{manifold with
corners!diffeomorphism}
\item[(e)] We call $f$ a {\it submersion\/} if for all $x\in
S^k(X)\subseteq Y$ with $f(x)=y\in S^l(Y)\subseteq Y$, then $\d
f\vert_x:T_xX\ra T_{f(x)}Y$ and $\d f\vert_x:T_x(S^k(X))\ra
T_{f(x)}(S^l(Y))$ are surjective. Submersions are automatically
semisimple. We call $f$ an {\it s-submersion\/} if it is a
simple submersion.
\item[(f)] We call $f$ an {\it immersion\/} if $\d f\vert_x:
T_xX\ra T_{f(x)}Y$ is injective for all $x\in X$. We call $f$ an
{\it s-immersion} (or {\it sf-immersion}) if $f$ is also simple
(or simple and flat). We call $f$ an {\it embedding\/} (or {\it
s-embedding}, or {\it sf-embedding\/}) if $f$ is an immersion
(or s-immersion, or sf-immersion), and $f:X\ra f(X)$ is a
homeomorphism with its image.
\end{itemize}
For manifolds without boundary, one considers immersed or embedded
submanifolds. Part (f) gives six different notions of submanifolds
$X$ of manifolds with corners $Y$: {\it immersed}, {\it s-immersed},
{\it sf-immersed}, {\it embedded}, {\it s-embedded\/} and {\it
sf-embedded submanifolds}.\I{manifold with corners!submanifold}
\label{ds5def3}
\end{dfn}

\begin{ex}{\bf(i)} The inclusion $i:[0,\iy)\hookra\R$ is an
embedding. It is semisimple and flat, but not simple, as
$s_i:S_i\ra\pd[0,\iy)$ maps $\es\ra\{0\}$, and is not surjective, so
$i$ is not an s- or sf-embedding. Thus $[0,\iy)$ is an embedded
submanifold of $\R$, but not an s- or sf-embedded submanifold.

\smallskip

\noindent{\bf(ii)} The map $f:[0,\iy)\ra[0,\iy)^2$ mapping
$f:x\mapsto(x,x)$ is an embedding. It is flat, but not semisimple,
as $s_f:S_f\ra\pd[0,\iy)$ maps two points to one point, and is not
injective. Hence $f$ is not an s- or sf-embedding, and
$\bigl\{(x,x):x\in[0,\iy)\bigr\}$ is an embedded submanifold of
$[0,\iy)^2$, but not s- or sf-embedded.

\smallskip

\noindent{\bf(iii)} The inclusion $i:\{0\}\hookra[0,\iy)$ has $\d
i\vert_0$ injective, so it is an embedding. It is simple, but not
flat, as $T_i=\bigl\{\bigl(0,(0,\{x=0\})\bigr)\bigr\} \ne\es$. Thus
$i$ is an s-embedding, but not an sf-embedding. Hence $\{0\}$ is an
s-embedded but not sf-embedded submanifold of~$[0,\iy)$.

\smallskip

\noindent{\bf(iv)} Let $X$ be a manifold with corners with $\pd
X\ne\es$. Then $i_X:\pd X\ra X$ is an immersion. Also
$s_{i_X}:S_{i_X}\ra\pd^2X$ is a bijection, so $i_X$ is simple,  but
$T_{i_X}\cong\pd X\ne\es$, so $i_X$ is not flat. Hence $i_X$ is an
s-immersion, but not an sf-immersion. If $\pd^2X=\es$ then $i_X$ is
an s-embedding, but not an sf-embedding.

\smallskip

\noindent{\bf(v)} Let $f:[0,\iy)\ra\R$ be smooth. Define
$g:[0,\iy)\ra[0,\iy)\t\R$ by $g(x)=(x,f(x))$. Then $g$ is an
sf-embedding, and $\Ga_f=\bigl\{(x,f(x):x\in[0,\iy)\bigr\}$ is an
sf-embedded submanifold of~$[0,\iy)\t\R$.
\label{ds5ex2}
\end{ex}

Simple and semisimple maps have a property of lifting to boundaries:

\begin{prop} Let\/ $f:X\ra Y$ be a semisimple map of manifolds with
corners. Then there exists a natural decomposition $\pd
X=\pd^f_+X\amalg\pd^f_-X$ with\/ $\pd^f_\pm X$ open and closed in
$\pd X,$ and semisimple maps\/ $f_+=f\ci
i_X\vert_{\smash{\pd_+^fX}}:\pd_+^fX\ra Y$ and\/ $f_-:\pd_-^fX\ra\pd
Y,$ such that the following commutes in
$\Manc\!:$\G[dfXa]{$\pd^f_\pm X$}{sets of decomposition $\pd
X=\pd^f_+X\amalg\pd^f_-X$ of boundary $\pd X$ induced by $f:X\ra Y$
in $\Manc$}
\e
\begin{split}
\xymatrix@C=120pt@R=10pt{ \pd_-^fX \ar[r]_{f_-} \ar[d]^(0.6){i_X
\vert_{\pd_-^fX}} & \pd Y \ar[d]_{i_Y} \\ X \ar[r]^{f} & Y.}
\end{split}
\label{ds5eq1}
\e
If\/ $f$ is also flat, then\/ \eq{ds5eq1} is a Cartesian square, so
that\/ $\pd_-^fX\cong X\t_Y\pd Y$. If\/ $f$ is simple then
$\pd^f_+X=\es$ and\/ $\pd^f_-X=\pd X$. If\/ $f$ is simple, flat, a
submersion, or an s-submersion, then $f_\pm$ are also simple,
\ldots, s-submersions, respectively.
\label{ds5prop1}
\end{prop}

In fact we define $\pd^f_-X=s_f(S_f)$, so that $s_f:S_f\ra\pd^f_-X$
is a bijection since $s_f$ is injective as $f$ is semisimple, and
then $f_-=u_f\ci s_f^{-1}$, using the notation of Definition
\ref{ds5def2}. If $f:X\ra Y$ is simple then $f_-:\pd X\ra\pd Y$ is
also simple, so $f_{-^k}:\pd^kX\ra\pd^kY$ is simple for
$k=1,2,\ldots.$ If $f$ is also flat then $f_{-^k}$ is flat and
$\pd^kX\cong X\t_Y\pd^kY$. A smooth map $f:X\ra Y$ is flat if and
only if $f(X^\ci)\subseteq Y^\ci$, or equivalently, if $f:X\ra Y$
and $i_Y:\pd Y\ra Y$ are transverse.

(S-)submersions are locally modelled on projections~$\pi_X:X\t Y\ra
X$:

\begin{prop}{\bf(a)} Let\/ $X,Y$ be manifolds with corners. Then the
projection $\pi_X:X\t Y\ra X$ is a submersion, and an s-submersion
if\/~$\pd Y=\es$.

\smallskip

\noindent{\bf(b)} Let\/ $f:X\ra Y$ be a submersion of manifolds with
corners, and\/ $x\in X$ with\/ $f(x)=y\in Y$. Then there exist open
neighbourhoods $V$ of\/ $x$ in $X$ and\/ $W$ of\/ $y$ in $Y$ with\/
$f(V)=W,$ a manifold with corners $Z,$ and a diffeomorphism $V\cong
W\t Z$ which identifies $f\vert_V:V\ra W$ with\/ $\pi_W:W\t Z\ra W$.
If\/ $f$ is an s-submersion then\/~$\pd Z=\es$.
\label{ds5prop2}
\end{prop}

S-immersions and sf-immersions are also locally modelled on
products:

\begin{prop}{\bf(a)} Let\/ $X$ be a manifold with corners and\/
$0\le k\le n$. Then $\id_X\t 0:X\ra X\t\bigl([0,\iy)^k\t\R^{n-k}
\bigr)$ mapping $x\mapsto(x,0)$ is an s-embedding, and an
sf-embedding if\/~$k=0$.
\smallskip

\noindent{\bf(b)} Let\/ $f:X\ra Y$ be an s-immersion of manifolds
with corners, and\/ $x\in X$ with\/ $f(x)=y\in Y$. Then there exist
open neighbourhoods $V$ of\/ $x$ in $X$ and\/ $W$ of\/ $y$ in $Y$
with\/ $f(V)\subseteq W,$ an open neighbourhood\/ $Z$ of\/ $0$ in
$[0,\iy)^k\t\R^{n-k},$ and a diffeomorphism $W\cong V\t Z$ which
identifies $f\vert_V:V\ra W$ with\/ $\id_V\t 0:V\ra V\t Z$. If\/ $f$
is an sf-immersion then\/~$k=0$.
\label{ds5prop3}
\end{prop}

Example \ref{ds5ex2}(ii) shows general immersions are not modelled
on products.\I{manifold with corners!simple map|)}\I{manifold with
corners!semisimple map|)}\I{manifold with corners!flat map|)}
\I{manifold with corners!submersion|)}\I{manifold with
corners!s-submersion|)}\I{manifold with corners!immersion|)}
\I{manifold with corners!s-immersion|)}\I{manifold with
corners!sf-immersion|)}\I{manifold with corners!embedding|)}
\I{manifold with corners!s-embedding|)}\I{manifold with
corners!sf-embedding|)}

\subsection{Corners and the corner functors}
\label{ds53}
\I{manifold with corners!k-corners C_k(X)@$k$-corners $C_k(X)$}

As in \cite[\S 2]{Joyc3}, \cite[\S 5.5]{Joyc6}, we define the
$k$-{\it corners\/} $C_k(X)$ of a manifold with corners~$X$.

\begin{dfn} Let $X$ be an $n$-manifold with corners.
Applying $\pd$ repeatedly gives manifolds with corners $\pd
X,\pd^2X,\ldots.$ There is a natural identification
\e
\begin{split}
\pd^kX\cong\bigl\{(x,\be_1,\ldots,\be_k):\text{$x\in X,$
$\be_1,\ldots,\be_k$ are distinct}&\\
\text{local boundary components for $X$ at $x$}&\bigr\}.
\end{split}
\label{ds5eq2}
\e
Using \eq{ds5eq2}, we see that the symmetric group $S_k$ of
permutations of $\{1,\ldots,k\}$ has a natural, free action on
$\pd^kX$ by diffeomorphisms, given by
\begin{equation*}
\si:(x,\be_1,\ldots,\be_k)\longmapsto
(x,\be_{\si(1)},\ldots,\be_{\si(k)}).
\end{equation*}

Define the $k$-{\it corners\/} of $X$, as a set, to be
\begin{align*}
C_k(X)=\bigl\{(x,\{\be_1,\ldots,\be_k\}):\text{$x\in X,$
$\be_1,\ldots,\be_k$ are distinct}&\\
\text{local boundary components for $X$ at $x$}&\bigr\}.
\end{align*}
Then $C_k(X)$ is naturally a manifold with corners of dimension
$n-k$, with $C_k(X)\cong\pd^kX/S_k$. The interior $C_k(X)^\ci$ is
naturally diffeomorphic to $S^k(X)$. We have natural diffeomorphisms
$C_0(X)\cong X$ and~$C_1(X)\cong\pd X$.
\label{ds5def4}
\end{dfn}

A surprising fact about manifolds with corners $X$ is that the
disjoint union $C(X):=\coprod_{k=0}^{\dim X}C_k(X)$ has strong
functorial properties. Since $C(X)$ is not a manifold with corners,
it is helpful to enlarge our category~$\Manc$:

\begin{dfn} Write $\cManc$ for the category whose objects are
disjoint unions $\coprod_{m=0}^\iy X_m$, where $X_m$ is a manifold
with corners of dimension $m$, and whose morphisms are continuous
maps $f:\coprod_{m=0}^\iy X_m\ra\coprod_{n=0}^\iy Y_n$, such that
$f\vert_{X_m\cap f^{-1}(Y_n)}:\bigl(X_m\cap f^{-1}(Y_n)\bigr)\ra
Y_n$ is a smooth map of manifolds with corners for all~$m,n\ge
0$.\I{manifold with corners!corner
functors|(}\G[Manc']{$\cManc$}{category of disjoint unions of
manifolds with corners of different dimensions}
\label{ds5def5}
\end{dfn}

\begin{dfn} Define {\it corner functors\/} $C,\hat C:\Manc\ra\cManc$
by $C(X)=\hat C(X)=\coprod_{k=0}^{\dim X}C_k(X)$ on objects, and on
morphisms $f:X\ra Y$ in $\Manc$,\G[CCa]{$C,\hat
C:\Manc\ra\cManc$}{`corner functors' for manifolds with corners}
\begin{align*}
&C(f):\bigl(x,\{\ti\be_1,\ldots,\ti\be_i\}\bigr)\longmapsto
\bigl(y,\{\be_1,\ldots,\be_j\}\bigr),\quad\text{where $y=f(x)$,}\\
&\{\be_1,\ldots,\be_j\}\!=\!\bigl\{\be:\bigl((x,\ti\be_l),(y,\be)\bigr)
\in S_f,\; \text{some $l=1,\ldots,i$}\bigr\},
\\
&\hat C(f):\bigl(x,\{\ti\be_1,\ldots,\ti\be_i\}\bigr)\longmapsto
\bigl(y,\{\be_1,\ldots,\be_j\}\bigr),\quad\text{where $y=f(x)$,}\\
&\{\be_1,\ldots,\be_j\}\!=\!\bigl\{\be:\bigl((x,\ti\be_l),(y,\be)\bigr)
\in S_f,\; l=1,\ldots,i\bigr\}\cup\bigl\{\be:\bigl(x,(y,\be)\bigr)
\!\in\! T_f\bigr\}.
\end{align*}

Write $C_j^{f,k}(X)=C_j(X)\cap C(f)^{-1}(C_k(Y))$ and
$C_j^k(f)=C(f)\vert_{C_j^{f,k}(X)}:C_j^{f,k}(X)\ab\ra C_k(Y)$ for
all $j,k$, and similarly for $\hat C{}_j^{f,k}(X),\hat C{}_j^k(f)$.
Then $C_j^k(f)$ and $\hat C{}_j^k(f)$ are smooth maps of manifolds
with corners. Note that $C_0^{f,0}(X)=C_0(X)\cong X$ and
$C_0(Y)\cong Y$, and these isomorphisms identify $C_0^0(f):C_0(X)\ra
C_0(Y)$ with~$f:X\ra Y$.
\label{ds5def6}
\end{dfn}

It turns out that $C,\hat C$ are both functors $\Manc\ra\cManc$.
Furthermore:
\begin{itemize}
\setlength{\itemsep}{0pt}
\setlength{\parsep}{0pt}
\item[(i)] For each $X\in\Manc$ we have a natural
diffeomorphism $C(\pd X)\cong \pd C(X)$ identifying
$C(i_X):C(\pd X)\ra C(X)$ with $i_{C(X)}:\pd C(X)\ra C(X)$
\item[(ii)] For all $X,Y$ in $\Manc$ we have a natural
diffeomorphism $C(X\t Y)\cong C(X)\t C(Y)$. These
diffeomorphisms commute with product morphisms and direct
product morphisms in $\Manc,\cManc$.
\item[(iii)] If $g:X\ra Z$ and $h:Y\ra Z$ are {\it strongly
transverse\/} maps in $\Manc$ then $C$ maps the fibre product
$X\t_{g,Z,h}Y$ in $\Manc$ to the fibre product
$C(X)\t_{C(g),C(Z),C(h)}C(Y)$ in $\cManc$.\I{manifold with
corners!strongly transverse maps}
\item[(iv)] If $f:X\ra Y$ is semisimple, then $C(f)$ maps
$C_k(X)\ra \coprod_{l=0}^kC_l(Y)$ for all $k\ge 0$. The natural
diffeomorphisms $C_1(X)\cong\pd X,$ $C_0(Y)\cong Y$ and
$C_1(Y)\cong\pd Y$ identify $C_1^{f,0}(X)\cong\pd^f_+X,$
$C_1^0(f)\cong f_+,$ $C_1^{f,1}(X)\cong\pd^f_-X$ and
$C_1^1(f)\cong f_-$. If $f$ is simple then $C(f)$ maps
$C_k(X)\ra C_k(Y)$ for all~$k\ge 0$.
\end{itemize}
The analogues hold for $\hat C$, except for (iv) and the last part
of~(i).\I{manifold with corners!corner functors|)}

\subsection{(Strong) transversality and fibre products}
\label{ds54}
\I{manifold with corners!strongly transverse maps|(}\I{manifold with
corners!transverse fibre products|(}

In \cite[\S 6]{Joyc3}, \cite[\S 5.6]{Joyc6} we discuss conditions
for fibre products to exist in~$\Manc$.

\begin{dfn} Let $g:X\ra Z$, $h:Y\ra Z$ be smooth maps of manifolds
with corners. We call $g,h$ {\it transverse\/} if whenever $x\in
S^j(X)\subseteq X$, $y\in S^k(Y)\subseteq Y$ and $z\in
S^l(Z)\subseteq Z$ with $g(x)=h(y)=z$, then $T_zZ=\d
g\vert_x(T_xX)+\d h\vert_y(T_yY)$ and~$T_z(S^l(Z))=\d
g\vert_x(T_x(S^j(X)))+\d h\vert_y(T_y(S^k(Y)))$.

We call $g,h$ {\it strongly transverse\/} if they are transverse,
and whenever there are points in $C_j(X),C_k(Y),C_l(Z)$ with
\begin{equation*}
C(g)(x,\{\be_1,\ldots,\be_j\})=C(h)(y,\{\ti\be_1,\ldots,
\ti\be_k\})=(z,\{\dot\be_1,\ldots,\dot\be_l\})
\end{equation*}
we have either $j+k>l$ or~$j=k=l=0$.
\label{ds5def7}
\end{dfn}

If one of $g,h$ is a submersion then $g,h$ are strongly transverse.
It is well known that transverse fibre products of manifolds without
boundary exist. Here is the (more difficult to prove) analogue for
manifolds with corners.

\begin{thm} Let\/ $g:X\ra Z,$ $h:Y\ra Z$ be transverse smooth
maps of manifolds with corners. Then a fibre product\/
$W\!=\!X\t_{g,Z,h}Y$ exists in~$\Manc$.
\label{ds5thm1}
\end{thm}

As a topological space, the fibre product in Theorem \ref{ds5thm1}
is just the topological fibre product $W=\bigl\{(x,y)\in X\t
Y:g(x)=h(y)\bigr\}$. In general, the boundary $\pd W$ is difficult
to describe explicitly: it is the quotient of a subset of $(\pd
X\t_ZY)\amalg(X\t_Z\pd Y)$ by an equivalence relation. Here are some
special cases in which we can give an explicit formula for~$\pd
W$.\I{manifold with corners!transverse fibre products!boundaries of}

\begin{prop} Let\/ $g:X\ra Z,$ $h:Y\ra Z$ be transverse smooth maps
in $\Manc,$ so that\/ $X\t_{g,Z,h}Y$ exists by Theorem\/
{\rm\ref{ds5thm1}}. Then:
\begin{itemize}
\setlength{\itemsep}{0pt}
\setlength{\parsep}{0pt}
\item[{\bf(a)}] If\/ $\pd Z=\es$ then
\e
\pd\bigl(X\t_{g,Z,h}Y\bigr)\cong \bigl(\pd X\t_{g\ci i_X,Z,h}Y\bigr)
\amalg \bigl(X\t_{g,Z,h\ci i_Y}\pd Y\bigr).
\label{ds5eq3}
\e
\item[{\bf(b)}] If\/ $g$ is semisimple then\I{manifold with
corners!semisimple map}
\e
\pd\bigl(X\t_{g,Z,h}Y\bigr)\cong \bigl(\pd_+^gX
\t_{g_+,Z,h}Y\bigr) \amalg \bigl(X\t_{g,Z,h\ci i_Y}\pd Y\bigr).
\label{ds5eq4}
\e
\item[{\bf(c)}] If both\/ $g,h$ are semisimple then
\end{itemize}
\e
\begin{split}
\pd\bigl(&X\t_{g,Z,h}Y\bigr)\cong \\
&\bigl(\pd_+^gX \t_{g_+,Z,h}Y\bigr)\amalg \bigl(X \t_{g,Z,h_+}
\pd_+^hY\bigr)\amalg\bigl(\pd_-^gX\t_{g_-,\pd
Z,h_-}\pd_-^hY\bigr).
\end{split}
\label{ds5eq5}
\e
Here all fibre products in \eq{ds5eq3}--\eq{ds5eq5} are transverse,
and so exist.
\label{ds5prop4}
\end{prop}

For {\it strongly\/} transverse smooth maps, fibre products commute
with the corner functors $C,\hat C:\Manc\ra\cManc$. Since
$C_1(W)\cong\pd W$, equation \eq{ds5eq6} with $i=1$ gives another
explicit description of $\pd W$ in this case.\I{manifold with
corners!corner functors}

\begin{thm} Let\/ $g:X\ra Z,$ $h:Y\ra Z$ be strongly transverse
smooth maps of manifolds with corners, and write\/ $W$ for the fibre
product\/ $X\t_{g,Z,h}Y$ given by Theorem\/ {\rm\ref{ds5thm1}}. Then
there is a canonical diffeomorphism
\e
C_i(W)\cong \coprod_{j,k,l\ge 0:i=j+k-l}
C_j^{g,l}(X)\t_{C_j^l(g),C_l(Z),C_k^l(h)}C_k^{h,l}(Y)
\label{ds5eq6}
\e
for all\/ $i\ge 0,$ where the fibre products are all transverse and
so exist. Hence\I{manifold with corners!strongly transverse
maps|)}\I{manifold with corners!transverse fibre products|)}
\begin{equation*}
C(W)\cong C(X)\t_{C(g),C(Z),C(h)}C(Y)\quad\text{in\/ $\cManc$.}
\end{equation*}
\label{ds5thm2}
\end{thm}

\subsection{Orientations on manifolds with corners}
\label{ds55}
\I{manifold with corners!orientations|(}

In \cite[\S 7]{Joyc3}, \cite[\S 5.8]{Joyc6} we discuss orientations
on manifolds with corners.

\begin{dfn} Let $X$ be an $n$-manifold with corners. An {\it
orientation\/} $\om$ on $X$ is an orientation on the fibres of the
real line bundle $\La^nT^*X$ over $X$. That is, $\om$ is an
equivalence class $[\tau]$ of isomorphisms $\tau:O_X\ra\La^nT^*X$,
where $O_X=\R\t X\ra X$ is the trivial line bundle on $X$, and
$\tau,\tau'$ are equivalent if $\tau'=\tau\cdot c$ for some smooth
$c:X\ra(0,\iy)$.

If $\om=[\tau]$ is an orientation, we write $-\om$ for the {\it
opposite orientation\/}~$[-\tau]$.

We call the pair $(X,\om)$ an {\it oriented manifold}. Usually we
suppress the orientation $\om$, and just refer to $X$ as an oriented
manifold. When $X$ is an oriented manifold, we write $-X$ for $X$
with the opposite orientation.
\label{ds5def8}
\end{dfn}

If $X,Y,Z$ are oriented manifolds with corners, then we can define
orientations on boundaries $\pd X$, products $X\t Y$, and transverse
fibre products $X\t_ZY$. To do this requires a choice of {\it
orientation convention}.\I{orientation convention} Our orientation
conventions are given in \cite[\S 5.8]{Joyc6}. Having fixed an
orientation convention, natural isomorphisms of manifolds with
corners such as $X\t_ZY\cong Y\t_ZX$ lift to isomorphisms of
oriented manifolds of corners, modified by signs depending on the
dimensions. For example, if $g:X\ra Z$ and $h:Y\ra Z$ are transverse
maps of oriented manifolds with corners then
\begin{equation*}
X\t_{g,Z,h}Y\cong (-1)^{(\dim X-\dim Z)(\dim Y-\dim Z)}Y\t_{h,Z,g}X,
\end{equation*}
and with orientations equations \eq{ds5eq3}--\eq{ds5eq5}
become\I{manifold with corners!orientations|)}
\begin{align*}
\pd\bigl(X\!\t_{g,Z,h}\!Y\bigr)\!&\cong \!\bigl(\pd X\!\t_{g\ci
i_X,Z,h}\!Y\bigr)\!\amalg\! (-\!1)^{\dim X+\dim
Z}\!\bigl(X\!\t_{g,Z,h\ci i_Y}\!\pd Y\bigr),\\
\pd\bigl(X\!\t_{g,Z,h}\!Y\bigr)\!&\cong\!\bigl(\pd_+^gX\!
\t_{g_+,Z,h}\!Y\bigr)\! \amalg\! (-\!1)^{\dim X+\dim
Z}\!\bigl(X\!\t_{g,Z,h\ci i_Y}\!\pd Y\bigr),\\
\pd\bigl(X\t_{g,Z,h}Y\bigr)\!&\cong\!\bigl(\pd_+^gX
\t_{g_+,Z,h}Y\bigr)\!\amalg\! (-\!1)^{\dim X+\dim Z}\!\bigl(X
\!\t_{h,Z,h_+}\!\pd_+^hY\bigr)\\
&\qquad \amalg\bigl(\pd_-^gX\t_{g_-,\pd Z,h_-}\pd_-^hY\bigr).
\end{align*}

\subsection{Fixed point loci in manifolds with corners}
\label{ds56}
\I{manifold with corners!fixed point loci|(}

In \cite[\S 5.5]{Joyc6} we study the fixed point locus $X^\Ga$ of a
group $\Ga$ acting on a manifold with corners $X$. These are related
to orbifold strata $\oX^\Ga$ of orbifolds with corners $\oX$, which
we will discuss in \S\ref{ds125}. Here is our main
result.\G[XGab]{$X^\Ga$}{fixed subset of a group $\Ga$ acting on a
manifold with corners $X$}

\begin{prop} Suppose $X$ is a manifold with corners, $\Ga$ a
finite group, and\/ $r:\Ga\ra\Aut(X)$ an action of\/ $\Ga$ on $X$ by
diffeomorphisms. Applying the corner functor $C$ of\/
{\rm\S\ref{ds53}} gives an action $C(r):\Ga\ra\Aut(C(X))$ of\/ $\Ga$
on $C(X)$ by diffeomorphisms. Write $X^\Ga,C(X)^\Ga$ for the subsets
of\/ $X,C(X)$ fixed by $\Ga,$ and\/ $j_{X,\Ga}:X^\Ga\ra X$ for the
inclusion. Then:\G[jXGab]{$j_{X,\Ga}:X^\Ga\hookra X$}{inclusion of
$\Ga$-fixed subset $X^\Ga$ in a manifold with corners $X$}
\begin{itemize}
\setlength{\itemsep}{0pt}
\setlength{\parsep}{0pt}
\item[{\bf(a)}] $X^\Ga$ has the structure of an object in
$\cManc$ {\rm(}a disjoint union of manifolds with corners of
different dimensions, as in {\rm\S\ref{ds53})} in a unique way,
such that\/ $j_{X,\Ga}:X^\Ga\ra X$ is an embedding. This\/
$j_{X,\Ga}$ is flat, but need not be (semi)simple.
\item[{\bf(b)}] By {\bf(a)} we have a smooth map
$C(j_{X,\Ga}):C(X^\Ga)\ra C(X)$. This $C(j_{X,\Ga})$ is a
diffeomorphism $C(X^\Ga)\ra C(X)^\Ga$. As\/ $j_{X,\Ga}$ need not
be simple, $C(j_{X,\Ga})$ need not map $C_k(X^\Ga)\ra C_k(X)$
for $k>0$.
\item[{\bf(c)}] By {\bf(b)}{\rm,} $C(j_{X,\Ga})$ identifies
$C_1(X^\Ga)\cong\pd(X^\Ga)$ with a subset of\/
$C(X)^\Ga\subseteq C(X)$. This gives the following description
of\/~$\pd(X^\Ga)\!:$
\begin{align*}
\pd(X^\Ga)\cong\bigl\{(x,\,&\{\be_1,\ldots,\be_k\})\in
C_k(X):\text{$x\in X^\Ga,$ $k\ge 1,$ $\be_1,\ldots,\be_k$}\\
&\text{are distinct local boundary components for $X$ at\/ $x,$}\\
&\text{and\/ $\Ga$ acts transitively on
$\{\be_1,\ldots,\be_k\}$}\bigr\}.
\end{align*}
\item[{\bf(d)}] Now suppose $Y$ is a manifold with corners with
an action of\/ $\Ga,$ and\/ $f:X\ra Y$ is a $\Ga$-equivariant
smooth map. Then $X^\Ga,Y^\Ga$ are objects in $\cManc$ by
{\bf(a)\rm,} and\/ $f^\Ga:=f\vert_{X^\Ga}:X^\Ga\ra Y^\Ga$ is a
morphism in\/~$\cManc$.
\end{itemize}
\label{ds5prop5}
\end{prop}

\begin{ex} Let $\Ga=\{1,\si\}$ with $\si^2=1$, so that
$\Ga\cong\Z_2$, and let $\Ga$ act on $X=[0,\iy)^2$ by
$\si:(x_1,x_2)\mapsto (x_2,x_1)$. Then
$X^\Ga=\bigl\{(x,x):x\in[0,\iy)\bigr\}\cong [0,\iy)$, a manifold
with corners, and the inclusion $j_{X,\Ga}:X^\Ga\ra X$ is
$j_{X,\Ga}:[0,\iy)\ra [0,\iy)^2$, $j_{X,\Ga}:x\mapsto(x,x)$, a
smooth, flat embedding, which is not semisimple. We have $\pd
X=\pd\bigl([0,\iy)^2\bigr)\cong [0,\iy)\amalg[0,\iy)$, where $\Ga$
acts freely on $\pd X$ by exchanging the two copies of $[0,\iy)$.
Hence $(\pd X)^\Ga=\es$, but $\pd(X^\Ga)$ is a point $*$, so in this
case $(\pd X)^\Ga\not\cong\pd(X^\Ga)$. Also
$C_2(X)=\bigl\{\bigl(0,\{\{x_1=0\}, \{x_2=0\}\}\bigr)\bigr\}$ is a
single point, which is $\Ga$-invariant, and
$C(j_{X,\Ga}):C(X^\Ga)\ra C(X)^\Ga$ identifies $(0,\{\{x=0\}\})\in
C_1(X^\Ga)\cong\pd X$ with this point in~$C_2(X)^\Ga$.
\label{ds5ex3}
\end{ex}

If a finite group $\Ga$ acts on a manifold with corners $X$ then as
in Proposition \ref{ds5prop5}(b) we have $C(X)^\Ga\cong C(X^\Ga)$,
but as in Example \ref{ds5ex3} in general we do not have $(\pd
X)^\Ga\cong\pd (X^\Ga)$, but only $(\pd X)^\Ga\subseteq\pd (X^\Ga)$.
Thus for fixed point loci, corners have more functorial behaviour
than boundaries.\I{manifold with boundary|)}\I{manifold with
corners|)}\I{manifold with corners!fixed point loci|)}

\section{D-spaces with corners}
\label{ds6}
\I{d-space with corners|(}

The goal of \cite[Chap.s 6 \& 7]{Joyc6} is to construct a
well-behaved 2-category $\dManc$ of {\it d-manifolds with corners},
a derived version of $\Manc$. It is tempting to define $\dManc$ as a
2-subcategory of d-spaces $\dSpa$, but this turns out not to be a
good idea. For example, the natural functor
$F_\Manc^\dSpa:\Manc\ra\dSpa$ is not full, as 1-morphisms $\bs
f:F_\Manc^\dSpa(X)\ra F_\Manc^\dSpa(Y)$ correspond to weakly smooth
rather than smooth maps $f:X\ra Y$, in the notation of~\S\ref{ds51}.

Therefore we begin in \cite[Chap.~6]{Joyc6} by defining a 2-category
$\dSpac$ of {\it d-spaces with corners}, and then define $\dManc$ in
\cite[Chap.~7]{Joyc6} as a 2-subcategory of $\dSpac$. Many
properties of manifolds with corners in \S\ref{ds5} work for
d-spaces with corners, e.g.\ boundaries $\pd X$, simple, semisimple
and flat maps $f:X\ra Y$, decompositions $\pd
X=\pd^f_+X\amalg\pd^f_-X$ and semisimple maps $f_+:\pd_+^fX\ra Y$
and $f_-:\pd_-^fX\ra\pd Y$ when $f$ is semisimple, and the corner
functors~$C,\hat C$.

\subsection{Outline of the definition of the 2-category $\dSpac$}
\label{ds61}
\I{d-space with corners!definition|(}\I{2-category|(}

The definition of the 2-category of d-spaces with corners $\dSpac$
in \cite[\S 6.1]{Joyc6} is long and complicated. So here we just
sketch the main ideas.\G[dSpac]{$\dSpac$}{2-category of d-spaces
with corners}

Let $X$ be a manifold with corners. Then it has a boundary $\pd X$
with a proper smooth map $i_X:\pd X\ra X$. On $\pd X$ we have an
exact sequence
\e
\smash{\xymatrix@C=20pt{ 0 \ar[r] & \cN_X \ar[rr] && i_X^*(T^*X)
\ar[rr]^{(\d i_X)^*} && T^*(\pd X) \ar[r] & 0,}}
\label{ds6eq1}
\e
where $\cN_X$ is the conormal bundle of $\pd X$ in $X$. The line
bundle $\cN_X$ has a natural orientation $\om_X$ induced by
outward-pointing normal vectors to $\pd X$ in~$X$.

Thus, for each manifold with corners $X$ we have a quadruple $(X,\pd
X,i_X,\om_X)$. D-spaces with corners are based on this idea. A {\it
d-space with corners\/} $\rX$ is a quadruple $\rX=(\bX,\bpX,\bs
i_\rX,\om_\rX)$ where $\bX,\bpX$ are d-spaces, and $\bs
i_\rX:\bpX\ra\bX$\G[iXb]{$\bs i_\rX:\pd\rX\ra\rX$}{inclusion of
boundary $\pd\rX$ into a d-space with corners $\rX$} is a proper
1-morphism, and we have an exact sequence in~$\qcoh(\upX)$:
\e
\smash{\xymatrix{ 0 \ar[r] & \cN_\rX \ar[rr]^{\nu_\rX} &&
\ui_\rX^*(\FX) \ar[rr]^{i_\rX^2} && \cF_{\pd X} \ar[r] & 0,}}
\label{ds6eq2}
\e
with $\cN_\rX$ a line bundle,\G[NXa]{$\cN_\rX$}{conormal line bundle
of $\pd\rX$ in $\rX$ for a d-space with corners $\rX$} and
$\om_\rX$\G[omXa]{$\om_\rX$}{orientation on line bundle $\cN_\rX$
for a d-space with corners $\rX$} is an orientation on $\cN_\rX$.
These $\bX,\bpX,\bs i_\rX,\om_\rX$ must satisfy some complicated
conditions in \cite[\S 6.1]{Joyc6}, that we will not give. They
require $\bpX$ to be locally equivalent to a fibre product
$\bX\t_{\bs{[0,\iy)}}\bs{*}$ in~$\dSpa$.

If $\rX=(\bX,\bpX,\bs i_\rX,\om_\rX)$ and $\rY=(\bY,\bpY,\bs
i_\rY,\om_\rY)$ are d-spaces with corners, a 1-{\it morphism\/} $\bs
f:\rX\ra\rY$ in $\dSpac$ is a 1-morphism $\bs f:\bX\ra\bY$ in
$\dSpa$ satisfying extra conditions over $\bpX,\bpY$, which are
analogous to the extra conditions for a weakly smooth map of
manifolds with corners $f:X\ra Y$ to be smooth in
Definition~\ref{ds5def2}.

If $\bs f:\rX\ra\rY$ is a 1-morphism in $\dSpac$, we can form the
$C^\iy$-scheme fibre products $\upX\t_{\uf\ci\ui_\rX,
\uY,\ui_\rY}\upY$ and $\uX\t_{\uf,\uY,\ui_\rY}\upY$. As for
$S_f,T_f$ in Definition \ref{ds5def2}, we can define open and closed
$C^\iy$-subschemes $\uS_{\bs f}\subseteq\upX\t_\uY\upY$ and
$\uT_{\bs f}\subseteq\uX\t_\uY\upY$, and define $C^\iy$-scheme
morphisms $\us_{\bs f}:\uS_{\bs f}\ra\upX$, $\ut_{\bs f}:\uT_{\bs
f}\ra\uX$, $\uu_{\bs f}:\uS_{\bs f}\ra\upY$ and $\uv_{\bs
f}:\uT_{\bs f}\ra\upY$ to be the projections from the fibre
products. Then $\us_{\bs f},\ut_{\bs f}$ are
\'etale.\G[Sfb]{$\uS_{\bs f}\subseteq\upX\t_\uY\upY$}{$C^\iy$-scheme
associated to 1-morphism $\bs f:\rX\ra\rY$ in
$\dSpac$}\G[Tfb]{$\uT_{\bs f}\subseteq\uX\t_\uY\upY$}{$C^\iy$-scheme
associated to 1-morphism $\bs f:\rX\ra\rY$ in $\dSpac$}

If $\bs f,\bs g:\rX\ra\rY$ are 1-morphisms in $\dSpac$, a 2-{\it
morphism\/} $\eta:\bs f\Ra\bs g$ in $\dSpac$ is a 2-morphism
$\eta:\bs f\Ra\bs g$ in $\dSpa$ such that $\uS_{\bs f}=\uS_{\bs g}$,
$\uT_{\bs f}=\uT_{\bs g}$ and extra vanishing conditions hold on
$\eta$ over $\uS_{\bs f},\uT_{\bs f}$. Identity 1- and 2-morphisms
in $\dSpac$, and the compositions of 1- and 2-morphisms in $\dSpac$,
are all given by identities and compositions
in~$\dSpa$.\I{2-category|)}

A d-space with corners $\rX=(\bX,\bpX,\bs i_\rX,\om_\rX)$ is called
a {\it d-space with boundary\/}\I{d-space with boundary} if $\bs
i_\rX:\bpX\ra\bX$ is injective, and a {\it d-space without
boundary\/} if $\bpX=\bs\es$. We write
$\dSpab$\G[dSpab]{$\dSpab$}{2-category of d-spaces with boundary}
for the full 2-subcategory of d-spaces with boundary, and
$\bdSpa$\G[dSpa']{$\bdSpa$}{2-subcategory of d-spaces with corners
equivalent to d-spaces} for the full 2-subcategory of d-spaces
without boundary, in $\dSpac$. There is an isomorphism of
2-categories $F_\dSpa^\dSpac:\dSpa\ra\bdSpa$ mapping
$\bX\mapsto\rX=(\bX,\bs\es,\bs\es,\bs\es)$ on objects, $\bs
f\mapsto\bs f$ on 1-morphisms and $\eta\mapsto\eta$ on 2-morphisms.
So we can consider d-spaces to be examples of d-spaces with
corners.\I{d-space with corners!definition|)}\I{d-space!as d-space
with corners}\I{d-space with corners!include d-spaces}

\begin{rem} If $X$ is a manifold with corners then the orientation
$\om_X$ on $\cN_X$ is determined uniquely by $X,\pd X,i_X$. But
there are examples of d-spaces with corners $\rX=(\bX,\bpX,\bs
i_\rX,\om_\rX)$ in which $\om_\rX$ is not determined by
$\bX,\bpX,\bs i_\rX$, and really is extra data. We include $\om_\rX$
in the definition so that orientations of d-manifolds with corners
behave well in relation to boundaries. If we had omitted $\om_\rX$
from the definition, then there would exist examples of oriented
d-manifolds with corners $\rX$ such that $\pd\rX$ is not orientable.
\label{ds6rem1}
\end{rem}

For each d-space with corners $\rX=(\bX,\bpX,\bs i_\rX,\om_\rX)$, in
\cite[\S 6.2]{Joyc6} we define a d-space with corners $\pd\rX=(\bpX,
\bs{\pd^2X},\bs i_{\pd\rX},\om_{\pd\rX})$ called the {\it
boundary\/}\I{d-space with corners!boundary}\I{boundary!of a d-space
with corners}\G[dXb]{$\pd\rX$}{boundary of a d-space with corners
$\rX$} of $\rX$, and show that $\bs i_\rX:\pd\rX\ra\rX$ is a
1-morphism in $\dSpac$. Motivated by \eq{ds5eq2} when $k=2$, the
d-space $\bs{\pd^2X}$ in $\pd\rX$ is given by
\e
\bs{\pd^2X}\simeq \bigl(\bpX\t_{\bs i_\rX,\bX,\bs i_\rX}\bpX\bigr)
\sm\bs\De_\bpX(\bpX),
\label{ds6eq3}
\e
where $\De_\bpX:\bpX\ra\bpX\t_\bX\bpX$ is the diagonal 1-morphism.
The 1-morphism $\bs i_{\pd\rX}:\bs{\pd^2X}\ra\bpX$ is projection to
the first factor in the fibre product. There is a natural
isomorphism $\cN_{\pd\rX}\cong\ui_\rX^*(\cN_\rX)$, and the
orientation $\om_{\pd\rX}$ on $\cN_{\pd\rX}$ is defined to
correspond to the orientation $\ui_\rX^*(\om_\rX)$
on~$\ui_\rX^*(\cN_\rX)$.

\subsection{Simple, semisimple and flat 1-morphisms}
\label{ds62}
\I{d-space with corners!simple 1-morphism|(}\I{d-space with
corners!semisimple 1-morphism|(}\I{d-space with corners!flat
1-morphism|(}

In \cite[\S 6.3]{Joyc6} we generalize the material on simple,
semisimple, and flat maps of manifolds with corners in \S\ref{ds52}
to d-spaces with corners. Here are the analogues of Definition
\ref{ds5def3}(a)--(c) and Proposition \ref{ds5prop1}.

\begin{dfn} Let $\bs f:\rX\ra\rY$ be a 1-morphism of d-spaces
with corners.
\begin{itemize}
\setlength{\itemsep}{0pt}
\setlength{\parsep}{0pt}
\item[(a)] We call $\bs f$ {\it simple\/} if $\us_{\bs
f}:\uS_{\bs f}\ra\upX$ is bijective.
\item[(b)] We call $\bs f$ {\it semisimple\/} if $\us_{\bs
f}:\uS_{\bs f}\ra\upX$ is injective.
\item[(c)] We call $\bs f$ {\it flat\/} if $\uT_{\bs f}=\es$.
\end{itemize}
\label{ds6def1}
\end{dfn}

\begin{thm} Let\/ $\bs f:\rX\ra\rY$ be a semisimple $1$-morphism
of d-spaces with corners. Then there exists a natural decomposition
$\pd\rX=\pd_+^{\bs f}\rX\amalg\pd_-^{\bs f}\rX$ with\/ $\pd_\pm^{\bs
f}\rX$ open and closed in $\pd\rX,$ such that:\G[dfXb]{$\pd^{\bs
f}_\pm\rX$}{sets of decomposition $\pd\rX=\pd_+^{\bs
f}\rX\amalg\pd_-^{\bs f}\rX$ of boundary $\pd\rX$ induced by
1-morphism $\bs f:\rX\ra\rY$ in $\dSpac$}
\begin{itemize}
\setlength{\itemsep}{0pt}
\setlength{\parsep}{0pt}
\item[{\bf(a)}] Define $\bs f_+=\bs f\ci \bs i_\rX\vert_{\pd^{\bs
f}_+\rX}:\pd_+^{\bs f}\rX\ra\rY$. Then $\bs f_+$ is semisimple.
If\/ $\bs f$ is flat then $\bs f_+$ is also flat.
\item[{\bf(b)}] There exists a unique, semisimple\/ $1$-morphism
$\bs f_-:\pd_-^{\bs f}\rX\ra\pd\rY$ with\/ $\bs f\ci\bs i_\rX
\vert_{\pd_-^{\bs f}\rX}=\bs i_\rY\ci\bs f_-$. If\/ $\bs f$ is
simple then $\pd_+^{\bs f}\rX=\bs\es,$ $\pd_-^{\bs
f}\rX=\pd\rX,$ and\/ $\bs f_-:\pd\rX\ra\pd\rY$ is also simple.
If\/ $\bs f$ is flat then $\bs f_-$ is flat, and the following
diagram is $2$-Cartesian\I{2-category!2-Cartesian square}
in\/~$\dSpac\!:$
\e
\begin{split}
\xymatrix@C=120pt@R=10pt{ \pd_-^{\bs f}\rX \ar[r]_(0.2){\bs f_-}
\ar[d]_{\bs i_\rX \vert_{\pd_-^{\bs f}\rX}}
\drtwocell_{}\omit^{}\omit{^{\id_{\bs i_\rY\ci
\bs f_-}\,\,\,\,\,\,\,\,\,\,\,\,\,\,{}}}
& \pd\rY \ar[d]^{\bs i_\rY} \\ \rX \ar[r]^(0.7){\bs f} & \rY.}
\end{split}
\label{ds6eq4}
\e
\item[{\bf(c)}] Let\/ $\bs g:\rX\ra\rY$ be another $1$-morphism,
and\/ $\eta:\bs f\Ra\bs g$ a $2$-morphism in $\dSpac$. Then $\bs
g$ is also semisimple, with\/ $\pd_-^{\bs g}\rX=\pd_-^{\bs
f}\rX$. If\/ $\bs f$ is simple, or flat, then $\bs g$ is simple,
or flat, respectively. Part\/ {\bf(b)} defines $1$-morphisms
$\bs f_-,\bs g_-:\pd_-^{\bs f}\rX\ra\pd\rY$. There is a unique
$2$-morphism $\eta_-:\bs f_-\Ra\bs g_-$ in\/ $\dSpac$ such
that\/ $\id_{\bs i_\rY}*\eta_-\!=\!\eta*\id_{\bs
i_\rX\vert_{\pd_-^{\bs f}\rX}}:\bs i_\rY\!\ci\!\bs f_- \Ra\bs
i_\rY\!\ci\!\bs g_-$.
\end{itemize}
\label{ds6thm1}
\end{thm}

We also show that the maps $\bs f\mapsto\bs f_-$,
$\eta\mapsto\eta_-$ in Theorem \ref{ds6thm1} are functorial, in that
they commute with compositions of 1- and 2-morphisms, and take
identities to identities. For simple 1-morphisms, this implies:

\begin{cor} Write $\mathop{\bf dSpa^c_{si}}$ for the
$2$-subcategory of\/ $\dSpac$ with arbitrary objects and\/
$2$-morphisms, but only simple $1$-morphisms. Then there is a
strict\/ $2$-functor\/ $\pd:\mathop{\bf dSpa^c_{si}}\ra \mathop{\bf
dSpa^c_{si}}$ mapping $\rX\mapsto\pd\rX$ on objects, $\bs
f\mapsto\bs f_-$ on (simple)\/ $1$-morphisms, and\/
$\eta\mapsto\eta_-$ on $2$-morphisms.
\label{ds6cor}
\end{cor}

Thus, boundaries in $\dSpac$ have strong functoriality properties.

\begin{rem} According to the general philosophy of working in
2-categories, when one constructs an object with some property in a
2-category, it is usually unique only up to equivalence. When one
constructs a 1-morphism with some property in a 2-category, it is
usually unique only up to 2-isomorphism. When one considers diagrams
of 1-morphisms in a 2-category, they usually commute only up to
(specified) 2-isomorphisms.

From this point of view, Theorem \ref{ds6thm1}(b) looks unnatural,
as it gives a 1-morphism $\bs f_-$ which is unique, not just up to
2-isomorphism, and a 1-morphism diagram \eq{ds6eq4} which commutes
strictly, not just up to 2-isomorphism.

In fact, this unnaturalness pervades our treatment of boundaries. In
our definition of d-space with corners $\rX=(\bX,\bpX,\bs
i_\rX,\om_\rX)$, the conditions on the 1-morphism $\bs
i_\rX:\bpX\ra\bX$ depend on $\bpX$ up to 1-isomorphism in $\dSpa$,
rather than up to equivalence, and depend on $\bs i_\rX$ up to
equality, not just up to 2-isomorphism. Boundaries $\pd\rX$ are
natural up to 1-isomorphism in $\dSpac$, not up to equivalence, and
1-morphisms $\bs i_\rX:\pd\rX\ra\rX$ natural up to equality.

The author chose this definition of $\dSpac$ for its (comparative!)
simplicity. In defining objects $\rX,\rY$, 1-morphisms $\bs f$, and
2-morphisms $\eta$ in $\dSpac$, we must impose extra conditions, and
possibly include extra data, over $\pd\rX,\pd\rY$. If these
conditions/extra data are imposed weakly, up to equivalence of
objects or 2-isomorphism of 1-morphisms, things rapidly become very
complicated and unwieldy. For instance, 1-morphisms in $\dSpac$
would comprise not just a 1-morphism $\bs f:\bX\ra\bY$ in $\dSpa$,
but also extra 2-morphism data over~$\uS_{\bs f},\uT_{\bs f}$.

So as a matter of policy, we generally do constructions involving
boundaries or corners in $\dSpac$ strictly, up to 1-isomorphism of
objects, and equality of 1-morphisms. One advantage of this is that
1-morphisms $\bs f:\rX\ra\rY$ and 2-morphisms $\eta:\bs f\Ra\bs g$
in $\dSpac$ are special examples of 1- and 2-morphisms in $\dSpa$ of
the underlying d-spaces $\bX,\bY$, rather than also containing
further data over $\pd\rX,\pd\rY$. Another advantage is that
boundaries in $\dSpac$ behave in a strictly functorial
way,\I{d-space with corners!boundary!strictly functorial} as in
Corollary \ref{ds6cor}, rather than weakly functorial.\I{d-space
with corners!simple 1-morphism|)}\I{d-space with corners!semisimple
1-morphism|)}\I{d-space with corners!flat 1-morphism|)}
\label{ds6rem2}
\end{rem}

\subsection{Manifolds with corners as d-spaces with corners}
\label{ds63}
\I{manifold with corners!as d-space with corners|(}\I{d-space with
corners!include manifolds with corners|(}

In \cite[\S 6.4]{Joyc6} we define a (2-)functor
$F_\Manc^\dSpac:\Manc\ra\dSpac$ from manifolds with corners to
d-spaces with corners.

\begin{dfn} Let $X$ be a manifold with corners. Then the
boundary $\pd X$ is a manifold with corners, with a smooth map
$i_{\pd X}:\pd X\ra X$. We will define a d-space with corners
$\rX=(\bX,\bpX,\bs i_\rX,\om_\rX)$. Set $\bX,\bpX,\bs
i_\rX=F_\Manc^\dSpa(X,\pd X,i_X)$. Then the conormal bundle
$\cN_\rX$ in \eq{ds6eq2} is the lift to the $C^\iy$-scheme $\upX$ of
the conormal line bundle $\cN_X$ of $\pd X$ in $X$, as in
\eq{ds6eq1}. Let $\om_\rX$ be the orientation on $\cN_\rX$
corresponding to that on $\cN_X$ induced by outward-pointing normal
vectors to $\pd X$ in $X$. Then $\rX$ is a d-space with corners.
Set~$F_\Manc^\dSpac(X)=\rX$.

Let $f:X\ra Y$ be a morphism in $\Manc$, and set
$\rX,\rY=F_\Manc^\dSpac(X,Y)$. Write $\bs
f=F_\Manc^\dSpa(f):\bX\ra\bY$, as a 1-morphism of d-spaces. Then
$\bs f:\rX\ra\rY$ is a 1-morphism of d-spaces with corners.
Define~$F_\Manc^\dSpac(f)=\bs f$.

The only 2-morphisms in $\Manc$, regarded as a 2-category, are
identity 2-morphisms $\id_f:f\Ra f$ for smooth $f:X\ra Y$. We
define~$F_\Manc^\dSpac(\id_f)=\id_{\bs f}$.

Define $F_\Man^\bdSpa:\Man\ra\bdSpa$ and $F_\Manb^\dSpab:
\Manb\ra\dSpab$ to be the restrictions of $F_\Manc^\dSpac$ to the
subcategories~$\Man,\Manb\subset\Manc$.

Write $\bMan,\bManb,\bManc$\G[Manc'']{$\bManc$}{2-subcategory of
d-spaces with corners equivalent to manifolds with corners} for the
full 2-subcategories of objects $\rX$ in $\dSpac$ equivalent to
$F_\Manc^\dSpac(X)$ for some manifold $X$ without boundary, or with
boundary, or with corners, respectively. Then $\bMan\subset\bdSpa$,
$\bManb\subset\dSpab$ and $\bManc\subset\dSpac$. When we say that a
d-space with corners $\rX$ {\it is a manifold}, we mean
that~$\rX\in\bManc$.\I{d-space with corners!is a manifold}
\label{ds6def2}
\end{dfn}

In \cite[\S 6.4]{Joyc6} we show that $F_\Man^\bdSpa: \Man\ra\bdSpa,$
$F_\Manb^\dSpab:\Manb\ra\dSpab$ and $F_\Manc^\dSpac: \Manc\ra\dSpac$
are full and faithful\I{functor!full}\I{functor!faithful} strict
2-functors.\I{2-category!strict 2-functor} We also prove that if $X$
is a manifold with corners, then there is a natural 1-isomorphism
$F_\Manc^\dSpac(\pd X)\cong \pd F_\Manc^\dSpac(X)$, and if $f:X\ra
Y$ is a smooth map of manifolds with corners and $\bs
f=F_\Manc^\dSpac(f)$, then $f$ is simple, semisimple or flat in
$\Manc$ if and only if $\bs f$ is simple, semisimple or flat in
$\dSpac$, respectively.\I{manifold with corners!as d-space with
corners|)}\I{d-space with corners!include manifolds with
corners|)}\I{d-space with corners!simple 1-morphism}\I{d-space with
corners!semisimple 1-morphism}\I{d-space with corners!flat
1-morphism}

\subsection{Equivalences, and gluing by equivalences}
\label{ds64}
\I{d-space with corners!equivalence|(}\I{d-space with corners!gluing
by equivalences|(}

In \cite[\S 6.5 \& \S 6.6]{Joyc6} we discuss {\it equivalences\/} in
$\dSpac$. First we characterize when a 1-morphism $\bs f:\rX\ra\rY$
in $\dSpac$ is an equivalence, in terms of the underlying 1-morphism
in~$\dSpa$:

\begin{prop}{\bf(a)} Suppose\/ $\bs f:\rX\ra\rY$ is an equivalence
in $\dSpac$. Then $\bs f$ is simple and flat, and\/ $\bs
f:\bX\ra\bY$ is an equivalence in\/ $\dSpa,$ where
$\rX=(\bX,\bpX,\bs i_\rX,\om_\rX)$ and\/ $\rY=(\bY,\bpY,\bs
i_\rY,\om_\rY)$. Also\/ $\bs f_-:\pd\rX\ra\pd\rY$ in Theorem\/
{\rm\ref{ds6thm1}(b)} is an equivalence in\/~$\dSpac$.
\smallskip

\noindent{\bf(b)} Let\/ $\bs f:\rX\ra\rY$ be a simple, flat\/
$1$-morphism in $\dSpac$ with\/ $\bs f:\bX\ra\bY$ an equivalence
in\/ $\dSpa$. Then $\bs f$ is an equivalence in\/~$\dSpac$.
\label{ds6prop1}
\end{prop}

Then we consider gluing d-spaces with corners by equivalences, as
for d-spaces in \S\ref{ds32}. The story is the same. Here is the
analogue of Definition~\ref{ds3def2}:

\begin{dfn} Let $\rX=(\bX,\bpX,\bs i_\rX,\om_\rX)$ be a d-space
with corners. Suppose $\bU\subseteq\bX$ is an open d-subspace in
$\dSpa$. Define $\bs{\pd U}=\bs i_\rX^{-1}(\bU)$, as an open
d-subspace of $\bpX$, and $\bs i_\rU:\bs{\pd U}\ra\bU$ by $\bs
i_\rU=\bs i_\rX\vert_{\bs{\pd U}}$. Then $\upU\subseteq\upX$ is an
open $C^\iy$-subscheme, and the conormal bundle of $\bs{\pd U}$ in
$\bU$ is $\cN_\rU=\cN_\rX\vert_\upU$ in $\qcoh(\upU)$. Define an
orientation $\om_\rU$ on $\cN_\rU$ by $\om_\rU=\om_\rX\vert_\upU$.
Write $\rU=(\bU,\bs{\pd U},\bs i_\rU,\om_\rU)$. Then $\rU$ is a
d-space with corners. We call $\rU$ an {\it open d-subspace\/} of
$\rX$. An {\it open cover\/} of $\rX$ is a family $\{\rU_a:a\in A\}$
of open d-subspaces $\rU_a$ of $\rX$ with~$\uX=\bigcup_{a\in
A}\uU_a$.\I{d-space with corners!open d-subspace}\I{d-space with
corners!open cover}
\label{ds6def3}
\end{dfn}

\begin{thm} Proposition\/ {\rm\ref{ds3prop}} and Theorems\/
{\rm\ref{ds3thm2}} and\/ {\rm\ref{ds3thm3}} hold without change in
the $2$-category\/ $\dSpac$ of d-spaces with corners.\I{d-space with
corners!equivalence|)}\I{d-space with corners!gluing by
equivalences|)}
\label{ds6thm2}
\end{thm}

\subsection{Corners and the corner functors}
\label{ds65}
\I{d-space with corners!k-corners C_k(X)@$k$-corners
$C_k(\rX)$|(}\I{d-space with corners!corner functors|(}

In \cite[\S 6.7]{Joyc6} we extend the material of \S\ref{ds53} on
corners and the corner functors from $\Manc$ to $\dSpac$. The next
theorem summarizes our results.

\begin{thm}{\bf(a)} Let\/ $\rX$ be a d-space with corners. Then
for each\/ $k=0,1,\ldots,$ we can define a d-space with corners
$C_k(\rX)$ called the $k$-\begin{bfseries}corners\end{bfseries} of\/
$\rX,$ and a\/ $1$-morphism $\bs\Pi^k_\rX:C_k(\rX)\ra\rX$ in
$\dSpac$. It has topological space
\e
\begin{split}
C_k(X)=\bigl\{&(x,\{x_1',\ldots,x_k'\}):x\in\bX,\;\>
x_1',\ldots,x_k'\in\bpX, \\
&\,\,\bs i_\rX(x_a')=x,\;\> a=1,\ldots,k,\;\> \text{$x_1',\ldots,x_k'$
are distinct\/}\bigr\}.
\end{split}
\label{ds6eq5}
\e
There is a natural, free action of the symmetric group $S_k$ on
$\pd^k\rX,$ and a $1$-isomorphism $C_k(\rX)\cong\pd^k\rX/S_k$. We
have $1$-isomorphisms\/ $C_0(\rX)\cong\rX$ and\/
$C_1(\rX)\cong\pd\rX$ in $\dSpac$. Write\/ $C(\rX)=\coprod_{k=0}^\iy
C_k(\rX)$ and\/ $\bs\Pi_\rX=\coprod_{k=0}^\iy\bs\Pi^k_\rX,$ so
that\/ $C(\rX)$ is a d-space with corners and\/
$\bs\Pi_\rX:C(\rX)\ra\rX$ is a $1$-morphism.
\smallskip

\noindent{\bf(b)} Let\/ $\bs f:\rX\ra\rY$ be a $1$-morphism of
d-spaces with corners. Then there is a unique $1$-morphism $C(\bs
f):C(\rX)\ra C(\rY)$ in $\dSpac$ such that\/ $\bs\Pi_\rY\ci C(\bs
f)=\bs f\ci\bs\Pi_\rX:C(\rX)\ra\rY,$ and\/ $C(\bs f)$ acts on points
as in \eq{ds6eq5} by
\begin{align*}
C(f):\,&\bigl(x,\{x_1',\ldots,x_k'\}\bigr)\longmapsto
\bigl(y,\{y_1',\ldots,y_l'\}\bigr),\quad\text{where}\\
&\{y_1',\ldots,y_l'\}\!=\!\bigl\{y':(x_i',y')\in\uS_{\bs f},\;
\text{some $i=1,\ldots,k$}\bigr\}.
\end{align*}
For all\/ $k,l\ge 0,$ write\/ $C_k^{\bs f,l}(\rX)=C_k(\rX)\cap C(\bs
f)^{-1}(C_l(\rY)),$ so that\/ $C_k^{\bs f,l}(\rX)$ is open and
closed in\/ $C_k(\rX)$ with\/ $C_k(\rX)=\coprod_{l=0}^\iy C_k^{\bs
f,l}(\rX),$ and write $C^l_k(\bs f)=C(\bs f)\vert_{C_k^{\bs
f,l}(\rX)},$ so that\/ $C^l_k(\bs f):C_k^{\bs f,l}(\rX)\ra C_l(\rY)$
is a $1$-morphism in~$\dSpac$.

\smallskip

\noindent{\bf(c)} Let\/ $\bs f,\bs g:\rX\ra\rY$ be $1$-morphisms
and\/ $\eta:\bs f\Ra\bs g$ a $2$-morphism in $\dSpac$. Then there
exists a unique $2$-morphism $C(\eta):C(\bs f)\Ra C(\bs g)$ in
$\dSpac,$ where $C(\bs f),C(\bs g)$ are as in {\bf(b)\rm,} such that
\begin{equation*}
\id_{\bs\Pi_\rY}*C(\eta)=\eta*\id_{\bs\Pi_\rX}:\bs\Pi_\rY\ci C(\bs f)
=\bs f\ci\bs\Pi_\rX\Longra\bs\Pi_\rY\ci C(\bs g)=\bs
g\ci\bs\Pi_\rX.
\end{equation*}

\noindent{\bf(d)} Define $C:\dSpac\ra\dSpac$ by $C:\rX\mapsto
C(\rX)$ on objects, $C:\bs f\mapsto C(\bs f)$ on\/ $1$-morphisms,
and\/ $C:\eta\mapsto C(\eta)$ on\/ $2$-morphisms, where
$C(\rX),C(\bs f), C(\eta)$ are as in {\bf(a)--(c)} above. Then $C$
is a strict\/ $2$-functor, called a
\begin{bfseries}corner functor\end{bfseries}.\G[CCb]{$C,\hat
C:\dSpac\ra\dSpac$}{`corner functors' for d-spaces with corners}

\smallskip

\noindent{\bf(e)} Let\/ $\bs f:\rX\ra\rY$ be semisimple. Then $C(\bs
f)$ maps $C_k(\rX)\ra\coprod_{l=0}^kC_l(\rY)$ for all\/ $k\ge 0$.
The natural\/ $1$-isomorphisms $C_1(\rX)\cong\pd\rX,$
$C_0(\rY)\cong\rY,$ $C_1(\rY)\cong\pd\rY$ identify $C_1^{\bs
f,0}(\rX)\cong\pd^{\bs f}_+\rX,$ $C_1^{\bs f,1}(\rX)\cong\pd^{\bs
f}_-\rX,$ $C_1^0(\bs f)\cong\bs f_+$ and\/~$C_1^1(\bs f)\cong\bs
f_-$.

If\/ $\bs f$ is simple then $C(\bs f)$ maps $C_k(\rX)\ra C_k(\rY)$
for all\/~$k\ge 0$.
\smallskip

\noindent{\bf(f)} Analogues of\/ {\bf(b)}--{\bf(d)} also hold for a
second corner functor\/ $\hat C:\dSpac\ra\dSpac,$ which acts on
objects by $\hat C:\rX\mapsto C(\rX)$ in {\bf(a)\rm,} and for
$1$-morphisms $\bs f:\rX\ra\rY$ in {\bf(b)\rm,} $\hat C(\bs
f):C(\rX)\ra C(\rY)$ acts on points by
\begin{align*}
&\hat C(f):\bigl(x,\{x_1',\ldots,x_k'\}\bigr)\longmapsto
\bigl(y,\{y_1',\ldots,y_l'\}\bigr),\quad\text{where}\\
&\{y_1',\ldots,y_l'\}\!=\!\bigl\{y':(x_i',y')\in\uS_{\bs f},\;
\text{some $i\!=\!1,\ldots,k$}\bigr\}\!\cup\!\bigl\{y':(x,y')\!\in
\!\uT_{\bs f}\bigr\}.
\end{align*}

If\/ $\bs f:\rX\ra\rY$ is flat then\/~$\hat C(\bs f)=C(\bs f)$.
\label{ds6thm3}
\end{thm}

The comments of Remark \ref{ds6rem2} also apply to Theorem
\ref{ds6thm3}: our construction characterizes $C_k(\rX)$ up to
1-isomorphism in $\dSpac$, not just up to equivalence, the
1-morphisms $C(\bs f),\hat C(\bs f)$ are characterized up to
equality, not just up to 2-isomorphism, and in $\bs\Pi_\rY\ci C(\bs
f)=\bs f\ci\bs\Pi_\rX$ we require the 1-morphisms to be equal, not
just 2-isomorphic. This may seem unnatural from a 2-category point
of view, but it has the advantage that corners are strictly
2-functorial rather than weakly 2-functorial.\I{d-space with
corners!k-corners C_k(X)@$k$-corners $C_k(\rX)$|)}\I{d-space with
corners!corner functors|)}

\subsection{Fibre products in $\dSpac$}
\label{ds66}
\I{d-space with corners!fibre products|(}\I{2-category!fibre
products in|(}

In \cite[\S 6.8--\S 6.9]{Joyc6} we study {\it fibre products\/} in
$\dSpac$. Here the situation is more complex than for d-spaces. As
in \S\ref{ds32}, all fibre products exist in $\dSpa$, but this fails
for $\dSpac$. The problem is that in a fibre product $\rW=\rX\t_{\bs
g,\rZ,\bs h}\rY$ in $\dSpac$, the boundary $\pd\rW$ depends in a
complicated way on $\rX,\rY,\rZ,\pd\rX,\ab \pd\rY,\ab\pd\rZ$, and
sometimes there is no good candidate for $\pd\rW$. Here is an
example.

\begin{ex} Let $X=Y=[0,\iy)\t\R$ and $Z=[0,\iy)^2\t\R$, as
manifolds with corners, and define smooth maps $g:X\ra Z$ and
$h:Y\ra Z$ by $g(u,v)=(u,u,v)$ and $h(u,v)=(u,e^vu,v)$.
Set~$\rX,\rY,\rZ,\bs g,\bs h=F_\Manc^\dSpac(X,Y,Z,g,h)$.

In \cite[\S 6.8.6]{Joyc6} we show that no fibre product
$\rW=\rX\t_{\bs g,\rZ,\bs h}\rY$ exists in $\dSpac$. We do this by
showing that $\pd\rW$ would have to have exactly one point, lying
over $(0,0)\in X$ and $(0,0)\in Y$, which is the only point in
$X\t_ZY$ where normal vectors to $\pd X,\pd Y$ in $X,Y$ project
under $\d g,\d h$ to parallel vectors in $TZ$. But this would
contradict other properties of~$\pd\rW$.\I{d-space with
corners!fibre products!may not exist}
\label{ds6ex}
\end{ex}

So, we would like to find useful sufficient conditions for existence
of fibre products $\rX\t_{\bs g,\rZ,\bs h}\rY$ in $\dSpac$; and
these conditions should be wholly to do with boundaries, since we
already know that fibre products exist in $\dSpa$. In \cite[\S
6.8.1]{Joyc6} we define two such sufficient conditions on $\bs g,\bs
h$, called {\it b-transversality\/} and {\it
c-transversality}.\I{b-transversality|(}\I{c-transversality|(}\I{d-space
with corners!b-transverse 1-morphisms|(}\I{d-space with
corners!c-transverse 1-morphisms|(}

\begin{dfn} Let $\bs g:\rX\ra\rZ$ and $\bs h:\rY\ra\rZ$ be
1-morphisms in $\dSpac$. As in \S\ref{ds61} we have line bundles
$\cN_\rX,\cN_\rZ$ over the $C^\iy$-schemes $\upX,\upZ$, and a
$C^\iy$-subscheme $\uS_{\bs g}\subseteq\upX\t_\uZ\upZ$. As in
\cite[\S 7.1]{Joyc6}, there is a natural isomorphism $\la_{\bs
g}:\uu_{\bs g}^*(\cN_\rZ)\ra \us_{\bs f}^*(\cN_\rX)$ in
$\qcoh(\uS_{\bs g})$. The same holds for~$\bs h$.

We say that $\bs g,\bs h$ are {\it b-transverse\/} if whenever
$x\in\uX$ and $y\in\uY$ with $\ug(x)=\uh(y)=z\in\uZ$, the following
morphism in $\qcoh(\ul{*})$ is injective:
\begin{align*}
&\bigop_{(x',z')\in\uS_{\bs g}:\ui_\rX(x')=x}\!\!\! \la_{\bs
g}\vert_{(x',z')}\op\bigop_{(y',z')\in\uS_{\bs
h}:\ui_\rY(y')=y}\!\!\! \la_{\bs h}\vert_{(y',z')}:\\
&\bigop_{z'\in\ui_\rZ^{-1}(z)}\cN_\rZ\vert_{z'}\longra
\bigop_{x'\in\ui_\rX^{-1}(x)}\cN_\rX\vert_{x'}\op
\bigop_{y'\in\ui_\rY^{-1}(y)}\cN_\rY\vert_{y'}.
\end{align*}
Roughly speaking, this says that the corners of $\rX,\rY$ are
transverse to the corners of $\rZ$. In Example \ref{ds6ex}, this
condition fails at $x=0\in\uX$ and $y=0\in\uY$, so $\bs g,\bs h$ are
not b-transverse.

We call $\bs g,\bs h$ {\it c-transverse\/} if the following two
conditions hold, using the notation of Theorem~\ref{ds6thm3}:
\begin{itemize}
\setlength{\itemsep}{0pt}
\setlength{\parsep}{0pt}
\item[(a)] whenever there are points in
$C_j(\rX),C_k(\rY),C_l(\rZ)$ with
\begin{equation*}
C(\bs g)(x,\{x_1',\ldots,x_j'\})=C(\bs h)(y,\{y_1',\ldots,
y_k'\})=(z,\{z_1',\ldots,z_l'\}),
\end{equation*}
we have either $j+k>l$ or $j=k=l=0;$ and
\item[(b)] whenever there are points in
$C_j(\rX),C_k(\rY),C_l(\rZ)$ with
\begin{equation*}
\hat C(\bs g)(x,\{x_1',\ldots,x_j'\})=\hat C(\bs
h)(y,\{y_1',\ldots, y_k'\})=(z,\{z_1',\ldots,z_l'\}),
\end{equation*}
we have~$j+k\ge l$.
\end{itemize}
\label{ds6def4}
\end{dfn}

Here b-transversality is a continuous condition on $\bs g,\bs h$,
and c-transversality is a discrete condition. Also c-transversality
implies b-transversality (though this is not obvious). Part (a)
corresponds to the condition in Definition \ref{ds5def7} for
transverse $g,h$ in $\Manc$ to be strongly transverse. We can show:

\begin{lem} Let\/ $\bs g:\rX\ra\rZ$ and\/ $\bs h:\rY\ra\rZ$ be\/
$1$-morphisms in $\dSpac$. The following are sufficient conditions
for\/ $\bs g,\bs h$ to be c-transverse, and hence b-transverse:
\begin{itemize}
\setlength{\itemsep}{0pt}
\setlength{\parsep}{0pt}
\item[{\bf(i)}] $\bs g$ or $\bs h$ is semisimple and flat;
or\I{d-space with corners!semisimple 1-morphism}\I{d-space with
corners!flat 1-morphism}
\item[{\bf(ii)}] $\rZ$ is a d-space without boundary.
\end{itemize}
\label{ds6lem}
\end{lem}

We summarize the main results of \cite[\S 6.8]{Joyc6} on fibre
products in~$\dSpac$:\I{fibre product!of d-spaces with corners}

\begin{thm}{\bf(a)} All b-transverse fibre products exist
in\/~$\dSpac$.\I{d-space with corners!fibre products!b-transverse}
\smallskip

\noindent{\bf(b)} The $2$-functor $F_\Manc^\dSpac$ of\/
{\rm\S\ref{ds63}} takes transverse fibre products in $\Manc$ to
b-transverse fibre products in $\dSpac$. That is, if
\begin{equation*}
\xymatrix@C=90pt@R=10pt{ W \ar[r]_{f} \ar[d]^{e}
& Y \ar[d]_{h} \\
X \ar[r]^{g} & Z}
\end{equation*}
is a Cartesian square\I{category!Cartesian square} in $\Manc$ with\/
$g,h$ transverse, and\/ $\rW,\rX,\rY,\rZ,\ab\bs e,\ab\bs f,\ab\bs
g,\ab\bs h=F_\Manc^\dSpac(W,X,Y,Z,e,f,g,h),$ then
\begin{equation*}
\xymatrix@C=90pt@R=10pt{ \rW \ar[r]_(0.2){\bs f} \ar[d]^{\bs e}
\drtwocell_{}\omit^{}\omit{^{\id_{\bs g\ci\bs
e}\,\,\,\,\,\,\,\,\,{}}} & \rY \ar[d]_{\bs h} \\ \rX
\ar[r]^(0.7){\bs g} & \rZ}
\end{equation*}
is $2$-Cartesian\I{2-category!2-Cartesian square} in $\dSpac,$ with
$\bs g,\bs h$ b-transverse. If also $g,h$ are strongly transverse
in\/ $\Manc,$ then $\bs g,\bs h$ are c-transverse in\/~$\dSpac$.
\smallskip

\noindent{\bf(c)} Suppose we are given a $2$-Cartesian diagram
in\/~{\rm$\dSpac$:}
\begin{equation*}
\xymatrix@C=60pt@R=10pt{ \rW \ar[r]_(0.25){\bs f} \ar[d]^{\bs e}
\drtwocell_{}\omit^{}\omit{^{\eta}}
 & \rY \ar[d]_{\bs h} \\ \rX \ar[r]^(0.7){\bs g} & \rZ,}
\end{equation*}
with\/ $\bs g,\bs h$ c-transverse. Then the following are also
$2$-Cartesian in\/~$\dSpac\!:$\I{d-space with corners!fibre
products!boundary and corners|(}\I{d-space with corners!corner
functors}
\ea
\begin{gathered}
\xymatrix@C=100pt@R=10pt{ C(\rW) \ar[r]_(0.25){C(\bs f)}
\ar[d]^{C(\bs e)} \drtwocell_{}\omit^{}\omit{^{C(\eta)
\,\,\,\,\,\,\,\,{}}} & C(\rY) \ar[d]_{C(\bs h)} \\ C(\rX)
\ar[r]^(0.7){C(\bs g)} & C(\rZ),}
\end{gathered}
\label{ds6eq6}\\
\begin{gathered}
\xymatrix@C=100pt@R=10pt{ C(\rW) \ar[r]_(0.25){\hat C(\bs f)}
\ar[d]^{\hat C(\bs e)} \drtwocell_{}\omit^{}\omit{^{\hat C(\eta)
\,\,\,\,\,\,\,\,{}}} & C(\rY) \ar[d]_{\hat C(\bs h)} \\ C(\rX)
\ar[r]^(0.7){\hat C(\bs g)} & C(\rZ).}
\end{gathered}
\label{ds6eq7}
\ea
Also \eq{ds6eq6}--\eq{ds6eq7} preserve gradings, in that they relate
points in $C_i(\rW),\ab C_j(\rX),\ab C_k(\rY),\ab C_k(\rZ)$ with\/
$i=j+k-l$. Hence \eq{ds6eq6} implies equivalences
in\/~$\dSpac\!:$\I{d-space with corners!fibre products!boundary and
corners|)}
\ea
C_i(\rW)&\simeq \coprod_{j,k,l\ge 0:i=j+k-l} C_j^{\bs
g,l}(\rX)\t_{C_j^l(\bs g),C_l(\rZ),C_k^l(\bs h)}C_k^{\bs h,l}(\rY),
\label{ds6eq8}\\
\pd\rW&\simeq \coprod_{j,k,l\ge 0:j+k=l+1} C_j^{\bs
g,l}(\rX)\t_{C_j^l(\bs g),C_l(\rZ),C_k^l(\bs h)}C_k^{\bs h,l}(\rY).
\label{ds6eq9}
\ea
\label{ds6thm4}
\end{thm}

Part (a) takes some work to prove. For fibre products in $\dSpa$, as
in \S\ref{ds33}, we gave an explicit global construction. But for
fibre products in $\dSpac$, we first prove that local fibre products
$\rX\t_{\bs g,\rZ,\bs h}\rY$ exist in $\dSpac$ near each $x\in\rX$,
$y\in\rY$ with $\bs g(x)=\bs h(y)\in\rZ$, and then we use the
results of \S\ref{ds64} to glue these local fibre products by
equivalences to get a global fibre product.\I{2-category!fibre
products in|)}

For general b-transverse fibre products $\rW=\rX\t_{\bs g,\rZ,\bs
h}\rY$ in $\dSpac$, the description of $\pd\rW$ can be complicated.
For c-transverse fibre products, we do at least have a (still
complicated) explicit formula \eq{ds6eq9} for $\pd\rW$. Here are
some cases when this formula simplifies, an analogue of
Proposition~\ref{ds5prop4}.

\begin{prop} Let\/ $\bs g:\rX\ra\rZ$ and\/ $\bs h:\rY\ra\rZ$ be
$1$-morphisms of d-spaces with corners. Then:
\begin{itemize}
\setlength{\itemsep}{0pt}
\setlength{\parsep}{0pt}
\item[{\bf(a)}] If\/ $\pd\rZ=\bs\es$ then there is an equivalence
\e
\pd\bigl(\rX\t_{\bs g,\rZ,\bs h}\rY\bigr)\simeq \bigl(\pd
\rX\t_{\bs g\ci\bs i_\rX,\rZ,\bs
h}\rY\bigr)\amalg\bigl(\rX\t_{\bs g,\rZ,\bs h \ci\bs
i_\rY}\pd\rY\bigr).
\label{ds6eq10}
\e
\item[{\bf(b)}] If\/ $\bs g$ is semisimple and flat then there is an
equivalence
\e
\pd\bigl(\rX\t_{\bs g,\rZ,\bs h}\rY\bigr)\simeq \bigl(\pd^{\bs
g}_+\rX\t_{\bs g_+,\rZ,\bs h}\rY\bigr)\amalg\bigl(\rX\t_{\bs
g,\rZ, \bs h\ci\bs i_\rY}\pd\rY\bigr).
\label{ds6eq11}
\e
\item[{\bf(c)}] If both\/ $\bs g$ and\/ $\bs h$ are semisimple and
flat then there is an equivalence
\e
\begin{split}
\pd\bigl(\rX\t_{\bs g,\rZ,\bs h}\rY\bigr)\simeq \bigl(\pd_+^{\bs
g}\rX \t_{\bs g_+,\rZ,\bs h}\rY\bigr)\amalg \bigl(\rX
\t_{\bs g,\rZ,\bs h_+}\pd_+^{\bs h}\rY\bigr)&\\
\amalg\bigl(\pd^{\bs g}_-\rX\t_{\bs g_-,\pd\rZ,\bs h_-} \pd^{\bs
h}_-\rY\bigr)&.
\end{split}
\label{ds6eq12}
\e
\end{itemize}
Here all fibre products in \eq{ds6eq10}--\eq{ds6eq12} are
c-transverse, and so
exist.\I{b-transversality|)}\I{c-transversality|)}\I{d-space with
corners!fibre products|)}\I{d-space with corners!b-transverse
1-morphisms|)}\I{d-space with corners!c-transverse 1-morphisms|)}
\label{ds6prop2}
\end{prop}

\subsection{Fixed point loci in d-spaces with corners}
\label{ds67}
\I{d-space with corners!fixed point loci|(}

Section \ref{ds34} discussed the fixed d-subspace $\bX^\Ga$ of a
finite group $\Ga$ acting on a d-space $\bX$, and \S\ref{ds56}
considered fixed point loci $X^\Ga$ of a finite group $\Ga$ acting
on a manifold with corners $X$. In \cite[\S 6.10]{Joyc6} we
generalize these to d-spaces with corners. Here is the analogue of
Theorem~\ref{ds3thm5}.

\begin{thm} Let\/ $\rX$ be a d-space with corners, $\Ga$ a finite
group, and\/ $\bs r:\Ga\ra\Aut(\rX)$ an action of\/ $\Ga$ on\/ $\rX$
by $1$-isomorphisms. Then we can define a d-space with corners\/
$\rX^\Ga$\G[XGac]{$\rX^\Ga$}{fixed d-subspace of group $\Ga$ acting
on a d-space with corners $\rX$} called the \begin{bfseries}fixed
d-subspace of\/ $\Ga$ in\/\end{bfseries} $\rX,$ with an inclusion
$1$-morphism\/ $\bs j_{\rX,\Ga}:\rX^\Ga\ra\rX$.\G[jXGac]{$\bs
j_{\rX,\Ga}:\rX^\Ga\hookra \rX$}{inclusion of $\Ga$-fixed d-subspace
$\rX^\Ga$ in a d-space with corners $\rX$} It has the following
properties:
\begin{itemize}
\setlength{\itemsep}{0pt}
\setlength{\parsep}{0pt}
\item[{\bf(a)}] Let\/ $\rX,\Ga,\bs r$ and\/ $\bs
j_{\rX,\Ga}:\rX^\Ga\ra\rX$ be as above. Suppose $\bs
f:\rW\ra\rX$ is a $1$-morphism in $\dSpac$. Then $\bs f$
factorizes as $\bs f=\bs j_{\rX,\Ga}\ci\bs g$ for some
$1$-morphism $\bs g:\rW\ra\rX^\Ga$ in $\dSpac,$ which must be
unique, if and only if\/ $\bs r(\ga)\ci\bs f=\bs f$ for
all\/~$\ga\in\Ga$.

\item[{\bf(b)}] Suppose $\rX,\rY$ are d-spaces with corners,
$\Ga$ is a finite group, $\bs r:\Ga\ra\Aut(\rX),$ $\bs
s:\Ga\ra\Aut(\rY)$ are actions of\/ $\Ga$ on $\rX,\rY,$ and\/
$\bs f:\rX\ra\rY$ is a $\Ga$-equivariant\/ $1$-morphism in
$\dSpac,$ that is, $\bs f\ci\bs r(\ga)=\bs s(\ga)\ci\bs f$ for
all\/ $\ga\in\Ga$. Then there exists a unique\/ $1$-morphism
$\bs f^\Ga:\rX^\Ga\ra\rY^\Ga$ such that\/~$\bs j_{\rY,\Ga}\ci\bs
f^\Ga=\bs f\ci\bs j_{\rX,\Ga}$.
\item[{\bf(c)}] Let\/ $\bs f,\bs g:\rX\ra\rY$ be
$\Ga$-equivariant\/ $1$-morphisms as in {\bf(b)\rm,} and\/
$\eta:\bs f\Ra\bs g$ be a $\Ga$-equivariant\/ $2$-morphism, that
is, $\eta*\id_{\bs r(\ga)}=\id_{\bs s(\ga)}*\eta$ for all\/
$\ga\in\Ga$. Then there exists a unique $2$-morphism
$\eta^\Ga:\bs f^\Ga\Ra\bs g^\Ga$ such that\/~$\id_{\bs
j_{\rY,\Ga}}*\eta^\Ga=\eta*\id_{\bs j_{\rX,\Ga}}$.
\end{itemize}
Note that\/ {\bf(a)} is a universal property that determines
$\bX^\Ga,\bs j_{\bX,\Ga}$ up to canonical\/ $1$-isomorphism.
\label{ds6thm5}
\end{thm}

As for manifolds with corners in \S\ref{ds56}, in general
$\pd(\rX^\Ga)\not\simeq(\pd\rX)^\Ga$, so fixed point loci do not
commute with boundaries. But the following analogue of Proposition
\ref{ds5prop5}(b) shows that fixed point loci do commute with
corners.

\begin{prop} Let\/ $\rX$ be a d-space with corners, $\Ga$ a finite
group, and\/ $\bs r:\Ga\ra\Aut(\rX)$ an action of\/ $\Ga$ on\/
$\rX$. Applying the corner functor $C$ of\/ {\rm\S\ref{ds65}} gives
an action $C(\bs r):\Ga\ra\Aut(C(\rX))$. Hence Theorem\/
{\rm\ref{ds6thm5}} defines fixed d-subspaces $\rX{}^\Ga,C(\rX)^\Ga$
and inclusion $1$-morphisms $\bs j_{\rX,\Ga}:\rX{}^\Ga\ra\rX,$ $\bs
j_{C(\rX),\Ga}:C(\rX)^\Ga\ra C(\rX)$. Applying $C$ to $\bs
j_{\rX,\Ga}$ also gives $C(\bs j_{\rX,\Ga}):C(\rX{}^\Ga)\ra C(\rX)$.

Then there exists a unique equivalence $\bs
k_{\rX,\Ga}:C(\rX{}^\Ga)\ra C(\rX)^\Ga$ in $\dSpac$ such
that\/~$C(\bs j_{\rX,\Ga})=\bs j_{C(\rX),\Ga}\ci\bs k_{\rX,\Ga}$.
\label{ds6prop3}
\end{prop}

We will use fixed d-subspaces $\rX^\Ga$ in \S\ref{ds137} below to
describe orbifold strata $\eX^\Ga$ of quotient d-stacks with
corners\I{d-stack with corners!quotients $[\rX/G]$} $\eX=[\rX/G]$.
If $\rX$ is a d-manifold with corners, as in \S\ref{ds7}, then in
general the fixed d-subspaces $\rX^\Ga$ are disjoint unions of
d-manifolds with corners of different dimensions, that is, $\rX^\Ga$
lies in~$\dcManc$.\I{d-space with corners|)}\I{d-space with
corners!fixed point loci|)}

\section{D-manifolds with corners}
\label{ds7}
\I{d-manifold with corners|(}

We can now define the 2-categories $\dManb$ of {\it d-manifolds with
boundary\/} and $\dManc$ of {\it d-manifolds with corners}, which
are derived versions of manifolds with boundary and with corners,
following~\cite[Chap.~7]{Joyc6}.

\subsection{The definition of d-manifolds with corners}
\label{ds71}
\I{d-manifold with corners!definition|(}

In \S\ref{ds41} we defined a d-manifold to be a d-space covered by
principal open d-submanifolds of fixed dimension, where Proposition
\ref{ds4prop1} gave three equivalent definitions of principal
d-manifolds,\I{d-manifold!principal} the first as a fibre product
$\bX\t_\bZ\bY$ in $\dSpa$ with $\bX,\bY,\bZ\in\hMan$, and the third
as a fibre product $\bV\t_{\bs s,\bs E,\bs 0}\bV$ in $\dSpa$, where
$V$ is a manifold, $E\ra V$ a vector bundle, and~$s\in C^\iy(E)$.

When we pass to d-spaces and d-manifolds with corners in \cite[\S
7.1]{Joyc6}, the analogues of Proposition \ref{ds4prop1}(a)--(c) are
no longer equivalent. So we have to choose which of them gives the
best idea of principal d-manifold with corners. Defining principal
d-manifolds with corners to be fibre products $\rX\t_\rZ\rY$ in
$\dSpac$ with $\rX,\rY,\rZ\in\bManc$ is unsatisfactory, since as in
\S\ref{ds66} fibre products $\rX\t_\rZ\rY$ may not exist in
$\dSpac$. So instead we define principal d-manifolds with corners to
be fibre products $\rV\t_{\bs s,\rE,\bs 0}\rV$ in~$\dSpac$.

\begin{dfn} A d-space with corners $\rW$ is called a {\it principal
d-manifold with corners\/}\I{d-manifold with
corners!principal}\I{principal d-manifold with corners} if is
equivalent in $\dSpac$ to a fibre product $\rV\t_{\bs s,\rE,\bs
0}\rV$, where $V$ is a manifold with corners, $E\ra V$ is a vector
bundle, $s:V\ra E$ is a smooth section of $E,$ $0:V\ra E$ is the
zero section, and $\rV,\rE,\bs s,\bs 0=F_\Manc^\dSpac(V,E,s,0)$.
Note that $s,0:V\ra E$ are simple, flat smooth maps in $\Manc$, so
$\bs s,\bs 0:\rV\ra\rE$ are simple, flat 1-morphisms in $\dSpac$,
and thus $\bs s,\bs 0$ are b-transverse by Lemma \ref{ds6lem}(a),
and the fibre product $\rV\t_{\bs s,\rE,\bs 0}\rV$ exists in
$\dSpac$ by Theorem~\ref{ds6thm4}(a).

If $\rW\simeq\rV\t_{\bs s,\rE,\bs 0}\rV$ then the virtual cotangent
sheaf $T^*\bW$ of the d-space $\bW$ is a virtual vector bundle with
$\rank T^*\bW=\dim V-\rank E$. Hence, if $\rW\ne\bs\es$ then the
integer $\dim V-\rank E$ depends only on $\bW$ up to equivalence in
$\dSpa$, and is independent of the choice of $V,E,s$ with
$\rW\simeq\rV\t_{\bs s,\rE,\bs 0}\rV$. Define the {\it virtual
dimension\/}\I{d-manifold with corners!virtual dimension} $\vdim\rW$
to be~$\vdim\rW=\rank T^*\bW=\dim V-\rank E$.

A d-space with corners $\rX$ is called a {\it d-manifold with
corners of virtual dimension\/} $n\in\Z$, written $\vdim\rX=n$, if
$\rX$ can be covered by open d-subspaces $\rW$ which are principal
d-manifolds with corners with $\vdim\rW=n$. A d-manifold with
corners $\rX$ is called a {\it d-manifold with
boundary\/}\I{d-manifold with boundary} if it is a d-space with
boundary, and a {\it d-manifold without boundary\/} if it is a
d-space without boundary.

Write $\bdMan,\dManb,\dManc$\G[dManc]{$\dManc$}{2-category of
d-manifolds with corners}\G[dManb]{$\dManb$}{2-category of
d-manifolds with boundary}\G[dMan']{$\bdMan$}{2-subcategory of
d-manifolds with corners equivalent to d-manifolds} for the full
2-subcategories of d-manifolds without boundary, and d-manifolds
with boundary, and d-manifolds with corners in $\dSpac$,
respectively. The 2-functor $F_\dSpa^\dSpac:\dSpa\ra\dSpac$ in
\S\ref{ds61} is an isomorphism of 2-categories $\dSpa\ra\bdSpa$, and
its restriction to $\dMan\subset\dSpa$ gives an isomorphism of
2-categories $F_\dMan^\dManc:\dMan\ra\bdMan\subset\dManc$. So we may
as well identify $\dMan$ with its image $\bdMan$, and consider
d-manifolds in \S\ref{ds4} as examples of d-manifolds with
corners.\I{d-manifold with corners!include
d-manifolds}\I{d-manifold!as d-manifold with corners}

If $\rX=(\bX,\bpX,\bs i_\rX,\om_\rX)$ is a d-manifold with corners,
then the virtual cotangent sheaf $T^*\bX$ of the d-space $\bX$ from
Definition \ref{ds4def3} is a virtual vector bundle\I{virtual vector
bundle} on $\uX$, of rank $\vdim\rX$. We will call
$T^*\bX\in\vvect(\uX)$ the {\it virtual cotangent
bundle\/}\I{d-manifold with corners!virtual cotangent
bundle}\I{virtual cotangent bundle} of $\rX$, and also write
it~$T^*\rX$.\G[T*Xb]{$T^*\rX$}{virtual cotangent sheaf of a d-space
with corners $\rX$}
\label{ds7def1}
\end{dfn}

Much of \S\ref{ds6} on d-spaces with corners applies immediately to
d-manifolds with corners. If $\rX$ is a d-manifold with corners with
$\vdim\rX=n$ then the boundary $\pd\rX$ as a d-space with corners
from \S\ref{ds61} is a d-manifold with corners, with
$\vdim\pd\rX=n-1$.\I{d-manifold with corners!boundary} The material
on simple, semisimple, and flat 1-morphisms in $\dSpac$ in
\S\ref{ds62} also holds in $\dManc$. The functor
$F_\Manc^\dSpac:\Manc\ra\dSpac$ in \S\ref{ds63} maps to
$\dManc\subset\dSpac$, so we write $F_\Manc^\dManc=
F_\Manc^\dSpac:\Manc\ra\dManc$. The 2-categories
$\bMan,\bManb,\bManc$ in Definition \ref{ds6def2} are
2-subcategories of $\bdMan,\ab\dManb,\ab\dManc$, respectively. When
we say that a d-manifold with corners $\rX$ {\it is a
manifold},\I{d-manifold with corners!is a manifold} we mean
that~$\rX\in\bManc$.

In \S\ref{ds64}, if we make a d-space with corners $\rY$ by gluing
together d-manifolds with corners $\rX_i$ for $i\in I$ by
equivalences, then $\rY$ is a d-manifold with corners with
$\vdim\rY=n$ provided $\vdim\rX_i=n$ for all~$i\in I$.

In \S\ref{ds65}, if $\rX$ is a d-manifold with corners with
$\vdim\rX=n$ then the $k$-corners $C_k(\rX)$ is a d-manifold with
corners, with $\vdim C_k(\rX)=n-k$. Note however that
$C(\rX)=\coprod_{k=0}^\iy C_k(\rX)$ in Theorem \ref{ds6thm3} is in
general {\it not\/} a d-manifold with corners, but only a disjoint
union of d-manifolds with corners with different dimensions. As for
$\cManc$ in \S\ref{ds53}, define
$\dcManc$\G[dManc']{$\dcManc$}{2-category of disjoint unions of
d-manifolds with corners of different dimensions} to be the full
2-subcategory of $\rX$ in $\dSpac$ which may be written as a
disjoint union $\rX=\coprod_{n\in\Z}\rX_n$ for $\rX_n\in\dManc$ with
$\vdim\rX_n=n$, where we allow $\rX_n=\bs\es$. We call such $\rX$ a
{\it d-manifold with corners of mixed dimension}.\I{d-manifold with
corners!of mixed dimension} Then $C,\hat C$ in Theorem \ref{ds6thm3}
restrict to strict 2-functors~$C,\hat C:\dManc\ra\dcManc$.\I{d-space
with corners!corner functors}\I{d-manifold with corners!corner
functors}\I{2-category!strict 2-functor}

Here are some examples. The fibre products we give all exist in
$\dManc$ by results in \S\ref{ds75} below.\I{d-manifold with
corners!fibre products}

\begin{ex}{\bf(i)} Let $\rX$ be the fibre product
$\bs{[0,\iy)}\t_{\bs i,\bR,\bs 0}\bs *$ in $\dManc$, where
$i:[0,\iy)\hookra\R$ is the inclusion. Then $\rX=(\bX,\bpX,\bs
i_\rX,\om_\rX)$ is `a point with point boundary', of virtual
dimension 0, and its boundary $\pd\rX$ is an `obstructed point', a
point with obstruction space $\R$, of virtual dimension~$-1$.

The conormal bundle $\cN_\rX$ of $\bpX$ in $\bX$ is the obstruction
space $\R$ of $\bpX$. In this case, the orientation $\om_\rX$ on
$\cN_\rX$ cannot be determined from $\bX,\bpX,\bs i_\rX$, in fact,
there is an automorphism of $\bX,\bpX,\bs i_\rX$ which reverses the
orientation of $\cN_\rX$. So $\om_\rX$ really is extra data. We
include $\om_\rX$ in the definition of d-manifolds with corners to
ensure that orientations of d-manifolds with corners are
well-behaved. If we omitted $\om_\rX$ from the definition, there
would exist oriented d-manifolds with corners $\rX$ whose boundaries
$\pd\rX$ are not orientable.
\smallskip

\noindent{\bf(ii)} The fibre product $\bs{[0,\iy)}\t_{\bs
i,\bs{[0,\iy)},\bs 0}\bs *$ is a point $\bs *$ without boundary. The
only difference with {\bf(i)} is that we have replaced the target
$\bR$ with $\bs{[0,\iy)}$, adding a boundary. So in a fibre product
$\rW=\rX\t_\rZ\rY$ in $\dManc$, the boundary of $\rZ$ affects the
boundary of $\rW$. This does not happen for fibre products
in~$\Manc$.

\smallskip

\noindent{\bf(iii)} Let $\rX'$ be the fibre product
$\bs{[0,\iy)}\t_{\bs i,\bR,\bs i}\bs{(-\iy,0]}$ in $\dManc$, that
is, the derived intersection of submanifolds $[0,\iy),(-\iy,0]$ in
$\R$. Topologically, $\rX'$ is just the point $\{0\}$, but as a
d-manifold with corners $\rX'$ has virtual dimension 1. The boundary
$\pd\rX'$ is the disjoint union of two copies of $\rX$ in {\bf(i)}.
The $C^\iy$-scheme $\uX'$ in $\rX'$ is the
spectrum\I{C-scheme@$C^\iy$-scheme!spectrum functor} of the
$C^\iy$-ring $C^\iy\bigl([0,\iy)^2\bigr)/(x+y)$, which is
infinite-dimensional, although its topological space is a
point.\I{d-manifold with corners!definition|)}
\label{ds7ex1}
\end{ex}

\subsection{`Standard model' d-manifolds with corners}
\label{ds72}
\I{d-manifold with corners!standard model|(}

In Examples \ref{ds4ex2} and \ref{ds4ex3} of \S\ref{ds42}, we
defined `standard model' d-manifolds $\bS_{V,E,s}$ and 1-morphisms
$\bS_{\smash{f,\hat f}}:\bS_{V,E,s}\ra \bS_{W,F,t}$. In \cite[\S
7.1--\S 7.2]{Joyc6} we show that all this extends to d-manifolds
with corners in a straightforward way.

\begin{ex} Let $V$ be a manifold with corners, $E\ra V$ a vector
bundle, and $s:V\ra E$ a smooth section of $E$. We will write down
an explicit principal d-manifold with corners~$\rS=(\bS,\bpS,\bs
i_\rS,\om_\rS)$.

Define a vector bundle $E_\pd\ra\pd V$ by $E_\pd=i_V^*(E)$, and a
section $s_\pd:\pd V\ra E_\pd$ by $s_\pd=i_V^*(s)$. Define d-spaces
$\bS=\bS_{V,E,s}$ and $\bpS=\bS_{\pd V,E_\pd,s_\pd}$ from the
triples $V,E,s$ and $\pd V,E_\pd,s_\pd$ exactly as in Example
\ref{ds4ex2}, although now $V,\pd V$ have corners. Define a
1-morphism $\bs i_\rS:\bpS\ra\bS$ in $\dSpa$ to be the `standard
model' 1-morphism $\bS_{i_V,\id_{E_\pd}}: \bS_{\pd
V,E_\pd,s_\pd}\ra\bS_{V,E,s}$ from Example~\ref{ds4ex3}.

Comparing the analogues of \eq{ds6eq1} for $i_V:\pd V\ra V$ and
\eq{ds6eq2} for $\bs i_\rS:\bpS\ra\bS$, we see that the conormal
bundle $\cN_\rS$ of $\bpS$ in $\bS$ is canonically isomorphic to the
lift to $\upS\subseteq\upV$ of the conormal bundle $\cN_V$ of $\pd
V$ in $V$. Define $\om_\rS$ to be the orientation on $\cN_\rS$
induced by the orientation on $\cN_V$ by outward-pointing normal
vectors to $\pd V$ in $V$. Then $\rS=(\bS,\bpS,\bs i_\rS,\om_\rS)$
is a d-space with corners. It is equivalent to $\rV\t_{\bs s,\rE,\bs
0}\rV$ in Definition \ref{ds7def1}, and so is a principal d-manifold
with corners. We call $\rS$ the {\it standard model\/} of $(V,E,s)$,
and write it~$\rS_{V,E,s}$.\G[SVEsb]{$\rS_{V,E,s}$}{`standard model'
d-manifold with corners}

There is a natural 1-isomorphism $\pd\rS_{V,E,s}\cong\rS_{\pd
V,E_\pd,s_\pd}$ in~$\dManc$.\I{d-manifold with corners!standard
model!boundary of}
\label{ds7ex2}
\end{ex}

\begin{ex}\I{d-manifold with corners!standard model!1-morphism|(}
Let $V,W$ be manifolds with corners, $E\ra V$, $F\ra W$ be vector
bundles, and $s:V\ra E$, $t:W\ra F$ be smooth sections. Then Example
\ref{ds7ex2} defines `standard model' principal d-manifolds with
corners $\rS_{V,E,s},\rS_{W,F,t}$, with underlying d-spaces
$\bS_{V,E,s},\bS_{W,F,t}$. Suppose $f:V\ra W$ is a smooth map, and
$\hat f:E\ra f^*(F)$ is a morphism of vector bundles on $V$
satisfying $\hat f\ci s= f^*(t)+O(s^2)$ in $C^\iy(f^*(F))$, where
$f^*(t)=t\ci f$, and $O(s^2)$ is as \S\ref{ds42}. Define a
1-morphism $\rS_{\smash{f,\hat f}}:\bS_{V,E,s}\ra\bS_{W,F,t}$ in
$\dSpa$ using $f,\hat f$ exactly as in Example \ref{ds4ex3}. Then
$\rS_{\smash{f,\hat f}}:\rS_{V,E,s}\ra\rS_{W,F,t}$ is a 1-morphism
in $\dManc$, which we call a `standard model'
1-morphism.\G[Sffb]{$\rS_{\smash{f,\hat f}}:\rS_{V,E,s}\ra
\rS_{W,F,t}$}{`standard model' 1-morphism in $\dManc$}

Suppose $\ti V\subseteq V$ is open, with inclusion $i_{\smash{\ti
V}}:\ti V\ra V$. Write $\ti E=E\vert_{\ti V}=i_{\ti V}^*(E)$ and
$\ti s=s\vert_{\ti V}$. Define $\bs i_{\smash{\ti
V,V}}=\rS_{i_{\smash{\ti V}},\id_{\ti E}}:\rS_{\smash{\ti V,\ti
E,\ti s}}\ra\rS_{V,E,s}$.\G[iVVb]{$\bs i_{\smash{\ti
V,V}}:\rS_{\smash{\ti V,\ti E,\ti s}}\ra\rS_{V,E,s}$}{inclusion of
open set in `standard model' d-manifold with corners} If
$s^{-1}(0)\subseteq\ti V$ then $\bs i_{\smash{\ti V,V}}$ is a
1-isomorphism, with inverse~$\bs i_{\ti V,V}^{-1}$.
\label{ds7ex3}
\end{ex}

In \cite[\S 7.2 \& \S 7.3]{Joyc6} we prove analogues of Theorems
\ref{ds4thm1} and~\ref{ds4thm2}:

\begin{thm} Let\/ $\rX$ be a d-manifold with corners, and\/
$x\in\rX$. Then there exists an open neighbourhood\/ $\rU$ of\/ $x$
in $\rX$ and an equivalence $\rU\simeq\rS_{V,E,s}$ in $\dManc$ for
some manifold with corners\/ $V,$ vector bundle $E\ra V$ and smooth
section $s:V\ra E$ which identifies $x\in\rU$ with a point\/ $v\in
S^k(V)\subseteq V,$ where $S^k(V)$ is as in\/ {\rm\S\ref{ds51},}
such that\/ $s(v)=\d s\vert_{S^k(V)}(v)=0$. Furthermore, $V,E,s$
and\/ $k$ are determined up to non-canonical isomorphism near $v$ by
$\rX$ near~$x$.
\label{ds7thm1}
\end{thm}

\begin{thm} Let\/ $V,W$ be manifolds with corners, $E\ra V,$ $F\ra
W$ be vector bundles, and\/ $s:V\ra E,$ $t:W\ra F$ be smooth
sections. Suppose $\bs g:\rS_{V,E,s}\ra\rS_{W,F,t}$ is a
$1$-morphism in $\dManc$. Then there exist an open neighbourhood\/
$\ti V$ of\/ $s^{-1}(0)$ in $V,$ a smooth map $f:\ti V\ra W,$ and a
morphism of vector bundles $\hat f:\ti E\ra f^*(F)$ with\/ $\hat
f\ci\ti s= f^*(t),$ where $\ti E=E\vert_{\ti V},$ $\ti s=s\vert_{\ti
V},$ such that\/ $\bs g=\rS_{\smash{f,\hat f}}\ci\smash{\bs
i_{\smash{\ti V,V}}^{-1}},$ using the notation of Examples\/
{\rm\ref{ds7ex2}} and\/~{\rm\ref{ds7ex3}}.\I{d-manifold with
corners!standard model|)}\I{d-manifold with corners!standard
model!1-morphism|)}
\label{ds7thm2}
\end{thm}

\subsection{Equivalences in $\dManc$, and gluing by equivalences}
\label{ds73}
\I{d-manifold with corners!equivalence|(}

In \cite[\S 7.4]{Joyc6} we study equivalences and gluing in
$\dManc$, as for $\dMan$ in \S\ref{ds44}. Here are the analogues of
Definition \ref{ds4def5} and Theorems
\ref{ds4thm3}--\ref{ds4thm5}.\I{d-manifold with corners!etale
1-morphism@\'etale 1-morphism|(}

\begin{dfn} Let $\bs f:\rX\ra\rY$ be a 1-morphism in $\dManc$. We
call $\bs f$ {\it \'etale\/} if it is a {\it local equivalence},
that is, if for each $x\in\rX$ there exist open
$x\in\rU\subseteq\rX$ and $\bs f(x)\in\rV\subseteq\rY$ such that
$\bs f(\rU)=\rV$ and $\bs f\vert_\rU:\rU\ra\rV$ is an equivalence.
\label{ds7def2}
\end{dfn}

\begin{thm} Suppose $\bs f:\rX\ra\rY$ is a $1$-morphism of
d-manifolds with corners. Then the following are equivalent:
\begin{itemize}
\setlength{\itemsep}{0pt}
\setlength{\parsep}{0pt}
\item[{\rm(i)}] $\bs f$ is \'etale;
\item[{\rm(ii)}] $\bs f$ is simple and flat, in the sense of\/
{\rm\S\ref{ds62},} and\/ $\Om_{\bs f}:\uf^* (T^*\rY)\ra T^*\rX$
is an equivalence in $\vqcoh(\uX);$ and
\item[{\rm(iii)}] $\bs f$ is simple and flat, and\/ \eq{ds4eq3} is
a split short exact sequence\I{split short exact
sequence}\I{abelian category!split short exact sequence}
in\/~$\qcoh(\uX)$.
\end{itemize}
If in addition $f:X\ra Y$ is a bijection, then $\bs f$ is an
equivalence in\/~$\dManc$.
\label{ds7thm3}
\end{thm}

\begin{thm} Let\/ $V,W$ be manifolds with corners, $E\ra V,$
$F\ra W$ be vector bundles, $s:V\ra E,$ $t:W\ra F$ be smooth
sections, $f:V\ra W$ be smooth, and\/ $\hat f:E\ra f^*(F)$ be a
morphism of vector bundles on $V$ with\/ $\hat f\ci s=
f^*(t)+O(s^2)$. Then Examples {\rm\ref{ds7ex2}} and\/
{\rm\ref{ds7ex3}} define principal d-manifolds with corners\/
$\rS_{V,E,s},\rS_{W,F,t}$ and a $1$-morphism\/ $\rS_{\smash{f,\hat
f}}:\rS_{V,E,s}\ra\rS_{W,F,t}$. This $\rS_{\smash{f,\hat f}}$ is
\'etale if and only if\/ $f$ is simple and flat near
$s^{-1}(0)\subseteq V,$ in the sense of\/ {\rm\S\ref{ds52},} and for
each\/ $v\in V$ with\/ $s(v)=0$ and\/ $w=f(v)\in W,$ equation
\eq{ds4eq4} is exact. Also $\rS_{\smash{f,\hat f}}$ is an
equivalence if and only if in addition\/
$f\vert_{s^{-1}(0)}:s^{-1}(0)\!\ra\! t^{-1}(0)$ is a bijection,
where $s^{-1}(0)\!=\!\{v\in V:s(v)\!=\!0\},$ $t^{-1}(0)\!=\!\{w\in
W:t(w)\!=\!0\}$.\I{d-manifold with corners!etale 1-morphism@\'etale
1-morphism|)}
\label{ds7thm4}
\end{thm}

\begin{thm} Suppose we are given the following data:\I{d-manifold
with corners!gluing by equivalences|(}
\begin{itemize}
\setlength{\itemsep}{0pt}
\setlength{\parsep}{0pt}
\item[{\rm(a)}] an integer $n;$
\item[{\rm(b)}] a Hausdorff, second countable topological space $X;$
\item[{\rm(c)}] an indexing set\/ $I,$ and a total order $<$ on $I;$
\item[{\rm(d)}] for each\/ $i$ in $I,$ a manifold with corners
$V_i,$ a vector bundle $E_i\ra V_i$ with\/ $\dim V_i-\rank
E_i=n,$ a smooth section $s_i:V_i\ra E_i,$ and a homeomorphism
$\psi_i:X_i\ra\hat X_i,$ where $X_i=\{v_i\in V_i:s_i(v_i)=0\}$
and\/ $\hat X_i\subseteq X$ is open; and
\item[{\rm(e)}] for all\/ $i<j$ in $I,$ an open submanifold\/
$V_{ij}\subseteq V_i,$ a simple, flat map $e_{ij}:V_{ij}\ra
V_j,$ and a morphism of vector bundles $\hat
e_{ij}:E_i\vert_{V_{ij}}\ra e_{ij}^*(E_j)$.
\end{itemize}
Let this data satisfy the conditions:
\begin{itemize}
\setlength{\itemsep}{0pt}
\setlength{\parsep}{0pt}
\item[{\rm(i)}] $X=\bigcup_{i\in I}\hat X_i;$
\item[{\rm(ii)}] if\/ $i<j$ in $I$ then $\hat e_{ij}\ci
s_i\vert_{V_{ij}}= e_{ij}^*(s_j),$ and\/ $\psi_i(X_i\cap
V_{ij})=\hat X_i\cap\hat X_j,$ and\/ $\psi_i\vert_{X_i\cap
V_{ij}}=\psi_j\ci e_{ij}\vert_{X_i\cap V_{ij}},$ and if\/
$v_i\in V_i$ with\/ $s_i(v_i)=0$ and\/ $v_j=e_{ij}(v_i)$ then
the following sequence of vector spaces is exact:
\begin{equation*}
\smash{\xymatrix@C=18pt{ 0 \ar[r] & T_{v_i}V_i
\ar[rrr]^(0.42){\d s_i(v_i)\op \,\d e_{ij}(v_i)} &&&
E_i\vert_{v_i}\!\op\! T_{v_j}V_j \ar[rrr]^(0.57){\hat
e_{ij}(v_i)\op\, -\d s_j(v_j)} &&& E_j\vert_{v_j} \ar[r] &
0;}}\!\!\!
\end{equation*}
\item[{\rm(iii)}] if\/ $i\!<\!j\!<\!k$ in $I$ then
$e_{ik}\vert_{V_{ik}\cap e_{ij}^{-1}(V_{jk})}= e_{jk}\ci
e_{ij}\vert_{V_{ik}\cap e_{ij}^{-1}(V_{jk})}+O(s_i^2)$ and\/
$\hat e_{ik}\vert_{V_{ik}\cap e_{ij}^{-1}(V_{jk})}=
e_{ij}\vert_{V_{ik}\cap e_{ij}^{-1}(V_{jk})}^*(\hat e_{jk})\ci
\hat e_{ij}\vert_{V_{ik}\cap e_{ij}^{-1}(V_{jk})}+O(s_i)$.
\end{itemize}

Then there exist a d-manifold with corners\/ $\rX$ with\/
$\vdim\rX=n$ and topological space $X,$ and a $1$-morphism
$\bs\psi_i:\rS_{V_i,E_i,s_i}\ra\rX$ in $\dManc$ with underlying
continuous map $\psi_i$ which is an equivalence with the open
d-submanifold\/ $\bs{\hat{\rm X}}_i\subseteq\rX$ corresponding to
$\hat X_i\subseteq X$ for all\/ $i\in I,$ such that for all\/ $i<j$
in $I$ there exists a $2$-morphism\/
$\eta_{ij}:\bs\psi_j\ci\rS_{e_{ij},\hat e_{ij}}\Ra\bs\psi_i\ci\bs
i_{V_{ij},V_i},$ where $\rS_{e_{ij},\hat
e_{ij}}:\rS_{V_{ij},E_i\vert_{V_{ij}},
s_i\vert_{V_{ij}}}\ra\rS_{V_j,E_j,s_j}$ and\/ $\bs i_{V_{ij},V_i}:
\rS_{V_{ij},E_i\vert_{V_{ij}}, s_i\vert_{V_{ij}}}\ra
\rS_{V_i,E_i,s_i}$ are as in Example\/ {\rm\ref{ds7ex3}}. This\/
$\rX$ is unique up to equivalence in~$\dManc$.

Suppose also that\/ $Y$ is a manifold with corners, and\/
$g_i:V_i\ra Y$ are smooth maps for all\/ $i\in I,$ and\/ $g_j\ci
e_{ij}=g_i\vert_{V_{ij}}+O(s_i^2)$ for all\/ $i<j$ in $I$. Then
there exist a $1$-morphism $\bs h:\rX\ra\rY$ unique up to
$2$-isomorphism, where $\rY=F_\Manc^\dManc(Y)=\rS_{Y,0,0},$ and\/
$2$-morphisms $\ze_i:\bs h\ci\bs\psi_i\Ra\rS_{g_i,0}$ for all\/
$i\in I$. Here $\rS_{Y,0,0}$ is from Example\/ {\rm\ref{ds7ex2}}
with vector bundle $E$ and section $s$ both zero, and\/
$\rS_{g_i,0}:\rS_{V_i,E_i,s_i}\ra\rS_{Y,0,0}=\rY$ is from Example
{\rm\ref{ds7ex3}} with\/~$\hat g_i=0$.
\label{ds7thm5}
\end{thm}

We can use Theorem \ref{ds7thm5} as a tool to prove the existence of
d-manifold with corner structures on spaces coming from other areas
of geometry.\I{d-manifold with corners!equivalence|)}\I{d-manifold
with corners!gluing by equivalences|)}

\subsection{Submersions, immersions and embeddings}
\label{ds74}
\I{d-manifold with corners!w-submersion|(}\I{d-manifold with
corners!sw-submersion|(}\I{d-manifold with
corners!submersion|(}\I{d-manifold with
corners!s-submersion|(}\I{d-manifold with
corners!w-immersion|(}\I{d-manifold with
corners!sw-immersion|(}\I{d-manifold with
corners!sfw-immersion|(}\I{d-manifold with
corners!immersion|(}\I{d-manifold with
corners!s-immersion|(}\I{d-manifold with
corners!sf-immersion|(}\I{d-manifold with
corners!w-embedding|(}\I{d-manifold with
corners!sw-embedding|(}\I{d-manifold with
corners!sfw-embedding|(}\I{d-manifold with
corners!embedding|(}\I{d-manifold with
corners!s-embedding|(}\I{d-manifold with corners!sf-embedding|(}

In \S\ref{ds45} we defined two kinds of submersions (submersions and
w-submersions), immersions, and embeddings for d-manifolds. In
\S\ref{ds52} we defined two kinds of submersions (submersions and
s-submersions), and three kinds of immersions (immersions, s- and
sf-immersions), and embeddings for manifolds with corners. In
\cite[\S 7.5]{Joyc6}, we combine both alternatives for d-manifolds
with corners, giving four types of submersions, and six types of
immersions and embeddings.

\begin{dfn} Let $\bs f:\rX\ra\rY$ be a 1-morphism in $\dManc$. As in
\S\ref{ds43} and \S\ref{ds71}, $T^*\rX,\uf^*(T^*\rY)$ are virtual
vector bundles on $\uX$ of ranks $\vdim\rX,\vdim\rY$, and $\Om_{\bs
f}:\uf^*(T^*\rY)\ra T^*\rX$ is a 1-morphism in $\vvect(\uX)$. Also
we have 1-morphisms $C(\bs f),\hat C(\bs f):C(\rX)\ra C(\rY)$ in
$\dcManc\subset\dSpac$ as in \S\ref{ds65} and \S\ref{ds71}, so we
can form $\Om_{C(\bs f)}:\ul{\smash{C(f)\!}}\,^*(T^*C(\rY))\ra
T^*C(\rX)$ and $\Om_{\smash{\hat C(\bs f)}}:\ul{\smash{\hat
C(f)\!}}\,^*(T^*C(\rY))\ra T^*C(\rX)$. Then:
\begin{itemize}
\setlength{\itemsep}{0pt}
\setlength{\parsep}{0pt}
\item[(a)] We call $\bs f$ a {\it w-submersion\/} if $\bs f$ is
semisimple and flat and $\Om_{\bs f}$ is weakly injective. We
call $\bs f$ an {\it sw-submersion\/} if it is also simple.
\item[(b)] We call $\bs f$ a {\it submersion\/} if $\bs f$ is
semisimple and flat and $\Om_{C(\bs f)}$ is injective. We call
$\bs f$ an {\it s-submersion\/} if it is also simple.
\item[(c)] We call $\bs f$ a {\it w-immersion\/} if $\Om_{\bs f}$ is
weakly surjective. We call $\bs f$ an {\it sw-immersion}, or
{\it sfw-immersion}, if $\bs f$ is also simple, or simple and
flat.
\item[(d)] We call $\bs f$ an {\it immersion\/} if
$\Om_{\smash{\hat C(\bs f)}}$ is surjective. We call $\bs f$ an
{\it s-immersion\/} if $\bs f$ is also simple, and an {\it
sf-immersion\/} if $\bs f$ is also simple and flat.
\item[(e)] We call $\bs f$ a {\it w-embedding, sw-embedding,
sfw-embedding, embedding, s-embedding}, or {\it sf-embedding},
if $\bs f$ is a w-immersion, \ldots, sf-immersion, respectively,
and $f:X\ra f(X)$ is a homeomorphism, so $f$ is injective.
\end{itemize}
Here (weakly) injective and (weakly) surjective 1-morphisms in
$\vvect(\uX)$ are defined in \S\ref{ds45}.

Parts (c)--(e) enable us to define {\it d-submanifolds\/} $\rX$ of a
d-manifold with corners $\rY$. {\it Open d-submanifolds\/} are open
d-subspaces $\rX$ in $\rY$. For more general d-submanifolds, we call
$\bs f:\rX\ra\rY$ a {\it w-immersed, sw-immersed, sfw-immersed,
immersed, s-immersed, sf-immersed, w-embedded, sw-embedded,
sfw-embedded, embedded, s-embedded}, or {\it sf-embedded
d-submanifold\/} of $\rY$ if $\rX,\rY$ are d-manifolds with corners
and $\bs f$ is a w-immersion, \ldots, sf-embedding,
respectively.\I{d-manifold with corners!d-submanifold}
\label{ds7def3}
\end{dfn}

Here is the analogue of Theorem \ref{ds4thm6}, proved in~\cite[\S
7.5]{Joyc6}.

\begin{thm}{\bf(i)} Any equivalence of d-manifolds with corners
is a w-sub\-mer\-sion, submersion, \ldots, sf-embedding.
\smallskip

\noindent{\bf(ii)} If\/ $\bs f,\bs g:\rX\ra\rY$ are $2$-isomorphic
$1$-morphisms of d-manifolds with corners then $\bs f$ is a
w-submersion, \ldots, sf-embedding, if and only if\/ $\bs g$ is.
\smallskip

\noindent{\bf(iii)} Compositions of w-submersions, \ldots,
sf-embeddings are of the same kind.
\smallskip

\noindent{\bf(iv)} The conditions that a $1$-morphism\/ $\bs
f:\rX\ra\rY$ in $\dManc$ is any kind of submersion or immersion are
local in $\rX$ and\/ $\rY$. The conditions that\/ $\bs f$ is any
kind of embedding are local in $\rY,$ but not in\/~$\rX$.
\smallskip

\noindent{\bf(v)} Let\/ $\bs f:\rX\ra\rY$ be a submersion in
$\dManc$. Then\/ $\vdim\rX\ge\vdim\rY,$ and if\/ $\vdim\rX=\vdim\rY$
then\/ $\bs f$ is \'etale.
\smallskip

\noindent{\bf(vi)} Let\/ $\bs f:\rX\ra\rY$ be an immersion in
$\dManc$. Then $\vdim\rX\le\vdim\rY$. If\/ $\bs f$ is an s-immersion
and\/ $\vdim\rX=\vdim\rY$ then\/ $\bs f$ is \'etale.
\smallskip

\noindent{\bf(vii)} Let\/ $f:X\ra Y$ be a smooth map of manifolds
with corners, and\/ $\bs f=F_\Manc^\dManc(f)$. Then $\bs f$ is a
submersion, s-submersion, immersion, s-immersion, sf-immersion,
embedding, s-embedding, or sf-embedding, in $\dManc$ if and only
if\/ $f$ is a submersion, \ldots, an sf-embedding in $\Manc,$
respectively. Also $\bs f$ is a w-immersion, sw-immersion,
sfw-immersion, w-embedding, sw-embedding, or sfw-embedding in
$\dManc$ if and only if\/ $f$ is an immersion, \ldots, sf-embedding
in $\Manc,$ respectively.
\smallskip

\noindent{\bf(viii)} Let\/ $\bs f:\rX\ra\rY$ be a $1$-morphism in
$\dManc,$ with\/ $\rY$ a manifold. Then $\bs f$ is a w-submersion if
and only if it is semisimple and flat, and\/ $\bs f$ is an
sw-submersion if and only if it is simple and flat.
\smallskip

\noindent{\bf(ix)} Let\/ $\rX,\rY$ be d-manifolds with corners,
with\/ $\rY$ a manifold. Then $\bs\pi_\rX:\rX\t\rY\ra\rX$ is a
submersion, and\/ $\bs\pi_\rX$ is an s-submersion
if\/~$\pd\rY=\bs\es$.
\smallskip

\noindent{\bf(x)} Suppose\/ $\bs f:\rX\ra\rY$ is a submersion in
$\dManc,$ and\/ $x\in\rX$ with\/ $\bs f(x)=y\in\rY$. Then there
exist open d-submanifolds $x\in\rU\subseteq\rX$ and\/
$y\in\rV\subseteq\rY$ with\/ $\bs f(\rU)=\rV,$ a manifold with
corners\/ $\rZ,$ and an equivalence $\bs i:\rU\ra\rV\t\rZ,$ such
that\/ $\bs f\vert_\rU:\rU\ra\rV$ is $2$-isomorphic to
$\bs\pi_\rV\ci\bs i,$ where $\bs\pi_\rV:\rV\t\rZ\ra\rV$ is the
projection. If\/ $\bs f$ is an s-submersion then\/~$\pd\rZ=\bs\es$.
\smallskip

\noindent{\bf(xi)} Let\/ $\bs f:\rX\ra\rY$ be a submersion of
d-manifolds with corners, with\/ $\rY$ a manifold with corners.
Then\/ $\rX$ is a manifold with corners.
\label{ds7thm6}
\end{thm}

Parts (ix)-(x) are a d-manifold analogue of
Proposition~\ref{ds5prop2}.\I{d-manifold with
corners!w-submersion|)}\I{d-manifold with
corners!sw-submersion|)}\I{d-manifold with
corners!submersion|)}\I{d-manifold with
corners!s-submersion|)}\I{d-manifold with
corners!w-immersion|)}\I{d-manifold with
corners!sw-immersion|)}\I{d-manifold with
corners!sfw-immersion|)}\I{d-manifold with
corners!immersion|)}\I{d-manifold with
corners!s-immersion|)}\I{d-manifold with
corners!sf-immersion|)}\I{d-manifold with
corners!w-embedding|)}\I{d-manifold with
corners!sw-embedding|)}\I{d-manifold with
corners!sfw-embedding|)}\I{d-manifold with
corners!embedding|)}\I{d-manifold with
corners!s-embedding|)}\I{d-manifold with corners!sf-embedding|)}

\subsection{Bd-transversality and fibre products}
\label{ds75}
\I{d-manifold with corners!fibre products|(}\I{d-manifold with
corners!bd-transverse 1-morphisms|(}\I{d-manifold with
corners!cd-transverse
1-morphisms|(}\I{bd-transversality|(}\I{cd-transversality|(}

In \cite[\S 7.6]{Joyc6} we extend \S\ref{ds46} to the corners case.
Here are the analogues of Definition \ref{ds4def8} and Theorems
\ref{ds4thm7}--\ref{ds4thm10}:

\begin{dfn} Let $\rX,\rY,\rZ$ be d-manifolds with corners and
$\bs g:\rX\ra\rZ,$ $\bs h:\rY\ra\rZ$ be 1-morphisms. We call $\bs
g,\bs h$ {\it bd-transverse\/} if they are both b-transverse in
$\dSpac$ in the sense of Definition \ref{ds6def4}, and d-transverse
in the sense of Definition \ref{ds4def8}. We call $\bs g,\bs h$ {\it
cd-transverse\/} if they are both c-transverse in $\dSpac$ in the
sense of Definition \ref{ds6def4}, and d-transverse. As in
\S\ref{ds66}, c-transverse implies b-transverse, so cd-transverse
implies bd-transverse.
\label{ds7def4}
\end{dfn}

\begin{thm} Suppose\/ $\rX,\rY,\rZ$ are d-manifolds with corners
and\/ $\bs g:\rX\ra\rZ,$ $\bs h:\rY\ra\rZ$ are bd-transverse
$1$-morphisms, and let\/ $\rW=\rX\t_{\bs g,\rZ,\bs h}\rY$ be the
fibre product in $\dSpac,$ which exists by Theorem\/
{\rm\ref{ds6thm4}(a)} as $\bs g,\bs h$ are b-transverse. Then $\rW$
is a d-manifold with corners, with\I{d-manifold with corners!fibre
products!bd-transverse}\I{fibre product!of d-manifolds with
corners}\I{d-manifold with corners!virtual dimension}
\e
\vdim\rW=\vdim\rX+\vdim\rY-\vdim\rZ.
\label{ds7eq1}
\e
Hence, all bd-transverse fibre products exist in\/~$\dManc$.
\label{ds7thm7}
\end{thm}

\begin{thm} Suppose\/ $\bs g:\rX\ra\rZ$ and\/ $\bs h:\rY\ra\rZ$ are
$1$-morphisms in $\dManc$. The following are sufficient conditions
for $\bs g,\bs h$ to be cd-transverse, and hence bd-transverse, so
that\/ $\rW=\rX\t_{\bs g,\rZ,\bs h}\rY$ is a d-manifold with corners
of virtual dimension\/ {\rm\eq{ds7eq1}:}
\begin{itemize}
\setlength{\itemsep}{0pt}
\setlength{\parsep}{0pt}
\item[{\bf(a)}] $\rZ$ is a manifold without boundary, that is,
$\bZ\in\bMan;$ or
\item[{\bf(b)}] $\bs g$ or $\bs h$ is a
w-submersion.\I{d-manifold with
corners!w-submersion}\I{d-manifold with corners!bd-transverse
1-morphisms|)}\I{bd-transversality|)}\I{d-manifold with
corners!cd-transverse 1-morphisms|)}\I{cd-transversality|)}
\end{itemize}
\label{ds7thm8}
\end{thm}

\begin{thm} Let\/ $\rX,\rY,\rZ$ be d-manifolds with corners with\/
$\rY$ a manifold, and\/ $\bs g:\rX\ra\rZ,$ $\bs h:\rY\ra\rZ$ be
$1$-morphisms with\/ $\bs g$ a submersion. Then\/ $\rW=\rX\t_{\bs
g,\rZ,\bs h}\rY$ is a manifold,
with\/~$\dim\rW=\vdim\rX+\dim\rY-\vdim\rZ$.\I{d-manifold with
corners!is a manifold}
\label{ds7thm9}
\end{thm}

\begin{thm}{\bf(i)} Let\/ $\rX$ be a d-manifold with corners
and\/ $\bs g:\rX\ra\bs{[0,\iy)^k}\t\bR^{\bs{n-k}}$ a semisimple,
flat\/ $1$-morphism in $\dManc$. Then the fibre product\/
$\rW=\rX\t_{\bs g,\bs{[0,\iy)^k}\t\bR^{\bs{n-k}},\bs 0}\bs *$ exists
in $\dManc,$ and\/ $\bs\pi_\rX:\rW\ra\rX$ is an s-embedding. When\/
$k=0,$ any $1$-morphism $\bs g:\rX\ra\bR^{\bs n}$ is semisimple and
flat, and\/ $\bs\pi_\rX:\rW\ra\rX$ is an sf-embedding.\I{d-manifold
with corners!semisimple 1-morphism}\I{d-manifold with corners!flat
1-morphism}\I{d-manifold with corners!s-embedding|(}
\smallskip

\noindent{\bf(ii)} Suppose\/ $\bs f:\rX\ra\rY$ is an
s-immersion\I{d-manifold with corners!s-immersion} of d-manifolds
with corners, and\/ $x\in\rX$ with\/ $\bs f(x)=y\in\rY$. Then there
exist open d-submanifolds $x\in\rU\subseteq\rX$ and\/
$y\in\rV\subseteq\rY$ with\/ $\bs f(\rU)\subseteq\rV$ and a
semisimple, flat\/ $1$-morphism $\bs
g:\rV\ra\bs{[0,\iy)^k}\t\bR^{\bs{n-k}}$ with\/ $\bs g(y)=0,$ where
$n=\vdim\rY-\vdim\rX\ge 0$ and\/ $0\le k\le n,$ fitting into a
$2$-Cartesian\I{2-category!2-Cartesian square} square in~$\dManc:$
\begin{equation*}
\xymatrix@C=170pt@R=10pt{*+[r]{\rU} \ar[d]^{\bs f\vert_\rU}
\ar[r]_(0.3){\bs\pi} \drtwocell_{}\omit^{}\omit{^{}} & *+[l]{\bs{*}}
\ar[d]_{\bs 0} \\
*+[r]{\rV} \ar[r]^(0.3){\bs g} & *+[l]{\bs{[0,\iy)^k}\t\bR^{\bs{n-k}}.} }
\end{equation*}
If\/ $\bs f$ is an sf-immersion\I{d-manifold with
corners!sf-immersion} then $k=0$. If\/ $\bs f$ is an s- or
sf-embedding then we may take\/~$\rU=\bs f^{-1}(\rV)$.\I{d-manifold
with corners!s-embedding|)}\I{d-manifold with corners!sf-embedding}
\label{ds7thm10}
\end{thm}

For ordinary manifolds, a submanifold $X$ in $Y$ may be described
locally either as the image of an embedding $X\hookra Y$, or
equivalently as the zeroes of a submersion $Y\ra\R^n$, where $n=\dim
Y-\dim X$. Theorem \ref{ds7thm10} is an analogue of this for
d-manifolds with corners. It should be compared with Proposition
\ref{ds5prop3} for manifolds with corners.\I{d-manifold with
corners!fibre products|)}

\subsection{Embedding d-manifolds with corners into manifolds}
\label{ds76}
\I{d-manifold with corners!embedding!into manifolds|(}

Section \ref{ds47} discussed embeddings of d-manifolds $\bX$ into
manifolds $Y$. Our two major results were Theorem \ref{ds4thm12},
which gave necessary and sufficient conditions on $\bX$ for
existence of embeddings $\bs f:\bX\hookra\bR^{\bs n}$ for $n\gg 0$,
and Theorem \ref{ds4thm13}, which showed that if an embedding $\bs
f:\bX\hookra\bY$ exists with $\bX$ a d-manifold and
$\bY=F_\Man^\dMan(Y)$, then $\bX\simeq\bS_{V,E,s}$ for open
$V\subseteq Y$, so $\bX$ is a principal d-manifold.

In \cite[\S 7.7]{Joyc6} we generalize these to d-manifolds with
corners. As in \S\ref{ds74}, we have three kinds of embeddings in
$\dManc$, embeddings, s-embeddings and sf-embeddings. The analogue
of Theorem \ref{ds4thm12} naturally holds for embeddings:

\begin{thm} Let\/ $\rX$ be a d-manifold with corners. Then there
exist immersions and/or embeddings $\bs f:\rX\ra\bR^{\bs n}$ for
some $n\gg 0$ if and only if there is an upper bound for\/ $\dim
T^*_x\uX$ for all\/ $x\in\uX$. If there is such an upper bound, then
immersions $\bs f:\rX\ra\bR^{\bs n}$ exist provided\/ $n\ge 2\dim
T_x^*\uX$ for all\/ $x\in\uX,$ and embeddings $\bs f:\rX\ra\bR^{\bs
n}$ exist provided\/ $n\ge 2\dim T_x^*\uX+1$ for all\/ $x\in\uX$.
For embeddings we may also choose $\bs f$ with\/ $f(X)$ closed
in\/~$\R^n$.\I{d-manifold with corners!immersion}\I{d-manifold with
corners!embedding}
\label{ds7thm11}
\end{thm}

Example \ref{ds4ex5} shows the hypotheses of Theorem \ref{ds7thm11}
need not hold, so there exist d-manifolds with corners $\rX$ with no
embedding into $\R^n$, or into any manifold with corners. The
analogue of Theorem \ref{ds4thm13} holds for
sf-embeddings:\I{d-manifold with
corners!sf-embedding|(}\I{d-manifold with
corners!principal|(}\I{principal d-manifold with
corners|(}\I{d-manifold with corners!standard model|(}

\begin{thm} Let\/ $\rX$ be a d-manifold with corners, $Y$ a
manifold with corners, and\/ $\bs f:\rX\ra\rY$ an sf-embedding, in
the sense of Definition\/ {\rm\ref{ds7def3}}. Then there exist an
open subset\/ $V$ in $Y$ with\/ $\bs f(\rX)\subseteq\rV,$ a vector
bundle $E\ra V,$ and a smooth section\/ $s:V\ra E$ of\/ $E$ fitting
into a $2$-Cartesian\I{2-category!2-Cartesian square} diagram in
$\dManc,$ where $0:V\ra E$ is the zero section
and\/~$\rY,\rV,\rE,\bs s,\bs 0=F_\Manc^\dManc(Y,V,E,s,0)\!:$
\begin{equation*}
\xymatrix@C=60pt@R=10pt{ \rX \ar[r]_(0.25){\bs f} \ar[d]^{\bs f}
\drtwocell_{}\omit^{}\omit{^{}}
 & \rV \ar[d]_{\bs 0} \\ \rV \ar[r]^(0.7){\bs s} & \rE.}
\end{equation*}
Hence $\rX$ is equivalent to the `standard model'\/ $\rS_{V,E,s}$ of
Example\/ {\rm\ref{ds7ex2},} and is a principal d-manifold with
corners.
\label{ds7thm12}
\end{thm}

Note that, unlike the d-manifolds case in \S\ref{ds47}, we cannot
immediately combine Theorems \ref{ds7thm11} and \ref{ds7thm12}: we
have first to bridge the gap between embeddings and sf-embeddings.
For d-manifolds with boundary, we can do this.

\begin{thm} Let\/ $\rX$ be a d-manifold with boundary. Then there
exist sf-immersions and/or sf-embeddings $\bs f:\rX\ra\bs{[0,\iy)}
\t\bR^{\bs{n-1}}$ for some $n\gg 0$ if and only if\/ $\dim T^*_x\uX$
is bounded above for all\/ $x\in\uX$. Such an upper bound always
exists if\/ $\rX$ is compact. If there is such an upper bound, then
sf-immersions $\bs f:\rX\ra\bs{[0,\iy)} \t\bR^{\bs{n-1}}$ exist
provided\/ $n\ge 2\dim T_x^*\uX+1$ for all\/ $x\in\uX,$ and
sf-embeddings $\bs f:\rX\ra\bs{[0,\iy)} \t\bR^{\bs{n-1}}$ exist
provided\/ $n\ge 2\dim T_x^*\uX+2$ for all\/ $x\in\uX$. For
sf-embeddings we may also choose $\bs f$ with\/ $f(X)$ closed
in\/~$[0,\iy)\t\R^{n-1}$.
\label{ds7thm13}
\end{thm}

Combining Theorems \ref{ds7thm12} and \ref{ds7thm13} shows that a
d-manifold with boundary $\rX$ is principal if and only if $\dim
T^*_x\uX$ is bounded above.

Since (nice) d-manifolds with boundary can be embedded into
$[0,\iy)\t\R^{n-1}$ for $n\gg 0$, one might guess that (nice)
d-manifolds with corners can be embedded into $[0,\iy)^k\t\R^{n-k}$
for $n\gg k\gg 0$. However, this is not true even for manifolds with
corners, as the following example from \cite[\S 5.7]{Joyc6} shows:

\begin{ex} Consider the {\it teardrop\/} $T=\bigl\{(x,y)\in\R^2:
x\ge 0$, $y^2\le x^2-x^4\bigr\}$, shown in Figure \ref{ds7fig}. It
is a compact 2-manifold with corners.

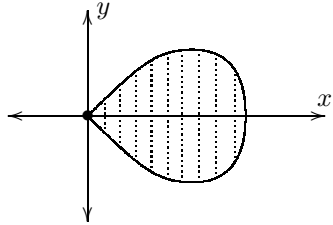
\begin{figure}[htb]
\begin{xy}
,(-1.5,0)*{}
,<6cm,-1.5cm>;<6.7cm,-1.5cm>:
,(3,.3)*{x}
,(-1.2,2)*{y}
,(-1.5,0)*{\bullet}
,(-1.5,0); (1.5,0) **\crv{(-.5,1)&(.1,1.4)&(1.5,1.2)}
?(.06)="a"
?(.12)="b"
?(.2)="c"
?(.29)="d"
?(.4)="e"
?(.5)="f"
?(.6)="g"
?(.7)="h"
?(.83)="i"
,(-1.5,0); (1.5,0) **\crv{(-.5,-1)&(.1,-1.4)&(1.5,-1.2)}
?(.06)="j"
?(.12)="k"
?(.2)="l"
?(.29)="m"
?(.4)="n"
?(.5)="o"
?(.6)="p"
?(.7)="q"
?(.83)="r"
,"a";"j"**@{.}
,"b";"k"**@{.}
,"c";"l"**@{.}
,"d";"m"**@{.}
,"e";"n"**@{.}
,"f";"o"**@{.}
,"g";"p"**@{.}
,"h";"q"**@{.}
,"i";"r"**@{.}
\ar (-1.5,0);(3,0)
\ar (-1.5,0);(-3,0)
\ar (-1.5,0);(-1.5,2)
\ar (-1.5,0);(-1.5,-2)
\end{xy}
\caption{The teardrop, a 2-manifold with corners.}
\label{ds7fig}
\end{figure}

Suppose that $f:T\ra[0,\iy)^k\t\R^{n-k}$ is an sf-embedding. As $f$
is simple and flat, it maps $S^j(T)\hookra
S^j\bigl([0,\iy)^k\t\R^{n-k}\bigr)$ for $j=0,1,2$, in the notation
of \S\ref{ds51}. The connected components of
$S^j\bigl([0,\iy)^k\t\R^{n-k}\bigr)$ correspond to subsets
$I\subseteq \{1,\ldots,k\}$ with $\md{I}=j$, with the component
corresponding to $I$ given by the equations $x_i=0$ for $i\in I$ and
$x_a>0$ for $a\in\{1,\ldots,k\}\sm I$. As $(0,0)\in S^2(T)$, we see
that $f(0,0)$ lies in the component of
$S^2\bigl([0,\iy)^k\t\R^{n-k}\bigr)$ given by $x_a=x_b=0$ for~$1\le
a<b\le k$.

Considering local models for $f$ near $(0,0)\in T$, we see that $f$
must map the two ends of $S^1(T)$ at $(0,0)$ into different
connected components $x_a=0$ and $x_b=0$ of
$S^1\bigl([0,\iy)^k\t\R^{n-k}\bigr)$. However, $S^1(T)\cong(0,1)$ is
connected, so $f$ maps $S^1(T)$ into a single connected component of
$S^1\bigl([0,\iy)^k\t\R^{n-k}\bigr)$, a contradiction. Hence there
do not exist sf-embeddings $f:T\ra [0,\iy)^k\t\R^{n-k}$ for
any~$n,k$.
\label{ds7ex4}
\end{ex}

Here are necessary and sufficient conditions for existence of
sf-embeddings from a d-manifold with corners $\rX$ into a manifold
with corners~$Y$.

\begin{thm} Let\/ $\rX$ be a d-manifold with corners. Then there
exist a manifold with corners $Y$ and an sf-embedding $\bs
f:\rX\ra\rY,$ where $\rY=F_\Manc^\dManc(Y),$ if and only if\/ $\dim
T^*_x\uX+\md{\ui_\rX^{-1}(x)}$ is bounded above for all\/ $x\in\uX$.
If such an upper bound exists, then we may take $Y$ to be an
embedded\/ $n$-dimensional submanifold of\/ $\R^n$ for any $n$
with\/ $n\ge 2\bigl(\dim T^*_x\uX+\md{\ui_\rX^{-1}(x)}\bigr)+1$ for
all\/~$x\in\uX$.

Such an upper bound always exists if\/ $\rX$ is compact. Thus, every
compact d-manifold with corners admits an sf-embedding into a
manifold with corners.
\label{ds7thm14}
\end{thm}

The idea of the proof of Theorem \ref{ds7thm14} is that we first
choose an embedding $\bs g:\rX\ra\bR^{\bs n}$ using Theorem
\ref{ds7thm11}, and then show that we can choose a submanifold
$Y\subseteq\R^n$ which is the set of points in an open neighbourhood
$U$ of $g(X)$ in $\R^n$ satisfying local transverse inequalities of
the form $c_i(x)\ge 0$ for $i=1,\ldots,k$, where $c_i:U\ra\R$ are
local smooth functions which lift under $\bs g$ to local boundary
defining functions for~$\pd\rX$.\I{d-manifold with
corners!sf-embedding|)}

Combining Theorems \ref{ds7thm12} and \ref{ds7thm14} yields:

\begin{cor} Let\/ $\rX$ be a d-manifold with corners. Then $\rX$ is
principal, that is, $\rX$ is equivalent in $\dManc$ to some
$\rS_{V,E,s}$ in Example\/ {\rm\ref{ds7ex2},} if and only if\/ $\dim
T^*_x\uX$ and\/ $\md{\ui_\rX^{-1}(x)}$ are bounded above for all\/
$x\in\uX$. This holds automatically if\/ $\rX$ is
compact.\I{d-manifold with corners!embedding!into
manifolds|)}\I{d-manifold with corners!principal|)}\I{principal
d-manifold with corners|)}\I{d-manifold with corners!standard
model|)}
\label{ds7cor}
\end{cor}

\subsection{Orientations}
\label{ds77}
\I{d-manifold with corners!orientations|(}

In \cite[\S 7.8]{Joyc6} we study orientations on d-manifolds with
corners, following the d-manifold case in \S\ref{ds48}. Here is the
analogue of Definition~\ref{ds4def9}:

\begin{dfn} Let $\rX$ be a d-manifold with corners. Then the
virtual cotangent bundle\I{virtual cotangent bundle}
$T^*\rX=(\EX,\FX,\phi_X)$ is a virtual vector bundle on $\uX$, so
Theorem \ref{ds4thm14} gives a line bundle $\cL_{T^*\rX}$ on $\uX$.
We call $\cL_{T^*\rX}$\G[LT*Xb]{$\cL_{T^*\rX}$}{orientation line
bundle of a d-manifold with corners $\rX$} the {\it orientation line
bundle\/}\I{d-manifold with corners!orientation line
bundle}\I{orientation line bundle} of~$\rX$.

An {\it orientation\/} $\om$ on $\rX$ is an orientation on
$\cL_{T^*\rX}$, in the sense of Definition \ref{ds4def9}. An {\it
oriented d-manifold with corners\/} is a pair $(\rX,\om)$ where
$\rX$ is a d-manifold with corners and $\om$ an orientation on
$\rX$. Usually we refer to $\rX$ as an oriented d-manifold, leaving
$\om$ implicit. We also write $-\rX$ for $\rX$ with the opposite
orientation, that is, $\rX$ is short for $(\rX,\om)$ and $-\rX$
short for~$(\rX,-\om)$.
\label{ds7def5}
\end{dfn}

Example \ref{ds4ex6}, Theorem \ref{ds4thm15} and Proposition
\ref{ds4prop4} now extend to d-manifolds with corners without
change. We can also orient boundaries of oriented d-manifolds with
corners. Theorem \ref{ds7thm15} is the main reason for including the
data $\om_\rX$ in a d-manifold with corners~$\rX=(\bX,\bpX,\bs
i_\rX,\om_\rX)$.\I{d-manifold with corners!boundary!orientation on}

\begin{thm} Let\/ $\rX$ be a d-manifold with corners. Then $\pd\rX$
is also a d-manifold with corners, so we have orientation line
bundles\I{orientation line bundle} $\cL_{T^*\rX}$ on $\uX$ and\/
$\cL_{T^*(\pd\rX)}$ on $\upX$. There is a canonical isomorphism of
line bundles on\/~$\upX\!:$
\e
\Psi:\cL_{T^*(\pd\rX)}\longra\ui_\rX^*(\cL_{T^*\rX})\ot\cN^*_\rX,
\label{ds7eq2}
\e
where $\cN_\rX$ is the conormal bundle of $\bpX$ in $\bX$
from\/~{\rm\S\ref{ds61}}.

Now\/ $\cN_\rX$ comes with an orientation $\om_\rX$ in
$\rX=(\bX,\bpX,\bs i_\rX,\om_\rX)$. Hence, if\/ $\rX$ is an oriented
d-manifold with corners, then $\pd\rX$ also has a natural
orientation, by combining the orientations on $\cL_{T^*\rX}$ and\/
$\cN^*_\rX$ to get an orientation on $\cL_{T^*(\pd\rX)}$
using\/~\eq{ds7eq2}.
\label{ds7thm15}
\end{thm}

As for Proposition \ref{ds4prop4}, natural equivalences of
d-manifolds with corners generally extend to natural equivalences of
oriented d-manifolds with corners, with some sign depending on the
orientation conventions.\I{orientation convention} Here are two such
results, which include signs in Theorem \ref{ds6thm1}(b) and
Proposition~\ref{ds6prop2}.\I{d-manifold with corners!fibre
products!orientations on|(}\I{d-manifold with corners!cd-transverse
1-morphisms|(}

\begin{prop} Suppose\/ $\rX,\rY$ are oriented d-manifolds with
corners, and\/ $\bs f:\rX\ra\rY$ is a semisimple, flat\/
$1$-morphism. Then the following holds in oriented d-manifolds with
corners, with fibre products cd-transverse:
\begin{equation*}
\pd_-^{\bs f}\rX\simeq\pd\rY\t_{\bs i_\rY,\rY,\bs f}\rX
\simeq(-1)^{\vdim\rX+\vdim\rY}\rX\t_{\bs f,\rY,\bs i_\rY}\pd\rY.
\end{equation*}
If\/ $\bs f$ is also simple then\/~$\pd_-^{\bs f}\rX=\pd\rX$.
\label{ds7prop1}
\end{prop}

\begin{prop} Let\/ $\bs g:\rX\ra\rZ$ and\/ $\bs h:\rY\ra\rZ$ be
$1$-morphisms of oriented d-manifolds with corners. Then the
following hold in oriented d-manifolds with corners, where all the
fibre products are cd-transverse, and so exist:
\begin{itemize}
\setlength{\itemsep}{0pt}
\setlength{\parsep}{0pt}
\item[{\bf(a)}] If\/ $\rZ$ is a manifold without boundary then
there is an equivalence
\begin{equation*}
\pd\bigl(\rX\t_{\bs g,\rZ,\bs h}\rY\bigr)\simeq\bigl(\pd
\rX\t_{\bs g\ci\bs i_\rX,\rZ,\bs h}\rY\bigr)\amalg
(-1)^{\vdim\rX+\dim\rZ}\bigl(\rX\t_{\bs g,\rZ,\bs h
\ci\bs i_\rY}\pd\rY\bigr).
\end{equation*}
\item[{\bf(b)}] If\/ $\bs g$ is a w-submersion then there is an
equivalence
\begin{equation*}
\pd\bigl(\rX\t_{\bs g,\rZ,\bs h}\rY\bigr)\simeq
\bigl(\pd^{\bs g}_+\rX\t_{\bs g_+,\rZ,\bs h}\rY\bigr)
\amalg(-1)^{\vdim\rX+\vdim\rZ}\bigl(\rX\t_{\bs g,\rZ, \bs
h\ci\bs i_\rY}\pd\rY\bigr).
\end{equation*}
\item[{\bf(c)}] If both\/ $\bs g$ and\/ $\bs h$ are
w-submersions then there is an equivalence
\begin{align*}
\pd\bigl(&\rX\t_{\bs g,\rZ,\bs h}\rY\bigr)\simeq
\bigl(\pd_+^{\bs g}\rX \t_{\bs g_+,\rZ,\bs h}\rY\bigr)\\
&\amalg(-1)^{\vdim\rX+\vdim\rZ}\bigl(\rX \t_{\bs g,\rZ,\bs
h_+}\pd_+^{\bs h}\rY\bigr) \amalg\bigl(\pd^{\bs g}_-\rX\t_{\bs
g_-,\pd\rZ,\bs h_-} \pd^{\bs h}_-\rY\bigr).
\end{align*}
\end{itemize}\I{d-manifold with corners|)}\I{d-manifold with
corners!orientations|)}\I{d-manifold with corners!fibre
products!orientations on|)}\I{d-manifold with corners!cd-transverse
1-morphisms|)}
\label{ds7prop2}
\end{prop}

\section{Deligne--Mumford $C^\iy$-stacks}
\label{ds8}
\I{Deligne--Mumford $C^\iy$-stack|(}\I{C-algebraic
geometry@$C^\iy$-algebraic geometry|(}

Next we discuss {\it Deligne--Mumford\/ $C^\iy$-stacks}, which are
related to $C^\iy$-schemes in the same way that Deligne--Mumford
stacks in algebraic geometry are related to schemes, and will be the
foundation of our theories of orbifolds, d-stacks and d-orbifolds.
$C^\iy$-stacks were introduced by the author in~\cite[\S 7--\S
11]{Joyc4}.

\subsection{$C^\iy$-stacks}
\label{ds81}
\I{C-stack@$C^\iy$-stack|(}
\I{C-stack@$C^\iy$-stack!Deligne--Mumford|see{Deligne-- \\ Mumford
$C^\iy$-stack}}

The next few definitions assume a lot of standard material from
stack theory, which is summarized in~\cite[\S 7]{Joyc4}.\I{stack|(}

\begin{dfn} Define a {\it Grothendieck topology\/}\I{Grothendieck
topology} $\cal J$ on the category $\CSch$ of $C^\iy$-schemes to
have coverings $\{\ui_a:\uU_a\ra\uU\}_{a\in A}$ where $V_a=i_a(U_a)$
is open in $U$ with $\ui_a:\uU_a\ra(V_a,\O_U\vert_{V_a})$ an
isomorphism for all $a\in A$, and $U=\bigcup_{a\in A}V_a$. Up to
isomorphisms of the $\uU_a$, the coverings
$\{\ui_a:\uU_a\ra\uU\}_{a\in A}$ of $\uU$ correspond exactly to open
covers $\{V_a:a\in A\}$ of $U$. Then $(\CSch,{\cal J})$ is a {\it
site}.

The stacks on $(\CSch,{\cal J})$ form a 2-category\I{2-category}
$\Sta_{(\CSch,{\cal J})}$, with all 2-morphisms invertible. As the
site $(\CSch,{\cal J})$ is subcanonical, there is a natural, fully
faithful\I{functor!full}\I{functor!faithful} functor $\CSch\ra
\Sta_{(\CSch,{\cal J})}$, defined explicitly below, which we write
as $\uX\mapsto\ul{\bar X\!}\,$ on objects and $\uf\mapsto\ul{\bar
f\!}\,$ on morphisms. A $C^\iy$-{\it stack\/} is a stack $\cX$ on
$(\CSch,{\cal J})$ such that the diagonal 1-morphism
$\De_\cX:\cX\ra\cX\t\cX$ is representable, and there exists a
surjective 1-morphism $\Pi:\bar\uU\ra\cX$ called an {\it
atlas\/}\I{atlas} for some $C^\iy$-scheme $\uU$. Write
$\CSta$\G[CSta]{$\CSta$}{2-category of $C^\iy$-stacks} for the full
2-subcategory of $C^\iy$-stacks in $\Sta_{(\CSch,{\cal J})}$. The
functor $\CSch\ra \Sta_{(\CSch,{\cal J})}$ above maps into $\CSta$,
so we also write it as~$F_\CSch^\CSta:\CSch\ra\CSta$.
\G[WXYZb]{$\cW,\cX,\cY,\cZ,\ldots$}{Deligne--Mumford $C^\iy$-stacks,
including orbifolds}

Formally, a $C^\iy$-stack is a category $\cX$ with a functor
$p_\cX:\cX\ra\CSch$, where $\cX,p_\cX$ must satisfy many complicated
conditions, including sheaf-like conditions for all open covers in
$\CSch$. A 1-{\it
morphism\/}\I{C-stack@$C^\iy$-stack!1-morphism}\I{2-category!1-morphism}
$f:\cX\ra\cY$ of $C^\iy$-stacks is a functor\I{category!functor}
$f:\cX\ra\cY$ with $p_\cY\ci f=p_\cX:\cX\ra\CSch$. Given 1-morphisms
$f,g:\cX\ra\cY$, a 2-{\it
morphism\/}\I{C-stack@$C^\iy$-stack!2-morphism}\I{2-category!2-morphism}
$\eta:f\Ra g$ is an isomorphism of
functors\I{category!functor!natural isomorphism} $\eta:f\Ra g$
with~$\id_{\smash{p_\cY}}
*\eta=\id_{\smash{p_\cX}}:p_\cY\ci f\Ra p_\cY\ci g$.

If $\uX$ is a $C^\iy$-scheme, the corresponding $C^\iy$-stack
$\ul{\bar X\!}\,=F_\CSch^\CSta(\uX)$\G[X]{$\bar\uX$}{$C^\iy$-stack
associated to a $C^\iy$-scheme $\uX$} is the category with objects
$(\uU,\uu)$ for $\uu:\uU\ra\uX$ a morphism in $\CSch$, and morphisms
$\uh:(\uU,\uu)\ra(\uV,\uv)$ for $\uh:\uU\ra\uV$ a morphism in
$\CSch$ with $\uv\ci\uh=\uu$. The functor $p_{\smash{\ul{\bar
X\!}\,}}: \ul{\bar X\!}\,\ra\CSch$ maps $p_{\smash{\ul{\bar
X\!}\,}}:(\uU,\uu)\mapsto\uU$ and~$p_{\smash{\ul{\bar
X\!}\,}}:\uh\mapsto\uh$.

If $\uf:\uX\ra\uY$ is a morphism of $C^\iy$-schemes, the
corresponding 1-morphism $\ul{\bar f\!}\,=F_\CSch^\CSta(\uf):
\ul{\bar X\!}\,\ra\ul{\bar Y\!\!}\,\,$ maps $\ul{\bar
f\!}\,:(\uU,\uu)\mapsto (\uU,\uf\ci\uu)$ on objects $(\uU,\uu)$ and
$\ul{\bar f\!}\,:\uh\mapsto\uh$ on morphisms $\uh$ in $\ul{\bar
X\!}\,$. This defines a functor $\ul{\bar f\!}\,:\ul{\bar X\!}\,
\ra\ul{\bar Y\!\!}\,\,$ with $p_{\smash{\ul{\bar Y\!\!}\,\,}}\ci
\ul{\bar f\!}\,=p_{\smash{\ul{\bar X\!}\,}}:\cX\ra \CSch$, so
$\ul{\bar f\!}\,$ is a 1-morphism $\ul{\bar f\!}\,:\ul{\bar X\!}\,
\ra\ul{\bar Y\!\!}\,\,$ in~$\CSta$.\I{stack|)}
\label{ds8def1}
\end{dfn}

We define some classes of morphisms of $C^\iy$-schemes:

\begin{dfn} Let $\uf:\uX\ra\uY$ be a
morphism in $\CSch$. Then:
\begin{itemize}
\setlength{\itemsep}{0pt}
\setlength{\parsep}{0pt}
\item We call $\uf$ an {\it open
embedding\/}\I{C-scheme@$C^\iy$-scheme!open embedding} if it is
an isomorphism with an open $C^\iy$-subscheme of $\uY$.
\item We call $\uf$ {\it
\'etale\/}\I{C-scheme@$C^\iy$-scheme!etale morphism@\'etale
morphism} if it is a local isomorphism (in the Zariski
topology).
\item We call $\uf$ {\it
proper\/}\I{C-scheme@$C^\iy$-scheme!proper morphism} if $f:X\ra
Y$ is a proper map of topological spaces, that is, if
$S\subseteq Y$ is compact then $f^{-1}(S)\subseteq X$ is
compact.
\item We call $\uf$ {\it universally
closed\/}\I{C-scheme@$C^\iy$-scheme!universally closed morphism}
if whenever $\ug:\uW\ra\uY$ is a morphism then
$\pi_W:X\t_{f,Y,g}W\ra W$ is a closed map of topological spaces.
\end{itemize}
Each one is invariant under base change and local in the target in
$(\CSch,{\cal J})$. Thus, they are also defined for representable
1-morphisms of $C^\iy$-stacks.\I{C-stack@$C^\iy$-stack!open
embedding}\I{C-stack@$C^\iy$-stack!etale 1-morphism@\'etale
1-morphism}\I{C-stack@$C^\iy$-stack!proper
1-morphism}\I{C-stack@$C^\iy$-stack!universally closed 1-morphism}
\label{ds8def2}
\end{dfn}

\begin{dfn} Let $\cX$ be a $C^\iy$-stack. We say that $\cX$
is {\it separated\/}\I{C-stack@$C^\iy$-stack!separated} if the
diagonal 1-morphism $\De_\cX:\cX\ra\cX\t\cX$ is universally closed.
If $\cX\simeq\ul{\bar X\!}\,$ for some $C^\iy$-scheme $\uX$ then
$\cX$ is separated if and only if $\uX$ is separated (Hausdorff).
\label{ds8def3}
\end{dfn}

\begin{dfn} Let $\cX$ be a $C^\iy$-stack. A $C^\iy$-{\it
substack\/}\I{C-stack@$C^\iy$-stack!C-substack@$C^\iy$-substack}
$\cY$ in $\cX$ is a strictly full subcategory $\cY$ in $\cX$ such
that $p_\cY:=p_\cX\vert_\cY:\cY\ra\CSch$ is also a $C^\iy$-stack. It
has a natural inclusion 1-morphism $i_\cY:\cY\hookra\cX$. We call
$\cY$ an {\it open\/ $C^\iy$-substack\/}
\I{C-stack@$C^\iy$-stack!C-substack@$C^\iy$-substack!open}of $\cX$
if $i_\cY$ is a representable open embedding. An {\it open
cover\/}\I{C-stack@$C^\iy$-stack!open cover} $\{\cY_a:a\in A\}$ of
$\cX$ is a family of open $C^\iy$-substacks $\cY_a$ in $\cX$ with
$\coprod_{a\in A}i_{\cY_a}:\coprod_{a\in A}\cY_a\ra\cX$ surjective.
\label{ds8def4}
\end{dfn}

\subsection{Topological spaces of $C^\iy$-stacks}
\label{ds82}
\I{C-stack@$C^\iy$-stack!underlying topological space $\cX_\top$|(}

By \cite[\S 8.4]{Joyc4}, a $C^\iy$-stack $\cX$ has an underlying
topological space~$\cX_\top$.

\begin{dfn} Let $\cX$ be a $C^\iy$-stack. Write $\ul{*}$ for the
point $\Spec\R$ in $\CSch$, and $\bar{\ul{*}}$ for the associated
point in $\CSta$. Define $\cX_\top$ to be the set of 2-isomorphism
classes $[x]$ of 1-morphisms $x:\bar{\ul{*}}\ra\cX$. If
$\cU\subseteq\cX$ is an open $C^\iy$-substack then any 1-morphism
$x:\bar{\ul{*}}\ra\cU$ is also a 1-morphism $x:\bar{\ul{*}}\ra\cX$,
and $\cU_\top$ is a subset of $\cX_\top$. Define
$\cT_{\cX_\top}=\bigl\{\cU_\top:\cU\subseteq\cX$ is an open
$C^\iy$-substack in $\cX\bigr\}$. Then $\cT_{\cX_\top}$ is a set of
subsets of $\cX_\top$ which is a topology on $\cX_\top$, so
$(\cX_\top,{\cal T}_{\cX_\top})$ is a topological space, which we
call the {\it underlying topological space\/} of $\cX$, and usually
write as $\cX_\top$. If $\uX=(X,\O_X)$ is a $C^\iy$-scheme, so that
$\ul{\bar X\!}\,$ is a $C^\iy$-stack, then $\ul{\bar X\!}\,_\top$ is
naturally homeomorphic to~$X$.\G[Xtop]{$\cX_\top$}{underlying
topological space of a $C^\iy$-stack $\cX$}

If $f:\cX\ra\cY$ is a 1-morphism of $C^\iy$-stacks then there is a
natural continuous map $f_\top:\cX_\top\ra\cY_\top$ defined by
$f_\top([x])=[f\ci x]$. If $f,g:\cX\ra\cY$ are 1-morphisms and
$\eta:f\Ra g$ is a 2-morphism then $f_\top=g_\top$. Mapping
$\cX\mapsto\cX_\top$, $f\mapsto f_\top$ and 2-morphisms to
identities defines a 2-functor $F_\CSta^\Top:\CSta\ra\Top$, where
the category of topological spaces $\Top$ is regarded as a
2-category with only identity 2-morphisms.
\label{ds8def5}
\end{dfn}

\begin{dfn} Let $\cX$ be a $C^\iy$-stack, and $[x]\in\cX_\top$.
Pick a representative $x$ for $[x]$, so that $x:\bar{\ul{*}}\ra\cX$
is a 1-morphism. Define the {\it orbifold
group\/}\I{C-stack@$C^\iy$-stack!orbifold group $\Iso_\cX([x])$} (or
{\it isotropy group}, or {\it stabilizer group\/}) $\Iso([x])$ or
$\Iso_\cX([x])$ of $[x]$ to be the group of 2-morphisms $\eta:x\Ra
x$. It is independent of the choice of $x\in[x]$ up to isomorphism,
which is canonical up to conjugation in~$\Iso_\cX([x])$.

If $f:\cX\ra\cY$ is a 1-morphism of $C^\iy$-stacks and
$[x]\in\cX_\top$ with $f_\top([x])=[y]\in\cY_\top$, for $y=f\ci x$,
then we define a group morphism $f_*:\Iso_\cX([x])\ra \Iso_\cY([y])$
by $f_*(\eta)=\id_f*\eta$. It is independent of choices of
$x\in[x]$, $y\in[y]$ up to conjugation
in~$\Iso_\cX([x]),\Iso_\cY([y])$.\I{C-stack@$C^\iy$-stack!underlying
topological space $\cX_\top$|)}
\label{ds8def6}
\end{dfn}

\subsection{Strongly representable 1-morphisms}
\label{ds83}
\I{C-stack@$C^\iy$-stack!strongly representable 1-morphism|(}

{\it Strongly representable\/} 1-morphisms, discussed in \cite[\S
8.6]{Joyc4}, will be important in the definitions of orbifolds,
d-stacks, and d-orbifolds with corners.

\begin{dfn} Let $\cY,\cZ$ be $C^\iy$-stacks, and $g:\cY\ra\cZ$ a
1-morphism. Then $\cY,\cZ$ are categories with functors
$p_\cY:\cY\ra\CSch$, $p_\cZ:\cZ\ra\CSch$, and $g:\cY\ra\cZ$ is a
functor with~$p_\cZ\ci g=p_\cY$.

We call $g$ {\it strongly representable\/} if whenever $A\in\cY$
with $p_\cY(A)=\uU\in\CSch$, so that $B=g(A)\in\cZ$ with
$p_\cZ(B)=\uU$, and $b:B\ra B'$ is an isomorphism in $\cZ$ with
$p_\cZ(B')=\uU$ and $p_\cZ(b)=\uid_\uU$, then there exist a unique
object $A'$ and isomorphism $a:A\ra A'$ in $\cY$ with $g(A')=B'$
and~$g(a)=b$.
\label{ds8def7}
\end{dfn}

Here are two important properties of strongly representable
1-morphisms. The first says that we may replace a representable
1-morphism $g:\cY\ra\cZ$ with a strongly representable 1-morphism
$g':\cY'\ra\cZ$ with $\cY'\simeq\cY$.

\begin{prop}{\bf(a)} Let\/ $g:\cY\ra\cZ$ be a strongly
representable $1$-morphism of\/ $C^\iy$-stacks. Then $g$ is
representable.
\smallskip

\noindent{\bf(b)} Suppose $g:\cY\ra\cZ$ is a representable
$1$-morphism of\/ $C^\iy$-stacks. Then there exist a\/ $C^\iy$-stack
$\cY',$ an equivalence $i:\cY\ra\cY',$ and a strongly representable
$1$-morphism $g':\cY'\ra\cZ$ with\/ $g=g'\ci i$. Also $\cY'$ is
unique up to canonical\/ $1$-isomorphism in\/~$\CSta$.
\label{ds8prop1}
\end{prop}

The second says that for some 2-commutative
diagrams\I{2-category!2-commutative diagram} involving strongly
representable morphisms, we can require the diagrams to commute {\it
up to equality}, not just up to 2-isomorphism.

\begin{prop} Suppose $\cX,\cY,\cZ$ are\/ $C^\iy$-stacks,
$f:\cX\ra\cY,$ $g:\cY\ra\cZ,$ $h:\cX\ra\cZ$ are $1$-morphisms with\/
$g$ strongly representable, and\/ $\eta:g\ci f\Ra h$ is a
$2$-morphism in $\CSta$. Then as in the diagram below there exist a
$1$-morphism $f':\cX\ra\cY$ with\/ $g\ci f'=h,$ and a $2$-morphism
$\ze:f\Ra f'$ with\/ $\id_g*\ze=\eta,$ and\/ $f',\ze$ are unique
under these conditions.
\begin{equation*}
\xymatrix@C=80pt@R=15pt{ & \cY \ar[dr]^g
\ar@{}@<-1ex>[d]^(0.55){\eta\ts\Downarrow} \\
\cX \ar@/^.8pc/[ur]^(0.4){f'} \ar@/_.5pc/[ur]_(0.8){f}
\ar@{}@<-1.3ex>[ur]^{\ze\ts\Uparrow} \ar@<-.3ex>[rr]^(0.7)h && \cZ.}
\end{equation*}
\label{ds8prop2}
\end{prop}

We will use strongly representable 1-morphisms to define orbifolds,
d-stacks, and d-orbifolds with corners so that boundaries\I{orbifold
with corners!boundary}\I{d-stack with corners!boundary}\I{d-orbifold
with corners!boundary} behave in a strictly functorial rather than
weakly functorial way, as for d-spaces with corners in Remark
\ref{ds6rem2}. Here is an explicit construction of fibre products
$\cX\t_{g,\cZ,h}\cY$\I{C-stack@$C^\iy$-stack!fibre products} in
$\CSta$ when $g$ is strongly representable, yielding a strictly
commutative 2-Cartesian square.

\begin{prop} Let\/ $g:\cX\ra\cZ$ and\/ $h:\cY\ra\cZ$ be
$1$-morphisms of\/ $C^\iy$-stacks with\/ $g$ strongly representable.
Define a category $\cW$ to have objects pairs $(A,B)$ for\/
$A\in\cX,$ $B\in\cY$ with\/ $g(A)=h(B)$ in $\cZ,$ so that\/
$p_\cX(A)=p_\cY(B)$ in $\CSch,$ and morphisms pairs
$(a,b):(A,B)\ra(A',B')$ with\/ $a:A\ra A',$ $b:B\ra B'$ morphisms in
$\cX,\cY$ with\/ $p_\cX(a)=p_\cY(b)$ in~$\CSch$.

Define functors $p_\cW:\cW\ra\CSch,$ $e:\cW\ra\cX,$ $f:\cW\ra\cY$ by
$p_\cW:(A,B)\mapsto p_\cX(A)=p_\cY(B),$ $e:(A,B)\mapsto A,$
$f:(A,B)\mapsto B$ on objects and\/ $p_\cW:(a,b)\mapsto
p_\cX(a)=p_\cY(b),$ $e:(a,b)\mapsto a,$ $f:(a,b)\mapsto b$ on
morphisms. Then $\cW$ is a\/ $C^\iy$-stack and\/ $e:\cW\ra\cX,$
$f:\cW\ra\cY$ are $1$-morphisms, with\/ $f$ strongly representable,
and\/ $g\ci e=h\ci f$. Furthermore, the following diagram in $\CSta$
is $2$-Cartesian:\I{2-category!2-Cartesian square}
\begin{equation*}
\xymatrix@C=140pt@R=10pt{ *+[r]{\cW} \ar[r]_(0.25){f} \ar[d]^{e}
\drtwocell_{}\omit^{}\omit{^{\id_{g\ci e}\,\,\,\,\,\,\,\,\,\,{}}} &
*+[l]{\cY} \ar[d]_{h} \\
*+[r]{\cX} \ar[r]^(0.7){g} & *+[l]{\cZ.} }
\end{equation*}
If also\/ $h$ is strongly representable, then $e$ is strongly
representable.\I{C-stack@$C^\iy$-stack!strongly representable
1-morphism|)}
\label{ds8prop3}
\end{prop}

\subsection{Quotient $C^\iy$-stacks}
\label{ds84}
\I{C-stack@$C^\iy$-stack!quotients $[\protect\uX/G]$|(}\I{quotient
C-stack@quotient $C^\iy$-stack|(}

An important class of examples of $C^\iy$-stacks $\cX$ are {\it
quotient\/ $C^\iy$-stacks\/} $[\uX/G]$, for $\uX$ a $C^\iy$-scheme
acted on by a finite group $G$. The next three examples define
quotient $C^\iy$-stacks $[\uX/G]$, quotient 1-morphisms
$[\uf,\rho]:[\uX/G]\ra[\uY/H]$, and quotient 2-morphisms
$[\de]:[\uf,\rho]\Ra [\ug,\si]$.

In fact Examples \ref{ds8ex1}--\ref{ds8ex3} are simplifications of
more complicated definitions given in \cite[\S 9.1]{Joyc4}. The
construction of \cite[\S 9.1]{Joyc4} gives equivalent $C^\iy$-stacks
$[\uX/G]$, but has the advantage of being strictly functorial, that
is, quotient 1-morphisms compose as
$[\ug,\si]\ci[\uf,\rho]=[\ug\ci\uf, \si\ci\rho]$, whereas in Example
\ref{ds8ex2} we only have a 2-isomorphism
$[\ug,\si]\ci[\uf,\rho]\cong[\ug\ci\uf,\si\ci\rho]$. We will
occasionally assume this strict functoriality below, for instance,
in Definition~\ref{ds11def7}.

\begin{ex} Let $\uX$ be a separated $C^\iy$-scheme, $G$ a finite
group, and $\ur:G\ra\Aut(\uX)$ an action of $G$ on $\uX$ by
isomorphisms. We will define the {\it quotient
$C^\iy$-stack\/}~$\cX=[\uX/G]$.

Define a category $\cX$ to have objects quintuples
$(\uT,\uU,\ut,\uu,\uv)$, where $\uT,\uU$ are $C^\iy$-schemes,
$\ut:G\ra\Aut(\uT)$ is a free action of $G$ on $\uT$ by
isomorphisms, $\uu:\uT\ra\uX$ is a morphism with $\uu\ci\ut(\ga)=
\ur(\ga)\ci\uu:\uT\ra\uX$ for all $\ga\in G$, and $\uv:\uT\ra\uU$ is
a morphism which makes $\uT$ into a principal $G$-bundle over $\uU$,
that is, $\uv$ is proper, \'etale and surjective, and its fibres are
$G$-orbits in $\uT$ under~$\ut$.

A morphism $(\ua,\ub):(\uT,\ab\uU,\ab\ut, \ab\uu,\ab\uv)
\ra(\uT',\uU',\ut',\uu',\uv')$ in $\cX$ is a pair of morphisms
$\ua:\uU\ra\uU'$ and $\ub:\uT\ra\uT'$ such that
$\ub\ci\ut(\ga)=\ut'(\ga)\ci\ub$ for $\ga\in G$, and
$\uu=\uu'\ci\ub$, and $\ua\ci\uv=\uv'\ci\ub$. Composition is
$(\uc,\ud)\ci(\ua,\ub)= (\uc\ci\ua,\ud\ci\ub)$, and identities
are~$\id_{(\uT,\ldots,\uv)}=(\uid_\uU, \uid_\uT)$.

This defines the category $\cX$. The functor $p_\cX:\cX\ra\CSch$
acts by $p_\cX:(\uT,\uU,\ut,\uu,\uv)\mapsto\uU$ on objects, and
$p_\cX:(\ua,\ub)\mapsto\ua$ on morphisms. Then $\cX$ is a
$C^\iy$-stack, which we write as~$[\uX/G]$.
\label{ds8ex1}
\end{ex}

\begin{ex} Let $\uX,\uY$ be separated $C^\iy$-schemes acted on by
finite groups $G,H$ with actions $\ur:G\ra\Aut(\uX)$,
$\us:H\ra\Aut(\uY)$, so that we have quotient $C^\iy$-stacks
$[\uX/G]$ and $[\uY/H]$ as in Example \ref{ds8ex1}. Suppose we have
morphisms $\uf:\uX\ra\uY$ of $C^\iy$-schemes and $\rho:G\ra H$ of
groups, with $\uf\ci \ur(\ga)=\us(\rho(\ga))\ci\uf$ for all $\ga\in
G$. Define a functor $[\uf,\rho]:[\uX/G]\ra[\uY/H]$ on objects in
$[\uX/G]$ by
\begin{equation*}
[\uf,\rho]:(\uT,\uU,\ut,\uu,\uv)
\longmapsto \bigl((\uT\t H)/G,\uU,
\ul{\ti t\kern -0.1em}\kern 0.1em,
\ul{\ti u\kern -0.1em}\kern 0.1em,
\ul{\ti v\kern -0.1em}\kern 0.1em\bigr).
\end{equation*}
Here for each $\de\in H$, write $L_\de,R_\de:H\ra H$ for left and
right multiplication by $\de$. Then to define $(\uT\t H)/G$, each
$\ga\in G$ acts by $\ur(\ga)\t R_{\rho(\ga)^{-1}}:\uT\t H\ra\uT\t
H$. For each $\de\in H$, the morphism $\ul{\ti t\kern -0.1em}\kern
0.1em(\de):(\uT\t H)/G\ra (\uT\t H)/G$ is induced by the morphism
$\uid_\uT\t L_\de:\uT\t H\ra\uT\t H$. The morphisms $\ul{\ti u\kern
-0.1em}\kern 0.1em:(\uT\t H)/G\ra\uY$ and $\ul{\ti v\kern
-0.1em}\kern 0.1em:(\uT\t H)/G\ra\uU$ are induced by
$\uf\ci\uu\ci\upi_\uT:\uT\t H\ra\uY$ and~$\uv\ci\upi_\uT:\uT\t
H\ra\uU$.

On morphisms $(\ua,\ub):(\uT,\uU,\ut,\uu,\uv) \ra(\uT',\uU',
\ut',\uu',\uv')$ in $[\uX/G]$, define $[\uf,\rho]$ to map
$(\ua,\ub)\mapsto(\ua,\ul{\ti b\kern -0.1em}\kern 0.1em)$, where
$\ul{\ti b\kern -0.1em}\kern 0.1em:(\uT\t H)/G\ra (\uT'\t H)/G$ is
induced by $\ub\t\id_H:\uT\t H\ra \uT'\t H$. Then $[\uf,\rho]:
[\uX/G]\ra[\uY/H]$ is a 1-morphism of $C^\iy$-stacks, which we call
a {\it quotient\/ $1$-morphism}.\I{C-stack@$C^\iy$-stack!quotients
$[\protect\uX/G]$!quotient 1-morphism}

If $\rho:G\ra H$ is injective, then $[\uf,\rho]: [\uX/G]\ra[\uY/H]$
is representable.
\label{ds8ex2}
\end{ex}

\begin{ex} Let $[\uf,\rho]:[\uX/G]\ra[\uY/H]$ and
$[\ug,\si]:[\uX/G]\ra[\uY/H]$ be quotient 1-morphisms, so that
$\uf,\ug:\uX\ra\uY$ and $\rho,\si:G\ra H$ are morphisms. Suppose
$\de\in H$ satisfies $\si(\ga)=\de\,\rho(\ga)\,\de^{-1}$ for all
$\ga\in G$, and $\ug=\us(\de)\ci\uf$.

For each object $(\uT,\uU,\ut,\uu,\uv)$ in $[\uX/G]$, define an
isomorphism in $[\uY/H]$:
\begin{align*}
[\de]\bigl((\uT,\uU,\ut,\uu,\uv)\bigr)\!=\!(\uid_\uU,\ui_\de)\!:\!
[\uf,\rho]\bigl((\uT,\uU,\ut,\uu,\uv)\bigr)\!=\!
\bigl((\uT\!\t\! H)/_{\ur\t R_{\rho^{-1}}} G,\uU,
\ul{\ti t\kern -0.1em}\kern 0.1em,
\ul{\ti u\kern -0.1em}\kern 0.1em,&
\ul{\ti v\kern -0.1em}\kern 0.1em\bigr)\\
\longra [\ug,\si]\bigl((\uT,\uU,\ut,\uu,\uv)\bigr)= \bigl((\uT\t
H)/_{\ur\t R_{\si^{-1}}} G,\uU, \ul{\dot t\kern -0.1em}\kern 0.1em,
\ul{\dot u\kern -0.1em}\kern 0.1em, \ul{\dot v\kern -0.1em}\kern
0.1em\bigr)&,
\end{align*}
where $\ui_\de:(\uT\!\t\! H)/_{\ur\t R_{\rho^{-1}}}G\!\ra\!(\uT\!\t
\! H)/_{\ur\t R_{\si^{-1}}}G$ is induced $\uid_\uT\!\t\!
R_{\de^{-1}}\!:\! \uT\t H\ra\uT\t H$. Then $[\de]:[\uf,\rho]\Ra
[\ug,\si]$ is a natural isomorphism of functors, and a 2-morphism of
$C^\iy$-stacks, which we call a {\it quotient\/
$2$-morphism}.\I{C-stack@$C^\iy$-stack!quotients
$[\protect\uX/G]$|)}\I{quotient C-stack@quotient
$C^\iy$-stack|)}\I{C-stack@$C^\iy$-stack!quotients
$[\protect\uX/G]$!quotient 2-morphism}
\label{ds8ex3}
\end{ex}

\subsection{Deligne--Mumford $C^\iy$-stacks}
\label{ds85}

Deligne--Mumford stacks in algebraic geometry are locally modelled
on quotient stacks $[X/G]$ for $X$ an affine scheme and $G$ a finite
group. This motivates:

\begin{dfn} A {\it Deligne--Mumford\/ $C^\iy$-stack\/} is a
$C^\iy$-stack $\cX$ which admits an open cover $\{\cY_a:a\in A\}$
with each $\cY_a$ equivalent to a quotient stack $[\uU_a/G_a]$ in
Example \ref{ds8ex1} for $\uU_a$ an affine $C^\iy$-scheme and $G_a$
a finite group. We call $\cX$ {\it locally
fair\/}\I{Deligne--Mumford $C^\iy$-stack!locally fair} if it has
such an open cover with each $\uU_a$ a fair affine $C^\iy$-scheme.

We call a Deligne--Mumford $C^\iy$-stack $\cX$ {\it second
countable, compact, locally compact,} or {\it paracompact}, if the
underlying topological space $\cX_\top$ from \S\ref{ds82} is second
countable, compact, locally compact, or paracompact, respectively.

Write $\DMCSta,\kern -.1em\DMCStalf,\kern -.1em\DMCStalfssc$%
\G[DMCSta]{$\DMCSta$}{2-category of Deligne--Mumford
$C^\iy$-stacks}\G[DMCStalf]{$\DMCStalf$}{2-category of locally fair
Deligne--Mumford
$C^\iy$-stacks}\G[DMCStalfssc]{$\DMCStalfssc$}{2-category of
separated, second countable, locally fair Deligne--Mumford
$C^\iy$-stacks} for the full 2-subcategories\I{2-category} of
Deligne--Mumford $C^\iy$-stacks, and locally fair Deligne--Mumford
$C^\iy$-stacks, and separated, second countable, locally fair
Deligne--Mumford $C^\iy$-stacks in $\CSta$,
respectively.\I{Deligne--Mumford $C^\iy$-stack!definition}
\label{ds8def8}
\end{dfn}

If $\cX$ is a Deligne--Mumford $C^\iy$-stack then $\Iso_\cX([x])$ is
finite for all $[x]$ in $\cX_\top$. If $f:\cX\ra\cY$ is a 1-morphism
of Deligne--Mumford $C^\iy$-stacks then $f$ is representable if and
only if the morphisms of orbifold groups $f_*:\Iso_\cX([x])\ra
\Iso_\cY([y])$ from Definition \ref{ds8def6} are injective for all
$[x]\in\cX_\top$ with $f_\top([x])=[x]\in\cY_\top$. From \cite[\S
8--\S 9]{Joyc4}, we have:

\begin{thm}{\bf(a)} All fibre products exist in the\/
$2$-category\/~$\CSta$.\I{C-stack@$C^\iy$-stack!fibre
products}\I{fibre product!of C-stacks@of
$C^\iy$-stacks}\I{2-category!fibre products in}
\smallskip

\noindent{\bf(b)} $\DMCSta,\DMCStalf$ and\/ $\DMCStalfssc$ are
closed under fibre products and under taking open\/
$C^\iy$-substacks in\/~$\CSta$.\I{Deligne--Mumford
$C^\iy$-stack!fibre products}
\label{ds8thm1}
\end{thm}

\begin{prop} Let\/ $\cX$ be a Deligne--Mumford\/ $C^\iy$-stack and\/
$[x]\in\cX_\top,$ so that\/ $\Iso_\cX([x])\cong
H$\I{C-stack@$C^\iy$-stack!orbifold group $\Iso_\cX([x])$} for some
finite group $H$. Then there exists an open $C^\iy$-substack\/ $\cU$
in $\cX$ with\/ $[x]\in\cU_\top\subseteq\cX_\top$ and an
equivalence\/ $\cU\simeq[\uY/H],$ where $\uY=(Y,\O_Y)$ is an affine
$C^\iy$-scheme with an action of\/ $H,$ and\/ $[x]\in\cU_\top\cong
Y/H$ corresponds to a fixed point\/ $y$ of\/ $H$ in\/~$Y$.
\label{ds8prop4}
\end{prop}

\begin{thm} Suppose $\cX$ is a Deligne--Mumford\/ $C^\iy$-stack
with\/ $\Iso_\cX([x])\cong\{1\}$ for all\/ $[x]\in\cX_\top$. Then\/
$\cX$ is equivalent to\/ $\ul{\bar X\!}\,$ for some\/
$C^\iy$-scheme\/~$\uX$.
\label{ds8thm2}
\end{thm}

In conventional algebraic geometry, a stack with all orbifold groups
trivial is (equivalent to) an {\it algebraic space},\I{algebraic
space} but may not be a scheme, so the category of algebraic spaces
is larger than the category of schemes. Here algebraic spaces are
spaces which are locally isomorphic to schemes in the \'etale
topology,\I{etale topology@\'etale topology} but not necessarily
locally isomorphic to schemes in the Zariski topology.\I{Zariski
topology}

In contrast, as Theorem \ref{ds8thm2} shows, in $C^\iy$-algebraic
geometry there is no difference between $C^\iy$-schemes and
$C^\iy$-algebraic spaces. This is because in $C^\iy$-geometry the
Zariski topology is already fine enough, as in Remark
\ref{ds2rem1}(iii), so we gain no extra generality by passing to the
\'etale topology.\I{C-stack@$C^\iy$-stack|)}

\subsection{Quasicoherent sheaves on $C^\iy$-stacks}
\label{ds86}
\I{Deligne--Mumford $C^\iy$-stack!quasicoherent sheaves on|(}

In \cite[\S 10]{Joyc4} we study sheaves on Deligne--Mumford
$C^\iy$-stacks.

\begin{dfn} Let $\cX$ be a Deligne--Mumford $C^\iy$-stack. Define
a category ${\cal C}_\cX$ to have objects pairs $(\uU,u)$ where
$\uU$ is a $C^\iy$-scheme and $u:\bar\uU\ra\cX$ is an \'etale
1-morphism, and morphisms $(\uf,\eta):(\uU,u)\ra(\uV,v)$ where
$\uf:\uU\ra\uV$ is an \'etale morphism of $C^\iy$-schemes, and
$\eta:u\Ra v\ci\ul{\bar f\!}\,$ is a 2-isomorphism. If
$(\uf,\eta):(\uU,u)\ra(\uV,v)$ and $(\ug,\ze):(\uV,v)\ra(\uW,w)$ are
morphisms in ${\cal C}_\cX$ then we define the composition
$(\ug,\ze)\ci(\uf,\eta)$ to be $(\ug\ci\uf,\th):(\uU,u)\ra(\uW,w)$,
where $\th$ is the composition of 2-morphisms across the diagram:
\begin{equation*}
\xymatrix@C=20pt@R=7pt{ \bar\uU \ar[dr]^(0.6){\ul{\bar f\!}\,}
\ar@/^/@<1ex>[drrr]_(0.4)u \ar[dd]_{\overline{\ug\ci\uf}}
\dduppertwocell_{}\omit^{}\omit{<-2.5>^{\id}} & \\
& \bar\uV \ar[rr]^(0.5){v} \ar[dl]^(0.45){\bar\ug} &
\ultwocell_{}\omit^{}\omit{^\eta} & \cX.  \\
\bar\uW \ar@/_/@<-1ex>[urrr]^(0.35)w && {}
\ultwocell_{}\omit^{}\omit{^\ze} }
\end{equation*}

Define an $\O_\cX$-{\it module\/} $\cE$ to assign an $\O_U$-module
$\cE(\uU,u)$ on $\uU=(U,\O_U)$ for all objects $(\uU,u)$ in ${\cal
C}_\cX$, and an isomorphism $\cE_{(\uf,\eta)}:\uf^*(\cE(\uV,v))
\ra\cE(\uU,u)$ for all morphisms $(\uf,\eta):(\uU,u)\ra(\uV,v)$ in
${\cal C}_\cX$, such that for all $(\uf,\eta),(\ug,\ze),
(\ug\ci\uf,\th)$ as above the following diagram of isomorphisms of
$\O_U$-modules commutes:
\e
\begin{gathered}
\xymatrix@C=7pt@R=7pt{ (\ug\ci\uf)^*\bigl(\cE(\uW,w)\bigr)
\ar[rrrrr]_(0.4){\cE_{(\ug\ci\uf,\th)}} \ar[dr]_{I_{\uf,\ug}(\cE
(\uW,w))\,\,\,\,\,\,\,\,\,\,\,\,{}} &&&&& \cE(\uU,u), \\
& \uf^*\bigl(\ug^*(\cE(\uW,w)\bigr)
\ar[rrr]^(0.55){\uf^*(\cE_{(\ug,\ze)})} &&&
\uf^*\bigl(\cE(\uV,v)\bigr) \ar[ur]_{{}\,\,\,\cE_{(\uf,\eta)}} }
\end{gathered}
\label{ds8eq1}
\e
for $I_{\uf,\ug}(\cE (\uW,w))$ as in Remark~\ref{ds2rem3}.

A {\it morphism of\/ $\O_\cX$-modules\/} $\phi:\cE\ra\cF$ assigns a
morphism of $\O_U$-modules $\phi(\uU,u):\cE(\uU,u)\ra\cF(\uU,u)$ for
each object $(\uU,u)$ in ${\cal C}_\cX$, such that for all morphisms
$(\uf,\eta):(\uU,u)\ra(\uV,v)$ in ${\cal C}_\cX$ the following
commutes:
\begin{equation*}
\xymatrix@C=40pt@R=11pt{ \uf^*\bigl(\cE(\uV,v)\bigr)
\ar[d]_{\uf^*(\phi(\uV,v))} \ar[r]_{\cE_{(\uf,\eta)}} & \cE(\uU,u)
\ar[d]^{\phi(\uU,u)} \\ \uf^*\bigl(\cF(\uV,v)\bigr)
\ar[r]^{\cF_{(\uf,\eta)}} & \cF(\uU,u). }
\end{equation*}

We call $\cE$ {\it quasicoherent}, or a {\it vector bundle of
rank\/}\I{Deligne--Mumford $C^\iy$-stack!vector bundles on} $n$, if
$\cE(\uU,u)$ is quasicoherent, or a vector bundle of rank $n$,
respectively, for all $(\uU,u)\in{\cal C}_\cX$. Write $\OcXmod$ for
the category of $\O_\cX$-modules, and $\qcoh(\cX)$, $\vect(\cX)$ for
the full subcategories of quasicoherent sheaves and vector bundles,
respectively. Then $\OcXmod$ is an abelian category,\I{abelian
category} and $\qcoh(\cX)$ an abelian subcategory of $\OcXmod$. If
$\cX$ is locally fair then~$\qcoh(\cX)=\OcXmod$.
\label{ds8def9}
\end{dfn}

Note that vector bundles $\cE$ on $\cX$ are locally trivial in the
\'etale topology, but need not be locally trivial in the Zariski
topology. In particular, the orbifold groups $\Iso_\cX([x])$ of
$\cX$ can act nontrivially on the fibres $\cE\vert_x$ of~$\cE$.

As in \cite[\S 10.5]{Joyc4}, as well as sheaves of $\OcX$-modules,
we can define other kinds of sheaves on Deligne--Mumford
$C^\iy$-stacks $\cX$ by the same method. In particular, to define
d-stacks in \S\ref{ds10}, we will need {\it sheaves of abelian
groups\/}\I{Deligne--Mumford $C^\iy$-stack!sheaves of abelian groups
on} and {\it sheaves of\/ $C^\iy$-rings\/}\I{Deligne--Mumford
$C^\iy$-stack!sheaves of C-rings on@sheaves of $C^\iy$-rings on} on
Deligne--Mumford $C^\iy$-stacks.

\begin{ex} Let $\cX$ be a Deligne--Mumford $C^\iy$-stack. Define
a quasicoherent sheaf $\OcX$ on $\cX$ called the {\it structure
sheaf\/}\I{Deligne--Mumford $C^\iy$-stack!structure sheaf $\O_\cX$}
of $\cX$ by $\OcX(\uU,u)=\O_U$ for all objects $(\uU,u)$ in ${\cal
C}_\cX$, and $(\OcX)_{(\uf,\eta)}:\uf^*(\O_V)\ra\O_U$ is the natural
isomorphism for all morphisms $(\uf,\eta):(\uU,u)\ra(\uV,v)$
in~${\cal C}_\cX$.

We may also consider $\OcX$ as a sheaf of $C^\iy$-rings on~$\cX$.
\label{ds8ex4}
\end{ex}

\begin{ex} Let $\cX$ be a Deligne--Mumford $C^\iy$-stack. Define
an $\O_\cX$-module $T^*\cX$ called the {\it cotangent
sheaf\/}\I{Deligne--Mumford $C^\iy$-stack!cotangent sheaf} of $\cX$
by $(T^*\cX)(\uU,u)=T^*\uU$ for all objects $(\uU,u)$ in ${\cal
C}_\cX$ and $(T^*\cX)_{(\uf,\eta)}=\Om_\uf:\uf^*(T^*\uV)\ra T^*\uU$
for all morphisms $(\uf,\eta):(\uU,u)\ra(\uV,v)$ in ${\cal C}_\cX$,
where $T^*\uU$ and $\Om_\uf$ are as in~\S\ref{ds24}.
\label{ds8ex5}
\end{ex}

\begin{ex} Let $\uX$ be a $C^\iy$-scheme. Then $\cX=\ul{\bar
X\!}\,$ is a Deligne--Mumford $C^\iy$-stack. We will define an {\it
inclusion functor\/} ${\cal I}_\uX:\OXmod\ra\OcXmod$. Let $\cE$ be
an object in $\OXmod$. If $(\uU,u)$ is an object in ${\cal C}_\cX$
then $u:\bar\uU\ra\cX=\ul{\bar X\!}\,$ is 2-isomorphic to
$\bar\uu:\bar\uU\ra\ul{\bar X\!}\,$ for some unique morphism
$\uu:\uU\ra\uX$. Define $\cE'(\uU,u)=\uu^*(\cE)$. If
$(\uf,\eta):(\uU,u)\ra(\uV,v)$ is a morphism in ${\cal C}_\cX$ and
$\uu,\uv$ are associated to $u,v$ as above, so that $\uu=\uv\ci\uf$,
then define
\begin{equation*}
\smash{\cE'_{(\uf,\eta)}= I_{\uf,\uv}(\cE)^{-1}:\uf^*(\cE'(\uV,v))=
\uf^*\bigl(\uv^*(\cE)\bigr)\longra (\uv\ci\uf)^*(\cE)=\cE'(\uU,u).}
\end{equation*}
Then \eq{ds8eq1} commutes for all $(\uf,\eta),(\ug,\ze)$, so $\cE'$
is an $\O_\cX$-module.

If $\phi:\cE\ra\cF$ is a morphism of $\O_X$-modules then we define a
morphism $\phi':\cE'\ra\cF'$ in $\OcXmod$ by
$\phi'(\uU,u)=\uu^*(\phi)$ for $\uu$ associated to $u$ as above.
Then defining ${\cal I}_\uX:\cE\mapsto\cE'$, ${\cal
I}_\uX:\phi\mapsto\phi'$ gives a functor $\OXmod\ra\OcXmod$, which
induces equivalences between the categories $\OXmod,\qcoh(\uX)$
defined in \S\ref{ds24} and $\OcXmod,\qcoh(\cX)$ above.
\label{ds8ex6}
\end{ex}

\begin{dfn} Let $f:\cX\ra\cY$ be a 1-morphism of
Deligne--Mumford $C^\iy$-stacks, and $\cF$ be an $\O_\cY$-module. A
{\it pullback\/}\I{Deligne--Mumford $C^\iy$-stack!quasicoherent
sheaves on!pullbacks|(} of $\cF$ to $\cX$ is an $\O_\cX$-module
$\cE$, together with the following data: if $\uU,\uV$ are
$C^\iy$-schemes and $u:\bar\uU\ra\cX$ and $v:\bar\uV\ra\cY$ are
\'etale 1-morphisms, then there is a $C^\iy$-scheme $\uW$ and
morphisms $\upi_\uU:\uW\ra\uU$, $\upi_\uV:\uW\ra\uV$ giving a
2-Cartesian\I{2-category!2-Cartesian square} diagram:
\e
\begin{gathered}
\xymatrix@C=60pt@R=10pt{ \bar\uW \ar[r]_(0.25){\bar\upi_\uV}
\ar[d]_{\bar\upi_\uU}
\drtwocell_{}\omit^{}\omit{^\ze} & \bar\uV \ar[d]^v \\
\bar\uU \ar[r]^(0.7){f\ci u} & \cY.}
\end{gathered}
\label{ds8eq2}
\e
Then an isomorphism $i(\cF,f,u,v,\ze):\ab\upi^*_\uU\bigl(
\cE(\uU,u)\bigr)\ra \upi^*_\uV\bigl(\cF(\uV,v)\bigr)$ of
$\O_W$-modules should be given, which is functorial in $(\uU,u)$ in
${\cal C}_\cX$ and $(\uV,v)$ in ${\cal C}_\cY$ and the 2-isomorphism
$\ze$ in \eq{ds8eq2}. We usually write pullbacks $\cE$ as
$f^*(\cF)$. Pullbacks $f^*(\cF)$ exist, and are unique up to unique
isomorphism. Using the Axiom of Choice, we choose a pullback
$f^*(\cF)$ for all such~$f,\cF$.

Let $f:\cX\ra\cY$ be a 1-morphism, and $\phi:\cE\ra\cF$ be a
morphism in $\OcYmod$. Then $f^*(\cE),f^*(\cF)\in\OcXmod$. The {\it
pullback morphism\/} $f^*(\phi):f^*(\cE)\ra f^*(\cF)$ is the unique
morphism in $\OcXmod$ such that whenever $u:\bar\uU\ra\cX$,
$v:\bar\uV\ra\cY$, $\uW,\upi_\uU,\upi_\uV$ are as above, the
following diagram in $\O_W$-mod commutes:
\begin{equation*}
\xymatrix@C=50pt@R=10pt{ \upi^*_\uU\bigl(f^*(\cE)(\uU,u)\bigr)
\ar[r]_{i(\cE,f,u,v,\ze)} \ar[d]_{\pi^*_\uU(f^*(\phi)(\uU,u))} &
\upi^*_\uV\bigl(\cE(\uV,v)\bigr) \ar[d]^{\pi_\uV^*(\phi(\uV,v))} \\
\upi^*_\uU\bigl(f^*(\cF)(\uU,u)\bigr) \ar[r]^{i(\cF,f,u,v,\ze)} &
\upi^*_\uV\bigl(\cF(\uV,v)\bigr).}
\end{equation*}
This defines a right exact functor $f^*:\OcYmod\ra\OcXmod$, which
also maps $\qcoh(\cY)\ra\qcoh(\cX)$.

Let $f,g:\cX\ra\cY$ be 1-morphisms of Deligne--Mumford
$C^\iy$-stacks, $\eta:f\Ra g$ a 2-morphism, and $\cE\in\OcYmod$.
Then we have $\O_\cX$-modules $f^*(\cE),g^*(\cE)$. Define
$\eta^*(\cE):f^*(\cE)\ra g^*(\cE)$ to be the unique isomorphism such
that whenever $\uU,\uV,\uW,u,v,\upi_\uU,\upi_\uV$ are as above, so
that we have 2-Cartesian diagrams
\begin{equation*}
\xymatrix@C=30pt@R=10pt{ \bar\uW \ar[rrr]_(0.8){\bar\upi_\uV}
\ar[d]_{\bar\upi_\uU}  &
\drrtwocell_{}\omit^{}\omit{^{\!\!\!\!\!\!\!\!\!\!\!
\!\!\!\!\!\!\!\!\!\!\!\!\!\ze\od(\eta*\id_{u\ci\bar\upi_\uU})
\,\,\,\,\,\,\,\,\,\,\,\,{}}} && \bar\uV \ar[d]^v & \bar\uW
\ar[rrr]_(0.7){\bar\upi_\uV} \ar[d]_{\bar\upi_\uU}
& \drtwocell_{}\omit^{}\omit{^\ze} && \bar\uV \ar[d]^v \\
\bar\uU \ar[rrr]^(0.8){f\ci u} &&& \cY, & \bar\uU
\ar[rrr]^(0.7){g\ci u} &&& \cY,}
\end{equation*}
as in \eq{ds8eq2}, then we have commuting isomorphisms of
$\O_W$-modules:
\begin{equation*}
\xymatrix@C=100pt@R=-3pt{ \upi^*_\uU\bigl(f^*(\cE)(\uU,u)\bigr)
\ar[dr]^{i(\cE,f,u,v,\ze\od(\eta*\id_{u\ci\bar\upi_\uU}))}
\ar[dd]_{\upi^*_\uU((\eta^*(\cE))(\uU,u))} \\
& \upi^*_\uV\bigl(\cE(\uV,v)\bigr). \\
\upi^*_\uU\bigl(g^*(\cE)(\uU,u)\bigr) \ar[ur]_{i(\cE,g,u,v,\ze)}}
\end{equation*}
This defines a natural isomorphism $\eta^*:f^*\Ra g^*$.

As in Remark \ref{ds2rem3}, if $f:\cX\ra\cY$ and $g:\cY\ra\cZ$ are
1-morphisms of Deligne--Mumford $C^\iy$-stacks and $\cE\in\OcZmod$,
then we have a canonical isomorphism $I_{f,g}(\cE): (g\ci
f)^*(\cE)\ra f^*(g^*(\cE))$. If $\cX$ is a Deligne--Mumford
$C^\iy$-stack and $\cE\in\OcXmod$, we have a canonical isomorphism
$\de_\cX(\cE):\id_\cX^*(\cE)\ra\cE$. These $I_{f,g},\de_\cX$ have
the same properties as in the $C^\iy$-scheme case.

In a similar way, we can define pullbacks $f^{-1}(\cE)$ for sheaves
of abelian groups and of $C^\iy$-rings $\cE$ on $\cY$, and
corresponding
isomorphisms~$I_{f,g}(\cE),\de_\cX(\cE)$.\I{Deligne--Mumford
$C^\iy$-stack!quasicoherent sheaves on!pullbacks|)}
\label{ds8def10}
\end{dfn}

\begin{ex} Let $f:\cX\ra\cY$ be a 1-morphism of Deligne--Mumford
$C^\iy$-stacks. Then Example \ref{ds8ex4} defines sheaves of
$C^\iy$-rings $\OcX,\OcY$ on $\cX,\cY$, so as in Definition
\ref{ds8def10} we have a pullback sheaf $f^{-1}(\OcY)$ of
$C^\iy$-rings on $\cX$. There is a natural morphism
$f^\sh:f^{-1}(\OcY)\ra\OcX$ of sheaves of $C^\iy$-rings on $\cX$,
characterized by the following property: for all
$(\uU,u),(\uV,v),\uW,\ze$ as in Definition \ref{ds8def10}, the
following diagram of sheaves of $C^\iy$-rings on $W$ commutes:
\begin{equation*}
\xymatrix@!0@C=70pt@R=30pt{
\pi^{-1}_U\bigl(f^{-1}(\OcY)(\uU,u)\bigr)
\ar@<-2ex>[d]^{i(\OcY,f,u,v,\ze)}_\cong
\ar[rr]_{\pi^{-1}_U(f^\sh(\uU,u))} &&
\pi^{-1}_U\bigl((\OcX)(\uU,u)\bigr)
\ar@{=}[r] & \pi^{-1}_U(\O_U) \ar[d]^\cong_{\pi_U^\sh} \\
\pi^{-1}_V\bigl(\OcY(\uV,v)\bigr) \ar@{=}[r] & \pi^{-1}_V(\O_V\bigr)
\ar[rr]^{\pi_V^\sh} && \O_W,}
\end{equation*}
where $\upi_\uU=(\pi_U,\pi_U^\sh)$ and $\upi_\uV=(\pi_V,\pi_V^\sh)$.
\label{ds8ex7}
\end{ex}

\begin{dfn} Let $f:\cX\ra\cY$ be a 1-morphism of Deligne--Mumford
$C^\iy$-stacks. Then $f^*(T^*\cY),T^*\cX$\I{Deligne--Mumford
$C^\iy$-stack!cotangent sheaf|(} are $\O_\cX$-modules, by Example
\ref{ds8ex5} and Definition \ref{ds8def10}. Define
$\Om_f:f^*(T^*\cY)\ra T^*\cX$ to be the unique morphism
characterized as follows. Let $u:\bar\uU\ra\cX$, $v:\bar\uV\ra\cY$,
$\uW,\upi_\uU,\upi_\uV$ be as in Definition \ref{ds8def10}, with
\eq{ds8eq2} 2-Cartesian. Then the following diagram commutes
in~$\O_W$-mod:
\begin{equation*}
\xymatrix@C=15pt@R=15pt{ \upi^*_\uU\bigl(f^*(T^*\cY)(\uU,u)\bigr)
\ar[d]_{\pi^*_\uU(\Om_f(\uU,u))} \ar[rrr]_{i(T^*\cY,f,u,v,\ze)}
&&& \upi^*_\uV\bigl((T^*\cY)(\uV,v)\bigr) \ar@{=}[r] &
\upi^*_\uV(T^*\uV) \ar[d]_{\Om_{\upi_\uV}} \\
\upi^*_\uU\bigl((T^*\cX)(\uU,u)\bigr)
\ar[rrr]^{(T^*\cX)_{(\upi_\uU,\id_{u\ci\upi_\uU})}} &&&
(T^*\cX)(\uW,u\ci\upi_\uU) \ar@{=}[r] & T^*\uW.}
\end{equation*}
\label{ds8def11}
\end{dfn}

Here \cite[Th.~10.15]{Joyc4} is the analogue of
Theorem~\ref{ds2thm2}.

\begin{thm}{\bf(a)} Let\/ $f:\cX\ra\cY$ and\/ $g:\cY\ra\cZ$ be\/
$1$-morphisms of Deligne--Mumford\/ $C^\iy$-stacks.
Then\/~$\Om_{g\ci f}=\Om_f\ci f^*(\Om_g)\ci I_{f,g}(T^*\cZ)$.
\smallskip

\noindent{\bf(b)} Let\/ $f,g:\cX\ra\cY$ be\/ $1$-morphisms of
Deligne--Mumford\/ $C^\iy$-stacks and\/ $\eta:f\Ra g$ a
$2$-morphism. Then $\Om_f=\Om_g\ci\eta^*(T^*\cY):f^*(T^*\cY)\ra
T^*\cX$.
\smallskip

\noindent{\bf(c)} Suppose\/ $\cW,\cX,\cY,\cZ$ are locally fair
Deligne--Mumford\/ $C^\iy$-stacks with a
$2$-Cartesian\I{2-category!2-Cartesian square} square
\begin{equation*}
\xymatrix@C=60pt@R=10pt{ \cW \ar[r]_(0.25)f \ar[d]^e
\drtwocell_{}\omit^{}\omit{^\eta}
 & \cY \ar[d]_h \\
\cX \ar[r]^(0.7)g & \cZ}
\end{equation*}
in $\DMCStalf,$ so that\/ $\cW\simeq\cX\t_\cZ\cY$. Then the
following is exact in $\qcoh(\cW)\!:$\I{Deligne--Mumford
$C^\iy$-stack!quasicoherent sheaves on|)}\I{Deligne--Mumford
$C^\iy$-stack!cotangent sheaf|)}
\begin{equation*}
\xymatrix@C=16pt{ (g\ci e)^*(T^*\cZ)
\ar[rrrrrr]^(0.51){\begin{subarray}{l}e^*(\Om_g)\ci I_{e,g}(T^*\cZ)\op\\
-f^*(\Om_h)\ci I_{f,h}(T^*\cZ)\ci\eta^*(T^*\cZ)
\end{subarray}} &&&&&&
{\raisebox{8pt}{$\displaystyle \begin{subarray}{l}\ts e^*(T^*\cX)\op\\
\ts f^*(T^*\cY)\end{subarray}$}} \ar[rr]^(0.58){\Om_e\op \Om_f} &&
T^*\cW \ar[r] & 0.}
\end{equation*}
\label{ds8thm3}
\end{thm}

\subsection{Orbifold strata of Deligne--Mumford $C^\iy$-stacks}
\label{ds87}
\I{Deligne--Mumford $C^\iy$-stack!orbifold strata|(}\I{orbifold
strata!of Deligne--Mumford $C^\iy$-stacks|(}

Let $\cX$ be a Deligne--Mumford $C^\iy$-stack, and $\Ga$ a finite
group. In \cite[\S 11.1]{Joyc4} we define six different notions of
{\it orbifold strata\/} of $\cX$, which are Deligne--Mumford
$C^\iy$-stacks written $\cX^\Ga,\tcX^\Ga,\hcX^\Ga$, and open
$C^\iy$-substacks $\cX{}^\Ga_\ci\subseteq\cX^\Ga$,
$\tcX^\Ga_\ci\subseteq\tcX^\Ga$, $\hcX^\Ga_\ci\subseteq
\hcX^\Ga$.\G[XGad]{$\cX^\Ga,\tcX^\Ga,\hcX^\Ga,
\cX{}^\Ga_\ci,\tcX^\Ga_\ci,\hcX^\Ga_\ci$}{orbifold strata of a
Deligne--Mumford $C^\iy$-stack $\cX$} The points and orbifold
groups\I{C-stack@$C^\iy$-stack!orbifold group $\Iso_\cX([x])$} of
$\cX^\Ga,\ldots,\hcX^\Ga_\ci$ are given by:
\begin{itemize}
\setlength{\itemsep}{0pt}
\setlength{\parsep}{0pt}
\item[(i)] Points of $\cX^\Ga$ are isomorphism classes
$[x,\rho]$, where $[x]\in\cX_\top$ and
$\rho:\Ga\ra\Iso_\cX([x])$ is an injective morphism, and
$\Iso_{\smash{\cX^\Ga}}([x,\rho])$ is the centralizer of
$\rho(\Ga)$ in $\Iso_\cX([x])$. Points of
$\cX{}^\Ga_\ci\subseteq\cX^\Ga$ are $[x,\rho]$ with $\rho$ an
isomorphism, and $\Iso_{\smash{ \cX{}^\Ga_\ci}}([x,\rho])\cong
C(\Ga)$, the centre of $\Ga$.
\item[(ii)] Points of $\tcX^\Ga$ are pairs $[x,\De]$, where
$[x]\in\cX_\top$ and $\De\subseteq\Iso_\cX([x])$ is isomorphic
to $\Ga$, and $\Iso_{\smash{\tcX^\Ga}}([x,\De])$ is the
normalizer of $\De$ in $\Iso_\cX([x])$. Points of
$\tcX^\Ga_\ci\subseteq\tcX^\Ga$ are $[x,\De]$ with
$\De=\Iso_\cX([x])$, and~$\Iso_{\smash{\tcX^\Ga_\ci}}
([x,\De])\cong\Ga$.
\item[(iii)] Points $[x,\De]$ of $\hcX^\Ga,\hcX^\Ga_\ci$
are the same as for $\tcX^\Ga,\tcX^\Ga_\ci$, but with orbifold
groups $\Iso_{\smash{\hcX^\Ga}}([x,\De])\cong
\Iso_{\smash{\tcX^\Ga}}([x,\De])/\De$
and~$\Iso_{\smash{\hcX^\Ga_\ci}}([x,\De])\cong\{1\}$.
\end{itemize}
Since the $C^\iy$-stack $\hcX^\Ga_\ci$ has trivial orbifold groups,
it is (equivalent to) a $C^\iy$-scheme. That is, there is a genuine
$C^\iy$-scheme $\hat\uX{}^\Ga_\ci$, unique up to isomorphism in
$\CSch$, such that $\hcX^\Ga_\ci\simeq \bar{\hat\uX}{}^\Ga_\ci$
in~$\CSta$.

There are 1-morphisms $O^\Ga(\cX),\ldots, \hat\Pi{}^\Ga_\ci(\cX)$
forming a strictly commutative diagram, where the columns are
inclusions of open $C^\iy$-substacks:\G[OGaXa]{$O^\Ga(\cX),\ti
O^\Ga(\cX),O{}^\Ga_\ci(\cX),\ti O{}^\Ga_\ci(\cX)$}{1-morphisms of
orbifold strata $\cX^\Ga,\ldots,\hcX^\Ga_\ci$ of a Deligne--Mumford
$C^\iy$-stack $\cX$}\G[PiGaXa]{$\ti\Pi^\Ga(\cX),\hat\Pi{}^\Ga(\cX),
\ti\Pi^\Ga_\ci(\cX),\hat\Pi{}^\Ga_\ci(\cX)$}{1-morphisms of orbifold
strata $\cX^\Ga,\ldots,\hcX^\Ga_\ci$ of a Deligne--Mumford
$C^\iy$-stack $\cX$}
\e
\begin{split}
\xymatrix@C=55pt@R=3pt{ \cX{}^\Ga_\ci
\ar[rr]^{\ti\Pi{}^\Ga_\ci(\cX)} \ar[dr]_(0.3){O{}^\Ga_\ci(\cX)}
\ar[dd]_\subset \ar@(ul,l)[]_(0.7){\Aut(\Ga)} && \tcX^\Ga_\ci
\ar[r]^{\hat\Pi{}^\Ga_\ci(\cX)} \ar[dl]^(0.3){\ti O^\Ga_\ci(\cX)}
\ar[dd]^\subset  &
*+[r]{\hcX^\Ga_\ci\simeq
\bar{\hat\uX}{}^\Ga_\ci\!\!\!\!\!\!}
\ar@<.5ex>[dd]^\subset \\ & \cX \\
\cX^\Ga \ar[rr]_{\ti\Pi^\Ga(\cX)} \ar[ur]^(0.3){O^\Ga(\cX)}
\ar@(dl,l)[]^(0.7){\Aut(\Ga)} && \tcX^\Ga
\ar[r]_{\hat\Pi{}^\Ga(\cX)} \ar[ul]_(0.3){\ti O^\Ga(\cX)} &
*+[r]{\hcX^\Ga.} }
\end{split}
\label{ds8eq3}
\e
Also $\Aut(\Ga)$ acts on $\cX^\Ga,\cX{}^\Ga_\ci$, with
$\tcX^\Ga\simeq[\cX^\Ga/\Aut(\Ga)]$,
$\tcX^\Ga_\ci\simeq[\cX^\Ga_\ci/\Aut(\Ga)]$. The topological space
$\cX_\top$\I{C-stack@$C^\iy$-stack!underlying topological space
$\cX_\top$} of $\cX$ from \S\ref{ds82} has stratifications
\begin{equation*}
\cX_\top\cong\coprod\nolimits_{\substack{\text{iso.\ classes of}\\
\text{finite groups $\Ga$}}} \cX_{\ci,\top}^\Ga/\Ga
\cong\coprod\nolimits_\Ga \tcX_{\ci,\top}^\Ga
\cong\coprod\nolimits_\Ga \hcX_{\ci,\top}^\Ga,
\end{equation*}
which is why $\cX^\Ga,\ldots,\hcX^\Ga_\ci$ are called orbifold
strata. The 1-morphisms $O^\Ga(\cX),\ab\ti O{}^\Ga(\cX)$ in
\eq{ds8eq3} are proper, and
$\hat\Pi{}^\Ga(\cX)_\top:\tcX^\Ga_\top\ra \hcX^\Ga_\top$ is a
homeomorphism of topological spaces. Hence, if $\cX$ is compact then
$\cX^\Ga,\tcX^\Ga,\hcX^\Ga$ are also compact.

\begin{ex} Let $\cX$ be a Deligne--Mumford $C^\iy$-stack. The {\it
inertia stack\/}\I{Deligne--Mumford $C^\iy$-stack!inertia stack}
$\cI_\cX$ of $\cX$ is the fibre product
$\cI_\cX=\cX\t_{\De_\cX,\cX\t\cX,\De_\cX}\cX$, where
$\De_\cX:\cX\ra\cX\t\cX$ is the diagonal 1-morphism. One can show
there is an equivalence
\begin{equation*}
\ts\cI_\cX\simeq\coprod_{k\ge 1}\cX^{\Z_k}.
\end{equation*}
Points of $\cI_\cX$ are isomorphism classes $[x,\eta]$, where
$[x]\in\cX_\top$ and $\eta\in\Iso_\cX([x])$. Each such
$\eta\in\Iso_\cX([x])$ has some finite order $k\ge 1$, and generates
an injective morphism $\rho:\Z_k\ra\Iso_\cX([x])$ mapping
$\rho:a\mapsto\eta^a$. We may identify $\cX^{\Z_k}$ with the open
and closed $C^\iy$-substack of $[x,\eta]$ in $\cI_\cX$ for which
$\eta$ has order~$k$.
\label{ds8ex8}
\end{ex}

Orbifold strata $\cX^\Ga$ are strongly functorial\I{Deligne--Mumford
$C^\iy$-stack!orbifold strata!functoriality|(} for representable
1-morphisms and their 2-morphisms. That is, if $f:\cX\ra\cY$ is a
representable 1-morphism of Deligne--Mumford $C^\iy$-stacks, we
define a unique representable 1-morphism $f^\Ga:\cX^\Ga\ra\cY^\Ga$
with $O^\Ga(\cY)\ci f^\Ga=f\ci O^\Ga(\cX)$. If $f,g:\cX\ra\cY$ are
representable and $\eta:f\Ra g$ is a 2-morphism, we define a unique
2-morphism $\eta^\Ga:f^\Ga\Ra g^\Ga$ with $\id_{O^\Ga(\cY)}*
\eta^\Ga=\eta*\id_{O^\Ga(\cX)}$. These $f^\Ga,\eta^\Ga$ are
compatible with compositions of 1- and 2-morphisms, and identities,
in the obvious way. Orbifold strata $\tcX^\Ga$ have the same kind of
functorial behaviour, and $\hcX^\Ga$ have a weaker functorial
behaviour, in that $\hat f{}^\Ga$ is only natural up to
2-isomorphism.\I{Deligne--Mumford $C^\iy$-stack!orbifold
strata!functoriality|)}

For $f:\cX\ra\cY$ and $\Ga$ as above, the restriction
$f^\Ga\vert_{\smash{\cX^\Ga_\ci}}$ need not map
$\cX^\Ga_\ci\ra\cY^\Ga_\ci$, but only $\cX^\Ga_\ci\ra\cY^\Ga$. So we
do not define 1-morphisms $f^\Ga_\ci:\cX^\Ga_\ci\ra\cY^\Ga_\ci$. The
same applies for the actions $\ti f{}^\Ga,\hat f{}^\Ga$ of $f$ on
orbifold strata~$\tcX^\Ga_\ci,\hcX^\Ga_\ci$.

In \cite[\S 11.3]{Joyc4} we describe the orbifold strata of a
quotient $C^\iy$-stack~$[\uX/G]$.

\begin{thm} Suppose\/ $\uX$ is a separated\/ $C^\iy$-scheme and\/
$G$ a finite group acting on $\uX$ by isomorphisms, and write\/
$\cX=[\uX/G]$ for the quotient\/
$C^\iy$-stack\I{C-stack@$C^\iy$-stack!quotients $[\protect\uX/G]$}
from Example {\rm\ref{ds8ex1},} which is a Deligne--Mumford\/
$C^\iy$-stack. Let\/ $\Ga$ be a finite group. Then there are
equivalences of\/ $C^\iy$-stacks
\ea
\cX^\Ga&\simeq\coprod_{\begin{subarray}{l}\text{conjugacy classes
$[\rho]$ of injective}\\ \text{group morphisms $\rho:\Ga\ra
G$}\end{subarray}\!\!\!\!\!\!\!\!\!\!\!\!\!\!\!\!\!\!\!
\!\!\!\!\!\!\!\!\!} \bigl[\,\uX^{\rho(\Ga)}/\bigl\{g\in
G:g\rho(\ga)=\rho(\ga)g\;\> \forall\ga\in\Ga\bigr\}\bigr],
\label{ds8eq4}\\
\cX^\Ga_\ci&\simeq\coprod_{\begin{subarray}{l}\text{conjugacy
classes $[\rho]$ of injective}\\ \text{group morphisms $\rho:\Ga\ra
G$}\end{subarray}\!\!\!\!\!\!\!\!\!\!\!\!\!\!\!\!\!\!\!
\!\!\!\!\!\!\!\!\!} \bigl[\,\uX^{\rho(\Ga)}_\ci/\bigl\{g\in
G:g\rho(\ga)=\rho(\ga)g\;\> \forall\ga\in\Ga\bigr\}\bigr],
\label{ds8eq5}\\
\tcX^\Ga&\simeq\coprod_{\text{conjugacy classes $[\De]$ of subgroups
$\De\subseteq G$ with $\De\cong\Ga$} \!\!\!\!\!\!\!\!\!
\!\!\!\!\!\!\!\!\!\!\!\!\!\!\!\!\!\!\!\!\!\!\!\!\!\!\!\!\!\!\!\!
\!\!\!\!\!\!\!\!\!\!\!\!\!\!\!\!\!\!\!\!\!\!\!}
\bigl[\,\uX^\De/\bigl\{g\in G:\De=g\De g^{-1}\bigr\}\bigr],
\label{ds8eq6}\\
\tcX^\Ga_\ci&\simeq\coprod_{\text{conjugacy classes $[\De]$ of
subgroups $\De\subseteq G$ with $\De\cong\Ga$} \!\!\!\!\!\!\!\!\!
\!\!\!\!\!\!\!\!\!\!\!\!\!\!\!\!\!\!\!\!\!\!\!\!\!\!\!\!\!\!\!\!
\!\!\!\!\!\!\!\!\!\!\!\!\!\!\!\!\!\!\!\!\!\!\!}
\bigl[\,\uX^\De_\ci/\bigl\{g\in G:\De=g\De g^{-1}\bigr\}\bigr].
\label{ds8eq7}\\
\hcX^\Ga&\simeq\coprod_{\text{conjugacy classes $[\De]$ of subgroups
$\De\subseteq G$ with $\De\cong\Ga$} \!\!\!\!\!\!\!\!\!
\!\!\!\!\!\!\!\!\!\!\!\!\!\!\!\!\!\!\!\!\!\!\!\!\!\!\!\!\!\!\!\!
\!\!\!\!\!\!\!\!\!\!\!\!\!\!\!\!\!\!\!\!\!\!\!}
\bigl[\,\uX^\De\big/\bigl(\{g\in G:\De=g\De
g^{-1}\}/\De\bigr)\bigr],
\label{ds8eq8}\\
\hcX^\Ga_\ci&\simeq\coprod_{\text{conjugacy classes $[\De]$ of
subgroups $\De\subseteq G$ with $\De\cong\Ga$} \!\!\!\!\!\!\!\!\!
\!\!\!\!\!\!\!\!\!\!\!\!\!\!\!\!\!\!\!\!\!\!\!\!\!\!\!\!\!\!\!\!
\!\!\!\!\!\!\!\!\!\!\!\!\!\!\!\!\!\!\!\!\!\!\!}
\bigl[\,\uX^\De_\ci\big/\bigl(\{g\in G:\De=g\De
g^{-1}\}/\De\bigr)\bigr].
\label{ds8eq9}
\ea
Here for each subgroup\/ $\De\subseteq G,$ we write\/ $\uX{}^\De$
for the closed\/ $C^\iy$-subscheme in\/ $\uX$ fixed by $\De$ in\/
$G,$ and\/ $\uX{}^\De_\ci$ for the open $C^\iy$-subscheme in
$\uX{}^\De$ of points in $\uX$ whose stabilizer group in $G$ is
exactly\/ $\De$. In {\rm\eq{ds8eq4}--\eq{ds8eq5},} morphisms
$\rho,\rho':\Ga\ra G$ are conjugate if\/ $\rho'=\Ad(g)\ci\rho$ for
some $g\in G,$ and subgroups $\De,\De'\subseteq G$ are conjugate
if\/ $\De=g\De'g^{-1}$ for some $g\in G$. In
\eq{ds8eq4}--\eq{ds8eq9} we sum over one representative $\rho$ or
$\De$ for each conjugacy class.
\label{ds8thm4}
\end{thm}

Let $\cX$ be a Deligne--Mumford $C^\iy$-stack and $\Ga$ a finite
group, so that as above we have an orbifold stratum $\cX^\Ga$ with a
1-morphism $O^\Ga(\cX):\cX^\Ga\ra\cX$. Let $\cE$ be a quasicoherent
sheaf on $\cX$,\I{Deligne--Mumford $C^\iy$-stack!quasicoherent
sheaves on!restriction to orbifold strata} so that
$\cE^\Ga:=O^\Ga(\cX)^*(\cE)$ is a quasicoherent sheaf on $\cX^\Ga$.
In \cite[\S 11.4]{Joyc4} we show that there is a natural
representation of $\Ga$ on $\cE^\Ga$ by isomorphisms. Also the
action of $\Aut(\Ga)$ on $\cX^\Ga$ lifts naturally to $\cE^\Ga$, so
that $\Aut(\Ga)\lt\Ga$ acts equivariantly on~$\cE^\Ga$.

Write $R_0,\ldots,R_k$ for the irreducible representations of $\Ga$
over $\R$ (that is, we choose one representative $R_i$ in each
isomorphism class of irreducible representations), with $R_0=\R$ the
trivial representation. Then the $\Ga$-representation on $\cE^\Ga$
induces a splitting
\e
\cE^\Ga\cong\ts\bigop_{i=0}^k\cE^\Ga_i\ot R_i\quad\text{for
$\cE_0^\Ga,\ldots,\cE^\Ga_k\in\qcoh(\cX^\Ga)$.}
\label{ds8eq10}
\e
We will be interested in splitting $\cE^\Ga$ into {\it trivial\/}
and {\it nontrivial\/} representations of $\Ga$, denoted by
subscripts `tr' and `nt'. So we write
\e
\cE^\Ga=\cE^\Ga_\tr\op\cE^\Ga_\nt,
\label{ds8eq11}
\e
where $\cE^\Ga_\tr,\cE^\Ga_\nt$ are the subsheaves of $\cE^\Ga$
corresponding to the factors $\cE^\Ga_0\ot R_0$ and
$\bigop_{i=1}^k\cE^\Ga_i\ot R_i$ respectively. The same applies for
the orbifold stratum~$\cX{}^\Ga_\ci\subseteq\cX^\Ga$.

We also have an orbifold stratum $\tcX^\Ga$ with a 1-morphism $\ti
O{}^\Ga(\cX):\tcX^\Ga\ra\cX$, so that $\ti{\cal E}{}^\Ga:=\ti
O^\Ga(\cX)^*(\cE)$ is a quasicoherent sheaf on $\tcX^\Ga$. In
general there is no natural $\Ga$-representation on $\ti{\cal
E}{}^\Ga$, as the quotient by $\Aut(\Ga)$ in
$\tcX^\Ga\simeq[\cX^\Ga/ \Aut(\Ga)]$ does not preserve the
$\Ga$-action. However, we do have a natural splitting
\e
\ti{\cal E}{}^\Ga=\ti{\cal E}{}^\Ga_\tr\op\ti{\cal E}{}^\Ga_\nt
\label{ds8eq12}
\e
corresponding to \eq{ds8eq11}. The same applies
for~$\tcX^\Ga_\ci\subseteq\tcX^\Ga$.

As in \eq{ds8eq3}, for the orbifold stratum $\hcX^\Ga$ we do not
have a natural 1-morphism $\hcX^\Ga\ra\cX$, so we cannot just pull
$\cE$ back to $\hcX^\Ga$. Instead, we push $\ti{\cal E}{}^\Ga$ down
to $\hcX^\Ga$ along the 1-morphism
$\hat\Pi{}^\Ga:\tcX^\Ga\ra\hcX^\Ga$. It turns out that in the
splitting \eq{ds8eq12}, the push down $\hat\Pi{}^\Ga_* (\ti{\cal
E}{}^\Ga_\nt)$ is zero, since $\hat\Pi{}^\Ga$ has fibre
$[\ul{*}/\Ga]$, and $\hat\Pi{}^\Ga_*$ essentially takes
$\Ga$-equivariant parts. So we define
$\hat\cE{}^\Ga_\tr=\hat\Pi{}^\Ga_* (\ti{\cal E}{}^\Ga_\tr)$, a
quasicoherent sheaf on $\hcX^\Ga$. The same applies
for~$\hcX^\Ga_\ci\subseteq\hcX^\Ga$.

When passing to orbifold strata, it is often natural to restrict to
the trivial parts $\cE{}^\Ga_\tr,\ti{\cal
E}{}^\Ga_\tr,\hat\cE{}^\Ga_\tr$ of the pullbacks of $\cE$. The next
theorem illustrates this.

\begin{thm} Let\/ $\cX$ be a Deligne--Mumford\/ $C^\iy$-stack
and\/ $\Ga$ a finite group, so that we have a\/ $1$-morphism
$O^\Ga(\cX):\cX^\Ga\ra\cX$. As in Example\/ {\rm\ref{ds8ex5}} we
have cotangent sheaves\I{Deligne--Mumford $C^\iy$-stack!cotangent
sheaf} $T^*\cX,T^*(\cX^\Ga)$ and a morphism
$\Om_{\smash{O^\Ga(\cX)}}:O^\Ga(\cX)^*(T^*\cX)\ab \ra T^*(\cX^\Ga)$
in $\qcoh(\cX^\Ga)$. But\/ $O^\Ga(\cX)^*(T^*\cX)=(T^*\cX)^\Ga,$ so
by \eq{ds8eq11} we have a splitting
$(T^*\cX)^\Ga=(T^*\cX)^\Ga_\tr\op(T^*\cX)^\Ga_\nt$. Then
$\Om_{\smash{O^\Ga(\cX)}}\vert_{\smash{(T^*\cX)^\Ga_\tr}}:
(T^*\cX)^\Ga_\tr\ra T^*(\cX^\Ga)$ is an isomorphism,
and\/~$\Om_{\smash{O^\Ga(\cX)}}\vert_{\smash{(T^*\cX)^\Ga_\nt}}=0$.

Similarly, using $\ti\O{}^\Ga(\cX): \tcX^\Ga\ra\cX$ and\/
\eq{ds8eq12} for ${}\,\,\,\,\,\widetilde{\!\!\!\!\!
(T^*\cX)\!\!\!\!\!} \,\,\,\,\,^\Ga$ we find that\/
$\Om_{\smash{\ti\O{}^\Ga(\cX)}} \vert_{\smash{\,\,\,\,\,
\widetilde{\!\!\!\!\!(T^*\cX)\!\!\!\!\!} \,\,\,\,\, ^\Ga_\tr}}:
\,\,\,\,\,\widetilde{\!\!\!\!\!(T^*\cX)
\!\!\!\!\!}\,\,\,\,\,^\Ga_\tr\!\ra\! T^*(\tcX^\Ga)$ is an
isomorphism, and\/ $\Om_{\smash{\ti\O{}^\Ga(\cX)}}
\vert_{\smash{\,\,\,\,\, \widetilde{\!\!\!\!\!(T^*\cX)\!\!\!\!\!}
\,\,\,\,\,^\Ga_\nt}}\!=\!0$. Also, there is a natural isomorphism
${}\,\,\,\,\,\widehat{ \!\!\!\!\!(T^*\cX)
\!\!\!\!\!}\,\,\,\,\,^\Ga_\tr\cong T^*(\hcX^\Ga)$
in\/~$\qcoh(\hcX^\Ga)$.\I{Deligne--Mumford
$C^\iy$-stack|)}\I{C-algebraic geometry@$C^\iy$-algebraic
geometry|)}\I{Deligne--Mumford $C^\iy$-stack!orbifold
strata|)}\I{orbifold strata!of Deligne--Mumford $C^\iy$-stacks|)}
\label{ds8thm5}
\end{thm}

\section{Orbifolds}
\label{ds9}
\I{orbifold|(}\I{orbifold!with corners|see{orbifold with \\
corners}}\I{orbifold!with boundary|see{orbifold with \\ boundary}}

We now summarize \cite[\S 8.1--\S 8.4]{Joyc6} on orbifolds. There is
already a substantial literature on orbifolds, and \S\ref{ds91}
indicates the main milestones in the field, and explains how our
definition of orbifolds relates to those by other authors.

\subsection{Different ways to define orbifolds}
\label{ds91}
\I{orbifold!different definitions|(}

Orbifolds are geometric spaces locally modelled on $\R^n/G$, for
$G\subset\GL(n,\R)$ a finite group. There are several {\it
nonequivalent\/} definitions of orbifolds in the literature, which
are reviewed in \cite[\S 8.1]{Joyc6}. They were first defined by
Satake \cite{Sata} (who called them `V-manifolds') and Thurston
\cite[\S 13]{Thur}. Satake and Thurston defined an orbifold to be a
Hausdorff topological space $X$ with an atlas of charts
$(U_i,\Ga_i,\phi_i)$ for $i\in I$, where $\Ga_i\subset\GL(n,\R)$ is
a finite subgroup, $U_i\subseteq\R^n$ a $\Ga_i$-invariant open
subset, and $\phi_i:U_i/\Ga_i\ra X$ a homeomorphism with an open set
in $X$, compatible on overlaps $\phi_i(U_i/\Ga_i)
\cap\phi_j(U_j/\Ga_j)$ in $X$. Smooth maps between orbifolds are
continuous maps $f:X\ra Y$, which lift locally to equivariant smooth
maps on the charts.

There is a problem with this notion of smooth maps: some
differential-geometric operations, such as pullbacks of vector
bundles by smooth maps, may not be well-defined. To fix this
problem, new definitions were needed. Moerdijk and Pronk
\cite{Moer,MoPr} defined orbifolds to be {\it proper \'etale Lie
groupoids\/}\I{orbifold!as groupoid in $\Man$|(} in $\Man$. Chen and
Ruan \cite[\S 4]{ChRu} gave an alternative theory more in the spirit
of \cite{Sata,Thur}. A book on orbifolds in the sense of
\cite{ChRu,Moer,MoPr} is Adem, Leida and Ruan~\cite{ALR}.

All of \cite{ALR,ChRu,Sata,Thur,Moer,MoPr} regard orbifolds as an
ordinary category. But orbifolds are differential-geometric
analogues of Deligne--Mumford stacks,\I{stack} which form a
2-category.\I{2-category} So it seems natural to define a 2-category
of orbifolds $\Orb$. Several important geometric constructions need
the extra structure of a 2-category to work properly. For example,
transverse fibre products\I{orbifold!transverse fibre
products}\I{2-category!fibre products in} exist in the 2-category
$\Orb$, where they satisfy a universal property involving
2-morphisms, as in \S\ref{dsA4}. In the homotopy category
$\Ho(\Orb)$,\I{2-category!homotopy category}\I{homotopy category}
`transverse fibre products' can be defined as an {\it ad hoc\/}
geometric construction, but they are not fibre products in the
category-theoretic sense, and do not satisfy a universal
property.\I{orbifold!a category or a 2-category?}

There are two main routes in the literature for defining a
2-category\I{2-category} of orbifolds $\Orb$. The first, as in Pronk
\cite{Pron} and Lerman \cite[\S 3.3]{Lerm}, is to define orbifolds
to be groupoids in $\Man$ as in \cite{Moer,MoPr}. But to define 1-
and 2-morphisms in $\Orb$ one must do more work: one makes proper
\'etale Lie groupoids into a 2-category $\bf Gpoid$, and then $\Orb$
is defined as a (weak) 2-category\I{2-category!weak} localization of
$\bf Gpoid$ at a suitable class of 1-morphisms.

The second route, as in Behrend and Xu \cite[\S 2]{BeXu}, Lerman
\cite[\S 4]{Lerm} and Metzler \cite[\S 3.5]{Metz}, is to define
orbifolds as a class of Deligne--Mumford stacks\I{orbifold!as stack
on $\Man$|(}\I{stack} on the site $(\Man,{\cal J}_\Man)$ of
manifolds with Grothendieck topology\I{Grothendieck topology} ${\cal
J}_\Man$ coming from open covers. The relationship between the two
routes is discussed in~\cite{BeXu,Lerm,Pron}.

In the `classical' approaches to orbifolds \cite{ALR,ChRu,Moer,
MoPr,Sata,Thur}, the objects, orbifolds, have a simple definition,
but the smooth maps, or 1- and 2-morphisms, are either badly
behaved, or very complicated to define. In contrast, in the `stacky'
approaches to orbifolds \cite{BeXu,Joyc4,Lerm,Metz}, the objects are
very complicated to define, but 1- and 2-morphisms are well-behaved
and easy to define --- 1-morphisms are just functors, and
2-morphisms are natural isomorphisms.

Our approach, described in \cite[\S 8.2]{Joyc6}, is similar to the
second route: we define orbifolds to be special examples of
Deligne--Mumford $C^\iy$-stacks,\I{orbifold!as Deligne--Mumford
$C^\iy$-stack} so that they are stacks\I{stack} on the site
$(\CSch,{\cal J})$. This will be convenient for our work on d-stacks
and d-orbifolds, which are also based on $C^\iy$-stacks.

\begin{dfn} An {\it orbifold of dimension\/} $n$ is a separated,
second countable Deligne--Mumford $C^\iy$-stack $\cX$ such that for
every $[x]\in\cX_\top$ there exist a linear action of
$G=\Iso_\cX([x])$ on $\R^n$, a $G$-invariant open neighbourhood $U$
of 0 in $\R^n$, and a 1-morphism $i:[\uU/G]\ra\cX$ which is an
equivalence with an open neighbourhood $\cU\subseteq\cX$ of $[x]$ in
$\cX$ with $i_\top([0])=[x]$, where $\uU=F_\Man^\CSch(U)$.

Write $\Orb$\G[Orb]{$\Orb$}{2-category of orbifolds} for the full
2-subcategory of orbifolds in $\DMCSta$. We may refer to 1-morphisms
$f:\cX\ra\cY$ in $\Orb$ as {\it smooth maps\/} of orbifolds. Define
a full and faithful\I{functor!full}\I{functor!faithful} functor
$F_\Man^\Orb:\Man\ra\Orb$ by~$F_\Man^\Orb=F_\CSch^\CSta\ci
F_\Man^\CSch$.
\label{ds9def1}
\end{dfn}

Here is \cite[Th.~9.26 \& Cor.~9.27]{Joyc4}. Since equivalent
(2-)categories are considered to be `the same', the moral of Theorem
\ref{ds9thm1} is that our orbifolds are essentially the same objects
as those considered by other recent authors.

\begin{thm} The $2$-category\/ $\Orb$ of orbifolds without boundary
defined above is equivalent to the $2$-categories of orbifolds
considered as stacks on $\Man$\I{orbifold!as stack on
$\Man$|)}\I{stack} defined in Metzler\/ {\rm\cite[\S 3.4]{Metz}} and
Lerman\/ {\rm\cite[\S 4]{Lerm},} and also equivalent as a weak\/
$2$-category\I{2-category!weak} to the weak\/ $2$-categories of
orbifolds regarded as proper \'etale Lie groupoids\I{orbifold!as
groupoid in $\Man$|)} defined in Pronk\/ {\rm\cite{Pron}} and
Lerman\/~{\rm\cite[\S 3.3]{Lerm}}.

Furthermore, the homotopy category\I{2-category!homotopy
category}\I{homotopy category} $\Ho(\Orb)$ of\/ $\Orb$ (that is, the
category whose objects are objects in\/ $\Orb,$ and whose morphisms
are $2$-isomorphism classes of\/ $1$-morphisms in\/ $\Orb$) is
equivalent to the category of orbifolds regarded as proper \'etale
Lie groupoids defined in Moerdijk\/ {\rm\cite{Moer}}. Transverse
fibre products in\/ $\Orb$ agree with the corresponding fibre
products in\/~$\CSta$.\I{orbifold!different definitions|)}
\label{ds9thm1}
\end{thm}

We define five classes of smooth maps:

\begin{dfn} Let $f:\cX\ra\cY$ be a smooth map (1-morphism) of
orbifolds.
\begin{itemize}
\setlength{\itemsep}{0pt}
\setlength{\parsep}{0pt}
\item[(i)] We call $f$ {\it representable\/} if it acts
injectively on orbifold groups, that is, $f_*:\Iso_\cX([x])\ra
\Iso_\cY\bigl(f_\top([x])\bigr)$ is an injective morphism for
all~$[x]\in\cX_\top$.\I{orbifold!representable 1-morphism}

Equivalently, $f$ is representable if it is a representable
1-morphism of $C^\iy$-stacks. This means that whenever $\uV$ is
a $C^\iy$-scheme and $\Pi:\bar\uV\ra\cY$ is a 1-morphism then
the $C^\iy$-stack fibre product $\cX\t_{f,\cY,\Pi}\bar\uV$ is a
$C^\iy$-scheme.
\item[(ii)] We call $f$ an {\it
immersion\/}\I{orbifold!immersion} if it is representable and
$\Om_f:f^*(T^*\cY)\ra T^*\cX$ is a surjective morphism of vector
bundles, i.e.\ has a right inverse in~$\qcoh(\cX)$.
\item[(iii)] We call $f$ an {\it
embedding\/}\I{orbifold!embedding} if it is an immersion, and
$f_*:\Iso_\cX([x])\ra\Iso_\cY \bigl(f_\top([x])\bigr)$ is an
isomorphism for all $[x]\in\cX_\top$, and
$f_\top:\cX_\top\ra\cY_\top$ is a homeomorphism with its image
(so in particular it is injective).
\item[(iv)] We call $f$ a {\it
submersion\/}\I{orbifold!submersion} if $\Om_f: f^*(T^*\cY)\ra
T^*\cX$ is an injective morphism of vector bundles, i.e.\ has a
left inverse in~$\qcoh(\cX)$.
\item[(v)] We call $f$ {\it \'etale\/}\I{orbifold!etale
1-morphism@\'etale 1-morphism} if it is representable and
$\Om_f: f^*(T^*\cY)\ra T^*\cX$ is an isomorphism, or
equivalently, if $f$ is \'etale as a 1-morphism of
$C^\iy$-stacks.
\end{itemize}
Note that submersions are not required to be representable.
\label{ds9def2}
\end{dfn}

\begin{dfn} An orbifold $\cX$ is called {\it
effective\/}\I{orbifold!effective|(} if $\cX$ is locally modelled
near each $[x]\in\cX_\top$ on $\R^n/G,$ where $G$ acts effectively
on $\R^n$, that is, every $1\ne\ga\in G$ acts nontrivially
on~$\R^n$.
\label{ds9def3}
\end{dfn}

In \cite[\S 8.4]{Joyc6} we prove a uniqueness property for
2-morphisms of effective orbifolds.

\begin{prop} Let\/ $\cX,\cY$ be effective orbifolds, and\/
$f,g:\cX\ra\cY$ be $1$-morphisms. Suppose that either:
\begin{itemize}
\setlength{\itemsep}{0pt}
\setlength{\parsep}{0pt}
\item[{\rm(i)}] $f$ is an embedding, a submersion, \'etale, or an
equivalence;
\item[{\rm(ii)}] $f_*:\Iso_\cX([x])\ra\Iso_\cY\bigl(f_\top([x])\bigr)$
is surjective for all\/~$[x]\in\cX_\top;$ or
\item[{\rm(iii)}] $\cY$ is a manifold.
\end{itemize}
Then there exists at most one $2$-morphism\/~$\eta:f\Ra g$.
\label{ds9prop1}
\end{prop}

Some authors include effectiveness in their definition of orbifolds.
The Satake--Thurston definitions are not as well-behaved for
noneffective orbifolds. One reason is that Proposition
\ref{ds9prop1} often allows us to treat effective orbifolds as if
they were a category rather than a 2-category, that is, one can work
in the homotopy category $\Ho(\Orb^{\bf eff})$\I{2-category!homotopy
category}\I{homotopy category} of the full 2-subcategory $\Orb^{\bf
eff}$ of effective orbifolds, because genuinely 2-categorical
behaviour comes from non-uniqueness of
2-morphisms.\I{orbifold!effective|)}

In \cite[\S 8.3]{Joyc6} we discuss {\it vector bundles\/} on
orbifolds.\I{orbifold!vector bundles on} Now an orbifold $\cX$ is an
example of a Deligne--Mumford $C^\iy$-stack, and in \S\ref{ds86} we
defined a category $\qcoh(\cX)$ of quasicoherent sheaves on $\cX$,
and a full subcategory $\vect(\cX)$ of vector bundles on $\cX$.
Unless we say otherwise, a {\it vector bundle $\cE$ on an
orbifold\/} $\cX$ will just mean an object in $\vect(\cX)$, a
special kind of quasicoherent sheaf on $\cX$, and a {\it smooth
section\/ $s$ of\/} $\cE$ will mean an element of $C^\iy(\cE)$, that
is, a morphism $s:\OcX\ra\cE$ in $\vect(\cX)$. The cotangent sheaf
$T^*\cX$ of an $n$-orbifold $\cX$ is a vector bundle on $\cX$ of
rank $n$, which we call the {\it cotangent
bundle}.\I{orbifold!cotangent bundle}

For some applications below, this point of view on vector bundles is
not ideal. If $E\ra X$ is a vector bundle on a manifold, then $E$ is
itself a manifold (with extra structure), with a submersion
$\pi:E\ra X$, and a section $s\in C^\iy(E)$ is a smooth map $s:X\ra
E$ with $\pi\ci s=\id_X$. In \S\ref{ds41}--\S\ref{ds42} we
considered d-space fibre products $\bV\t_{\bs s,\bE,\bs 0}\bV$ where
$\bV,\bE,\bs s,\bs 0=F_\Man^\dSpa(V,E,s,0)$. For the d-orbifold
analogue of this, we would like to regard a vector bundle $\cE$ over
an orbifold $\cX$ as being an orbifold in its own right, rather than
just a quasicoherent sheaf, and a section $s\in C^\iy(\cE)$ as being
a 1-morphism $s:\cX\ra\cE$ in~$\Orb$.

To get round this, in \cite[\S 8.3]{Joyc6} we define a {\it total
space functor\/}\I{orbifold!vector bundles on!total space functor
$\Tot$} $\Tot$, which to each $\cE$ in $\vect(\cX)$ associates an
orbifold $\Tot(\cE)$, called the {\it total space\/} of $\cE$, and
to each section $s\in C^\iy(\cE)$ associates a 1-morphism
$\Tot(s):\cX\ra\Tot(\cE)$ in $\Orb$. Then the d-orbifold analogue of
$\bV\t_{\bs s,\bE,\bs 0}\bV$ in Proposition \ref{ds4prop1}(c) is
$\bcV\t_{\bs s,\bcE,\bs 0}\bcV$, where~$\bcV,\bcE,\bs s,\bs
0=F_\Orb^\dSta\bigl(\cV,\Tot(\cE),\Tot(s),\Tot(0)\bigr)$.

Many other standard ideas in differential geometry extend simply to
orbifolds, such as submanifolds, transverse fibre products, and
orientations, and we will generally use these without
comment.\I{orbifold!suborbifolds}\I{orbifold!transverse fibre
products}\I{orbifold!orientations}

\subsection{Orbifold strata of orbifolds}
\label{ds92}
\I{orbifold!orbifold strata|(}\I{orbifold strata!of orbifolds|(}

Section \ref{ds87} discussed orbifold strata $\cX^\Ga,\ldots,
\hcX^\Ga_\ci$ of a Deligne--Mumford $C^\iy$-stack $\cX$. In \cite[\S
8.4]{Joyc6} we work these ideas out for orbifolds. If $\cX$ is an
orbifold, then $\cX^\Ga,\ldots, \hcX^\Ga_\ci$ need not be orbifolds,
as the next example shows, but are disjoint unions of orbifolds of
different dimensions.

\begin{ex} Let the real projective plane $\RP^2$ have homogeneous
coordinates $[x_0,x_1,x_2]$, and let $\Z_2=\{1,\si\}$ act on $\RP^2$
by $\si:[x_0,x_1,x_2]\mapsto[x_0,x_1,-x_2]$. The fixed point locus
of $\si$ in $\RP^2$ is the disjoint union of the circle
$\bigl\{[x_0,x_1,0]:[x_0,x_1]\in\RP^1\bigr\}$ and the
point~$\bigl\{[0,0,1]\}$.

Write $\ul{\RP\!}\,^2=F_\Man^\CSch(\RP^2)$, and form the quotient
orbifold $\cX=[\ul{\RP\!}\,^2/\Z_2]$. Then \eq{ds8eq4} shows that
the orbifold stratum $\cX^{\Z_2}$ is the disjoint union of orbifolds
$\ul{\RP\!}\,^1\t[\ul */\Z_2]$ and $[\ul */\Z_2]$ of dimensions 1
and 0, respectively. Note that $\cX^{\Z_2}$ is not an orbifold, as
it does not have pure dimension, and nor are~$\tcX^{\Z_2},\ldots,
\hcX^{\Z_2}_\ci$.
\label{ds9ex1}
\end{ex}

So that our constructions remain within the world of orbifolds, we
will find it useful to define a decomposition
$\cX^\Ga=\coprod_{\smash{\la\in\La^\Ga_+}}\cX^{\Ga,\la}$ of
$\cX^\Ga$ such that each $\cX^{\Ga,\la}$ is an orbifold of
dimension~$\dim\cX-\dim\la$.

\begin{dfn} Let $\Ga$ be a finite group. Consider representations
$(V,\rho)$ of $\Ga$, where $V$ is a finite-dimensional real vector
space and $\rho:\Ga\ra\Aut(V)$ a group morphism. We call $(V,\rho)$
{\it nontrivial\/} if $V^{\rho(\Ga)}=\{0\}$. Write $\Rep_\nt(\Ga)$
for the abelian category\I{abelian category} of nontrivial
$(V,\rho)$, and $K_0(\Rep_\nt(\Ga))$ for its Grothendieck group.
Then any $(V,\rho)$ in $\Rep_\nt(\Ga)$ has a class
$\bigl[(V,\rho)\bigr]$ in $K_0(\Rep_\nt(\Ga))$. For brevity, we will
use the notation
$\La^\Ga=K_0\bigl(\Rep_\nt(\Ga)\bigr)$\G[LaGa]{$\La^\Ga$}{lattice
generated by nontrivial representations of a finite group $\Ga$} and
$\La^\Ga_+=\bigl\{\bigl[(V,\rho)\bigr]:(V,\rho)\in\Rep_\nt(\Ga)\bigr\}
\subseteq\La^\Ga$.\G[LaGa+]{$\La^\Ga_+$}{`positive cone' of classes
of $\Ga$-representations in lattice $\La^\Ga$} We think of
$\La^\Ga_+$ as the `positive cone' in~$\La^\Ga$.

By elementary representation theory, up to isomorphism $\Ga$ has
finitely many irreducible representations. Let $R_0,R_1,\ldots,R_k$
be choices of irreducible representations in these isomorphism
classes, with $R_0=\R$ the trivial irreducible representation, so
that $R_1,\ldots,R_k$ are nontrivial. Then $\La^\Ga$ is freely
generated over $\Z$ by $[R_1],\ldots,[R_k]$, so that
\begin{align*}
\La^\Ga&=\bigl\{a_1[R_1]+\cdots+a_k[R_k]:a_1,\ldots,a_k\in\Z\bigr\},
\qquad\text{and}\\
\La^\Ga_+&=\bigl\{a_1[R_1]+\cdots+a_k[R_k]:a_1,\ldots,a_k\in\N\bigr\}
\subseteq\La^\Ga,
\end{align*}
where $\N=\{0,1,2,\ldots\}\subset\Z$. Hence $\La^\Ga\cong\Z^k$
and~$\La^\Ga_+\cong\N^k$.

Define a group morphism $\dim:\La^\Ga\ra\Z$ by
$\dim:a_1[R_1]+\cdots+a_k[R_k]\mapsto a_1\dim R_1+\cdots+a_k\dim
R_k$, so that $\dim:[(V,\rho)]\mapsto\dim V$.
Then~$\dim(\La^\Ga_+)\subseteq\N$.

Now let $\cX$ be an orbifold. As in \eq{ds8eq10}--\eq{ds8eq11} we
have decompositions $O^\Ga(\cX)^*(T^*\cX)\!=\!(T^*\cX)^\Ga_\tr\!
\op\! (T^*\cX)^\Ga_\nt$ with
$(T^*\cX)^\Ga_\tr\!\cong\!(T^*\cX)^\Ga_0\!\ot\! R_0$ and
$(T^*\cX)^\Ga_\nt \!\cong\!\bigop_{i=1}^k(T^*\cX)^\Ga_i\!\ot\! R_i$,
where $(T^*\cX)^\Ga_0,\ldots,\ab(T^*\cX)^\Ga_k\in\qcoh(\cX^\Ga)$.
Since $T^*\cX$ is a vector bundle, $O^\Ga(\cX)^*(T^*\cX)$ is a
vector bundle, and so the $(T^*\cX)^\Ga_i$ are {\it vector bundles
of mixed rank}, that is, locally they are vector bundles, but their
ranks may vary on different connected components of~$\cX^\Ga$.

For each $\la\in\La^\Ga_+$, define
$\cX^{\Ga,\la}$\G[XGae]{$\cX^{\Ga,\la},\tcX^{\Ga,\mu},\hcX^{\Ga,\mu},
\cX{}^{\Ga,\la}_\ci,\tcX^{\Ga,\mu}_\ci,\hcX^{\Ga,\mu}_\ci$}{orbifold
strata of an orbifold $\cX$} to be the open and closed
$C^\iy$-substack in $\cX^\Ga$ with
$\rank\bigl((T^*\cX)^\Ga_1\bigr)[R_1]+\cdots+\rank\bigl((T^*\cX)^\Ga_k
\bigr)[R_k]=\la$ in $\La^\Ga_+$. Then
$(T^*\cX)^\Ga_\nt\vert_{\cX^{\Ga,\la}}$ is a vector bundle of rank
$\dim\la$, so $(T^*\cX)^\Ga_\tr\vert_{\smash{\cX^{\Ga,\la}}}$ is a
vector bundle of dimension $\dim\cX-\dim\la$ on $\cX^{\Ga,\la}$. But
$(T^*\cX)^\Ga_\tr\cong T^*(\cX^\Ga)$ by Theorem \ref{ds8thm5}. Hence
$T^*(\cX^{\Ga,\la})$ is a vector bundle of rank $\dim\cX-\dim\la$.
Since $\cX^\Ga$ is a disjoint union of orbifolds of different
dimensions, we see that $\cX^{\Ga,\la}$ is an orbifold, with
$\dim\cX^{\Ga,\la}=\dim\cX-\dim\la$.
Then~$\cX^\Ga=\coprod_{\la\in\La^\Ga_+}\cX^{\Ga,\la}$.

Write $O^{\Ga,\la}(\cX)=O^\Ga(\cX)\vert_{\smash{\cX^{\Ga,\la}}}:
\cX^{\Ga,\la}\ra\cX$. It is a proper, representable immersion of
orbifolds. We interpret $(T^*\cX)^\Ga_\nt
\vert_{\smash{\cX^{\Ga,\la}}}$ as the {\it conormal bundle\/} of
$\cX^{\Ga,\la}$ in $\cX$. It carries a nontrivial
$\Ga$-representation of class $\la\in\La^\Ga_+$, so we refer to
$\la$ as the {\it conormal\/ $\Ga$-representation}
of~$\cX^{\Ga,\la}$.

Define $\cX^{\Ga,\la}_\ci=\cX^\Ga_\ci\cap\cX^{\Ga,\la}$, and
$O^{\Ga,\la}_\ci(\cX)=O^\Ga_\ci(\cX)\vert_{\smash{\cX^{\Ga,\la}_\ci}}:
\cX^{\Ga,\la}_\ci\ra\cX$. Then $\cX^{\Ga,\la}_\ci$ is an orbifold
with $\dim\cX^{\Ga,\la}_\ci=\dim\cX-\dim\la$, and
$\cX^\Ga_\ci=\ts\coprod_{\smash{\la\in\La^\Ga_+}}\cX^{\Ga,\la}_\ci$.

As in \S\ref{ds87}, we have $\tcX^\Ga\simeq[\cX^\Ga/ \Aut(\Ga)]$.
Now $\Aut(\Ga)$ acts on the right on $\Rep_\nt(\Ga)$ by
$\al:(V,\rho)\mapsto (V,\rho\ci\al)$ for $\al\in\Aut(\Ga)$, and this
induces right actions of $\Aut(\Ga)$ on
$\La^\Ga=K_0\bigl(\Rep_\nt(\Ga)\bigr)$ and
$\La^\Ga_+\subseteq\La^\Ga$. Write these actions as
$\al:\la\mapsto\la\cdot\al$. Then the action of $\al\in\Aut(\Ga)$ on
$\cX^\Ga$ maps $\cX^{\Ga,\la}\ra \cX^{\Ga,\la\cdot\al}$. Write
$\La^\Ga_+/\Aut(\La)$ for the set of $\Aut(\Ga)$-orbits
$\mu=\la\cdot\Aut(\Ga)$ in $\La^\Ga_+$. The map $\dim:\La^\Ga\ra\Z$
is $\Aut(\Ga)$-invariant, and so descends
to~$\dim:\La^\Ga/\Aut(\Ga)\ra\Z$.

Then $\coprod_{\la\in\mu}\cX^{\Ga,\la}$ is an open and closed
$\Aut(\Ga)$-invariant $C^\iy$-substack in $\cX^\Ga$ for each
$\mu\in\La^\Ga_+/\Aut(\La)$, so we may define $\tcX^{\Ga,\mu}\simeq
\bigl[\bigl(\coprod_{\la\in\mu}
\cX^{\Ga,\la}\bigr)\big/\Aut(\Ga)\bigr]$, an open and closed
$C^\iy$-substack of $\tcX^\Ga\simeq[\cX^\Ga/ \Aut(\Ga)]$. Write
$\tcX^{\Ga,\mu}_\ci=\tcX^\Ga_\ci\cap \tcX^{\Ga,\mu}$. Then
$\tcX^{\Ga,\mu},\tcX {}^{\Ga,\mu}_\ci$ are orbifolds of dimension
$\dim\cX-\dim\mu$, with
\begin{equation*}
\tcX^\Ga=\ts\coprod_{\mu\in\La^\Ga_+/\Aut(\Ga)}\tcX^{\Ga,\mu}
\quad\text{and}\quad \tcX^\Ga_\ci=
\ts\coprod_{\mu\in\La^\Ga_+/\Aut(\Ga)} \tcX^{\Ga,\mu}_\ci.
\end{equation*}
Set $\ti O{}^{\Ga,\mu}(\cX)\!=\!\ti O{}^\Ga(\cX)\vert_{\smash{
\tcX^{\Ga,\mu}}}:\tcX^{\Ga,\mu}\!\ra\!\cX$ and $\ti
O{}^{\Ga,\mu}_\ci(\cX)\!=\!\ti O{}^\Ga_\ci(\cX)\vert_{\smash{\tcX
{}^{\Ga,\mu}_\ci}}:\tcX^{\Ga,\mu}_\ci\!\ra\!\cX$. Then $\ti
O^{\Ga,\mu}(\cX),\ti O{}^{\Ga,\mu}_\ci(\cX)$ are representable
immersions, with $\ti O^{\Ga,\mu}(\cX)$ proper.

The 1-morphism $\hat\Pi{}^\Ga(\cX):\tcX^\Ga\ra\hcX^\Ga$ maps open
and closed $C^\iy$-substacks of $\tcX^\Ga$ to open and closed
$C^\iy$-substacks of $\hcX^\Ga$. Let
$\hcX^{\Ga,\mu}=\hat\Pi{}^\Ga(\cX) (\tcX^{\Ga,\mu})$ for each
$\mu\in\La^\Ga_+/\Aut(\La)$, and write
$\hcX^{\Ga,\mu}_\ci=\hcX^\Ga_\ci\cap \hcX^{\Ga,\mu}$. Then
$\hcX^{\Ga,\mu}, \hcX^{\Ga,\mu}_\ci$ are orbifolds of dimension
$\dim\cX-\dim\mu$, with
\begin{equation*}
\hcX^\Ga=\ts\coprod_{\mu\in\La^\Ga_+/\Aut(\Ga)}\hcX^{\Ga,\mu}
\quad\text{and}\quad \hcX^\Ga_\ci=
\ts\coprod_{\mu\in\La^\Ga_+/\Aut(\Ga)} \hcX^{\Ga,\mu}_\ci.
\end{equation*}
\label{ds9def4}
\end{dfn}

If $f:\cX\ra\cY$ is a representable 1-morphism of Deligne--Mumford
$C^\iy$-stacks and $\Ga$ a finite group, then as in \S\ref{ds87} we
have a representable 1-morphism of orbifold strata
$f^\Ga:\cX^\Ga\ra\cY^\Ga$. Note that if $\cX,\cY$ are orbifolds,
then $f^\Ga$ need not map $\cX^{\Ga,\la}\ra\cY^{\Ga,\la}$, or map
$\cX^\Ga_\ci\ra\cY^\Ga_\ci$. The analogue applies for~$\ti
f{}^\Ga,\hat f{}^\Ga$.

Some important properties of orbifolds can be characterized by the
vanishing of certain orbifold strata $\cX^{\Ga,\la}$. For example:
\begin{itemize}
\setlength{\itemsep}{0pt}
\setlength{\parsep}{0pt}
\item An orbifold $\cX$ is {\it locally
orientable\/}\I{orbifold!locally orientable} if and only if
$\cX^{\Z_2,\la}=\es$ for all odd $\la\in\La^{\Z_2}_+
\cong\N=\{0,1,2,\ldots\}$.
\item An orbifold $\cX$ is {\it
effective\/}\I{orbifold!effective} in the sense of Definition
\ref{ds9def3} if and only if $\cX^{\Ga,0}=\es$ for all
nontrivial finite groups~$\Ga$.
\end{itemize}

In \cite[\S 8.4]{Joyc6} we consider the question: if $\cX$ is an
oriented orbifold,\I{orbifold!orientations}\I{orbifold!orbifold
strata!orientations on|(} can we define orientations on the orbifold
strata $\cX^{\Ga,\la},\ldots,\hcX^{\Ga,\mu}_\ci$? Here is an
example:

\begin{ex} Let ${\cal S}^4=\bigl\{(x_1,\ldots,x_5)\in\R^5:
x_1^2+\cdots+x_5^2=1\bigr\}$, an oriented 4-manifold. Let
$G=\{1,\si,\tau,\si\tau\}\cong\Z_2^2$ act on ${\cal S}^4$ preserving
orientations by
\begin{align*}
\si:(x_1,\ldots,x_5)&\longmapsto (x_1,x_2,x_3,-x_4,-x_5),\\
\tau:(x_1,\ldots,x_5)&\longmapsto (-x_1,-x_2,-x_3,-x_4,x_5),\\
\si\tau:(x_1,\ldots,x_5)&\longmapsto (-x_1,-x_2,-x_3,x_4,-x_5).
\end{align*}
Then $\cX=[\ul{\cal S}^4/G]$ is an oriented 4-orbifold. The orbifold
groups $\Iso_\cX([x])$ for $[x]\in\cX_\top$ are all $\{1\}$ or
$\Z_2$. The singular locus of $\cX$ is the disjoint union of a copy
of $\RP^2$ from the fixed points $\pm(x_1,x_2,x_3,0,0)$ of $\si$,
and two isolated points $\{\pm(0,0,0,0,1)\}$ and
$\{\pm(0,0,0,1,0)\}$ from the fixed points of $\tau$ and~$\si\tau$.

Identifying $\La_+^{\Z_2}$ and $\La_+^{\Z_2}/\Aut(\Z_2)$ with $\N$,
it follows that
\begin{align*}
\cX^{\Z_2,2}&=\cX^{\Z_2,2}_\ci\cong \tcX^{\Z_2,2}=
\tcX^{\Z_2,2}_\ci\cong\RP^2\t[\ul */\Z_2],&
\hcX^{\Z_2,2}&=\hcX^{\Z_2,2}_\ci\cong\RP^2,\\
\cX^{\Z_2,4}&=\cX^{\Z_2,4}_\ci\cong \tcX^{\Z_2,4}
=\tcX^{\Z_2,4}_\ci\cong[\ul */\Z_2]\amalg [\ul */\Z_2],&
\hcX^{\Z_2,4}&=\hcX^{\Z_2,4}_\ci\cong *\amalg *.
\end{align*}
Since $\RP^2$ is not orientable, we see that $\cX$ is an oriented
orbifold, but none of $\cX^{\Z_2,2},\tcX^{\Z_2,2},
\hcX^{\Z_2,2},\cX^{\Z_2,2}_\ci,\tcX^{\Z_2,2}_\ci, \hcX^{\Z_2,2}_\ci$
are orientable.
\label{ds9ex2}
\end{ex}

Thus, we can only orient $\cX^{\Ga,\la},\ldots, \hcX^{\Ga,\mu}_\ci$
for all oriented orbifolds $\cX$ under some conditions on
$\Ga,\la,\mu$. The next proposition sets out these conditions:

\begin{prop}{\bf(a)} Suppose $\Ga$ is a finite group and\/
$(V,\rho)$ a nontrivial\/ $\Ga$-representation which has no
odd-dimensional subrepresentations, and write
$\la=[(V,\rho)]\in\La^\Ga_+$. Choose an orientation on $V$. Then for
all oriented orbifolds\/ $\cX$ we can define natural orientations on
the orbifold strata\/~$\cX^{\Ga,\la},\cX^{\Ga,\la}_\ci$.

If\/ $\md{\Ga}$ is odd then all nontrivial\/ $\Ga$-representations
are even-dimensional, so we can orient\/
$\cX^{\Ga,\la},\cX^{\Ga,\la}_\ci$ for all\/~$\la\in\La^\Ga_+$.
\smallskip

\noindent{\bf(b)} Let\/ $\Ga,(V,\rho),\la$ be as in {\bf(a)\rm,} and
set\/ $\mu=\la\cdot\Aut(\Ga)\in\La^\Ga_+/\Aut(\Ga)$. Write $H$ for
the subgroup of\/ $\Aut(\Ga)$ fixing\/ $\la$ in\/ $\La^\Ga_+$. Then
for each\/ $\de\in H$ there exists an isomorphism of\/
$\Ga$-representations $i_\de:(V,\rho\ci\de)\ra(V,\rho)$. Suppose
$i_\de:V\ra V$ is orientation-preserving for all\/ $\de\in H$. If\/
$\la\in 2\La^\Ga_+$ this holds automatically.

Then for all oriented orbifolds\/ $\cX$ we can define orientations
on the orbifold strata\/ $\tcX^{\Ga,\mu},\hcX^{\Ga,\mu},
\tcX^{\Ga,\mu}_\ci,\hcX^{\Ga,\mu}_\ci$. For $\tcX^{\Ga,\mu}$ this
works as\/ $\tcX^{\Ga,\mu}\simeq [\cX^{\Ga,\la}/H],$ where
$\cX^{\Ga,\la}$ is oriented by {\bf(a)\rm,} and the $H$-action on
$\cX^{\Ga,\la}$ preserves orientations, so the orientation on
$\cX^{\Ga,\la}$ descends to an orientation
on\/~$\tcX^{\Ga,\mu}\simeq [\cX^{\Ga,\la}/H]$.
\smallskip

\noindent{\bf(c)} Suppose that\/ $\Ga$ and\/ $\la\in\La^\Ga_+$ do
not satisfy the conditions in {\bf(a)\rm,} or\/ $\Ga$ and\/
$\mu\in\La^\Ga_+/\Aut(\Ga)$ do not satisfy the conditions in
{\bf(b)}. Then as in Example\/ {\rm\ref{ds9ex2}} we can find
examples of oriented orbifolds\/ $\cX$ such that\/
$\cX^{\Ga,\la},\cX^{\Ga,\la}_\ci$ are not orientable, or
$\tcX^{\Ga,\mu},\hcX^{\Ga,\mu},\tcX^{\Ga,\mu}_\ci,
\hcX^{\Ga,\mu}_\ci$ are not orientable, respectively. That is, the
conditions on $\Ga,\la,\mu$ in {\bf(a)\rm,\bf(b)} are necessary as
well as sufficient to be able to orient orbifold strata\/
$\cX^{\Ga,\la},\ldots, \hcX^{\Ga,\mu}_\ci$ of all oriented
orbifolds\/~$\cX$.\I{orbifold|)}\I{orbifold!orbifold
strata|)}\I{orbifold strata!of orbifolds|)}\I{orbifold!orbifold
strata!orientations on|)}
\label{ds9prop2}
\end{prop}

Note that Proposition \ref{ds9prop2}(a),(b) do not apply in Example
\ref{ds9ex2}, since the nontrivial representation of $\Z_2$ on
$\R^2$ has an odd-dimensional subrepresentation.

\section{The 2-category of d-stacks}
\label{ds10}
\I{d-stack|(}\I{d-stack!with boundary|see{d-stack with \\
boundary}}\I{d-stack!with corners|see{d-stack with \\ corners}}

In \cite[Chap.~9]{Joyc6} we define and study the 2-category of {\it
d-stacks\/} $\dSta$, which are orbifold versions of d-spaces in
\S\ref{ds3}. Broadly, to go from d-spaces
$\bX=(\uX,\OXp,\EX,\im_X,\jm_X)$ to d-stacks we just replace the
$C^\iy$-scheme $\uX$ by a Deligne--Mumford $C^\iy$-stack~$\cX$.

One might expect that combining the 2-categories $\DMCSta$ and
$\dSpa$ should result in a 3-category $\dSta$, but in fact a
2-category is sufficient. For 1-morphisms $\bs f,\bs g:\bcX\ra\bcY$
in $\dSta$, a 2-morphism $\bs\eta:\bs f\Ra\bs g$ in $\dSta$ is a
pair $(\eta,\eta')$, where $\eta:f\Ra g$ is a 2-morphism in $\CSta$,
and $\eta': f^*(\cF_\cY)\ra\cE_\cX$ is as for 2-morphisms in
$\dSpa$. These $\eta,\eta'$ do not interact very much.

\subsection{The definition of d-stacks}
\label{ds101}
\I{d-stack!definition|(}\I{2-category|(}

\begin{dfn} A {\it d-stack\/}
$\bcX$\G[WXYZd]{$\bcW,\bcX,\bcY,\bcZ,\ldots$}{d-stacks, including
d-orbifolds} is a quintuple $\bcX=(\cX,\OcXp,\EcX,\im_\cX,\jm_\cX)$,
where $\cX$ is a separated, second countable, locally fair
Deligne--Mumford $C^\iy$-stack in the sense of \S\ref{ds8}, and
$\OcXp,\EcX,\im_\cX,\jm_\cX$ fit into an exact sequence of sheaves
of abelian groups on $\cX$, in the sense of~\S\ref{ds86}
\begin{equation*}
\smash{\xymatrix@C=25pt{ \EcX \ar[rr]^(0.45){\jm_\cX} && \OcXp
\ar[rr]^(0.55){\im_\cX} && \OcX \ar[r] & 0,}}
\end{equation*}
satisfying the conditions:
\begin{itemize}
\setlength{\itemsep}{0pt}
\setlength{\parsep}{0pt}
\item[(a)] $\OcXp$ is a sheaf of $C^\iy$-rings on $\cX$, and
 $\im_\cX:\OcXp\ra\OcX$ is a morphism of sheaves of
$C^\iy$-rings on $\cX$, where $\OcX$ is the structure sheaf of
$\cX$ as in Example \ref{ds8ex4}, such that for all $(\uU,u)$ in
${\cal C}_\cX$, $(U,\OcXp(\uU,u))$ is a $C^\iy$-scheme, and
$\im_\cX(\uU,u):\OcXp(\uU,u)\ra\OcX(\uU,u)=\O_U$ is a surjective
morphism of sheaves of $C^\iy$-rings on $U$, whose kernel is a
sheaf of square zero ideals.\I{square zero ideal}

We call $\im_\cX:\OcXp\ra\OcX$ satisfying these conditions a
{\it square zero extension}.\I{square zero extension}
\item[(b)] As $\im_\cX:\OcXp\ra\OcX$ is a square zero extension, its
kernel $\IcX$ is a quasi\-coherent sheaf on $\cX$. We require
that $\EcX$ is also a quasicoherent sheaf on $\cX$, and
$\jm_\cX:\EcX\ra\IcX$ is a surjective morphism in~$\qcoh(\cX)$.
\end{itemize}

The sheaf of $C^\iy$-rings $\OcXp$ has a sheaf of cotangent modules
$\Om_{\OcXp}$, which is an $\OcXp$-module with exterior derivative
$\d:\OcXp\ra\Om_{\OcXp}$. Define $\FcX=\Om_{\OcXp}\ot_\OcXp\OcX$ to
be the associated $\OcX$-module, a quasicoherent sheaf on $\cX$, and
set $\psi_\cX=\Om_{\im_\cX}\ot\id:\FcX\ra T^*\cX$, a morphism in
$\qcoh(\cX)$. Define $\phi_\cX:\EcX\ra\FcX$ to be the composition of
morphisms of sheaves of abelian groups on~$\cX$:
\begin{equation*}
\smash{\xymatrix@C=11pt{ \EcX  \ar[rr]^{\jm_\cX} && \IcX
\ar[rr]^{\d\vert_{\IcX}} && \Om_{\OcXp} \ar@{=}[r]^(0.31)\sim &
\Om_{\OcXp}\ot_\OcXp\OcXp \ar[rrr]^{\id\ot\im_\cX} &&&
\Om_{\OcXp}\ot_\OcXp\OcX \ar@{=}[r] & \FcX. }}
\end{equation*}
Then $\phi_\cX$ is a morphism in $\qcoh(\cX)$, and the following
sequence is exact:
\e
\smash{\xymatrix@C=20pt{ \EcX \ar[rr]^{\phi_\cX} && \FcX
\ar[rr]^{\psi_\cX} && T^*\cX \ar[r] & 0.}}
\label{ds10eq1}
\e
The morphism $\phi_\cX:\EcX\ra\FcX$ will be called the {\it virtual
cotangent sheaf\/}\I{d-stack!virtual cotangent sheaf} of $\bcX$. It
is a d-stack analogue of the cotangent complex\I{cotangent complex}
in algebraic geometry.

Let $\bcX,\bcY$ be d-stacks. A 1-{\it
morphism\/}\I{2-category!1-morphism} $\bs f:\bcX\ra\bcY$ is a triple
$\bs f=(f,f',f'')$, where $f:\cX\ra\cY$ is a 1-morphism of
$C^\iy$-stacks, $f':f^{-1}(\OcYp)\ra\OcXp$ a morphism of sheaves of
$C^\iy$-rings on $\cX$, and $f'':f^*(\EcY)\ra\EcX$ a morphism in
$\qcoh(\cX)$, such that the following diagram of sheaves on $\cX$
commutes:
\begin{equation*}
\xymatrix@C=11pt@R=1pt{
f^{-1}\!(\EcY)\ot_{f^{-1}(\OcY)}^{\id}\!f^{-1}\!(\OcY) \ar@{=}[r]
\ar[dd]^(0.4){{}\,\id\ot f^\sh} & f^{-1}\!(\EcY)
\ar[rr]_(0.45){\raisebox{-9pt}{$\scriptstyle f^{-1}(\jm_\cY)$}}
 && f^{-1}\!(\OcYp)
\ar[rr]_(0.52){\raisebox{-9pt}{$\scriptstyle f^{-1}(\im_\cY)$}}
\ar[ddd]^{f'} && f^{-1}\!(\OcY) \ar[r] \ar[ddd]_{f^\sh} & 0 \\ \\
{\begin{subarray}{l}\ts f^*(\EcY)=\\
\ts f^{-1}(\EcY) \ot_{f^{-1}(\OcY)}^{f^\sh}\OcX\end{subarray}}
\ar[dr]^(0.7){f''} \\  & \EcX \ar[rr]^{\jm_\cX} && \OcXp
\ar[rr]^(0.55){\im_\cX} && \OcX \ar[r] &  {0.\kern -.28em} }
\end{equation*}
Define morphisms $f^2=\Om_{f'}\ot\id:f^*(\FcY)\ra\FcX$ and
$f^3=\Om_f:f^*(T^*\cY)\ra T^*\cX$ in $\qcoh(\cX)$. Then the
following commutes in $\qcoh(\cX)$, with exact rows:
\e
\begin{gathered}
\xymatrix@C=20pt@R=12pt{ f^*(\EcY) \ar[rr]_{f^*(\phi_\cY)}
\ar[d]^{f''} && f^*(\FcY) \ar[rr]_{f^*(\psi_\cY)} \ar[d]^{f^2} &&
f^*(T^*\cY) \ar[r] \ar[d]^{f^3} & 0 \\
\EcX \ar[rr]^{\phi_\cX} && \FcX \ar[rr]^{\psi_\cX} && T^*\cX \ar[r]
& {0.\kern -.28em} }
\end{gathered}
\label{ds10eq2}
\e

If $\bcX$ is a d-stack, the {\it identity\/ $1$-morphism\/}
$\bs\id_\bcX:\bcX\ra\bcX$ is $\bs\id_\bcX=\bigl(\id_\cX,\ab
\de_\cX(\OcXp),\ab\de_\cX(\EcX)\bigr)$, with $\de_\cX(*)$ the
canonical isomorphisms of Definition~\ref{ds8def10}.

Let $\bcX,\bcY,\bcZ$ be d-stacks, and $\bs f:\bcX\ra\bcY$, $\bs
g:\bcY\ra\bcZ$ be 1-morphisms. As in \eq{ds3eq2} define the {\it
composition\I{2-category!1-morphism!composition} of\/
$1$-morphisms\/} $\bs g\ci\bs f:\bcX\ra\bcZ$ to be
\begin{equation*}
\bs g\ci\bs f=\bigl(g\ci f,f'\ci f^{-1}(g')\ci I_{f,g}(\OcZp),f''\ci
f^*(g'')\ci I_{f,g}(\EcZ)\bigr),
\end{equation*}
where $I_{*,*}(*)$ are the canonical isomorphisms of
Definition~\ref{ds8def10}.

Let $\bs f,\bs g:\bcX\ra\bcY$ be 1-morphisms of d-stacks, where $\bs
f=(f,f',f'')$ and $\bs g=(g,g',g'')$. A 2-{\it
morphism\/}\I{2-category!2-morphism} $\bs\eta:\bs f\Ra\bs g$ is a
pair $\bs\eta=(\eta,\eta')$, where $\eta:f\Ra g$ is a 2-morphism in
$\CSta$ and $\eta':f^*(\FcY)\ra\EcX$ a morphism in $\qcoh(\cX)$,
with
\begin{gather*}
g'\ci\eta^{-1}(\OcYp)=f'+\!\ka_\cX\!\ci\!\jm_\cX
\!\ci\!\eta'\!\ci\!\bigl(\id\!\ot(f^\sh\!\ci\!f^{-1}(\im_\cY))
\bigr)\!\ci \!\bigl(f^{-1}(\d)\bigr),\\
\text{and}\quad g''\ci\eta^*(\EcY)=f''+\eta'\ci f^*(\phi_\cY).
\end{gather*}
Then $g^2\ci\eta^*(\FcY)=f^2+\phi_\cX\ci\eta'$ and
$g^3\ci\eta^*(T^*\cY)=f^3$, so \eq{ds10eq2} for $\bs f,\bs g$
combine to give a commuting diagram (except $\eta'$) in
$\qcoh(\cX)$, with exact rows:
\begin{equation*}
\xymatrix@C=12pt@R=12pt{ f^*(\EcY)
\ar[ddr]_(0.6){\begin{subarray}{l}f''+ \\ \eta'\ci
f^*(\phi_\cY)\end{subarray}\!\!\!\!\!} \ar[dr]^(0.6){\eta^*(\EcY)}
\ar[rr]^{f^*(\phi_\cY)} && f^*(\FcY) \ar@{.>}[ddl]_(0.33){\eta'}
\ar[ddr]_(0.6){\begin{subarray}{l}f^2+
\\ \phi_\cX\ci\eta'\end{subarray}\!\!\!\!\!} \ar[dr]^(0.6){\eta^*(\FcY)}
\ar[rr]^{f^*(\psi_\cY)} && f^*(T^*\cY) \ar[ddr]_(0.65){f^3}
\ar[dr]^(0.6){\eta^*(T^*\cY)} \ar[r] & 0 \\
& g^*(\EcY) \ar[d]^(0.4){g''} \ar[rr]_(0.45){g^*(\phi_\cY)} &&
g^*(\FcY) \ar[d]^(0.4){g^2} \ar[rr]_(0.4){g^*(\psi_\cY)} &&
g^*(T^*\cY) \ar[d]^(0.4){g^3} \ar[r] & 0 \\
& \EcX \ar[rr]^(0.4){\phi_\cX} && \FcX \ar[rr]^(0.4){\psi_\cX} &&
T^*\cX \ar[r] & {0.\!\!} }
\end{equation*}

If $\bs f=(f,f',f''):\bcX\ra\bcY$ is a 1-morphism, the {\it
identity\/ $2$-morphism\/} $\bs\id_{\bs f}:\bs f\Ra\bs f$
is~$\bs\id_{\bs f}=(\id_f,0)$.

Let $\bs f,\bs g,\bs h:\bcX\ra\bcY$ be 1-morphisms and $\bs\eta:\bs
f\Ra\bs g$, $\bs\ze:\bs g\Ra\bs h$ 2-morphisms. Define the {\it
vertical composition\I{2-category!2-morphism!vertical composition}
of\/ $2$-morphisms\/} $\bs\ze\od\bs\eta:\bs f\Ra\bs h$ to be
\begin{equation*}
\bs\ze\od\bs\eta=\bigl(\ze\od\eta,\ze'\ci\eta^*(\FcY)+\eta'\bigr).
\end{equation*}

Suppose $\bcX,\bcY,\bcZ$ are d-stacks, $\bs f,\bs{\ti
f}:\bcX\ra\bcY$ and $\bs g,\bs{\ti g}:\bcY\ra\bcZ$ are 1-morphisms,
and $\bs\eta:\bs f\Ra\bs{\ti f}$, $\bs\ze:\bs g\Ra\bs{\ti g}$ are
2-morphisms. Define the {\it
horizontal\I{2-category!2-morphism!horizontal composition}
composition of\/ $2$-morphisms\/} $\bs\ze*\bs\eta:\bs g\ci\bs f\Ra
\bs{\ti g}\ci\bs{\ti f}$ to be
\begin{equation*}
\bs\ze*\bs\eta=\bigl(\ze*\eta,\bigl[\eta'\ci f^*(g^2)+f''\ci
f^*(\ze')+\eta'\ci f^*(\phi_\cY)\ci f^*(\ze')\bigr]\ci
I_{f,g}(\FcZ)\bigr).
\end{equation*}
This completes the definition of the 2-category of
d-stacks~$\dSta$.\G[dSta]{$\dSta$}{2-category of
d-stacks}\I{2-category|)}

Write $\DMCStalfssc$ for the 2-category of separated, second
countable, locally fair Deligne--Mumford $C^\iy$-stacks. Define a
strict 2-functor\I{2-category!strict 2-functor}
$F_\CSta^\dSta:\DMCStalfssc\ra\dSta$ to map objects $\cX$ to
$\bcX=(\cX,\OcX,0,\id_{\OcX},0)$, to map 1-morphisms $f$ to $\bs
f=(f,f^\sh,0)$, and to map 2-morphisms $\eta$ to $\bs\eta=(\eta,0)$.
Write $\hDMCStalfssc$ for the full 2-subcategory of $\bcX\in\dSta$
equivalent to $F_\CSta^\dSta(\cX)$ for $\cX\in\DMCStalfssc$. When we
say that a d-stack $\bcX$ {\it is a\/ $C^\iy$-stack}, we mean
that~$\bcX\in\hDMCStalfssc$.

Define a strict 2-functor $F_\Orb^\dSta:\Orb\ra\dSta$ by
$F_\Orb^\dSta=F_\CSta^\dSta\vert_\Orb$, noting that $\Orb$ is a full
2-subcategory of $\DMCStalfssc$. Write $\hOrb$ for the full
2-subcategory of objects $\bcX$ in $\dSta$ equivalent to
$F_\Orb^\dSta(\cX)$ for some orbifold $\cX$. When we say that a
d-stack $\bcX$ {\it is an orbifold},\I{d-stack!is an orbifold} we
mean that~$\bcX\in\hOrb$.\G[Orb']{$\hOrb$}{2-subcategory of d-stacks
equivalent to orbifolds}

Recall from \S\ref{ds81} that there is a natural (2-)functor
$F_\CSch^\CSta:\CSch\ra\CSta$ mapping $\uX\mapsto\ul{\bar X\!}\,$ on
objects and $\uf\mapsto\ul{\bar f\!}\,$ on morphisms. Also, if $\uX$
is a $C^\iy$-scheme and $\ul{\bar X\!}\,$ the corresponding
$C^\iy$-stack then Example \ref{ds8ex6} defines a functor ${\cal
I}_\uX:\OXmod\ra \O_{\smash{\ul{\bar X\!}\,}}$-mod. In the same way,
we can define functors from the category of sheaves of abelian
groups on $X$ to the category of sheaves of abelian groups on
$\ul{\bar X\!}\,$, and from the category of sheaves of $C^\iy$-rings
on $X$ to the category of sheaves of $C^\iy$-rings on $\ul{\bar
X\!}\,$, both of which we also denote by~${\cal I}_\uX$.

With this notation, define a strict 2-functor\I{2-category!strict
2-functor} $F_\dSpa^\dSta:\dSpa \ra\dSta$ to map
$\bX=(\uX,\OXp,\ab\EX,\ab\im_X,\ab\jm_X)$ to $\bcX=\bigl(\ul{\bar
X\!}\,,{\cal I}_\uX(\OXp),{\cal I}_\uX(\EX),{\cal I}_\uX(\im_X),
{\cal I}_\uX(\jm_X)\bigr)$ on objects, and to map $\bs
f=(\uf,f',f'')$ to $\bs{\hat f}=\bigl(\ul{\bar f\!}\,,{\cal
I}_\uX(f'),{\cal I}_\uX(f'')\bigr)$ on 1-morphisms, and to map
$\eta$ to $\bs\eta=\bigl(\id_{\ul{\bar f\!}\,},{\cal
I}_\uX(\eta)\bigr)$ on 2-morphisms. Write
$\hdSpa$\G[dSpa'']{$\hdSpa$}{2-subcategory of d-stacks equivalent to
d-spaces} for the full 2-subcategory of $\bcX$ in $\dSta$ equivalent
to $F_\dSpa^\dSta(\bX)$ for some $\bX$
in~$\dSpa$.\I{d-stack!definition|)}
\label{ds10def1}
\end{dfn}

In \cite[\S 9.2]{Joyc6} we prove:

\begin{thm}{\bf(a)} Definition\/ {\rm\ref{ds10def1}} defines a
strict\/ $2$-category $\dSta,$ in which all\/ $2$-morphisms are
$2$-isomorphisms.
\smallskip

\noindent{\bf(b)} $F_\CSta^\dSta,F_\Orb^\dSta$ and\/ $F_\dSpa^\dSta$
are full and faithful\I{functor!full}\I{functor!faithful} strict\/
$2$-functors. Hence $\DMCStalfssc,$ $\Orb,$ $\dSpa$ and\/
$\hDMCStalfssc,$ $\hOrb,$ $\hdSpa$ are equivalent\/ $2$-categories,
respectively.
\label{ds10thm1}
\end{thm}

\subsection{D-stacks as quotients of d-spaces}
\label{ds102}
\I{d-stack!quotients $[\bX/G]$|(}

Section \ref{ds84} defined quotient Deligne--Mumford $C^\iy$-stacks
$[\uX/G]$, quotient 1-morphisms $[\uf,\rho]:[\uX/G]\ra[\uY/H]$, and
quotient 2-morphisms $[\de]:[\uf,\rho]\Ra[\ug,\si]$. In \cite[\S
9.3]{Joyc6} we generalize all this to d-stacks. The next two
theorems summarize our results.

\begin{thm}{\bf(i)} Let\/ $\bX$ be a d-space, $G$ a finite group,
and\/ $\bs r:G\ra\Aut(\bX)$ a (strict) action of\/ $G$ on $\bX$ by\/
$1$-isomorphisms. Then we can define a \begin{bfseries}quotient
d-stack\end{bfseries}\/ $\bcX=[\bX/G],$ which is natural up to
$1$-isomorphism in $\dSta$. The underlying $C^\iy$-stack\/ $\cX$ is
$[\uX/G]$ from Example\/~{\rm\ref{ds8ex1}}.
\smallskip

\noindent{\bf(ii)} Let\/ $\bX,\bY$ be d-spaces, $G,H$ finite groups,
and\/ $\bs r:G\ra\Aut(\bX),$ $\bs s:H\ra\Aut(\bY)$ be actions of\/
$G,H$ on $\bX,\bY,$ so that by\/ {\bf(i)} we have quotient
d-stacks\/ $\bcX=[\bX/G]$ and\/ $\bcY=[\bY/H]$. Suppose $\bs
f:\bX\ra\bY$ is a\/ $1$-morphism in $\dSpa$ and $\rho:G\ra H$ is a
group morphism, satisfying $\bs f\ci \bs r(\ga)=\bs
s(\rho(\ga))\ci\bs f$ for all\/ $\ga\in G$ (this is an equality of\/
$1$-morphisms in $\dSpa,$ not just a\/ $2$-isomorphism). Then we can
define a \begin{bfseries}quotient\/
$1$-morphism\end{bfseries}\I{d-stack!quotients $[\bX/G]$!quotient
1-morphism} $\bs{\ti f}:\bcX\ra\bcY$ in $\dSta,$ which we will also
write as~$[\bs f,\rho]:[\bX/G]\ra[\bY/H]$.
\smallskip

\noindent{\bf(iii)} Let\/ $\bs{\ti f}=[\bs f,\rho]:[\bX/G]
\ra[\bY/H]$ and\/ $\bs{\ti g}=[\bs g,\si]:[\bX/G]\ra[\bY/H]$ be two
quotient\/ $1$-morphisms as in {\bf(ii)}. Suppose $\de\in H$
satisfies $\de^{-1}\,\si(\ga)=\rho(\ga)\,\de^{-1}$ for all\/ $\ga\in
G,$ and\/ $\eta:\bs f\Ra\bs s(\de^{-1})\ci\bs g$ is a $2$-morphism
in $\dSpa$ such that\/ $\eta*\id_{\bs r(\ga)}=\id_{\bs
s(\si(\ga))}*\eta$ for all\/ $\ga\in G,$ using the diagram:
\begin{equation*}
\xymatrix@C=55pt@R=15pt{
*+[r]{\bs f\ci\bs r(\ga)} \ar@{=>}[d]^{\eta*\id_{\bs r(\ga)}} \ar@{=}[rr]
&& *+[l]{\bs s(\rho(\ga))\ci\bs f}
\ar@{=>}[d]_{\id_{\bs s(\si(\ga))}*\eta} \\
*+[r]{\bs s(\de^{-1})\ci\bs g\ci\bs r(\ga)} \ar@{=}[r]
& \bs s(\de^{-1})\ci\bs s(\si(\ga))\ci \bs g \ar@{=}[r]
& *+[l]{\bs s(\rho(\ga))\ci\bs s(\de^{-1})\ci \bs g.}  }
\end{equation*}
Then we can define a \begin{bfseries}quotient\/
$2$-morphism\end{bfseries}\I{d-stack!quotients $[\bX/G]$!quotient
2-morphism} $\bs\ze:\bs{\ti f}\Ra\bs{\ti g}$ in $\dSta,$ which we
also write as\/~$[\eta,\de]:[\bs f,\rho]\Ra[\bs g,\si]$.
\label{ds10thm2}
\end{thm}

\begin{thm}{\bf(a)} Let\/ $\bcX$ be a d-stack and\/
$[x]\in\cX_\top,$ and write\/ $G=\Iso_\cX([x])$. Then there exist a
quotient d-stack\/ $[\bU/G],$ as in Theorem\/
{\rm\ref{ds10thm2}(i),} and an equivalence $\bs i:[\bU/G]\ra\bcX$
with an open d-substack\/ $\bcU$ in $\bcX,$ with\/
$i_\top:[u]\mapsto[x]\in\cU_\top\subseteq\cX_\top$ for some fixed
point\/ $u$ of\/ $G$ in\/~$U$.
\smallskip

\noindent{\bf(b)} Let\/ $\bs{\ti f}:\bcX\ra\bcY$ be a $1$-morphism
in $\dSta,$ and\/ $[x]\in\cX_\top$ with\/ $\ti
f_\top:[x]\mapsto[y]\in\cY_\top,$ and write\/ $G=\Iso_\cX([x])$
and\/ $H=\Iso_\cY([y])$. Part\/ {\bf(a)} gives $1$-morphisms $\bs
i:[\bU/G]\ra\bcX,$ $\bs j:[\bV/H]\ra\bcY$ which are equivalences
with open $\bcU\subseteq\bcX,$ $\bcV\subseteq\bcY,$ such that\/
$i_\top:[u]\mapsto[x]\in \cU_\top\subseteq\cX_\top,$
$j_\top:[v]\mapsto[y]\in \cV_\top\subseteq\cY_\top$ for $u,v$ fixed
points of\/ $G,H$ in~$U,V$.

Then there exist a $G$-invariant open d-subspace\/ $\bU'$ of\/ $u$
in $\bU$ and a quotient\/ $1$-morphism $[\bs
f,\rho]:[\bU'/G]\ra[\bV/H],$ as in Theorem\/
{\rm\ref{ds10thm2}(ii),} such that\/ $\bs f(u)=v,$ and\/ $\rho:G\ra
H$ is $\ti f_*:\Iso_\cX([x])\ra \Iso_\cY([y]),$ fitting into a
$2$-commutative diagram:
\begin{equation*}
\xymatrix@C=120pt@R=10pt{ *+[r]{[\bU'/G]} \ar[r]_(0.3){[\bs f,\rho]}
\ar[d]^{\bs i\vert_{[\bU'/G]}}
\drtwocell_{}\omit^{}\omit{^{\bs\ze}} & *+[l]{[\bV/H]} \ar[d]_{\bs j} \\
*+[r]{\bcX} \ar[r]^(0.6){\bs{\ti f}} & *+[l]{\bcY.} }
\end{equation*}

\noindent{\bf(c)} Let\/ $\bs{\ti f},\bs{\ti g}:\bcX\ra\bcY$ be
$1$-morphisms in $\dSta$ and\/ $\bs\eta:\bs{\ti f}\Ra\bs{\ti g}$ a
$2$-morphism, let\/ $[x]\in\cX_\top$ with\/ $\ti
f_\top:[x]\mapsto[y]\in\cY_\top,$ and write\/ $G=\Iso_\cX([x])$
and\/ $H=\Iso_\cY([y])$. Part\/ {\bf(a)} gives $\bs
i:[\bU/G]\ra\bcX,$ $\bs j:[\bV/H]\ra\bcY$ which are equivalences
with open $\bcU\subseteq\bcX,$ $\bcV\subseteq\bcY$ and map
$i_\top:[u]\mapsto[x],$ $j_\top:[v]\mapsto[y]$ for $u,v$ fixed
points of\/~$G,H$.

By making $\bU'$ smaller, we can take the same $\bU'$ in {\bf(b)}
for both $\bs{\ti f},\bs{\ti g}$. Thus part\/ {\bf(b)} gives a
$G$-invariant open $\bU'\subseteq\bU,$ quotient $1$-morphisms $[\bs
f,\rho]:[\bU'/G]\ra[\bV/H]$ and\/ $[\bs g,\si]:[\bU'/G]\ra[\bV/H]$
with\/ $\bs f(u)=\bs g(u)=v$ and\/ $\rho=\ti f_*:\Iso_\cX([x])\ra
\Iso_\cY([y]),$ $\si=\ti g_*:\Iso_\cX([x])\ra\Iso_\cY([y]),$ and\/
$2$-morphisms $\bs\ze:\bs{\ti f}\ci\bs i\vert_{[\bU'/G]}\Ra \bs
j\ci[\bs f,\rho],$ $\bs\th:\bs{\ti g}\ci\bs i\vert_{[\bU'/G]}\Ra\bs
j\ci[\bs g,\si]$.

Then there exist a $G$-invariant open neighbourhood\/ $\bU''$ of\/
$u$ in $\bU'$ and a quotient\/ $2$-morphism $[\la,\de]:[\bs
f\vert_{\smash{\bU''}},\rho] \Ra[\bs g\vert_{\smash{\bU''}},\si],$
as in Theorem\/ {\rm\ref{ds10thm2}(iii),} such that the following
diagram of\/ $2$-morphisms in $\dSta$ commutes:
\begin{equation*}
\xymatrix@C=120pt@R=13pt{ *+[r]{\bs{\ti f}\ci\bs i\vert_{[\bU''/G]}}
\ar@{=>}[r]_{\bs\eta*\id_{\bs i\vert_{[\bU''/G]}}}
\ar@{=>}[d]^{{}\,\,\bs\ze\vert_{[\bU''/G]}} & *+[l]{\bs{\ti g}\ci\bs
i\vert_{[\bU''/G]}} \ar@{=>}[d]_{\bs\th\vert_{[\bU''/G]}\,\,{}} \\
*+[r]{\bs j\ci[\bs f\vert_{\bU''},\rho]} \ar@{=>}[r]^{
\id_{\bs j}*[\la,\de]} & *+[l]{\bs j\ci[\bs g\vert_{\bU''},\si].} }
\end{equation*}
\label{ds10thm3}
\end{thm}

Effectively, this says that d-stacks and their 1-morphisms and
2-morphisms are Zariski locally modelled on quotient d-stacks,
quotient 1-morphisms, and quotient 2-morphisms, up to equivalence
in~$\dSta$.

In \cite[\S 9.2]{Joyc6} we define when a 1-morphism of d-stacks $\bs
f:\bcX\ra\bcY$ is {\it \'etale}.\I{d-stack!etale 1-morphism@\'etale
1-morphism} Essentially, $\bs f$ is \'etale if it is an equivalence
locally in the \'etale topology.\I{etale topology@\'etale topology}
It implies that the $C^\iy$-stack 1-morphism $f:\cX\ra\cY$ in $\bs
f$ is \'etale, and so representable.

We can characterize \'etale 1-morphisms in $\dSta$ using Theorem
\ref{ds10thm3}: a 1-morphism $\bs{\ti f}:\bcX\ra\bcY$ in $\dSta$ is
\'etale if and only if for all $[\bs f,\rho]:[\bU'/G]\ra[\bV/H]$ in
Theorem \ref{ds10thm3}(b), $\bs f:\bU'\ra\bV$ is an \'etale
1-morphism in $\dSpa$ (that is, a local equivalence in the Zariski
topology), and $\rho:G\ra H$ is injective.\I{d-stack!quotients
$[\bX/G]$|)}

\subsection{Gluing d-stacks by equivalences}
\label{ds103}
\I{d-stack!equivalence|(}\I{d-stack!gluing by equivalences|(}

Section \ref{ds32} discussed gluing d-spaces by equivalences in
$\dSpa$. In \cite[\S 9.4]{Joyc6} we generalize this to $\dSta$. Here
are the analogues of Definition \ref{ds3def2}, Proposition
\ref{ds3prop}, and Theorems \ref{ds3thm2} and~\ref{ds3thm3}.

\begin{dfn} Let $\bcX=(\cX,\OcXp,\EcX,\im_\cX,\jm_\cX)$ be a d-stack.
Suppose $\cU\subseteq\cX$ is an open $C^\iy$-substack, in the
Zariski topology, with inclusion 1-morphism $i_\cU:\cU\ra\cX$. Then
$\bcU=\bigl(\cU,\im_\cU^{-1}(\OcXp),i_\cU^*(\EcX),i_\cU^\sh\ci
i_\cU^{-1}(\im_\cX),i_\cU^*(\jm_\cX)\bigr)$ is a d-stack, where
$i_\cU^\sh:i_\cU^{-1}(\OcX)\ra\O_\cU$ is as in Example \ref{ds8ex7},
and is an isomorphism as $i_\cU$ is \'etale. We call $\bcU$ an {\it
open d-substack\/}\I{d-stack!open d-substack} of $\bcX$. An {\it
open cover\/}\I{d-stack!open cover} of a d-stack $\bcX$ is a family
$\{\bcU_a:a\in A\}$ of open d-substacks $\bcU_a$ of $\bcX$ such that
$\{\cU_a:a\in A\}$ is an open cover of $\cX$, in the Zariski
topology.\I{Zariski topology}
\label{ds10def2}
\end{dfn}

\begin{prop} Let\/ $\bcX,\bcY$ be d-stacks,
$\bcU,\bcV\!\subseteq\!\bcX$ be open d-substacks with\/
$\bcX=\bcU\cup\bcV,$ $\bs f:\bcU\ra\bcY$ and\/ $\bs g:\bcV\ra\bcY$
be $1$-morphisms, and\/ $\bs\eta:\bs f\vert_{\bcU\cap\bcV}\Ra\bs
g\vert_{\bcU\cap\bcV}$ a $2$-morphism. Then there exist a
$1$-morphism $\bs h:\bcX\ra\bcY$ and\/ $2$-morphisms $\bs\ze:\bs
h\vert_\bcU\Ra\bs f,$ $\bs\th:\bs h\vert_\bcV\Ra\bs g$ in $\dSta$
such that\/ $\bs\th\vert_{\bcU\cap\bcV}=
\bs\eta\od\bs\ze\vert_{\bcU\cap\bcV}:\bs h\vert_{\bcU\cap\bcV}\Ra
\bs g\vert_{\bcU\cap\bcV}$. This $\bs h$ is unique up to
$2$-isomorphism.
\label{ds10prop}
\end{prop}

\begin{thm} Suppose $\bcX,\bcY$ are d-stacks, $\bcU\subseteq\bcX,$
$\bcV\subseteq\bcY$ are open d-substacks, and\/ $\bs f:\bcU\ra\bcV$
is an equivalence in $\dSta$. At the level of topological spaces, we
have open $\cU_\top\subseteq\cX_\top,$ $\cV_\top\subseteq\cY_\top$
with a homeomorphism $f_\top:\cU_\top\ra\cV_\top,$ so we can form
the quotient topological space $\cZ_\top:=\cX_\top\amalg_{f_\top}
\cY_\top=(\cX_\top\amalg\cY_\top)/\sim,$ where the equivalence
relation $\sim$ on $\cX_\top\amalg\cY_\top$ identifies
$[u]\in\cU_\top\subseteq\cX_\top$ with\/~$f_\top([u])\in
\cV_\top\subseteq\cY_\top$.

Suppose $\cZ_\top$ is Hausdorff. Then there exist a d-stack\/
$\bcZ,$ open d-substacks\/ $\bs{\hat{\cal X}},\bs{\hat{\cal Y}}$ in
$\bcZ$ with\/ $\bcZ=\bs{\hat{\cal X}}\cup\bs{\hat{\cal Y}},$
equivalences $\bs g:\bcX\ra\bs{\hat{\cal X}}$ and\/ $\bs
h:\bcY\ra\bs{\hat{\cal Y}}$ such that\/ $\bs g\vert_\bcU$ and\/ $\bs
h\vert_\bcV$ are both equivalences with\/ $\bs{\hat{\cal
X}}\cap\bs{\hat{\cal Y}},$ and a $2$-morphism $\bs\eta:\bs
g\vert_\bcU\Ra\bs h\ci\bs f$. Furthermore, $\bcZ$ is independent of
choices up to equivalence.
\label{ds10thm4}
\end{thm}

\begin{thm} Suppose\/ $I$ is an indexing set, and\/ $<$ is a total
order on $I,$ and\/ $\bcX_i$ for $i\in I$ are d-stacks, and for
all\/ $i<j$ in $I$ we are given open d-substacks\/
$\bcU_{ij}\subseteq\bcX_i,$ $\bcU_{ji}\subseteq\bcX_j$ and an
equivalence $\bs e_{ij}:\bcU_{ij}\ra\bcU_{ji},$ satisfying the
following properties:
\begin{itemize}
\setlength{\itemsep}{0pt}
\setlength{\parsep}{0pt}
\item[{\rm(a)}] For all\/ $i<j<k$ in $I$ we have a
$2$-commutative diagram
\begin{equation*}
\xymatrix@C=70pt@R=10pt{ & \bcU_{ji}\cap\bcU_{jk} \ar@<.5ex>[dr]^{\bs
e_{jk}\vert_{\bcU_{ji}\cap\bU_{jk}}} \ar@{=>}[d]^{\bs\eta_{ijk}} \\
\bcU_{ij}\cap\bcU_{ik} \ar@<.5ex>[ur]^{\bs
e_{ij}\vert_{\bcU_{ij}\cap\bcU_{ik}}} \ar@<-.25ex>[rr]^(0.37){\bs
e_{ik}\vert_{\bcU_{ij}\cap\bcU_{ik}}} && \bcU_{ki}\cap\bcU_{kj}}
\end{equation*}
for some $\bs\eta_{ijk},$ where all three $1$-morphisms are
equivalences; and
\item[{\rm(b)}] For all\/ $i<j<k<l$ in $I$ the components $\eta_{ijk}$
in $\bs\eta_{ijk}=(\eta_{ijk},\eta_{ijk}')$ satisfy
\e
\eta_{ikl}\od(\id_{f_{kl}}*\eta_{ijk})\vert_{\cU_{ij}\cap
\cU_{ik}\cap\cU_{il}}=\eta_{ijl}\od (\eta_{jkl}*\id_{f_{ij}})
\vert_{\cU_{ij}\cap \cU_{ik}\cap\cU_{il}}.
\label{ds10eq3}
\e
\end{itemize}

On the level of topological spaces, define the quotient topological
space $\cY_\top=(\coprod_{i\in I}\cX_{i,\top})/\sim,$ where $\sim$
is the equivalence relation generated by $[x_i]\sim[x_j]$ if\/
$[x_i]\in \cU_{ij,\cX_i,\top}\subseteq\cX_{i,\top}$ and\/ $[x_j]\in
\cU_{ji,\top}\subseteq\cX_{j,\top}$ with\/
$e_{ij,\top}([x_i])=[x_j]$. Suppose $\cY_\top$ is Hausdorff and
second countable. Then there exist a d-stack\/ $\bcY$ and a
$1$-morphism $\bs f_i:\bcX_i\ra\bcY$ which is an equivalence with an
open d-substack\/ $\bs{\hat{\cal X}}_i\subseteq\bcY$ for all\/ $i\in
I,$ where $\bcY=\bigcup_{i\in I}\bs{\hat{\cal X}}_i,$ such that\/
$\bs f_i\vert_{\bcU_{ij}}$ is an equivalence
$\bcU_{ij}\ra\bs{\hat{\cal X}}_i\cap\bs{\hat{\cal X}}_j$ for all\/
$i<j$ in $I,$ and there exists a $2$-morphism\/ $\bs\eta_{ij}:\bs
f_j\ci\bs e_{ij}\Ra\bs f_i\vert_{\bcU_{ij}}$. The d-stack\/ $\bcY$
is unique up to equivalence.

Suppose also that\/ $\bcZ$ is a d-stack, and\/ $\bs
g_i:\bcX_i\ra\bcZ$ are $1$-morphisms for all\/ $i\in I,$ and there
exist\/ $2$-morphisms $\bs\ze_{ij}:\bs g_j\ci\bs e_{ij}\Ra\bs
g_i\vert_{\bcU_{ij}}$ for all\/ $i<j$ in $I,$ such that for all\/
$i<j<k$ in $I$ the components $\ze_{ij},\eta_{ijk}$ in
$\bs\ze_{ij},\bs\eta_{ijk}$ satisfy
\e
\bigl(\ze_{ij}\vert_{\cU_{ij}\cap\cU_{ik}}\bigr)\od
\bigl(\ze_{jk}*\id_{e_{ij}}\vert_{\cU_{ij}\cap\cU_{ik}}\bigr)=
\bigl(\ze_{ik}\vert_{\cU_{ij}\cap\cU_{ik}}\bigr)\od
\bigl(\id_{g_k}*\eta_{ijk}\vert_{\cU_{ij}\cap\cU_{ik}}\bigr).
\label{ds10eq4}
\e
Then there exist a $1$-morphism $\bs h:\bcY\ra\bcZ$ and\/
$2$-morphisms $\bs\ze_i:\bs h\ci\bs f_i\Ra\bs g_i$ for all\/ $i\in
I$. The $1$-morphism $\bs h$ is unique up to $2$-isomorphism.
\label{ds10thm5}
\end{thm}

\begin{rem} Note that in Proposition \ref{ds3prop} for d-spaces,
$\bs h$ is independent of $\eta$ up to 2-isomorphism, but in
Proposition \ref{ds10prop} for d-stacks, $\bs h$ may depend on
$\bs\eta$. Similarly, in Theorem \ref{ds3thm3} for d-spaces, we
impose no conditions on 2-morphisms $\eta_{ijk}$ on quadruple
overlaps or $\ze_{ij}$ on triple overlaps, but in Theorem
\ref{ds10thm5} for d-stacks, we do impose extra conditions
\eq{ds10eq3} on the 2-morphisms $\bs\eta_{ijk}$ on quadruple
overlaps and \eq{ds10eq4} on the 2-morphisms $\bs\ze_{ij}$ on triple
overlaps. Thus, the d-stack versions of these results are
weaker.\I{d-stack!gluing by equivalences!conditions on overlaps}

The reason for this is that 2-morphisms $\eta:\bs f\Ra\bs g$ of
d-space 1-morphisms $\bs f,\bs g:\bX\ra\bY$ are morphisms
$\eta:\uf^*(\FY)\ra\EX$ in $\qcoh(\uX)$. We can interpolate between
such morphisms using partitions of unity on $\uX$, and in Remark
\ref{ds3rem2} we explained why this enables us to prove $\bs h$ is
independent of $\eta$ in Proposition \ref{ds3prop}, and to do
without overlap conditions on $\eta_{ijk},\ze_{ij}$ in
Theorem~\ref{ds3thm3}.

In contrast, for 2-morphisms $\bs\eta=(\eta,\eta'):\bs f\Ra\bs g$ in
$\dSta$, the $C^\iy$-stack 2-morphisms $\eta:f\Ra g$ are discrete
objects, and we cannot join them using partitions of unity. So $\bs
h$ may depend on $\eta$ in Proposition \ref{ds10prop}, and we need
overlap conditions on the components $\eta_{ijk},\ze_{ij}$ in
$\bs\eta_{ijk},\bs\ze_{ij}$ in Theorem~\ref{ds10thm5}.

If $f,g:\cX\ra\cY$ are 1-morphisms of Deligne--Mumford
$C^\iy$-stacks, we can make extra assumptions on $\cX,\cY$ or $f,g$
which imply that there is at most one 2-morphism $\eta:f\Ra g$, as
in Proposition \ref{ds9prop1} for orbifolds. Such assumptions can
make \eq{ds10eq3} or \eq{ds10eq4} hold automatically, as both sides
of \eq{ds10eq3} or \eq{ds10eq4} are 2-morphisms $f\Ra g$. So, for
instance, if the $C^\iy$-stacks $\cX_i$ are all effective then
\eq{ds10eq3} holds, and if the d-stack $\bcZ$ is a d-space then
\eq{ds10eq4} holds.\I{d-stack!gluing by
equivalences|)}\I{d-stack!equivalence|)}
\label{ds10rem}
\end{rem}

\subsection{Fibre products of d-stacks}
\label{ds104}
\I{d-stack!fibre products|(}

Section \ref{ds33} discussed fibre products of d-spaces. In \cite[\S
9.5]{Joyc6} we generalize this to d-stacks. Here is the analogue of
Theorem~\ref{ds3thm4}:

\begin{thm}{\bf(a)} All fibre products exist in the
$2$-category~$\dSta$.\I{fibre product!of
d-stacks}\I{2-category!fibre products in}
\smallskip

\noindent{\bf(b)} The\/ $2$-functor\/ $F_\dSpa^\dSta:\dSpa\ra\dSta$
preserves fibre products.
\smallskip

\noindent{\bf(c)} Let\/ $g:\cX\ra\cZ$ and\/ $h:\cY\ra\cZ$ be smooth
maps (\/$1$-morphisms) of orbifolds, and write\/
$\bcX=F_\Orb^\dSta(\cX),$ and similarly for\/ $\bcY,\bcZ,\bs g,\bs
h$. If\/ $g,h$ are transverse,\I{orbifold!transverse fibre products}
so that a fibre product $\cX\t_{g,\cZ,h}\cY$ exists in $\Orb,$ then
the fibre product\/ $\bcX\t_{\bs g,\bcZ,\bs h}\bcY$ in $\dSta$ is
equivalent in $\dSta$ to\/ $F_\Orb^\dSta(\cX\t_{g,\cZ,h}\cY)$. If\/
$g,h$ are not transverse then $\bcX\t_{\bs g,\bcZ,\bs h}\bcY$ exists
in $\dSta,$ but is not an orbifold.
\label{ds10thm6}
\end{thm}

As for d-spaces, we prove (a) by explicitly constructing a d-stack
$\bcW=\bcX\t_{\bs g,\bcZ,\bs h}\bcY$ and showing it satisfies the
universal property to be a fibre product in the 2-category $\dSta$.
The proof follows that of Theorem \ref{ds3thm4} closely, inserting
extra terms for 2-morphisms of $C^\iy$-stacks.\I{d-stack!fibre
products|)}

\subsection{Orbifold strata of d-stacks}
\label{ds105}
\I{d-stack!orbifold strata|(}\I{orbifold strata!of d-stacks|(}

Section \ref{ds87} discussed orbifold strata of Deligne--Mumford
$C^\iy$-stacks. In \cite[\S 9.6]{Joyc6} we generalize this to
d-stacks. The next theorems summarize the results.

\begin{thm} Let\/ $\bcX$ be a d-stack, and\/ $\Ga$ a finite group.
Then we can define d-stacks\/
$\bcX^\Ga,\btcX^\Ga,\bhcX^\Ga,$\G[XGaf]{$\bcX^\Ga,\btcX{}^\Ga,\bhcX{}^\Ga,
\bcX{}^\Ga_\ci,\btcX{}^\Ga_\ci,\bhcX{}^\Ga_\ci$}{orbifold strata of
a d-stack $\bcX$} and open d-substacks\/
$\bcX^\Ga_\ci\subseteq\bcX^\Ga,$ $\btcX^\Ga_\ci\subseteq\btcX^\Ga,$
$\bhcX^\Ga_\ci\subseteq \bhcX^\Ga,$ all natural up to
$1$-isomorphism in $\dSta,$ a d-space\/ $\bs{\hat X}{}^\Ga_\ci$
natural up to $1$-isomorphism in $\dSpa,$ and\/ $1$-morphisms $\bs
O^\Ga(\bcX),\bs{\ti\Pi}{}^\Ga (\bcX),\ldots$ fitting into a strictly
commutative diagram in\/~$\dSta\!:$\G[OGaXb]{$\bs
O{}^\Ga(\bcX),\bs{\ti O}{}^\Ga(\bcX),\bs O{}^\Ga_\ci(\bcX),\bs{\ti
O}{}^\Ga_\ci(\bcX)$}{1-morphisms of orbifold strata
$\bcX^\Ga,\ldots,\bhcX{}^\Ga_\ci$ of a d-stack
$\bcX$}\G[PiGaXb]{$\bs{\ti\Pi}{}^\Ga(\bcX),\bs{\hat\Pi}{}^\Ga(\bcX),
\bs{\ti\Pi}{}^\Ga_\ci(\bcX),\bs{\hat\Pi}{}^\Ga_\ci(\bcX)$}{1-morphisms
of orbifold strata $\bcX^\Ga,\ldots,\bhcX{}^\Ga_\ci$ of a d-stack
$\bcX$}
\e
\begin{gathered}
\xymatrix@C=48pt@R=7pt{ \bcX^\Ga_\ci
\ar[rr]^{\bs{\ti\Pi}{}^\Ga_\ci(\bcX)} \ar[dr]_(0.3){\bs
O{}^\Ga_\ci(\bcX)} \ar[dd]_\subset \ar@(ul,l)[]_(0.8){\Aut(\Ga)} &&
\btcX^\Ga_\ci \ar[r]^(0.35){\bs{\hat\Pi}{}^\Ga_\ci(\bcX)}
\ar[dl]^(0.3){\bs{\ti O}{}^\Ga_\ci(\bcX)} \ar[dd]^\subset &
{\bhcX{}^\Ga_\ci\simeq F_\dSpa^\dSta(\bs{\hat
X}{}^\Ga_\ci)\!\!\!\!\!\!\!\!\!\!\!\!}
\ar@<.5ex>[dd]^\subset \\ & \bcX \\
\bcX^\Ga \ar[rr]_{\bs{\ti\Pi}{}^\Ga(\bcX)} \ar[ur]^(0.3){\bs
O^\Ga(\bcX)} \ar@(dl,l)[]^(0.8){\Aut(\Ga)} && \btcX^\Ga
\ar[r]_{\bs{\hat\Pi}{}^\Ga(\bcX)} \ar[ul]_(0.3){\bs{\ti
O}^\Ga(\bcX)}  &
*+[r]{\bhcX^\Ga.} }\!\!\!\!
\end{gathered}
\label{ds10eq5}
\e
We will call\/ $\bcX^\Ga,\btcX^\Ga,\bhcX^\Ga,\bcX^\Ga_\ci,
\btcX^\Ga_\ci,\bhcX^\Ga_\ci,\bs{\hat X}{}^\Ga_\ci$ the
\begin{bfseries}orbifold strata\end{bfseries} of\/~$\bcX$.

The underlying\/ $C^\iy$-stacks of\/ $\bcX^\Ga,\ldots,\bhcX^\Ga_\ci$
are the orbifold strata\/ $\cX^\Ga,\ldots,\hcX^\Ga_\ci$ from\/
{\rm\S\ref{ds87}} of the\/ $C^\iy$-stack\/ $\cX$ in\/ $\bcX$. The
$C^\iy$-stack\/ $1$-morphisms underlying the d-stack\/ $1$-morphisms
in \eq{ds10eq5} are those given in\/~\eq{ds8eq3}.
\label{ds10thm7}
\end{thm}

\begin{thm}{\bf(a)} Let\/ $\bcX,\bcY$ be d-stacks, $\Ga$ a finite
group, and\/ $\bs f:\bcX\ra\bcY$ a representable\/
$1$-morphism\I{d-stack!representable 1-morphism}\I{Deligne--Mumford
$C^\iy$-stack!representable 1-morphism} in $\dSta,$ that is, the
underlying\/ $C^\iy$-stack\/ $1$-morphism $f:\cX\ra\cY$ is
representable. Then there is a unique representable\/ $1$-morphism\/
$\bs f^\Ga:\bcX^\Ga\ra\bcY^\Ga$ in $\dSta$ with\/ $\bs
O^\Ga(\bcY)\ci\bs f^\Ga=\bs f\ci\bs O^\Ga(\bcX)$. Here
$\bcX^\Ga,\bcY^\Ga,\bs O^\Ga(\bcX), \bs O^\Ga(\bcY)$ are as in
Theorem\/~{\rm\ref{ds10thm7}}.
\smallskip

\noindent{\bf(b)} Let\/ $\bs f,\bs g:\bcX\ra\bcY$ be representable\/
$1$-morphisms and\/ $\bs\eta:\bs f\Ra\bs g$ a $2$-morphism in
$\dSta,$ and\/ $\bs f^\Ga,\bs g^\Ga:\bcX^\Ga\ra\bcY^\Ga$ be as in\/
{\bf(a)}. Then there is a unique $2$-morphism $\bs\eta^\Ga:\bs
f^\Ga\Ra\bs g^\Ga$ in $\dSta$ with\/ $\bs\id_{\bs
O^\Ga(\bcY)}*\bs\eta^\Ga=\bs\eta*\bs\id_{\bs O^\Ga(\bcX)}$.
\smallskip

\noindent{\bf(c)} Write\/ $\dSta^{\bf re}$ for the $2$-subcategory
of\/ $\dSta$ with only representable\/ $1$-morphisms. Then mapping
$\bcX\mapsto F^\Ga(\bcX)=\bcX^\Ga$ on objects, $\bs f\mapsto
F^\Ga(\bs f)=\bs f^\Ga$ on (representable)\/ $1$-morphisms, and\/
$\bs\eta\mapsto F^\Ga(\bs\eta)=\bs\eta^\Ga$ on\/ $2$-morphisms
defines a strict\/ $2$-functor\/~$F^\Ga:\dSta^{\bf re}\ra\dSta^{\bf
re}$.
\smallskip

\noindent{\bf(d)} Analogues of\/ {\bf(a)}--{\bf(c)} hold for the
orbifold strata\/ $\btcX^\Ga,$ yielding a strict\/ $2$-functor $\ti
F{}^\Ga:\dSta^{\bf re}\ra\dSta^{\bf re}$. Weaker analogues of\/
{\bf(a)}--{\bf(c)} also hold for the orbifold strata\/ $\bhcX^\Ga$.
In {\bf(a)\rm,} the $1$-morphism\/ $\bs{\hat f}{}^\Ga:\bhcX^\Ga
\ra\bhcY^\Ga$ is natural only up to $2$-isomorphism, and in\/
{\bf(c)} we get a weak\/ $2$-functor\I{2-category!weak 2-functor}
$\hat F{}^\Ga:\dSta^{\bf re}\ra\dSta^{\bf re}$.
\label{ds10thm8}
\end{thm}

Since equivalences in $\dSta$ are automatically representable, and
(strict or weak) 2-functors take equivalences to equivalences, we
deduce:\I{d-stack!equivalence}

\begin{cor} Suppose $\bcX,\bcY$ are equivalent d-stacks, and\/
$\Ga$ is a finite group. Then $\bcX^\Ga$ and\/ $\bcY^\Ga$ are
equivalent in $\dSta,$ and similarly for $\btcX^\Ga,\ab\bhcX^\Ga,
\ab\bcX^\Ga_\ci,\ab\btcX^\Ga_\ci,\ab\bhcX^\Ga_\ci$ and\/
$\btcY^\Ga,\bhcY^\Ga,\bcY^\Ga_\ci,\btcY^\Ga_\ci, \bhcY^\Ga_\ci$.
Also $\bs{\hat X}{}^\Ga_\ci,\bs{\hat Y}{}^\Ga_\ci$ are equivalent
in\/~$\dSpa$.
\label{ds10cor}
\end{cor}

Here are the d-stack analogues of Theorems \ref{ds8thm4}
and~\ref{ds8thm5}:

\begin{thm} Let\/ $\bX$ be a d-space and\/ $G$ a finite group
acting on $\bX$ by $1$-isomorphisms, and write\/ $\bcX=[\bX/G]$ for
the quotient d-stack, from Theorem\/ {\rm\ref{ds10thm2}}. Let\/
$\Ga$ be a finite group. Then there are equivalences of
d-stacks\I{d-stack!quotients $[\bX/G]$}
\ea
\bcX^\Ga&\simeq\coprod_{\begin{subarray}{l}\text{conjugacy classes
$[\rho]$ of injective}\\ \text{group morphisms $\rho:\Ga\ra
G$}\end{subarray}\!\!\!\!\!\!\!\!\!\!\!\!\!\!\!\!\!\!\!
\!\!\!\!\!\!\!\!\!} \bigl[\,\bX^{\rho(\Ga)}/\bigl\{g\in
G:g\rho(\ga)=\rho(\ga)g\;\> \forall\ga\in\Ga\bigr\}\bigr],
\label{ds10eq6}
\allowdisplaybreaks\\
\bcX^\Ga_\ci&\simeq\coprod_{\begin{subarray}{l}\text{conjugacy
classes $[\rho]$ of injective}\\ \text{group morphisms $\rho:\Ga\ra
G$}\end{subarray}\!\!\!\!\!\!\!\!\!\!\!\!\!\!\!\!\!\!\!
\!\!\!\!\!\!\!\!\!} \bigl[\,\bX^{\rho(\Ga)}_\ci/\bigl\{g\in
G:g\rho(\ga)=\rho(\ga)g\;\> \forall\ga\in\Ga\bigr\}\bigr],
\label{ds10eq7}
\allowdisplaybreaks\\
\btcX^\Ga&\simeq\coprod_{\text{conjugacy classes $[\De]$ of
subgroups $\De\subseteq G$ with $\De\cong\Ga$} \!\!\!\!\!\!\!\!\!
\!\!\!\!\!\!\!\!\!\!\!\!\!\!\!\!\!\!\!\!\!\!\!\!\!\!\!\!\!\!\!\!
\!\!\!\!\!\!\!\!\!\!\!\!\!\!\!\!\!\!\!\!\!\!\!}
\bigl[\,\bX^\De/\bigl\{g\in G:\De=g\De g^{-1}\bigr\}\bigr],
\label{ds10eq8}
\allowdisplaybreaks\\
\btcX^\Ga_\ci&\simeq\coprod_{\text{conjugacy classes $[\De]$ of
subgroups $\De\subseteq G$ with $\De\cong\Ga$} \!\!\!\!\!\!\!\!\!
\!\!\!\!\!\!\!\!\!\!\!\!\!\!\!\!\!\!\!\!\!\!\!\!\!\!\!\!\!\!\!\!
\!\!\!\!\!\!\!\!\!\!\!\!\!\!\!\!\!\!\!\!\!\!\!}
\bigl[\,\bX^\De_\ci/\bigl\{g\in G:\De=g\De g^{-1}\bigr\}\bigr],
\label{ds10eq9}
\allowdisplaybreaks\\
\bhcX^\Ga&\simeq\coprod_{\text{conjugacy classes $[\De]$ of
subgroups $\De\subseteq G$ with $\De\cong\Ga$} \!\!\!\!\!\!\!\!\!
\!\!\!\!\!\!\!\!\!\!\!\!\!\!\!\!\!\!\!\!\!\!\!\!\!\!\!\!\!\!\!\!
\!\!\!\!\!\!\!\!\!\!\!\!\!\!\!\!\!\!\!\!\!\!\!}
\bigl[\,\bX^\De\big/\bigl(\{g\in G:\De=g\De
g^{-1}\}/\De\bigr)\bigr],
\label{ds10eq10}
\allowdisplaybreaks\\
\bhcX^\Ga_\ci&\simeq\coprod_{\text{conjugacy classes $[\De]$ of
subgroups $\De\subseteq G$ with $\De\cong\Ga$} \!\!\!\!\!\!\!\!\!
\!\!\!\!\!\!\!\!\!\!\!\!\!\!\!\!\!\!\!\!\!\!\!\!\!\!\!\!\!\!\!\!
\!\!\!\!\!\!\!\!\!\!\!\!\!\!\!\!\!\!\!\!\!\!\!}
\bigl[\,\bX^\De_\ci\big/\bigl(\{g\in G:\De=g\De
g^{-1}\}/\De\bigr)\bigr].
\label{ds10eq11}
\ea
Here for each subgroup\/ $\De\subseteq G,$ we write\/ $\bX^\De$ for
the closed d-subspace in\/ $\bX$ fixed by $\De$ in\/ $G,$ as in\/
{\rm\S\ref{ds34},} and\/ $\bX^\De_\ci$ for the open d-subspace in
$\bX^\De$ of points in $\bX$ whose stabilizer group in $G$ is
exactly\/ $\De$. In {\rm\eq{ds10eq6}--\eq{ds10eq7},} morphisms
$\rho,\rho':\Ga\ra G$ are conjugate if\/ $\rho'=\Ad(g)\ci\rho$ for
some $g\in G,$ and subgroups $\De,\De'\subseteq G$ are conjugate
if\/ $\De=g\De'g^{-1}$ for some $g\in G$. In
\eq{ds10eq6}--\eq{ds10eq11} we sum over one representative $\rho$ or
$\De$ for each conjugacy class.
\label{ds10thm9}
\end{thm}

\begin{thm} Let\/ $\bcX$ be a d-stack and\/ $\Ga$ a finite group,
so that Theorem\/ {\rm\ref{ds10thm7}} gives a d-stack\/ $\bcX^\Ga$
and a $1$-morphism $\bs O{}^\Ga(\bcX):\bcX^\Ga\ra\bcX$. Equation
\eq{ds10eq2} for $\bs O{}^\Ga(\bcX)$ becomes:
\e
\begin{gathered}
\xymatrix@C=11pt@R=20pt{ {\begin{subarray}{l}\ts O^\Ga(\cX)^*
(\EcX)\!=
\\ \ts (\EcX)^\Ga_\tr\!\op\!(\EcX)^\Ga_\nt\end{subarray}}
\ar@<-1.6ex>[rr]_{\raisebox{-10pt}{$\sst O^\Ga(\cX)^*(\phi_\cX)$}}
\ar[d]^(0.67){O^\Ga(\cX)''} &&
{\begin{subarray}{l}\ts O^\Ga(\cX)^*(\FcX)\!= \\
\ts (\FcX)^\Ga_\tr\!\op\!(\FcX)^\Ga_\nt\end{subarray}}
\ar@<-1.6ex>[rr]_{\raisebox{-10pt}{$\sst O^\Ga(\cX)^*(\psi_\cX)$}}
\ar[d]^(0.67){O^\Ga(\cX)^2} &&
{\begin{subarray}{l}\ts O^\Ga(\cX)^*(T^*\cX)\!= \\
\ts (T^*\cX)^\Ga_\tr\!\op\!(T^*\cX)^\Ga_\nt\end{subarray}}
\ar@<-1.6ex>[r] \ar[d]^(0.6){\begin{subarray}{l} O{}^\Ga(\cX)^3=\\
\Om_{O^\Ga(\cX)}\end{subarray}} &
{\raisebox{-20pt}{$0$}} \\
\cE_{\smash{\cX^\Ga}} \ar[rr]^{\phi_{\smash{\cX^\Ga}}} &&
\cF_{\smash{\cX^\Ga}} \ar[rr]^{\psi_{\smash{\cX^\Ga}}} &&
T^*(\cX^\Ga) \ar[r] & {0.\kern -.28em} }\!\!\!\!\!\!
\end{gathered}
\label{ds10eq12}
\e

Then the columns $O^\Ga(\cX)'',$ $O^\Ga(\cX)^2,$ $O^\Ga(\cX)^3$ of\/
\eq{ds10eq12} are isomorphisms when restricted to the `trivial'
summands $(\EcX)^\Ga_\tr,(\FcX)^\Ga_\tr,(T^*\cX)^\Ga_\tr,$ and are
zero when restricted to the `nontrivial' summands
$(\EcX)^\Ga_\nt,(\FcX)^\Ga_\nt,(T^*\cX)^\Ga_\nt$. In particular,
this implies that the virtual cotangent sheaf\/
$\phi_{\smash{\cX^\Ga}}:\cE_{\smash{\cX^\Ga}}\ra
\cF_{\smash{\cX^\Ga}}$ of\/ $\bcX^\Ga$ is $1$-isomorphic in
$\vqcoh(\cX^\Ga)$ to $(\phi_\cX)^\Ga_\tr:(\EcX)^\Ga_\tr
\ra(\FcX)^\Ga_\tr,$ the `trivial' part of the pullback to $\cX^\Ga$
of the virtual cotangent sheaf\/\I{d-stack!virtual cotangent sheaf}
$\phi_\cX:\EcX\ra\FcX$ of\/~$\bcX$.

The analogous results also hold for $\btcX^\Ga,\bhcX^\Ga,
\bcX^\Ga_\ci,\btcX^\Ga_\ci$
and\/~$\bhcX^\Ga_\ci$.\I{d-stack|)}\I{d-stack!orbifold
strata|)}\I{orbifold strata!of d-stacks|)}
\label{ds10thm10}
\end{thm}

\section{The 2-category of d-orbifolds}
\label{ds11}
\I{d-orbifold|(}

In \cite[Chap.~10]{Joyc6} we discuss {\it d-orbifolds}, orbifold
versions of d-manifolds. They are related to {\it Kuranishi
spaces\/}\I{Kuranishi space!and d-orbifolds} (without boundary) in
the work of Fukaya, Oh, Ohta and Ono \cite{FuOn,FOOO} on symplectic
geometry. As we explain briefly in \S\ref{ds16}, and in more detail
in \cite[\S 14.3]{Joyc6}, although Kuranishi spaces are similar to
d-orbifolds in many ways, the theory of Kuranishi spaces in
\cite{FuOn,FOOO} is incomplete --- for instance, there is no notion
of morphism of Kuranishi spaces, so they do not form a category. We
argue in \cite[\S 14.3]{Joyc6} that the `right' way to define
Kuranishi spaces is as d-orbifolds, or d-orbifolds with corners.

\subsection{Definition of d-orbifolds}
\label{ds111}
\I{d-orbifold!definition|(}\I{virtual quasicoherent
sheaf|(}\I{virtual vector bundle|(}

In \S\ref{ds43} we discussed {\it virtual quasicoherent sheaves\/}
and {\it virtual vector bundles\/} on $C^\iy$-schemes $\uX$. The
next remark, drawn from \cite[\S 10.1.1]{Joyc6}, explains how these
generalize to Deligne--Mumford $C^\iy$-stacks~$\cX$.

\begin{rem} In the $C^\iy$-stack analogue of Definition
\ref{ds4def2}, the 2-categories
$\vqcoh(\cX)$\G[vqcoh(X)b]{$\vqcoh(\cX)$}{2-category of virtual
quasicoherent sheaves on a Deligne--Mumford $C^\iy$-stack $\cX$} and
$\vvect(\cX)$\G[vvect(X)b]{$\vvect(\cX)$}{2-category of virtual
vector bundles on a Deligne--Mumford $C^\iy$-stack $\cX$} for a
Deligne--Mumford $C^\iy$-stack $\cX$ are defined exactly as for
$C^\iy$-schemes.\I{virtual quasicoherent sheaf!on Deligne--Mumford
$C^\iy$-stack}\I{virtual vector bundle!on a Deligne--Mumford
$C^\iy$-stack} For $\cX\ne\es$, virtual vector bundles
$(\cE^\bu,\phi)$ have a well-defined {\it rank\/}
$\rank(\cE^\bu,\phi)\in\Z$. If $f:\cX\ra\cY$ is a 1-morphism of
Deligne--Mumford $C^\iy$-stacks then pullback $f^*$ defines strict
2-functors\I{2-category!strict 2-functor} $f^*:\vqcoh(\cY)\ra
\vqcoh(\cX)$ and $f^*:\vvect(\cY)\ra \vvect(\cX)$, as for
$C^\iy$-schemes. If $f,g:\cX\ra\cY$ are 1-morphisms and $\eta:f\Ra
g$ a 2-morphism then $\eta^*:f^*\Ra g^*$ is a 2-natural
transformation.

In the d-stack version of Definition \ref{ds4def3}, we define the
{\it virtual cotangent sheaf\/} $T^*\bcX$\G[T*Xc]{$T^*\bcX$}{virtual
cotangent sheaf of a d-stack $\bcX$}\I{d-stack!virtual cotangent
sheaf} of a d-stack $\bcX$ to be the morphism $\phi_\cX:\EcX\ra\FcX$
in $\qcoh(\cX)$ from Definition \ref{ds10def1}. If $\bs
f:\bcX\ra\bcY$ is a 1-morphism in $\dSta$ then $\Om_{\bs
f}:=(f'',f^2)$ is a 1-morphism $f^*(T^*\bcY)\ra T^*\bcX$ in
$\vqcoh(\cX)$. If $\bs f,\bs g:\bcX\ra\bcY$ are 1-morphisms and
$\bs\eta=(\eta,\eta'):\bs f\Ra\bs g$ is a 2-morphism in $\dSta$,
then we have 1-morphisms $\Om_{\bs f}:f^*(T^*\bcY)\ra T^*\bcX$,
$\Om_{\bs g}:g^*(T^*\bcY)\ra T^*\bcX$, and
$\eta^*(T^*\bcY):f^*(T^*\bcY)\ra g^*(T^*\bcY)$ in $\qcoh(\cX)$, and
$\eta':\Om_{\bs f}\Ra\Om_{\bs g}\ci\eta^*(T^*\bcY)$ is a 2-morphism
in~$\vqcoh(\cX)$.\I{virtual quasicoherent sheaf|)}\I{virtual vector
bundle|)}
\label{ds11rem}
\end{rem}

We can now define d-orbifolds.

\begin{dfn} A d-stack $\bcW$ is called a {\it principal
d-orbifold\/}\I{d-orbifold!principal} if is equivalent in $\dSta$ to
a fibre product $\bcX\t_{\bs g,\bcZ,\bs h}\bcY$ with
$\bcX,\bcY,\bcZ\in\hOrb$. If $\bcW$ is a nonempty principal
d-orbifold then as in Proposition \ref{ds4prop2}, the virtual
cotangent sheaf $T^*\bcW$ is a virtual vector bundle on $\cW$, in
the sense of Remark \ref{ds11rem}. We define the {\it virtual
dimension\/} of $\bcW$ to be $\vdim\bcW=\rank T^*\bcW\in\Z$. If
$\bcW\simeq\bcX\t_{\bcZ}\bcY$ for orbifolds $\bcX,\bcY,\bcZ$
then~$\vdim\bcW=\dim\cX+\dim\cY-\dim\cZ$.

A d-stack $\bcX$ is called a {\it d-orbifold\/} ({\it without
boundary\/}) {\it of virtual dimension\/}\I{d-orbifold!virtual
dimension} $n\in\Z$, written $\vdim\bcX=n$, if $\bcX$ can be covered
by open d-substacks $\bcW$ which are principal d-orbifolds with
$\vdim\bcW=n$. The virtual cotangent sheaf
$T^*\bcX=(\EcX,\FcX,\phi_\cX)$ of $\bcX$ is a virtual vector bundle
of rank $\vdim\bcX=n$, so we call it the {\it virtual cotangent
bundle\/}\I{d-orbifold!virtual cotangent bundle}\I{virtual cotangent
bundle} of~$\bcX$.

Let $\dOrb$\G[dOrb]{$\dOrb$}{2-category of d-orbifolds} be the full
2-subcategory of d-orbifolds in $\dSta$. The 2-functor
$F_\Orb^\dSta:\Orb\ra\dSta$ in Definition \ref{ds10def1} maps into
$\dOrb$, and we will write $F_\Orb^\dOrb=F_\Orb^\dSta:
\Orb\ra\dOrb$. Also $\hOrb$ is a 2-subcategory of $\dOrb$. We say
that a d-orbifold $\bcX$ {\it is an orbifold\/} if it lies in
$\hOrb$.\I{d-orbifold!is an orbifold} The 2-functor $F_\dSpa^\dSta$
maps $\dMan\ra\dOrb$, and we will write $F_\dMan^\dOrb=F_\dSpa^\dSta
\vert_{\dMan}:\dMan\ra\dOrb$. Then~$F_\dMan^\dOrb\ci F_\Man^\dMan=
F_\Orb^\dOrb\ci F_\Man^\Orb:\Man\ra\dOrb$.

Write $\hdMan$\G[dMan'']{$\hdMan$}{2-subcategory of d-orbifolds
equivalent to d-manifolds} for the full 2-subcategory of objects
$\bcX$ in $\dOrb$ equivalent to $F_\dMan^\dOrb(\bX)$ for some
d-manifold $\bX$. When we say that a d-orbifold $\bcX$ {\it is a
d-manifold},\I{d-orbifold!is a d-manifold} we mean
that~$\bcX\in\hdMan$.\I{d-orbifold!definition|)}
\label{ds11def1}
\end{dfn}

The orbifold analogue of Proposition \ref{ds4prop1} holds. Using
Theorem \ref{ds8thm2} we can deduce:

\begin{lem} Let\/ $\bcX$ be a d-orbifold. Then $\bcX$ is a d-manifold,
that is, $\bcX$ is equivalent to\/ $F_\dMan^\dOrb(\bX)$ for some
d-manifold\/ $\bX,$ if and only if\/
$\Iso_\cX([x])\cong\{1\}$\I{C-stack@$C^\iy$-stack!orbifold group
$\Iso_\cX([x])$} for all\/ $[x]$ in\/~$\cX_\top$.
\label{ds11lem}
\end{lem}

\subsection{Local properties of d-orbifolds}
\label{ds112}
\I{d-orbifold!local properties|(}\I{d-orbifold!standard model
a@standard model $\bcS_{\cV,\cE,s}$|(}

Following Examples \ref{ds4ex2} and \ref{ds4ex3}, we define
`standard model' d-orbifolds $\bcS_{\cV,\cE,s}$ and
1-morphisms~$\bcS_{\smash{f,\hat f}}$.

\begin{ex} Let $\cV$ be an orbifold, $\cE\in\vect(\cV)$ a vector
bundle on $\cV$ as in \S\ref{ds86}, and $s\in C^\iy(\cE)$ a smooth
section, that is, $s:\O_\cV\ra\cE$ is a morphism in $\vect(\cV)$. We
will define a principal d-orbifold $\bcS_{\cV,\cE,s}=(\cS,\OcSp,
\EcS,\im_\cS,\jm_\cS)$, which we call a `standard model' d-orbifold.

Let the Deligne--Mumford $C^\iy$-stack $\cS$ be the $C^\iy$-substack
in $\cV$ defined by the equation $s=0$, so that informally
$\cS=s^{-1}(0)\subset\cV$. Explicitly, as in \S\ref{ds8}, a
$C^\iy$-stack $\cV$ consists of a category $\cV$ and a functor
$p_\cV:\cV\ra\CSch$, where there is a 1-1 correspondence between
objects $u$ in $\cV$ with $p_\cV(u)=\uU$ in $\CSch$ and 1-morphisms
$\ti u:\bar\uU\ra\cV$ in $\CSta$. Define $\cS$ to be the full
subcategory of objects $u$ in $\cV$ such that the morphism $\ti
u{}^*(s):\ti u{}^*(\OcV)\ra\ti u{}^*(\cE)$ in $\qcoh(\bar\uU)$ is
zero, and define~$p_\cS=p_\cV\vert_\cS:\cS\ra\CSch$.

Since $i_\cV:\cS\ra\cV$ is the inclusion of a $C^\iy$-substack,
$i_\cV^\sh:i_\cV^{-1}(\OcV)\ra\OcS$ is a surjective morphism of
sheaves of $C^\iy$-rings on $\cS$. Write $\cI_s$ for the kernel of
$i_\cV^\sh$, as a sheaf of ideals in $i_\cV^{-1}(\OcV)$, and
$\cI_s^2$ for the corresponding sheaf of squared ideals, and
$\OcSp=i_\cV^{-1}(\OcV)/\cI_s^2$ for the quotient sheaf of
$C^\iy$-rings, and $\im_\cS:\OcSp\ra\OcS$ for the natural projection
$i_\cV^{-1}(\OcV)/\cI_s^2\twoheadrightarrow
i_\cV^{-1}(\OcV)/\cI_s\cong\OcS$ induced by the
inclusion~$\cI_s^2\subseteq \cI_s$.

Write $\cE^*\in\vect(\cV)$ for the dual vector bundle of $\cE$, and
set $\EcS=i_\cV^*(\cE{}^*)$. There is a natural, surjective morphism
$\jm_\cS:\EcS\ra\cI_\cS=\cI_s/\cI_s^2$ in $\qcoh(\cS)$ which locally
maps $\al+(\cI_s\cdot C^\iy(\cE^*))\mapsto \al\cdot s+\cI_s^2$. Then
$\bcS_{\cV,\cE,s}=(\cS,\OcSp,\EcS,\im_\cS,\jm_\cS)$ is a d-stack. As
in the d-manifold case, we can show that $\bcS_{\cV,\cE,s}$ is
equivalent in $\dSta$ to $\bcV\t_{\bs s,\bcE,\bs 0}\bcV$, where
$\bcV,\bcE,\bs s,\bs 0=F_\Orb^\dSta\bigl(\cV,\Tot(\cE),
\Tot(s),\Tot(0)\bigr)$, using the notation of \S\ref{ds91}. Thus
$\bcS_{\cV,\cE,s}$ is a principal d-orbifold. Every principal
d-orbifold $\bcW$ is equivalent in $\dSta$ to
some~$\bcS_{\cV,\cE,s}$.\G[SVEsc]{$\bcS_{\cV,\cE,s}$}{`standard
model' d-orbifold}

Sometimes it is useful to take $\cV$ to be an {\it effective\/}
orbifold,\I{orbifold!effective} as in~\S\ref{ds91}.
\label{ds11ex1}
\end{ex}

\begin{ex} Let \I{d-orbifold!standard model a@standard model
$\bcS_{\cV,\cE,s}$!1-morphism|(} $\cV,\cW$ be orbifolds, $\cE,\cF$
be vector bundles on $\cV,\cW$, and $s\in C^\iy(\cE)$, $t\in
C^\iy(\cF)$ be smooth sections, so that Example \ref{ds11ex1}
defines `standard model' principal d-orbifolds
$\bcS_{\cV,\cE,s},\bcS_{\cW,\cF,t}$. Write
$\bcS_{\cV,\cE,s}=\bcS=(\cS,\OcSp,\EcS,\im_\cS,\jm_\cS)$ and
$\bcS_{\cW,\cF,t}=\bcT=(\cT,\O_\cT',\cE_\cT,\im_\cT,\jm_\cT)$.
Suppose $f:\cV\ra\cW$ is a 1-morphism, and $\hat f:\cE\ra f^*(\cF)$
is a morphism in $\vect(\cV)$ satisfying
\e
\hat f\ci s=f^*(t).
\label{ds11eq1}
\e
We will define a 1-morphism $\bs g=(g,g',g''):\bcS\ra\bcT$ in
$\dSta$, which we write as $\bcS_{\smash{f,\hat f}}:
\bcS_{\cV,\cE,s}\ra \bcS_{\cW,\cF,t}$, and call a `standard model'
1-morphism.

As in Example \ref{ds11ex1}, $\cV,\cW$ are categories,
$\cS\subseteq\cV$, $\cT\subseteq\cW$ are full subcategories, and
$f:\cV\ra\cW$ is a functor. Using \eq{ds11eq1} one can show that
$f(\cS)\subseteq\cT\subseteq\cW$. Define $g=f\vert_\cS:\cS\ra\cT$.
Then $g:\cS\ra\cT$ is a 1-morphism of Deligne--Mumford
$C^\iy$-stacks, with~$i_\cW\ci g=f\ci i_\cV:\cS\ra\cW$.

To define $g':g^{-1}(\O_\cT')\ra\OcSp$, consider the commutative
diagram:
\begin{equation*}
\xymatrix@C=10pt@R=15pt{
g^{-1}(\cI_t^2) \ar[r] \ar@{.>}[d] &
g^{-1}(i_\cW^{-1}(\OcW)) \ar[rr] \ar[d]^{\begin{subarray}{l}
i_\cV^{-1}(f^\sh)\ci I_{i_\cV,f}(\OcW)\ci\\
I_{g,i_\cW}(\OcW)^{-1}\end{subarray}} &&
g^{-1}(\O_\cT')\!=\!g^{-1}(i_\cW^{-1}(\OcW)/\cI_t^2)
\ar@{.>}[d]^{g'} \ar[r] & 0 \\
\cI_s^2 \ar[r] & i_\cV^{-1}(\OcV)
\ar[rr] && \OcSp=i_\cV^{-1}(\OcV)/\cI_s^2 \ar[r] & {0.\kern -.28em} }
\end{equation*}
The rows are exact. Using \eq{ds11eq1}, we see the central column
maps $g^{-1}(\cI_t)\ra\cI_s$, and so maps $g^{-1}(\cI_t^2)\ra
\cI_s^2$, and the left column exists. Thus by exactness there is a
unique morphism $g'$ making the diagram commute.

We have $\EcS=i_\cV^*(\cE^*)$ and $\cE_\cT=i_\cW^*(\cF^*)$, and
$\hat f:\cE\ra f^*(\cF)$ induces $\hat f{}^*:f^*(\cF{}^*)\ra
\cE{}^*$. Define $g''=i_\cV^*(\hat f{}^*)\ci I_{i_\cV,f}(\cF{}^*)\ci
I_{g,i_\cW}(\cF{}^*)^{-1}:g^*(\cE_\cT)\ra\EcS$ in $\qcoh(\cS)$. Then
$\bs g=(g,g',g''):\bcS\ra\bcT$ is a 1-morphism in $\dSta$, which we
write as~$\bcS_{\smash{f,\hat f}}: \bcS_{\cV,\cE,s}\ra
\bcS_{\cW,\cF,t}$.\G[Sffc]{$\bcS_{\smash{f,\hat
f}}:\bcS_{\cV,\cE,s}\ra \bcS_{\cW,\cF,t}$}{`standard model'
1-morphism in $\dOrb$}

Suppose now that $\ti\cV\subseteq\cV$ is open, with inclusion
1-morphism $i_{\smash{\ti\cV}}:\ti\cV\ra\cV$. Write $\ti{\cal
E}=\cE\vert_{\smash{\ti\cV}}= i_{\smash{\ti\cV}}^*(\cE)$ and $\ti
s=s\vert_{\smash{\ti\cV}}$. Define $\bs
i_{\smash{\ti\cV,\cV}}=\bcS_{i_{\smash{\ti\cV}},\id_{\ti{\cal E}}}:
\bcS_{\smash{\ti\cV,\ti{\cal E},\ti
s}}\ra\bcS_{\cV,\cE,s}$.\G[iVVc]{$\bs
i_{\smash{\ti\cV,\cV}}:\bcS_{\smash{\ti\cV,\ti{\cal E},\ti
s}}\ra\bcS_{\cV,\cE,s}$}{inclusion of open set in `standard model'
d-orbifold} If $s^{-1}(0)\subseteq\ti\cV$ then $\bs
i_{\smash{\ti\cV,\cV}}:\bcS_{\smash{\ti\cV,\ti{\cal E},\ti
s}}\ra\bcS_{\cV,\cE,s}$ is a 1-isomorphism.\I{d-orbifold!standard
model a@standard model $\bcS_{\cV,\cE,s}$!1-morphism|)}
\label{ds11ex2}
\end{ex}

We do not define `standard model' 2-morphisms in $\dOrb$, as in
Example \ref{ds4ex4} for d-manifolds, to avoid inconvenience in
combining the $O(s),O(s^2)$ notation with 2-morphisms of orbifolds.
But see Example \ref{ds11ex5} below for a different form of
`standard model' 2-morphism.

Any d-orbifold $\bcX$ is locally equivalent near a point $[x]$ to a
principal d-orbifold, and so to a standard model d-orbifold
$\bcS_{\cV,\cE,s}$. The next theorem, the analogue of Theorem
\ref{ds4thm1}, shows that $\cV,\cE,s$ are locally determined
essentially uniquely if $\dim\cV$ is chosen to be minimal (which
corresponds to the condition~$\d s(v)=0$).

\begin{thm} Suppose\/ $\bcX$ is a d-orbifold, and\/
$[x]\in\cX_\top$. Then there exists an open neighbourhood\/ $\bcU$
of\/ $[x]$ in $\bcX$ and an equivalence\I{d-orbifold!equivalence}
$\bcU\simeq\bcS_{\cV,\cE,s}$ in $\dOrb$ for $\bcS_{\cV,\cE,s}$ as in
Example\/ {\rm\ref{ds11ex1},} such that the equivalence identifies
$[x]$ with $[v]\in\cV_\top$ with\/ $s(v)=\d s(v)=0$. Furthermore,
$\cV,\cE,s$ are determined up to non-canonical equivalence near
$[v]$ by $\bcX$ near $[x]$. In fact, they depend only on the
$C^\iy$-stack $\cX,$ the point\/ $[x]\in\cX_\top,$ and the
representation of\/ $\Iso_\cX([x])$\I{C-stack@$C^\iy$-stack!orbifold
group $\Iso_\cX([x])$} on the finite-dimensional vector
space\/~$\Ker\bigl(x^*(\phi_\cX):x^*(\EcX)\ra x^*(\FcX)\bigr)$.
\label{ds11thm1}
\end{thm}

In a d-orbifold $\bcX=(\cX,\OcXp,\EcX,\im_\cX,\jm_\cX)$, we think of
$\cX$ as `classical' and $\OcXp,\EcX,\im_\cX,\jm_\cX$ as `derived'.
The extra information in the `derived' data is like a vector bundle
$\cE$ over $\cX$. A vector bundle $\cE$ on a Deligne--Mumford
$C^\iy$-stack $\cX$ is determined locally near $[x]\in\cX_\top$ by
the representation of $\Iso_\cX([x])$ on the fibre $x^*(\cE)$ of
$\cE$ at $[x]$. Thus, it is reasonable that $\bcX$ should be
determined up to equivalence near $[x]$ by $\cX$ and a
representation of~$\Iso_\cX([x])$.\I{d-orbifold!standard model
a@standard model $\bcS_{\cV,\cE,s}$|)}

Here are alternative forms of `standard model' d-orbifolds,
1-morphisms and 2-morphisms, using the quotient d-stack notation
of~\S\ref{ds102}.\I{d-orbifold!standard model b@standard model
$[\bS_{V,E,s}/\Ga]$|(}

\begin{ex} Let $V$ be a manifold, $E\ra V$ a vector bundle, $\Ga$ a
finite group acting smoothly on $V,E$ preserving the vector bundle
structure, and $s:V\ra E$ a smooth, $\Ga$-equivariant section of
$E$. Write the $\Ga$-actions on $V,E$ as $r(\ga):V\ra V$ and $\hat
r(\ga):E\ra r(\ga)^*(E)$ for $\ga\in\Ga$. Then Examples \ref{ds4ex2}
and \ref{ds4ex3} give an explicit principal d-manifold
$\bS_{V,E,s}$, and 1-morphisms $\bS_{\smash{r(\ga),\hat
r(\ga)}}:\bS_{V,E,s}\ra\bS_{V,E,s}$ for $\ga\in\Ga$ which are an
action of $\Ga$ on $\bS_{V,E,s}$. Hence Theorem \ref{ds10thm2}(i)
gives a quotient
d-stack~$[\bS_{V,E,s}/\Ga]$.\G[SVEse]{$[\bS_{V,E,s}/\Ga]$}{alternative
`standard model' d-orbifold}

In fact $[\bS_{V,E,s}/\Ga]\simeq\bcS_{\smash{\ti\cV,\ti{\cal E},\ti
s}}$ for $\ti\cV,\ti{\cal E},\ti s$ defined using $V,E,s,\Ga$, with
$\ti\cV=[\uV/\Ga]$. Thus, $[\bS_{V,E,s}/\Ga]$ is a principal
d-orbifold. But not all principal d-orbifolds $\bcW$ have
$\bcW\simeq[\bS_{V,E,s}/\Ga]$, as not all orbifolds $\cV$ have
$\cV\simeq[\uV/\Ga]$ for some manifold $V$ and finite group~$\Ga$.
\label{ds11ex3}
\end{ex}

\begin{ex} Let $[\bS_{V,E,s}/\Ga]$, $[\bS_{W,F,t}/\De]$ be quotient
d-orbifolds as in Example \ref{ds11ex3}, where $\Ga$ acts on $V,E$
by $q(\ga):V\ra V$ and $\hat q(\ga):E\ra q(\ga)^*(E)$ for
$\ga\in\Ga$, and $\De$ acts on $W,F$ by $r(\de):W\ra W$ and $\hat
r(\de):F\ra r(\de)^*(F)$ for $\de\in\De$. Suppose $f:V\ra W$ is a
smooth map, and $\hat f:E\ra f^*(F)$ is a morphism of vector bundles
on $V$ satisfying $\hat f\ci s= f^*(t)+O(s^2)$, as in \eq{ds4eq2},
and $\rho:\Ga\ra\De$ is a group morphism satisfying $f\ci
q(\ga)=r(\rho(\ga))\ci f:V\ra W$ and $q(\ga)^*(\hat f)\ci\hat
q(\ga)=f^*(\hat r(\rho(\ga)))\ci\hat f:E\ra(f\ci q(\ga))^*(F)$ for
all $\ga\in\Ga$, so that $f,\hat f$ are equivariant under
$\Ga,\De,\rho$. Then Example \ref{ds4ex2} defines a 1-morphism
$\bS_{\smash{f,\hat f}}:\bS_{V,E,s}\ra \bS_{W,F,t}$ in $\dMan$. The
equivariance conditions on $f,\hat f$ imply that $\bS_{\smash{f,\hat
f}}\ci\bS_{\smash{q(\ga),\hat q(\ga)}}= \bS_{\smash{r(\rho(\ga)),
\hat r(\rho(\ga))}}\ci\bS_{\smash{f,\hat f}}$ for $\ga\in\Ga$. Hence
Theorem \ref{ds10thm2}(ii) gives a quotient 1-morphism~$[\bS_{
\smash{f,\hat f}},\rho]:[\bS_{V,E,s}/\Ga]\ra
[\bS_{W,F,t}/\De]$.\G[Sffe]{$[\bS_{ \smash{f,\hat
f}},\rho]:[\bS_{V,E,s}/\Ga]\ra [\bS_{W,F,t}/\De]$}{`standard model'
1-morphism in $\dOrb$}\I{d-orbifold!standard model b@standard model
$[\bS_{V,E,s}/\Ga]$!1-morphism}
\label{ds11ex4}
\end{ex}

\begin{ex} Suppose $[\bS_{\smash{f,\hat f}},\rho],
[\bS_{\smash{g,\hat g}},\si] :[\bS_{V,E,s}/\Ga]\ra
[\bS_{W,F,t}/\De]$ are two 1-morphisms as in Example \ref{ds11ex4},
and write $q,\hat q$ for the actions of $\Ga$ on $V,E$ and $r,\hat
r$ for the actions of $\De$ on $W,F$. Then $\rho,\si:\Ga\ra\De$ are
group morphisms. Suppose $\de\in\De$ satisfies
$\si(\ga)=\de\,\rho(\ga)\,\de^{-1}$ for all $\ga\in\Ga$, and
$\La:E\ra f^*(TW)$ is a morphism of vector bundles on $V$ which
satisfies
\ea
r(\de^{-1})\ci g= f+\La\cdot s+O(s^2)\;\>\text{and}\;\> g^*(\hat
r(\de^{-1}))\ci\hat g=\hat f+\La\cdot\d t+O(s),
\label{ds11eq2}\\
f^*(\d r(\rho(\ga)))\ci\La=q(\ga)^*(\La)\ci\hat q(\ga):E\longra
(f\ci q(\ga))^*(TW),\quad \forall\ga\in\Ga,
\label{ds11eq3}
\ea
where $\d r(\rho(\ga)):TW\ra r(\rho(\ga))^*(TW)$ is the derivative
of $r(\rho(\ga))$. Here \eq{ds11eq2} is the conditions for Example
\ref{ds4ex4} to define a `standard model' 2-morphism
$\bS_\La:\bS_{\smash{f,\hat f}}\Ra\bS_{\smash{r(\de^{-1})\ci g,
g^*(\hat r(\de^{-1}))\ci\hat g }}=\bS_{\smash{r(\de^{-1}), \hat
r(\de^{-1})}}\ci\bS_{\smash{g,\hat g}}$ in $\dMan$. Then
\eq{ds11eq3} implies that $\bS_\La*\id_{ \bS_{\smash{q(\ga),\hat
q(\ga)}}}=\id_{\bS_{\smash{r(\rho(\ga)), \hat
r(\rho(\ga))}}}*\bS_\La$ for all $\ga\in\Ga$. Hence Theorem
\ref{ds10thm2}(iii) gives a quotient 2-morphism~$[\bS_\La,\de]:
[\bS_{\smash{f,\hat f}},\rho]\Ra[\bS_{\smash{g,\hat
g}},\si]$.\G[SLa']{$[\bS_\La,\de]: [\bS_{\smash{f,\hat
f}},\rho]\Ra[\bS_{\smash{g,\hat g}},\si]$}{`standard model'
2-morphism in $\dOrb$}\I{d-orbifold!standard model b@standard model
$[\bS_{V,E,s}/\Ga]$!2-morphism}
\label{ds11ex5}
\end{ex}

Here is an analogue of Theorem \ref{ds11thm1} for the alternative
form~$[\bS_{V,E,s}/\Ga]$.

\begin{prop} A d-stack\/ $\bcX$ is a d-orbifold of virtual
dimension\/ $n\in\Z$ if and only if each\/ $[x]\in\cX_\top$ has an
open neighbourhood\/ $\bcU$ equivalent to some $[\bS_{V,E,s}/\Ga]$
in Example\/ {\rm\ref{ds11ex3}} with\/ $\dim V-\rank E=n,$ where
$\Ga=\Iso_\cX([x])$ and\/ $[x]\in\cX_\top$ is identified with a
fixed point\/ $v$ of\/ $\Ga$ in\/ $V$ with\/ $s(v)=0$ and\/ $\d
s(v)=0$. Furthermore, $V,E,s,\Ga$ are determined up to non-canonical
isomorphism near $v$ by $\bcX$ near~$[x]$.\I{d-orbifold!local
properties|)}\I{d-orbifold!standard model b@standard model
$[\bS_{V,E,s}/\Ga]$|)}
\label{ds11prop1}
\end{prop}

\subsection{Equivalences in $\dOrb$, and gluing by equivalences}
\label{ds113}
\I{d-orbifold!equivalence|(}\I{d-orbifold!gluing by equivalences|(}

Next we summarize the results of \cite[\S 10.2]{Joyc6}, the analogue
of \S\ref{ds44}. Section \ref{ds102} discussed \'etale 1-morphisms
in $\dSta$. We characterize when 1-morphisms $\bs f:\bcX\ra\bcY$ and
$\bcS_{\smash{f,\hat f}}:\bcS_{\cV,\cE,s}\ra \bcS_{\cW,\cF,t}$ in
$\dOrb$ are \'etale, or equivalences.\I{d-orbifold!etale
1-morphism@\'etale 1-morphism|(}

\begin{thm} Suppose $\bs f:\bcX\ra\bcY$ is a $1$-morphism of
d-orbifolds, and\/ $f:\cX\ra\cY$ is representable. Then the
following are equivalent:
\begin{itemize}
\setlength{\itemsep}{0pt}
\setlength{\parsep}{0pt}
\item[{\rm(i)}] $\bs f$ is \'etale;
\item[{\rm(ii)}] $\Om_{\bs f}:f^*(T^*\bcY)\ra T^*\bcX$ is an
equivalence in $\vqcoh(\cX);$ and
\item[{\rm(iii)}] The following is a split short
exact sequence\I{split short exact sequence}\I{abelian
category!split short exact sequence} in\/~$\qcoh(\cX)\!:$
\e
\smash{{}\!\!\!\!\xymatrix@C=20pt{ 0 \ar[r] & f^*(\EcY)
\ar[rr]^(0.45){f'' \op -f^*(\phi_\cY)} && \EcX\op f^*(\FcY)
\ar[rr]^(0.6){\phi_\cX\op f^2} && \FcX \ar[r] & 0.}}
\label{ds11eq4}
\e
\end{itemize}
If in addition $f_*:\Iso_\cX([x])\ra\Iso_\cY(f_\top([x]))$ is an
isomorphism for all\/ $[x]\in\cX_\top,$ and\/
$f_\top:\cX_\top\ra\cY_\top$ is a bijection, then $\bs f$ is an
equivalence in\/~$\dOrb$.
\label{ds11thm2}
\end{thm}

\begin{thm} Suppose $\bcS_{\smash{f,\hat f}}:\bcS_{\cV,\cE,s}\ra
\bcS_{\cW,\cF,t}$\I{d-orbifold!standard model a@standard model
$\bcS_{\cV,\cE,s}$!1-morphism} is a `standard model'\/ $1$-morphism,
in the notation of Examples\/ {\rm\ref{ds11ex1}} and\/
{\rm\ref{ds11ex2},} with\/ $f:\cV\ra\cW$ representable. Then
$\bcS_{\smash{f,\hat f}}$ is \'etale if and only if for each\/
$[v]\in\cV_\top$ with\/ $s(v)=0$ and\/ $[w]=f_\top([v])\in\cW_\top,$
the following sequence of vector spaces is exact:
\begin{equation*}
\smash{\xymatrix@C=18pt{ 0 \ar[r] & T_v\cV \ar[rrr]^(0.42){\d
s(v)\op \,\d f(v)} &&& \cE_v\op T_w\cW \ar[rrr]^(0.57){\hat
f(v)\op\, -\d t(w)} &&& \cF_w \ar[r] & 0.}}
\end{equation*}
Also $\bcS_{\smash{f,\hat f}}$ is an equivalence if and only if in
addition\/ $f_\top\vert_{s^{-1}(0)}:s^{-1}(0)\!\ra\! t^{-1}(0)$ is a
bijection, where $s^{-1}(0)\!=\!\{[v]\in\cV_\top:s(v)\!=\!0\},$
$t^{-1}(0)\!=\!\{[w]\in\cW_\top:t(w)\!=\!0\},$ and\/
$f_*:\Iso_\cV([v])\ra\Iso_\cW(f_\top([v]))$ is an isomorphism for
all\/~$[v]\in s^{-1}(0)\subseteq\cV_\top$.\I{d-orbifold!etale
1-morphism@\'etale 1-morphism|)}\I{d-orbifold!equivalence|)}
\label{ds11thm3}
\end{thm}

Here is an analogue of Theorem \ref{ds4thm5} for d-orbifolds, taken
from \cite[\S 10.2]{Joyc6}. It is proved by applying Theorem
\ref{ds10thm5} to glue together the `standard model' d-orbifolds
$\bcS_{\cV_i,\cE_i,s_i}$ by equivalences. Now Theorem \ref{ds10thm5}
includes extra conditions \eq{ds10eq3}--\eq{ds10eq4} on the
2-morphisms $\eta_{ijk},\ze_{jk}$. But by taking the $\cV_i,\cY$ to
be effective orbifolds and the $g_i$ to be submersions, the
$\eta_{ijk},\ze_{jk}$ are unique by Proposition \ref{ds9prop1}, and
so \eq{ds10eq3}--\eq{ds10eq4} hold
automatically.\I{d-orbifold!standard model a@standard model
$\bcS_{\cV,\cE,s}$|(}

\begin{thm} Suppose we are given the following data:
\begin{itemize}
\setlength{\itemsep}{0pt}
\setlength{\parsep}{0pt}
\item[{\rm(a)}] an integer $n;$
\item[{\rm(b)}] a Hausdorff, second countable topological space $X;$
\item[{\rm(c)}] an indexing set\/ $I,$ and a total order $<$ on $I;$
\item[{\rm(d)}] for each\/ $i$ in $I,$ an effective orbifold\/
$\cV_i$ in the sense of Definition\/ {\rm\ref{ds9def3},} a
vector bundle $\cE_i$ on $\cV_i$ with\/
$\dim\cV_i-\rank\cE_i=n,$ a section $s_i\in C^\iy(\cE_i),$ and a
homeomorphism $\psi_i:s_i^{-1}(0) \ra\hat X_i,$ where
$s_i^{-1}(0)=\{[v_i]\in\cV_{i,\top}:s_i(v_i)=0\}$ and\/ $\hat
X_i\subseteq X$ is open; and
\item[{\rm(e)}] for all\/ $i<j$ in $I,$ an open suborbifold\/
$\cV_{ij}\subseteq\cV_i,$ a $1$-morphism
$e_{ij}:\cV_{ij}\ra\cV_j,$ and a morphism of vector bundles
$\hat e_{ij}:\cE_i\vert_{\cV_{ij}}\ra e_{ij}^*(\cE_j)$.
\end{itemize}
Let this data satisfy the conditions:
\begin{itemize}
\setlength{\itemsep}{0pt}
\setlength{\parsep}{0pt}
\item[{\rm(i)}] $X=\bigcup_{i\in I}\hat X_i;$
\item[{\rm(ii)}] if\/ $i<j$ in $I$ then $(e_{ij})_*:
\Iso_{\cV_{ij}}([v])\ra\Iso_{\cV_j}(e_{ij,\top}([v]))$ is an
isomorphism for all\/ $[v]\in\cV_{ij,\top},$ and\/ $\hat
e_{ij}\ci s_i\vert_{\cV_{ij}}=e_{ij}^*(s_j)\ci\io_{ij}$ where
$\io_{ij}:\O_{\cV_{ij}}\ra e_{ij}^*(\O_{\cV_j})$ is the natural
isomorphism, and\/ $\psi_i(s_i\vert_{\cV_{ij}}^{-1}(0))=\hat
X_i\cap\hat X_j,$ and\/ $\psi_i\vert_{s_i\vert_{\cV_{ij}}^{-1}
(0)}=\psi_j\ci e_{ij,\top}\vert_{s_i\vert_{\cV_{ij}}^{-1}(0)},$
and if\/ $[v_i]\in\cV_{ij,\top}$ with\/ $s_i(v_i)=0$ and\/
$[v_j]=e_{ij,\top}([v_i])$ then the following sequence is exact:
\begin{equation*}
\smash{\xymatrix@C=19pt{ 0 \ar[r] & T_{v_i}\cV_i \ar[rrr]^(0.42){\d
s_i(v_i)\op \,\d e_{ij}(v_i)} &&& \cE_i\vert_{v_i}\!\op\!
T_{v_j}\cV_j \ar[rrr]^(0.57){\hat e_{ij}(v_i)\op\, -\d s_j(v_j)} &&&
\cE_j\vert_{v_j} \ar[r] & 0;}}\!\!\!
\end{equation*}
\item[{\rm(iii)}] if\/ $i<j<k$ in $I$ then there exists a
$2$-morphism $\eta_{ijk}:e_{jk}\ci e_{ij}\vert_{\cV_{ik}\cap
e_{ij}^{-1}(\cV_{jk})}\ab\Ra e_{ik}\vert_{\cV_{ik}\cap
e_{ij}^{-1}(\cV_{jk})}$ in $\Orb$ with
\begin{equation*}
{}\!\!\!\!\!\!\!\!\!\!\!
\hat e_{ik}\vert_{\cV_{ik}\cap e_{ij}^{-1}(\cV_{jk})}\!=\!
\eta_{ijk}^*(\cE_k)\!\ci\! I_{e_{ij},e_{jk}}(\cE_k)^{-1}\!\ci\!
e_{ij}\vert_{\cV_{ik}\cap e_{ij}^{-1}(\cV_{jk})}^*(\hat e_{jk})\!\ci\!\hat
e_{ij}\vert_{\cV_{ik}\cap e_{ij}^{-1}(\cV_{jk})}.
\end{equation*}
Note that\/ $\eta_{ijk}$ is unique by
Proposition\/~{\rm\ref{ds9prop1}}.
\end{itemize}

Then there exist a d-orbifold\/ $\bcX$ with\/ $\vdim\bcX=n$ and
underlying topological space $\cX_\top\cong X,$ and a $1$-morphism
$\bs\psi_i:\bcS_{\cV_i,\cE_i,s_i}\ra\bcX$ with underlying continuous
map $\psi_i$ which is an equivalence with the open d-suborbifold\/
$\bs{\hcX}_i\subseteq\bcX$ corresponding to $\hat X_i\subseteq X$
for all\/ $i\in I,$ such that for all\/ $i<j$ in $I$ there exists a
$2$-morphism\/ $\bs\eta_{ij}:\bs\psi_j\ci\bcS_{e_{ij},\hat
e_{ij}}\Ra\bs\psi_i\ci\bs i_{\cV_{ij},\cV_i},$ where
$\bcS_{e_{ij},\hat e_{ij}}:
\bcS_{\cV_{ij},\cE_i\vert_{\cV_{ij}},s_i\vert_{\cV_{ij}}}\ra
\bcS_{\cV_j,\cE_j,s_j}$ and\/ $\bs i_{\cV_{ij},\cV_i}:
\bcS_{\cV_{ij},\cE_i\vert_{\cV_{ij}},s_i\vert_{\cV_{ij}}}\ra
\bcS_{\cV_i,\cE_i,s_i},$ using the notation of Examples\/
{\rm\ref{ds11ex1}} and\/ {\rm\ref{ds11ex2}}. This d-orbifold\/
$\bcX$ is unique up to equivalence in~$\dOrb$.

Suppose also that\/ $\cY$ is an effective orbifold, and\/
$g_i:\cV_i\ra\cY$ are submersions for all\/ $i\in I,$ and there
are\/ $2$-morphisms $\ze_{ij}:g_j\ci e_{ij}\Ra g_i\vert_{\cV_{ij}}$
in $\Orb$ for all\/ $i<j$ in $I$. Then there exist a $1$-morphism
$\bs h:\bcX\ra\bcY$ in $\dOrb$ unique up to $2$-isomorphism, where
$\bcY=F_\Orb^\dOrb(\cY)=\bcS_{\cY,0,0},$ and\/ $2$-morphisms
$\bs\ze_i:\bs h\ci\bs\psi_i\Ra\bcS_{g_i,0}$ for all\/~$i\in
I$.\I{d-orbifold!standard model a@standard model
$\bcS_{\cV,\cE,s}$|)}
\label{ds11thm4}
\end{thm}

Here is another version of the same result using the alternative
form of `standard model' d-orbifolds
in~\S\ref{ds111}.\I{d-orbifold!standard model b@standard model
$[\bS_{V,E,s}/\Ga]$|(}

\begin{thm} Suppose we are given the following data:
\begin{itemize}
\setlength{\itemsep}{0pt}
\setlength{\parsep}{0pt}
\item[{\rm(a)}] an integer $n;$
\item[{\rm(b)}] a Hausdorff, second countable topological space $X;$
\item[{\rm(c)}] an indexing set\/ $I,$ and a total order $<$ on $I;$
\item[{\rm(d)}] for each\/ $i$ in $I,$ a manifold\/ $V_i,$ a vector
bundle $E_i\ra V_i$ with\/ $\dim V_i-\rank E_i=n,$ a finite
group\/ $\Ga_i,$ smooth, locally effective actions\/
$r_i(\ga):V_i\ra V_i,$ $\hat r_i(\ga):E_i\ra r(\ga)^*(E_i)$ of\/
$\Ga_i$ on $V_i,E_i$ for $\ga\in\Ga_i,$ a smooth,
$\Ga_i$-equivariant section $s_i:V_i\ra E_i,$ and a
homeomorphism $\psi_i:X_i\ra\hat X_i,$ where $X_i=\{v_i\in
V_i:s_i(v_i)=0\}/\Ga_i$ and\/ $\hat X_i\subseteq X$ is an open
set; and
\item[{\rm(e)}] for all\/ $i<j$ in $I,$ an open submanifold\/
$V_{ij}\subseteq V_i,$ invariant under $\Ga_i,$ a group morphism
$\rho_{ij}:\Ga_i\ra\Ga_j$, a smooth map $e_{ij}:V_{ij}\ra V_j,$
and a morphism of vector bundles $\hat
e_{ij}:E_i\vert_{V_{ij}}\ra e_{ij}^*(E_j)$.
\end{itemize}

Let this data satisfy the conditions:
\begin{itemize}
\setlength{\itemsep}{0pt}
\setlength{\parsep}{0pt}
\item[{\rm(i)}] $X=\bigcup_{i\in I}\hat X_i;$
\item[{\rm(ii)}] if\/ $i<j$ in $I$ then $\hat e_{ij}\ci s_i
\vert_{V_{ij}}= e_{ij}^*(s_j)+O(s_i^2),$ and for all\/
$\ga\in\Ga$ we have
\begin{align*}
e_{ij}\ci r_i(\ga)&= r_j(\rho_{ij}(\ga))\ci e_{ij}:V_{ij}
\longra V_j,\\
r_i(\ga)^*(\hat e_{ij})\ci\hat r_i(\ga)&= e_{ij}^*(\hat
r_j(\rho_{ij}(\ga)))\ci\hat e_{ij}:E_i\vert_{V_{ij}}
\longra(e_{ij}\ci r_i(\ga))^*(E_j),
\end{align*}
and\/ $\psi_i(X_i\cap (V_{ij}/\Ga_i))=\hat X_i\cap\hat X_j,$
and\/ $\psi_i\vert_{X_i\cap V_{ij}/\Ga_i}=\psi_j\ci
(e_{ij})_*\vert_{X_i\cap V_{ij}/\Ga_j},$ and if\/ $v_i\in
V_{ij}$ with\/ $s_i(v_i)=0$ and\/ $v_j=e_{ij}(v_i)$ then\/
$\rho\vert_{\Stab_{\Ga_i}(v_i)}:\Stab_{\Ga_i}(v_i)\ra
\Stab_{\Ga_j}(v_j)$ is an isomorphism, and the following
sequence of vector spaces is exact:
\begin{equation*}
\smash{\xymatrix@C=19pt{ 0 \ar[r] & T_{v_i}V_i \ar[rrr]^(0.42){\d
s_i(v_i)\op \,\d e_{ij}(v_i)} &&& E_i\vert_{v_i}\!\op\! T_{v_j}V_j
\ar[rrr]^(0.57){\hat e_{ij}(v_i)\op\, -\d s_j(v_j)} &&&
E_j\vert_{v_j} \ar[r] & 0;}}
\end{equation*}
\item[{\rm(iii)}] if\/ $i<j<k$ in $I$ then there exists\/
$\ga_{ijk}\in\Ga_k$ satisfying
\begin{align*}
\rho_{ik}(\ga)&=
\ga_{ijk}\,\ab\rho_{jk}(\rho_{ij}(\ga))\ab\,\ga_{ijk}^{-1}
\quad\text{for all\/ $\ga\in\Ga_i,$}\\
e_{ik}\vert_{V_{ik}\cap e_{ij}^{-1}(V_{jk})}&= r_k(\ga_{ijk})\ci
e_{jk}\ci e_{ij}\vert_{V_{ik}\cap e_{ij}^{-1}(V_{jk})},\quad\text{and}\\
\hat e_{ik}\vert_{V_{ik}\cap e_{ij}^{-1}(V_{jk})}
&= \bigl(e_{ij}^*(e_{jk}^*(\hat r_k(\ga_{ijk})))\ci e_{ij}^*
(\hat e_{jk})\ci\hat e_{ij}\bigr)\vert_{V_{ik}\cap e_{ij}^{-1}(V_{jk})}.
\end{align*}
\end{itemize}

Then there exist a d-orbifold\/ $\bcX$ with\/ $\vdim\bcX=n$ and
underlying topological space $\cX_\top\cong X,$ and a $1$-morphism
$\bs\psi_i:[\bS_{V_i,E_i,s_i}/\Ga_i]\ra\bcX$ with underlying
continuous map $\psi_i$ which is an equivalence with the open
d-suborbifold\/ $\bhcX_i\subseteq\bcX$ corresponding to $\hat
X_i\subseteq X$ for all\/ $i\in I,$ such that for all\/ $i<j$ in $I$
there exists a $2$-morphism\/
$\bs\eta_{ij}:\bs\psi_j\ci[\bS_{e_{ij},\hat
e_{ij}},\rho_{ij}]\Ra\bs\psi_i\ci[\bs i_{V_{ij},V_i},\id_{\Ga_i}],$
where $[\bS_{V_i,E_i,s_i}/\Ga_i]$ is as in Example\/
{\rm\ref{ds11ex3},} and\/ $[\bS_{e_{ij},\hat e_{ij}},\rho_{ij}]:
[\bS_{V_{ij},E_i\vert_{V_{ij}},s_i\vert_{V_{ij}}}/\Ga_i]\ra
[\bS_{V_j,E_j,s_j}/\Ga_j]$ and\/ $[\bs
i_{V_{ij},V_i},\id_{\Ga_i}]:[\bS_{V_{ij},
E_i\vert_{V_{ij}},s_i\vert_{V_{ij}}}/\Ga_i]\ra[\bS_{V_i,E_i,s_i}/\Ga_j]$
as in Example\/ {\rm\ref{ds11ex4}}. This d-orbifold\/ $\bcX$ is
unique up to equivalence in~$\dOrb$.

Suppose also that\/ $Y$ is a manifold, and\/ $g_i:V_i\ra Y$ are
smooth maps for all\/ $i\in I$ with\/ $g_i\ci r_i(\ga)= g_i$ for
all\/ $\ga\in\Ga_i,$ and\/ $g_j\ci e_{ij}=g_i\vert_{V_{ij}}$ for
all\/ $i<j$ in $I$. Then there exist a $1$-morphism $\bs
h:\bcX\ra\bcY$ unique up to $2$-isomorphism, where
$\bcY=F_\Man^\dOrb(Y)=[\bS_{Y,0,0}/\{1\}],$ and\/ $2$-morphisms
$\bs\ze_i:\bs h\ci\bs\psi_i\Ra[\bS_{g_i,0},\pi_{\{1\}}]$ for all\/
$i\in I$. Here $[\bS_{Y,0,0}/\{1\}]$ is from Example\/
{\rm\ref{ds11ex3}} with\/ $E,s$ both zero and\/ $\Ga=\{1\}$, and\/
$[\bS_{g_i,0},\pi_{\{1\}}]:[\bS_{V_i,E_i,s_i}/\Ga_i]\ra
[\bS_{Y,0,0}/\{1\}]=\bcY$ is from Example\/ {\rm\ref{ds11ex4}}
with\/ $\hat g_i=0$
and\/~$\rho=\pi_{\{1\}}:\Ga_i\ra\{1\}$.\I{d-orbifold!standard model
b@standard model $[\bS_{V,E,s}/\Ga]$|)}
\label{ds11thm5}
\end{thm}

The importance of Theorems \ref{ds11thm4} and \ref{ds11thm5} is that
all the ingredients are described wholly in differential-geometric
or topological terms. So we can use these theorems as tools to prove
the existence of d-orbifold structures on spaces coming from other
areas of geometry, such as moduli spaces of $J$-holomorphic
curves.\I{moduli space!of J-holomorphic curves@of $J$-holomorphic
curves} The theorems are used to define functors to d-orbifolds from
other geometric structures, as discussed
in~\S\ref{ds16}.\I{d-orbifold!gluing by equivalences|)}

\subsection{Submersions, immersions, and embeddings}
\label{ds114}

Section \ref{ds45} discussed (w-)submersions, (w-)immersions, and
(w-)embeddings for d-manifolds. Following \cite[\S 10.3]{Joyc6},
here are the analogues for d-orbifolds.

\begin{dfn} Let $\cX$ be a Deligne--Mumford $C^\iy$-stack, so that
as in Remark \ref{ds11rem} we have a 2-category $\vvect(\cX)$ of
virtual vector bundles on $\cX$. We define when a 1-morphism
$f^\bu:(\cE^\bu,\phi)\ra(\cF^\bu,\psi)$ in $\vvect(\cX)$ is {\it
weakly injective}, {\it injective}, {\it weakly surjective\/} or
{\it surjective\/} exactly as in Definition~\ref{ds4def6}.\I{virtual
vector bundle!weakly injective 1-morphism} \I{virtual vector
bundle!injective 1-morphism} \I{virtual vector bundle!weakly
surjective 1-morphism} \I{virtual vector bundle!surjective
1-morphism}\I{virtual vector bundle!on a Deligne--Mumford
$C^\iy$-stack}

Let $\bs f:\bcX\ra\bcY$ be a 1-morphism of d-orbifolds. Then
$\Om_{\bs f}:f^*(T^*\bcY)\ra T^*\bcX$ is a 1-morphism in
$\vvect(\cX)$.\I{d-orbifold!w-submersion|(}\I{d-orbifold!submersion|(}
\I{d-orbifold!w-immersion|(}\I{d-orbifold!immersion|(}
\I{d-orbifold!w-embedding|(}\I{d-orbifold!embedding|(}
\begin{itemize}
\setlength{\itemsep}{0pt}
\setlength{\parsep}{0pt}
\item[(a)] We call $\bs f$ a {\it w-submersion\/} if $\Om_{\bs
f}$ is weakly injective.
\item[(b)] We call $\bs f$ a {\it submersion\/} if $\Om_{\bs
f}$ is injective.
\item[(c)] We call $\bs f$ a {\it w-immersion\/} if
$f:\cX\!\ra\!\cY$ is representable, i.e.\ $f_*:\Iso_\cX([x])
\!\ra\!\Iso_\cY(f_\top([x]))$ is injective for all
$[x]\in\cX_\top$, and $\Om_{\bs f}$ is weakly surjective.
\item[(d)] We call $\bs f$ an {\it immersion\/} if
$f:\cX\ra\cY$ is representable and $\Om_{\bs f}$ is surjective.
\item[(e)] We call $\bs f$ a {\it w-embedding\/} or {\it
embedding\/} if it is a w-immersion or immersion, respectively,
and $f_*:\Iso_\cX([x])\ra\Iso_\cY(f_\top([x]))$ is an
isomorphism for all $[x]\in\cX_\top$, and
$f_\top:\cX_\top\ra\cY_\top$ is a homeomorphism with its image,
so in particular $f_\top$ is injective.
\end{itemize}

Parts (c)--(e) enable us to define {\it
d-suborbifolds\/}\I{d-orbifold!d-suborbifold} of d-orbifolds. {\it
Open d-suborbifolds\/} are (Zariski) open d-substacks of a
d-orbifold. For more general d-suborbifolds, we call $\bs
i:\bcX\ra\bcY$ a {\it w-immersed d-suborbifold}, or {\it immersed
d-suborbifold}, or {\it w-embedded d-suborbifold}, or {\it embedded
d-suborbifold\/} of $\bcY$, if $\bcX,\bcY$ are d-orbifolds and $\bs
i$ is a w-immersion, \ldots, embedding, respectively.
\label{ds11def2}
\end{dfn}

Theorem \ref{ds4thm6} in \S\ref{ds45} holds with orbifolds and
d-orbifolds in place of manifolds and d-manifolds, except part (v),
when we need also to assume $f:\cX\ra\cY$
representable\I{Deligne--Mumford $C^\iy$-stack!representable
1-morphism} to deduce $\bs f$ is \'etale,\I{d-orbifold!etale
1-morphism@\'etale 1-morphism} and part (x), which is false for
d-orbifolds (in the Zariski topology,\I{Zariski topology} at
least).\I{d-orbifold!w-submersion|)}\I{d-orbifold!submersion|)}
\I{d-orbifold!w-immersion|)}\I{d-orbifold!immersion|)}
\I{d-orbifold!w-embedding|)}\I{d-orbifold!embedding|)}

\subsection{D-transversality and fibre products}
\label{ds115}
\I{d-transversality|(}\I{d-orbifold!d-transverse 1-morphisms|(}
\I{d-orbifold!fibre products|(}

Section \ref{ds46} discussed d-transversality and fibre products for
d-manifolds. This is extended to d-orbifolds in \cite[\S
10.4]{Joyc6}, with little essential change. Here are the analogues
of Definition \ref{ds4def8} and
Theorems~\ref{ds4thm7}--\ref{ds4thm10}.

\begin{dfn} Let $\bcX,\bcY,\bcZ$ be d-orbifolds and $\bs
g:\bcX\ra\bcZ,$ $\bs h:\bcY\ra\bcZ$ be 1-morphisms. Let
$\cW=\cX\t_{g,\cZ,h}\cY$ be the $C^\iy$-stack fibre product, and
write $e:\cW\ra\cX$, $f:\cW\ra\cY$ for the projection 1-morphisms,
and $\eta:g\ci e\Ra h\ci f$ for the 2-morphism from the fibre
product. Consider the morphism
\begin{equation*}
\al=\begin{pmatrix} e^*(g'')\ci
I_{e,g}(\EcZ) \\ -f^*(h'')\ci I_{f,h}(\EcZ)\ci\eta^*(\EcZ) \\
(g\ci e)^*(\phi_\cZ)\end{pmatrix} :\begin{aligned}[t]
&(g\ci e)^*(\EcZ)\longra \\
&\quad e^*(\EcX)\op f^*(\EcY)\op (g\ci e)^*(\FcZ)\end{aligned}
\end{equation*}
in $\qcoh(\cW)$. We call $\bs g,\bs h$ {\it d-transverse\/} if $\al$
has a left inverse.
\label{ds11def3}
\end{dfn}

\begin{thm} Suppose\/ $\bcX,\bcY,\bcZ$ are d-orbifolds and\/ $\bs
g:\bcX\ra\bcZ,$ $\bs h:\bcY\ra\bcZ$ are d-transverse $1$-morphisms,
and let\/ $\bcW=\bcX\t_{\bs g,\bcZ,\bs h}\bcY$ be the d-stack fibre
product, which exists by Theorem\/ {\rm\ref{ds10thm6}(a)}. Then
$\bcW$ is a d-orbifold, with\I{d-orbifold!virtual dimension}\I{fibre
product!of d-orbifolds}
\e
\vdim\bcW=\vdim\bcX+\vdim\bcY-\vdim\bcZ.
\label{ds11eq5}
\e
\label{ds11thm6}
\end{thm}

\begin{thm} Suppose\/ $\bs g:\bcX\ra\bcZ,$ $\bs h:\bcY\ra\bcZ$ are
$1$-morphisms of d-orbifolds. The following are sufficient
conditions for $\bs g,\bs h$ to be d-transverse, so that\/
$\bcW=\bcX\t_{\bs g,\bcZ,\bs h}\bcY$ is a d-orbifold of virtual
dimension\/ {\rm\eq{ds11eq5}:}
\begin{itemize}
\setlength{\itemsep}{0pt}
\setlength{\parsep}{0pt}
\item[{\bf(a)}] $\bcZ$ is an orbifold, that is, $\bcZ\in\hOrb;$
or\I{d-orbifold!is an orbifold}
\item[{\bf(b)}] $\bs g$ or $\bs h$ is a
w-submersion.\I{d-orbifold!w-submersion}
\end{itemize}
\label{ds11thm7}
\end{thm}

\begin{thm} Let\/ $\bcX,\bcZ$ be d-orbifolds, $\bcY$ an orbifold,
and\/ $\bs g:\bcX\ra\bcZ,$ $\bs h:\bcY\ra\bcZ$ be $1$-morphisms
with\/ $\bs g$ a submersion. Then\/ $\bcW=\bcX\t_{\bs g,\bcZ,\bs
h}\bcY$ is an orbifold,\I{d-orbifold!is an orbifold}
with\/~$\dim\bcW=\vdim\bcX+\dim\bcY-
\vdim\bcZ$.\I{d-orbifold!virtual dimension}
\label{ds11thm8}
\end{thm}

\begin{thm}{\bf(i)} Let\/ $\rho:G\ra H$ be a morphism of finite
groups, and\/ $H$ act linearly on $\R^n$. Then as in
{\rm\S\ref{ds102}} we have quotient d-orbifolds\I{d-stack!quotients
$[\bX/G]$} $[\bs */G],$ $[\bR^{\bs n}/H]$ and a quotient\/
$1$-morphism $[\bs 0,\rho]:[\bs
*/G]\ra[\kern .1em\bR^{\bs n}/H]$. Suppose $\bcX$ is a d-orbifold
and\/ $\bs g:\bcX\ra[\kern .1em\bR^{\bs n}/H]$ a $1$-morphism in
$\dOrb$. Then the fibre product\/ $\bcW=\bcX\t_{\bs g,[\kern
.1em\bR^{\bs n}/H],[\bs 0,\rho]}[\bs */G]$ exists in $\dOrb$ by
Theorem\/ {\rm\ref{ds11thm7}(a)}. The projection
$\bs\pi_\bcX:\bcW\ra\bcX$ is an immersion\I{d-orbifold!immersion}
if\/ $\rho$ is injective, and an embedding if\/ $\rho$ is an
isomorphism.
\smallskip

\noindent{\bf(ii)} Suppose\/ $\bs f:\bcX\ra\bcY$ is an immersion of
d-orbifolds, and\/ $[x]\in\cX_\top$ with\/
$f_\top([x])=[y]\in\cY_\top$. Write $\rho:G\ra H$ for
$f_*:\Iso_\cX([x])\ra\Iso_\cY([y])$. Then $\rho$ is injective, and
there exist open neighbourhoods $\bcU\subseteq\bcX$ and\/
$\bcV\subseteq\bcY$ of\/ $[x],[y]$ with\/ $\bs
f(\bcU)\subseteq\bcV,$ a linear action of\/ $H$ on $\R^n$ where
$n=\vdim\bcY-\vdim\bcX\ge 0,$ and a $1$-morphism $\bs
g:\bcV\ra[\kern .1em\bR^{\bs n}/H]$ with\/ $g_\top([y])=[0],$
fitting into a $2$-Cartesian\I{2-category!2-Cartesian square} square
in~$\dOrb:$
\begin{equation*}
\xymatrix@C=80pt@R=10pt{*+[r]{\bcU} \ar[d]^{\bs f\vert_\bcU} \ar[r]
\drtwocell_{}\omit^{}\omit{^{}} & *+[l]{[\bs */G]} \ar[d]_{[\bs 0,\rho]} \\
*+[r]{\bcV} \ar[r]^(0.3){\bs g} & *+[l]{[\kern .1em\bR^{\bs n}/H].} }
\end{equation*}
If\/ $\bs f$ is an embedding\I{d-orbifold!embedding} then $\rho$ is
an isomorphism, and we may take\/~$\bcU\!=\!\bs
f^{-1}(\bcV)$.\I{d-transversality|)}\I{d-orbifold!d-transverse
1-morphisms|)}\I{d-orbifold!fibre products|)}

\label{ds11thm9}
\end{thm}

\subsection{Embedding d-orbifolds into orbifolds}
\label{ds116}
\I{principal
d-orbifold|(}\I{d-orbifold!principal|(}\I{d-orbifold!embedding!into
orbifolds|(}

Section \ref{ds47} discussed embeddings of d-manifolds into
manifolds. Theorem \ref{ds4thm12} gave necessary and sufficient
conditions for the existence of embeddings $\bs f:\bX\ra\bR^{\bs n}$
for any d-manifold $\bX$, and Theorem \ref{ds4thm13} showed that if
a d-manifold $\bX$ has an embedding $\bs f:\bX\ra\bY$ for a manifold
$Y$ then $\bX\simeq\bS_{V,E,s}$ for open $f(X)\subset V\subseteq Y$.
Combining these proves that large classes of d-manifolds --- all
compact d-manifolds, for instance --- are principal d-manifolds.

In \cite[\S 10.5]{Joyc6} we consider how to generalize all this to
d-orbifolds. The proof of Theorem \ref{ds4thm13} extends to
(d-)orbifolds, giving:

\begin{thm} Suppose $\bcX$ is a d-orbifold, $\cY$ an orbifold,
and\/ $\bs f:\bcX\ra\bcY$ an embedding, in the sense of Definition\/
{\rm\ref{ds11def2}}. Then there exist an open suborbifold\/
$\cV\subseteq\cY$ with $\bs f(\bcX)\subseteq\bcV,$ a vector bundle
$\cE$ on $\cV,$ and a smooth section\/ $s\in C^\iy(\cE)$ fitting
into a $2$-Cartesian\I{2-category!2-Cartesian square} diagram in
$\dOrb,$ where\/ $\bcY,\bcV,\bcE,\bs s,\bs
0=F_\Orb^\dOrb\bigl(\cY,\cV,\ab\Tot(\cE),\ab\Tot(s),\ab\Tot(0)
\bigr),$ in the notation of\/~{\rm\S\ref{ds91}:}
\begin{equation*}
\xymatrix@C=60pt@R=10pt{ \bcX \ar[r]_(0.25){\bs f} \ar[d]^{\bs f}
\drtwocell_{}\omit^{}\omit{^{}}
 & \bcV \ar[d]_{\bs 0} \\ \bcV \ar[r]^(0.7){\bs s} & \bcE.}
\end{equation*}
Hence $\bcX$ is equivalent to the `standard model' d-orbifold\/
$\bcS_{\cV,\cE,s}$\I{d-orbifold!standard model a@standard model
$\bcS_{\cV,\cE,s}$} of Example\/ {\rm\ref{ds11ex1},} and is a
principal d-orbifold.
\label{ds11thm10}
\end{thm}

However, we do not presently have a good analogue of Theorem
\ref{ds4thm12} for d-orbifolds, so we cannot state useful necessary
and sufficient conditions for when a d-orbifold $\bcX$ can be
embedded into an orbifold, or is a principal d-orbifold.\I{principal
d-orbifold|)}\I{d-orbifold!principal|)}\I{d-orbifold!embedding!into
orbifolds|)}

\subsection{Orientations of d-orbifolds}
\label{ds117}
\I{d-orbifold!orientations|(}

Section \ref{ds48} discusses orientations on d-manifolds. As in
\cite[\S 10.6]{Joyc6}, all this material generalizes easily to
d-orbifolds, so we will give few details.

If $\cX$ is a Deligne--Mumford $C^\iy$-stack and $(\cE^\bu,\phi)$ a
virtual vector bundle on $\cX$, then we define a line bundle
$\cL_{\smash{(\cE^\bu,\phi)}}$ on $\cX$ called the {\it orientation
line bundle\/}\I{virtual vector bundle!orientation line bundle
of}\I{orientation line bundle} of $(\cE^\bu,\phi)$. It has
functorial properties as in Theorem \ref{ds4thm14}(a)--(f). If
$\bcX$ is a d-orbifold, the virtual cotangent bundle\I{virtual
cotangent bundle} $T^*\bcX=(\EcX,\FcX,\phi_\cX)$ is a virtual vector
bundle on $\cX$. We define an {\it orientation\/} $\om$ on $\bcX$ to
be an orientation on the orientation line bundle
$\cL_{T^*\bcX}$.\G[LT*Xc]{$\cL_{T^*\bcX}$}{orientation line bundle
of a d-orbifold $\bcX$}\I{d-orbifold!orientation line bundle} The
analogues of Theorem \ref{ds4thm15} and Proposition \ref{ds4prop4}
hold for d-orbifolds.

One difference between (d-)manifolds and (d-)orbifolds is that line
bundles $\cL$ on Deligne--Mumford $C^\iy$-stacks $\cX$ (such as
orientation line bundles)\I{orientation line bundle} need only be
locally trivial in the \'etale topology,\I{etale topology@\'etale
topology} not in the Zariski topology.\I{Zariski topology} Because
of this, orbifolds and d-orbifolds need not be (Zariski) locally
orientable. For example, the orbifold $[\ul\R^{2n+1}/\{\pm 1\}]$ is
not locally orientable near~0.\I{d-orbifold!orientations|)}

\subsection{Orbifold strata of d-orbifolds}
\label{ds118}
\I{orbifold strata!of d-orbifolds|(}\I{d-orbifold!orbifold strata|(}

Section \ref{ds87} discussed the {\it orbifold strata\/}
$\cX^\Ga,\tcX^\Ga, \ldots,\hcX^\Ga_\ci$ of a Deligne--Mumford
$C^\iy$-stack $\cX$. When $\cX$ is an orbifold, \S\ref{ds92}
explained that $\cX^\Ga$ decomposes as
$\cX^\Ga=\coprod_{\smash{\la\in\La^\Ga_+}} \cX^{\Ga,\la}$, where
each $\cX^{\Ga,\la}$ is an orbifold of dimension $\dim\cX-\dim\la$,
and similarly for $\tcX^\Ga,\ldots,\hcX^\Ga_\ci$. Section
\ref{ds105} discussed the orbifold strata
$\bcX^\Ga,\btcX^\Ga,\ldots,\bhcX{}^\Ga_\ci$ of a d-stack $\bcX$. In
\cite[\S 10.7]{Joyc6} we show that for a d-orbifold $\bcX$, the
orbifold strata decompose as $\bcX^\Ga=\coprod_{
\smash{\la\in\La^\Ga}} \bcX^{\Ga,\la}$, where $\bcX^{\Ga,\la}$ is a
d-orbifold of virtual dimension $\vdim\bcX-\dim\la$, and similarly
for $\btcX^\Ga,\ldots,\bhcX{}^\Ga_\ci$.

\begin{dfn} Let $\Ga$ be a finite group, and use the notation
$\Rep_\nt(\Ga)$, $\La^\Ga=K_0\bigl(\Rep_\nt(\Ga)\bigr)$,
$\La^\Ga_+\subseteq\La^\Ga$ and $\dim:\La^\Ga\ra\Z$ of Definition
\ref{ds9def4}. Let $R_0,R_1,\ldots,R_k$ be the irreducible
$\Ga$-representations up to isomorphism, with $R_0=\R$ the trivial
representation, so that $\La^\Ga\cong\Z^k$ and~$\La^\Ga_+\cong\N^k$.

Suppose $\bcX$ is a d-orbifold. Theorem \ref{ds10thm7} gives a
d-stack $\bcX^\Ga$ and a 1-morphism $\bs
O{}^\Ga(\bcX):\bcX^\Ga\ra\bcX$. The virtual cotangent
bundle\I{virtual cotangent bundle} of $\bcX$ is
$T^*\bcX=(\EcX,\FcX,\phi_\cX)$, a virtual vector bundle of rank
$\vdim\bcX$ on $\cX$. So $O^\Ga(\cX)^*
(T^*\bcX)=\bigl(O^\Ga(\cX)^*(\EcX),O^\Ga(\cX)^*(\FcX),O^\Ga(\cX)^*
(\phi_\cX)\bigr)$ is a virtual vector bundle on $\cX^\Ga$. As in
\S\ref{ds87}, $O^\Ga(\cX)^*(\EcX),O^\Ga(\cX)^*(\FcX)$ have natural
$\Ga$-representations inducing decompositions of the form
\eq{ds8eq10}--\eq{ds8eq11}, and $O^\Ga(\cX)^*(\phi_\cX)$ is
$\Ga$-equivariant and so preserves these splittings. Hence we have
decompositions in $\vqcoh(\cX^\Ga)$:
\e
\begin{gathered}
O^\Ga(\cX)^* (T^*\bcX)\cong\ts\bigop_{i=0}^k (T^*\bcX)^\Ga_i\ot
R_i\;\>\text{for}\;\> (T^*\bcX)^\Ga_i\in\vqcoh(\cX^\Ga), \\
\text{and}\;\> O^\Ga(\cX)^*(T^*\bcX)=(T^*\bcX)^\Ga_\tr
\op(T^*\bcX)^\Ga_\nt,\;\>\text{with} \\
(T^*\bcX)^\Ga_\tr\cong (T^*\bcX)^\Ga_0\ot R_0\;\>\text{and}\;\>
(T^*\bcX)^\Ga_\nt\cong\ts \bigop_{i=1}^k (T^*\bcX)^\Ga_i\ot R_i.
\end{gathered}
\label{ds11eq6}
\e
Also Theorem \ref{ds10thm10} shows that $T^*(\bcX^\Ga)\cong
(T^*\bcX)^\Ga_\tr$.

As $O^\Ga(\cX)^*(T^*\bcX)$ is a virtual vector bundle, \eq{ds11eq6}
implies the $(T^*\bcX)^\Ga_i$ are {\it virtual vector bundles of
mixed rank},\I{virtual vector bundle!of mixed rank} whose ranks may
vary on different connected components of $\bcX^\Ga$. For each
$\la\in\La^\Ga$, define
$\bcX^{\Ga,\la}$\G[XGag]{$\bcX^{\Ga,\la},\btcX{}^{\Ga,\mu},\bhcX{}^{\Ga,\mu},
\bcX{}^{\Ga,\la}_\ci,\btcX{}^{\Ga,\mu}_\ci,\bhcX{}^{\Ga,\mu}_\ci$}{orbifold
strata of a d-orbifold $\bcX$} to be the open and closed d-substack
in $\bcX^\Ga$ with $\rank\bigl((T^*\bcX)^\Ga_1
\bigr)[R_1]+\cdots+\rank\bigl((T^*\bcX)^\Ga_k\bigr)[R_k]=\la$ in
$\La^\Ga$. Then $\bcX^{\Ga,\la}$ is a d-orbifold, with
$\vdim\bcX^{\Ga,\la}=\vdim\bcX-\dim\la$. Also we have a
decomposition $\bcX^\Ga=\coprod_{\la\in\La^\Ga}\bcX^{\Ga,\la}$
in~$\dSta$.

Note that in the d-orbifold case $\dim\la$ may be negative, so we
can have $\vdim\bcX^{\Ga,\la}>\vdim\bcX$. This is counterintuitive:
the (w-immersed) d-suborbifold $\bcX^{\Ga,\la}$ has larger dimension
than the d-orbifold $\bcX$ that contains it.

Write $\bs O^{\Ga,\la}(\bcX)=\bs O^\Ga(\bcX)\vert_{\smash{
\bcX^{\Ga,\la}}}:\bcX^{\Ga,\la}\ra\bcX$. Then $\bs
O^{\Ga,\la}(\bcX)$ is a proper w-immersion of d-orbifolds, in the
sense of \S\ref{ds114}. Define
$\bcX^{\Ga,\la}_\ci=\bcX^\Ga_\ci\cap\bcX^{\Ga,\la}$, and $\bs
O^{\Ga,\la}_\ci(\bcX)=\bs O^\Ga_\ci(\bcX)\vert_{\smash{
\bcX^{\Ga,\la}_\ci}}:\bcX^{\Ga,\la}_\ci\ra\bcX$. Then
$\bcX^{\Ga,\la}_\ci$ is a d-orbifold with
$\vdim\bcX^{\Ga,\la}_\ci=\vdim\bcX-\dim\la$, and
$\bcX^\Ga_\ci=\ts\coprod_{\smash{\la\in\La^\Ga}}\bcX^{\Ga,\la}_\ci$.

As for $\tcX^{\Ga,\mu},\ldots,\hcX^{\Ga,\mu}_\ci$ in \S\ref{ds92},
for each $\mu\in\La^\Ga/\Aut(\Ga)$ we define
$\btcX^{\Ga,\mu}\simeq\bigl[\bigl(\coprod_{\la\in\mu}\bcX^{\Ga,\la}
\bigr)\big/\Aut(\Ga)\bigr]$ in $\btcX^\Ga\simeq
[\bcX^\Ga/\Aut(\Ga)]$, and
$\btcX^{\Ga,\mu}_\ci=\btcX^\Ga_\ci\cap\btcX^{\Ga,\mu}$, and
$\bhcX^{\Ga,\mu}=\bs{\hat\Pi}{}^\Ga(\bcX)(\btcX^{\Ga,\mu})$, and
$\bhcX^{\Ga,\mu}_\ci= \bhcX^\Ga_\ci\cap \bhcX^{\Ga,\mu}$. Then
$\btcX^{\Ga,\mu},\ldots,\bhcX^{\Ga,\mu}_\ci$ are d-orbifolds with
$\vdim\btcX^{\Ga,\mu}=\cdots=\vdim\bhcX^{\Ga,\mu}_\ci=\vdim\bcX-\dim\mu$,
with
\begin{equation*}
\btcX^\Ga=\ts\coprod_\mu\btcX^{\Ga,\mu},\;\>
\btcX^\Ga_\ci=\ts\coprod_\mu \btcX^{\Ga,\mu}_\ci,\;\>
\bhcX^\Ga=\ts\coprod_\mu \bhcX^{\Ga,\mu},\;\>
\bhcX^\Ga_\ci=\ts\coprod_\mu \bhcX^{\Ga,\mu}_\ci.
\end{equation*}
Also $\bhcX^{\Ga,\mu}_\ci$ is a d-manifold, that is, it lies
in~$\hdMan$.\I{d-orbifold!is a d-manifold}
\label{ds11def4}
\end{dfn}

In \cite[\S 10.7]{Joyc6}\I{d-orbifold!orbifold strata!orientations
on|(}\I{d-orbifold!orientations|(} we also consider the question: if
$\bcX$ is an oriented d-orbifold, under what conditions on
$\Ga,\la,\mu$ do the orbifold strata
$\bcX^{\Ga,\la},\ldots,\bhcX{}^{\Ga,\mu}_\ci$ have natural
orientations? Here is the analogue of Proposition~\ref{ds9prop2}:

\begin{prop}{\bf(a)} Let\/ $\Ga$ be a finite group with\/
$\md{\Ga}$ odd, and\/ $\bcX$ an oriented d-orbifold. Then we may
define orientations on $\bcX^{\Ga,\la},\bcX^{\Ga,\la}_\ci$ for
all\/~$\la\in\La^\Ga$.
\smallskip

\noindent{\bf(b)} Let\/ $\Ga$ be a finite group with\/ $\md{\Ga}$
odd, $\la\in\La^\Ga$ and\/ $\mu=\la\cdot\Aut(\Ga)$ in $\La^\Ga/
\Aut(\Ga)$. We may write $\la=[(V^+,\rho^+)]-[(V^-,\rho^-)]$ for
nontrivial $\Ga$-representations $(V^\pm,\rho^\pm)$ with no common
subrepresentation, and then\/ $(V^\pm,\rho^\pm)$ are unique up to
isomorphism. Define $H$ to be the subgroup of\/ $\Aut(\Ga)$ fixing\/
$\la$ in\/ $\La^\Ga$. Then for each\/ $\de\in H$ there exist
isomorphisms of\/ $\Ga$-representations
$i_\de^\pm:(V^\pm,\rho^\pm\ci\de) \ra(V^\pm,\rho^\pm)$. Suppose
$i_\de^+\op i_\de^-:V^+\op V^-\ra V^+\op V^-$ is
orientation-preserving for all\/ $\de\in H$. If\/ $\la\in 2\La^\Ga$
this holds automatically.

Then for all oriented d-orbifolds\/ $\bcX$ we can define
orientations on the orbifold strata\/ $\btcX{}^{\Ga,\mu},
\ab\btcX{}^{\Ga,\mu}_\ci,\ab
\bhcX{}^{\Ga,\mu},\ab\bhcX{}^{\Ga,\mu}_\ci$. For $\btcX{}^{\Ga,\mu}$
this works as\/ $\btcX{}^{\Ga,\mu}\simeq [\bcX^{\Ga,\la}/H],$ where
$\bcX^{\Ga,\la}$ is oriented by {\bf(a)\rm,} and the $H$-action on
$\bcX^{\Ga,\la}$ preserves orientations, so the orientation on
$\bcX^{\Ga,\la}$ descends to an orientation
on\/~$\btcX{}^{\Ga,\mu}\simeq [\bcX^{\Ga,\la}/H]$.
\smallskip

\noindent{\bf(c)} Suppose that\/ $\Ga$ and\/ $\la\in\La^\Ga$ do not
satisfy the conditions in {\bf(a)} (i.e.\ $\md{\Ga}$ is even), or\/
$\Ga$ and\/ $\mu\in\La^\Ga/\Aut(\Ga)$ do not satisfy the conditions
in {\bf(b)}. Then we can find examples of oriented d-orbifolds\/
$\bcX$ such that\/ $\bcX^{\Ga,\la},\bcX^{\Ga,\la}_\ci$ are not
orientable, or $\btcX{}^{\Ga,\mu}, \ab\btcX{}^{\Ga,\mu}_\ci,\ab
\bhcX{}^{\Ga,\mu},\ab\bhcX{}^{\Ga,\mu}_\ci$ are not orientable,
respectively. That is, the conditions on $\Ga,\la,\mu$ in
{\bf(a)\rm,\bf(b)} are necessary as well as sufficient to be able to
orient orbifold strata\/ $\bcX^{\Ga,\la},\ldots,
\bhcX{}^{\Ga,\mu}_\ci$ of all oriented d-orbifolds\/~$\bcX$.
\label{ds11prop2}
\end{prop}

Note that Proposition \ref{ds11prop2} for d-orbifolds is weaker than
Proposition \ref{ds9prop2} for orbifolds. That is, if $\Ga$ is a
finite group with $\md{\Ga}$ even then for some choices of $\la,\mu$
we can orient $\cX^{\Ga,\la},\ldots,\hcX^{\Ga,\mu}_\ci$ for all
oriented orbifolds $\cX$, but we cannot orient $\bcX^{\Ga,\la},
\ldots,\bhcX{}^{\Ga,\mu}_\ci$ for all oriented
d-orbifolds~$\bcX$.\I{orbifold strata!of
d-orbifolds|)}\I{d-orbifold!orbifold strata|)}\I{d-orbifold!orbifold
strata!orientations on|)}\I{d-orbifold!orientations|)}

\subsection[Kuranishi neighbourhoods and good coordinate
systems]{Kuranishi neighbourhoods, good coordinate systems}
\label{ds119}
\I{d-orbifold!Kuranishi neighbourhood|(}\I{d-orbifold!good
coordinate system|(}\I{good coordinate
system|(}\I{d-orbifold!standard model b@standard model
$[\bS_{V,E,s}/\Ga]$|(}

We now explain the main ideas of \cite[\S 10.8]{Joyc6}, which are
based on parallel material about Kuranishi spaces\I{Kuranishi space}
due to Fukaya, Oh, Ohta and Ono~\cite{FuOn,FOOO}.

\begin{dfn} Let $\bcX$ be a d-orbifold. A {\it type A Kuranishi
neighbourhood\/} on $\bcX$ is a quintuple $(V,E,\Ga,s,\bs\psi)$
where $V$ is a manifold, $E\ra V$ a vector bundle, $\Ga$ a finite
group acting smoothly and locally effectively on $V,E$ preserving
the vector bundle structure, and $s:V\ra E$ a smooth,
$\Ga$-equivariant section of $E$. Write the $\Ga$-actions on $V,E$
as $r(\ga):V\ra V$ and $\hat r(\ga):E\ra r(\ga)^*(E)$ for
$\ga\in\Ga$. Then Example \ref{ds11ex3} defines a principal
d-orbifold $[\bS_{V,E,s}/\Ga]$. We require that $\bs\psi:
[\bS_{V,E,s}/\Ga]\ra\bcX$ is a 1-morphism of d-orbifolds which is an
equivalence with a nonempty open d-suborbifold
$\bs\psi([\bS_{V,E,s}/\Ga])\subseteq\bcX$.
\label{ds11def5}
\end{dfn}

\begin{dfn} Suppose $(V_i,E_i,\Ga_i,s_i,\bs\psi_i),(V_j,E_j,\Ga_j,
s_j,\bs\psi_j)$ are type A Kuranishi neighbourhoods on a d-orbifold
$\bcX$, with\I{d-orbifold!Kuranishi neighbourhood!coordinate
change|(}
\begin{equation*}
\bs\es\ne\bs\psi_i([\bS_{V_i,E_i,s_i}/\Ga_i])\cap
\bs\psi_j([\bS_{V_j,E_j,s_j}/\Ga_j])\subseteq\bcX.
\end{equation*}
A {\it type A coordinate change from\/
$(V_i,E_i,\Ga_i,s_i,\bs\psi_i)$ to\/}
$(V_j,E_j,\Ga_j,s_j,\bs\psi_j)$ is a quintuple $(V_{ij},e_{ij},\hat
e_{ij},\rho_{ij},\bs\eta_{ij})$, where:
\begin{itemize}
\setlength{\itemsep}{0pt}
\setlength{\parsep}{0pt}
\item[(a)] $\es\ne V_{ij}\subseteq V_i$ is a $\Ga_i$-invariant
open submanifold, with
\begin{equation*}
\bs\psi_i\bigl([\bS_{V_{ij},E_i\vert_{V_{ij}},s_i\vert_{V_{ij}}}/
\Ga_i]\bigr)=\bs\psi_i([\bS_{V_i,E_i,s_i}/\Ga_i])\cap
\bs\psi_j([\bS_{V_j,E_j,s_j}/\Ga_j])\subseteq\bcX.
\end{equation*}
\item[(b)] $\rho_{ij}:\Ga_i\ra\Ga_j$ is an injective group morphism.
\item[(c)] $e_{ij}:V_{ij}\ra V_j$ is an embedding of manifolds with
$e_{ij}\ci r_i(\ga)= r_j(\rho_{ij}(\ga))\ci e_{ij}:V_{ij} \ra
V_j$ for all $\ga\in\Ga_i$. If $v_i,v_i'\in V_{ij}$ and
$\de\in\Ga_j$ with $r_j(\de)\ci e_{ij}(v_i')=e_{ij}(v_i)$, then
there exists $\ga\in\Ga_i$ with $\rho_{ij}(\ga)=\de$ and
$r_i(\ga)(v_i')=v_i$.
\item[(d)] $\hat e_{ij}:E_i\vert_{V_{ij}}\ra
e_{ij}^*(E_j)$ is an embedding of vector bundles (that is, $\hat
e_{ij}$ has a left inverse), such that $\hat e_{ij}\ci
s_i\vert_{V_{ij}}=e_{ij}^*(s_j)$ and $r_i(\ga)^*(\hat
e_{ij})\ci\hat r_i(\ga)= e_{ij}^*(\hat
r_j(\rho_{ij}(\ga)))\ci\hat e_{ij}:E_i\vert_{V_{ij}}
\ra(e_{ij}\ci r_i(\ga))^*(E_j)$ for all $\ga\in\Ga_i$. Thus
Example \ref{ds11ex4} defines a quotient 1-morphism
\e
[\bS_{\smash{e_{ij},\hat e_{ij}}},\rho_{ij}]:[\bS_{V_{ij},
E_i\vert_{V_{ij}},s_i \vert_{V_{ij}}}/\Ga_i] \longra
[\bS_{V_j,E_j,s_j}/\Ga_j],
\label{ds11eq7}
\e
where $[\bS_{V_{ij},E_i\vert_{V_{ij}},s_i\vert_{V_{ij}}}/
\Ga_i]$ is an open d-suborbifold in $[\bS_{V_i,E_i,s_i}/\Ga_i]$.
\item[(e)] If $v_i\in V_{ij}$ with $s_i(v_i)=0$ and
$v_j=e_{ij}(v_i)\in V_j$ then the following linear map is an
isomorphism:
\begin{equation*}
{}\!\!\bigl(\d s_j(v_j)\bigr){}_*\!:\!\bigl(T_{v_j}V_j\bigr)\big/
\bigl(\d e_{ij}(v_i)[T_{v_i}V_i]\bigr)\!\ra\!\bigl(E_j\vert_{ v_j}
\bigr)\big/\bigl(\hat e_{ij}(v_i)[E_i\vert_{v_i}]\bigr).\!\!{}
\end{equation*}
Theorem \ref{ds11thm3} then implies that $[\bS_{\smash{e_{ij},
\hat e_{ij}}},\rho_{ij}]$ in \eq{ds11eq7} is an equivalence with
an open d-suborbifold of~$[\bS_{V_j,E_j,s_j}/\Ga_j]$.
\item[(f)] $\bs\eta_{ij}:\bs\psi_j\ci[\bS_{e_{ij},\hat e_{ij}},
\rho_{ij}]\Ra\bs\psi_i\vert_{[\bS_{V_{ij},E_i \vert_{V_{ij}},
s_i\vert_{V_{ij}}}/\Ga_i]}$ is a 2-morphism in~$\dOrb$.
\item[(g)] The quotient topological space $V_i\amalg_{V_{ij}}V_j=
(V_i\amalg V_j)/\sim$ is Hausdorff, where the equivalence
relation $\sim$ identifies $v\in V_{ij}\subseteq V_i$
with~$e_{ij}(v)\in V_j$.\I{d-orbifold!Kuranishi
neighbourhood!coordinate change|)}
\end{itemize}
\label{ds11def6}
\end{dfn}

\begin{dfn} Let $\bcX$ be a d-orbifold. A {\it type A good
coordinate system\/} on $\bcX$ consists of the following data
satisfying conditions~(a)--(e):
\begin{itemize}
\setlength{\itemsep}{0pt}
\setlength{\parsep}{0pt}
\item[(a)] We are given a countable indexing set $I,$ and a
total order $<$ on $I$ making $(I,<)$ into a well-ordered set.
\item[(b)] For each $i\in I$ we are given a Kuranishi
neighbourhood $(V_i,E_i,\Ga_i,s_i,\bs\psi_i)$ of type A on
$\bcX$. Write $\bcX_i=\bs\psi_i([\bS_{V_i,E_i,s_i}/\Ga_i])$, so
that $\bcX_i\subseteq\bcX$ is an open d-suborbifold, and
$\bs\psi_i:[\bS_{V_i,E_i,s_i}/\Ga_i]\ra\bcX_i$ is an
equivalence. We require that $\bigcup_{i\in I}\bcX_i=\bcX$, so
that $\{\bcX_i:i\in I\}$ is an open cover of~$\bcX$.
\item[(c)] For all $i<j$ in $I$ with $\bcX_i\cap\bcX_j \ne\bs\es$
we are given a type A coordinate change $(V_{ij},e_{ij},\hat
e_{ij},\rho_{ij},\bs\eta_{ij})$ from
$(V_i,E_i,\Ga_i,s_i,\bs\psi_i)$ to
$(V_j,E_j,\Ga_j,s_j,\bs\psi_j)$.
\item[(d)] For all $i<j<k$ in $I$ with $\bcX_i\cap\bcX_j
\cap\bcX_k\ne\bs\es$, we are given $\ga_{ijk}\in\Ga_k$
satisfying $\rho_{ik}(\ga)=\ga_{ijk}\,\rho_{jk}(\rho_{ij}(\ga))
\,\ga_{ijk}^{-1}$ for all $\ga\in\Ga_i,$ and
\e
\begin{split}
{}\!\!\!\!\!\! e_{ik}\vert_{V_{ik}\cap e_{ij}^{-1}(V_{jk})}&=
r_k(\ga_{ijk})\ci e_{jk}\ci
e_{ij}\vert_{V_{ik}\cap e_{ij}^{-1}(V_{jk})},\\
{}\!\!\!\!\!\! \hat e_{ik}\vert_{V_{ik}\cap
e_{ij}^{-1}(V_{jk})}&= \bigl(e_{ij}^*(e_{jk}^*(\hat
r_k(\ga_{ijk})))\ci e_{ij}^* (\hat e_{jk})\ci\hat
e_{ij}\bigr)\vert_{V_{ik}\cap e_{ij}^{-1}(V_{jk})}.
\end{split}
\label{ds11eq8}
\e
Combining the first equation of \eq{ds11eq8} with Definition
\ref{ds11def6}(c) for $e_{ik}$ and $\Ga_i$ acting effectively on
$V_{ik}\cap e_{ij}^{-1}(V_{jk})$ shows that $\ga_{ijk}$ is
unique. Example \ref{ds11ex5} with $\de=\ga_{ijk}$ and $\La=0$
then gives a 2-morphism in~$\dOrb$:
\begin{align*}
\bs\eta_{ijk}=[\bS_0,\ga_{ijk}]:[\bS_{\smash{e_{jk},\hat e_{jk}}},
\rho_{jk}]\ci[\bS_{\smash{e_{ij}, \hat e_{ij}}},\rho_{ij}]
\vert_{[\bS_{V_{ik}\cap e_{ij}^{-1}(V_{jk}),E_i, s_i}/\Ga_i]}&\\
\Longra [\bS_{\smash{e_{ik},\hat e_{ik}}},\rho_{ik}]
\vert_{[\bS_{V_{ik}\cap e_{ij}^{-1}(V_{jk}),E_i, s_i}/\Ga_i]}&.
\end{align*}
\item[(e)] For all $i<j<k$ in $I$ with
$\bcX_i\cap\bcX_k\ne\bs\es$ and $\bcX_j\cap\bcX_k\ne\bs\es$, we
require that if $v_i\in V_{ik}$, $v_j\in V_{jk}$ and
$\de\in\Ga_k$ with $e_{jk}(v_j)=r_k(\de)\ci e_{ik}(v_i)$ in
$V_k$, then $\bcX_i\cap\bcX_j\cap\bcX_k\ne\bs\es$, and $v_i\in
V_{ij}$, and there exists $\ga\in\Ga_j$ with
$\rho_{jk}(\ga)=\de\,\ga_{ijk}$ and~$v_j=r_j(\ga)\ci
e_{ij}(v_i)$.
\end{itemize}

Suppose now that $Y$ is a manifold, and $\bs h:\bcX\ra\bcY$ is a
1-morphism in $\dOrb$, where $\bcY=F_\Man^\dOrb(Y)$. A {\it type A
good coordinate system for\/} $\bs h:\bcX\ra\bcY$ consists of a type
A good coordinate system $\bigl(I,<,\ldots,\ga_{ijk}\bigr)$ for
$\bcX$ as in (a)--(e) above, together with the following data
satisfying conditions~(f)--(g):
\begin{itemize}
\setlength{\itemsep}{0pt}
\setlength{\parsep}{0pt}
\item[(f)] For each $i\in I$, we are given a smooth map
$g_i:V_i\ra Y$ with $g_i\ci r_i(\ga)= g_i$ for all
$\ga\in\Ga_i$, so that Example \ref{ds11ex4} defines a quotient
1-morphism
\begin{equation*}
[\bS_{g_i,0},\pi]:[\bS_{V_i, E_i,s_i}/\Ga_i] \longra
[\bS_{Y,0,0}/\{1\}]=\bcY,
\end{equation*}
where $\pi:\Ga_i\ra\{1\}$ is the projection. We are given a
2-morphism $\bs\ze_i:\bs h\ci\bs\psi_i\Ra[\bS_{g_i,0},\pi]$ in
$\dOrb$. Sometimes we require $g_i$ to be a submersion.
\item[(g)] For all $i<j$ in $I$ with $\bcX_i\cap\bcX_j \ne\es$,
we require that $g_j\ci e_{ij}= g_i\vert_{V_{ij}}$. This implies
that
\begin{align*}
[\bS_{g_j,0},\pi]\ci [\bS_{\smash{e_{ij},\hat
e_{ij}}},\rho_{ij}]=[\bS_{g_i,0},\pi]\vert_{[\bS_{V_{ij},
E_i\vert_{V_{ij}},s_i\vert_{V_{ij}}}/\Ga_i]}:&\\
[\bS_{V_{ij},E_i\vert_{V_{ij}},s_i \vert_{V_{ij}}}/\Ga_i]
\longra [\bS_{Y,0,0}/\{1\}]&=\bcY.
\end{align*}
\end{itemize}
\label{ds11def7}
\end{dfn}

Here is the main result of \cite[\S 10.8]{Joyc6}, which is proved
in~\cite[App.~D]{Joyc6}.

\begin{thm} Suppose $\bcX$ is a d-orbifold. Then there exists a
type A good coordinate system
$\bigl(I,<,(V_i,E_i,\Ga_i,s_i,\bs\psi_i),(V_{ij},e_{ij},\hat
e_{ij},\rho_{ij},\bs\eta_{ij}), \ga_{ijk}\bigr)$ for $\bcX$. If\/
$\bcX$ is compact, we may take $I$ to be finite. If\/ $\{\bcU_j:j\in
J\}$ is an open cover of\/ $\bcX,$ we may take
$\bcX_i=\bs\psi_i([\bS_{V_i,E_i,s_i}/\Ga_i]) \subseteq
\bcU_{\smash{j_i}}$ for each\/ $i\in I$ and some\/~$j_i\in J$.

Now let\/ $Y$ be a manifold and\/ $\bs h:\bcX\ra\bcY=
F_\Man^\dOrb(Y)$ a $1$-morphism in $\dOrb$. Then all the above
extends to type A good coordinate systems for $\bs
h\!:\!\bcX\!\ra\!\bcY,$ and we may take the $g_i$ in Definition\/
{\rm\ref{ds11def7}(f)} to be submersions.\I{d-orbifold!standard
model b@standard model $[\bS_{V,E,s}/\Ga]$|)}
\label{ds11thm11}
\end{thm}

In \cite[\S 10.8]{Joyc6} we also give `type B' versions of
Definitions \ref{ds11def5}--\ref{ds11def7} and Theorem
\ref{ds11thm11} using the standard model d-orbifolds
$\bcS_{\cV,\cE,s}$\I{d-orbifold!standard model a@standard model
$\bcS_{\cV,\cE,s}$} and 1-morphisms $\bcS_{\smash{e_{ij},\hat
e_{ij}}}$ of Examples \ref{ds11ex1} and \ref{ds11ex2} in place of
$[\bS_{V,E,s}/\Ga]$ and $[\bS_{\smash{e_{ij},\hat
e_{ij}}},\rho_{ij}]$ from Examples \ref{ds11ex3} and~\ref{ds11ex4}.

Observe that Definition \ref{ds11def7} is similar to the hypotheses
of Theorem \ref{ds11thm5}. Given a good coordinate system
$I,<,(V_i,E_i,\Ga_i,s_i,\bs\psi_i),\ldots$ on $\bcX$, Theorem
\ref{ds11thm5} reconstructs $\bcX$ up to equivalence in $\dOrb$ from
the data $I,\ab <,\ab V_i,\ab E_i,\ab\Ga_i,s_i,V_{ij},e_{ij},\hat
e_{ij},\rho_{ij}, \ga_{ijk}$. Thus, we can regard Theorem
\ref{ds11thm11} as a kind of converse to Theorem \ref{ds11thm5}.
Combining the two, we see that every d-orbifold $\bcX$ can be
described up to equivalence by a collection of
differential-geometric data $I,<,V_i,\ldots, \ga_{ijk}$. The `type
B' version of Theorem \ref{ds11thm11} is a kind of converse to
Theorem~\ref{ds11thm4}.

Fukaya and Ono \cite[\S 5]{FuOn} and Fukaya, Oh, Ohta and Ono
\cite[\S A1]{FOOO} define {\it Kuranishi spaces},\I{Kuranishi
space|(} the geometric structure they put on moduli spaces of
$J$-holomorphic curves\I{moduli space!of J-holomorphic curves@of
$J$-holomorphic curves} in symplectic geometry.\I{symplectic
geometry} We argue in \cite[\S 14.3]{Joyc6} that their definition is
not really satisfactory, and that the `right' way to define
Kuranishi spaces is as d-orbifolds, or d-orbifolds with corners.

A {\it Kuranishi space\/} in \cite[\S A1]{FOOO} is a topological
space $X$ with a cover by `Kuranishi neighbourhoods'
$(V,E,\Ga,s,\psi)$, which are as in Definition \ref{ds11def5} except
that $\psi$ is a homeomorphism with an open set in $X$, rather than
an equivalence with an open d-suborbifold. On overlaps between
(images of) Kuranishi neighbourhoods in $X$ we are given `coordinate
changes', roughly as in Definition \ref{ds11def6} except for the
2-morphisms $\bs\eta_{ij}$. Fukaya et al.\ define `good coordinate
systems' for Kuranishi spaces, roughly as in Definition
\ref{ds11def7}. They state without proof in \cite[Lem.~A1.11]{FOOO}
that good coordinate systems exist for any (compact) Kuranishi
space, the analogue of Theorem~\ref{ds11thm11}.

Good coordinate systems are used in \cite{FuOn,FOOO} in some kinds
of proof involving Kuranishi spaces, in particular, in the
construction of virtual classes\I{virtual class} and virtual
chains.\I{virtual chain} The proofs involve choosing data (such as a
multi-valued perturbation of $s_i$) on each Kuranishi neighbourhood
$(V_i,E_i,\Ga_i,s_i,\psi_i)$, by induction on $i$ in $I$ in the
order $<$, where the data must satisfy compatibility conditions with
coordinate changes~$(V_{ij},e_{ij},\hat e_{ij},\rho_{ij})$.

In fact we have already met the problem good coordinate systems are
designed to solve in \S\ref{ds116}: in contrast to the d-manifold
case, we do not have useful criteria for when a d-orbifold $\bcX$ is
principal.\I{d-orbifold!principal}\I{principal d-orbifold} The
parallel issue for Kuranishi spaces is that we cannot cover a
general Kuranishi space $\bcX$ with a single Kuranishi neighbourhood
$(V,E,\Ga,s,\psi)$. So we cover (compact) $\bcX$ with (finitely)
many Kuranishi neighbourhoods $(V_i,E_i,\Ga_i,s_i,\psi_i)$ with
particularly well-behaved coordinate changes on overlaps, and then
carry out the construction we want on each
$(V_i,E_i,\Ga_i,s_i,\psi_i)$, compatibly with coordinate changes.

The material above is used in \cite[\S 14.3]{Joyc6} to explain the
relations between d-orbifolds and Kuranishi spaces.\I{d-orbifold!and
Kuranishi spaces}\I{Kuranishi space!and d-orbifolds} As for
Kuranishi spaces, it is also helpful for some proofs involving
d-orbifolds, for instance, in constructing virtual classes\I{virtual
class!for d-orbifolds} for compact oriented d-orbifolds, and in
studying d-orbifold bordism.\I{bordism!d-orbifold
bordism}\I{d-orbifold!bordism}\I{d-orbifold!Kuranishi
neighbourhood|)}\I{d-orbifold!good coordinate system|)}\I{Kuranishi
space|)}\I{good coordinate system|)}

\subsection{Semieffective and effective d-orbifolds}
\label{ds1110}
\I{d-orbifold!semieffective|(}\I{d-orbifold!effective|(}

In \cite[\S 10.9]{Joyc6} we define {\it semieffective\/} and {\it
effective\/} d-orbifolds, which are related to the notion of
effective orbifold in Definition~\ref{ds9def3}.

\begin{dfn} Let $\bcX$ be a d-orbifold. For $[x]\in\cX_\top$, so that
$x:\bar{\ul *}\ra\cX$ is a $C^\iy$-stack 1-morphism, applying
pullback $x^*$ to \eq{ds10eq1} gives an exact sequence in
$\qcoh(\bar{\ul *})$, where $K_{[x]}=\Ker(x^*(\phi_\cX))$:
\begin{equation*}
\xymatrix@C=15pt{ 0 \ar[r] & K_{[x]} \ar[r] & x^*(\EcX)
\ar[rr]^{x^*(\phi_\cX)} && x^*(\FcX) \ar[rr]^(0.4){x^*(\psi_\cX)} &&
x^*(T^*\cX)\cong T^*_x\cX \ar[r] & 0.}
\end{equation*}
We may think of this as an exact sequence of real vector spaces,
where $K_{[x]},T^*_x\cX$ are finite-dimensional with~$\dim
T^*_x\cX-\dim K_{[x]}=\vdim\bcX$.

The orbifold group $\Iso_\cX([x])$\I{C-stack@$C^\iy$-stack!orbifold
group $\Iso_\cX([x])$} is the group of 2-morphisms $\eta:x\Ra x$.
Definition \ref{ds8def10} defines isomorphisms
$\eta^*(\EcX):x^*(\EcX)\ra x^*(\EcX)$ in $\qcoh(\bar{\ul *})$, which
make $x^*(\EcX)$ into a representation of $\Iso_\cX([x])$. The same
holds for $x^*(\FcX),x^*(T^*\cX)$, and $x^*(\phi_\cX),x^*(\psi_\cX)$
are equivariant. Hence $K_{[x]},T^*_x\cX$ are also
$\Iso_\cX([x])$-representations.

We call $\bcX$ a {\it semieffective d-orbifold\/} if $K_{[x]}$ is a
trivial representation of $\Iso_\cX([x])$ for all $[x]\in\cX_\top$.
We call $\bcX$ an {\it effective d-orbifold\/} if it is
semieffective, and $T^*_x\cX$ is an effective representation of
$\Iso_\cX([x])$ for all~$[x]\in\cX_\top$.
\label{ds11def8}
\end{dfn}

That is, $\bcX$ is semieffective if the orbifold groups
$\Iso_\cX([x])$ act trivially on the obstruction spaces of $\bcX$,
and effective if the $\Iso_\cX([x])$ also act effectively on the
tangent spaces of $\bcX$. One useful property of (semi)effective
d-orbifolds is that generic perturbations of semieffective (or
effective) d-orbifolds are (effective)
orbifolds.\I{d-orbifold!perturbing to orbifolds|(} We state this for
`standard model' d-orbifolds $\bcS_{\cV,\cE,s}$.

\begin{prop} Let\/ $\cV$ be an orbifold, $\cE$ a vector bundle on
$\cV,$ and\/ $s\in C^\iy(\cE),$ and let\/ $\bcS_{\cV,\cE,s}$ be as
in Example\/ {\rm\ref{ds11ex1}}. Suppose $\bcS_{\cV,\cE,s}$ is a
semieffective d-orbifold. Then for any generic perturbation $\ti s$
of\/ $s$ in $C^\iy(\cE)$ with $\ti s-s$ sufficiently small in $C^1$
locally on $\cV,$ the d-orbifold $\bcS_{\cV,\cE,\ti s}$ is an
orbifold, that is, it lies in $\hOrb\subset\dOrb$. If\/
$\bcS_{\cV,\cE,s}$ is an effective d-orbifold, then
$\bcS_{\cV,\cE,\ti s}$ is an effective
orbifold.\I{orbifold!effective}\I{d-orbifold!perturbing to
orbifolds|)}
\label{ds11prop3}
\end{prop}

Here are some other good properties of (semi)effective d-orbifolds:
\begin{itemize}
\setlength{\itemsep}{0pt}
\setlength{\parsep}{0pt}
\item If $\cX$ is an orbifold then $\bcX=F_\Orb^\dOrb(\cX)$ is
a semieffective d-orbifold, and if $\cX$ is effective then
$\bcX$ is effective.\I{orbifold!effective}
\item Let $\bcX$ be a semieffective d-orbifold, $\Ga$
a finite group, and $\la\in\La^\Ga$. Then the orbifold stratum
$\bcX^{\Ga,\la}=\bs\es$ unless $\la\in\La^\Ga_+\subset\La^\Ga$.
If $\bcX$ is effective then $\bcX^{\Ga,\la}=\bs\es$ unless
$\la=[R]$ for $R$ an effective
$\Ga$-representation.\I{d-orbifold!semieffective!orbifold strata
of}\I{d-orbifold!effective!orbifold strata of}
\item If $\bcX,\bcY$ are (semi)effective d-orbifolds, then
the product $\bcX\t\bcY$ is also (semi)effective. More
generally, any fibre product $\bcX\t_\bcZ\bcY$ in $\dOrb$ with
$\bcX,\bcY$ (semi)effective and $\bcZ$ a manifold is also
(semi)effective.
\item Proposition \ref{ds11prop2} says that if $\bcX$ is an
oriented d-orbifold, then when $\md{\Ga}$ is odd we can define
orientations on the orbifold strata $\bcX^{\Ga,\la},
\bcX^{\Ga,\la}_\ci$, and under extra conditions on $\mu$ we can
also orient $\btcX{}^{\Ga,\mu},\ab\btcX{}^{\Ga,\mu}_\ci,\ab
\bhcX{}^{\Ga,\mu},\ab\bhcX{}^{\Ga,\mu}_\ci$.\I{d-orbifold!orbifold
strata!orientations on}

For general d-orbifolds $\bcX$, this is the best we can do. But
for semieffective d-orbifolds $\bcX$ the analogue of Proposition
\ref{ds9prop2} for orbifolds holds. This is stronger, as it
orients $\bcX^{\Ga,\la},\ldots, \bhcX{}^{\Ga,\mu}_\ci$ under
weaker conditions on $\Ga,\la,\mu$, which allow $\md{\Ga}$ even
for some~$\la,\mu$.\I{d-orbifold|)}%
\I{d-orbifold!semieffective|)}\I{d-orbifold!effective|)}
\end{itemize}

\section{Orbifolds with corners}
\label{ds12}
\I{orbifold with corners|(}

In \cite[\S 8.5--\S 8.9]{Joyc6} we discuss 2-categories $\Orbb$ and
$\Orbc$ of {\it orbifolds with boundary\/} and {\it orbifolds with
corners}, which are orbifold versions of manifolds with boundary and
with corners in \S\ref{ds5}. This is new material, and the author
knows of no other foundational work on orbifolds with corners.

\subsection{The definition of orbifolds with corners}
\label{ds121}
\I{orbifold with corners!definition|(}\I{2-category|(}

\begin{dfn} An {\it orbifold with corners\/} $\oX$ of dimension
$n\ge 0$ is a triple $\oX=(\cX,\cpX,i_\oX)$ where $\cX,\cpX$ are
separated, second countable Deligne--Mumford $C^\iy$-stacks, and
$i_\oX:\cpX\ra\cX$\G[iXc]{$\bs i_\oX:\pd\oX\ra\oX$}{inclusion of
boundary $\pd\oX$ into an orbifold with corners $\oX$} is a proper,
strongly representable\I{C-stack@$C^\iy$-stack!strongly
representable 1-morphism} 1-morphism of $C^\iy$-stacks, in the sense
of \S\ref{ds83}, such that for each $[x]\in\cX_\top$ there exists a
2-Cartesian\I{2-category!2-Cartesian square} diagram in $\CSta$:
\begin{equation*}
\xymatrix@C=150pt@R=10pt{ *+[r]{\bar\upU} \ar[r]_(0.3){u_\pd}
\ar[d]^{\bar\ui_U} \drtwocell_{}\omit^{}\omit{^{\id\,\,{}}} &
*+[l]{\cpX} \ar[d]_{i_\oX} \\
*+[r]{\bar\uU} \ar[r]^(0.65){u} & *+[l]{\cX.} }
\end{equation*}
Here $U$ is an $n$-manifold with corners, so that $i_U:\pd U\ra U$
is smooth, and $\uU,\upU,\ui_U=F_\Manc^\CSch(U,\pd U,i_U)$, and
$u,u_\pd$ are \'etale 1-morphisms, and $u_\top([p])=[x]$ for some
$p\in U$. We call $\oX$ an {\it orbifold with boundary},\I{orbifold
with boundary} or an {\it orbifold without boundary},\I{orbifold!as
orbifold with corners} if the above condition holds with $U$ a
manifold with boundary, or a manifold without boundary,
respectively, for
each~$[x]\in\cX_\top$.\G[WXYZd]{$\oW,\oX,\oY,\oZ,\ldots$}{orbifolds
with corners}

Now suppose $\oX=(\cX,\cpX,i_\oX)$ and $\oY=(\cY,\cpY,i_\oY)$ are
orbifolds with corners. A 1-{\it
morphism\/}\I{2-category!1-morphism} $f:\oX\ra\oY$, or {\it smooth
map}, is a 1-morphism of $C^\iy$-stacks $f:\cX\ra\cY$ such that for
each $[x]\in\cX_\top$ with $f_\top([x])=[y]\in\cY_\top$ there exists
a 2-commutative diagram\I{2-category!2-commutative diagram} in
$\CSta$:
\begin{equation*}
\xymatrix@C=150pt@R=10pt{ *+[r]{\bar\uU} \ar[r]_(0.3){u}
\ar[d]^{\bar\uh} \drtwocell_{}\omit^{}\omit{_{\,\,\eta}} &
*+[l]{\cX} \ar[d]_{f} \\
*+[r]{\bar\uV} \ar[r]^(0.65){v} & *+[l]{\cY.\!{}} }
\end{equation*}
Here $U,V$ are manifolds with corners, $h:U\ra V$ is a smooth map,
$\uU,\uV,\uh=F_\Manc^\CSch(U,V,h)$, and $u,v$ are \'etale, and
$u_\top([p])=[x]$ for some~$p\in U$.

Let $f,g:\oX\ra\oY$ be 1-morphisms of orbifolds with corners. A
2-{\it morphism\/}\I{2-category!2-morphism} $\eta:f\Ra g$ is a
2-morphism of 1-morphisms $f,g:\cX\ra\cY$ in~$\CSta$.

{\it Composition of\/ $1$-morphisms\/} $g\ci f$, {\it identity
$1$-morphisms\/} $\id_\oX$, {\it vertical\/} and {\it horizontal
composition of\/ $2$-morphisms\/} $\ze\od\eta$, $\ze*\eta$, and {\it
identity\/ $2$-morphisms\/} for orbifolds with corners, are all
given by the corresponding compositions and identities in $\CSta$.
This defines the 2-category $\Orbc$\G[Orbc]{$\Orbc$}{2-category of
orbifolds with corners} of orbifolds with corners.\I{2-category|)}
Write $\Orbb$\G[Orbb]{$\Orbb$}{2-category of orbifolds with
boundary} and $\dotOrb$\G[Orb'']{$\dotOrb$}{2-subcategory of
orbifolds with corners equivalent to orbifolds} for the full
2-subcategories of orbifolds with boundary, and orbifolds without
boundary, in~$\Orbc$.

If $\cX$ is an orbifold in the sense of Definition \ref{ds9def1},
then $\oX=(\cX,\es,\es)$ is an orbifold without boundary in this
sense, and vice versa. Thus the 2-functor
$F_\Orb^\Orbc:\Orb\ra\Orbc$ mapping $\cX\mapsto\oX=(\cX,\es,\es)$ on
objects, $f\mapsto f$ on 1-morphisms, and $\eta\mapsto\eta$ on
2-morphisms, is an isomorphism of 2-categories~$\Orb\ra\dotOrb$.

Define $F_\Manc^\Orbc:\Manc\ra\Orbc$ by $F_\Manc^\Orbc:
X\mapsto\oX=(\ul{\bar X\!}\,,\bar\upX,\bar\ui_X)$ on objects $X$ in
$\Manc$, where $\uX,\upX,\ui_X=F_\Manc^\CSch(X,\pd X,i_X)$, and
$F_\Manc^\Orbc:f\mapsto\ul{\bar f\!}\,$ on morphisms $f:X\ra Y$ in
$\Manc$, where $\uf=F_\Manc^\CSch(f)$. Then $F_\Manc^\Orbc$ is a
full and faithful\I{functor!full}\I{functor!faithful} strict
2-functor.\I{2-category!strict 2-functor}

Let $\oX=(\cX,\cpX,i_\oX)$ be an orbifold with corners, and
$\cV\subseteq\cX$ an open $C^\iy$-substack. Define
$\cpV=i_\oX^{-1}(\cV)$, as an open $C^\iy$-substack of $\cpX$, and
$i_\oV:\cpV\ra\cV$ by $i_\oV=i_\oX\vert_\cpV$. Then
$\oV=(\cV,\cpV,i_\eV)$ is an orbifold with corners. We call $\oV$ an
{\it open suborbifold\/}\I{orbifold with corners!open suborbifold}
of $\oX$. An {\it open cover\/}\I{orbifold with corners!open cover}
of $\oX$ is a family $\{\oV_a:a\in A\}$ of open suborbifolds $\oV_a$
of $\oX$ with~$\cX=\bigcup_{a\in A}\cV_a$.\I{orbifold with
corners!definition|)}
\label{ds12def1}
\end{dfn}

\begin{ex} Suppose $X$ is a manifold with corners, $G$ a finite
group, and $r:G\ra\Aut(X)$ an action of $G$ on $X$ by
diffeomorphisms. Since $r(\ga):X\ra X$ is simple for each $\ga\in
G$, as in \S\ref{ds52} we have $r_-(\ga):\pd X\ra\pd X$, which is
also a diffeomorphism. Then $r_-:G\ra\Aut(\pd X)$ is an action of
$G$ on $\pd X$, and $i_X:\pd X\ra X$ is $G$-equivariant. Set
$\uX,\upX,\ui_X,\ur,\ur_-=F_\Manc^\CSch(X,\pd X,i_X,r,r_-)$. Then
$\uX,\upX$ are $C^\iy$-schemes with $G$-actions $\ur,\ur_-$, and
$\ui_X:\upX\ra\uX$ is $G$-equivariant, so Examples \ref{ds8ex1} and
\ref{ds8ex2} define Deligne--Mumford $C^\iy$-stacks $[\uX/G],
[\upX/G]$ and a 1-morphism $[\ui_X,\id_G]:[\upX/G]\ra[\uX/G]$, which
turns out to be strongly representable. One can show that
$\oX=\bigl([\uX/G],\ab [\upX/G],\ab[\ui_X,\id_G]\bigr)$ is an
orbifold with corners, which we will write as~$[X/G]$.\I{orbifold
with corners!quotients $[X/G]$}
\label{ds12ex1}
\end{ex}

\begin{rem}{\bf(a)} We could have defined $\Orbc$ equivalently and
more simply as a (non-full) 2-subcategory of $\DMCSta$, so that an
orbifold with corners would be a $C^\iy$-stack $\cX$ rather than a
triple $\oX=(\cX,\cpX,i_\oX)$. We chose the set-up of Definition
\ref{ds12def1} partly for its compatibility with the definitions of
d-stacks and d-orbifolds with corners $\eX=(\bcX,\bcpX,\bs
i_\eX,\om_\eX)$ in \S\ref{ds13}--\S\ref{ds14}, and partly because,
to make several important constructions more functorial, it is
useful to have a particular choice of boundary $\cpX$ for $\cX$
already made.
\smallskip

\noindent{\bf(b)} In Remark \ref{ds6rem2} we noted that boundaries
in $\dSpac$ are strictly functorial. One sign of this is that for a
semisimple 1-morphism $\bs f:\rX\ra\rY$ in $\dSpac$, the 1-morphism
$\bs f_-:\pd_-^{\bs f}\rX\ra\pd\rY$ is unique, not just unique up to
2-isomorphism, with an equality of 1-morphisms $\bs f\ci\bs i_\rX
\vert_{\smash{\pd_-^{\bs f}\rX}}=\bs i_\rY\ci\bs f_-$, not just a
2-isomorphism. By the general philosophy of 2-categories, this may
seem unnatural.

We will arrange that boundaries in $\Orbc$ and also in
$\dStac,\dOrbc$ are strictly functorial\I{orbifold with
corners!boundary!strictly functorial} in the same way. This is our
reason for taking $i_\oX:\cpX\ra\cX$ in Definition \ref{ds12def1} to
be {\it strongly representable},\I{C-stack@$C^\iy$-stack!strongly
representable 1-morphism}
 in the sense of \S\ref{ds83}.
Proposition \ref{ds8prop1}(b) shows that this is no real
restriction: $i_\oX:\cpX\ra\cX$ is naturally representable, and we
can make it strongly representable by replacing $\cpX$ by an
equivalent $C^\iy$-stack. Then Proposition \ref{ds8prop2} applied to
$i_\oY:\cpY\ra\cY$ is what we need to show that a semisimple
1-morphism $f:\oX\ra\oY$ in $\dOrbc$ lifts to a unique 1-morphism
$f_-:\pd_-^f\oX\ra\pd\oY$ with~$f\ci i_\oX
\vert_{\smash{\pd_-^f\oX}}=i_\oY\ci f_-$.
\smallskip

\noindent{\bf(c)} An orbifold with corners $\oX$ of dimension $n$ is
locally modelled near each point $[x]\in\cX_\top$ on
$\bigl([0,\iy)^k\t\R^{n-k} \bigr)/G$ near $0$, where $G$ is a finite
group acting linearly on $\R^n$ preserving the subset
$[0,\iy)^k\t\R^{n-k}$. Note that $G$ is allowed to permute the
coordinates $x_1,\ldots,x_k$ in $[0,\iy)^k$. So, for example, we
allow 2-dimensional orbifolds with corners modelled on
$[0,\iy)^2/\Z_2$, where $\Z_2=\an{\si}$ acts on $[0,\iy)^2$
by~$\si:(x_1,x_2)\mapsto(x_2,x_1)$.

This implies that the 1-morphism $i_\oX:\cpX\ra\cX$ induces
morphisms of orbifold groups $(i_\oX)_*:\Iso_\cpX([x'])\ra
\Iso_\cX([x])$\I{C-stack@$C^\iy$-stack!orbifold group
$\Iso_\cX([x])$} which are injective (so that $i_\oX$ is
representable), but need not be isomorphisms. We will call an
orbifold with corners $\oX$ {\it straight\/}\I{orbifold with
corners!straight} if the morphisms $(i_\oX)_*:\Iso_\cpX([x'])\ra
\Iso_\cX([x])$ are isomorphisms for all $[x']\in\cpX_\top$ with
$i_{\oX,\top}([x'])=[x]$. That is, straight orbifolds with corners
are locally modelled on $[0,\iy)^k\t(\R^{n-k}/G)$. Orbifolds with
boundary,\I{orbifold with boundary} with $k=0$ or 1, are
automatically straight. Boundaries of orbifold strata behave better
for straight orbifolds with corners.
\label{ds12rem}
\end{rem}

In \S\ref{ds91} we explained that a  vector bundle $\cE$ on an
orbifold $\cX$ is a vector bundle on $\cX$ as a Deligne--Mumford
$C^\iy$-stack, in the sense of \S\ref{ds86}. But sometimes it is
convenient to regard $\cE$ as an orbifold in its own right, so we
define a `total space functor' mapping vector bundles $\cE$ to
orbifolds $\Tot(\cE)$.

In the same way, if $\oX=(\cX,\cpX,i_\oX)$ is an orbifold with
corners, in \cite[\S 8.5]{Joyc6} we define a {\it vector bundle\/
$\cE$ on\/}\I{orbifold with corners!vector bundles on} $\oX$ to be a
vector bundle on $\cX$ as a Deligne--Mumford $C^\iy$-stack. To
regard $\cE$ as an orbifold with corners in its own right, we define
a `total space functor'\I{orbifold with corners!vector bundles
on!total space functor $\Totc$} $\Totc:\vect(\cX)\ra\Orbc$, which
maps a vector bundle $\cE$ on $\oX$ to an orbifold with corners
$\Totc(\cE)$, and maps a section $s\in C^\iy(\cE)$ to a simple, flat
1-morphism $\Totc(s):\oX\ra\Totc(\cE)$ in~$\Orbc$.

\begin{dfn} An orbifold with corners $\oX$ is called {\it
effective\/}\I{orbifold with corners!effective} if $\oX$ is locally
modelled near each $[x]\in\cX_\top$ on $([0,\iy)^k\t\R^{n-k})/G,$
where $G$ acts effectively on $\R^n$ preserving
$[0,\iy)^k\t\R^{n-k}$, that is, every $1\ne\ga\in G$ acts
nontrivially.
\label{ds12def2}
\end{dfn}

The analogue of Proposition \ref{ds9prop1} holds for effective
orbifolds with corners.

\subsection{Boundaries of orbifolds with corners, and \\ simple,
semisimple and flat 1-morphisms}
\label{ds122}
\I{orbifold with corners!boundary|(}\I{orbifold with corners!simple
1-morphism|(}\I{orbifold with corners!semisimple
1-morphism|(}\I{orbifold with corners!flat 1-morphism|(}

In \cite[\S 8.6]{Joyc6} we define boundaries of orbifolds with
corners.

\begin{dfn} Let $\oX=(\cX,\cpX,i_\oX)$ be an orbifold with
corners. We will define an orbifold with corners $\pd\oX=(\cpX,
\cptX,i_{\pd\oX})$,\G[dXc]{$\pd\oX$}{boundary of an orbifold with
corners $\oX$} called the {\it boundary\/}\I{boundary!of an orbifold
with corners}\I{orbifold with corners!boundary} of $\oX$, such that
$i_\oX:\pd\oX\ra\oX$ is a 1-morphism in $\Orbc$. Here $\cpX$ and
$i_\oX$ are given in $\oX$, so the new data we have to construct
is~$\cptX,i_{\pd\oX}$.

As $i_\oX:\cpX\ra\cX$ is strongly
representable\I{C-stack@$C^\iy$-stack!strongly representable
1-morphism} by Definition \ref{ds12def1}, Proposition \ref{ds8prop3}
defines an explicit fibre product $\cpX\t_{i_\oX,\cX,i_\oX}\cpX$
with strongly representable projection morphisms
$\pi_1,\pi_2:\cpX\t_\cX\cpX\ra\cpX$ such that
$i_\oX\ci\pi_1=i_\oX\ci\pi_2$. We will use this explicit fibre
product throughout. There is a unique diagonal 1-morphism
$\De_\cpX:\cpX\ra\cpX\t_\cX\cpX$ with
$\pi_1\ci\De_\cpX=\pi_2\ci\De_\cpX=\id_\cpX$. It is an equivalence
with an open and closed $C^\iy$-substack
$\De_\cpX(\cpX)\subseteq\cpX\t_\cX\cpX$. Define $\cptX=
\cpX\t_\cX\cpX\sm\De_\cpX(\cpX)$. Then $\cptX$ is also an open and
closed $C^\iy$-substack in $\cpX\t_\cX\cpX$. Define
$i_{\pd\oX}=\pi_1\vert_\cptX: \cptX\ra\cpX$. Then
$\pd\oX=(\cpX,\cptX,i_{\pd\oX})$ is an orbifold with corners, with
$\dim(\pd\oX)=\dim\oX-1$. Also $i_\oX:\cpX\ra\cX$ in $\oX$ is a
1-morphism $i_\oX:\pd\oX\ra\oX$ in~$\Orbc$.\I{orbifold with
corners!boundary|)}
\label{ds12def3}
\end{dfn}

Here is the orbifold analogue of parts of
\S\ref{ds51}--\S\ref{ds52}.

\begin{dfn} Let $\oX=(\cX,\cpX,i_\oX)$ and $\oY=(\cY,\cpY,i_\oY)$ be
orbifolds with corners, and $f:\oX\ra\oY$ a 1-morphism in $\Orbc$.
Consider the $C^\iy$-stack fibre products $\cpX\t_{f\ci i_\oX,\cY,
i_\oY}\cpY$ and $\cX\t_{f,\cY, i_\oY}\cpY$. Since $i_\oY$ is
strongly representable, we may define these using the explicit
construction of Proposition~\ref{ds8prop3}.

The topological space $(\cpX\t_\cY\cpY)_\top$ associated to the
$C^\iy$-stack $\cpX\t_\cY\cpY$ may be written explicitly as
\e
\begin{split}
(\cpX\t_\cY\cpY)_\top\cong\bigl\{&[x',y']:\text{$x':\bar{\ul
*}\ra\cpX$ and $y':\bar{\ul *}\ra\cpY$ are}\\
&\text{1-morphisms with $f\!\ci\! i_\oX\!\ci\! x'\!=\! i_\oY\!\ci
\!y':\bar{\ul *}\!\ra\!\cY$}\bigr\},
\end{split}
\label{ds12eq1}
\e
where $[x',y']$ in \eq{ds12eq1} denotes the $\sim$-equivalence class
of pairs $(x',y')$, with $(x',y')\sim(\ti x',\ti y')$ if there exist
2-morphisms $\eta:x'\Ra\ti x'$ and $\ze:y'\Ra\ti y'$ with $\id_{f\ci
i_\oX}*\eta=\id_{i_\oY}*\ze$. There is a natural open and closed
$C^\iy$-substack $\cS_f\subseteq\cpX\t_\cY\cpY$,\G[Sfc]{$\cS_f
\subseteq\cpX\t_\cY\cpY$}{$C^\iy$-stack associated to 1-morphism
$f:\oX\ra\oY$ in $\Orbc$} the analogue of $S_f$ in \S\ref{ds51},
such that $[x',y']$ in \eq{ds12eq1} lies in $\cS_{f,\top}$ if and
only if we can complete the following commutative diagram in
$\qcoh(\ul{\bar *})$ with morphisms `$\dashra$' as shown:
\begin{equation*}
\xymatrix@C=13pt@R=25pt{ 0 \ar[r] & (y')^*(\cN_\oY)
\ar[rr]^(0.45){(y')^*(\nu_\oY)} \ar@{.>}[d]_\cong && (y')^*\!\ci\!
i_\oY^*(T^*\cY) \ar[rr]^{(y')^*(\Om_{i_\oY})}
\ar@<-7ex>[d]^{\begin{subarray}{l}
I_{x',i_\oX}(T^*\cX)\ci (i_\oX\ci x')^*(\Om_f) \ci \\
I_{i_\oX\ci x',f}(T^*\cY)\ci
I_{y',i_\oY}(T^*\cY)^{-1}\end{subarray}} &&
(y')^*(T^*(\cpY)) \ar@{.>}[d] \ar[r] & 0 \\
0 \ar[r] & (x')^*(\cN_\oX) \ar[rr]_(0.45){(x')^*(\nu_\oX)} &&
(x')^*\!\ci\! i_\oX^*(T^*\cX) \ar[rr]_{(x')^*(\Om_{i_\oX})} &&
(x')^*(T^*(\cpX)) \ar[r] & 0, }
\end{equation*}
where $\cN_\oX,\cN_\oY$ are the {\it conormal line bundles\/} of
$\cpX,\cpY$ in~$\cX,\cY$.

Similarly, the topological space $(\cX\t_\cY\cpY)_\top$ may be
written explicitly as
\e
\begin{split}
(\cX\t_\cY\cpY)_\top\cong\bigl\{&[x,y']:\text{$x:\bar{\ul
*}\ra\cX$ and $y':\bar{\ul *}\ra\cpY$ are}\\
&\text{1-morphisms with $f\!\ci\! x\!=\! i_\oY\!\ci \!y':\bar{\ul
*}\!\ra\!\cY$}\bigr\},
\end{split}
\label{ds12eq2}
\e
where $[x,y']$ in \eq{ds12eq2} denotes the $\approx$-equivalence
class of $(x,y')$, with $(x,y')\approx(\ti x,\ti y')$ if there exist
$\eta:x\Ra\ti x$ and $\ze:y'\Ra\ti y'$ with
$\id_f*\eta=\id_{i_\oY}*\ze$. There is a natural open and closed
$C^\iy$-substack $\cT_f\subseteq\cX\t_\cY\cpY$,\G[Tfc]{$\cT_f
\subseteq\cX\t_\cY\cpY$}{$C^\iy$-stack associated to 1-morphism
$f:\oX\ra\oY$ in $\Orbc$} the analogue of $T_f$ in \S\ref{ds51},
such that $[x,y']$ in \eq{ds12eq2} lies in $\cT_{f,\top}$ if and
only if we can complete the following commutative diagram in
$\qcoh(\ul{\bar *})$:
\begin{equation*}
\xymatrix@C=13pt@R=20pt{ 0 \ar[r] & (y')^*(\cN_\oY)
\ar[rr]_(0.45){(y')^*(\nu_\oY)} && (y')^*\!\ci\! i_\oY^*(T^*\cY)
\ar[rr]_{(y')^*(\Om_{i_\oY})} \ar@<2ex>[d]_(0.6){\begin{subarray}{l}
x^*(\Om_f) \ci I_{x,f}(T^*\cY)\ci
I_{y',i_\oY}(T^*\cY)^{-1}\end{subarray}} &&
(y')^*(T^*(\cpY)) \ar@<1ex>@{.>}[dll] \ar[r] & 0 \\
&&& x^*(T^*\cX). }
\end{equation*}

Define $s_f=\pi_\cpX\vert_{\cS_f}: \cS_f\ra\cpX$,
$u_f=\pi_\cpY\vert_{\cS_f}:\cS_f\ra\cpY$,
$t_f=\pi_\cX\vert_{\cT_f}:\cT_f\ra\cX$, and
$v_f=\pi_\cpY\vert_{\cT_f}:\cT_f\ra\cpY$. Then $s_f,t_f$ are proper,
\'etale 1-morphisms. We call $f$ {\it simple\/} if
$s_f:\cS_f\ra\cpX$ is an equivalence, and we call $f$ {\it
semisimple\/} if $s_f:\cS_f\ra\cpX$ is injective as a 1-morphism of
Deligne--Mumford $C^\iy$-stacks, and we call $f$ {\it flat\/} if
$\cT_f=\es$. Simple implies semisimple.
\label{ds12def4}
\end{dfn}

The condition that $i_\oX$ is strongly
representable\I{C-stack@$C^\iy$-stack!strongly representable
1-morphism} in Definition \ref{ds12def1} is essential in
constructing $f_-,\eta_-$ in parts (b),(c) of the next theorem.

\begin{thm} Let\/ $f:\oX\ra\oY$ be a semisimple $1$-morphism of
orbifolds with corners. Then there is a natural decomposition
$\pd\oX=\pd^f_+\oX\amalg\pd^f_-\oX,$ where $\pd^f_\pm\oX$ are open
and closed suborbifolds in $\pd\oX,$ such that:\G[dfXc]{$\pd^f_\pm
\oX$}{sets of decomposition $\pd\oX=\pd^f_+\oX\amalg\pd^f_-\oX$ of
boundary $\pd\oX$ induced by 1-morphism $f:\oX\ra\oY$ in $\Orbc$}
\begin{itemize}
\setlength{\itemsep}{0pt}
\setlength{\parsep}{0pt}
\item[{\bf(a)}] Define $f_+=f\ci i_\oX\vert_{\pd^f_+\oX}:
\pd_+^f\oX\ra\oY$. Then $f_+$ is semisimple. If\/ $f$ is flat
then $f_+$ is also flat.
\item[{\bf(b)}] There exists a unique, semisimple\/ $1$-morphism
$f_-:\pd_-^f\oX\ra\pd\oY$ in $\Orbc$ with\/ $f\ci
i_\oX\vert_{\pd_-^f\oX}=i_\oY\ci f_-$. If\/ $f$ is simple then\/
$\pd^f_+\oX=\es,$ $\pd^f_-\oX=\pd\oX$ and\/ $f_-:\pd\oX\ra
\pd\oY$ is simple. If\/ $f$ is flat then $f_-$ is flat.
\item[{\bf(c)}] Let\/ $g:\oX\ra\oY$ be another $1$-morphism
and\/ $\eta:f\Ra g$ a $2$-morphism in $\Orbc$. Then $g$ is also
semisimple, with\/ $\pd_-^g\oX=\pd_-^f\oX$. If\/ $f$ is simple,
or flat, then $g$ is too. Part\/ {\bf(b)} defines $1$-morphisms
$f_-,g_-:\pd_-^f\oX\ra\pd\oY$. There is a unique $2$-morphism
$\eta_-:f_-\Ra g_-$ in\/ $\Orbc$ such that\I{orbifold with
corners!simple 1-morphism|)}\I{orbifold with corners!semisimple
1-morphism|)}\I{orbifold with corners!flat 1-morphism|)}
\end{itemize}\vskip -5pt
\begin{equation*}
\id_{i_\oY}*\eta_-\!=\!\eta*\id_{i_\oX\vert_{\pd_-^f\oX}}:
f\!\ci\!i_\oX\vert_{\pd_-^f\oX}\!=\!i_\oY\!\ci\!f_-\!\Longra\!
g\!\ci\!i_\oX\vert_{\pd_-^f\oX}\!=\!i_\oY\!\ci\!g_-.
\end{equation*}
\label{ds12thm1}
\end{thm}

\subsection{Corners $C_k(\oX)$ and the corner functors $C,\hat C$}
\label{ds123}

In \cite[\S 8.7]{Joyc6} we extend \S\ref{ds53} to orbifolds. Here is
the orbifold analogue of the category $\cManc$ in
Definition~\ref{ds5def5}.

\begin{dfn} We will define a 2-category $\cOrbc$ whose objects are
disjoint unions $\coprod_{m=0}^\iy \oX_m$, where $\oX_m$ is a
(possibly empty) orbifold with corners of dimension $m$. In more
detail, objects of $\cOrbc$ are triples $\oX=(\cX,\cpX,i_\oX)$ with
$i_\oX:\cpX\ra\cX$ a strongly representable 1-morphism of
Deligne--Mumford $C^\iy$-stacks, such that there exists a
decomposition $\cX=\coprod_{m=0}^\iy\cX_m$ with each
$\cX_m\subseteq\cX$ an open and closed $C^\iy$-substack, for which
$\oX_m:=\bigl(\cX_m,i_\oX^{-1}(\cX_m),i_\oX\vert_{\smash{i_\oX^{-1}
(\cX_m)}}\bigr)$ is an orbifold with corners of dimension $m$.

A 1-morphism $f:\oX\ra\oY$ in $\cOrbc$ is a 1-morphism $f:\cX\ra\cY$
in $\CSta$ such that $f\vert_{\cX_m\cap
f^{-1}(\cY_n)}:\bigl(\cX_m\cap f^{-1}(\cY_n)\bigr)\ra \cY_n$ is a
1-morphism in $\Orbc$ for all $m,n\ge 0$. For 1-morphisms
$f,g:\oX\ra\oY$, a 2-morphism $\eta:f\Ra g$ is a 2-morphism
$\eta:f\Ra g$ in $\CSta$. Then $\Orbc$ is a full 2-subcategory
of~$\cOrbc$.\G[Orbc']{$\cOrbc$}{2-category of disjoint unions of
orbifolds with corners of different dimensions}
\label{ds12def5}
\end{dfn}

The next theorem summarizes our results on corners functors
in~$\Orbc$.\I{orbifold with corners!k-corners C_k(X)@$k$-corners
$C_k(\oX)$|(}\I{orbifold with corners!corner functors|(}

\begin{thm}{\bf(a)} Suppose\/ $\oX$ is an orbifold with corners.
Then for each\/ $k=0,1,\ldots,\dim\oX$ we can define an orbifold
with corners $C_k(\oX)$ of dimension $\dim\oX-k$ called the
$k$-\begin{bfseries}corners\end{bfseries} of\/ $\oX,$ and a\/
$1$-morphism $\Pi^k_\oX:C_k(\oX)\ra\oX$ in $\Orbc$. It has
topological space
\e
\begin{split}
C_k(\cX)_\top\!\cong\!\bigl\{&[x,\{x_1',\ldots,x_k'\}]:
\text{$x:\bar{\ul{*}}\!\ra\!\cX,$ $x_i':\bar{\ul{*}}\!\ra\!\cpX$
are $1$-morphisms}\\
&\text{with\/ $x_1',\ldots,x_k'$ distinct and\/ $x=i_\oX\ci
x_1'=\cdots=i_\oX\ci x_k'$}\bigr\}.
\end{split}
\label{ds12eq3}
\e
There is a natural action of the symmetric group $S_k$ on $\pd^k\oX$
by $1$-isomorphisms, and an equivalence
$C_k(\oX)\simeq\pd^k\oX/S_k$. We have $1$-isomorphisms\/
$C_0(\oX)\cong\oX$ and\/ $C_1(\oX)\cong\pd\oX$ in $\Orbc$. Write\/
$C(\oX)=\coprod_{k=0}^{\dim\oX} C_k(\oX)$ and\/
$\Pi_\oX=\coprod_{k=0}^{\dim\oX}\Pi^k_\oX,$ so that\/ $C(\oX)$ is an
object and\/ $\Pi_\oX:C(\oX)\ra\oX$ a $1$-morphism in\/~$\cOrbc$.
\smallskip

\noindent{\bf(b)} Let\/ $f:\oX\ra\oY$ be a $1$-morphism of orbifolds
with corners. Then there is a unique $1$-morphism $C(f):C(\oX)\ra
C(\oY)$ in $\cOrbc$ such that\/ $\Pi_\oY\ci C(f)=f\ci\Pi_\oX:
C(\oX)\ra\oY,$ and\/ $C(f)$ acts on points as in \eq{ds12eq3} by
\e
\begin{split}
C(f)_\top:\bigl[x,\{x_1',\ldots,x_k'\}\bigr]\longmapsto
\bigl[y,\{y_1',\ldots,y_l'\}\bigr],\quad\text{where $y=f\ci x,$}&\\
\text{and\/}\quad\{y_1',\ldots,y_l'\}=\bigl\{y':[x_i',y']\in
\cS_{f,\top},\; \text{some\/ $i=1,\ldots,k\bigr\},$}&
\end{split}
\label{ds12eq4}
\e
where $\cS_f$ is as in Definition\/~{\rm\ref{ds12def4}}.

For all\/ $k,l\ge 0,$ write\/ $C_k^{f,l}(\oX)=C_k(\oX)\cap
C(f)^{-1}(C_l(\oY)),$ so that\/ $C_k^{f,l}(\oX)$ is open and closed
in\/ $C_k(\oX)$ with\/ $C_k(\oX)=\coprod_{l=0}^{\dim\oY} C_k^{f,l}
(\oX),$ and write $C^l_k(f)=C(f)\vert_{C_k^{f,l}(\oX)},$ so that\/
$C^l_k(f):C_k^{f,l}(\oX)\ra C_l(\oY)$ is a $1$-morphism in~$\Orbc$.

\noindent{\bf(c)} Let\/ $f,g:\oX\ra\oY$ be $1$-morphisms and\/
$\eta:f\Ra g$ a $2$-morphism in $\Orbc$. Then there exists a unique
$2$-morphism $C(\eta):C(f)\Ra C(g)$ in $\cOrbc,$ where $C(f),C(g)$
are as in {\bf(b)\rm,} such that
\begin{equation*}
\id_{\Pi_\oY}*
C(\eta)=\eta*\id_{\Pi_\oX}:\Pi_\oY\ci C(f)=
f\ci\Pi_\oX\Longra \Pi_\oY\ci C(g)=g\ci\Pi_\oX.
\end{equation*}

\noindent{\bf(d)} Define $C:\Orbc\ra\cOrbc$ by $C:\oX\mapsto C(\oX)$
on objects, $C:f\mapsto C(f)$ on\/ $1$-morphisms, and\/ $C:\eta
\mapsto C(\eta)$ on $2$-morphisms, where $C(\oX),C(f),C(\eta)$ are
as in {\bf(a)--(c)} above. Then $C$ is a strict\/ $2$-functor,
called a \begin{bfseries}corner
functor\end{bfseries}.\G[CCc]{$C,\hat C:\Orbc\ra\cOrbc$}{`corner
functors' for orbifolds with corners}
\smallskip

\noindent{\bf(e)} Let\/ $f:\oX\ra\oY$ be semisimple.\I{orbifold with
corners!semisimple 1-morphism} Then $C(f)$ maps
$C_k(\oX)\ra\coprod_{l=0}^kC_l(\oY)$ for all\/ $k\ge 0$. The
natural\/ $1$-isomorphisms $C_1(\oX)\cong\pd\oX,$
$C_0(\oY)\cong\oY,$ $C_1(\oY)\cong\pd\oY$ identify
$C_1^{f,0}(\oX)\cong\pd^f_+\oX,$ $C_1^{f,1}(\oX)\cong\pd^f_-\oX,$
$C_1^0(f)\cong f_+$ and\/~$C_1^1(f)\cong f_-$.

If\/ $f$ is simple\I{orbifold with corners!simple 1-morphism} then
$C(f)$ maps $C_k(\oX)\ra C_k(\oY)$ for all\/~$k\ge 0$.
\smallskip

\noindent{\bf(f)} Analogues of\/ {\bf(b)}--{\bf(d)} also hold for a
second corner functor\/ $\hat C:\Orbc\ra\cOrbc,$ which acts on
objects by $\hat C:\oX\mapsto C(\oX)$ in {\bf(a)\rm,} and for
$1$-morphisms $f:\oX\ra\oY$ in {\bf(b)\rm,} $\hat C(f):C(\oX)\ra
C(\oY)$ acts on points by
\begin{align*}
&\hat C(f)_\top:\bigl[x,\{x_1',\ldots,x_k'\}\bigr]\longmapsto
\bigl[y,\{y_1',\ldots,y_l'\}\bigr],\quad\text{where $y=f\ci x,$}\\
&\{y_1',\ldots,y_l'\}\!=\!\bigl\{y':[x_i',y']\!\in\!\cS_{f,\top},\;
i\!=\!1,\ldots,k\bigr\}\!\cup\!\bigl\{y':[x,y']\!\in\!
\cT_{f,\top}\bigr\}.
\end{align*}

If\/ $f$ is flat\I{orbifold with corners!flat 1-morphism}
then\/~$\hat C(f)=C(f)$.\I{orbifold with corners!k-corners
C_k(X)@$k$-corners $C_k(\oX)$|)}\I{orbifold with corners!corner
functors|)}
\label{ds12thm2}
\end{thm}

\begin{ex} Suppose $\oX$ is a quotient $[X/G]$\I{orbifold
with corners!quotients $[X/G]$} as in Example \ref{ds12ex1}, where
$X$ is a manifold with corners and $G$ is a finite group. Then the
action $r:G\ra\Aut(X)$ lifts to $C(r):G\ra\Aut(C(X))$, and there is
an equivalence $C([X/G])\simeq [C(X)/G]$ in $\cOrbc$, where to
define $[C(X)/G]$ we note that Example \ref{ds12ex1} also works with
$X$ in $\cManc$ rather than $\Manc$, yielding~$[X/G]\in\cOrbc$.
\label{ds12ex2}
\end{ex}

Section \ref{ds52} defined (s-)submersions, (s- or sf-)immersions
and \hbox{(s- or sf-)}\ab embeddings in $\Manc$. Section \ref{ds91}
defined submersions, immersions and embeddings in $\Orb$. We combine
the two definitions.

\begin{dfn} Let $f:\oX\ra\oY$ be a 1-morphism of orbifolds with
corners.
\begin{itemize}
\setlength{\itemsep}{0pt}
\setlength{\parsep}{0pt}
\item[(i)] We call $f$ a {\it submersion\/}\I{orbifold with
corners!submersion} if $\Om_{C(f)}:C(f)^*(T^*C(\cY))\ra
T^*C(\cX)$ is an injective morphism of vector bundles, i.e.\ has
a left inverse in $\qcoh(C(\cX))$, and $f$ is semisimple and
flat. We call $f$ an {\it s-submersion\/}\I{orbifold with
corners!s-submersion} if $f$ is also simple.
\item[(ii)] We call $f$ an {\it immersion\/}\I{orbifold with
corners!immersion} if it is representable and
$\Om_f:f^*(T^*\cY)\ra T^*\cX$ is a surjective morphism of vector
bundles, i.e.\ has a right inverse in $\qcoh(\cX)$. We call $f$
an {\it s-immersion\/}\I{orbifold with corners!s-immersion} if
$f$ is also simple, and an {\it sf-immersion\/}\I{orbifold with
corners!sf-immersion} if $f$ is also simple and flat.
\item[(iii)] We call $f$ an {\it embedding, s-embedding}, or
{\it sf-embedding},\I{orbifold with
corners!embedding}\I{orbifold with
corners!s-embedding}\I{orbifold with corners!sf-embedding} if it
is an immersion, s-immersion, or sf-immersion, respectively, and
$f_*:\Iso_\cX([x])\ra\Iso_\cY \bigl(f_\top([x])\bigr)$ is an
isomorphism for all $[x]\in\cX_\top$, and
$f_\top:\cX_\top\ra\cY_\top$ is a homeomorphism with its image
(so in particular it is injective).
\end{itemize}
\label{ds12def6}
\end{dfn}

Then submersions, \ldots, sf-embeddings in $\Orbc$ are \'etale
locally modelled on submersions, \ldots, sf-embeddings in~$\Manc$.

\subsection{Transversality and fibre products}
\label{ds124}
\I{orbifold with corners!strongly transverse
1-morphisms|(}\I{orbifold with corners!transverse fibre products|(}

Section \ref{ds54} discussed transversality and fibre products for
manifolds with corners. In \cite[\S 8.8]{Joyc6} we generalize this
to orbifolds with corners.

\begin{dfn} Let $\oX,\oY,\oZ$ be orbifolds with corners and
$g:\oX\ra\oZ$, $h:\oY\ra\oZ$ be 1-morphisms. Then as in \S\ref{ds13}
we have 1-morphisms $C(g):C(\oX)\ra C(\oZ)$ and $C(h):C(\oY)\ra
C(\oZ)$ in $\cOrbc$, and hence 1-morphisms $C(g):C(\cX)\ra C(\cZ)$
and $C(h):C(\cY)\ra C(\cZ)$ in $\CSta$. We call $g,h$ {\it
transverse\/} if the following holds. Suppose $x:\bar{\ul *}\ra
C(\cX)$ and $y:\bar{\ul *}\ra C(\cY)$ are 1-morphisms in $\CSta$,
and $\eta:C(g)\ci x\Ra C(h)\ci y$ a 2-morphism. Then the following
morphism in $\qcoh(\bar{\ul*})$ should be injective:
\begin{align*}
\bigl(x^*(\Om_{C(g)})\ci &I_{x,C(g)}(T^*C(\cZ))\bigr)\!\op\!
\bigl(y^*(\Om_{C(h)})\!\ci\! I_{y,C(h)}(T^*C(\cZ))\!\ci\!
\eta^*(T^*C(\cZ))\bigr):\\
&(C(g)\ci x)^*(T^*C(\cZ)) \longra x^*(T^*C(\cX))\op y^*(T^*C(\cY)).
\end{align*}

Now identify $C_k(\cX)_\top\subseteq C(\cX)_\top$ with the right
hand of \eq{ds12eq3}, and similarly for $C(\cY)_\top,C(\cZ)_\top$.
Then $C(g)_\top,C(h)_\top$ act as in \eq{ds12eq4}. We call $g,h$
{\it strongly transverse\/} if they are transverse, and whenever
there are points in $C_j(\cX)_\top,C_k(\cY)_\top,C_l(\cZ)_\top$ with
\begin{equation*}
C(g)_\top\bigl([x,\{x_1',\ldots,x_j'\}]\bigr)\!=\!
C(h)_\top\bigl([y,\{y_1',\ldots,y_k'\}]\bigr)\!=\!
[z,\{z_1',\ldots,z_l'\}],
\end{equation*}
we have either $j+k>l$ or $j=k=l=0$.

One can show that $g,h$ are (strongly) transverse if and only if
they are \'etale locally equivalent to (strongly) transverse smooth
maps in~$\Manc$.
\label{ds12def7}
\end{dfn}

Here is the analogue of Theorem~\ref{ds5thm1}:

\begin{thm} Suppose $g:\oX\ra\oZ$ and\/ $h:\oY\ra\oZ$ are
transverse $1$-morphisms in $\Orbc$. Then a fibre product\/
$\oW=\oX\t_{g,\oZ,h}\oY$ exists in the $2$-category~$\Orbc$.\I{fibre
product!of orbifolds with corners}\I{2-category!fibre products in}
\label{ds12thm3}
\end{thm}

Proposition \ref{ds5prop4} and Theorem \ref{ds5thm2} also extend to
$\Orbc$, with equivalences natural up to 2-isomorphism rather than
canonical diffeomorphisms.\I{orbifold with corners!strongly
transverse 1-morphisms|)}\I{orbifold with corners!transverse fibre
products|)}

\subsection{Orbifold strata of orbifolds with corners}
\label{ds125}
\I{orbifold with corners!orbifold strata|(}\I{orbifold strata!of
orbifolds with corners|(}

Sections \ref{ds87} and \ref{ds92} discussed orbifold strata of
Deligne--Mumford $C^\iy$-stacks and orbifolds, respectively. In
\cite[\S 8.9]{Joyc6} we extend this to orbifolds with corners. This
is also related to the material on fixed points of finite group
actions on manifolds with corners in~\S\ref{ds56}.

\begin{thm} Let\/ $\oX$ be an orbifold with corners, and\/ $\Ga$
a finite group. Then we can define objects\/ $\oX^\Ga,\ti\oX{}^\Ga,
\hat\oX{}^\Ga$ in\/ $\cOrbc,$ and open subobjects\/
$\oX^\Ga_\ci\subseteq\oX^\Ga,$ $\ti\oX{}^\Ga_\ci\subseteq
\ti\oX^\Ga,$ $\hat\oX{}^\Ga_\ci\subseteq\hat\oX^\Ga,$ all natural up
to $1$-isomorphism in $\cOrbc,$ and\/ $1$-morphisms
$O^\Ga(\oX),\ti\Pi^\Ga(\oX),\ldots$ fitting into a strictly
commutative diagram in\/~$\cOrbc\!:$\G[OGaXc]{$O^\Ga(\oX),\ti
O^\Ga(\oX),O{}^\Ga_\ci(\oX),\ti O{}^\Ga_\ci(\oX)$}{1-morphisms of
orbifold strata $\oX^\Ga,\ldots,\hat\oX{}^\Ga_\ci$ of an orbifold
with corners $\oX$}\G[PiGaXc]{$\ti\Pi^\Ga(\oX),\hat\Pi{}^\Ga(\oX),
\ti\Pi^\Ga_\ci(\oX),\hat\Pi{}^\Ga_\ci(\oX)$}{1-morphisms of orbifold
strata $\oX^\Ga,\ldots,\hat\oX{}^\Ga_\ci$ of an orbifold with
corners $\oX$}
\e
\begin{gathered}
\xymatrix@C=48pt@R=7pt{ \oX^\Ga_\ci \ar[rr]^{\ti\Pi{}^\Ga_\ci(\oX)}
\ar[dr]_(0.3){O^\Ga_\ci(\oX)} \ar[dd]_\subset
\ar@(ul,l)[]_(0.8){\Aut(\Ga)} && \ti\oX{}^\Ga_\ci
\ar[r]^{\hat\Pi{}^\Ga_\ci(\oX)} \ar[dl]^(0.3){\ti O{}^\Ga_\ci(\oX)}
\ar[dd]^\subset & {\hat\oX{}^\Ga_\ci}
\ar@<.5ex>[dd]^\subset \\ & \oX \\
\oX^\Ga \ar[rr]_{\ti\Pi{}^\Ga(\oX)} \ar[ur]^(0.3){O^\Ga(\oX)}
\ar@(dl,l)[]^(0.8){\Aut(\Ga)} && \ti\oX^\Ga
\ar[r]_{\hat\Pi{}^\Ga(\oX)} \ar[ul]_(0.3){\ti O^\Ga(\oX)} &
*+[r]{\hat\oX^\Ga.} }\!\!\!\!
\end{gathered}
\label{ds12eq5}
\e
The underlying\/ $C^\iy$-stacks of\/ $\oX^\Ga,\ldots,
\hat\oX^\Ga_\ci$ are the orbifold strata\/
$\cX^\Ga,\ldots,\hcX^\Ga_\ci$ from\/ {\rm\S\ref{ds87}} of the\/
$C^\iy$-stack\/ $\cX$ in\/ $\oX,$ and the\/ $1$-morphisms in
{\rm\eq{ds12eq5},} as $C^\iy$-stack\/ $1$-morphisms, are those given
in\/~\eq{ds8eq3}.

Use the notation of Definition\/ {\rm\ref{ds9def4}}. Then there are
natural decompositions
\begin{align*}
\oX^\Ga&=\ts\coprod_{\la\in\La^\Ga_+}\oX^{\Ga,\la},&\ti\oX{}^\Ga
&=\ts\coprod_{\mu\in\La^\Ga_+/\Aut(\Ga)}\ti\oX{}^{\Ga,\mu},&
\hat\oX{}^\Ga &=\ts\coprod_{\mu\in\La^\Ga_+/\Aut(\Ga)}\hat\oX{}^{\Ga,\mu},\\
\oX^\Ga_\ci&=\ts\coprod_{\la\in\La^\Ga_+}\oX^{\Ga,\la}_\ci,&
\ti\oX{}^\Ga_\ci &=\ts\coprod_{\mu\in\La^\Ga_+/\Aut(\Ga)}
\ti\oX{}^{\Ga,\mu}_\ci,& \hat\oX{}^\Ga_\ci
&=\ts\coprod_{\mu\in\La^\Ga_+/\Aut(\Ga)}\hat\oX{}^{\Ga,\mu}_\ci,
\end{align*}
where\/
$\oX^{\Ga,\la},\ldots,\hat\oX{}^{\Ga,\mu}_\ci$\G[XGah]{$\oX^{\Ga,\la},
\ti\oX{}^{\Ga,\mu},\hat\oX{}^{\Ga,\mu},
\oX{}^{\Ga,\la}_\ci,\ti\oX{}^{\Ga,\mu}_\ci,\hat\oX{}^{\Ga,\mu}_\ci$}{orbifold
strata of an orbifold with corners $\oX$} are orbifolds with
corners, open and closed in $\oX^\Ga,\ldots,\hat\oX{}^\Ga_\ci,$ and
of dimensions $\dim\oX-\dim\la,\dim\oX-\dim\mu$. All of\/
$\oX^\Ga,\ti\oX{}^\Ga,\hat\oX{}^\Ga,\oX^\Ga_\ci,\ti\oX{}^\Ga_\ci,
\hat\oX{}^\Ga_\ci,\ab\oX^{\Ga,\la},\ab\ti\oX{}^{\Ga,\mu},\ab
\hat\oX{}^{\Ga,\mu},\ab\oX_\ci^{\Ga,\la},\ab\ti\oX{}^{\Ga,\mu}_\ci,
\hat\oX{}^{\Ga,\mu}_\ci$ will be called \begin{bfseries}orbifold
strata\end{bfseries} of\/~$\oX$.
\label{ds12thm4}
\end{thm}

The definitions of $\oX^\Ga,\ti\oX{}^\Ga,\ldots, \hat\oX{}^\Ga_\ci$
also make sense if $\oX$ lies in $\cOrbc$ rather than $\Orbc$. We
will not use notation $\oX^{\Ga,\la},\ldots,\hat\oX{}^{\Ga,\mu}_\ci$
for~$\oX\in\cOrbc\sm\Orbc$.

As for Deligne--Mumford $C^\iy$-stacks in \S\ref{ds87}, orbifold
strata $\oX^\Ga$ are strongly functorial for representable
1-morphisms in $\Orbc$ and their 2-morphisms. That is, if
$f:\oX\ra\oY$ is a representable 1-morphism in $\Orbc$, there is a
unique representable 1-morphism $f^\Ga:\oX^\Ga\ra\oY^\Ga$ in
$\cOrbc$ with $O^\Ga(\oY)\ci f^\Ga=f\ci O^\Ga(\oX)$, which is just
the 1-morphism $f^\Ga$ from \S\ref{ds87} for the $C^\iy$-stack
1-morphism $f:\cX\ra\cY$. Note however that $f^\Ga$ need not map
$\oX^{\Ga,\la}\ra\oY^{\Ga,\la}$ for~$\la\in\La^\Ga_+$.

If $f,g:\oX\ra\oY$ are representable and $\eta:f\Ra g$ is a
2-morphism in $\Orbc$, there is a unique 2-morphism
$\eta^\Ga:f^\Ga\Ra g^\Ga$ in $\cOrbc$ with
$\id_{O^\Ga(\oY)}*\eta^\Ga=\eta*\id_{O^\Ga(\oX)}$, which is just the
$C^\iy$-stack 2-morphism $\eta^\Ga$ from \S\ref{ds87}. These
$f^\Ga,\eta^\Ga$ are compatible with compositions of 1- and
2-morphisms, and identities, in the obvious way. Orbifold strata
$\ti\oX{}^\Ga$ have the same strong functorial behaviour, and
orbifold strata $\hat\oX{}^\Ga$ a weaker functorial behaviour.

We also investigate the relationship between orbifold strata and
corners.

\begin{thm} Let\/ $\oX$ be an orbifold with corners, and\/ $\Ga$ a
finite group. The corners $C(\oX)$ lie in $\cOrbc$ as in
{\rm\S\ref{ds123},} so we have orbifold strata $\oX^\Ga,C(\oX)^\Ga$
and\/ $1$-morphisms $O^\Ga(\oX):\oX^\Ga\ra\oX,$
$O^\Ga(C(\oX)):C(\oX)^\Ga\ra C(\oX)$. Applying the corner functor
$C$ from {\rm\S\ref{ds123}} gives a $1$-morphism
$C(O^\Ga(\oX)):C(\oX^\Ga)\ra C(\oX)$. Then there exists a unique
equivalence $K^\Ga(\oX):C(\oX^\Ga)\ra C(\oX)^\Ga$ such that\/
$O^\Ga(C(\oX))\ci K^\Ga(\oX)=C(O^\Ga(\oX)):C(\oX^\Ga)\ra C(\oX)$. It
restricts to an equivalence~$K^\Ga_\ci(\oX):=K^\Ga(\oX)
\vert_{\smash{C(\oX_\ci^\Ga) }}:C(\oX^\Ga_\ci)\ra C(\oX)^\Ga_\ci$.

Similarly, there is a unique equivalence $\ti K{}^\Ga(\oX):
C(\ti\oX{}^\Ga)\ra\,\,\,\widetilde{\!\!\!C(\oX)\!\!\!}\,\,\,^\Ga$
with\/ $\ti O{}^\Ga(C(\oX))\ci\ti K{}^\Ga(\oX)= C(\ti O{}^\Ga(\oX))$
and\/ $\ti\Pi{}^\Ga(C(\oX))\ci K{}^\Ga(\oX)=\ti K{}^\Ga(\oX)\ci
C(\ti\Pi{}^\Ga(\oX))$. There is an equivalence\/ $\hat
K{}^\Ga(\oX):C(\hat\oX{}^\Ga)\ra \,\,\,\widehat{\!\!\!
C(\oX)\!\!\!}\,\,\,^\Ga,$ unique up to $2$-isomorphism, with a
$2$-morphism $\hat\Pi{}^\Ga(C(\oX))\ci\ti K{}^\Ga(\oX)\Ra\hat
K{}^\Ga(\oX)\ci C(\hat\Pi{}^\Ga(\oX))$. They both restrict to
equivalences $\ti K{}^\Ga_\ci(\oX):C(\ti\oX{}^\Ga_\ci)\ra
\,\,\,\widetilde{\!\!\!C(\oX)\!\!\!}\,\,\,^\Ga_\ci$ and\/~$\hat
K{}^\Ga_\ci(\oX):C(\hat\oX{}^\Ga_\ci)\ra \,\,\,\widehat{\!\!\!
C(\oX)\!\!\!}\,\,\,^\Ga_\ci$.
\label{ds12thm5}
\end{thm}

Here is an example:

\begin{ex} Let $\Z_2=\{1,\si\}$ with $\si^2=1$ act on $X=[0,\iy)^2$ by
$\si:(x_1,x_2)\mapsto(x_2,x_1)$. Then $\oX=\bigl[
[0,\iy)^2/\Z_2\bigr]$ is an orbifold with corners. We have
$\pd\oX\cong[0,\iy)$ and $\pd^2\oX\cong *$, so that
$C_2(\oX)\simeq[*/S_2]=[*/\Z_2]$. Hence $C(\oX)=C_0(\oX)\amalg
C_1(\oX)\amalg C_2(\oX)$ with $C_0(\oX)\simeq\bigl[
[0,\iy)^2/\Z_2\bigr]$, $C_1(\oX)\simeq[0,\iy)$
and~$C_2(\oX)\simeq[*/\Z_2]$. The orbifold strata
$\oX^\Ga,\ldots,\hat\oX{}^\Ga_\ci$ are given by
\begin{equation*}
\oX^{\Z_2}=\oX^{\Z_2}_\ci\simeq\ti\oX{}^{\Z_2}=
\ti\oX^{\Z_2}_\ci\simeq[0,\iy)\t[*/\Z_2],\qquad
\hat\oX{}^{\Z_2}=\hat\oX^{\Z_2}_\ci\simeq[0,\iy).
\end{equation*}
Therefore
\begin{align*}
C_0(\oX^{\Z_2})&\simeq[0,\iy)\t[*/\Z_2], & C_1(\oX^{\Z_2})&\simeq [*/\Z_2],
& C_2(\oX^{\Z_2})&=\es,\\
C_0(\oX)^{\Z_2}&\simeq [0,\iy)\t[*/\Z_2], & C_1(\oX)^{\Z_2}&=\es,
& C_2(\oX)^{\Z_2}&\simeq [*/\Z_2].
\end{align*}
We see from this that $K^{\Z_2}(\oX):C(\oX^{\Z_2})\ra C(\oX)^{\Z_2}$
identifies $C_1(\oX^{\Z_2})$ with $C_2(\oX)^{\Z_2}$, so $K^\Ga(\oX)$
need not map $C_k(\oX^\Ga)$ to $C_k(\oX)^\Ga$ for $k>0$. The same
applies to~$\ti K{}^\Ga(\oX),\hat K{}^\Ga(\oX)$.
\label{ds12ex3}
\end{ex}

The construction of $K^\Ga(\oX)$ in Theorem \ref{ds12thm5} implies
that it maps $C_k(\oX^\Ga)$ into $\coprod_{l\ge k}C_l(\oX)^\Ga$ for
$k>0$. This implies that $C_1(\oX)^\Ga\simeq(\pd\oX)^\Ga$ is
equivalent to an open and closed subobject of $C_1(\oX^\Ga)\simeq
\pd(\oX^\Ga)$. Hence we can choose a 1-morphism
$J^\Ga(\oX):(\pd\oX)^\Ga\ra \pd(\oX^\Ga)$ identified with a
quasi-inverse for $K^\Ga(\oX)\vert_{\cdots}:
K^\Ga(\oX)^{-1}(C_1(\oX)^\Ga)\ra C_1(\oX)^\Ga$ by the equivalences
$C_1(\oX)^\Ga\simeq(\pd\oX)^\Ga$ and $C_1(\oX^\Ga)\simeq
\pd(\oX^\Ga)$, and $J^\Ga(\oX)$ is an equivalence between
$(\pd\oX)^\Ga$ and an open and closed subobject of $\pd(\oX^\Ga)$.
We then deduce:

\begin{cor} Let\/ $\oX$ be an orbifold with corners, and\/ $\Ga$ a
finite group. Then there exist\/ $1$-morphisms
$J^\Ga(\oX):(\pd\oX)^\Ga\ra \pd(\oX^\Ga),$ $\ti J{}^\Ga(\oX):
\,\,\,\,\widetilde{\!\!\!\!(\pd\oX)\!\!\!\!}\,\,\,\,^\Ga\ra
\pd(\ti\oX{}^\Ga),$ $\hat J{}^\Ga(\oX):\,\,\,\,\widehat{\!\!\!\!
(\pd\oX)\!\!\!\!}\,\,\,\,^\Ga\ra \pd(\hat\oX{}^\Ga)$ in\/ $\cOrbc,$
natural up to\/ $2$-isomorphism, such that\/ $J^\Ga(\oX)$ is an
equivalence from $(\pd\oX)^\Ga$ to an open and closed subobject of\/
$\pd(\oX^\Ga),$ and similarly for\/~$\ti J{}^\Ga(\oX),\hat
J{}^\Ga(\oX)$.

For $\la\in\La^\Ga_+,$ $\mu\in\La^\Ga_+/\Aut(\La)$ these restrict to
$1$-morphisms $J^{\Ga,\la}(\oX):(\pd\oX)^{\Ga,\la} \ra
\pd(\oX^{\Ga,\la}),$ $\ti J{}^{\Ga,\mu}(\oX):
\,\,\,\,\widetilde{\!\!\!\!(\pd\oX)\!\!\!\!}\,\,\,\,^{\Ga,\mu}\ra
\pd(\ti\oX{}^{\Ga,\mu}),$ $\hat
J{}^{\Ga,\mu}(\oX):\,\,\,\,\widehat{\!\!\!\!
(\pd\oX)\!\!\!\!}\,\,\,\,^{\Ga,\mu}\ra \pd(\hat\oX{}^{\Ga,\mu})$
in\/ $\Orbc,$ which are equivalences with open and closed
suborbifolds. Hence, if\/ $\oX^{\Ga,\la}=\es$ then
$(\pd\oX)^{\Ga,\la}=\es,$ and similarly for $\ti\oX{}^{\Ga,\mu},
\,\,\,\,\widetilde{\!\!\!\!(\pd\oX)\!\!\!\!}\,\,\,\,^{\Ga,\mu},
\hat\oX{}^{\Ga,\mu},\,\,\,\,\widehat{\!\!\!\!
(\pd\oX)\!\!\!\!}\,\,\,\,^{\Ga,\mu}$.
\label{ds12cor}
\end{cor}

As in Remark \ref{ds12rem}(c), an orbifold with corners $\oX$ is
called {\it straight\/}\I{orbifold with corners!straight} if
$(i_\oX)_*:\Iso_\cpX([x'])\ra \Iso_\cX([x])$ is an isomorphism for
all $[x']\in\cpX_\top$ with $i_{\oX,\top}([x'])=[x]$. If $\oX$ is
straight then $K^\Ga(\oX)$ in Theorem \ref{ds12thm5} is an
equivalence $C_k(\oX^\Ga)\ra C_k(\oX)^\Ga$ for all $k\ge 0$, and so
$J^\Ga(\oX)$ in Corollary \ref{ds12cor} is an equivalence
$(\pd\oX)^\Ga\ra \pd(\oX^\Ga)$. The same applies for~$\ti
J{}^\Ga(\oX),\hat J{}^\Ga(\oX),\ti K{}^\Ga(\oX),\hat K{}^\Ga(\oX)$.

Proposition \ref{ds9prop2} on orientations of orbifold strata
$\cX^{\Ga,\la},\ldots,\hcX^{\Ga,\mu}_\ci$ of oriented orbifolds
$\cX$ also holds without change for orbifolds with
corners~$\oX$.\I{orbifold with corners|)}\I{orbifold with
corners!orbifold strata|)}\I{orbifold strata!of orbifolds with
corners|)}

\section{D-stacks with corners}
\label{ds13}
\I{d-stack with corners|(}

In \cite[Chap.~11]{Joyc6} we define and discuss the 2-category
$\dStac$ of {\it d-stacks with corners}. There are few new issues
here: almost all the material just combines ideas we have seen
already on d-spaces with corners from \S\ref{ds6}, and on d-stacks
from \S\ref{ds10}, and on orbifolds with corners from \S\ref{ds12}.
So we will be brief.

\subsection{Outline of the definition of the 2-category $\dStac$}
\label{ds131}
\I{d-stack with corners!definition|(}\I{2-category|(}

The definition of the 2-category
$\dStac$\G[dStac]{$\dStac$}{2-category of d-stacks with corners} in
\cite[\S 11.1]{Joyc6} is long and complicated. So as for $\dSpac$ in
\S\ref{ds61}, we will just sketch the main ideas.

A {\it d-stack with corners\/} is a quadruple $\eX=(\bcX,\bcpX,\bs
i_\eX,\om_\eX)$, where $\bcX,\bcpX$ are d-stacks and $\bs
i_\eX:\bcpX\ra\bcX$\G[iXd]{$\bs i_\eX:\pd\eX\ra\eX$}{inclusion of
boundary $\pd\eX$ into a d-stack with corners $\eX$} is a 1-morphism
of d-stacks with $i_\eX:\cpX\ra\cX$ a proper, strongly representable
1-morphism of Deligne--Mumford $C^\iy$-stacks, as in \S\ref{ds83}.
We should have an exact sequence in $\qcoh(\cpX)$:
\e
\xymatrix{ 0 \ar[r] & \cN_\eX \ar[rr]^{\nu_\eX} && i_\eX^*(\FcX)
\ar[rr]^{i_\eX^2} && \cF_{\cpX} \ar[r] & 0,}
\label{ds13eq1}
\e
where $\cN_\eX$ is a line bundle on $\cpX$, the {\it conormal
bundle\/} of $\bcpX$ in $\bcX$,\G[NXb]{$\cN_\eX$}{conormal line
bundle of $\pd\eX$ in $\eX$ for a d-stack with corners $\eX$} and
$\om_\eX$\G[omXb]{$\om_\eX$}{orientation on line bundle $\cN_\eX$
for a d-stack with corners $\eX$} is an orientation on $\cN_\eX$.
These $\bcX,\bcpX,\bs i_\eX,\om_\eX$ must satisfy some complicated
conditions in \cite[\S 11.1]{Joyc6}, that we will not give. They
require $\bcpX$ to be locally equivalent to a fibre product
$\bcX\t_{\bs{[0,\iy)}}\bs{*}$
in~$\dSta$.\G[WXYZe]{$\eW,\eX,\eY,\eZ,\ldots$}{d-stacks with
corners, including d-orbifolds with corners}

If $\eX=(\bcX,\bcpX,\bs i_\eX,\om_\eX)$ and $\eY=(\bcY,\bcpY,\bs
i_\eY,\om_\eY)$ are d-stacks with corners, a 1-{\it morphism\/} $\bs
f:\eX\ra\eY$ in $\dStac$ is a 1-morphism $\bs f:\bcX\ra\bcY$ in
$\dSta$ satisfying extra conditions over $\bcpX,\bcpY$. If $\bs
f,\bs g:\eX\ra\eY$ be 1-morphisms in $\dStac$, so $\bs f,\bs
g:\bcX\ra\bcY$ are 1-morphisms in $\dSta$, a 2-morphism $\bs\eta:\bs
f\Ra\bs g$ in $\dStac$ is a 2-morphism $\bs\eta:\bs f\Ra\bs g$ in
$\dSta$ satisfying extra conditions over $\bcpX,\bcpY$. In both
cases, 1- and 2-morphisms in $\dStac$ are \'etale locally modelled
on 1- and 2-morphisms in $\dSpac$. Identity 1- and 2-morphisms in
$\dStac$, and the compositions of 1- and 2-morphisms in $\dStac$,
are all given by identities and compositions
in~$\dSta$.\I{2-category|)}

A d-stack with corners $\eX=(\bcX,\bcpX,\bs i_\eX,\om_\eX)$ is
called a {\it d-stack with boundary\/}\I{d-stack with boundary} if
$i_\eX:\cpX\ra\cX$ is injective as a representable 1-morphism of
$C^\iy$-stacks, and a {\it d-stack without boundary\/} if
$\bcpX=\bs\es$. We write $\dStab$\G[dStab]{$\dStab$}{2-category of
d-stacks with boundary} for the full 2-subcategory of d-stacks with
boundary, and $\bdSta$\G[dSta']{$\bdSta$}{2-subcategory of d-stacks
with corners equivalent to d-stacks} for the full 2-subcategory of
d-stacks without boundary, in $\dStac$. There is an isomorphism of
2-categories $F_\dSta^\dStac:\dSta\ra\bdSta$ mapping
$\bcX\mapsto\eX=(\bcX,\bs\es,\bs\es,\bs\es)$ on objects, $\bs
f\mapsto\bs f$ on 1-morphisms and $\bs\eta\mapsto\bs\eta$ on
2-morphisms. So we can consider d-stacks to be examples of d-stacks
with corners.

Define a strict 2-functor\I{2-category!strict 2-functor}
$F_\dSpac^\dStac:\dSpac\ra\dStac$ as follows. If
$\rX=(\bX,\ab\bpX,\ab\bs i_\rX,\ab\om_\rX)$ is an object in
$\dSpac$, set $F_\dSpac^\dStac(\rX)=\eX=(\bcX,\bcpX,\bs
i_\eX,\om_\eX)$, where $\bcX,\bcpX,\bs i_\eX=F_\dSpa^\dSta
(\bX,\bpX,\bs i_\rX)$. Then comparing equations \eq{ds6eq2} and
\eq{ds13eq1}, we find there is a natural isomorphism of line bundles
$\cN_\eX\cong\cI_\upX(\cN_\rX)$, where ${\cal I}_\upX:\qcoh(\upX)
\ra\qcoh(\bar\upX)$ is the equivalence of categories from Example
\ref{ds8ex6}. We define $\om_\eX$ to be the orientation on $\cN_\eX$
identified with the orientation $\cI_\upX(\om_\rX)$ on
$\cI_\upX(\cN_\rX)$ by this isomorphism. On 1- and 2-morphisms $\bs
f,\eta$ in $\dSpac$, we define $F_\dSpac^\dStac(\bs
f)=F_\dSpa^\dSta(\bs f)$ and~$F_\dSpac^\dStac(\eta)=
F_\dSpa^\dSta(\eta)$.

Write $\hdSpac$\G[dSpac']{$\hdSpac$}{2-subcategory of d-stacks with
corners equivalent to d-spaces with corners} for the full
2-subcategory of objects $\eX$ in $\dStac$ equivalent to
$F_\dSpac^\dStac(\rX)$ for some d-space with corners $\rX$. When we
say that a d-stack with corners $\eX$ {\it is a d-space},\I{d-stack
with corners!is a d-space} we mean that~$\eX\in\hdSpac$.

Define a strict 2-functor $F_\Orbc^\dStac:\Orbc\ra\dStac$ as
follows. If $\oX=(\cX,\cpX,i_\oX)$ is an orbifold with corners, as
in \S\ref{ds121}, define $F_\Orbc^\dStac(\oX)=\eX=(\bcX,\bcpX,\bs
i_\eX,\om_\eX)$, where $\bcX,\bcpX,\bs i_\eX=F_\CSta^\dSta
(\cX,\cpX,i_\oX)$. Then $\cN_\eX$ in \eq{ds13eq1} is isomorphic to
the conormal line bundle of $\cpX$ in $\cX$, and we define $\om_\eX$
to be the orientation on $\cN_\eX$ induced by `outward-pointing'
normal vectors to $\cpX$ in $\cX$. Then $\eX=(\bcX,\bcpX,\bs
i_\eX,\om_\eX)$ is a d-orbifold with corners. On 1- and 2-morphisms
$f,\eta$ in $\Orbc$, we define $F_\Orbc^\dStac(f)=F_\CSta^\dSta(f)$
and~$F_\Orbc^\dStac(\eta)=F_\CSta^\dSta(\eta)$.

Write $\bOrb,\bOrbb,\bOrbc$\G[Orbc']{$\bOrbc$}{2-subcategory of
d-stacks with corners equivalent to orbifolds with corners} for the
full 2-subcategories of objects $\eX$ in $\dStac$ equivalent to
$F_\Orbc^\dStac(\oX)$ for some orbifold $\oX$ without boundary, or
with boundary, or with corners, respectively. Then
$\bOrb\subset\bdSta$, $\bOrbb\subset\dStab$ and
$\bOrbc\subset\dStac$. When we say that a d-stack with corners $\eX$
{\it is an orbifold},\I{d-stack with corners!is an orbifold} we mean
that~$\eX\in\bOrbc$.

\begin{rem} As discussed for orbifolds with corners in Remark
\ref{ds12rem}(b), in a d-stack with corners $\eX=(\bcX,\bcpX,\bs
i_\eX,\om_\eX)$ we require $i_\eX:\cpX\ra\cX$ to be {\it strongly
representable},\I{C-stack@$C^\iy$-stack!strongly representable
1-morphism} in the sense of \S\ref{ds83}, so that we can make
boundaries and corners in $\dStac$ strictly functorial,\I{d-stack
with corners!boundary!strictly functorial} as in Remark
\ref{ds6rem2} for~$\dSpac$.
\label{ds13rem}
\end{rem}

For each d-stack with corners $\eX=(\bcX,\bcpX,\bs i_\eX,\om_\eX)$,
in \cite[\S 11.3]{Joyc6} we define a d-stack with corners
$\pd\eX=(\bcpX, \bcptX,\bs
i_{\pd\eX},\om_{\pd\eX})$\G[dXd]{$\pd\eX$}{boundary of a d-stack
with corners $\eX$} called the {\it boundary\/}\I{d-stack with
corners!boundary}\I{boundary!of a d-stack with corners} of $\eX$,
and show that $\bs i_\eX:\pd\eX\ra\eX$ is a 1-morphism in $\dStac$.
As for d-spaces with corners in \eq{ds6eq3}, the d-stack $\bcptX$ in
$\pd\eX$ satisfies
\begin{equation*}
\bcptX\simeq\bigl(\bcpX\t_{\bs i_\eX,\bcX,\bs i_\eX}\bcpX\bigr)
\sm\bs\De_\bcpX(\bcpX),
\end{equation*}
where $\De_\bcpX:\bcpX\ra\bcpX\t_\bcX\bcpX$ is the diagonal
1-morphism. The 1-morphism $\bs i_{\pd\eX}:\bcptX\ra\bcpX$ is
projection to the first factor in the fibre product. There is a
natural isomorphism $\cN_{\pd\eX}\cong i_\eX^*(\cN_\eX)$, and the
orientation $\om_{\pd\eX}$ on $\cN_{\pd\eX}$ corresponds to the
orientation $i_\eX^*(\om_\eX)$ on~$i_\eX^*(\cN_\eX)$.\I{d-stack with
corners!definition|)}

\subsection[D-stacks with corners as quotients of d-spaces with
corners]{D-stacks with corners as quotients of d-spaces}
\label{ds132}
\I{d-stack with corners!quotients $[\rX/G]$|(}

Section \ref{ds102} discussed quotient d-stacks $[\bX/G]$, for $\bX$
a d-space and $\bs r:G\ra\Aut(\bX)$ an action of $G$ on $\bX$ by
1-isomorphisms. In \cite[\S 11.2]{Joyc6} we extend this to d-spaces
with corners and d-stacks with corners, and prove:

\begin{thm} Theorems\/ {\rm\ref{ds10thm2}} and\/
{\rm\ref{ds10thm3}} hold unchanged in\/~$\dStac$.
\label{ds13thm1}
\end{thm}

Here if $\rX=(\bX,\bpX,\bs i_\rX,\om_\rX)$ is a d-space with corners
and $\bs r:G\ra\Aut(\rX)$ an action of $G$ on $\rX$ then each $\bs
r(\ga):\rX\ra\rX$ for $\ga\in G$ is simple, so Theorem
\ref{ds6thm1}(b) gives a lift $\bs r_-(\ga):\pd\rX\ra\pd\rX$,
defining an action $\bs r_-:G\ra\Aut(\pd\rX)$ of $G$ on $\pd\rX$.
Then $\bs r:G\ra\Aut(\bX)$ and $\bs r_-:G\ra\Aut(\bpX)$ are actions
of $G$ on the d-spaces $\bX,\bpX$, and $\bs i_\rX:\bpX\ra\bX$ is
$G$-equivariant. So Theorem \ref{ds10thm2}(a),(b) give quotient
d-stacks $[\bX/G],[\bpX/G]$ and a quotient 1-morphism $[\bs
i_\rX,\id_G]: [\bpX/G]\ra[\bX/G]$. The quotient d-stack with corners
$[\rX/G]$ given by the analogue of Theorem \ref{ds10thm2} is defined
to be $[\rX/G]=\bigl([\bX/G],[\bpX/G],[\bs
i_\rX,\id_G],\om_{[\rX/G]}\bigr)$, for a natural orientation
$\om_{[\rX/G]}$ on $\cN_{[\rX/G]}$ constructed from~$\om_\rX$.

In \cite[\S 11.4]{Joyc6} we define when a 1-morphism of d-stacks
with corners $\bs f:\eX\ra\eY$ is {\it \'etale}.\I{d-stack with
corners!etale 1-morphism@\'etale 1-morphism} Essentially, $\bs f$ is
\'etale if it is an equivalence locally in the \'etale topology. It
implies that the $C^\iy$-stack 1-morphism $f:\cX\ra\cY$ in $\bs f$
is \'etale, and so representable. As for d-stacks in \S\ref{ds102},
we can characterize \'etale 1-morphisms in $\dStac$ using the
corners analogue of Theorem \ref{ds10thm3}(b) and the definition of
\'etale 1-morphisms in $\dSpac$ as (Zariski) local
equivalences.\I{d-stack with corners!quotients $[\rX/G]$|)}

\subsection{Simple, semisimple and flat 1-morphisms}
\label{ds133}

In \cite[\S 11.3]{Joyc6} we generalize \S\ref{ds62} to d-stacks with
corners. Here is the analogue of Definition~\ref{ds12def4}.

\begin{dfn} Let $\eX=(\bcX,\bcpX,\bs i_\eX,\om_\eX)$ and
$\eY=(\bcY,\bcpY,\bs i_\eY,\om_\eY)$ be d-stacks with corners, and
$\bs f:\eX\ra\eY$ a 1-morphism in $\dStac$. Consider the
$C^\iy$-stack fibre products $\cpX\t_{f\ci i_\eX,\cY, i_\eY}\cpY$
and $\cX\t_{f,\cY, i_\eY}\cpY$. Since $i_\eY$ is strongly
representable, we may define these using the construction of
Proposition~\ref{ds8prop3}.

As in \eq{ds12eq1}, we may write $(\cpX\t_\cY\cpY)_\top$ explicitly
as
\e
\begin{split}
(\cpX\t_\cY\cpY)_\top\cong\bigl\{&[x',y']:\text{$x':\bar{\ul
*}\ra\cpX$ and $y':\bar{\ul *}\ra\cpY$ are}\\
&\text{1-morphisms with $f\!\ci\! i_\eX\!\ci\! x'\!=\! i_\eY\!\ci
\!y':\bar{\ul *}\!\ra\!\cY$}\bigr\},
\end{split}
\label{ds13eq2}
\e
where $[x',y']$ in \eq{ds13eq2} denotes the $\sim$-equivalence class
of pairs $(x',y')$, with $(x',y')\sim(\ti x',\ti y')$ if there exist
2-morphisms $\eta:x'\Ra\ti x'$ and $\ze:y'\Ra\ti y'$ with $\id_{f\ci
i_\oX}*\eta=\id_{i_\eY}*\ze$. There is a natural open and closed
$C^\iy$-substack $\cS_{\bs
f}\subseteq\cpX\t_\cY\cpY$\G[Sfd]{$\cS_{\bs
f}\subseteq\cpX\t_\cY\cpY$}{$C^\iy$-stack associated to 1-morphism
$\bs f:\eX\ra\eY$ in $\dStac$} such that $[x',y']$ in \eq{ds13eq2}
lies in $\cS_{\bs f,\top}$ if and only if we can complete the
following commutative diagram in $\qcoh(\ul{\bar
*})$ with morphisms~`$\dashra$':
\begin{equation*}
\xymatrix@C=13pt@R=25pt{ 0 \ar[r] & (y')^*(\cN_\eY)
\ar[rr]^(0.45){(y')^*(\nu_\eY)} \ar@{.>}[d]_\cong && (y')^*\!\ci\!
i_\eY^*(\FcY) \ar[rr]^{(y')^*(i_\eY^2)}
\ar@<-7ex>[d]^{\begin{subarray}{l}
I_{x',i_\eX}(\FcX)\ci (i_\eX\ci x')^*(f^2) \ci \\
I_{i_\eX\ci x',f}(\FcY)\ci I_{y',i_\eY}(\FcY)^{-1}\end{subarray}} &&
(y')^*(\cF_\cpY) \ar@{.>}[d] \ar[r] & 0 \\
0 \ar[r] & (x')^*(\cN_\eX) \ar[rr]_(0.45){(x')^*(\nu_\eX)} &&
(x')^*\!\ci\! i_\eX^*(\FcX) \ar[rr]_{(x')^*(i_\eX^2)} &&
(x')^*(\cF_\cpX) \ar[r] & 0. }
\end{equation*}

Similarly, as in \eq{ds12eq2} we may write $(\cX\t_\cY\cpY)_\top$
explicitly as
\e
\begin{split}
(\cX\t_\cY\cpY)_\top\cong\bigl\{&[x,y']:\text{$x:\bar{\ul
*}\ra\cX$ and $y':\bar{\ul *}\ra\cpY$ are}\\
&\text{1-morphisms with $f\!\ci\! x\!=\! i_\eY\!\ci \!y':\bar{\ul
*}\!\ra\!\cY$}\bigr\},
\end{split}
\label{ds13eq3}
\e
where $[x,y']$ in \eq{ds13eq3} denotes the $\approx$-equivalence
class of $(x,y')$, with $(x,y')\approx(\ti x,\ti y')$ if there exist
$\eta:x\Ra\ti x$ and $\ze:y'\Ra\ti y'$ with
$\id_f*\eta=\id_{i_\eY}*\ze$. There is a natural open and closed
$C^\iy$-substack $\cT_{\bs
f}\subseteq\cX\t_\cY\cpY$\G[Tfd]{$\cT_{\bs
f}\subseteq\cX\t_\cY\cpY$}{$C^\iy$-stack associated to 1-morphism
$\bs f:\eX\ra\eY$ in $\dStac$} with $[x,y']$ in \eq{ds13eq3} lies in
$\cT_{\bs f,\top}$ if and only if we can complete the following
commutative diagram:
\begin{equation*}
\xymatrix@C=13pt@R=20pt{ 0 \ar[r] & (y')^*(\cN_\eY)
\ar[rr]_(0.45){(y')^*(\nu_\eY)} && (y')^*\!\ci\! i_\eY^*(\FcY)
\ar[rr]_{(y')^*(i_\eY^2)} \ar@<2ex>[d]_(0.6){\begin{subarray}{l}
x^*(f^2) \ci I_{x,f}(\FcY)\ci I_{y',i_\eY}(\FcY)^{-1}\end{subarray}}
&& (y')^*(\cF_\cpY) \ar@<1ex>@{.>}[dll] \ar[r] & 0 \\
&&& x^*(\FcX). }
\end{equation*}

Define $s_{\bs f}=\pi_\cpX\vert_{\cS_{\bs f}}:\cS_{\bs f}\ra\cpX$,
$u_{\bs f}=\pi_\cpY \vert_{\cS_{\bs f}}:\cS_{\bs f}\ra\cpY$, $t_{\bs
f}=\pi_\cX\vert_{\cT_{\bs f}}:\cT_{\bs f}\ra\cX$, and $v_{\bs
f}=\pi_\cpY\vert_{\cT_{\bs f}}:\cT_{\bs f}\ra\cpY$. Then $s_{\bs
f},t_{\bs f}$ are proper, \'etale 1-morphisms. We call $\bs f$ {\it
simple\/}\I{d-stack with corners!simple 1-morphism|(} if $s_{\bs
f}:\cS_{\bs f}\ra\cpX$ is an equivalence, and we call $\bs f$ {\it
semisimple\/}\I{d-stack with corners!semisimple 1-morphism|(} if
$s_{\bs f}:\cS_{\bs f}\ra\cpX$ is injective, as a 1-morphism of
Deligne--Mumford $C^\iy$-stacks, and we call $\bs f$ {\it
flat\/}\I{d-stack with corners!flat 1-morphism|(} if $\cT_{\bs
f}=\es$. Simple implies semisimple.
\label{ds13def1}
\end{dfn}

\begin{thm} Let\/ $\bs f:\eX\ra\eY$ be a semisimple $1$-morphism
of d-stacks with corners. Then there exists a natural decomposition
$\pd\eX=\pd_+^{\bs f}\eX\amalg\pd_-^{\bs f}\eX$ with\/ $\pd_\pm^{\bs
f}\eX$ open and closed in $\pd\eX,$ such that:\G[dfXd]{$\pd^{\bs
f}_\pm\eX$}{sets of decomposition $\pd\eX=\pd^{\bs
f}_+\eX\amalg\pd^{\bs f}_-\eX$ of boundary $\pd\eX$ induced by
1-morphism $\bs f:\eX\ra\eY$ in $\dStac$}
\begin{itemize}
\setlength{\itemsep}{0pt}
\setlength{\parsep}{0pt}
\item[{\bf(a)}] Define $\bs f_+=\bs f\ci \bs i_\eX\vert_{\pd^{\bs
f}_+\eX}:\pd_+^{\bs f}\eX\ra\eY$. Then $\bs f_+$ is semisimple.
If\/ $\bs f$ is flat then $\bs f_+$ is also flat.
\item[{\bf(b)}] There exists a unique, semisimple\/ $1$-morphism
$\bs f_-:\pd_-^{\bs f}\eX\ra\pd\eY$ with\/ $\bs f\ci\bs i_\eX
\vert_{\pd_-^{\bs f}\eX}=\bs i_\eY\ci\bs f_-$. If\/ $\bs f$ is
simple then $\pd_+^{\bs f}\eX=\bs\es,$ $\pd_-^{\bs
f}\eX=\pd\eX,$ and\/ $\bs f_-:\pd\eX\ra\pd\eY$ is also simple.
If\/ $\bs f$ is flat then $\bs f_-$ is flat, and the following
diagram is $2$-Cartesian\I{2-category!2-Cartesian square}
in\/~$\dStac\!:$
\begin{align*}
\xymatrix@C=120pt@R=10pt{ \pd_-^{\bs f}\eX \ar[r]_(0.2){\bs f_-}
\ar[d]_{\bs i_\eX \vert_{\pd_-^{\bs f}\eX}}
\drtwocell_{}\omit^{}\omit{^{\bs\id_{\bs i_\eY\ci \bs
f_-}\,\,\,\,\,\,\,\,\,\,\,\,\,\,{}}} & \pd\eY \ar[d]^{\bs i_\eY}
\\ \eX \ar[r]^(0.7){\bs f} & \eY.}
\end{align*}
\item[{\bf(c)}] Let\/ $\bs g:\eX\ra\eY$ be another $1$-morphism
and\/ $\bs\eta:\bs f\Ra\bs g$ a $2$-morphism in $\dStac$. Then
$\bs g$ is also semisimple, with\/ $\pd_-^{\bs g}\eX=\pd_-^{\bs
f}\eX$. If\/ $\bs f$ is simple, or flat, then $\bs g$ is simple,
or flat, respectively. Part\/ {\bf(b)} defines $1$-morphisms
$\bs f_-,\bs g_-:\pd_-^{\bs f}\eX\ra\pd\eY$. There is a unique
$2$-morphism $\bs\eta_-:\bs f_-\Ra\bs g_-$ in\/ $\dSpac$ such
that\/\I{d-stack with corners!flat 1-morphism|)}\I{d-stack with
corners!semisimple 1-morphism|)}\I{d-stack with corners!simple
1-morphism|)} $\bs\id_{\bs
i_\eY}*\bs\eta_-\!=\!\bs\eta*\bs\id_{\bs i_\eX\vert_{\pd_-^{\bs
f}\eX}}:\bs i_\eY\!\ci\!\bs f_- \Ra\bs i_\eY\!\ci\!\bs g_-$.
\end{itemize}
\label{ds13thm2}
\end{thm}

\subsection{Equivalences in $\dStac$, and gluing by equivalences}
\label{ds134}
\I{d-stack with corners!equivalence|(}\I{d-stack with corners!gluing
by equivalences|(}

Sections \ref{ds32}, \ref{ds64} and \ref{ds103} discussed
equivalences and gluing for d-spaces, d-spaces with corners, and
d-stacks. In \cite[\S 11.4]{Joyc6} we generalize these to~$\dStac$.

\begin{prop}{\bf(a)} Suppose\/ $\bs f:\eX\ra\eY$ is an equivalence
in $\dStac$. Then $\bs f$ is simple\I{d-stack with corners!simple
1-morphism} and flat,\I{d-stack with corners!flat 1-morphism} and\/
$\bs f:\bcX\ra\bcY$ is an equivalence in\/ $\dSta,$ where
$\eX=(\bcX,\bcpX,\bs i_\eX,\om_\eX)$ and\/ $\eY=(\bcY,\bcpY,\bs
i_\eY,\om_\eY)$. Also\/ $\bs f_-:\pd\eX\ra\pd\eY$ in Theorem\/
{\rm\ref{ds13thm2}(b)} is an equivalence in\/~$\dStac$.
\smallskip

\noindent{\bf(b)} Let\/ $\bs f:\eX\ra\eY$ be a simple, flat\/
$1$-morphism in $\dStac$ with\/ $\bs f:\bcX\ra\bcY$ an equivalence
in\/ $\dSta$. Then $\bs f$ is an equivalence in\/~$\dStac$.
\label{ds13prop1}
\end{prop}

Here is the analogue of Definition~\ref{ds10def2}:

\begin{dfn} Let $\eX=(\bcX,\bcpX,\bs i_\eX,\om_\eX)$ be a d-stack
with corners. Suppose $\bcV\subseteq\bcX$ is an open d-substack in
$\dSta$. Define $\bs{\pd\cV}=\bs i_\eX^{-1}(\bcV)$, as an open
d-substack of $\bcpX$, and $\bs i_\eV:\bs{\pd\cV}\ra\bcV$ by $\bs
i_\eV=\bs i_\eX\vert_{\bs{\pd\cV}}$. Then $\cpV\subseteq\cpX$ is
open, and the conormal bundle of $\bs{\pd\cV}$ in $\bcV$ is
$\cN_\eV=\cN_\eX\vert_\cpV$ in $\qcoh(\cpV)$. Define an orientation
$\om_\eV$ on $\cN_\eV$ by $\om_\eV=\om_\eX\vert_\cpV$. Write
$\eV=(\bcV,\bs{\pd\cV},\bs i_\eV,\om_\eV)$. Then $\eV$ is a d-stack
with corners. We call $\eV$ an {\it open d-substack\/} of $\eX$. An
{\it open cover\/} of $\eX$ is a family $\{\eV_a:a\in A\}$ of open
d-substacks $\eV_a$ of $\eX$ with~$\cX=\bigcup_{a\in
A}\cV_a$.\I{d-stack with corners!open d-substack}\I{d-stack with
corners!open cover}
\label{ds13def2}
\end{dfn}

\begin{thm} Proposition\/ {\rm\ref{ds10prop}} and Theorems\/
{\rm\ref{ds10thm4}} and\/ {\rm\ref{ds10thm5}} hold without change in
the $2$-category\/ $\dStac$ of d-stacks with corners.\I{d-stack with
corners!equivalence|)}\I{d-stack with corners!gluing by
equivalences|)}
\label{ds13thm3}
\end{thm}

\subsection{Corners $C_k(\eX)$, and the corner functors $C,\hat C$}
\label{ds135}
\I{d-stack with corners!k-corners C_k(X)@$k$-corners
$C_k(\eX)$|(}\I{d-stack with corners!corner functors|(}

In \cite[\S 11.5]{Joyc6} we generalize the material of \S\ref{ds53},
\S\ref{ds65}, and \S\ref{ds123} to d-stacks with corners. Here are
the main results.

\begin{thm}{\bf(a)} Let\/ $\eX$ be a d-stack with corners. Then
for each\/ $k\ge 0$ we can define a d-stack with corners\/
$C_k(\eX)$ called the $k$-\begin{bfseries}corners\end{bfseries} of\/
$\eX,$ and a $1$-morphism $\bs\Pi^k_\eX:C_k(\eX)\ra\eX,$ such that\/
$C_k(\eX)$ is equivalent to a quotient d-stack\/ $[\pd^k\eX/S_k]$
for a natural action of\/ $S_k$ on $\pd^k\eX$ by $1$-isomorphisms.
The construction of\/ $C_k(\eX)$ is unique up to canonical\/
$1$-isomorphism.

We can describe the topological space $C_k(\cX)_\top$ as follows.
Consider pairs $(x,\{x_1',\ldots,x_k'\}),$ where $x:\bar{\ul{*}}\ra
\cX$ and\/ $x_i':\bar{\ul{*}}\ra\cpX$ for $i=1,\ldots,k$ are
$1$-morphisms in $\CSta$ with\/ $x_1',\ldots,x_k'$ distinct and\/
$x=i_\eX\ci x_1'=\cdots=i_\eX\ci x_k'$. Define an equivalence
relation $\approx$ on such pairs by
$(x,\{x_1',\ldots,x_k'\})\approx(\ti x,\{\ti x_1',\ldots,\ti
x_k'\})$ if there exist\/ $\si\in S_k$ and\/ $2$-morphisms
$\eta:x\Ra\ti x$ and\/ $\eta_i':x_i'\Ra\ti x_{\si(i)}'$ for
$i=1,\ldots,k$ with\/ $\eta=\id_{i_\eX}*\eta_1'=
\cdots=\id_{i_\eX}*\eta_k'$. Write $[x,\{x_1',\ldots,x_k'\}]$ for
the $\approx$-equivalence class of\/ $(x,\{x_1',\ldots,x_k'\})$.
Then
\e
\begin{split}
C_k(\cX)_\top\!\cong\!\bigl\{&[x,\{x_1',\ldots,x_k'\}]:\text{$x:\bar{\ul{*}}
\!\ra\!\cX,$ $x_i':\bar{\ul{*}}\!\ra\!\cpX$ $1$-morphisms}\\
&\text{with\/ $x_1',\ldots,x_k'$ distinct and\/ $x\!=\!i_\eX\!\ci\!
x_1'\!=\!\cdots\!=\!i_\eX\!\ci\!x_k'$}\bigr\}.
\end{split}
\label{ds13eq4}
\e
We have $1$-isomorphisms $C_0(\eX)\cong\eX$ and\/
$C_1(\eX)\cong\pd\eX$. We write $C(\eX)=\coprod_{k\ge 0}C_k(\eX),$
so that\/ $C(\eX)$ is a d-stack with corners, called the
\begin{bfseries}corners\end{bfseries} of\/~$\eX$.

\smallskip

\noindent{\bf(b)} Let\/ $\bs f:\eX\ra\eY$ be a $1$-morphism of
d-stacks with corners. Then there are unique $1$-morphisms $C(\bs
f):C(\eX)\ra C(\eY)$ and\/ $\hat C(\bs f):C(\eX)\ra C(\eY)$ in
$\dStac$ such that\/ $\bs\Pi_\eY\ci C(\bs f)=\bs f\ci\bs\Pi_\eX
=\bs\Pi_\eY\ci\hat C(\bs f):C(\eX)\ra\eY,$ with maps
$C(f)_\top:C(\cX)_\top\ra C(\cY)_\top,$ $\hat
C(f)_\top:C(\cX)_\top\ra C(\cY)_\top$ characterized as follows.
Identify $C_k(\cX)_\top\subseteq C(\cX)_\top$ with the right hand
side of\/ {\rm\eq{ds13eq4},} and similarly for\/ $C_l(\cY)_\top,$
and identify\/ $\cS_{\bs f,\top},\cT_{\bs f,\top}$ with subsets of\/
\eq{ds13eq2}--\eq{ds13eq3} as in\/ {\rm\S\ref{ds133}}. Then
$C(f)_\top$ and\/ $\hat C(f)_\top$ act by
\ea
\begin{split}
&C(f)_\top:\bigl[x,\{x_1',\ldots,x_k'\}\bigr]\longmapsto
\bigl[y,\{y_1',\ldots,y_l'\}\bigr],\;\>\text{where $y=f\ci x$,}\\
&\{y_1',\ldots,y_l'\}\!=\!\bigl\{y':[x_i',y']\in \cS_{\bs f,\top},\;
\text{some\/ $i=1,\ldots,k\bigr\},$ and}
\end{split}
\label{ds13eq5}\\
\begin{split}
 &\hat C(f)_\top:\bigl[x,\{x_1',\ldots,x_k'\}\bigr]\longmapsto
\bigl[y,\{y_1',\ldots,y_l'\}\bigr],\;\>\text{where $y=f\ci x$,}\\
&\{y_1',\ldots,y_l'\}\!=\!\bigl\{y':[x_i',y']\!\in\!\cS_{\bs
f,\top},\; i\!=\!1,\ldots,k\bigr\}\!\cup\!\bigl\{y':[x,y']\!\in\!
\cT_{\bs f,\top}\bigr\}.
\end{split}
\label{ds13eq6}
\ea

For all\/ $k,l\ge 0,$ write\/ $C_k^{\bs f,l}(\eX)=C_k(\eX)\cap C(\bs
f)^{-1}(C_l(\eY)),$ so that\/ $C_k^{\bs f,l}(\eX)$ is an open and
closed d-substack of\/ $C_k(\eX)$ with\/ $C_k(\eX)=\coprod_{l=0}^\iy
C_k^{\bs f,l}(\eX),$ and write $C^l_k(\bs f)=C(\bs f)\vert_{C_k^{\bs
f,l}(\eX)}:C_k^{\bs f,l}(\eX)\ra C_l(\eY)$. If\/ $\bs f$ is
simple\I{d-stack with corners!simple 1-morphism} then $C(\bs f)$
maps $C_k(\eX)\ra C_k(\eY)$ for all\/ $k\ge 0$. If\/ $\bs f$ is
flat\I{d-stack with corners!flat 1-morphism} then~$C(\bs f)=\hat
C(\bs f)$.
\smallskip

\noindent{\bf(c)} Let\/ $\bs f,\bs g:\eX\ra\eY$ be $1$-morphisms
and\/ $\bs\eta:\bs f\Ra\bs g$ a $2$-morphism in $\dStac$. Then there
exist unique $2$-morphisms $C(\bs\eta):C(\bs f)\Ra C(\bs g),$ $\hat
C(\bs\eta):\hat C(\bs f)\Ra\hat C(\bs g)$ in $\dStac,$ where $C(\bs
f),C(\bs g),\hat C(\bs f),\hat C(\bs g)$ are as in {\bf(b)\rm,} such
that
\begin{align*}
&\id_{\bs\Pi_\eY}*
C(\bs\eta)\!=\!\bs\eta*\id_{\bs\Pi_\eX}:\bs\Pi_\eY\!\ci\! C(\bs
f)\!=\!\bs f\!\ci\!\bs\Pi_\eX\Longra \bs\Pi_\eY\!\ci\!C(\bs
g)\!=\!\bs g\!\ci\!\bs\Pi_\eX,\\
&\id_{\bs\Pi_\eY}*
\hat C(\bs\eta)\!=\!\bs\eta*\id_{\bs\Pi_\eX}:\bs\Pi_\eY\!\ci\!\hat C(\bs
f)\!=\!\bs f\!\ci\!\bs\Pi_\eX\Longra \bs\Pi_\eY\!\ci\!\hat C(\bs
g)\!=\!\bs g\!\ci\!\bs\Pi_\eX.
\end{align*}
If\/ $\bs f,\bs g$ are flat\I{d-stack with corners!flat 1-morphism}
then $C(\bs\eta)=\hat C(\bs\eta)$.
\smallskip

\noindent{\bf(d)} Define\/ $C:\dStac\ra\dStac$ by\/ $C:\eX\mapsto
C(\eX),$ $C:\bs f\mapsto C(\bs f),$ $C:\bs\eta\mapsto C(\bs\eta)$ on
objects, $1$- and\/ $2$-morphisms, where $C(\eX),C(\bs
f),C(\bs\eta)$ are as in {\bf(a){\rm --}(c)} above. Similarly,
define\/ $\hat C:\dStac\ra\dStac$ by $\hat C:\eX\mapsto C(\eX),$
$\hat C:\bs f\mapsto\hat C(\bs f),$ $\hat C:\bs\eta\mapsto\hat
C(\bs\eta)$. Then $C,\hat C$ are strict\/ $2$-functors, called
\begin{bfseries}corner functors\end{bfseries}.\G[CCd]{$C,\hat
C:\dStac\ra\dStac$}{`corner functors' for d-stacks with
corners}\I{d-stack with corners!k-corners C_k(X)@$k$-corners
$C_k(\eX)$|)}\I{d-stack with corners!corner functors|)}
\label{ds13thm4}
\end{thm}

\subsection{Fibre products in $\dStac$}
\label{ds136}
\I{d-stack with corners!fibre
products|(}\I{b-transversality|(}\I{c-transversality|(}\I{d-stack
with corners!b-transverse 1-morphisms|(}\I{d-stack with
corners!c-transverse 1-morphisms|(}\I{2-category!fibre products
in|(}

In \cite[\S 11.6]{Joyc6} we generalize \S\ref{ds66} and
\S\ref{ds104} to d-stacks with corners. Here are the analogues of
Definition \ref{ds6def4}, Lemma \ref{ds6lem} and
Theorem~\ref{ds6thm4}:

\begin{dfn} Let $\bs g:\eX\ra\eZ$ and $\bs h:\eY\ra\eZ$ be
1-morphisms in $\dStac$. As in \S\ref{ds131} we have line bundles
$\cN_\eX,\cN_\eZ$ over the $C^\iy$-stacks $\cpX,\cpZ$, and
\S\ref{ds133} defines a $C^\iy$-substack $\cS_{\bs
g}\subseteq\cpX\t_\cZ\cpZ$. As in \cite[\S 11.1]{Joyc6}, there is a
natural isomorphism $\la_{\bs g}:u_{\bs g}^*(\cN_\eZ)\ra s_{\bs
f}^*(\cN_\eX)$ in $\qcoh(\cS_{\bs g})$. The same holds for~$\bs h$.

We call $\bs g,\bs h$ {\it b-transverse\/} if the following holds.
Suppose $x:\bar{\ul *}\ra\cX$ and $y:\bar{\ul *}\ra\cY$ are
1-morphisms in $\CSta$, and $\eta:g\ci x\Ra h\ci y$ is a 2-morphism.
Since $i_\eX:\cpX\ra\cX$ is finite and strongly representable, there
are finitely many 1-morphisms $x':\bar{\ul *}\ra\cpX$ with
$x=i_\eX\ci x'$. Write these $x'$ as $x_1',\ldots,x_j'$. Similarly,
write $y_1',\ldots,y_k'$ for the 1-morphisms $y':\bar{\ul *}\ra\cpY$
with $y=i_\eY\ci y'$. Write $z=g\ci x$ and $\ti z=h\ci y$, so that
$z,\ti z:\bar{\ul *}\ra\cZ$ and $\eta:z\Ra\ti z$. Write
$z_1',\ldots,z_l'$ for the 1-morphisms $z':\bar{\ul *}\ra\cpZ$ with
$z=i_\eZ\ci z'$. Then by Proposition \ref{ds8prop2}, for each
$c=1,\ldots,l$ there are unique $\ti z_c':\bar{\ul *}\ra\cpZ$ and
$\eta_c':z_c'\Ra\ti z_c'$ with $i_\eZ\ci\ti z_c'=\ti z$
and~$\id_{i_\eZ}*\eta_c'=\eta$.

Definition \ref{ds13def1} defined $\cS_{\bs g}\subseteq
\cpX\t_\cZ\cpZ$ in terms of points $[x',z']$ in \eq{ds13eq2}; write
$(x',z'):\ul{\bar *}\ra\cS_{\bs g}$ for the corresponding
1-morphisms. Then we require that for all such $x,y,\eta$, the
following morphism in $\qcoh(\bar{\ul *})$ is injective:
\begin{align*}
&\bigop_{a=1,\ldots,j,\; c=1,\ldots,l:\; [x_a',z_c']\in\cS_{\bs
g,\top}\!\!\!\!\!\!\!\!\!\!\!\!\!\!\!\!\!\!\!\!\!\!\!\!\!\!
\!\!\!\!\!\!\!\!\!\!\!\!\!\!\!\!\!\!\!\!\!\!\!\!\!\!\!
\!\!\!\!\!\!\!\!\!\!\!\!\!\!\!\!\!\!\!\!\!\!\!\!}
I_{(x_a',z_c'),s_{\bs g}}(\cN_\eX)^{-1}\ci (x_a',z_c')^*(\la_{\bs
g})\ci I_{(x_a',z_c'),u_{\bs g}}(\cN_\eZ)\op{}\\
&\bigop_{b=1,\ldots,k,\; c=1,\ldots,l:\; [y_b',\ti z_c']\in\cS_{\bs
h,\top}\!\!\!\!\!\!\!\!\!\!\!\!\!\!\!\!\!\!\!\!\!\!\!\!\!\!
\!\!\!\!\!\!\!\!\!\!\!\!\!\!\!\!\!\!\!\!\!\!\!\!\!\!\!
\!\!\!\!\!\!\!\!\!\!\!\!\!\!\!\!\!\!\!\!\!\!\!\!} I_{(y_b',\ti
z_c'),s_{\bs h}}(\cN_\eY)^{-1}\!\ci\! (y_b',\ti z_c')^*(\la_{\bs
h})\!\ci\!I_{(y_b',\ti z_c'),u_{\bs h}}(\cN_\eZ)\!\ci\!
(\eta_c')^*(\cN_\eZ):\\
&\bigop\nolimits_{c=1}^l(z_c')^*(\cN_\eZ)\longra
\bigop\nolimits_{a=1}^j(x_a')^*(\cN_\eX)\op
\bigop\nolimits_{b=1}^k(y_b')^*(\cN_\eY).
\end{align*}

We call $\bs g,\bs h$ {\it c-transverse\/} if the following holds.
Identify $C_k(\cX)_\top\subseteq C(\cX)_\top$ with the right hand of
\eq{ds13eq4}, and similarly for $C(\cY)_\top,C(\cZ)_\top$. Then
$C(g)_\top,\ab C(h)_\top,\ab\hat C(g)_\top,\ab\hat C(h)_\top$ act as
in \eq{ds13eq5}--\eq{ds13eq6}. We require that:
\begin{itemize}
\setlength{\itemsep}{0pt}
\setlength{\parsep}{0pt}
\item[(a)] whenever there are points in
$C_j(\cX)_\top,C_k(\cY)_\top,C_l(\cZ)_\top$ with
\begin{equation*}
C(g)_\top\bigl([x,\{x_1',\ldots,x_j'\}]\bigr)\!=\!
C(h)_\top\bigl([y,\{y_1',\ldots,y_k'\}]\bigr)\!=\!
[z,\{z_1',\ldots,z_l'\}],
\end{equation*}
we have either $j+k>l$ or $j=k=l=0;$ and
\item[(b)] whenever there are points in
$C_j(\cX)_\top,C_k(\cY)_\top,C_l(\cZ)_\top$ with
\begin{equation*}
\hat C(g)_\top\bigl([x,\{x_1',\ldots,x_j'\}]\bigr)\!=\!
\hat C(h)_\top\bigl([y,\{y_1',\ldots, y_k'\}]\bigr)\!=\!
[z,\{z_1',\ldots,z_l'\}],
\end{equation*}
we have~$j+k\ge l$.
\end{itemize}
Then $\bs g,\bs h$ c-transverse implies $\bs g,\bs h$ b-transverse.
\label{ds13def3}
\end{dfn}

\begin{lem} Let\/ $\bs g:\eX\ra\eZ$ and\/ $\bs h:\eY\ra\eZ$ be\/
$1$-morphisms in $\dStac$. The following are sufficient conditions
for\/ $\bs g,\bs h$ to be c-transverse, and hence b-transverse:
\begin{itemize}
\setlength{\itemsep}{0pt}
\setlength{\parsep}{0pt}
\item[{\bf(i)}] $\bs g$ or $\bs h$ is semisimple\I{d-stack with
corners!semisimple 1-morphism} and flat;\I{d-stack with
corners!flat 1-morphism} or
\item[{\bf(ii)}] $\eZ$ is a d-stack without boundary.
\end{itemize}
\label{ds13lem}
\end{lem}

\begin{thm}{\bf(a)} All b-transverse fibre products exist
in\/~$\dStac$.\I{d-stack with corners!fibre
products!b-transverse}\I{fibre product!of d-stacks with corners}
\smallskip

\noindent{\bf(b)} The $2$-functor\/ $F_\dSpac^\dStac:\dSpac\ra
\dStac$ of\/ {\rm\S\ref{ds131}} takes b- and c-transverse fibre
products in\/ $\dSpac$ to b- and c-transverse fibre products
in\/~$\dStac$.
\smallskip

\noindent{\bf(c)} The $2$-functor $F_\Orbc^\dStac$ of\/
{\rm\S\ref{ds131}} takes transverse fibre products\I{orbifold with
corners!transverse fibre products} in $\Orbc$ to b-transverse fibre
products in $\dStac$. That is, if
\begin{equation*}
\xymatrix@C=90pt@R=10pt{ \oW \ar[r]_(0.2){f} \ar[d]^{e}
\drtwocell_{}\omit^{}\omit{^{\eta}\,\,\,{}} & \oY \ar[d]_{h} \\
\oX \ar[r]^(0.7){g} & \oZ}
\end{equation*}
is a $2$-Cartesian\I{2-category!2-Cartesian square} square in
$\Orbc$ with\/ $g,h$ transverse, and\/ $\eW,\eX,\eY,\eZ,\ab\bs
e,\ab\bs f,\ab\bs g,\ab\bs
h,\bs\eta=F_\Orbc^\dStac(\oW,\oX,\oY,\oZ,e,f,g,h,\eta),$ then
\begin{equation*}
\xymatrix@C=90pt@R=10pt{ \eW \ar[r]_(0.2){\bs f} \ar[d]^{\bs e}
\drtwocell_{}\omit^{}\omit{^{\bs\eta}\,\,\,{}} & \eY \ar[d]_{\bs h} \\
\eX \ar[r]^(0.7){\bs g} & \eZ}
\end{equation*}
is $2$-Cartesian in $\dStac,$ with $\bs g,\bs h$ b-transverse. If
also $g,h$ are strongly transverse\I{orbifold with corners!strongly
transverse 1-morphisms} in\/ $\Orbc,$ then $\bs g,\bs h$ are
c-transverse in\/~$\dStac$.
\smallskip

\noindent{\bf(d)} Suppose we are given a
$2$-Cartesian\I{2-category!2-Cartesian square} diagram
in\/~{\rm$\dStac$:}\I{d-stack with corners!fibre products!boundary
and corners|(}
\begin{equation*}
\xymatrix@C=60pt@R=10pt{ \eW \ar[r]_(0.25){\bs f} \ar[d]^{\bs e}
\drtwocell_{}\omit^{}\omit{^{\bs\eta}}
 & \eY \ar[d]_{\bs h} \\ \eX \ar[r]^(0.7){\bs g} & \eZ,}
\end{equation*}
with\/ $\bs g,\bs h$ c-transverse. Then the following are also
$2$-Cartesian in\/~$\dStac\!:$\I{d-stack with corners!corner
functors}
\ea
\begin{gathered}
\xymatrix@C=100pt@R=10pt{ C(\eW) \ar[r]_(0.25){C(\bs f)}
\ar[d]^{C(\bs e)} \drtwocell_{}\omit^{}\omit{^{C(\bs\eta)
\,\,\,\,\,\,\,\,{}}} & C(\eY) \ar[d]_{C(\bs h)} \\
C(\eX) \ar[r]^(0.7){C(\bs g)} & C(\eZ),}
\end{gathered}
\label{ds13eq7}\\
\begin{gathered}
\xymatrix@C=100pt@R=10pt{ C(\eW) \ar[r]_(0.25){\hat C(\bs f)}
\ar[d]^{\hat C(\bs e)} \drtwocell_{}\omit^{}\omit{^{\hat C(\bs\eta)
\,\,\,\,\,\,\,\,{}}} & C(\eY) \ar[d]_{\hat C(\bs h)} \\
C(\eX) \ar[r]^(0.7){\hat C(\bs g)} & C(\eZ).}
\end{gathered}
\label{ds13eq8}
\ea
Also \eq{ds13eq7}--\eq{ds13eq8} preserve gradings, in that they
relate points in $C_i(\eW),\ab C_j(\eX),\ab C_k(\eY),\ab C_k(\eZ)$
with\/ $i=j+k-l$. Hence \eq{ds13eq7} implies equivalences
in\/~$\dStac\!:$\I{d-stack with corners!fibre products!boundary and
corners|)}
\begin{align*}
C_i(\eW)&\simeq \coprod_{j,k,l\ge 0:i=j+k-l} C_j^{\bs
g,l}(\eX)\t_{C_j^l(\bs g),C_l(\eZ),C_k^l(\bs h)}C_k^{\bs h,l}(\eY),\\
\pd\eW&\simeq \coprod_{j,k,l\ge 0:j+k=l+1} C_j^{\bs
g,l}(\eX)\t_{C_j^l(\bs g),C_l(\eZ),C_k^l(\bs h)}C_k^{\bs h,l}(\eY).
\end{align*}
\label{ds13thm5}
\end{thm}

The analogue of Proposition \ref{ds6prop2} also holds
in~$\dStac$.\I{d-stack with corners!fibre
products|)}\I{b-transversality|)}\I{c-transversality|)}\I{d-stack
with corners!b-transverse 1-morphisms|)}\I{d-stack with
corners!c-transverse 1-morphisms|)}\I{2-category!fibre products
in|)}

\subsection{Orbifold strata of d-stacks with corners}
\label{ds137}
\I{d-stack with corners!orbifold strata|(}\I{orbifold strata!of
d-stacks with corners|(}

In \cite[\S 11.7]{Joyc6} we combine material in \S\ref{ds105} and
\S\ref{ds125} on orbifold strata of d-stacks and of orbifolds with
corners. It is also related to \S\ref{ds67} on fixed loci in
d-spaces with corners. Here is the analogue of
Theorem~\ref{ds10thm7}.

\begin{thm} Let\/ $\eX$ be a d-stack with corners, and\/ $\Ga$ a
finite group. Then we can define d-stacks with corners\/
$\eX^\Ga,\etX^\Ga,\ehX^\Ga,$\G[XGai]{$\eX^\Ga,\etX{}^\Ga,\ehX{}^\Ga,
\eX{}^\Ga_\ci,\etX{}^\Ga_\ci,\ehX{}^\Ga_\ci$}{orbifold strata of a
d-stack with corners $\eX$} and open d-substacks\/
$\eX^\Ga_\ci\subseteq\eX^\Ga,$ $\etX^\Ga_\ci\subseteq\etX^\Ga,$
$\ehX^\Ga_\ci\subseteq \ehX^\Ga,$ all natural up to $1$-isomorphism
in $\dStac,$ a d-space with corners\/ $\bs{\hat{\rm X}}{}^\Ga_\ci$
natural up to $1$-isomorphism in $\dSpac,$ and\/ $1$-morphisms $\bs
O^\Ga(\eX),\bs{\ti\Pi}{}^\Ga (\eX),\ldots$ fitting into a strictly
commutative diagram in\/~$\dStac\!:$\G[OGaXd]{$\bs
O{}^\Ga(\eX),\bs{\ti O}{}^\Ga(\eX),\bs O{}^\Ga_\ci(\eX),\bs{\ti
O}{}^\Ga_\ci(\eX)$}{1-morphisms of orbifold strata
$\eX^\Ga,\ldots,\ehX{}^\Ga_\ci$ of a d-stack with corners
$\eX$}\G[PiGaXd]{$\bs{\ti\Pi}{}^\Ga(\eX),\bs{\hat\Pi}{}^\Ga(\eX),
\bs{\ti\Pi}{}^\Ga_\ci(\eX),\bs{\hat\Pi}{}^\Ga_\ci(\eX)$}{1-morphisms
of orbifold strata $\eX^\Ga,\ldots,\ehX{}^\Ga_\ci$ of a d-stack with
corners $\eX$}
\e
\begin{gathered}
\xymatrix@C=48pt@R=7pt{ \eX^\Ga_\ci
\ar[rr]^{\bs{\ti\Pi}{}^\Ga_\ci(\eX)} \ar[dr]_(0.3){\bs
O{}^\Ga_\ci(\eX)} \ar[dd]_\subset \ar@(ul,l)[]_(0.8){\Aut(\Ga)} &&
\etX^\Ga_\ci \ar[r]^(0.35){\bs{\hat\Pi}{}^\Ga_\ci(\eX)}
\ar[dl]^(0.3){\bs{\ti O}{}^\Ga_\ci(\eX)} \ar[dd]^\subset &
{\ehX{}^\Ga_\ci\simeq F_\dSpac^\dStac(\bs{\hat{\rm
X}}{}^\Ga_\ci)\!\!\!\!\!\!\!\!\!\!\!\!}
\ar@<.5ex>[dd]^\subset \\ & \eX \\
\eX^\Ga \ar[rr]_{\bs{\ti\Pi}{}^\Ga(\eX)} \ar[ur]^(0.3){\bs
O^\Ga(\eX)} \ar@(dl,l)[]^(0.8){\Aut(\Ga)} && \etX^\Ga
\ar[r]_{\bs{\hat\Pi}{}^\Ga(\eX)} \ar[ul]_(0.3){\bs{\ti O}^\Ga(\eX)}
&
*+[r]{\ehX^\Ga.} }\!\!\!\!
\end{gathered}
\label{ds13eq9}
\e
We will call\/ $\eX^\Ga,\etX^\Ga,\ehX^\Ga,\eX^\Ga_\ci,
\etX^\Ga_\ci,\ehX^\Ga_\ci,\bs{\hat{\rm X}}{}^\Ga_\ci$ the
\begin{bfseries}orbifold strata\end{bfseries} of\/~$\eX$.

The underlying d-stacks of\/ $\eX^\Ga,\ldots,\ehX^\Ga_\ci$ are the
orbifold strata\/ $\bcX^\Ga,\ldots,\bhcX{}^\Ga_\ci$ from\/
{\rm\S\ref{ds105}} of the d-stack\/ $\bcX$ in\/ $\eX$. The
$1$-morphisms {\rm\eq{ds13eq9},} as $1$-morphisms in $\dSta,$ are
those given in\/~\eq{ds10eq5}.
\label{ds13thm6}
\end{thm}

The rest of \S\ref{ds105} also extends to~$\dStac$:

\begin{thm} Theorems\/ {\rm\ref{ds10thm8}, \ref{ds10thm9},
\ref{ds10thm10}} and Corollary {\rm\ref{ds10cor}} hold without
change in\/ $\dStac,\dSpac$ rather than\/~$\dSta,\dSpa$.
\label{ds13thm7}
\end{thm}

Here are analogues of Theorem \ref{ds12thm5} and
Corollary~\ref{ds12cor}.

\begin{thm} Let\/ $\eX$ be a d-stack with corners, and\/ $\Ga$ a
finite group. The corners $C(\eX)$ from {\rm\S\ref{ds135}} lie in
$\dStac,$ so Theorem\/ {\rm\ref{ds13thm6}} gives orbifold strata
$\eX^\Ga,C(\eX)^\Ga$ and\/ $1$-morphisms $\bs
O^\Ga(\eX):\eX^\Ga\ra\eX,$ $\bs O^\Ga(C(\eX)):C(\eX)^\Ga\ra C(\eX)$.
Applying the corner functor $C$ from {\rm\S\ref{ds135}} gives a
$1$-morphism $C(\bs O^\Ga(\eX)):C(\eX^\Ga)\ra C(\eX)$. There exists
a unique equivalence $\bs K^\Ga(\eX):C(\eX^\Ga)\ra C(\eX)^\Ga$ in
$\dStac$ with\/ $\bs O^\Ga(C(\eX))\ci \bs K^\Ga(\eX)=C(\bs
O^\Ga(\eX)):C(\eX^\Ga)\ra C(\eX)$. It restricts to an equivalence
$\bs K^\Ga_\ci(\eX):=\bs K^\Ga(\eX)\vert_{\smash{C(\eX_\ci^\Ga)
}}:C(\eX^\Ga_\ci)\ra C(\eX)^\Ga_\ci$.

Similarly, there is a unique equivalence $\bs{\ti K}{}^\Ga(\eX):
C(\etX{}^\Ga)\ra\,\,\,\widetilde{\!\!\!C(\eX)\!\!\!}\,\,\,^\Ga$
with\/ $\bs{\ti O}{}^\Ga(C(\eX))\!\ci\!\bs{\ti K}{}^\Ga(\eX)\!=\!
C(\bs{\ti O}{}^\Ga(\eX))$ and\/ $\bs{\ti\Pi}{}^\Ga(C(\eX))\!\ci\!\bs
K{}^\Ga(\eX)\!=\!\bs{\ti K}{}^\Ga(\eX)\!\ci\!
C(\bs{\ti\Pi}{}^\Ga(\eX))$. There is an equivalence\/ $\bs{\hat
K}{}^\Ga(\eX):C(\hat\eX{}^\Ga) \ra \,\,\,\widehat{\!\!\!
C(\eX)\!\!\!}\,\,\,^\Ga,$ unique up to $2$-isomorphism, with a
$2$-morphism $\bs{\hat\Pi}{}^\Ga(C(\eX))\ci\bs{\ti
K}{}^\Ga(\eX)\Ra\bs{\hat K}{}^\Ga(\eX)\ci
C(\bs{\hat\Pi}{}^\Ga(\eX))$. They restrict to equivalences $\bs{\ti
K}{}^\Ga_\ci(\eX):C(\etX{}^\Ga_\ci)\ra\,\,\,\widetilde{\!\!\!
C(\eX)\!\!\!}\,\,\,^\Ga_\ci$ and\/ $\bs{\hat K}{}^\Ga_\ci(\eX):
C(\ehX{}^\Ga_\ci)\ra \,\,\,\widehat{\!\!\!
C(\eX)\!\!\!}\,\,\,^\Ga_\ci$.
\label{ds13thm8}
\end{thm}

\begin{cor} Let\/ $\eX$ be a d-stack with corners, and\/ $\Ga$ a
finite group. Then there exist\/ $1$-morphisms $\bs
J{}^\Ga(\eX):(\pd\eX)^\Ga\ra \pd(\eX^\Ga),$ $\bs{\ti J}{}^\Ga(\eX):
\,\,\,\,\widetilde{\!\!\!\!(\pd\eX)\!\!\!\!}\,\,\,\,^\Ga\ra
\pd(\etX{}^\Ga),$ $\bs{\hat J}{}^\Ga(\eX):\,\,\,\,\widehat{\!\!\!\!
(\pd\eX)\!\!\!\!}\,\,\,\,^\Ga\ra\pd(\ehX{}^\Ga)$ in\/ $\dStac,$
natural up to\/ $2$-isomorphism, such that\/ $\bs J{}^\Ga(\eX)$ is
an equivalence from $(\pd\eX)^\Ga$ to an open and closed d-substack
of\/ $\pd(\eX^\Ga),$ and similarly for\/~$\bs{\ti
J}{}^\Ga(\eX),\bs{\hat J}{}^\Ga(\eX)$.
\label{ds13cor}
\end{cor}

A d-stack with corners $\eX$ is called {\it straight\/}\I{d-stack
with corners!straight} if
$(i_\oX)_*:\Iso_\cpX([x'])\ra\Iso_\cX([x])$ is an isomorphism for
all $[x']\in\cpX_\top$ with $i_{\oX,\top}([x'])=[x]$. D-stacks with
boundary\I{d-stack with boundary} are automatically straight. If
$\eX$ is straight then $\pd\eX$ is straight, so by induction
$\pd^k\eX$ is also straight for all $k\ge 0$. If $\eX$ is straight
then $\bs K^\Ga(\eX)$ in Theorem \ref{ds13thm8} is an equivalence
$C_k(\eX^\Ga)\ra C_k(\eX)^\Ga$ for all $k\ge 0$, and so $\bs
J{}^\Ga(\eX)$ in Corollary \ref{ds13cor} is an equivalence
$(\pd\eX)^\Ga\ra \pd(\eX^\Ga)$. The same applies for~$\bs{\ti
J}{}^\Ga(\eX),\bs{\hat J}{}^\Ga(\eX),\bs{\ti K}{}^\Ga(\eX),\bs{\hat
K}{}^\Ga(\eX)$.\I{d-stack with corners!orbifold strata|)}\I{orbifold
strata!of d-stacks with corners|)}\I{d-stack with corners|)}

\section{D-orbifolds with corners}
\label{ds14}
\I{d-orbifold with corners|(}

In \cite[Chap.~12]{Joyc6} we discuss the 2-category $\dOrbc$ of {\it
d-orbifolds with corners}. Again, there are few new issues here:
almost all the material combines ideas we have seen already on
d-manifolds with corners from \S\ref{ds7}, on orbifolds with corners
from \S\ref{ds12}, and on d-stacks with corners from \S\ref{ds13}.
So we will be brief.

As we explain in \S\ref{ds16} and \cite[\S 14.3]{Joyc6}, d-orbifolds
with corners are related to Kuranishi spaces\I{Kuranishi space!and
d-orbifolds} in the work of Fukaya, Oh, Ohta and Ono~\cite{FOOO}.

\subsection{Definition of d-orbifolds with corners}
\label{ds141}
\I{d-orbifold with corners!definition|(}

In \cite[\S 12.1]{Joyc6} we define d-orbifolds with corners,
following \S\ref{ds71} and~\S\ref{ds111}.

\begin{dfn} A d-stack with corners $\eW$ is called a {\it principal
d-orbifold with corners\/}\I{principal d-orbifold with
corners}\I{d-orbifold with corners!principal} if is equivalent in
$\dStac$ to a fibre product $\eV\t_{\bs s,\eE,\bs 0}\eV$, where
$\oV$ is an orbifold with corners, $\cE$ is a vector bundle on
$\oV$, $s\in C^\iy(\cE)$, and $\eV,\eE,\bs s,\bs
0=F_\Orbc^\dStac\bigl(\oV,\Totc(\cE), \Totc(s),\Totc(0)\bigr)$, for
$\Totc$\I{orbifold with corners!vector bundles on!total space
functor $\Totc$} as in \S\ref{ds121}. Note that
$\Totc(s),\Totc(0):\oV\ra\Totc(\cE)$ are simple, flat 1-morphisms in
$\Orbc$, so $\bs s,\bs 0:\eV\ra\eE$ are simple, flat 1-morphisms in
$\dStac$. Thus $\bs s,\bs 0$ are b-transverse\I{d-stack with
corners!b-transverse 1-morphisms} by Lemma \ref{ds13lem}(i), and
$\eV\t_{\bs s,\eE,\bs 0}\eV$ exists in $\dStac$ by
Theorem~\ref{ds13thm5}(a).

If $\eW$ is a nonempty principal d-orbifold with corners, then
$T^*\bcW$ is a virtual vector bundle. We define the {\it virtual
dimension\/} of $\eW$ to be $\vdim\eW=\rank T^*\eW\in\Z$. If
$\eW\simeq\eV\t_{\bs s,\eE,\bs 0}\eV$
then~$\vdim\eW=\dim\oV-\rank\cE$.

A d-stack with corners $\eX$ is called a {\it d-orbifold with
corners of virtual dimension\/}\I{d-orbifold with corners!virtual
dimension} $n\in\Z$, written $\vdim\eX=n$, if $\eX$ can be covered
by open d-substacks $\eW$ which are principal d-orbifolds with
corners with $\vdim\eW=n$. A d-orbifold with corners $\eX$ is called
a {\it d-orbifold with boundary\/}\I{d-orbifold with boundary} if it
is a d-stack with boundary, and a {\it d-orbifold without
boundary\/} if it is a d-stack without boundary.

Write $\bdOrb,\dOrbb,\dOrbc$\G[dOrbc]{$\dOrbc$}{2-category of
d-orbifolds with corners}\G[dOrbb]{$\dOrbb$}{2-category of
d-orbifolds with boundary}\G[dOrb']{$\bdOrb$}{2-subcategory of
d-orbifolds with corners equivalent to d-orbifolds} for the full
2-subcategories of d-orbifolds without boundary, with boundary, and
with corners, in $\dStac$, respectively.\I{d-orbifold with
corners!include d-orbifolds}\I{d-orbifold!as d-orbifold with
corners} Then $\bOrb,\bOrbb,\bOrbc$ in \S\ref{ds131} are full
2-subcategories of $\bdOrb,\ab\dOrbb,\ab \dOrbc$. When we say that a
d-orbifold with corners $\eX$ {\it is an orbifold},\I{d-orbifold
with corners!is an orbifold} we mean that $\eX$ lies in $\bOrbc$.
Define full and faithful\I{functor!full}\I{functor!faithful} strict
2-functors\I{2-category!strict 2-functor}
\begin{gather*}
\begin{aligned}
F_\dOrb^\dOrbc&:\dOrb\ra\bdOrb\subset\dOrbc, &
F_\Orbc^\dOrbc&:\Orbc\ra\dOrbc, \\
F_\Orbb^\dOrbc&:\Orbb\ra\dOrbb\!\subset\!\dOrbc, &
F_\Orb^\dOrbc&:\Orb\ra\bdOrb\subset\dOrbc, \\
F_\dManc^\dOrbc&:\dManc\ra\dOrbc, & F_\dManb^\dOrbc&:\dManb\ra\dOrbb
\!\subset\!\dOrbc,
\end{aligned}\\
\text{and}\qquad F_\dMan^\dOrbc:\dMan\ra\bdOrb\subset\dOrbc,
\qquad \text{by} \\
\begin{aligned}
F_\dOrb^\dOrbc\!&=F_\dSta^\dStac\vert_{\dOrb}, &
F_\Orbc^\dOrbc&=F_\Orbc^\dStac, &
F_\Orbb^\dOrbc&=F_\Orbc^\dStac\vert_{\smash{\Orbb}}, \\
F_\Orb^\Orbc\!&=\!F_\dSta^\dStac\!\ci\!F_\Orb^\dSta,\!\!\! &
F_\dManc^\dOrbc\!&=\!F_\dSpac^\dStac\vert_{\dManc},\!\!\! &
F_\dManb^\dOrbc\!&=\!F_\dSpac^\dStac\vert_{\dManb},
\end{aligned}\\
\text{and}\qquad F_\dMan^\dOrbc=F_\dSpac^\dStac\ci
F_\dMan^\dManc=F_\dOrb^\dOrbc\ci F_\dMan^\dOrb,
\end{gather*}
where $F_\Orb^\Orbc$, $F_\Orb^\dSta$, $F_\dSta^\dStac$,
$F_\dSpac^\dStac$, $F_\dMan^\dManc$, $F_\dMan^\dOrb$,
$F_\Orbc^\dStac$ are as in \S\ref{ds71}, \S\ref{ds111},
\S\ref{ds121}, and \S\ref{ds131}. Here
$F_\dOrb^\dOrbc:\dOrb\ra\bdOrb$ is an isomorphism of 2-categories.
So we may as well identify $\dOrb$ with its image $\bdOrb$, and
consider d-orbifolds in \S\ref{ds11} as examples of d-orbifolds with
corners.

Write $\hdManc$\G[dManc']{$\hdManc$}{2-subcategory of d-orbifolds
with corners equivalent to d-manifolds with corners} for the full
2-subcategory of objects $\eX$ in $\dOrbc$ equivalent to
$F_\dManc^\dOrbc(\rX)$ for some d-manifold with corners $\rX$. When
we say that a d-orbifold with corners $\eX$ {\it is a
d-manifold},\I{d-orbifold with corners!is a d-manifold} we mean
that~$\eX\in\hdManc$.

These 2-categories lie in a commutative diagram:
\begin{equation*}
\xymatrix@C=36pt@R=21pt{ \dSpa \ar@<1.5ex>@/^1pc/[ddd]^(0.5){F_\dSpa^\dSta}
& \Man \ar[dl] \ar[d]^{F_\Man^\dManc}
\ar[dl]_{F_\Man^\dMan} \ar[r]_\subset
\ar@<1.5ex>@/^1pc/[ddd]^(0.45){F_\Man^\Orb}
& \Manb \ar[d]_{F_\Manb^\dManc} \ar[r]_\subset
\ar@<1.5ex>@/^1pc/[ddd]^(0.45){F_\Manb^\Orbc} & \Manc
\ar@<1.5ex>@/^1pc/[ddd]^(0.45){F_\Manc^\Orbc}
\ar[d]_{F_\Manc^\dManc} \ar[dr]^{F_\Manc^\dSpac} \\
\dMan \ar[u]^{\subset} \ar[r]^(0.65){F_\dMan^\dManc}_(0.6)\cong
\ar[d]_{F_\dMan^\dOrb} &
\bdMan \ar[r]^\subset \ar[d] & \dManb \ar[r]^\subset
\ar[d]_(0.6){F_\dManb^\dOrbc} & \dManc \ar[d]_(0.6){F_\dManc^\dOrbc}
\ar[r]^\subset & \dSpac \ar[d]_(0.6){F_\dSpac^\dStac} \\
\dOrb \ar[d]_{\subset} \ar[r]_(0.65){F_\dOrb^\dOrbc}^(0.6)\cong &
\bdOrb \ar[r]_\subset & \dOrbb \ar[r]_\subset
\ar[d] & \dOrbc \ar[r]_\subset & \dStac \\
\dSta & \Orb \ar[u]_{F_\Orb^\dOrbc} \ar[ul]^{F_\Orb^\dOrb}
\ar[r]^{F_\Orb^\Orbc} & \Orbb \ar[u]^{F_\Orbb^\dOrbc}
\ar[r]^\subset & \Orbc. \ar[u]^{F_\Orbc^\dOrbc} \ar[ur]_{F_\Orbc^\dStac} }
\end{equation*}

If $\eX=(\bcX,\bcpX,\bs i_\eX,\om_\eX)$ is a d-orbifold with
corners, then the virtual cotangent sheaf $T^*\bcX$ of the d-stack
$\bcX$ from Remark \ref{ds11rem} is a virtual vector bundle on
$\cX$, of rank $\vdim\eX$. We will call $T^*\bcX\in\vvect(\cX)$ the
{\it virtual cotangent bundle\/}\I{d-orbifold with corners!virtual
cotangent bundle}\I{virtual cotangent bundle} of $\eX$, and also
write it~$T^*\eX$.\G[T*Xd]{$T^*\eX$}{virtual cotangent sheaf of a
d-stack with corners $\eX$}\I{d-orbifold with corners!definition|)}
\label{ds14def1}
\end{dfn}

Here is the analogue of Lemma~\ref{ds11lem}:

\begin{lem} Let\/ $\eX$ be a d-orbifold with corners. Then $\eX$ is
a d-manifold,\I{d-orbifold with corners!is a d-manifold} that is,
$\eX\simeq F_\dManc^\dOrbc(\rX)$ for some d-manifold with corners\/
$\rX,$ if and only if\/ $\Iso_\cX([x])\cong\{1\}$ for all\/ $[x]$
in\/~$\cX_\top$.
\label{ds14lem1}
\end{lem}

D-orbifolds with corners are preserved by boundaries and corners.

\begin{prop} Suppose $\eX$ is a d-orbifold with corners. Then
$\pd\eX$ in\/ {\rm\S\ref{ds131}} and\/ $C_k(\eX)$ in\/
{\rm\S\ref{ds135}} are d-orbifolds with corners, with\/
$\vdim\pd\eX=\vdim\eX-1$ and\/ $\vdim C_k(\eX)=\vdim\eX-k$ for
all\/~$k\ge 0$.
\label{ds14prop1}
\end{prop}

\begin{dfn} As for $\dcManc$ in \S\ref{ds71}, define
$\dcOrbc$\G[dOrbc']{$\dcOrbc$}{2-category of disjoint unions of
d-orbifolds with corners of different dimensions} to be the full
2-sub\-cat\-eg\-ory of $\eX$ in $\dStac$ which may be written as a
disjoint union $\eX=\coprod_{n\in\Z}\eX_n$ for $\eX_n\in\dOrbc$ with
$\vdim\eX_n=n$, where we allow $\eX_n=\bs\es$. Then
$\dOrbc\subset\dcOrbc\subset\dStac$, and the corner functors $C,\hat
C:\dStac\ra\dStac$\I{d-stack with corners!corner
functors}\I{d-orbifold with corners!corner functors} in
\S\ref{ds135} restrict to strict 2-functors~$C,\hat
C:\dOrbc\ra\dcOrbc$.\I{2-category!strict 2-functor}
\label{ds14def2}
\end{dfn}

\subsection{Local properties of d-orbifolds with corners}
\label{ds142}
\I{d-orbifold with corners!local properties|(}\I{d-orbifold with
corners!standard model a@standard model $\eS_{\oV,\cE,s}$|(}

In \cite[\S 12.2]{Joyc6} we combine \S\ref{ds72} and \S\ref{ds112}.
Here are analogues of Examples \ref{ds7ex2}, \ref{ds7ex3} and
Theorem \ref{ds7thm1}, and of Examples \ref{ds11ex1}, \ref{ds11ex2}
and Theorem~\ref{ds11thm1}.

\begin{ex} Let $\oV=(\cV,\cpV,i_\oV)$ be an orbifold with
corners, $\cE$ a vector bundle on $\oV$ as in \S\ref{ds121}, and
$s\in C^\iy(\cE)$. We will define an explicit principal d-orbifold
with corners~$\eS=(\bcS,\bcpS,\bs i_\eS,\om_\eS)$

Define a vector bundle $\cE_\pd$ on $\cpV$ by
$\cE_\pd=i_\oV^*(\cE)$, and a section $s_\pd=i_\oV^*(s)\in
C^\iy(\cE_\pd)$. Define d-stacks $\bcS=\bcS_{\cV,\cE,s}$ and
$\bcpS=\bS_{\cpV,\cE_\pd,s_\pd}$ from the triples $\cV,\cE,s$ and
$\cpV,\cE_\pd,s_\pd$ exactly as in Example \ref{ds11ex1}, although
now $\cV,\cpV$ have corners. Define a 1-morphism $\bs
i_\bcS:\bcpS\ra\bS$ in $\dSta$ to be the `standard model' 1-morphism
$\bcS_{i_\oV,\id_{\cE_\pd}}:\bcS_{\cpV,\cE_\pd,s_\pd}\ra\bcS_{\cV,\cE,s}$
from Example~\ref{ds11ex2}.

As in Example \ref{ds7ex2}, the conormal bundle $\cN_\eS$ of $\bcpS$
in $\bcS$ is canonically isomorphic to the lift to
$\cpS\subseteq\cpV$ of the conormal bundle $\cN_\oV$ of $\cpV$ in
$\cV$. Define $\om_\eS$ to be the orientation on $\cN_\eS$ induced
by the orientation on $\cN_\oV$ by outward-pointing normal vectors
to $\cpV$ in $\cV$. Then $\eS=(\bcS,\bcpS,\bs i_\eS,\om_\eS)$ is a
d-stack with corners. It is equivalent in $\dStac$ to $\eV\t_{\bs
s,\eE,\bs 0}\eV$ in Definition \ref{ds14def1}. We call $\eS$ a
`standard model' d-orbifold with corners, and write
it~$\eS_{\oV,\cE,s}$.\G[SVEsd]{$\eS_{\oV,\oE,s}$}{`standard model'
d-orbifold with corners}

Every principal d-orbifold with corners\I{principal d-orbifold with
corners}\I{d-orbifold with corners!principal} $\eW$ is equivalent in
$\dOrbc$ to some $\eS_{\oV,\cE,s}$. Sometimes it is useful to take
$\oV$ to be an {\it effective\/} orbifold with corners,\I{orbifold
with corners!effective} as in \S\ref{ds121}. There is a natural
1-isomorphism $\pd\eS_{\oV,\cE,s}\cong\eS_{\pd
\oV,\cE_\pd,s_\pd}$\I{d-orbifold with corners!standard model
a@standard model $\eS_{\oV,\cE,s}$!boundary} in~$\dOrbc$.
\label{ds14ex1}
\end{ex}

\begin{ex} Let $\oV,\oW$ be orbifolds with corners, $\cE,\cF$ be
vector bundles on $\oV,\oW$, and $s\in C^\iy(\cE)$, $t\in
C^\iy(\cF)$, so that Example \ref{ds14ex1} defines d-orbifolds with
corners $\eS_{\oV,\cE,s},\eS_{\oW,\cF,t}$. Suppose $f:\oV\ra\oW$ is
a 1-morphism in $\Orbc$, and $\hat f:\cE\ra f^*(\cF)$ is a morphism
in $\vect(\cV)$ satisfying $\hat f\ci s= f^*(t)$, as
in~\eq{ds11eq1}.

The d-stacks $\bcS_{\oV,\cE,s},\bcS_{\oW,\cF,t}$ in
$\eS_{\oV,\cE,s},\eS_{\oW,\cF,t}$ are defined as for `standard
model' d-orbifolds $\bcS_{\cV,\cE,s}$ in Example \ref{ds11ex1}. Thus
we can follow Example \ref{ds11ex2} to define a 1-morphism
$\eS_{\smash{f,\hat f}}:\bcS_{\oV,\cE,s}\ra\bcS_{\oW,\cF,t}$ in
$\dSta$. Then $\eS_{\smash{f,\hat f}}:\eS_{\oV,\cE,s}\ra
\eS_{\oW,\cF,t}$ is a 1-morphism in $\dOrbc$. We call it a
`standard model' 1-morphism.\G[Sffd]{$\eS_{\smash{f,\hat
f}}:\eS_{\oV,\cE,s}\ra \eS_{\oW,\cF,t}$}{`standard model' 1-morphism
in $\dOrbc$}\I{d-orbifold with corners!standard model a@standard
model $\eS_{\oV,\cE,s}$!1-morphism}

Suppose now that $\ti\oV\subseteq\oV$ is open, with inclusion
1-morphism $i_{\smash{\ti\oV}}:\ti\oV\ra\oV$. Write $\ti{\cal
E}=\cE\vert_{\smash{\ti\oV}}=i_{\smash{\ti\oV}}^*(\cE)$ and $\ti
s=s\vert_{\smash{\ti\oV}}$. Define $\bs
i_{\smash{\ti\oV,\oV}}=\eS_{i_{\smash{\ti\oV}},\id_{\ti{\cal E}}}:
\eS_{\smash{\ti\oV,\ti{\cal E},\ti s}}\ra\eS_{\oV,\cE,s}$.\G[iVVd]{$\bs
i_{\smash{\ti\oV,\oV}}:\eS_{\smash{\ti\oV,\ti{\cal E},\ti
s}}\ra\eS_{\oV,\cE,s}$}{inclusion of open set in `standard model'
d-orbifold with corners} If $s^{-1}(0)\subseteq\ti\oV$ then $\bs
i_{\smash{\ti\oV,\oV}}:\eS_{\smash{\ti\oV,\ti{\cal E},\ti
s}}\ra\eS_{\oV,\cE,s}$ is a 1-isomorphism.
\label{ds14ex2}
\end{ex}

\begin{thm} Let\/ $\eX$ be a d-orbifold with corners, and\/
$[x]\in\cX_\top$. Then there exists an open neighbourhood\/ $\eU$
of\/ $[x]$ in $\eX$ and an equivalence $\eU\simeq\eS_{\oV,\cE,s}$ in
$\dOrbc$ for some orbifold with corners\/ $\oV,$ vector bundle $\cE$
over $\oV$ and $s\in C^\iy(\cE)$ which identifies $[x]\in\cU_\top$
with a point\/ $[v]\in S^k(\cV)_\top\subseteq\cV_\top$ such that\/
$s(v)=\d s\vert_{S^k(\cV)}(v)=0,$ where $S^k(\cV)\subseteq\cV$ is
the locally closed\/ $C^\iy$-substack of\/ $[v]\in\cV_\top$ such
that\/ $\bar{\ul *}\t_{v,\cV,i_\oV}\cpV$ is $k$ points, for $k\ge
0$. Furthermore, $\oV,\cE,s,k$ are determined up to non-canonical
equivalence near $[v]$ by $\eX$ near~$[x]$.
\label{ds14thm1}
\end{thm}

As in Examples \ref{ds11ex3}--\ref{ds11ex4} for d-orbifolds, we can
combine the `standard model' d-manifolds with corners $\rS_{V,E,s}$
and 1-morphisms $\rS_{ \smash{f,\hat f}}:\rS_{V,E,s} \ra\rS_{W,F,t}$
of Examples \ref{ds7ex2}--\ref{ds7ex3} with quotient d-stacks with
corners\I{d-stack with corners!quotients $[\rX/G]$} of \S\ref{ds132}
to define an alternative form of `standard model' d-orbifolds with
corners\I{d-orbifold with corners!standard model b@standard model
$[\rS_{V,E,s}/\Ga]$}
$[\rS_{V,E,s}/\Ga]$\G[SVEsf]{$[\rS_{V,E,s}/\Ga]$}{alternative
`standard model' d-orbifold with corners} and `standard model'
1-morphisms~$[\rS_{ \smash{f,\hat
f}},\rho]:[\rS_{V,E,s}/\Ga]\ra[\rS_{W,F,t}/\De]$.\I{d-orbifold with
corners!standard model b@standard model
$[\rS_{V,E,s}/\Ga]$!1-morphism}\G[Sfff]{$[\rS_{ \smash{f,\hat
f}},\rho]:[\rS_{V,E,s}/\Ga]\ra [\rS_{W,F,t}/\De]$}{`standard model'
1-morphism in $\dOrbc$}\I{d-orbifold with corners!local
properties|)}\I{d-orbifold with corners!standard model a@standard
model $\eS_{\oV,\cE,s}$|)}

\subsection{Equivalences in $\dOrbc$, and gluing by equivalences}
\label{ds143}
\I{d-orbifold with corners!equivalence|(}

In \cite[\S 12.3]{Joyc6} we combine \S\ref{ds73} and \S\ref{ds113}.
Here are the analogues of Theorems \ref{ds11thm2}--\ref{ds11thm5}.
\'Etale 1-morphisms in $\dStac$ were discussed
in~\S\ref{ds132}.\I{d-orbifold with corners!etale 1-morphism@\'etale
1-morphism|(}

\begin{thm} Suppose $\bs f:\eX\ra\eY$ is a $1$-morphism of
d-orbifolds with corners, and\/ $f:\cX\ra\cY$ is representable. Then
the following are equivalent:
\begin{itemize}
\setlength{\itemsep}{0pt}
\setlength{\parsep}{0pt}
\item[{\rm(i)}] $\bs f$ is \'etale;
\item[{\rm(ii)}] $\bs f$ is simple\I{d-orbifold with
corners!simple 1-morphism} and flat,\I{d-orbifold with
corners!flat 1-morphism} in the sense of\/ {\rm\S\ref{ds133},}
and\/ $\Om_{\bs f}:f^*(T^*\eY) \ra T^*\eX$ is an equivalence in
$\vqcoh(\cX);$ and
\item[{\rm(iii)}] $\bs f$ is simple and flat, and\/
\eq{ds11eq4} is a split short exact\I{split short exact
sequence}\I{abelian category!split short exact sequence}
sequence in\/~$\qcoh(\cX)$.
\end{itemize}
If in addition $f_*:\Iso_\cX([x])\ra\Iso_\cY(f_\top([x]))$ is an
isomorphism for all\/ $[x]\in\cX_\top,$ and\/
$f_\top:\cX_\top\ra\cY_\top$ is a bijection, then $\bs f$ is an
equivalence in\/~$\dOrbc$.

\label{ds14thm2}
\end{thm}

\begin{thm} Suppose $\eS_{\smash{f,\hat f}}:\eS_{\oV,\cE,s}\ra
\eS_{\oW,\cF,t}$ is a `standard model'\/ $1$-morphism\I{d-orbifold
with corners!standard model a@standard model
$\eS_{\oV,\cE,s}$!1-morphism} in $\dOrbc,$ in the notation of
Examples\/ {\rm\ref{ds14ex1}} and\/ {\rm\ref{ds14ex2},} with\/
$f:\cV\ra\cW$ representable. Then $\eS_{\smash{f,\hat f}}$ is
\'etale if and only if\/ $f$ is simple and flat near
$s^{-1}(0)\subseteq\oV,$ in the sense of\/ {\rm\S\ref{ds122},} and
for each\/ $[v]\in\cV_\top$ with\/ $s(v)=0$ and\/
$[w]=f_\top([v])\in\cW_\top,$ the following sequence is exact:
\begin{equation*}
\smash{\xymatrix@C=18pt{ 0 \ar[r] & T_v\cV \ar[rrr]^(0.42){\d
s(v)\op \,\d f(v)} &&& \cE_v\op T_w\cW \ar[rrr]^(0.57){\hat
f(v)\op\, -\d t(w)} &&& \cF_w \ar[r] & 0.}}
\end{equation*}
Also\/ $\eS_{\smash{f,\hat f}}$ is an equivalence if and only if in
addition\/ $f_\top\vert_{s^{-1}(0)}:s^{-1}(0)\!\ra\! t^{-1}(0)$ is a
bijection, where $s^{-1}(0)\!=\!\{[v]\in\cV_\top:s(v)\!=\!0\},$
$t^{-1}(0)\!=\!\{[w]\in\cW_\top:t(w)\!=\!0\},$ and\/
$f_*:\Iso_\cV([v])\ra\Iso_\cW(f_\top([v]))$ is an isomorphism for
all\/~$[v]\in s^{-1}(0)\subseteq\cV_\top$.\I{d-orbifold with
corners!etale 1-morphism@\'etale 1-morphism|)}
\label{ds14thm3}
\end{thm}

\begin{thm} Suppose we are given the following data:\I{d-orbifold
with corners!gluing by equivalences|(}
\begin{itemize}
\setlength{\itemsep}{0pt}
\setlength{\parsep}{0pt}
\item[{\rm(a)}] an integer $n;$
\item[{\rm(b)}] a Hausdorff, second countable topological space $X;$
\item[{\rm(c)}] an indexing set\/ $I,$ and a total order $<$ on $I;$
\item[{\rm(d)}] for each\/ $i$ in $I,$ an effective orbifold
with corners\/ $\oV_i,$ a vector bundle $\cE_i$ on $\oV_i$
with\/ $\dim\oV_i-\rank\cE_i=n,$ a section $s_i\in
C^\iy(\cE_i),$ and a homeomorphism $\psi_i:s_i^{-1}(0) \ra\hat
X_i,$ where $s_i^{-1}(0)=\bigl\{[v_i]\in\cV_{i,\top}:
s_i(v_i)=0\bigr\}$ and\/ $\hat X_i\subseteq X$ is open; and
\item[{\rm(e)}] for all\/ $i<j$ in $I,$ an open suborbifold\/
$\oV_{ij}\subseteq\oV_i,$ a simple, flat\/ $1$-morphism
$e_{ij}:\oV_{ij}\ra\oV_j,$ and a morphism of vector bundles
$\hat e_{ij}:\cE_i\vert_{\oV_{ij}}\ra e_{ij}^*(\cE_j)$.
\end{itemize}
Let this data satisfy the conditions:
\begin{itemize}
\setlength{\itemsep}{0pt}
\setlength{\parsep}{0pt}
\item[{\rm(i)}] $X=\bigcup_{i\in I}\hat X_i;$
\item[{\rm(ii)}] if\/ $i<j$ in $I$ then $(e_{ij})_*:
\Iso_{\cV_{ij}}([v])\ra\Iso_{\cV_j}(e_{ij,\top}([v]))$ is an
isomorphism for all\/ $[v]\in\cV_{ij,\top},$ and\/ $\hat
e_{ij}\ci s_i\vert_{\oV_{ij}}=e_{ij}^*(s_j)\ci\io_{ij}$ where
$\io_{ij}:\O_{\oV_{ij}}\ra e_{ij}^*(\O_{\oV_j})$ is the natural
isomorphism, and\/ $\psi_i(s_i\vert_{\oV_{ij}}^{-1}(0))=\hat
X_i\cap\hat X_j,$ and\/ $\psi_i\vert_{s_i\vert_{\oV_{ij}}^{-1}
(0)}=\psi_j\ci e_{ij,\top}\vert_{s_i\vert_{\oV_{ij}}^{-1}(0)},$
and if\/ $[v_i]\in\cV_{ij,\top}$ with\/ $s_i(v_i)=0$ and\/
$[v_j]=e_{ij,\top}([v_i])$ then the following sequence is exact:
\begin{equation*}
\smash{\xymatrix@C=16pt{ 0 \ar[r] & T_{v_i}\cV_i \ar[rrr]^(0.42){\d
s_i(v_i)\op \,\d e_{ij}(v_i)} &&& \cE_i\vert_{v_i}\!\op\!
T_{v_j}\cV_j \ar[rrr]^(0.57){\hat e_{ij}(v_i)\op\, -\d s_j(v_j)} &&&
\cE_j\vert_{v_j} \ar[r] & 0;}}\!\!\!
\end{equation*}
\item[{\rm(iii)}] if\/ $i<j<k$ in $I$ then there exists a
$2$-morphism $\eta_{ijk}:e_{jk}\ci e_{ij}\vert_{\oV_{ik}\cap
e_{ij}^{-1}(\oV_{jk})}\ab\Ra e_{ik}\vert_{\oV_{ik}\cap
e_{ij}^{-1}(\oV_{jk})}$ in $\Orbc$ with\/
\begin{equation*}
{}\!\!\!\!\!\!\!\!\!\!\!
\hat e_{ik}\vert_{\oV_{ik}\cap e_{ij}^{-1}(\oV_{jk})}\!=\!
\eta_{ijk}^*(\cE_k)\!\ci\! I_{e_{ij},e_{jk}}(\cE_k)^{-1}\!\ci\!
e_{ij}\vert_{\oV_{ik}\cap e_{ij}^{-1}(\oV_{jk})}^*(\hat e_{jk})\!\ci\!
\hat e_{ij}\vert_{\oV_{ik}\cap e_{ij}^{-1}(\oV_{jk})}.
\end{equation*}
Note that\/ $\eta_{ijk}$ is unique by the corners analogue of
Proposition\/~{\rm\ref{ds9prop1}(ii)}.
\end{itemize}

Then there exist a d-orbifold with corners\/ $\eX$ with\/
$\vdim\eX=n$ and underlying topological space $\cX_\top\cong X,$ and
a $1$-morphism $\bs\psi_i:\eS_{\oV_i,\cE_i,s_i}\ra\eX$ in $\dOrbc$
with underlying continuous map $\psi_i$ which is an equivalence with
the open d-suborbifold\/ $\ehX_i\subseteq\eX$ corresponding to $\hat
X_i\subseteq X$ for all\/ $i\in I,$ such that for all\/ $i<j$ in $I$
there exists a $2$-morphism\/
$\bs\eta_{ij}:\bs\psi_j\ci\eS_{e_{ij},\hat e_{ij}}\Ra\bs\psi_i\ci\bs
i_{\oV_{ij},\oV_i},$ where $\eS_{e_{ij},\hat e_{ij}}:
\eS_{\oV_{ij},\cE_i\vert_{\oV_{ij}},s_i\vert_{\oV_{ij}}}\ra
\eS_{\oV_j,\cE_j,s_j}$ and\/ $\bs i_{\oV_{ij},\oV_i}:
\eS_{\oV_{ij},\cE_i\vert_{\oV_{ij}},s_i\vert_{\oV_{ij}}}\ra
\eS_{\oV_i,\cE_i,s_i}$ are as in Examples\/
{\rm\ref{ds14ex1}--\ref{ds14ex2}}.\I{d-orbifold with
corners!standard model a@standard model $\eS_{\oV,\cE,s}$} This
$\eX$ is unique up to equivalence in~$\dOrbc$.

Suppose also that\/ $\oY$ is an effective orbifold with corners,
and\/ $g_i:\oV_i\ra\oY$ are submersions for all\/ $i\in I,$ and
there are\/ $2$-morphisms $\ze_{ij}:g_j\ci e_{ij}\Ra
g_i\vert_{\cV_{ij}}$ in $\Orbc$ for all\/ $i<j$ in $I$. Then there
exist a $1$-morphism $\bs h:\eX\ra\eY$ in $\dOrbc$ unique up to
$2$-isomorphism, where $\eY=F_\Orbc^\dOrbc(\oY)=\eS_{\oY,0,0},$
and\/ $2$-morphisms $\bs\ze_i:\bs h\ci\bs\psi_i\Ra\eS_{g_i,0}$ for
all\/~$i\in I$.
\label{ds14thm4}
\end{thm}

\begin{thm} Suppose we are given the following data:
\begin{itemize}
\setlength{\itemsep}{0pt}
\setlength{\parsep}{0pt}
\item[{\rm(a)}] an integer $n;$
\item[{\rm(b)}] a Hausdorff, second countable topological space $X;$
\item[{\rm(c)}] an indexing set\/ $I,$ and a total order $<$ on $I;$
\item[{\rm(d)}] for each\/ $i$ in $I,$ a manifold with corners\/
$V_i,$ a vector bundle $E_i\ra V_i$ with\/ $\dim V_i-\rank
E_i=n,$ a finite group\/ $\Ga_i,$ smooth, locally effective
actions\/ $r_i(\ga):V_i\ra V_i,$ $\hat r_i(\ga):E_i\ra
r(\ga)^*(E_i)$ of\/ $\Ga_i$ on $V_i,E_i$ for $\ga\in\Ga_i,$ a
smooth, $\Ga_i$-equivariant section $s_i:V_i\ra E_i,$ and a
homeomorphism $\psi_i:X_i\ra\hat X_i,$ where $X_i=\{v_i\in
V_i:s_i(v_i)=0\}/\Ga_i$ and\/ $\hat X_i\subseteq X$ is an open
set; and
\item[{\rm(e)}] for all\/ $i<j$ in $I,$ an open submanifold\/
$V_{ij}\subseteq V_i,$ invariant under $\Ga_i,$ a group morphism
$\rho_{ij}:\Ga_i\ra\Ga_j$, a simple, flat, smooth map
$e_{ij}:V_{ij}\ra V_j,$ and a morphism of vector bundles $\hat
e_{ij}:E_i\vert_{V_{ij}}\ra e_{ij}^*(E_j)$.
\end{itemize}

Let this data satisfy the conditions:
\begin{itemize}
\setlength{\itemsep}{0pt}
\setlength{\parsep}{0pt}
\item[{\rm(i)}] $X=\bigcup_{i\in I}\hat X_i;$
\item[{\rm(ii)}] if\/ $i<j$ in $I$ then $\hat e_{ij}\ci s_i
\vert_{V_{ij}}= e_{ij}^*(s_j)+O(s_i^2),$ and for all\/
$\ga\in\Ga$ we have
\begin{align*}
e_{ij}\ci r_i(\ga)&= r_j(\rho_{ij}(\ga))\ci e_{ij}:V_{ij}
\longra V_j,\\
r_i(\ga)^*(\hat e_{ij})\ci\hat r_i(\ga)&= e_{ij}^*(\hat
r_j(\rho_{ij}(\ga)))\ci\hat e_{ij}:E_i\vert_{V_{ij}}
\longra(e_{ij}\ci r_i(\ga))^*(E_j),
\end{align*}
and\/ $\psi_i(X_i\cap (V_{ij}/\Ga_i))=\hat X_i\cap\hat X_j,$
and\/ $\psi_i\vert_{X_i\cap V_{ij}/\Ga_i}=\psi_j\ci
(e_{ij})_*\vert_{X_i\cap V_{ij}/\Ga_j},$ and if\/ $v_i\in
V_{ij}$ with\/ $s_i(v_i)=0$ and\/ $v_j=e_{ij}(v_i)$ then\/
$\rho\vert_{\Stab_{\Ga_i}(v_i)}:\Stab_{\Ga_i}(v_i)\ra
\Stab_{\Ga_j}(v_j)$ is an isomorphism, and the following
sequence of vector spaces is exact:
\begin{equation*}
\smash{\xymatrix@C=19pt{ 0 \ar[r] & T_{v_i}V_i \ar[rrr]^(0.42){\d
s_i(v_i)\op \,\d e_{ij}(v_i)} &&& E_i\vert_{v_i}\!\op\! T_{v_j}V_j
\ar[rrr]^(0.57){\hat e_{ij}(v_i)\op\, -\d s_j(v_j)} &&&
E_j\vert_{v_j} \ar[r] & 0;}}
\end{equation*}
\item[{\rm(iii)}] if\/ $i<j<k$ in $I$ then there exists\/
$\ga_{ijk}\in\Ga_k$ satisfying
\begin{align*}
\rho_{ik}(\ga)&=
\ga_{ijk}\,\ab\rho_{jk}(\rho_{ij}(\ga))\ab\,\ga_{ijk}^{-1}
\quad\text{for all\/ $\ga\in\Ga_i,$}\\
e_{ik}\vert_{V_{ik}\cap  e_{ij}^{-1}(V_{jk})}&= r_k(\ga_{ijk})\ci e_{jk}\ci
e_{ij}\vert_{V_{ik}\cap  e_{ij}^{-1}(V_{jk})},\quad\text{and}\\
\hat e_{ik}\vert_{V_{ik}\cap  e_{ij}^{-1}(V_{jk})}
&= \bigl(e_{ij}^*(e_{jk}^*(\hat r_k(\ga_{ijk})))\ci e_{ij}^*
(\hat e_{jk})\ci\hat e_{ij}\bigr)\vert_{V_{ik}\cap  e_{ij}^{-1}(V_{jk})}.
\end{align*}
\end{itemize}

Then there exist a d-orbifold with corners\/ $\eX$ with\/
$\vdim\eX=n$ and underlying topological space $\cX_\top\cong X,$ and
a $1$-morphism $\bs\psi_i:[\rS_{V_i,E_i,s_i}/\Ga_i]\ra\eX$ in
$\dOrbc$ with underlying continuous map $\psi_i$ which is an
equivalence with the open d-suborbifold\/ $\ehX_i\subseteq\eX$
corresponding to $\hat X_i\subseteq X$ for all\/ $i\in I,$ such that
for all\/ $i<j$ in $I$ there exists a $2$-morphism\/
$\bs\eta_{ij}:\bs\psi_j\ci[\rS_{e_{ij},\hat
e_{ij}},\rho_{ij}]\Ra\bs\psi_i\ci[\bs i_{V_{ij},V_i},\id_{\Ga_i}],$
where $[\rS_{e_{ij},\hat e_{ij}},\rho_{ij}]:
[\rS_{V_{ij},E_i\vert_{V_{ij}},s_i\vert_{V_{ij}}}/\Ga_i]\ra
[\rS_{V_j,E_j,s_j}/\Ga_j]$ and\/ $[\bs i_{V_{ij},V_i},
\id_{\Ga_i}]:[\rS_{V_{ij}, E_i\vert_{V_{ij}},s_i\vert_{V_{ij}}}
/\Ga_i]\ra[\rS_{V_i,E_i,s_i}/\Ga_j]$ combine the notation of
Examples\/ {\rm\ref{ds7ex2}--\ref{ds7ex3}} and\/
{\rm\S\ref{ds132}}.\I{d-orbifold with corners!standard model
b@standard model $[\rS_{V,E,s}/\Ga]$} This $\eX$ is unique up to
equivalence in~$\dOrbc$.

Suppose also that\/ $Y$ is a manifold with corners, and\/
$g_i:V_i\ra Y$ are smooth maps for all\/ $i\in I$ with\/ $g_i\ci
r_i(\ga)= g_i$ for all\/ $\ga\in\Ga_i,$ and\/ $g_j\ci
e_{ij}=g_i\vert_{V_{ij}}$ for all\/ $i<j$ in $I$. Then there exist a
$1$-morphism $\bs h:\eX\ra\eY$ in $\dOrbc$ unique up to
$2$-isomorphism, where $\eY=F_\Manc^\dOrbc(Y)=[\rS_{Y,0,0}/\{1\}],$
and\/ $2$-morphisms $\bs\ze_i:\bs
h\ci\bs\psi_i\Ra[\rS_{g_i,0},\pi_{\{1\}}]$ for all\/ $i\in I$. Here
$[\rS_{g_i,0},\pi_{\{1\}}]:[\rS_{V_i,E_i,s_i}/\Ga_i]\ra
[\rS_{Y,0,0}/\{1\}]=\eY$ with\/ $\hat g_i=0$
and\/~$\rho=\pi_{\{1\}}:\Ga_i\ra\{1\}$.
\label{ds14thm5}
\end{thm}

We can use Theorems \ref{ds14thm4} and \ref{ds14thm5} to prove the
existence of d-orbifold with corners structures on spaces coming
from other areas of geometry, such as moduli spaces of
$J$-holomorphic curves.\I{moduli space!of J-holomorphic curves@of
$J$-holomorphic curves}\I{d-orbifold with
corners!equivalence|)}\I{d-orbifold with corners!gluing by
equivalences|)}

\subsection{Submersions, immersions and embeddings}
\label{ds144}
\I{d-orbifold with corners!w-submersion|(}\I{d-orbifold with
corners!sw-submersion|(}\I{d-orbifold with
corners!submersion|(}\I{d-orbifold with
corners!s-submersion|(}\I{d-orbifold with
corners!w-immersion|(}\I{d-orbifold with
corners!sw-immersion|(}\I{d-orbifold with
corners!sfw-immersion|(}\I{d-orbifold with
corners!immersion|(}\I{d-orbifold with
corners!s-immersion|(}\I{d-orbifold with
corners!sf-immersion|(}\I{d-orbifold with
corners!w-embedding|(}\I{d-orbifold with
corners!sw-embedding|(}\I{d-orbifold with
corners!sfw-embedding|(}\I{d-orbifold with
corners!embedding|(}\I{d-orbifold with
corners!s-embedding|(}\I{d-orbifold with corners!sf-embedding|(}

In \cite[\S 12.4]{Joyc6} we extend \S\ref{ds74} and \S\ref{ds114} to
d-orbifolds with corners.

\begin{dfn} Let $\bs f:\eX\ra\eY$ be a 1-morphism in $\dOrbc$. Then
$T^*\eX$ and $f^*(T^*\eY)$ are virtual vector bundles on $\cX$ of
ranks $\vdim\eX,\vdim\eY$, and $\Om_{\bs f}:f^*(T^*\eY)\ra T^*\eX$
is a 1-morphism in $\vvect(\cX)$, as in Remark \ref{ds11rem} and
Definition \ref{ds14def1}. `Weakly injective', \ldots, below are as
in Definition~\ref{ds11def2}.
\begin{itemize}
\setlength{\itemsep}{0pt}
\setlength{\parsep}{0pt}
\item[(a)] We call $\bs f$ a {\it w-submersion\/} if $\bs f$ is
semisimple and flat and $\Om_{\bs f}$ is weakly injective. We
call $\bs f$ an {\it sw-submersion\/} if it is also simple.
\item[(b)] We call $\bs f$ a {\it submersion\/} if $\bs f$ is
semisimple and flat and $\Om_{C(\bs f)}$ is injective, for
$C(\bs f)$ as in \S\ref{ds135}. We call $\bs f$ an {\it
s-submersion\/} if it is also simple.
\item[(c)] We call $\bs f$ a {\it w-immersion\/} if
$f:\cX\!\ra\!\cY$ is representable and $\Om_{\bs f}$ is weakly
surjective. We call $\bs f$ an {\it sw-immersion}, or {\it
sfw-immersion}, if $\bs f$ is also simple, or simple and flat.
\item[(d)] We call $\bs f$ an {\it immersion\/} if
$f:\cX\ra\cY$ is representable and $\Om_{\smash{\hat C(\bs f)}}$
is surjective, for $\hat C(\bs f)$ as in \S\ref{ds135}. We call
$\bs f$ an {\it s-immersion\/} if $\bs f$ is also simple, and an
{\it sf-immersion\/} if $\bs f$ is also simple and flat.
\item[(e)] We call $\bs f$ a {\it w-embedding, sw-embedding,
sfw-embedding, embedding, s-embedding}, or {\it sf-embedding},
if $\bs f$ is a w-immersion, \ldots, sf-immersion, respectively,
and $f_*:\Iso_\cX([x])\ra\Iso_\cY(f_\top([x]))$ is an
isomorphism for all $[x]\in\cX_\top$, and $f_\top:\cX_\top
\ra\cY_\top$ is a homeomorphism with its image, so in particular
$f_\top$ is injective.
\end{itemize}

Parts (c)--(e) enable us to define {\it
d-suborbifolds\/}\I{d-orbifold with corners!d-suborbifold} $\eX$ of
a d-orbifold with corners $\eY$. {\it Open d-suborbifolds\/} are
open d-substacks $\eX$ in $\eY$. For more general d-suborbifolds, we
call $\bs f:\eX\ra\eY$ a {\it w-immersed, sw-immersed, sfw-immersed,
immersed, s-immersed, sf-immersed, w-embedded, sw-embedded,
sfw-embedded, embedded, s-embedded}, or {\it sf-embedded
suborbifold\/} of $\eY$ if $\eX,\eY$ are d-orbifolds with corners
and $\bs f$ is a w-immersion, \ldots, sf-embedding, respectively.
\label{ds14def3}
\end{dfn}

Theorem \ref{ds7thm6} in \S\ref{ds74} holds with orbifolds and
d-orbifolds with corners in place of manifolds and d-manifolds with
corners, except part (v), when we need also to assume $f:\cX\ra\cY$
representable\I{Deligne--Mumford $C^\iy$-stack!representable
1-morphism} to deduce $\bs f$ is \'etale,\I{d-orbifold with
corners!etale 1-morphism@\'etale 1-morphism} and part (x), which is
false for d-orbifolds with corners (in the Zariski
topology,\I{Zariski topology} at least).\I{d-orbifold with
corners!w-submersion|)}\I{d-orbifold with
corners!sw-submersion|)}\I{d-orbifold with
corners!submersion|)}\I{d-orbifold with
corners!s-submersion|)}\I{d-orbifold with
corners!w-immersion|)}\I{d-orbifold with
corners!sw-immersion|)}\I{d-orbifold with
corners!sfw-immersion|)}\I{d-orbifold with
corners!immersion|)}\I{d-orbifold with
corners!s-immersion|)}\I{d-orbifold with
corners!sf-immersion|)}\I{d-orbifold with
corners!w-embedding|)}\I{d-orbifold with
corners!sw-embedding|)}\I{d-orbifold with
corners!sfw-embedding|)}\I{d-orbifold with
corners!embedding|)}\I{d-orbifold with
corners!s-embedding|)}\I{d-orbifold with corners!sf-embedding|)}

\subsection{Bd-transversality and fibre products}
\label{ds145}
\I{d-orbifold with corners!fibre products|(}\I{d-orbifold with
corners!bd-transverse 1-morphisms|(}\I{d-orbifold with
corners!cd-transverse
1-morphisms|(}\I{bd-transversality|(}\I{cd-transversality|(}\I{2-category!fibre
products in|(}

In \cite[\S 12.5]{Joyc6} we generalize \S\ref{ds75} and
\S\ref{ds115} to $\dOrbc$. Here are the analogues of Definition
\ref{ds7def4} and Theorems \ref{ds7thm7}--\ref{ds7thm10}.

\begin{dfn} Let $\eX,\eY,\eZ$ be d-orbifolds with corners and
$\bs g:\eX\ra\eZ,$ $\bs h:\eY\ra\eZ$ be 1-morphisms. We call $\bs
g,\bs h$ {\it bd-transverse\/} if they are both b-transverse in
$\dStac$ in the sense of Definition \ref{ds13def3}, and d-transverse
in the sense of Definition \ref{ds11def3}. We call $\bs g,\bs h$
{\it cd-transverse\/} if they are both c-transverse in $\dStac$ in
the sense of Definition \ref{ds13def3}, and d-transverse. As in
\S\ref{ds136}, c-transverse implies b-transverse, so cd-transverse
implies bd-transverse.
\label{ds14def4}
\end{dfn}

\begin{thm} Suppose\/ $\eX,\eY,\eZ$ are d-orbifolds with corners
and\/ $\bs g:\eX\ra\eZ,$ $\bs h:\eY\ra\eZ$ are bd-transverse
$1$-morphisms, and let\/ $\eW=\eX\t_{\bs g,\eZ,\bs h}\eY$ be the
fibre product in $\dStac,$ which exists by Theorem\/
{\rm\ref{ds13thm5}(a)} as $\bs g,\bs h$ are b-transverse. Then $\eW$
is a d-orbifold with corners, with\I{d-orbifold with corners!fibre
products!bd-transverse}\I{fibre product!of d-orbifolds with
corners}\I{d-orbifold with corners!virtual dimension}
\e
\vdim\eW=\vdim\eX+\vdim\eY-\vdim\eZ.
\label{ds14eq1}
\e
Hence, all bd-transverse fibre products exist
in\/~$\dOrbc$.\I{2-category!fibre products in|)}
\label{ds14thm6}
\end{thm}

\begin{thm} Suppose\/ $\bs g:\eX\ra\eZ$ and\/ $\bs h:\eY\ra\eZ$ are
$1$-morphisms in $\dOrbc$. The following are sufficient conditions
for $\bs g,\bs h$ to be cd-transverse, and hence bd-transverse, so
that\/ $\eW=\eX\t_{\bs g,\eZ,\bs h}\eY$ is a d-orbifold with corners
of virtual dimension\/ {\rm\eq{ds14eq1}:}
\begin{itemize}
\setlength{\itemsep}{0pt}
\setlength{\parsep}{0pt}
\item[{\bf(a)}] $\eZ$ is an orbifold without boundary, that is,
$\eZ\in\bOrb;$ or
\item[{\bf(b)}] $\bs g$ or $\bs h$ is a w-submersion.\I{d-orbifold with
corners!w-submersion}\I{d-orbifold with corners!bd-transverse
1-morphisms|)}\I{d-orbifold with corners!cd-transverse
1-morphisms|)}\I{bd-transversality|)}\I{cd-transversality|)}
\end{itemize}
\label{ds14thm7}
\end{thm}

\begin{thm} Let\/ $\eX,\eY,\eZ$ be d-orbifolds with corners with\/
$\eY$ an orbifold, and\/ $\bs g:\eX\ra\eZ,$ $\bs h:\eY\ra\eZ$ be
$1$-morphisms with\/ $\bs g$ a submersion.\I{d-orbifold with
corners!submersion} Then\/ $\eW=\eX\t_{\bs g,\eZ,\bs h}\eY$ is an
orbifold,\I{d-orbifold with corners!is an orbifold}
with\/~$\dim\eW=\vdim\eX+\dim\eY-\vdim\eZ$.
\label{ds14thm8}
\end{thm}

\begin{thm}{\bf(i)} Let\/ $\rho:G\ra H$ be an injective morphism
of finite groups, and\/ $H$ act linearly on $\R^n$ preserving\/
$[0,\iy)^k\t\R^{n-k}$. Then\/ {\rm\S\ref{ds132}} gives a quotient\/
$1$-morphism $[\bs 0,\rho]:[\bs */G]\ra\bigl[\bs{[0,\iy)^k}
\t\bR^{\bs{n-k}}/H\bigr]$ in $\dOrbc$. Suppose $\eX$ is a d-orbifold
with corners and\/ $\bs g:\eX\ra\bigl[\bs{[0,\iy)^k}
\t\bR^{\bs{n-k}}/H\bigr]$ is a semisimple, flat\/ $1$-morphism in
$\dOrbc$. Then the fibre product\/ $\eW=\eX\t_{\bs
g,[\bs{[0,\iy)^k}\t\bR^{\bs{n-k}}/H],[\bs 0,\rho]}[\bs */G]$ exists
in $\dOrbc$. The projection $\bs\pi_\eX:\eW\ra\eX$ is an
s-immersion, and an s-embedding if\/ $\rho$ is an isomorphism.

When\/ $k=0,$ any $1$-morphism $\bs g:\eX\ra[\bR^{\bs n}/H]$ is
semisimple and flat, and\/ $\bs\pi_\eX:\eW\!\ra\!\eX$ is an
sf-immersion, and an sf-embedding if\/ $\rho$ is an
isomorphism.\I{d-orbifold with corners!semisimple
1-morphism}\I{d-orbifold with corners!flat 1-morphism}\I{d-orbifold
with corners!s-immersion|(}\I{d-orbifold with
corners!s-embedding|(}\I{d-orbifold with corners!sf-embedding|(}

\smallskip

\noindent{\bf(ii)} Suppose\/ $\bs f:\eX\ra\eY$ is an s-immersion of
d-orbifolds with corners, and\/ $[x]\in\cX_\top$ with\/
$f_\top([x])=[y]\in\cY_\top$. Write $\rho:G\ra H$ for
$f_*:\Iso_\cX([x])\ra\Iso_\cY([y])$. Then $\rho$ is injective, and
there exist open neighbourhoods $\eU\subseteq\eX$ and\/
$\eV\subseteq\eY$ of\/ $[x],[y]$ with\/ $\bs f(\eU)\subseteq\eV,$ a
linear action of\/ $H$ on $\R^n$ preserving\/ $[0,\iy)^k\t\R^{n-k}
\subseteq\R^n$ where $n=\vdim\eY-\vdim\eX\ge 0$ and\/ $0\le k\le n,$
and a $1$-morphism $\bs g:\eV\ra\bigl[\bs{[0,\iy)^k}\t
\bR^{\bs{n-k}}/H\bigr]$ with\/ $g_\top([y])=[0],$ fitting into a
$2$-Cartesian\I{2-category!2-Cartesian square} square in~$\dOrbc\!:$
\begin{equation*}
\xymatrix@C=80pt@R=10pt{*+[r]{\eU} \ar[d]^{\bs f\vert_\eU} \ar[r]
\drtwocell_{}\omit^{}\omit{^{}} & *+[l]{[\bs */G]} \ar[d]_{[\bs 0,\rho]} \\
*+[r]{\eV} \ar[r]^(0.3){\bs g} &
*+[l]{\bigl[\bs{[0,\iy)^k}\t\bR^{\bs{n-k}}/H\bigr].} }
\end{equation*}
If\/ $\bs f$ is an sf-immersion then $k=0$. If\/ $\bs f$ is an s- or
sf-embedding then $\rho$ is an isomorphism, and we may
take\/~$\eU=\bs f^{-1}(\eV)$.\I{d-orbifold with
corners!s-immersion|)}\I{d-orbifold with
corners!s-embedding|)}\I{d-orbifold with
corners!sf-embedding|)}\I{d-orbifold with corners!fibre products|)}
\label{ds14thm9}
\end{thm}

\subsection{Embedding d-orbifolds with corners into orbifolds}
\label{ds146}
\I{d-orbifold with corners!embedding!into orbifolds|(}

Section \ref{ds47} discussed embeddings of d-manifolds $\bX$ into
manifolds $Y$. Our two main results were Theorem \ref{ds4thm12},
which gave necessary and sufficient conditions on $\bX$ for
existence of embeddings $\bs f:\bX\hookra\bR^{\bs n}$ for $n\gg 0$,
and Theorem \ref{ds4thm13}, which showed that if an embedding $\bs
f:\bX\hookra\bY$ exists with $\bX$ a d-manifold and
$\bY=F_\Man^\dMan(Y)$, then $\bX\simeq\bS_{V,E,s}$ for open
$V\subseteq Y$.

Section \ref{ds76} generalized \S\ref{ds47} to d-manifolds with
corners, requiring $\bs f:\rX\hookra\rY$ to be an sf-embedding for
the analogue of Theorem \ref{ds4thm13}. Section \ref{ds116}
explained that while Theorem \ref{ds4thm13} generalizes to
d-orbifolds, we do not have a good d-orbifold generalization of
Theorem \ref{ds4thm12}. Thus, we do not have a useful necessary and
sufficient criterion for when a d-orbifold is principal.

As in \cite[\S 12.6]{Joyc6}, the situation is the same for
d-orbifolds with corners as for d-orbifolds. Here is the analogue of
Theorem~\ref{ds4thm13}:\I{d-orbifold with
corners!sf-embedding|(}\I{d-orbifold with
corners!principal|(}\I{principal d-orbifold with corners|(}
\I{d-orbifold with corners!standard model a@standard model
$\eS_{\oV,\cE,s}$|(}

\begin{thm} Suppose $\eX$ is a d-orbifold with corners, $\oY$ an
orbifold with corners, and\/ $\bs f:\eX\ra\eY$ an sf-embedding, in
the sense of\/ {\rm\S\ref{ds144}}. Then there exist an open
suborbifold\/ $\oV\subseteq\oY$ with $\bs f(\eX)\subseteq\eV,$ a
vector bundle $\cE$ on $\oV,$ and a section\/ $s\in C^\iy(\cE)$
fitting into a $2$-Cartesian\I{2-category!2-Cartesian square}
diagram in\/~$\dOrbc\!:$
\begin{equation*}
\xymatrix@C=60pt@R=10pt{ \eX \ar[r]_(0.25){\bs f} \ar[d]^{\bs f}
\drtwocell_{}\omit^{}\omit{^{}}
 & \eV \ar[d]_{\bs 0} \\ \eV \ar[r]^(0.7){\bs s} & \eE,}
\end{equation*}
where\/ $\eY,\eV,\eE,\bs s,\bs 0=F_\Orbc^\dOrbc\bigl(\oY,
\oV,\ab\Totc(\cE),\ab\Totc(s),\ab\Totc(0)\bigr),$\I{orbifold with
corners!vector bundles on!total space functor $\Totc$} in the
notation of\/ {\rm\S\ref{ds121}}. Thus $\eX$ is equivalent to the
`standard model'\/ $\eS_{\oV,\cE,s}$ of Example\/
{\rm\ref{ds14ex1},} and is a principal d-orbifold with corners.
\label{ds14thm10}
\end{thm}

Again, in contrast to d-manifolds with corners, the author does not
know useful necessary and sufficient conditions for when a
d-orbifold with corners admits an sf-embedding into an orbifold, or
is a principal d-orbifold with \I{d-orbifold with
corners!sf-embedding|)}\I{d-orbifold with
corners!principal|)}\I{principal d-orbifold with corners|)}
\I{d-orbifold with corners!standard model a@standard model
$\eS_{\oV,\cE,s}$|)}\I{d-orbifold with corners!embedding!into
orbifolds|)}corners.

\subsection{Orientations on d-orbifolds with corners}
\label{ds147}
\I{d-orbifold with corners!orientations|(}

Section \ref{ds48} discussed orientations on d-manifolds. This was
extended to d-manifolds with corners in \S\ref{ds77}, and to
d-orbifolds in \S\ref{ds117}. As in \cite[\S 12.7]{Joyc6}, all this
generalizes easily to d-orbifolds with corners, so we will give few
details.

If $\eX$ is a d-orbifold with corners, the virtual cotangent
bundle\I{d-orbifold with corners!virtual cotangent bundle}\I{virtual
cotangent bundle} $T^*\eX$ is a virtual vector bundle\I{virtual
vector bundle} on $\cX$. We define an {\it orientation\/} $\om$ on
$\eX$ to be an orientation on the orientation line
bundle\I{orientation line bundle}
$\cL_{T^*\eX}$.\G[LT*Xd]{$\cL_{T^*\eX}$}{orientation line bundle of
a d-orbifold with corners $\eX$}\I{d-orbifold with
corners!orientation line bundle} The analogues of Example
\ref{ds4ex6}, Theorem \ref{ds4thm15}, Proposition \ref{ds4prop4},
Theorem \ref{ds7thm15}, and Propositions \ref{ds7prop1} and
\ref{ds7prop2} all hold for d-orbifolds with corners.\I{d-orbifold
with corners!orientations|)}

\subsection{Orbifold strata of d-orbifolds with corners}
\label{ds148}
\I{orbifold strata!of d-orbifolds with corners|(}\I{d-orbifold with
corners!orbifold strata|(}

Sections \ref{ds87}, \ref{ds92}, \ref{ds118}, \ref{ds125}, and
\ref{ds137} discussed orbifold strata for Deligne--Mumford
$C^\iy$-stacks, orbifolds, d-orbifolds, orbifolds with corners, and
d-stacks with corners. In \cite[\S 12.8]{Joyc6} we extend this to
d-orbifolds with corners.

Let $\eX$ be a d-orbifold with corners and $\Ga$ a finite group, so
that \S\ref{ds137} gives orbifold strata
$\eX^\Ga,\etX^\Ga,\ehX^\Ga,\eX^\Ga_\ci,\etX^\Ga_\ci,\ehX^\Ga_\ci$,
which are d-stacks with corners. Use the notation
$\La^\Ga,\La^\Ga/\Aut(\Ga)$ of Definition \ref{ds9def4}. Exactly as
in the d-orbifold case in \S\ref{ds118}, there are natural
decompositions
\begin{align*}
\eX^\Ga&=\ts\coprod_{\la\in\La^\Ga}\eX^{\Ga,\la}, &
\etX^\Ga&=\ts\coprod_{\mu\in\La^\Ga/\Aut(\Ga)}\etX^{\Ga,\mu},&
\ehX^\Ga&=\ts\coprod_{\mu\in\La^\Ga/\Aut(\Ga)}\ehX^{\Ga,\mu},\\
\eX^\Ga_\ci&=\ts\coprod_{\smash{\la\in\La^\Ga}}\eX^{\Ga,\la}_\ci, &
\etX^\Ga_\ci&=
\ts\coprod_{\mu\in\La^\Ga/\Aut(\Ga)} \etX^{\Ga,\mu}_\ci, &
\ehX^\Ga_\ci&=\ts\coprod_{\mu\in\La^\Ga/\Aut(\Ga)} \ehX^{\Ga,\mu}_\ci,
\end{align*}
where
$\eX^{\Ga,\la},\ldots,\etX^{\Ga,\mu}_\ci$\G[XGaj]{$\eX^{\Ga,\la},
\etX{}^{\Ga,\mu},\ehX{}^{\Ga,\mu},\eX{}^{\Ga,\la}_\ci,
\etX{}^{\Ga,\mu}_\ci,\ehX{}^{\Ga,\mu}_\ci$}{orbifold strata of a
d-orbifold with corners $\eX$} are d-orbifolds with corners with
$\vdim\eX^{\Ga,\la}=\vdim\eX^{\Ga,\la}_\ci\ab= \vdim\eX-\dim\la$ and
$\vdim\etX^{\Ga,\mu}=\cdots=\vdim\ehX^{\Ga,\mu}_\ci=\vdim\eX-
\dim\mu$.\I{d-orbifold with corners!virtual dimension}

The analogue of Proposition \ref{ds11prop2} on orientations of
orbifold strata $\bcX^{\Ga,\la},\ab\ldots,\ab\bhcX^{\Ga,\mu}_\ci$
for oriented d-orbifolds $\bcX$ also holds for d-orbifolds with
corners.\I{d-orbifold with corners!orbifold strata!orientations on}
In an analogue of Corollary \ref{ds12cor}, we can relate boundaries
of orbifold strata to orbifold strata of boundaries:\I{d-orbifold
with corners!orbifold strata!boundaries of}

\begin{prop} Let\/ $\eX$ be a d-orbifold with corners, and\/ $\Ga$ a
finite group. Then Corollary\/ {\rm\ref{ds13cor}} gives\/
$1$-morphisms $\bs J{}^\Ga(\eX):(\pd\eX)^\Ga\ra \pd(\eX^\Ga),$
$\bs{\ti J}{}^\Ga(\eX): {{}\,\,\,\,\widetilde{\!\!\!\!(\pd\eX)
\!\!\!\!}\,\,\,\,^\Ga}\ra \pd(\etX{}^\Ga),$ $\bs{\hat
J}{}^\Ga(\eX):\,\,\,\,\widehat{\!\!\!\!
(\pd\eX)\!\!\!\!}\,\,\,\,^\Ga\ra\pd(\ehX{}^\Ga)$ in\/ $\dcOrbc,$
which are equivalences with open and closed subobjects
in\/~$\pd(\eX^\Ga),\pd(\etX{}^\Ga),\pd(\ehX{}^\Ga)$.

These restrict to $1$-morphisms $\bs J^{\Ga,\la}(\eX):
(\pd\eX)^{\Ga,\la} \ra \pd(\eX^{\Ga,\la})$ in\/ $\dOrbc$ for
$\la\in\La^\Ga$ and\/ $\bs{\ti J}{}^{\Ga,\mu}(\eX):
\,\,\,\,\widetilde{\!\!\!\!(\pd\eX)\!\!\!\!}\,\,\,\,^{\Ga,\mu}\ra
\pd(\etX{}^{\Ga,\mu}),$ $\bs{\hat J}{}^{\Ga,\mu}(\eX):
\,\,\,\,\widehat{\!\!\!\! (\pd\eX)\!\!\!\!}\,\,\,\,^{\Ga,\mu}\ra
\pd(\ehX{}^{\Ga,\mu})$ for $\mu\in\La^\Ga/\Aut(\La),$ which are
equivalences with open and closed d-suborbifolds. Hence, if\/
$\eX^{\Ga,\la}=\bs\es$ then $(\pd\eX)^{\Ga,\la}=\bs\es,$ and
similarly for\/~$\etX{}^{\Ga,\mu},\ehX{}^{\Ga,\mu}$.

Now suppose $\eX$ is straight,\I{d-orbifold with
corners!straight}\I{d-stack with corners!straight} in the sense of\/
{\rm\S\ref{ds137},} for instance $\eX$ could be a d-orbifold with
boundary.\I{d-orbifold with boundary} Then $\bs
J{}^\Ga(\eX),\ab\ldots,\ab\bs{\hat J}{}^{\Ga,\mu}(\eX)$ are
equivalences, so that\/ $(\pd\eX)^\Ga\simeq \pd(\eX^\Ga),$
$(\pd\eX)^{\Ga,\la}\simeq \pd(\eX^{\Ga,\la}),$ and so on.\I{orbifold
strata!of d-orbifolds with corners|)}\I{d-orbifold with
corners!orbifold strata|)}
\label{ds14prop2}
\end{prop}

\subsection[Kuranishi neighbourhoods and good coordinate
systems]{Kuranishi neighbourhoods, good coordinate systems}
\label{ds149}
\I{d-orbifold with corners!Kuranishi neighbourhood|(}\I{d-orbifold
with corners!good coordinate system|(}\I{d-orbifold with
corners!standard model b@standard model
$[\rS_{V,E,s}/\Ga]$|(}\I{good coordinate system|(}

Section \ref{ds119} defined type A Kuranishi neighbourhoods,
coordinate changes, and good coordinate systems, on d-orbifolds. In
\cite[\S 12.9]{Joyc6} we generalize these to d-orbifolds with
corners. The definitions in the corners case are obtained by
replacing $\Man,\ab\Orb,\ab\dMan,\ab\dOrb$ by
$\Manc,\ab\Orbc,\ab\dManc, \ab\dOrbc$ throughout, and making a few
other easy changes such as taking the $e_{ij}$ to be sf-embeddings
in Definitions \ref{ds11def6}(c). For brevity we will not write the
definitions out again, but just indicate the differences.

\begin{dfn} Let $\eX$ be a d-orbifold with corners. Define a {\it
type A Kuranishi neighbourhood\/} $(V,E,\Ga,s,\bs\psi)$ on $\eX$
following Definition \ref{ds11def5}, but taking $V$ to be a manifold
with corners, and defining the principal d-orbifold with corners
$[\rS_{V,E,s}/\Ga]$ by combining Example \ref{ds7ex2} and
\S\ref{ds132}, as in~\S\ref{ds142}.

If $(V_i,E_i,\Ga_i,s_i,\bs\psi_i),(V_j,E_j,\Ga_j, s_j,\bs\psi_j)$
are type A Kuranishi neighbourhoods on $\eX$ with
$\bs\es\ne\bs\psi_i([\rS_{V_i,E_i,s_i}/\Ga_i])\cap
\bs\psi_j([\rS_{V_j,E_j,s_j}/\Ga_j])\subseteq\eX$, define a {\it
type A coordinate change\I{d-orbifold with corners!Kuranishi
neighbourhood!coordinate change} $(V_{ij},e_{ij},\hat
e_{ij},\rho_{ij}, \bs\eta_{ij})$ from\/
$(V_i,E_i,\Ga_i,s_i,\bs\psi_i)$ to\/}
$(V_j,E_j,\Ga_j,s_j,\bs\psi_j)$ following Definition \ref{ds11def6},
but taking $e_{ij}:V_{ij}\ra V_j$ to be an sf-embedding of manifolds
with corners,\I{manifold with corners!sf-embedding} and defining the
quotient 1-morphism $[\rS_{\smash{e_{ij},\hat e_{ij}}},\rho_{ij}]$
by combining Example \ref{ds7ex3} and \S\ref{ds132}, as
in~\S\ref{ds142}.

Define a {\it type A good coordinate system\/} on $\eX$ following
Definition \ref{ds11def7}, defining quotient 2-morphisms
$\bs\eta_{ijk}=[0,\ga_{ijk}]$ in (d) using \S\ref{ds132}. Let $Y$ be
a manifold with corners, and $\bs h:\eX\ra\eY$ a 1-morphism in
$\dOrbc$, where $\eY=F_\Manc^\dOrbc(Y)$. Define a {\it type A good
coordinate system for\/} $\bs h:\eX\ra\eY$ following
Definition~\ref{ds11def7}.
\label{ds14def5}
\end{dfn}

Here is the analogue of Theorem \ref{ds11thm11}. It is proved
in~\cite[App.~D]{Joyc6}.

\begin{thm} Suppose $\eX$ is a d-orbifold with corners. Then there
exists a type A good coordinate system $\bigl(I,\kern -.1em <,\kern
-.1em (V_i,\kern -.1em E_i,\kern -.1em \Ga_i,\kern -.1em s_i,\kern
-.1em \bs\psi_i),\kern -.1em (V_{ij},\kern -.1em e_{ij},\kern -.1em
\hat e_{ij},\kern -.1em\rho_{ij},\kern -.1em \bs\eta_{ij}),\kern
-.1em \ga_{ijk}\bigr)$ for $\eX,$ in the sense of Definition\/
{\rm\ref{ds14def5}}. If\/ $\eX$ is compact, we may take $I$ to be
finite. If\/ $\{\eU_j:j\in J\}$ is an open cover of\/ $\eX,$ we may
take $\eX_i=\bs\psi_i([\rS_{V_i,E_i,s_i}/\Ga_i]) \subseteq
\eU_{\smash{j_i}}$ for each\/ $i\in I$ and some\/ $j_i\in J$.
If\/ $\eX$ is a d-orbifold with boundary,\I{d-orbifold with
boundary} we may take the $V_i$ to be manifolds with boundary.

Now let\/ $Y$ be a manifold with corners and\/ $\bs h:\eX\ra\eY=
F_\Manc^\dOrbc(Y)$ a semisimple,\I{d-orbifold with
corners!semisimple 1-morphism} flat\/\I{d-orbifold with corners!flat
1-morphism} $1$-morphism in $\dOrbc$. Then all the above extends to
type A good coordinate systems for $\bs h:\eX\ra\eY,$ and we may
take the $g_i:V_i\ra Y$ to be submersions in\/~$\Manc$.
\label{ds14thm11}
\end{thm}

We can regard Theorem \ref{ds14thm11} as a kind of converse to
Theorem \ref{ds14thm5}. Note that we make the extra assumption that
$\bs h$ is semisimple and flat in the last part. This happens
automatically if $Y$ is without boundary. Since
submersions\I{manifold with corners!submersion} in $\Manc$ are
automatically semisimple and flat, $\bs h$ being semisimple and flat
is a necessary condition for the $g_i:V_i\ra Y$ to be submersions.
In \cite[\S 12.9]{Joyc6} we also give `type B' versions of
Definition \ref{ds14def5} and Theorem \ref{ds14thm11} using the
`standard model' d-orbifolds with corners
$\eS_{\oV,\cE,s}$\I{d-orbifold with corners!standard model
a@standard model $\eS_{\oV,\cE,s}$} and 1-morphisms
$\eS_{\smash{e_{ij},\hat e_{ij}}}$ of Examples
\ref{ds14ex1}--\ref{ds14ex2} instead of $[\rS_{V_i,E_i,s_i}/\Ga_i]$
and~$[\rS_{\smash{e_{ij},\hat e_{ij}}},\rho_{ij}]$.\I{d-orbifold
with corners!Kuranishi neighbourhood|)}\I{d-orbifold with
corners!good coordinate system|)}\I{d-orbifold with corners!standard
model b@standard model $[\rS_{V,E,s}/\Ga]$|)}\I{good coordinate
system|)}

\subsection{Semieffective and effective d-orbifolds with corners}
\label{ds1410}
\I{d-orbifold with corners!semieffective|(}\I{d-orbifold with
corners!effective|(}

Section \ref{ds1110} discussed semieffective and effective
d-orbifolds. As in \cite[\S 12.10]{Joyc6}, all this material extends
to d-orbifolds with corners essentially unchanged. We define {\it
semieffective\/} and {\it effective\/} d-orbifolds with corners
$\eX$ following Definition \ref{ds11def8}. The analogues of
Proposition \ref{ds11prop3} and the rest of \S\ref{ds119} then hold,
with (d-)orbifolds replaced by (d-)orbifolds with corners
throughout.

\begin{prop} Let\/ $\eX$ be an effective (or semieffective) d-orbifold
with corners. Then $\pd^k\eX$ is also effective (or semieffective),
for all\/~$k\ge 0$.
\label{ds14prop3}
\end{prop}

However, $\eX$ (semi)effective does not imply $C_k(\eX)$
(semi)effective.\I{d-orbifold with
corners!semieffective|)}\I{d-orbifold with
corners!effective|)}\I{d-orbifold with corners|)}

\section{D-manifold and d-orbifold bordism}
\label{ds15}
\I{bordism|(}

In \cite[Chap.~13]{Joyc6} we discuss {\it bordism groups\/} of
manifolds and orbifolds, defined using manifolds, d-manifolds,
orbifolds, and d-orbifolds. We can use these to prove that compact,
oriented d-manifolds admit virtual cycles, and so can be used in
enumerative invariant problems. The same applies for compact,
oriented d-orbifolds, although the direct proof using bordism no
longer works.

\subsection{Classical bordism groups for manifolds}
\label{ds151}
\I{bordism!classical bordism|(}

In \cite[\S 13.1]{Joyc6} we review background material on bordism
from the literature. Classical bordism groups $MSO_k(Y)$ were
defined by Atiyah \cite{Atiy} for topological spaces $Y$, using
continuous maps $f:X\ra Y$ for $X$ a compact, oriented manifold.
Conner \cite[\S I]{Conn} gives a good introduction. We define
bordism $B_k(Y)$ only for manifolds $Y$, using smooth $f:X\ra Y$,
following Conner's {\it differential bordism groups\/} \cite[\S
I.9]{Conn}. By \cite[Th.~I.9.1]{Conn}, the natural projection
$B_k(Y)\ra MSO_k(Y)$ is an isomorphism, so our notion of bordism
agrees with the usual definition.

\begin{dfn} Let $Y$ be a manifold without boundary, and $k\in\Z$.
Consider pairs $(X,f)$, where $X$ is a compact, oriented manifold
without boundary with $\dim X=k$, and $f:X\ra Y$ is a smooth map.
Define an equivalence relation $\sim$ on such pairs by $(X,f)\sim
(X',f')$ if there exists a compact, oriented \ab $(k+1)$-manifold
with boundary\I{manifold with boundary} $W$, a smooth map $e:W\ra
Y$, and a diffeomorphism of oriented manifolds $j:-X\amalg X'\ra\pd
W$, such that $f\amalg f'=e\ci i_W\ci j$, where $-X$ is $X$ with the
opposite orientation.

Write $[X,f]$ for the $\sim$-equivalence class ({\it bordism
class\/}) of a pair $(X,f)$. For each $k\in\Z$, define the $k^{\it
th}$ {\it bordism group\/} $B_k(Y)$\G[BkYa]{$B_k(Y)$}{classical
bordism group of manifold $Y$} of $Y$ to be the set of all such
bordism classes $[X,f]$ with $\dim X=k$. We give $B_k(Y)$ the
structure of an abelian group, with zero element $0_Y=[\es,\es]$,
and addition given by $[X,f]+[X',f']=[X\amalg X',f\amalg f']$, and
additive inverses $-[X,f]=[-X,f]$.

Define $\Pi_\bo^\hom:B_k(Y)\ra H_k(Y;\Z)$ by
$\Pi_\bo^\hom:[X,f]\mapsto f_*([X])$, where $H_*(-;\Z)$ is singular
homology,\I{homology}\I{bordism!projection to homology} and $[X]\in
H_k(X;\Z)$ is the fundamental class.

If $Y$ is oriented and of dimension $n$, there is a biadditive,
associative, supercommutative {\it intersection
product\/}\I{bordism!classical bordism!intersection product}
$\bu:B_k(Y)\t B_l(Y)\ra B_{k+l-n}(Y)$, such that if $[X,f],[X',f']$
are classes in $B_*(Y)$, with $f,f'$ transverse, then the fibre
product $X\t_{f,Y,f'}X'$ exists as a compact oriented manifold, and
\begin{equation*}
[X,f]\bu[X',f']=[X\t_{f,Y,f'}X',f\ci\pi_X].
\end{equation*}
\label{ds15def1}
\end{dfn}

As in \cite[\S I.5]{Conn}, bordism is a generalized homology
theory.\I{generalized homology theory} Results of Thom, Wall and
others in \cite[\S I.2]{Conn} compute the bordism groups $B_k(*)$ of
the point $*$. This partially determines the bordism groups of
general manifolds $Y$, as there is a spectral sequence\I{spectral
sequence} $H_i\bigl(Y;B_j(*)\bigr)\Ra
B_{i+j}(Y)$.\I{bordism!classical bordism|)}

\subsection{D-manifold bordism groups}
\label{ds152}
\I{bordism!d-manifold bordism|(}\I{d-manifold!bordism|(}

In \cite[\S 13.2]{Joyc6} we define {\it d-manifold bordism\/} by
replacing manifolds $X$ in pairs $[X,f]$ in \S\ref{ds151} by
d-manifolds $\bX$. For simplicity, we identify the 2-category
$\dMan$ of d-manifolds $\bX$ in \S\ref{ds41}, and the 2-subcategory
$\bdMan$ of d-manifolds without boundary $\rX=(\bX,\bs\es,\bs\es,
\bs\es)$ in $\dManc$ in \S\ref{ds71}, writing both as~$\bX$.

\begin{dfn} Let $Y$ be a manifold without boundary, and
$k\in\Z$. Consider pairs $(\bX,\bs f)$, where $\bX\in\dMan$ is a
compact, oriented d-manifold without boundary with $\vdim\bX=k$, and
$\bs f:\bX\ra\bY$ is a 1-morphism in $\dMan$,
where~$\bY=F_\Man^\dMan(Y)$.

Define an equivalence relation $\sim$ between such pairs by
$(\bX,\bs f)\sim (\bX',\bs f')$ if there exists a compact, oriented
d-manifold with boundary\I{d-manifold with boundary} $\rW$ with
$\vdim\rW=k+1$, a 1-morphism $\bs e:\rW\ra\bY$ in $\dManb$, an
equivalence of oriented d-manifolds $\bs j:-\bX\amalg\bX'\ra\pd
\rW$, and a 2-morphism~$\eta:\bs f\amalg\bs f'\Ra\bs e\ci\bs
i_\rW\ci \bs j$.

Write $[\bX,\bs f]$ for the $\sim$-equivalence class ({\it d-bordism
class\/}) of a pair $(\bX,\bs f)$. For each $k\in\Z$, define the
$k^{\it th}$ {\it d-manifold bordism group}, or {\it d-bordism
group}, $dB_k(Y)$\G[dBkYa]{$dB_k(Y)$}{d-manifold bordism group of
manifold $Y$} of $Y$ to be the set of all such d-bordism classes
$[\bX,\bs f]$ with $\vdim\bX=k$. As for $B_k(Y)$, we give $dB_k(Y)$
the structure of an abelian group, with zero element
$0_Y=[\bs\es,\bs\es]$, addition $[\bX,\bs f]+[\bX',\bs
f']=[\bX\amalg\bX',\bs f\amalg\bs f']$, and additive
inverses~$-[\bX,\bs f]=[-\bX,\bs f]$.

If $Y$ is oriented and of dimension $n$, we define a biadditive,
associative, supercommutative {\it intersection
product\/}\I{bordism!d-manifold bordism!intersection
product}\I{d-manifold!bordism!intersection product} $\bu:dB_k(Y)\t
dB_l(Y)\ra dB_{k+l-n}(Y)$ by
\e
[\bX,\bs f]\bu[\bX',\bs f']=[\bX\t_{\bs f,\bY,\bs f'}\bX',
\bs f\ci\bs\pi_\bX].
\label{ds15eq1}
\e
Here $\bX\t_{\bs f,\bY,\bs f'}\bX'$ exists as a d-manifold by
Theorem \ref{ds4thm8}(a), and is oriented by Theorem \ref{ds4thm15}.
Note that we do not need to restrict to $[X,f],[X',f']$ with $f,f'$
transverse as in Definition \ref{ds15def1}. Define a morphism
$\Pi_\bo^\dbo:B_k(Y)\ra dB_k(Y)$ for $k\ge 0$ by~$\Pi_\bo^\dbo:
[X,f]\mapsto\bigl[F_\Man^\dMan(X),F_\Man^\dMan(f)\bigr]$.
\label{ds15def2}
\end{dfn}

In \cite[\S 13.2]{Joyc6} we prove that $B_*(Y)$ and $dB_*(Y)$ are
isomorphic. See Spivak \cite[Th.~2.6]{Spiv} for the analogous result
for Spivak's derived manifolds.\I{Spivak's derived manifolds}

\begin{thm} For any manifold\/ $Y,$ we have $dB_k(Y)=0$ for $k<0,$
and\/ $\Pi_\bo^\dbo:B_k(Y)\ra dB_k(Y)$ is an isomorphism for $k\ge
0$. When $Y$ is oriented, $\Pi_\bo^\dbo$ identifies the intersection
products $\bu$ on\/ $B_*(Y)$ and\/~$dB_*(Y)$.
\label{ds15thm1}
\end{thm}

Here is the main idea in the proof of Theorem \ref{ds15thm1}. Let
$[\bX,\bs f]\in dB_k(Y)$. By Corollary \ref{ds4cor1} there exists an
embedding $\bs g:\bX\ra\bR^{\bs n}$ for $n\gg 0$. Then the direct
product $(\bs f,\bs g):\bX\ra\bY\t\bR^{\bs n}$ is also an embedding.
Theorem \ref{ds4thm13} shows that there exist an open set
$V\subseteq Y\t\R^n$, a vector bundle $E\ra V$ and $s\in C^\iy(E)$
such that $\bX\simeq\bS_{V,E,s}$. Let $\ti s\in C^\iy(E)$ be a
small, generic perturbation of $s$. As $\ti s$ is generic, the graph
of $\ti s$ in $E$ intersects the zero section transversely. Hence
$\ti X=\ti s^{-1}(0)$ is a $k$-manifold for $k\ge 0$, which is
compact and oriented for $\ti s-s$ small, and $\ti X=\es$ for $k<0$.
Set $\ti f=\pi_Y\vert_{\smash{\ti X}}:\ti X\ra Y$. Then
$\Pi_\bo^\dbo\bigl([\ti X,\ti f]\bigr)=[\bX,\bs f]$, so that
$\Pi_\bo^\dbo$ is surjective. A similar argument for $\rW,\bs e$ in
Definition \ref{ds15def2} shows that $\Pi_\bo^\dbo$ is injective.

Theorem \ref{ds15thm1} implies that we may define
projections\I{homology}\I{bordism!projection to homology}
\e
\Pi_\dbo^\hom:dB_k(Y)\longra H_k(Y;\Z)\;\> \text{by}\;\>
\Pi_\dbo^\hom=\Pi_\bo^\hom\ci(\Pi_\bo^\dbo)^{-1}.
\label{ds15eq2}
\e
We think of $\Pi_\dbo^\hom$ as a {\it virtual class map}.\I{virtual
class!for d-manifolds}\I{d-manifold!virtual class} Virtual classes
(or virtual cycles, or virtual chains)\I{virtual chain} are used in
several areas of geometry to construct enumerative invariants using
moduli spaces. In algebraic geometry, Behrend and Fantechi
\cite{BeFa} construct virtual classes for schemes with obstruction
theories.\I{scheme
with obstruction theory}\I{virtual class!for schemes with obstruction \\
theory} In symplectic geometry,\I{symplectic geometry} there are
many versions --- see for example Fukaya et al.\ \cite[\S 6]{FuOn},
\cite[\S A1]{FOOO}, Hofer et al.\ \cite{HWZ4}, and
McDuff~\cite{McDu}.

The main message we want to draw from this is that {\it compact
oriented d-manifolds admit virtual classes\/}\I{virtual class} (or
virtual cycles, or virtual chains, as appropriate). Thus, we can use
d-manifolds (and d-orbifolds) as the geometric structure on moduli
spaces\I{moduli space} in enumerative invariant problems such as
Gromov--Witten invariants,\I{Gromov--Witten invariants} Lagrangian
Floer cohomology,\I{Lagrangian Floer cohomology} Donaldson--Thomas
invariants,\I{Donaldson--Thomas invariants} \ldots, as this
structure is strong enough to contain all the `counting'
information.\I{bordism!d-manifold bordism|)}\I{d-manifold!bordism|)}

\subsection{Classical bordism for orbifolds}
\label{ds153}
\I{bordism!orbifold bordism|(}\I{orbifold!orbifold bordism|(}

In \cite[\S 13.3]{Joyc6} we generalize \S\ref{ds151} to orbifolds.
Here we will be brief, much more information is given in \cite[\S
13.3]{Joyc6}. We use the notation of \S\ref{ds9} on orbifolds $\Orb$
and \S\ref{ds12} on orbifolds with boundary $\Orbb$ freely. For
simplicity we do not distinguish between the 2-categories $\Orb$ in
\S\ref{ds91} and $\dotOrb$ in~\S\ref{ds121}.

\begin{dfn} Let $\cY$ be an orbifold, and $k\in\Z$. Consider
pairs $(\cX,f)$, where $\cX$ is a compact, oriented orbifold
(without boundary) with $\dim\cX=k$, and $f:\cX\ra\cY$ is a
1-morphism in $\Orb$. Define an equivalence relation $\sim$ between
such pairs by $(\cX,f)\sim (\cX',f')$ if there exists a compact,
oriented $(k\!+\!1)$-orbifold with boundary\I{orbifold with
boundary} $\oW$, a 1-morphism $e:\oW\ra\cY$ in $\Orbb$, an
orientation-preserving equivalence $j:-\cX\amalg\cX'\ra\pd\oW$, and
a 2-morphism $\eta:f\amalg f'\Ra e\ci i_\oW\ci j$ in~$\Orbb$.

Write $[\cX,f]$ for the $\sim$-equivalence class ({\it bordism
class\/}) of a pair $(\cX,f)$. For each $k\in\Z$, define the $k^{\it
th}$ {\it orbifold bordism group\/}
$B_k^\orb(\cY)$\G[BkYb]{$B_k^\orb(\cY)$}{orbifold bordism group of
orbifold $\cY$} of $\cY$ to be the set of all such bordism classes
$[\cX,f]$ with $\dim\cX=k$. It is an abelian group, with zero
$0_\cY=[\es,\es]$, addition $[\cX,f]+[\cX',f']=[\cX\amalg
\cX',f\amalg f']$, and additive inverses $-[\cX,f]=[-\cX,f]$. If
$k<0$ then $B_k^\orb(\cY)=0$.

Define {\it effective orbifold bordism\/}\I{bordism!orbifold
bordism!effective}\I{orbifold!orbifold bordism!effective}
$B_k^\eff(\cY)$\G[BkYc]{$B_k^\eff(\cY)$}{effective orbifold bordism
group of orbifold $\cY$} in the same way, but requiring both
orbifolds $\cX$ and orbifolds with boundary $\oW$ to be effective
(as in \S\ref{ds91} and \S\ref{ds121}) in pairs $(\cX,f)$ and the
definition of~$\sim$.

If $\cY$ is an orbifold, define group
morphisms\I{homology}\I{bordism!projection to homology}
\begin{gather*}
\Pi_\eff^\orb:B_k^\eff(\cY)\longra B_k^\orb(\cY),\;\>
\Pi_\orb^\hom:B_k^\orb(\cY)\longra H_k(\cY_\top;\Q)\\
\text{and}\qquad \Pi_\eff^\hom:B_k^\eff(\cY)\longra H_k(\cY_\top;\Z)
\end{gather*}
by $\Pi_\eff^\orb:[\cX,f]\mapsto[\cX,f]$ and $\Pi_\orb^\hom,
\Pi_\eff^\hom:[\cX,f]\mapsto (f_\top)_*([\cX])$, where $[\cX]$ is
the fundamental class of the compact, oriented $k$-orbifold $\cX$,
which lies in $H_k(\cX_\top;\Q)$ for general $\cX$, and in
$H_k(\cX_\top;\Z)$ for effective $\cX$. The morphisms
$\Pi_\eff^\orb:B_k^\eff(\cY)\ra B_k^\orb(\cY)$ are injective.

Suppose $\cY$ is an oriented orbifold of dimension $n$ which is a
manifold, that is, the orbifold groups $\Iso_\cY([y])$ are trivial
for all $[y]\in\cY_\top$. Define biadditive, associative,
supercommutative {\it intersection products\/}\I{bordism!orbifold
bordism!intersection product}\I{orbifold!orbifold
bordism!intersection product} $\bu:B_k^\orb(\cY)\t B_l^\orb(\cY)\ra
B_{k+l-n}^\orb(\cY)$ and $\bu:B_k^\eff(\cY)\t B_l^\eff(\cY)\ra
B_{k+l-n}^\eff(\cY)$ as follows. Given classes $[\cX,f],[\cX',f']$,
we perturb $f,f'$ in their bordism classes to make $f:\cX\ra\cY$ and
$f':\cX'\ra\cY$ transverse 1-morphisms, and then as in \eq{ds15eq1}
we set
\begin{equation*}
[\cX,f]\bu[\cX',f']=[\cX\t_{f,\cY,f'}\cX',f\ci\pi_\cX].
\end{equation*}
If $\cY$ is not a manifold then $f,f'$ may not admit transverse
perturbations.
\label{ds15def3}
\end{dfn}

Again, orbifold bordism is a generalized homology
theory.\I{generalized homology theory} Results of Druschel
\cite{Drus1,Drus2} and Angel \cite{Ange1,Ange2,Ange3} give a
complete description of the rational effective orbifold bordism ring
$B_*^\eff(*)\ot_\Z\Q$ when $\cY$ is the point $*$, and some
information on the full ring $B_*^\eff(*)$. It is much more
complicated than bordism $B_*(*)$ for manifolds in \S\ref{ds151},
because of contributions from orbifold strata.\I{bordism!orbifold
bordism!and orbifold strata|(}\I{orbifold!orbifold bordism!and
orbifold strata|(}\I{orbifold strata!of orbifolds|(}

As in \S\ref{ds92}, if $\cX$ is an orbifold and $\Ga$ a finite group
then we may define orbifold strata\I{orbifold!orbifold strata}
$\cX^{\Ga,\la}$ for $\la\in\La^\Ga_+$ and $\tcX^{\Ga,\mu}$ for
$\mu\in\La^\Ga_+/\Aut(\Ga)$, which are orbifolds, with proper
1-morphisms $O^{\Ga,\la}(\cX): \cX^{\Ga,\la}\ra\cX$ and $\ti
O{}^{\Ga,\mu}(\cX):\tcX^{\Ga,\mu}\ra\cX$. Hence, if $\cX$ is compact
then $\cX^{\Ga,\la},\tcX^{\Ga,\mu}$ are compact. If $\cX$ is
oriented then under extra conditions on $\Ga,\la,\mu$, which hold
automatically for $\cX^{\Ga,\la}$ if $\md{\Ga}$ is odd, we can
define natural orientations on $\cX^{\Ga,\la}, \tcX^{\Ga,\mu}$.
Using these ideas, under the assumptions on $\Ga,\la,\mu$ needed to
orient $\cX^{\Ga,\la},\tcX^{\Ga,\mu}$ we define morphisms
\ea
\Pi_\orb^{\Ga,\la}:B_k^\orb(\cY)\ra B_{k-\dim\la}^\orb(\cY)
\;\>\text{by}\;\> \Pi_\orb^{\Ga,\la}:[\cX,f]\mapsto
[\cX^{\Ga,\la},f\ci O^{\Ga,\la}(\cX)],
\label{ds15eq3}\\
\ti\Pi{}_\orb^{\Ga,\mu}:B_k^\orb(\cY)\ra B_{k-\dim\mu}^\orb(\cY)
\;\>\text{by}\;\> \ti\Pi{}_\orb^{\Ga,\mu}:[\cX,f]\mapsto
[\tcX^{\Ga,\mu},f\ci\ti O{}^{\Ga,\mu}(\cX)].
\label{ds15eq4}
\ea

One moral is that orbifold bordism groups $B_*^\orb(\cY),
B_*^\eff(\cY)$ are generally much bigger than manifold bordism
groups $B_*(Y)$, because in elements $[\cX,f]$ of orbifold bordism
groups, extra information is contained in the orbifold strata of
$\cX$. The morphisms $\Pi_\orb^{\Ga,\la},\ti\Pi{}_\orb^{\Ga,\mu}$
recover some of this extra information.\I{bordism!orbifold
bordism|)}\I{orbifold!orbifold bordism|)}\I{bordism!orbifold
bordism!and orbifold strata|)}\I{orbifold!orbifold bordism!and
orbifold strata|)}\I{orbifold strata!of orbifolds|)}

\subsection{Bordism for d-orbifolds}
\label{ds154}
\I{d-orbifold!bordism|(}\I{bordism!d-orbifold bordism|(}

In \cite[\S 13.4]{Joyc6} we combine the ideas of \S\ref{ds152} and
\S\ref{ds153} to define bordism for d-orbifolds. For simplicity we
do not distinguish between the 2-categories $\dOrb$ in \S\ref{ds111}
and $\bdOrb\subset\dOrbc$ in~\S\ref{ds141}.

\begin{dfn} Let $\cY$ be an orbifold, and $k\in\Z$. Consider
pairs $(\bcX,\bs f)$, where $\bcX\in\dOrb$ is a compact, oriented
d-orbifold without boundary with $\vdim\bcX=k$, and $\bs
f:\bcX\ra\bcY$ is a 1-morphism in $\dOrb$,
where~$\bcY=F_\Orb^\dOrb(\cY)$.

Define an equivalence relation $\sim$ between such pairs by
$(\bcX,\bs f) \sim(\bcX',\bs f')$ if there exists a compact,
oriented d-orbifold with boundary\I{d-orbifold with boundary} $\eW$
with $\vdim\eW=k+1$, a 1-morphism $\bs e:\eW\ra\bcY$ in $\dOrbb$, an
equivalence of oriented d-orbifolds $\bs
j:-\bcX\amalg\bcX'\ra\pd\eW$, and a 2-morphism $\eta:\bs f\amalg\bs
f'\Ra\bs e\ci\bs i_\eW\ci \bs j$.

Write $[\bcX,\bs f]$ for the $\sim$-equivalence class ({\it
d-bordism class\/}) of a pair $(\bcX,\bs f)$. For each $k\in\Z$,
define the $k^{\it th}$ {\it d-orbifold bordism group\/}
$dB_k^\orb(\cY)$\G[dBkYb]{$dB_k^\orb(\cY)$}{d-orbifold bordism group
of orbifold $\cY$} of $\cY$ to be the set of all such d-bordism
classes $[\bcX,\bs f]$ with $\vdim\bcX=k$. We give $dB_k^\orb(\cY)$
the structure of an abelian group, with zero element
$0_\cY=[\bs\es,\bs\es]$, addition $[\bcX,\bs f]+[\bcX',\bs
f']=[\bcX\amalg\bcX',\bs f\amalg\bs f']$, and additive
inverses~$-[\bcX,\bs f]=[-\bcX,\bs f]$.

Similarly, define the {\it semieffective d-orbifold bordism
group\/}\I{d-orbifold!bordism!semieffective|(}\I{bordism!d-orbifold
bordism!semieffective|(}
$dB_k^\sef(\cY)$\G[dBkYc]{$dB_k^\sef(\cY)$}{semieffective d-orbifold
bordism group of orbifold $\cY$} and the {\it effective d-orbifold
bordism
group\/}\I{d-orbifold!bordism!effective|(}\I{bordism!d-orbifold
bordism!effective|(}
$dB_k^\eff(\cY)$\G[dBkYd]{$dB_k^\eff(\cY)$}{effective d-orbifold
bordism group of orbifold $\cY$} as above, but taking $\bcX$ and
$\eW$ to be semieffective, or effective, respectively, in the sense
of \S\ref{ds1110} and~\S\ref{ds1410}.

If $\cY$ is oriented and of dimension $n$, we define a biadditive,
associative, supercommutative {\it intersection
product\/}\I{d-orbifold!bordism!intersection
product}\I{bordism!d-orbifold bordism!intersection product}
$\bu:dB_k^\orb(\cY)\!\t\! dB_l^\orb(\cY)\!\ra\!
dB_{k+l-n}^\orb(\cY)$~by
\begin{equation*}
[\bcX,\bs f]\bu[\bcX',\bs f']=[\bcX\t_{\bs f,\bcY,\bs f'}\bcX',\bs
f\ci\bs\pi_\bcX],
\end{equation*}
as in \eq{ds15eq1}. Here $\bcX\t_{\bs f,\bcY,\bs f'}\bcX'$ exists in
$\dOrb$ by Theorem \ref{ds11thm7}(a), and is oriented by the
d-orbifold analogue of Theorem~\ref{ds4thm15}.

If $\cY$ is an orbifold, define group morphisms
\e
\begin{gathered}
\begin{aligned}
\Pi_\orb^\sef&:B_k^\orb(\cY)\longra dB_k^\sef(\cY),&
\Pi_\eff^\deff&:B_k^\eff(\cY)\longra dB_k^\eff(\cY),\\
\Pi_\deff^\sef&:dB_k^\eff(\cY)\longra dB_k^\sef(\cY),&
\Pi_\deff^\dorb&:dB_k^\eff(\cY)\longra dB_k^\orb(\cY),
\end{aligned}\\
\text{and}\qquad\Pi_\sef^\dorb:dB_k^\sef(\cY)\longra dB_k^\orb(\cY)
\end{gathered}
\label{ds15eq5}
\e
by $\Pi_\orb^\sef,\Pi_\eff^\deff:[\cX,f]\mapsto[\bcX,\bs f]$, where
$\bcX,\bs f=F_\Orb^\dOrb(\cX,f)$, and~$\Pi_\deff^\sef,\ab
\Pi_\deff^\dorb,\ab\Pi_\sef^\dorb:[\bcX,\bs f]\mapsto[\bcX,\bs f]$.
\label{ds15def4}
\end{dfn}

Here is the main result of \cite[\S 13.4]{Joyc6}, an orbifold
analogue of Theorem \ref{ds15thm1}. The key idea is that
semieffective (or effective)
d-orbifolds\I{d-orbifold!semieffective}\I{d-orbifold!effective}
$\bcX$ can be perturbed\I{d-orbifold!perturbing to orbifolds} to
(effective) orbifolds, as in \S\ref{ds1110}; to make this rigorous,
we use good coordinate systems\I{good coordinate
system}\I{d-orbifold!good coordinate system} on $\bcX$, as
in~\S\ref{ds119}.

\begin{thm} For any orbifold\/ $\cY,$ the maps
$\Pi_\orb^\sef:B_k^\orb(\cY)\ra dB_k^\sef(\cY)$ and\/
$\Pi_\eff^\deff:B_k^\eff(\cY)\ra dB_k^\eff(\cY)$ in \eq{ds15eq5} are
isomorphisms for all\/ $k\in\Z$.
\label{ds15thm2}
\end{thm}

As for \eq{ds15eq2}, the theorem implies that we may define
projections
\begin{align*}
\Pi_\sef^\hom:dB_k^\sef(\cY)\ra H_k(\cY_\top;\Q),\;\>
\Pi_\deff^\hom:dB_k^\eff(\cY)\ra H_k(\cY_\top;\Z)\\
\text{by}\;\> \Pi_\sef^\hom=\Pi_\orb^\hom\ci(\Pi_\orb^\sef)^{-1}\;\>
\text{and}\;\> \Pi_\deff^\hom=\Pi_\eff^\hom\ci(\Pi_\eff^\deff)^{-1}.
\end{align*}
We think of these $\Pi_\sef^\hom,\Pi_\deff^\hom$ as {\it virtual
class maps\/}\I{virtual class!for d-orbifolds}\I{d-orbifold!virtual
class} on $dB_*^\sef(\cY),dB_*^\eff(\cY)$. In fact, with more work,
one can also define virtual class maps
on~$dB_*^\orb(\cY)$:\I{homology}\I{bordism!projection to homology}
\e
\Pi_\dorb^\hom:dB_k^\orb(\cY)\longra H_k(\cY_\top;\Q),
\label{ds15eq6}
\e
satisfying $\Pi_\dorb^\hom\ci\Pi_\sef^\dorb=\Pi_\sef^\hom$, for
instance following the method of Fukaya et al.\ \cite[\S 6]{FuOn},
\cite[\S A1]{FOOO} for virtual classes of Kuranishi
spaces\I{Kuranishi space}\I{virtual class!for Kuranishi spaces}
using `multisections'.

In future work the author intends to define a virtual
chain\I{virtual chain} construction for d-manifolds and d-orbifolds,
expressed in terms of new (co)homology theories whose (co)chains are
built from d-manifolds or d-orbifolds, as for the `Kuranishi
(co)homology'\I{Kuranishi (co)homology} described
in~\cite{Joyc1,Joyc2}.

As in \S\ref{ds118}, if $\bcX$ is a d-orbifold and $\Ga$ a finite
group then we may define orbifold strata\I{d-orbifold!orbifold
strata}\I{bordism!d-orbifold bordism!and orbifold
strata|(}\I{d-orbifold!bordism!and orbifold strata|(}\I{orbifold
strata!of d-orbifolds|(} $\bcX^{\Ga,\la}$ for $\la\in\La^\Ga$ and
$\btcX{}^{\Ga,\mu}$ for $\mu\in\La^\Ga/\Aut(\Ga)$, which are
d-orbifolds, with proper 1-morphisms $\bs O^{\Ga,\la}(\bcX):
\bcX^{\Ga,\la}\ra\bcX$ and $\bs{\ti
O}{}^{\Ga,\mu}(\bcX):\btcX{}^{\Ga,\mu}\ra\bcX$. Hence, if $\bcX$ is
compact then $\bcX^{\Ga,\la},\btcX{}^{\Ga,\mu}$ are compact. If
$\bcX$ is oriented and $\Ga$ is odd then we under extra conditions
on $\mu$ can define natural orientations on
$\bcX^{\Ga,\la},\btcX{}^{\Ga,\mu}$. As in
\eq{ds15eq3}--\eq{ds15eq4}, for such $\Ga,\la,\mu$ we define
morphisms
\begin{align*}
\Pi_\dorb^{\Ga,\la}:dB_k^\orb(\cY)\ra dB_{k-\dim\la}^\orb(\cY)
\;\>\text{by}\;\> \Pi_\dorb^{\Ga,\la}:[\bcX,\bs f]\mapsto
[\bcX^{\Ga,\la},\bs f\ci\bs O^{\Ga,\la}(\bcX)],\\
\ti\Pi{}_\dorb^{\Ga,\mu}:dB_k^\orb(\cY)\ra dB_{k-\dim\mu}^\orb(\cY)
\;\>\text{by}\;\> \ti\Pi{}_\dorb^{\Ga,\mu}:[\bcX,\bs f]\mapsto
[\btcX{}^{\Ga,\mu},\bs f\ci\bs{\ti O}{}^{\Ga,\mu}(\bcX)].
\end{align*}

We can use these operators $\Pi_\dorb^{\Ga,\la}$ to study the
d-orbifold bordism ring $dB_*^\orb(*)$. Let $\Ga$ be a finite group
with $\md{\Ga}$ odd, and $R$ be a nontrivial $\Ga$-representation.
Define an element $[\bs*\t_{\bs 0,\bs R,\bs 0}\bs */\Ga,\bs\pi]$ in
$dB_{-\dim R}^\orb(*)$, where $\bs R=F_\Man^\dMan(R)$, and set
$\la=[R]\in \La^\Ga_+$. Then $\Pi_\dorb^{\Ga,-\la}\bigl([\bs*\t_{\bs
0,\bs R,\bs 0}\bs */\Ga,\bs\pi]\bigr)\in dB_0^\orb(*)$, so
$\Pi_\dorb^\hom\ci \Pi_\dorb^{\Ga,-\la}\bigl([\bs*\t_{\bs 0,\bs
R,\bs 0}\bs */\Ga,\bs\pi]\bigr)$ lies in $H_0^\orb(*;\Q)\cong\Q$ by
\eq{ds15eq6}. Calculation shows that $\Pi_\dorb^\hom\ci
\Pi_\dorb^{\Ga,-\la} \bigl([\bs*\t_{\bs 0,\bs R,\bs 0}\bs
*/\Ga,\bs\pi]\bigr)$ is either $\md{\Aut(\Ga)}/\md{\Ga}$ or 0,
depending on $\la$. In the first case, $[\bs*\t_{\bs 0,\bs R,\bs
0}\bs */\Ga,\bs\pi]$ has infinite order in $dB_{-\dim R}^\orb(*)$.
Extending this argument, we can show that $dB_{4k}^\orb(*)$ has
infinite rank for all $k\le 0$. In contrast,
$dB_k^\sef(*)=dB_k^\eff(*)=0$ for all $k<0$ by
Theorem~\ref{ds15thm2}.\I{bordism|)}\I{d-orbifold!bordism|)}
\I{bordism!d-orbifold
bordism|)}\I{d-orbifold!bordism!semieffective|)}\I{bordism!d-orbifold
bordism!semieffective|)}\I{d-orbifold!bordism!effective|)}\I{bordism!d-orbifold
bordism!effective|)}\I{bordism!d-orbifold bordism!and orbifold
strata|)}\I{d-orbifold!bordism!and orbifold strata|)}\I{orbifold
strata!of d-orbifolds|)}

\section[Relation to other classes of spaces in
mathematics]{Relation to other classes of spaces}
\label{ds16}
\I{derived manifold|see{Spivak's derived manifolds}}

In \cite[Chap.~14]{Joyc6} we study the relationships between
d-manifolds and d-orbifolds and other classes of geometric spaces in
the literature. The next theorem summarizes our
results:\I{truncation functor|(}\I{functor!truncation|(}

\begin{thm} We may construct `truncation functors' from various
classes of geometric spaces to d-manifolds and d-orbifolds, as
follows:
\begin{itemize}
\setlength{\itemsep}{0pt}
\setlength{\parsep}{0pt}
\item[{\bf(a)}] There is a functor\/ $\Pi_{\bf BManFS}^{\bf
dMan}:\mathop{\bf BManFS}\ra\Ho(\dMan),$ where $\mathop{\bf
BManFS}$ is a category whose objects are triples $(V,E,s)$ of a
Banach manifold $V,$ Banach vector bundle $E\ra V,$ and smooth
section $s:V\ra E$ whose linearization $\d s\vert_x:T_xV\ra
E\vert_x$ is Fredholm with index $n\in\Z$ for each $x\in V$
with\/ $s\vert_x=0,$ and\/ $\Ho(\dMan)$ is the homotopy
category\I{2-category!homotopy category|(}\I{homotopy
category|(} of the $2$-category
of d-manifolds\/~$\dMan$.\I{d-manifold!and Banach manifolds with \\
Fredholm sections}\I{Banach manifold}

There is also an orbifold version\I{d-orbifold!and Banach
orbifolds with \\ Fredholm sections} $\Pi_{\bf BOrbFS}^{\bf
dOrb}:\Ho(\mathop{\bf BOrbFS})\ra\Ho(\dOrb)$ of this using
Banach orbifolds $\cV,$ and `corners' versions of
both.\I{d-manifold with corners!and Banach manifolds with \\
Fredholm sections}\I{d-orbifold with corners!and Banach
orbifolds with \\ Fredholm sections}
\item[{\bf(b)}] There is a functor\/ $\Pi_{\bf MPolFS}^{\bf
dMan}:\mathop{\bf MPolFS}\ra\Ho(\dMan),$ where $\mathop{\bf
MPolFS}$ is a category whose objects are triples $(V,E,s)$ of an
\begin{bfseries}M-polyfold\end{bfseries}\I{d-manifold!and
M-polyfolds} without boundary\/ $V$ as in Hofer, Wysocki and
Zehnder\/ {\rm \cite[\S 3.3]{HWZ1},}\I{polyfold!and d-manifolds}
a fillable strong M-polyfold bundle $E$ over\/ $V$ {\rm \cite[\S
4.3]{HWZ1},} and an sc-smooth Fredholm section $s$ of\/ $E$ {\rm
\cite[\S 4.4]{HWZ1}} whose linearization $\d s\vert_x:T_xV\ra
E\vert_x$ {\rm \cite[\S 4.4]{HWZ1}} has Fredholm index\/
$n\in\Z$ for all\/ $x\in V$ with\/~$s\vert_x=0$.

There is also an orbifold version\I{d-orbifold!and
polyfolds}\I{polyfold!and d-orbifolds} $\Pi_{\bf PolFS}^{\bf
dOrb}:\Ho(\mathop{\bf PolFS})\ra\Ho(\dOrb)$ of this using
\begin{bfseries}polyfolds\end{bfseries} $\cV,$ and `corners'
versions of both.\I{d-manifold with corners!and
M-polyfolds}\I{d-orbifold with corners!and polyfolds}
\item[{\bf(c)}] Given a d-orbifold with corners $\eX,$ we can
construct a \begin{bfseries}Kuranishi
space\end{bfseries}\I{d-orbifold!and Kuranishi
spaces}\I{d-orbifold with corners!and Kuranishi
spaces}\I{Kuranishi space!and d-orbifolds} $(X,\ka)$ in the
sense of Fukaya, Oh, Ohta and Ono\/ {\rm\cite[\S A]{FOOO},} with
the same underlying topological space $X$. Conversely, given a
Kuranishi space $(X,\ka),$ we can construct a d-orbifold with
corners $\eX'$. Composing the two constructions, $\eX$ and\/
$\eX'$ are equivalent in\/~$\dOrbc$.

Very roughly speaking, this means that the `categories' of
d-orbifolds with corners, and Kuranishi spaces, are equivalent.
However, Fukaya et al.\ {\rm\cite{FOOO}} do not define morphisms
of Kuranishi spaces, so we have no category of Kuranishi spaces.
\item[{\bf(d)}] There is a functor\/ $\Pi_{\bf SchObs}^{\bf
dMan}:\mathop{\bf Sch_\C Obs}\ra\Ho(\dMan),$\I{d-manifold!and
schemes with obstruction \\ theories}\I{scheme with obstruction
theory!and d-manifolds} where $\mathop{\bf Sch_\C Obs}$ is a
category whose objects are triples $(X,E^\bu,\phi),$ for $X$ a
separated, second countable $\C$-scheme and\/ $\phi:E^\bu\ra
\tau_{\ge -1}(L_X)$ a perfect obstruction theory on $X$ with
constant virtual dimension, in the sense of Behrend and
Fantechi\/ {\rm\cite{BeFa}}. We may define a natural
orientation\I{d-manifold!orientations} on $\Pi_{\bf SchObs}^{\bf
dMan}(X,E^\bu,\phi)$ for each\/~$(X,E^\bu,\phi)$.

There is also an orbifold version\/ $\Pi_{\bf StaObs}^{\bf
dOrb}:\Ho(\mathop{\bf Sta_\C Obs})\ra\Ho(\dOrb),$ taking $\cX$
to be a Deligne--Mumford\/ $\C$-stack.\I{d-orbifold!and
Deligne--Mumford stacks with obstruction
theories}\I{Deligne--Mumford stack with \\
obstruction theory!and d-orbifolds}\I{stack}
\item[{\bf(e)}] There is a functor\/ $\Pi_{\bf QsDSch}^{\bf
dMan}:\Ho(\mathop{\bf QsDSch_\C})\longra\Ho(\dMan),$ where
$\mathop{\bf QsDSch_\C}$ is the $\iy$-category\I{$\iy$-category}
of separated, second countable, quasi-smooth\I{quasi-smooth}
derived\/ $\C$-schemes\I{d-manifold!and quasi-smooth derived
schemes}\I{derived scheme!quasi-smooth!and d-manifolds}  $X$ of
constant dimension, as in To\"en and Vezzosi\/
{\rm\cite{Toen,ToVe1,ToVe2}}. We may define a natural
orientation on $\Pi_{\bf QsDSch}^{\bf dMan}(X)$ for each\/~$X$.

There is also an orbifold version\/ $\Pi_{\bf QsDSta}^{\bf
dOrb}:\Ho(\mathop{\bf QsDSta_\C})\ra\Ho(\dOrb),$ taking $\cX$ to
be a derived Deligne--Mumford\/~$\C$-stack.\I{d-orbifold!and
quasi-smooth derived Deligne--Mumford stacks}\I{derived
Deligne--Mumford stack!quasi-smooth!and d-orbifolds}\I{stack}
\item[{\bf(f)}]{\bf (Borisov \cite{Bori})} There is a natural
functor\/ $\Pi_{\bf DerMan}^{\bf dMan}:\Ho(\mathop{\bf
DerMan}^{\bf pd}_{\smash{\bf ft}})\ra\Ho(\dMan_{\bf pr})$ from
the homotopy category of the $\iy$-category $\smash{\DerMan^{\bf
pd}_{\smash{\bf ft}}}$ of \begin{bfseries}derived
manifolds\end{bfseries}\I{d-manifold!and Spivak's derived
manifolds}\I{Spivak's derived manifolds!and d-manifolds} of
finite type with pure dimension, in the sense of Spivak\/
{\rm\cite{Spiv},} to the homotopy category of the full\/
$2$-subcategory $\dMan_{\smash{\bf pr}}$ of principal
d-manifolds in\/ $\dMan$. This functor induces a bijection
between isomorphism classes of objects in $\Ho(\mathop{\bf
DerMan}^{\bf pd}_{\smash{\bf ft}})$ and\/ $\Ho(\dMan_{\bf pr})$.
It is full, but not faithful. If\/ $[\bs f]$ is a morphism in
$\Ho(\mathop{\bf DerMan}^{\bf pd}_{\smash{\bf ft}}),$ then $[\bs
f]$ is an isomorphism if and only if\/ $\Pi_{\bf DerMan}^{\bf
dMan}([\bs f])$ is an isomorphism.\I{truncation
functor|)}\I{functor!truncation|)}
\end{itemize}
\label{ds16thm1}
\end{thm}

Here, as in \S\ref{dsA3}, if $\bs{\cal C}$ is a 2-category (or
$\iy$-category), the {\it homotopy category\/} $\Ho(\bs{\cal C})$ of
$\bs{\cal C}$ is the category whose objects are objects of $\bs{\cal
C}$, and whose morphisms are 2-isomorphism classes of 1-morphisms in
$\bs{\cal C}$. Then equivalences in $\bs{\cal C}$ become
isomorphisms in $\Ho(\bs{\cal C})$, 2-commutative diagrams in
$\bs{\cal C}$ become commutative diagrams in $\Ho(\bs{\cal C})$, and
so on.\I{2-category!homotopy category|)}\I{homotopy category|)}

One moral of Theorem \ref{ds16thm1} is that essentially every
geometric structure on moduli spaces\I{moduli space|(} which is used
to define enumerative invariants, either in differential geometry,
or in algebraic geometry over $\C$, has a truncation functor to
d-manifolds or d-orbifolds. Combining Theorem \ref{ds16thm1} with
proofs from the literature of the existence on moduli spaces of the
geometric structures listed in Theorem \ref{ds16thm1}, in
\cite[Chap.~14]{Joyc6} we deduce:

\begin{thm}{\bf(i)} Any solution set of a smooth nonlinear elliptic
equation with fixed topological invariants on a compact manifold
naturally has the structure of a d-manifold, uniquely up to
equivalence in\/~$\dMan$.\I{d-manifold!and solutions of elliptic
equations}\I{elliptic equations}\I{moduli space!of solutions of
nonlinear elliptic \\ equations}

For example, let\/ $(M,g),(N,h)$ be Riemannian manifolds, with\/ $M$
compact. Then the family of \begin{bfseries}harmonic
maps\end{bfseries}\I{harmonic maps}\I{moduli space!of harmonic maps}
$f:M\ra N$ is a d-manifold\/ $\bs\cH_{M,N}$ with\/
$\vdim\bs\cH_{M,N}=0$. If\/ $M={\cal S}^1,$ then $\bs\cH_{M,N}$ is
the moduli space of
\begin{bfseries}parametrized closed geodesics\end{bfseries}
in\/~$(N,h)$.
\smallskip

\noindent{\bf(ii)} Let\/ $(X,\om)$\I{moduli space!of J-holomorphic
curves@of $J$-holomorphic curves|(}\I{symplectic geometry|(} be a
compact symplectic manifold of dimension $2n,$ and\/ $J$ an almost
complex structure on\/ $X$ compatible with\/ $\om$. For\/ $\be\in
H_2(X,\Z)$ and\/ $g,m\ge 0$, write $\oM_{g,m}(X,J,\be)$ for the
moduli space of stable triples $(\Si,\vec z,u)$ for $\Si$ a genus
$g$ prestable Riemann surface with\/ $m$ marked points $\vec
z=(z_1,\ldots,z_m)$ and\/ $u:\Si\ra X$ a $J$-holomorphic map with\/
$[u(\Si)]=\be$ in\/ $H_2(X,\Z)$. Using results of Hofer, Wysocki and
Zehnder\/ {\rm\cite{HWZ6}} involving their theory of
polyfolds,\I{polyfold} we can make $\oM_{g,m}(X,J,\be)$ into a
compact, oriented d-orbifold\/~$\bs\oM_{g,m}(X,J,\be)$.

\smallskip

\noindent{\bf(iii)} Let\/ $(X,\om)$ be a compact symplectic
manifold, $J$ an almost complex structure on\/ $X$ compatible with\/
$\om,$ and\/ $Y$ a compact, embedded Lagrangian
submanifold\I{Lagrangian submanifold} in\/ $X$. For\/ $\be\in
H_2(X,Y;\Z)$ and\/ $k\ge 0,$ write $\oM_k(X,Y,J,\be)$ for the moduli
space of $J$-\begin{bfseries}holomorphic stable maps\end{bfseries}
$(\Si,\vec z,u)$ to $X$ from a prestable holomorphic disc $\Si$
with\/ $k$ boundary marked points $\vec z=(z_1,\ldots,z_k),$ with\/
$u(\pd\Si)\subseteq Y$ and\/ $[u(\Si)]=\be$ in $H_2(X,Y;\Z)$. Using
results of Fukaya, Oh, Ohta and Ono\/ {\rm\cite[\S 7--\S 8]{FOOO}}
involving their theory of Kuranishi spaces,\I{Kuranishi space} we
can make $\oM_k(X,Y,J,\be)$ into a compact d-orbifold with corners\/
$\bs\oM_k(X,Y,J,\be)$. Given a relative spin structure for $(X,Y),$
we may define an orientation on\/~$\bs\oM_k(X,Y,J,\be)$.\I{moduli
space!of J-holomorphic curves@of $J$-holomorphic
curves|)}\I{symplectic geometry|)}\I{d-orbifold!orientations}
\smallskip

\noindent{\bf(iv)} Let\/ $X$ be a complex projective manifold, and\/
$\oM_{g,m}(X,\be)$ the Deligne--Mumford moduli\/ $\C$-stack\I{stack}
of stable triples $(\Si,\vec z,u)$ for $\Si$ a genus $g$ prestable
Riemann surface with\/ $m$ marked points $\vec z=(z_1,\ldots,z_m)$
and\/ $u:\Si\ra X$ a morphism with\/ $u_*([\Si])=\be\in H_2(X;\Z)$.
Then Behrend\/ {\rm\cite{Behr}} defines a perfect obstruction
theory\I{Deligne--Mumford stack with \\ obstruction theory} on
$\oM_{g,m}(X,\be),$ so we can make $\oM_{g,m}(X,\be)$ into a
compact, oriented d-orbifold\/~$\bs\oM_{g,m}(X,\be)$.\I{moduli
space!of algebraic curves}\I{d-orbifold!orientations}
\smallskip

\noindent{\bf(v)} Let\/ $X$ be a complex algebraic surface, and\/
$\cM$ a stable moduli\/ $\C$-scheme of vector bundles or coherent
sheaves $E$ on $X$ with fixed Chern character. Then Mochizuki
{\rm\cite{Moch}} defines a perfect obstruction theory\I{scheme with
obstruction theory|(} on $\cM,$ so we can make $\cM$ into an
oriented d-manifold\/~$\bs\cM$.\I{moduli space!of coherent sheaves
on a surface}\I{d-manifold!orientations|(}
\smallskip

\noindent{\bf(vi)} Let\/ $X$ be a complex Calabi--Yau
$3$-fold\I{Calabi--Yau 3-fold} or smooth Fano $3$-fold,\I{Fano
3-fold} and\/ $\cM$ a stable moduli\/ $\C$-scheme of coherent
sheaves $E$ on $X$ with fixed Hilbert polynomial.\I{moduli space!of
coherent sheaves on a 3-fold} Then Thomas\/ {\rm\cite{Thom}} defines
a perfect obstruction theory on $\cM,$ so we can make $\cM$ into an
oriented d-manifold\/~$\bs\cM$.
\smallskip

\noindent{\bf(vii)} Let\/ $X$ be a smooth complex projective
$3$-fold, and\/ $\cM$ a moduli\/ $\C$-scheme of `stable PT pairs'
$(C,D)$ in $X,$ where $C\subset X$ is a curve and $D\subset C$ is a
divisor. Then Pandharipande and Thomas\/ {\rm\cite{PaTh}} define a
perfect obstruction theory on $\cM,$ so we can make $\cM$ into a
compact, oriented d-manifold\/~$\bs\cM$.\I{moduli space!of PT pairs
on a 3-fold}
\smallskip

\noindent{\bf(ix)} Let\/ $X$ be a complex Calabi--Yau
$3$-fold,\I{Calabi--Yau 3-fold} and\/ $\cM$ a separated moduli\/
$\C$-scheme of simple perfect complexes in the derived
category\I{derived category} $D^b\coh(X)$. Then Huybrechts and
Thomas\/ {\rm\cite{HuTh}} define a perfect obstruction theory on
$\cM,$ so we can make $\cM$ into an oriented
d-manifold\/~$\bs\cM$.\I{moduli space!of perfect complexes on a
3-fold}\I{moduli space|)}\I{scheme with obstruction
theory|)}\I{d-manifold!orientations|)}
\label{ds16thm2}
\end{thm}

A consequence of Theorems \ref{ds16thm1} and \ref{ds16thm2} is that
we can use d-manifolds and d-orbifolds in many areas of differential
and algebraic geometry to do with moduli spaces, enumerative
invariants, Floer homology theories, and so on. Working with
d-manifolds and d-orbifolds instead of the geometric structures
currently in use may lead to new results, or simplifications of
known proofs. We discuss some applications in symplectic geometry:

\begin{rem}{\bf(a)}\I{Gromov--Witten invariants|(}\I{symplectic
geometry|(}\I{d-orbifold!bordism!and Gromov--Witten
invariants|(}\I{moduli space!of J-holomorphic curves@of
$J$-holomorphic curves|(}\I{virtual class|(} Suppose $(X,\om)$ is a
compact symplectic manifold, and choose an almost complex structure
$J$ on $X$ compatible with $\om$. Let $g,m$ be nonnegative integers,
and $\be\in H_2(X;\Z)$. Consider the moduli space $\oM_{g,m}(X,\be)$
of genus $g$ stable $J$-holomorphic curves $(\Si,\vec z,u)$ with $m$
marked points in homology class $\be$ in $X$. Using Hofer, Wysocki
and Zehnder \cite{HWZ6}, Theorem \ref{ds16thm2}(ii) makes
$\oM_{g,m}(X,\be)$ into a compact, oriented
d-orbifold~$\bs\oM_{g,m}(X,J,\be)$.

There are evaluation maps $\ev_i:\oM_{g,m}(X,J,\be)\ra X$ for
$i=1,\ldots,m$ mapping $\ev_i:(\Si,\vec z,u)\mapsto u(z_i)$. Also,
if $2g+m\ge 3$ there is a projection
$\pi_{g,m}:\oM_{g,m}(X,J,\be)\ra\oM_{g,m}$ mapping $\ev_i:(\Si,\vec
z,u)\mapsto\bar\Si$, where $\oM_{g,m}$ is the moduli space of
Deligne--Mumford stable Riemann surfaces of genus $g$ with $m$
marked points, a compact orbifold of real dimension $2(m+3g-3)$, and
$\bar\Si$ is the stabilization of $\Si$. As in \cite[\S
14.2]{Joyc6}, these lift naturally to 1-morphisms
$\bs\ev_i:\bs\oM_{g,m}(X,J,\be)\ra \bcX$ and $\bs\pi_{g,m}:
\bs\oM_{g,m}(X,J,\be)\ra\bs\oM_{g,m}$ in $\dOrb$, where
$\bcX=F_\Man^\dOrb(X)$ and~$\bs\oM_{g,m}=F_\Orb^\dOrb(\oM_{g,m})$.

Define the {\it Gromov--Witten d-orbifold bordism invariant\/}
$GW_{g,m}^\dorb(X,\om,\be)$ by
\begin{align*}
&GW_{g,m}^\dorb(X,\om,\be)=\\
&\begin{cases} \bigl[\bs\oM_{g,m}(X,J,\be),
\bs{\rm ev}_1\!\t\!\cdots\!\t\!\bs{\rm ev}_m\bigr]\!\in\!
dB^\orb_{2k}(X^m), & 2g\!+\!m\!<\!3, \\
\bigl[\bs\oM_{g,m}(X,J,\be),\bs{\rm ev}_1\!\t\!\cdots\!\t\!
\bs{\rm ev}_m\!\t\!\bs\pi_{g,m}\bigr]\!\in\! dB^\orb_{2k}
(X^m\!\t\!\oM_{g,m}), & 2g\!+\!m\!\ge\!3,
\end{cases}
\end{align*}
where $k=c_1(X)\cdot\be+(n-3)(1-g)+m$, with $\dim X=2n$. As in
\cite[\S 14.2]{Joyc6}, we can prove this is {\it independent of the
choice of almost complex structure\/ $J$ on\/} $X$. This is because
any two choices $J_0,J_1$ for $J$ may be joined by a smooth path of
almost complex structures $J_t$, $t\in[0,1]$ compatible with $\om$.
One can then make $\coprod_{t\in[0,1]}\oM_{g,m}(X, J_t,\be)$ into a
compact oriented d-orbifold with boundary, whose boundary is
$\bs\oM_{g,m}(X,J_1,\be) \amalg -\bs\oM_{g,m}(X,J_0,\be)$, and this
defines a bordism between the moduli spaces for $J_0$ and~$J_1$.

Applying the virtual class maps $\Pi_\dorb^\hom$ in \eq{ds15eq6}
gives homology classes
\begin{equation*}
\Pi_\dorb^\hom\bigl(GW_{g,m}^\dorb(X,\om,\be)\bigr)\in
\begin{cases} H_{2k}(X^m;\Q), & 2g+m<3, \\
H_{2k}\bigl(X^m\t(\oM_{g,m})_\top;\Q\bigr), & 2g+m\ge 3.\end{cases}
\end{equation*}
These homology classes are essentially the {\it Gromov--Witten
invariants\/} of $(X,\om)$, as in Hofer, Wysocki and Zehnder
\cite{HWZ6} or Fukaya and Ono~\cite{FuOn}.

As in \S\ref{ds154}, d-orbifold bordism groups are far larger than
homology groups, and can have infinite rank even in negative
degrees. Therefore the new invariants $GW_{g,m}^\dorb(X,\om,\be)$
{\it contain more information than conventional Gromov--Witten
invariants}, in particular, on the orbifold strata of the moduli
spaces $\bs\oM_{g,m}(X,J,\be)$, which can be recovered by applying
the functors $\Pi_\dorb^{\Ga,\la}$ of \S\ref{ds154}. Since the
$dB_*^\orb(-)$ are defined over $\Z$ rather than $\Q$, these new
invariants may be good tools for studying {\it integrality
properties\/} of Gromov--Witten invariants.\I{d-orbifold!bordism!and
Gromov--Witten invariants|)}
\smallskip

\noindent{\bf(b)} Let $(X,J)$ be a projective complex manifold, with
K\"ahler form $\om$. Consider the moduli space\I{moduli space!of
J-holomorphic curves@of $J$-holomorphic curves} $\oM_{g,m}(X,\be)$
of genus $g$ stable $J$-holomorphic curves with $m$ marked points in
class $\be\in H_2(X;\Z)$ in $X$. As in {\bf(a)}, by constructing a
virtual class\I{virtual class} $[\oM_{g,m}(X,\be)]_{\rm virt}$, we
define the {\it Gromov--Witten invariants\/} of~$X$.

We can do this in two different ways. In symplectic geometry, by
Hofer, Wysocki and Zehnder \cite{HWZ6} we can write
$\oM_{g,m}(X,\be)$ as the zeroes of a Fredholm section of a
polyfold\I{polyfold} bundle over a polyfold, as in Theorem
\ref{ds16thm2}(ii). In algebraic geometry, as in Behrend \cite{Behr}
we can make $\oM_{g,m}(X,\be)$ into a
Deligne--Mumford $\C$-stack\I{Deligne--Mumford stack with \\
obstruction theory}\I{stack} with obstruction theory, as in Theorem
\ref{ds16thm2}(iv). We then apply the virtual class construction for
polyfolds \cite{HWZ4} or stacks with obstruction theories
\cite{BeFa} to define the virtual class, and the Gromov--Witten
invariants.

A priori it is not clear that the symplectic and algebraic routes
give the same Gromov--Witten invariants (though it is generally
accepted that they do). Comparing the two is difficult, as the
techniques used are so different.

Now by Theorem \ref{ds16thm1}(b),(d) we can make $\oM_{g,m}(X,\be)$
into a compact oriented d-orbifold $\bs\oM_{g,m}(X,J,\be)_{\rm sym}$
or $\bs\oM_{g,m}(X,J,\be)_{\rm alg}$ in two different ways, by
applying the truncation functors\I{truncation
functor}\I{functor!truncation} from polyfolds and Deligne--Mumford
$\C$-stacks\I{stack} with obstruction theories. Then we could apply
a virtual class construction for d-orbifolds, as discussed in
\S\ref{ds154}, to define Gromov--Witten invariants.

Thus, d-orbifolds provide a bridge between symplectic and algebraic
geometry, and a means to compare symplectic and algebraic
Gromov--Witten theory. To use d-orbifolds to prove that symplectic
and algebraic Gromov--Witten invariants coincide, we would have to
do two things. Firstly, we would have to show that the virtual class
constructions for polyfolds and Deligne--Mumford stacks\I{stack}
with obstruction theories, yield the same answer as applying the
truncation functors to d-orbifolds and then applying the virtual
class construction for d-orbifolds. This does not look difficult.

Secondly, we would have to compare the d-orbifolds
$\bs\oM_{g,m}(X,J,\be)_{\rm sym}$ and $\bs\oM_{g,m}(X,J,\be)_{\rm
alg}$. The author does {\it not\/} expect these to be equivalent in
$\dOrb$. This is because the smooth structure near singular curves
depends on a choice of {\it gluing profile\/}\I{gluing
profile}\I{polyfold!gluing profile} $\varphi:(0,1]\ra[0,\iy)$, in
the language of Hofer et al.\ \cite[\S 4.2]{HWZ5}, \cite[\S
2.1]{HWZ6}. As in \cite[\S 2.1]{HWZ6}, the gluing profiles used to
construct $\bs\oM_{g,m}(X,J,\be)_{\rm sym}$ and
$\bs\oM_{g,m}(X,J,\be)_{\rm alg}$ are $\varphi(r)=e^{1/r}-e$ and
$\varphi(r)=-\frac{1}{2\pi}\ln r$, respectively. However, the author
expects that $\bs\oM_{g,m}(X,J,\be)_{\rm sym}$ and
$\bs\oM_{g,m}(X,J,\be)_{\rm alg}$ are naturally bordant, so they
define the same Gromov--Witten d-orbifold bordism invariants in
{\bf(a)}, and the same Gromov--Witten invariants.\I{Gromov--Witten
invariants|)}
\smallskip

\noindent{\bf(c)} Some areas of symplectic geometry --- Lagrangian
Floer cohomology,\I{Lagrangian Floer cohomology} Fukaya
categories,\I{Fukaya categories} contact homology,\I{contact
homology} and Symplectic Field Theory\I{Symplectic Field Theory} ---
involve compact moduli spaces $\oM$ of $J$-holomorphic curves $\Si$
which either have boundary in a Lagrangian submanifold,\I{Lagrangian
submanifold} or have ends asymptotic to cones on a Reeb orbit. Then
$\oM$ naturally becomes a d-orbifold {\it with corners\/} $\bs\oM$,
as in Theorem \ref{ds16thm2}(iii) above, where the boundary
$\pd\bs\oM$ parametrizes real codimension one bubbling behaviour of
the curves $\Si$. This situation does not arise in (complex)
algebraic geometry, where real codimension one singular behaviour
does not occur.

For a moduli space $\bs\oM$ with corners, we cannot define a virtual
class $[\bs\oM]_{\rm virt}$ in some homology group $H_*(X)$.
Instead, we can only construct a (nonunique, and depending on
choices) {\it virtual chain\/}\I{virtual chain} $\bs\oM_{\rm virt}$
in the chains $\bigl(C_*(X),\pd\bigr)$ of a homology theory
$H_*(X)$, where $\pd\bigl(\bs\oM_{\rm virt}\bigr)\in C_{*-1}(X)$ is
a virtual chain for $\pd\bs\oM$. In applications (as in the
Lagrangian Floer cohomology of Fukaya et al.\ \cite{FOOO}, for
instance), one wants to choose $\bs\oM_{\rm virt}$ so that
$\pd\bigl(\bs\oM_{\rm virt}\bigr)$ is compatible with choices of
virtual chains for other moduli spaces. This leads to difficult
technical issues about virtual chains, and some ugly homological
algebra.

The author believes that these areas of symplectic geometry
involving `moduli spaces with corners' can be made simpler, more
elegant, and less technical, by writing them in terms of d-orbifolds
with corners. There are two parts to this. Firstly, by developing a
well-behaved 2-category $\dOrbc$ of d-orbifolds with corners, we
have made it easy to state precise relationships between boundaries
of moduli spaces and other moduli spaces. For example, for the
moduli spaces $\bs\oM_k(X,Y,J,\be)$ of Theorem \ref{ds16thm2}(iii),
the Kuranishi space boundary formula of Fukaya et al.\
\cite[Prop.~8.3.3]{FOOO} translates to an equivalence in $\dOrbc$:
\begin{equation*}
\pd\bs\oM_k(X,Y,J,\be)\simeq \coprod_{i+j=k,\;\>\be=\ga+\de\!\!\!\!
\!\!\!\!\!\!\!\!\!\!\!\!\!\!\!\!\!\!\!\!\!\!}
\bs\oM_{i+1}(X,Y,J,\ga)\t_{\bs\ev_{i+1},\bcY,\bs\ev_{j+1}}
\bs\oM_{j+1}(X,Y,J,\de),
\end{equation*}
where the fibre products exist in $\dOrbc$ by
Theorem~\ref{ds14thm7}(a).

Secondly, in future work the author intends to define a virtual
chain\I{virtual chain} construction for d-manifolds and d-orbifolds
with corners, expressed in terms of new (co)homology theories whose
(co)chains are built from d-manifolds or d-orbifolds with corners,
as for the `Kuranishi (co)homology'\I{Kuranishi (co)homology}
described in \cite{Joyc1,Joyc2}. Defining virtual chains for moduli
spaces in this homology theory would be almost trivial, and would
not involve perturbing the moduli space. Many issues to do with
transversality would also become \I{moduli space!of J-holomorphic
curves@of $J$-holomorphic curves|)}\I{symplectic
geometry|)}\I{virtual class|)}trivial.
\label{ds16rem}
\end{rem}

Finally, we note that d-manifolds should not be confused with {\it
differential graded manifolds}, or {\it
dg-manifolds}.\I{dg-manifold|(}\I{d-manifold!and dg-manifolds} This
term is used in two senses, in algebraic geometry to mean a special
kind of dg-scheme,\I{dg-scheme} as in Ciocan-Fontanine and Kapranov
\cite[Def.~2.5.1]{CiKa}, and in differential geometry to mean a
supermanifold with extra structure, as in Cattaneo and Sch\"atz
\cite[Def.~3.6]{CaSc}. In both cases, a dg-manifold $\mathfrak E$ is
roughly the total space of a graded vector bundle $E^\bu$ over a
manifold $V$, with a vector field $Q$ of degree 1
satisfying~$[Q,Q]=0$.

For example, if $E$ is a vector bundle over $V$ and $s\in C^\iy(E)$,
we can make $E$ into a dg-manifold $\mathfrak E$ by giving $E$ the
grading $-1$, and taking $Q$ to be the vector field on $E$
corresponding to $s$. To this $\mathfrak E$ we can associate the
d-manifold $\bS_{V,E,s}$\I{d-manifold!standard model} from Example
\ref{ds4ex2}. Note that $\bS_{V,E,s}$ only knows about an
infinitesimal neighbourhood of $s^{-1}(0)$ in $V$, but $\mathfrak E$
remembers all of~$V,E,s$.\I{dg-manifold|)}

\appendix

\section{Categories and 2-categories}
\label{dsA}
\I{category|(}\I{category!abelian|see{abelian category}}

We now explain the background in category theory we need. Some good
references are Behrend et al.\ \cite[App.~B]{BEFF}, and MacLane
\cite{MacL} for~\S\ref{dsA1}--\S\ref{dsA2}.

\subsection{Basics of category theory}
\label{dsA1}

For completeness, here are the basic definitions in category theory,
as in~\cite[\S I]{MacL}.

\begin{dfn} A {\it category\/} (or 1-{\it category\/}) $\cC$
consists of a proper class of {\it objects\/} $\Obj(\cC)$, and for
all $X,Y\in\Obj(\cC)$ a set $\Hom(X,Y)$ of {\it
morphisms\/}\I{category!morphism} $f$ from $X$ to $Y$, written
$f:X\ra Y$, and for all $X,Y\in\Obj(\cC)$ a {\it composition map\/}
$\ci:\Hom(X,Y)\t\Hom(Y,Z)\ra\Hom(X,Z)$, written $(f,g)\mapsto g\ci
f$. Composition must be {\it associative}, that is, if $f:W\ra X$,
$g:X\ra Y$ and $h:Y\ra Z$ are morphisms in $\cC$ then $(h\ci g)\ci
f=h\ci(g\ci f)$. For each $X\in \Obj(\cC)$ there must exist an {\it
identity morphism\/} $\id_X:X\ra X$ such that $f\ci\id_X=f=\id_Y\ci
f$ for all $f:X\ra Y$ in~$\cC$.

A morphism $f:X\ra Y$ is an {\it isomorphism\/} if there exists
$f^{-1}:Y\ra X$ with $f^{-1}\ci f=\id_X$ and $f\ci f^{-1}=\id_Y$. A
category $\cC$ is called a {\it
groupoid\/}\I{groupoid}\I{category!groupoid} if every morphism is an
isomorphism. In a (small) groupoid $\cC$, for each $X\in\Obj(\cC)$
the set $\Hom(X,X)$ of morphisms $f:X\ra X$ form a group.

If $\cC$ is a category, the {\it opposite
category\/}\I{category!opposite} ${\cC}^{\rm op}$ is $\cC$ with the
directions of all morphisms reversed. That is, we define
$\Obj(\cC^{\rm op})=\Obj(\cC)$, and for all $X,Y,Z\in\Obj(\cC)$ we
define $\Hom_{\cC^{\rm op}}(X,Y)=\Hom_\cC(Y,X)$, and for $f:X\ra Y$,
$g:Y\ra Z$ in $\cC$ we define $f\ci_{\cC^{\rm op}}g=g\ci_\cC f$, and
$\id_{\cC^{\rm op}}X=\id_\cC X$.

Given categories $\cC,\cD$, the {\it product category\/} ${\cal
C}\t\cD$ has objects $(W,X)$ in $\Obj(\cC)\t\Obj(\cD)$ and morphisms
$f\t g:(W,X)\ra (Y,Z)$ when $f:W\ra Y$ is a morphism in $\cC$ and
$g:X\ra Z$ is a morphism in $\cD$, in the obvious way.

We call $\cD$ a {\it subcategory\/}\I{category!subcategory} of $\cC$
if $\Obj(\cD)\subseteq \Obj(\cC)$, and
$\Hom_\cD(X,Y)\subseteq\Hom_\cC(X,Y)$ for all $X,Y\in\Obj(\cD)$. We
call $\cD$ a {\it full\/}\I{category!subcategory!full} subcategory
if also $\Hom_\cD(X,Y)=\Hom_\cC(X,Y)$ for all~$X,Y$.
\label{dsAdef1}
\end{dfn}

\begin{dfn} Let $\cC,\cD$ be categories. A ({\it
covariant\/}) {\it functor}\I{functor}\I{category!functor}
$F:\cC\ra\cD$ gives, for all objects $X$ in $\cC$ an object $F(X)$
in $\cD$, and for all morphisms $f:X\ra Y$ in $\cC$ a morphism
$F(f):F(X)\ra F(Y)$, such that $F(g\ci f)=F(g)\ci F(f)$ for all
$f:X\ra Y$, $g:Y\ra Z$ in $\cC$, and $F(\id_X)=\id_{F(X)}$ for all
$X\in\Obj(\cC)$. A {\it contravariant functor\/} $F:\cC\ra\cD$ is a
covariant functor $F:\cC^{\rm op}\ra\cD$.

Functors compose in the obvious way. Each category $\cC$ has an
obvious {\it identity functor\/} $\id_\cC:\cC\ra\cC$ with $\id_{\cal
C}(X)=X$ and $\id_\cC(f)=f$ for all $X,f$. A functor $F:\cC\ra\cD$
is called {\it full\/}\I{functor!full}\I{category!functor!full} if
the maps $\Hom_\cC(X,Y)\ra\Hom_\cD (F(X),F(Y))$, $f\mapsto F(f)$ are
surjective for all $X,Y\in\Obj(\cC)$, and {\it
faithful\/}\I{functor!faithful}\I{category!functor!faithful} if the
maps $\Hom_\cC(X,Y)\ra\Hom_\cD(F(X),F(Y))$ are injective for
all~$X,Y\in\Obj(\cC)$.

Let $\cC,\cD$ be categories and $F,G:\cC\ra\cD$ be functors. A {\it
natural transformation\/}\I{functor!natural
transformation}\I{category!functor!natural transformation}
$\eta:F\Ra G$ gives, for all objects $X$ in $\cC$, a morphism
$\eta(X):F(X)\ra G(X)$ such that if $f:X\ra Y$ is a morphism in
$\cC$ then $\eta(Y)\ci F(f)=G(f)\ci\eta(X)$ as a morphism $F(X)\ra
G(Y)$ in $\cD$. We call $\eta$ a {\it natural
isomorphism\/}\I{functor!natural
isomorphism}\I{category!functor!natural isomorphism} if $\eta(X)$ is
an isomorphism for all~$X\in\Obj(\cC)$.

An {\it equivalence\/}\I{category!equivalence of} between categories
$\cC,\cD$ consists of functors $F:\cC\ra\cD$ and $G:\cD\ra\cC$ with
natural isomorphisms $\eta:G\ci F\Ra\id_\cC$, $\ze:F\ci
G\Ra\id_\cD$.
\label{dsAdef2}
\end{dfn}

It is a fundamental principle of category theory that equivalent
categories $\cC,\cD$ should be thought of as being `the same', and
naturally isomorphic functors $F,G:\cC\ra\cD$ should be thought of
as being `the same'. Note that equivalence of categories $\cC,{\cal
D}$ is much weaker than strict isomorphism: isomorphism classes of
objects in $\cC$ are naturally in bijection with isomorphism classes
of objects in $\cD$, but there is no relation between the sizes of
the isomorphism classes, so that $\cC$ could have many more objects
than $\cD$, for instance.

\subsection{Limits, colimits and fibre products in categories}
\label{dsA2}
\I{category!limit|(}\I{category!colimit|(}

We shall be interested in various kinds of {\it limits\/} and {\it
colimits\/} in our categories of spaces. These are objects in the
category with a universal property\I{category!universal property}
with respect to some class of diagrams. For an introduction to
limits and colimits in category theory, see MacLane \cite[\S
III]{MacL}. Here are the basic definitions:

\begin{dfn} Let $\cC$ be a category. A {\it diagram\/} $\De$ in $\cal
C$ is a class of objects $S_i$ in $\cC$ for $i\in I$, and a class of
morphisms $\rho_j:S_{b(j)}\ra S_{e(j)}$ in $\cC$ for $j\in J$, where
$b,e:J\ra I$. The diagram is called {\it finite\/} if $I,J$ are
finite sets.

A {\it limit\/} of the diagram $\De$ is an object $L$ in $\cC$ and
morphisms $\pi_i:L\ra S_i$ for $i\in I$ such that
$\rho_j\ci\pi_{b(j)}=\pi_{e(j)}$ for all $j\in J$, which has the
universal property that given $L'\in\cC$ and $\pi_i':L'\ra S_i$ for
$i\in I$ with $\rho_j\ci\pi_{b(j)}'=\pi_{e(j)}'$ for all $j\in J$,
there is a unique morphism $\la:L'\ra L$ with $\pi_i'=\pi_i\ci \la$
for all~$i\in I$.

A {\it fibre product\/}\I{category!fibre product}\I{fibre
product!definition} is a limit of a diagram $X\,{\buildrel
g\over\longra}\,Z\,{\buildrel h\over\longleftarrow}\,Y$. The limit
object $W$ is often written $X\t_{g,Z,h}Y$ or $X\t_ZY$, and the
diagram
\begin{equation*}
\xymatrix@R=10pt@C=40pt{ W \ar[r]_{\pi_Y} \ar[d]^{\pi_X} & Y
\ar[d]_h \\ X \ar[r]^g & Z}
\end{equation*}
is called a {\it Cartesian square\/}\I{Cartesian
square}\I{category!Cartesian square} in the category $\cC$. A {\it
terminal object\/}\I{category!terminal object} is a limit of the
empty diagram.

A {\it colimit\/} of the diagram $\De$ is an object $L$ in $\cC$ and
morphisms $\la_i:S_i\ra L$ for $i\in I$ such that
$\la_{b(j)}=\la_{e(j)}\ci\rho_j$ for all $j\in J$, which has the
universal property that given $L'\in\cC$ and $\la_i':S_i\ra L'$ for
$i\in I$ with $\la_{b(j)}'=\la_{e(j)}'\ci\rho_j$ for all $j\in J$,
there is a unique morphism $\pi:L\ra L'$ with $\la_i'=\pi\ci\la_i$
for all~$i\in I$.

A {\it pushout\/}\I{pushout}\I{category!pushout} is a colimit of a
diagram $X\,{\buildrel e\over\longleftarrow}\,W\,{\buildrel
f\over\longra}\,Y$.

If a limit or colimit exists, it is unique up to unique isomorphism
in $\cC$. We say that {\it all finite limits}, or {\it all fibre
products exist in\/} $\cC$, if limits exist for all finite diagrams,
or for all diagrams $X\,{\buildrel g\over\longra}\,Z\,{\buildrel
h\over\longleftarrow}\,Y$
respectively.\I{category|)}\I{category!limit|)}\I{category!colimit|)}
\label{dsAdef3}
\end{dfn}

\subsection{2-categories}
\label{dsA3}
\I{2-category|(}\I{category!2-category|see{2-category}}

Next we discuss 2-{\it categories}. A good reference for our
purposes is Behrend et al.\ \cite[App.~B]{BEFF}, and Kelly and
Street \cite{KeSt} is also helpful.

\begin{dfn} A 2-{\it category\/} $\bs{\cal C}$ (also called a {\it strict\/
$2$-category\/})\I{2-category!strict} consists of a proper class of
{\it objects\/} $\Obj(\bs{\cal C})$, for all $X,Y\in\Obj(\bs{\cal
C})$ a category $\Hom(X,Y)$, for all $X\in\Obj(\bs{\cal C})$ an
object $\id_X$ in $\Hom(X,X)$ called the {\it identity
$1$-morphism}, and for all $X,Y,Z\in\Obj(\bs{\cal C})$ a functor
$\mu_{X,Y,Z}:\Hom(X,Y)\t\Hom(Y,Z)\ra\Hom(X,Z)$. These must satisfy
the {\it identity property}, that $\mu_{X,X,Y}(\id_X,-)\!=\!
\mu_{X,Y,Y}(-,\id_Y)\!=\!\id_{\Hom(X,Y)}$ as functors
$\Hom(X,Y)\ra\Hom(X,Y)$, and the {\it associativity property}, that
$\mu_{W,Y,Z}\ci(\mu_{W,X,Y}\t\id_{\Hom(Y,Z)})
=\mu_{W,X,Z}\ab\ci\ab(\id_{\Hom(W,X)}\ab\t\mu_{X,Y,Z})$ as functors
$\Hom(W,X)\t\Hom(X,Y)\t\Hom(Y,Z)\ra\Hom(W,X)$, for all~$W,X,Y,Z$.

Objects $f$ of $\Hom(X,Y)$ are called 1-{\it
morphisms},\I{2-category!1-morphism} written $f:X\ra Y$. For
1-morphisms $f,g:X\ra Y$, morphisms $\eta\in \Hom_{\Hom(X,Y)}(f,g)$
are called 2-{\it morphisms},\I{2-category!2-morphism} written
$\eta:f\Ra g$. Thus, a 2-category has objects $X$, and two kinds of
morphisms, 1-morphisms $f:X\ra Y$ between objects, and 2-morphisms
$\eta:f\Ra g$ between 1-morphisms. In many examples, all 2-morphisms
are 2-isomorphisms (i.e.\ have an inverse), so that the categories
$\Hom(X,Y)$ are groupoids.\I{category!groupoid} Such 2-categories
are called (2,1)-{\it categories}.
\label{dsAdef4}
\end{dfn}

This is quite a complicated structure. There are three kinds of
composition in a 2-category, satisfying various associativity
relations. If $f:X\ra Y$ and $g:Y\ra Z$ are 1-morphisms then
$\mu_{X,Y,Z}(f,g)$ is the {\it composition of\/
$1$-morphisms},\I{2-category!1-morphism!composition} written $g\ci
f:X\ra Z$. If $f,g,h:X\ra Y$ are 1-morphisms and $\eta:f\Ra g$,
$\ze:g\Ra h$ are 2-morphisms then composition of $\eta,\ze$ in the
category $\Hom(X,Y)$ gives the {\it vertical composition of\/
$2$-morphisms}\I{2-category!2-morphism!vertical composition} of
$\eta,\ze$, written $\ze\od\eta:f\Ra h$, as a diagram
\e
\begin{gathered}
\xymatrix@C=25pt{ X \rruppertwocell^f{\eta} \rrlowertwocell_h{\ze}
\ar[rr]_(0.35)g && Y & \ar@{~>}[r] && X
\rrtwocell^f_h{{}\,\,\,\,\ze\od\eta\!\!\!\!\!} && Y.}
\end{gathered}
\label{dsAeq1}
\e
And if $f,\ti f:X\ra Y$ and $g,\ti g:Y\ra Z$ are 1-morphisms and
$\eta:f\Ra\ti f$, $\ze:g\Ra\ti g$ are 2-morphisms then
$\mu_{X,Y,Z}(\eta,\ze)$ is the {\it horizontal composition of\/
$2$-morphisms},\I{2-category!2-morphism!horizontal composition}
written $\ze*\eta:g\ci f\Ra\ti g\ci\ti f$, as a diagram
\e
\begin{gathered}
\xymatrix@C=20pt{ X \rrtwocell^f_{\ti f}{\eta} && Y
\rrtwocell^g_{\ti g}{\ze} && Z & \ar@{~>}[r] && X \rrtwocell^{g\ci
f}_{\ti g\ci\ti f}{{}\,\,\,\ze*\eta\!\!\!\!\!} && Z. }
\end{gathered}
\label{dsAeq2}
\e
There are also two kinds of identity: {\it identity\/
$1$-morphisms\/} $\id_X:X\ra X$ and {\it identity\/
$2$-morphisms\/}~$\id_f:f\Ra f$.

A basic example is the {\it $2$-category of categories\/}
$\mathfrak{Cat}$, with objects categories $\cC$, 1-morphisms
functors $F:\cC\ra\cD$, and 2-morphisms natural transformations
$\eta:F\Ra G$ for functors $F,G:\cC\ra\cD$. Orbifolds naturally form
a 2-category, as do Deligne--Mumford stacks\I{stack} and Artin
stacks in algebraic geometry.

In a 2-category $\bs{\cal C}$, there are three notions of when
objects $X,Y$ in $\bs{\cal C}$ are `the same': {\it equality\/}
$X=Y$, and {\it isomorphism}, that is we have 1-morphisms $f:X\ra
Y$, $g:Y\ra X$ with $g\ci f=\id_X$ and $f\ci g=\id_Y$, and {\it
equivalence},\I{2-category!equivalence in} that is we have
1-morphisms $f:X\ra Y$, $g:Y\ra X$ and 2-isomorphisms $\eta:g\ci
f\Ra\id_X$ and $\ze:f\ci g\Ra\id_Y$. Usually equivalence is the most
useful.

Let $\bs{\cal C}$ be a 2-category. The {\it homotopy
category\/}\I{homotopy category} $\Ho(\bs{\cal C})$ of $\bs{\cal C}$
is the category whose objects are objects of $\bs{\cal C}$, and
whose morphisms $[f]:X\ra Y$ are 2-isomorphism classes $[f]$ of
1-morphisms $f:X\ra Y$ in $\bs{\cal C}$. Then equivalences in
$\bs{\cal C}$ become isomorphisms in $\Ho(\bs{\cal C})$,
2-commutative diagrams in $\bs{\cal C}$ become commutative diagrams
in $\Ho(\bs{\cal C})$, and so on.

As in Borceux \cite[\S 7.7]{Borc}, there is also a second kind of
2-category, called a {\it weak\/ $2$-category\/}\I{2-category!weak}
(or {\it bicategory\/}), which we will not define in detail. In a
weak 2-category, compositions of 1-morphisms need only be
associative up to (specified) 2-isomorphisms. That is, part of the
data of a weak 2-category $\bs{\cal C}$ is a 2-isomorphism
$\al(f,g,h):(h\ci g)\ci f\Ra h\ci(g\ci f)$ for all 1-morphisms
$f:W\ra X$, $g:X\ra Y$, $h:Y\ra Z$ in $\bs{\cal C}$. A strict
2-category $\bs{\cal C}$ can be made into a weak 2-category by
putting $\al(f,g,h)=\id_{h\ci g\ci f}$ for all~$f,g,h$.

Some categorical constructions naturally yield weak 2-categories
rather than strict 2-categories, e.g.\ the weak 2-categories of
orbifolds defined by Pronk \cite{Pron} and Lerman \cite[\S
3.3]{Lerm} mentioned in \S\ref{ds91}. Every weak 2-category is
equivalent as a weak 2-category to a strict 2-category (that is,
weak 2-categories can be strictified), so we lose little by working
only with strict 2-categories.

\subsection{Fibre products in 2-categories}
\label{dsA4}
\I{2-category!fibre products in|(}

{\it Commutative diagrams\/} in 2-categories should in general only
commute {\it up to (specified)\/ $2$-isomorphisms}, rather than
strictly. Then we say the diagram 2-{\it commutes}. A simple example
of a commutative diagram in a 2-category $\bs{\cal C}$
is\I{2-category!2-commutative diagram}
\begin{equation*}
\xymatrix@C=50pt@R=8pt{ & Y \ar[dr]^g \ar@{=>}[d]^\eta \\
X \ar[ur]^f \ar[rr]_h && Z, }
\end{equation*}
which means that $X,Y,Z$ are objects of $\bs{\cal C}$, $f:X\ra Y$,
$g:Y\ra Z$ and $h:X\ra Z$ are 1-morphisms in $\bs{\cal C}$, and
$\eta:g\ci f\Ra h$ is a 2-isomorphism.

The generalizations of {\it limit\/} and {\it colimit\/} to
2-categories turn out to be rather
complicated.\I{2-category!limit}\I{2-category!colimit} As in
\cite[\S 7]{Borc} there are many different kinds --- 2-limits,
bilimits, pseudolimits, lax limits, and weighted limits (or indexed
limits), depending on whether one considers diagrams to commute on
the nose, up to 2-isomorphism, or up to 2-morphisms, and what kind
of universal property one requires. Our definition of {\it fibre
product}, following Behrend et al.\ \cite[Def.~B.13]{BEFF}, is
actually an example of a pseudolimit.

\begin{dfn} Let $\bs{\cal C}$ be a 2-category and $g:X\ra Z$,
$h:Y\ra Z$ be 1-morphisms in $\bs{\cal C}$. A {\it fibre product\/}
$X\t_ZY$ in $\bs{\cal C}$ consists of an object $W$, 1-morphisms
$e:W\ra X$ and $f:W\ra Y$ and a 2-isomorphism $\eta:g\ci e\Ra h\ci
f$ in $\bs{\cal C}$, so that we have a 2-commutative diagram
\e
\begin{gathered}
\xymatrix@C=80pt@R=10pt{ W \ar[r]_(0.25){f} \ar[d]^{e}
\drtwocell_{}\omit^{}\omit{^{\eta}}
 & Y \ar[d]_{h} \\ X \ar[r]^(0.7){g} & Z}
\end{gathered}
\label{dsAeq3}
\e
with the following universal property: suppose $e':W'\ra X$ and
$f':W'\ra Y$ are 1-morphisms and $\eta':g\ci e'\Ra h\ci f'$ is a
2-isomorphism in $\bs{\cal C}$. Then there should exist a 1-morphism
$b:W'\ra W$ and 2-isomorphisms $\ze:e\ci b\Ra e'$, $\th:f\ci b\Ra
f'$ such that the following diagram of 2-isomorphisms commutes:
\begin{equation*}
\xymatrix@C=50pt@R=11pt{ g\ci e\ci b \ar@{=>}[r]_{\eta*\id_b}
\ar@{=>}[d]_{\id_g*\ze} & h\ci f\ci b \ar@{=>}[d]^{{}\,\id_h*\th}
\\ g\ci e' \ar@{=>}[r]^{\eta'} & h\ci f'.}
\end{equation*}
Furthermore, if $\ti b,\ti\ze,\ti\th$ are alternative choices of
$b,\ze,\th$ then there should exist a unique 2-isomorphism $\ep:\ti
b\Ra b$ with
\begin{equation*}
\ti\ze=\ze\od(\id_{e}*\ep)\quad\text{and}\quad
\ti\th=\th\od(\id_{f}*\ep).
\end{equation*}
We call such a fibre product diagram \eq{dsAeq3} a 2-{\it Cartesian
square}.\I{2-category!2-Cartesian square}

If a fibre product $X\t_ZY$ in $\bs{\cal C}$ exists then it is
unique up to equivalence in $\bs{\cal C}$. If $\bs{\cal C}$ is a
category, that is, all 2-morphisms are identities $\id_f:f\Ra f$,
this definition of fibre products in $\bs{\cal C}$ coincides with
that in~\S\ref{dsA2}.
\label{dsAdef5}
\end{dfn}

Orbifolds,\I{orbifold!transverse fibre products} and stacks\I{stack}
in algebraic geometry, form 2-categories, and Definition
\ref{dsAdef5} is the right way to define fibre products of orbifolds
or stacks, as in \cite{BEFF}.\I{2-category!fibre products
in|)}\I{2-category|)} One can also define {\it pushouts\/} in a
2-category $\bs{\cal C}$\I{2-category!pushout}\I{pushout} in a dual
way to Definition \ref{dsAdef5}, reversing directions of morphisms.

\clearpage
\printnomenclature[1.3cm]
\clearpage
\addcontentsline{toc}{section}{Index}
\printindex

\end{document}